\documentclass[english]{article}
\usepackage[T1]{fontenc}
\usepackage[latin9]{inputenc}
\usepackage{geometry}
\usepackage{babel}
\usepackage{verbatim}
\usepackage{float}
\usepackage{mathrsfs}
\usepackage{mathtools}
\usepackage{bm}
\usepackage{amsmath}
\usepackage{amssymb}
\usepackage{graphicx}
\usepackage[unicode=true,
 bookmarks=false,
 breaklinks=false,pdfborder={0 0 1},colorlinks=false]
 {hyperref}
\hypersetup{
 colorlinks,linkcolor=red,anchorcolor=blue,citecolor=blue}

\makeatletter

\providecommand{\tabularnewline}{\\}
\floatstyle{ruled}
\newfloat{algorithm}{tbp}{loa}
\providecommand{\algorithmname}{Algorithm}
\floatname{algorithm}{\protect\algorithmname}

\usepackage{babel}
\usepackage{babel}
\usepackage{babel}

\usepackage{mathtools}

\usepackage{cite}\usepackage{amsthm}\usepackage{dsfont}\usepackage{array}\usepackage{mathrsfs}\usepackage{comment}\onecolumn
\usepackage{natbib}
\usepackage{color}\usepackage{babel}

\allowdisplaybreaks

\usepackage{enumitem}
\setlist[itemize]{leftmargin=1.5em}
\setlist[enumerate]{leftmargin=1.5em}

\usepackage{babel}
\usepackage{algorithm}
\usepackage{algorithmic}
\usepackage{arydshln}

\DeclareMathOperator{\ind}{\mathds{1}}  

\numberwithin{equation}{section}

\definecolor{yxc}{RGB}{255,0,0}
\definecolor{yjc}{RGB}{125,0,0}
\definecolor{cm}{RGB}{0,0,200}
\definecolor{yly}{RGB}{0,150,0}

\newcommand{\yxc}[1]{\textcolor{yxc}{[YXC: #1]}}

\makeatother

\begin{document}
\theoremstyle{plain} \newtheorem{lemma}{\textbf{Lemma}} \newtheorem{prop}{\textbf{Proposition}}\newtheorem{theorem}{\textbf{Theorem}}\setcounter{theorem}{0}
\newtheorem{corollary}{\textbf{Corollary}} \newtheorem{assumption}{\textbf{Assumption}}
\newtheorem{example}{\textbf{Example}} \newtheorem{definition}{\textbf{Definition}}
\newtheorem{fact}{\textbf{Fact}} \newtheorem{condition}{\textbf{Condition}}\theoremstyle{definition}

\theoremstyle{remark}\newtheorem{remark}{\textbf{Remark}}\newtheorem{claim}{\textbf{Claim}}\newtheorem{conjecture}{\textbf{Conjecture}}
\title{Inference for  Heteroskedastic PCA with Missing Data}
\author{Yuling Yan\thanks{Institute for Data, Systems, and Society, MIT, Cambridge, MA 02142, USA; Email: \texttt{yulingy@mit.edu}.}\and Yuxin Chen\thanks{Department of Statistics and Data Science, Wharton School, University
of Pennsylvania, Philadelphia, PA, 19104, USA; Email: \texttt{yuxinc@wharton.upenn.edu}.} \and Jianqing Fan\thanks{Department of Operations Research and Financial Engineering, Princeton
University, Princeton, NJ 08544, USA; Email: \texttt{jqfan@princeton.edu}.}}

\date{July 2021; ~Revised: February 2024}
\maketitle
\begin{abstract}
This paper studies how to construct confidence regions for principal
component analysis (PCA) in high dimension, a problem that has been
vastly under-explored. While computing measures of uncertainty for
nonlinear/nonconvex estimators is in general difficult in high dimension,
the challenge is further compounded by the prevalent presence of missing
data and heteroskedastic noise. We propose a novel approach to
	performing valid inference on the principal subspace under a spiked covariance model with missing data, 
on the basis of an estimator called \textsf{HeteroPCA} \citep{zhang2018heteroskedastic}.
We develop non-asymptotic distributional guarantees for \textsf{HeteroPCA},
and demonstrate how these can be invoked to compute both confidence
regions for the principal subspace and entrywise confidence intervals
for the spiked covariance matrix. Our inference procedures are fully
data-driven and adaptive to heteroskedastic random noise, without
requiring prior knowledge about the noise levels. 
\end{abstract}

\noindent \textbf{Keywords: }principal component analysis, confidence
regions, missing data, uncertainty quantification, heteroskedastic
data, subspace estimation 

\def\r#1{\textcolor{red}{\bf #1}}
\def\b#1{\textcolor{blue}{\bf #1}}

\setcounter{tocdepth}{2}

\tableofcontents{}

\section{Introduction}

The applications of modern data science frequently ask for succinct
representations of high-dimensional data. At the core of this pursuit
lies principal component analysis (PCA), which serves as an effective
means of dimension reduction and has been deployed across a broad
range of domains \citep{jolliffe1986principal,johnstone2018pca,vaswani2018rethinking,fan2021robust}.
In reality, data collection could often be far from ideal --- for
instance, the acquired data might be subject to random contamination
and contain incomplete observations --- which inevitably affects
the fidelity of PCA and calls for additional care when interpreting
the results. To enable informative assessment of the influence of
imperfect data acquisition, it would be desirable to accompany the
PCA estimators in use with valid measures of uncertainty or ``confidence''.

\subsection{Problem formulation \label{subsec:intro-model}}

To allow for concrete and precise studies, the present paper concentrates
on a tractable model that captures the effects of random heteroskedastic
noise and missing data in PCA. In what follows, we start by formulating
the problem, in the hope of facilitating more precise discussions.

\paragraph{Model.}

Imagine we are interested in $n$ independent random vectors $\bm{x}_{j}=[x_{1,j},\cdots,x_{d,j}]^{\top}\in\mathbb{R}^{d}$
drawn from the following distribution\footnote{	All results in this paper continue to hold 
	if the sample vectors are generated such that $\bm{x}_j=\bm{U}^\star (\bm{\Lambda}^\star)^{1/2}\bm{f}_j$, where $\{\bm{f}_j\}_{j=1}^n$ are independent {\em sub-Gaussian} random vectors in $\mathbb{R}^r$ satisfying $\bm{E}[\bm{f}_j]=\bm{0}$, $\mathbb{E}[\bm{f}_j\bm{f}_j^\top]=\bm{I}_r$ and $\Vert\bm{f}_j\Vert_{\psi_2}=O(1)$. Here, $\Vert\cdot\Vert_{\psi_2}$ denotes the sub-Gaussian norm 
	 \citep{vershynin2016high}.} 
\begin{equation}
\bm{x}_{j}\overset{\text{ind.}}{\sim}\mathcal{N}\left(\bm{0},\bm{S}^{\star}\right),\qquad1\leq j\leq n,\label{eq:xi-data-samples-noiseless}
\end{equation}
where the unknown covariance matrix $\bm{S}^{\star}\in\mathbb{R}^{d\times d}$
is assumed to be rank-$r$ $(r<n)$ with eigen-decomposition 
\begin{equation}
\bm{S}^{\star}=\bm{U}^{\star}\bm{\Lambda}^{\star}\bm{U}^{\star\top}.\label{eq:covariance-spike}
\end{equation}
Here, the orthonormal columns of $\bm{U}^{\star}\in\mathbb{R}^{d\times r}$
constitute the $r$ leading eigenvectors of $\bm{S}^{\star}$, whereas
$\bm{\Lambda}^{\star}\in\mathbb{R}^{r\times r}$ is a diagonal matrix
whose diagonal entries are composed of the non-zero eigenvalues of
$\bm{S}^{\star}$. In other words, these vectors $\{\bm{x}_{j}\}_{1\leq j\leq n}$
are randomly drawn from a low-dimensional subspace when $r$ is small.
What we have available are partial and randomly corrupted observations
of the entries of the above vectors. Specifically, suppose that we
only get to observe 
\begin{equation}
y_{l,j}=x_{l,j}+\eta_{l,j}\qquad\text{for}\text{ all }\left(l,j\right)\in\Omega\label{eq:observed-data}
\end{equation}
over a subsampled index set $\Omega\subseteq[d]\times[n]$ (with $[n]\coloneqq\{1,\cdots,n\}$),
where $\eta_{l,j}$ represents the noise that contaminates the observation
in this location. Throughout this paper, we focus on the following random sampling
and random noise models.  
\begin{itemize}
\item \emph{Random sampling}: each index $\left(l,j\right)$ is contained
in $\Omega$ independently with probability $p$; 
\item \emph{Heteroskedastic random noise with unknown variance}: the noise
components $\{\eta_{l,j}\}$ are independently generated sub-Gaussian
random variables obeying 
\[
\mathbb{E}[\eta_{l,j}]=0,\qquad\mathbb{E}[\eta_{l,j}^{2}]=\omega_{l}^{\star2},\qquad\text{and}\qquad\|\eta_{l,j}\|_{\psi_{2}}=O(\omega_{l}^{\star}),
\]
where $\{\omega_{l}^{\star}\}_{1\leq l\leq d}$ denote the standard
deviations that are \emph{a priori} unknown, and $\|\cdot\|_{\psi_{2}}$
stands for the sub-Gaussian norm of a random variable \citep{vershynin2016high}.
The noise levels $\{\omega_{l}^{\star}\}_{1\leq l\leq d}$ are allowed
to vary across locations, so as to model the so-called\emph{ heteroskedasticity}
of noise. 
\end{itemize}
This model can be viewed as a generalization of the spiked covariance model \citep{johnstone2001distribution,baik2005phase,paul2007asymptotics,donoho2018optimal,nadler2008finite,cai2019subspace,bao2022statistical} to account for missing data and heteroskedastic noise. With the observed data $\{y_{l,j}\mid(l,j)\in\Omega\}$ in hand, can
we perform statistical inference on the orthonormal matrix $\bm{U}^{\star}$
--- which embodies the ground-truth $r$-dimensional principal subspace
underlying the vectors $\{\bm{x}_{j}\}_{1\leq j\leq n}$ --- and
make inference on the underlying covariance matrix $\bm{S}^{\star}$.
Mathematically, the task can often be phrased as constructing valid
confidence intervals/regions for both $\bm{U}^{\star}$ and $\bm{S}^{\star}$
based on the incomplete and corrupted observations $\{y_{l,j}\mid(l,j)\in\Omega\}$.
Noteworthily, this model is frequently studied in econometrics and
financial modeling under the name of factor models \citep{fan2017elements,fan2021recent, fan2021robust, bai2016econometric,gagliardini2019estimation},
and is closely related to the noisy matrix completion problem where we also quantify uncertainty of missing entries  \citep{ExactMC09,CanPla10,Se2010Noisy,chi2018nonconvex}.



\paragraph{Inadequacy of prior works.}

While methods for estimating principal subspace are certainly not
in shortage (e.g., \citet{lounici2014high,zhang2018heteroskedastic,cai2018rate,balzano2018streaming,zhu2019high,cai2019subspace,li2021minimax}),
methods for constructing confidence regions for principal subspace
remain vastly under-explored. The fact that the estimators in use
for PCA are typically nonlinear and nonconvex presents a substantial
challenge in the development of a distributional theory, let alone
uncertainty quantification. As some representative recent attempts,
\citet{xia2019normal,bao2018singular} established normal approximations
of the distance between the true subspace and its estimate for the
matrix denoising setting, while \citet{koltchinskii2020efficient}
further established asymptotic normality of some debiased estimator
for linear functions of principal components. These distributional
guarantees pave the way for the development of statistical inference
procedures for PCA. However, it is noteworthy that these results required
the noise components to either be i.i.d.~Gaussian or at least exhibit
matching moments (up to the 4th order), which fell short of accommodating
heteroskedastic noise. The challenge is further compounded when statistical
inference needs to be conducted in the face of missing data, a scenario
that is beyond the reach of these prior works.

\subsection{Our contributions}

In light of the insufficiency of prior results, this paper takes a
step towards developing data-driven inference and uncertainty quantification
procedures for PCA, in the hope of accommodating both heteroskedastic
noise and missing data. Our inference procedures are built on an
estimator called \textsf{HeteroPCA} recently proposed by \citet{zhang2018heteroskedastic}, which is an iterative algorithm in nature and will be detailed in Section~\ref{sec:Background:-two-estimation}.
The main contributions of this paper are summarized as follows. 
\begin{itemize}
\item \textit{Distributional theory for PCA and covariance estimation.}
We derive, in a non-asymptotic manner, row-wise distributional characterizations
of the principal subspace estimate returned by \textsf{HeteroPCA} (see Theorem~\ref{thm:pca}),  
as well as entrywise distributional guarantees of the estimate for the covariance
matrix estimate of $\{\bm{x}_{l}\}_{1\leq l\leq n}$ (see Theorem~\ref{thm:ce}). These distributional
characterizations take the form of tractable Gaussian approximations
centered at the ground truth. 

\item \textit{Fine-grained confidence regions and intervals.} Our distributional 
theory in turns allows for construction of row-wise confidence region
for the subspace $\bm{U}^{\star}$ (see Algorithm~\ref{alg:PCA-HeteroPCA-CR} and Theorem~\ref{thm:pca-cr}) as well
as entrywise confidence intervals for the matrix $\bm{S}^{\star}$ (see Algorithm~\ref{alg:CE-HeteroPCA-CI} and Theorem~\ref{thm:ce-CI}).
The proposed inference procedures are fully data-driven and do not
require prior knowledge of the noise levels. 
\end{itemize}
Along the way, we have significantly strengthened the estimation guarantees for \textsf{HeteroPCA} in the presence of missing data. 
It is noteworthy that all of our theory allows the observed data to be highly incomplete
and covers heteroskedastic noise, which is previously unavailable.

%
%

\subsection{Paper organization}

The remainder of the paper is organized as follows. In Section~\ref{sec:Background:-two-estimation},
we introduce the estimation algorithms available in prior literature.
Section~\ref{sec:Distributional-theory-inference} develops a suite
of distributional theory for \textsf{HeteroPCA} and demonstrates how
to use it to construct fine-grained confidence regions and confidence
intervals for the unknowns; the detailed proofs of our theorems are
deferred to the appendices. In Section~\ref{sec:Numerical-experiments},
we carry out a series of numerical experiments to confirm the validity
and applicability of our theoretical findings. Section~\ref{sec:Related-works}
gives an overview of several related works. Section~\ref{sec:detour-subspace}
takes a detour to analyze two intimately related problems, which will
then be utilized to establish our main results. We conclude the paper
with a discussion of future directions in Section~\ref{sec:Discussion}.
Most of the proof details are deferred to the appendices.

\subsection{Notation\label{subsec:Paper-organization-notation}}

Before proceeding, we introduce several notation that will be useful
throughout. We let $f(n)\lesssim g(n)$ or $f(n)=O(g(n))$ represent
the condition that $\vert f(n)\vert\leq Cg(n)$ for some constant
$C>0$ when $n$ is sufficiently large; we use $f(n)\gtrsim g(n)$
to denote $f(n)\geq C\vert g(n)\vert$ for some constant $C>0$ when
$n$ is sufficiently large; and we let $f(n)\asymp g(n)$ indicate
that $f(n)\lesssim g(n)$ and $f(n)\gtrsim g(n)$ hold simultaneously.
The notation $f(n)\gg g(n)$ (resp.~$f(n)\ll g(n)$) means that there
exists some sufficiently large (resp.~small) constant $c_{1}>0$
(resp.~$c_{2}>0$) such that $f(n)\geq c_{1}g(n)$ (resp.~$f(n)\leq c_{2}g(n)$). We also let $f(n)=o(g(n))$ 
indicate that $\lim_{n\to\infty} f(n)/g(n)=0$.
For any real number $a,b\in\mathbb{R}$, we shall define $a\land b\coloneqq\min\{a,b\}$
and $a\lor b\coloneqq\max\{a,b\}$.

For any matrix $\bm{M}=[M_{i,j}]_{1\leq i\leq n_{1},1\leq j\leq n_{2}}$,
we let $\bm{M}_{i,\cdot}$ and $\bm{M}_{\cdot,j}$ stand for the $i$-th
row and the $j$-th column of $\bm{M}$, respectively. We shall also
let $\|\bm{M}\|$, $\|\bm{M}\|_{\mathrm{F}}$, $\|\bm{M}\|_{2,\infty}$
and $\|\bm{M}\|_{\infty}$ denote the spectral norm, the Frobenius
norm, the $\ell_{2,\infty}$ norm (i.e., $\|\bm{M}\|_{2,\infty}\coloneqq\max_{i}\|\bm{M}_{i,\cdot}\|_{2}$),
and the entrywise $\ell_{\infty}$ norm ($\|\bm{M}\|_{\infty}\coloneqq\max_{i,j}|M_{i,j}|$)
of $\bm{M}$, respectively. For any index set $\Omega$, the notation
$\mathcal{P}_{\Omega}(\bm{M})$ represents the Euclidean projection
of a matrix $\bm{M}$ onto the subspace of matrices supported on $\Omega$,
and define $\mathcal{P}_{\Omega^{\mathrm{c}}}(\bm{M})\coloneqq\bm{M}-\mathcal{P}_{\Omega}(\bm{M})$
as well. In addition, we denote by $\mathcal{P}_{\mathsf{diag}}(\bm{G})$
the Euclidean projection of a square matrix $\bm{G}$ onto the subspace
of matrices that vanish outside the diagonal, and define $\mathcal{P}_{\mathsf{off}\text{-}\mathsf{diag}}(\bm{G})\coloneqq\bm{G}-\mathcal{P}_{\mathsf{diag}}(\bm{G})$.
For a non-singular matrix $\bm{H}\in\mathbb{R}^{k\times k}$ with
SVD $\bm{U}_{H}\bm{\Sigma}_{H}\bm{V}_{H}^{\top}$, we denote by $\mathsf{sgn}(\bm{H})$
the following orthogonal matrix 
\begin{equation}
\mathsf{sgn}(\bm{H})\coloneqq\bm{U}_{H}\bm{V}_{H}^{\top}.\label{eq:defn-sgn-H-UV}
\end{equation}

Finally, we denote by $\mathscr{C}^{d}$ the set of all convex sets
in $\mathbb{R}^{d}$. For any Lebesgue measurable set $\mathcal{A}\subseteq\mathbb{R}^{d}$,
we adopt the shorthand notation $\mathcal{N}(\bm{\mu},\bm{\Sigma})\{\mathcal{A}\}\coloneqq\mathbb{P}(\bm{z}\in\mathcal{A})$,
where $\bm{z}\sim\mathcal{N}(\bm{\mu},\bm{\Sigma})$. Throughout this
paper, we let $\Phi(\cdot)$ (resp.~$\phi(\cdot)$) represent the
cumulative distribution function (resp.~probability distribution
function) of the standard Gaussian distribution. We also denote by
$\chi^{2}_k$ the chi-square distribution with $k$ degrees of freedom. 

\section{Background: the estimation algorithm \textsf{HeteroPCA} \label{sec:Background:-two-estimation}}

In order to conduct statistical inference for PCA, the first step
lies in selecting an algorithm to estimate the principal subspace
and the covariance matrix of interest, which we discuss in this
section. Before continuing, we introduce several useful matrix notation
as follows\begin{subequations} 
\begin{align}
\bm{X} & \coloneqq[\bm{x}_{1},\cdots,\bm{x}_{n}]\in\mathbb{R}^{d\times n},\label{eq:definition-X-matrix}\\
\bm{Y} & \coloneqq\mathcal{P}_{\Omega}(\bm{X}+\bm{N})\in\mathbb{R}^{d\times n},\label{eq:defn-Y-matrix}
\end{align}
\end{subequations}where $\mathcal{P}_{\Omega}$ has been defined
in Section~\ref{subsec:Paper-organization-notation}, and $\bm{N}\in\mathbb{R}^{d\times n}$
represents the noise matrix such that the $(l,j)$-th entry of $\bm{N}$
is given by $\eta_{l,j}$. In other words, $\bm{Y}$ encapsulates
all the observed data $\{y_{l,j}\mid(l,j)\in\Omega\}$, with any entry
outside $\Omega$ taken to be zero. If one has full access to the
noiseless data matrix $\bm{X}$, then a natural strategy to estimate
$\bm{U}^{\star}$ would be to return the top-$r$ eigenspace of the
sample covariance matrix $n^{-1}\bm{X}\bm{X}^{\top}$, or equivalently,
the top-$r$ left singular subspace of $\bm{X}$. In practice, however,
one needs to extract information from the corrupted and incomplete
data matrix $\bm{Y}$.

\paragraph{A vanilla SVD-based approach. }

Given that $p^{-1}\bm{Y}=p^{-1}\mathcal{P}_{\Omega}(\bm{X}+\bm{N})$
is an unbiased estimate of $\bm{X}$ (conditional on $\bm{X}$), a
natural idea that comes into mind is to resort to the top-$r$ left
singular subspace of $p^{-1}\bm{Y}$ when estimating $\bm{U}^{\star}$.
This simple procedure is summarized in Algorithm \ref{alg:PCA-SVD}.

\begin{algorithm}[h]
\caption{A vanilla SVD-based approach.}

\label{alg:PCA-SVD}\begin{algorithmic}

\STATE \textbf{{Input}}: data matrix $\bm{Y}$
(cf.~(\ref{eq:defn-Y-matrix})), sampling rate $p$, rank $r$.

\STATE \textbf{{Compute} }the truncated rank-$r$ SVD $\bm{U}\bm{\Sigma}\bm{V}^{\top}$
of $p^{-1}\bm{Y}/\sqrt{n}$, where $\bm{U}\in\mathbb{R}^{d\times r}$,
$\bm{\Sigma}\in\mathbb{R}^{r\times r}$ and $\bm{V}\in\mathbb{R}^{n\times r}$.

\STATE \textbf{{Output}}: $\bm{U}$ as the subspace estimate, $\bm{\Sigma}$
as an estimate of $(\bm{\Lambda}^{\star})^{1/2}$, and $\bm{S}=\bm{U}\bm{\Sigma}^{2}\bm{U}^{\top}$
as the covariance matrix estimate of $\bm{x}$.

\end{algorithmic} 
\end{algorithm}

\paragraph{An improved iterative estimator: \textsf{HeteroPCA}. }

While Algorithm~\ref{alg:PCA-SVD} returns reliable estimates of
$\bm{U}^{\star}$ and $\bm{S}^{\star}$ in the regime of moderate-to-high
signal-to-noise ratio (SNR), it might fail to be effective if either
the missing rate $1-p$ or the noise levels are too large. To offer
a high-level explanation, we find it helpful to compute the expectation of a properly rescaled sample covariance matrix:
\begin{equation}
\frac{1}{p^{2}}\mathbb{E}\left[\bm{Y}\bm{Y}^{\top}\,\big|\,\bm{X}\right]=\bm{X}\bm{X}^{\top}+\left(\frac{1}{p}-1\right)\mathcal{P}_{\mathsf{diag}}\left(\bm{X}\bm{X}^{\top}\right)+\frac{n}{p}\mathsf{diag}\left\{ \left[\omega_{l}^{\star2}\right]_{1\leq l\leq d}\right\} ,\label{eq:PCA-E-gram-decompose}
\end{equation}
where for any vector $\bm{z}=[z_{l}]_{1\leq l\leq d}$ we denote by
$\mathsf{diag}(\bm{z})\in\mathbb{R}^{d\times d}$ a diagonal matrix
whose $(l,l)$-th entry equals $z_{l}$. Here, we rescale the sample covariance matrix by $p^{-2}$ on the left-hand side, 
given that $p^{-1}\bm{Y}$ is an unbiased estimate for $\bm{X}$ and therefore we expect $p^{-2}\bm{Y}\bm{Y}^\top$ to be close to $\bm{X}\bm{X}^\top$. If the sampling rate $p$
is overly small and/or if the noise is of large size but heteroskedastic,
then the second and the third terms on the right-hand side of (\ref{eq:PCA-E-gram-decompose})
might result in significant bias on the diagonal of the matrix $\mathbb{E}\left[\bm{Y}\bm{Y}^{\top}\,\big|\,\bm{X}\right]$,
thus hampering the statistical accuracy of the eigenspace of $p^{-2}\bm{Y}\bm{Y}^{\top}$
(or equivalently, the left singular space of $p^{-1}\bm{Y}$) when
employed to estimate $\bm{U}^{\star}$. Viewed in this light, a more
effective estimator would include procedures that properly handle
the diagonal components of $p^{-2}\bm{Y}\bm{Y}^{\top}$. 

To remedy this issue, several previous works (e.g., \citet{florescu2016spectral,cai2019subspace})
adopted a spectral method with diagonal deletion, which essentially
discards any diagonal entry of $p^{-2}\bm{Y}\bm{Y}^{\top}$ before
computing its top-$r$ eigenspace. However, diagonal deletion comes
at a price: while this operation mitigates the significant bias due
to heteroskedasticity and missing data, it introduces another type
of bias that might be non-negligible if the goal is to enable efficient
fine-grained inference. To address this bias issue, \citet{zhang2018heteroskedastic}
proposed an iterative refinement scheme --- termed \textsf{HeteroPCA}
--- that copes with the diagonal entries in a more refined manner.
Informally, \textsf{HeteroPCA} starts by computing the rank-$r$ eigenspace
of the diagonal-deleted version of $p^{-2}\bm{Y}\bm{Y}^{\top}$, and
then alternates between imputing the diagonal entries of $\bm{X}\bm{X}^{\top}$
and estimating the eigenspace of $p^{-2}\bm{Y}\bm{Y}^{\top}$ with
the aid of the imputed diagonal. A precise description of this procedure
is summarized in Algorithm \ref{alg:PCA-HeteroPCA}; here, $\mathcal{P}_{\mathsf{off}\text{-}\mathsf{diag}}$
and $\mathcal{P}_{\mathsf{diag}}$ have been defined in Section~\ref{subsec:Paper-organization-notation}.

\begin{algorithm}[h]
\caption{\textsf{HeteroPCA} (by \citet{zhang2018heteroskedastic}).}

\label{alg:PCA-HeteroPCA}\begin{algorithmic}

\STATE \textbf{{Input}}:  data matrix $\bm{Y}$
(cf.~(\ref{eq:defn-Y-matrix})), sampling rate $p$, rank $r$, maximum
number of iterations $t_{0}$.

\STATE \textbf{{Initialization}}: set $\bm{G}^{0}=\frac{1}{np^{2}}\mathcal{P}_{\mathsf{off}\text{-}\mathsf{diag}}(\bm{Y}\bm{Y}^{\top})$.

\STATE \textbf{{Updates}}: \textbf{for }$t=0,1,\ldots,t_{0}$ \textbf{do}

\STATE \vspace{-1em}
 \begin{subequations} 
\begin{align*}
\left(\bm{U}^{t},\bm{\Lambda}^{t}\right) & =\mathsf{eigs}\left(\bm{G}^{t},r\right);\\  
\bm{G}^{t+1} & =\mathcal{P}_{\mathsf{off}\text{-}\mathsf{diag}}\left(\bm{G}^{t}\right)+\mathcal{P}_{\mathsf{diag}}\left(\bm{U}^{t}\bm{\Lambda}^{t}\bm{U}^{t\top}\right)=\frac{1}{np^{2}}\mathcal{P}_{\mathsf{off}\text{-}\mathsf{diag}}\left(\bm{Y}\bm{Y}^{\top}\right)+\mathcal{P}_{\mathsf{diag}}\left(\bm{U}^{t}\bm{\Lambda}^{t}\bm{U}^{t\top}\right).  
\end{align*}
\end{subequations} Here, for any symmetric matrix $\bm{G}\in\mathbb{R}^{d\times d}$
and $1\leq r\leq d$, $\mathsf{eigs}(\bm{G},r)$ returns $(\bm{U},\bm{\Lambda})$,
where $\bm{U}\bm{\Lambda}\bm{U}^{\top}$ is the top-$r$ eigen-decomposition
of $\bm{G}$.

\STATE \textbf{{Output}}: $\bm{U}=\bm{U}^{t_{0}}$ as the subspace
estimate, $\bm{\Sigma}=(\bm{\Lambda}^{t_{0}})^{1/2}$ as an estimate
of $(\bm{\Lambda}^{\star})^{1/2}$, and $\bm{S}=\bm{U}^{t_{0}}\bm{\Lambda}^{t_{0}}\bm{U}^{t_{0}\top}$
as the covariance matrix estimate.

\end{algorithmic} 
\end{algorithm}

\section{Distributional theory and inference procedures\label{sec:Distributional-theory-inference}}

In this section, we augment the \textsf{HeteroPCA} estimator introduced
in Section~\ref{sec:Background:-two-estimation} by a suite of distributional
theory, and demonstrate how to employ our distributional characterizations
to perform inference on both the principal subspace represented by
$\bm{U}^{\star}$ and the covariance matrix $\bm{S}^{\star}$.

\subsection{Key quantities and assumptions}

Before continuing, we introduce several additional notation and assumptions
that play a key role in our theoretical development. Recall that the
eigen-decomposition of the covariance matrix $\bm{S}^{\star}\in\mathbb{R}^{d\times d}$
(see (\ref{eq:covariance-spike})) is assumed to be $\bm{U}^{\star}\bm{\Lambda}^{\star}\bm{U}^{\star\top}.$
We assume the diagonal matrix $\bm{\Lambda}^{\star}$ to be $\bm{\Lambda}^{\star}=\mathsf{diag}\{\lambda_{1}^{\star},\ldots,\lambda_{r}^{\star}\}$,
where the diagonal entries are given by the non-zero eigenvalues of
$\bm{S}^{\star}$ obeying 
\[
\lambda_{1}^{\star}\geq\cdots\geq\lambda_{r}^{\star}>0.
\]
The condition number of $\bm{S}^{\star}$ is denoted by 
\begin{equation}
\kappa\coloneqq\lambda_{1}^{\star}/\lambda_{r}^{\star}.\label{eq:defn-condition-number}
\end{equation}
We also find it helpful to introduce the square root of $\bm{\Lambda}^{\star}$
as follows 
\begin{equation}
\bm{\Sigma}^{\star}=\mathsf{diag}\{\sigma_{1}^{\star},\ldots,\sigma_{r}^{\star}\}=(\bm{\Lambda}^{\star})^{1/2},\qquad\text{where }\sigma_{i}^{\star}=(\lambda_{i}^{\star})^{1/2},\ 1\leq i\leq r.\label{eq:defn-Sigma-star-sqrt-Lambda}
\end{equation}
Furthermore, we introduce an incoherence parameter commonly employed
in prior literature \citep{candes2014mathematics,chi2018nonconvex}.

\begin{definition}[\bf Incoherence]\label{assumption:incoherence}The
rank-$r$ matrix $\bm{S}^{\star}\in\mathbb{R}^{d\times d}$ defined
in (\ref{eq:covariance-spike}) is said to be $\mu$-incoherent if
the following condition holds: 
\begin{align}
\left\Vert \bm{U}^{\star}\right\Vert _{2,\infty} & \leq\sqrt{\frac{\mu r}{d}}.\label{eq:Ustar-incoherence-condition}
\end{align}
Here, we recall that $\Vert\bm{U}^{\star}\Vert_{2,\infty}$ denotes
the largest $\ell_{2}$ norm of all rows of the matrix $\bm{U}^{\star}$.
\end{definition} \begin{remark}When $\mu$ is small (e.g., $\mu\asymp1$),
this condition essentially ensures that the energy of $\bm{U}^{\star}$
is nearly evenly dispersed across all of its rows. As a worthy note,
the theory developed herein allows the incoherence parameter $\mu$
to grow with the problem dimension. \end{remark}

In light of a global rotational ambiguity issue (i.e., for any $r\times r$
rotation matrix $\bm{R}$, the matrices $\bm{U}\in\mathbb{R}^{d\times r}$
and $\bm{U}\bm{R}\in\mathbb{R}^{d\times r}$ share the same column
space), in general we can only hope to estimate $\bm{U}^{\star}$
up to global rotation (unless additional eigenvalue separation conditions
are imposed). Consequently, our theoretical development focuses on
characterizing the error distribution $\bm{U}\bm{R}-\bm{U}^{\star}$
of an estimator $\bm{U}$ when accounting for a proper rotation matrix
$\bm{R}$. In particular, we shall pay particular attention to a specific
way of rotation as follows 
\[
\bm{U}\mathsf{sgn}\left(\bm{U}^{\top}\bm{U}^{\star}\right)-\bm{U}^{\star},
\]
where we recall that for any non-singular matrix $\bm{H}\in\mathbb{R}^{k\times k}$
with SVD $\bm{U}_{H}\bm{\Sigma}_{H}\bm{V}_{H}^{\top}$, the matrix
$\mathsf{sgn}(\bm{H})$ is defined to be the rotation matrix $\bm{U}_{H}\bm{V}_{H}^{\top}$.
This particular choice aligns $\bm{U}$ and $\bm{U}^{\star}$ in the
following sense 
\[
\mathsf{sgn}\big(\bm{U}^{\top}\bm{U}^{\star}\big)=\arg\min_{\bm{R}\in\mathcal{O}^{r\times r}}\left\Vert \bm{U}\bm{R}-\bm{U}^{\star}\right\Vert _{\mathrm{F}},
\]
where $\mathcal{O}^{r\times r}$ indicates the set of all $r\times r$
rotation matrices; see \citet[Appendix D.2.1]{ma2017implicit}. 
This is often referred to as Whaba's problem \citep{wahba1965least} or the orthogonal  Procrustes problem \citep{schonemann1966generalized}. 

The last assumption is concerned with the noise levels, which are
allowed to vary across different locations.

\begin{assumption}[\bf Noise levels]\label{assumption:noise}The
noise levels $\{\omega_{i}^{\star}\}_{1\leq i\leq d}$ obey 
\begin{equation}
\frac{\omega_{\max}^{2}}{\omega_{\min}^{2}}\leq\kappa_{\omega}\qquad\text{with}\quad\omega_{\max}\coloneqq\max_{1\leq i\leq d}\omega_{i}^{\star}\quad\text{and}\quad\omega_{\min}\coloneqq\min_{1\leq i\leq d}\omega_{i}^{\star}.\label{eq:defn-sigma-min-sigma-max}
\end{equation}
\end{assumption}

\subsection{Inferential procedure and theory for \textsf{HeteroPCA}\label{subsec:Inference-for-HeteroPCA}}

We are now positioned to investigate how to assess the uncertainty
of the estimator \textsf{HeteroPCA}. For simplicity of presentation,
we shall abuse some notation (e.g., $\bm{\Sigma}_{U,l}^{\star}$ and
$v_{i,j}^{\star}$) whenever it is clear from the context.

\subsubsection{Distributional theory and inference for the principal subspace $\bm{U}^{\star}$} \label{subsec:Inference-for-Subspace}

In this subsection, we shall begin by establishing a distributional theory for the subspace estimate $\bm{U}$ returned by \textsf{HeteroPCA} (see Theorem~\ref{thm:pca}), followed by a data-driven and provably valid  method to construct fine-grained confidence regions for $\bm{U}^{\star}$ (see Algorithm~\ref{alg:PCA-HeteroPCA-CR} and Theorem~\ref{thm:pca-cr}). We shall also briefly discuss how our results improve upon prior estimation guarantees for \textsf{HeteroPCA} in the presence of missing data. 

\paragraph{Distributional guarantees.}

As it turns out, the subspace estimate returned by Algorithm \ref{alg:PCA-HeteroPCA}
is approximately unbiased and Gaussian under milder conditions, as
posited in the following theorem. The general result beyond the case
with $\kappa,\mu,r,\kappa_{\omega}\asymp1$ is postponed to Theorem~\ref{thm:pca-complete}
in Appendix~\ref{sec:A-list-of-general-thms}.

\begin{theorem}\label{thm:pca} Assume that each column of the ground truth $\bm{X}$ (cf.~\eqref{eq:definition-X-matrix}) is independently generated from $\mathcal{N}(\bm{0},\bm{S}^\star)$, and that the sampling set $\Omega$ follows the random sampling model in Section~\ref{subsec:intro-model}. Suppose that $p<1-\delta$   for some
arbitrary constant $0<\delta<1$ or $p=1$, and $\kappa,\mu,r,\kappa_{\omega}\asymp1$.
Assume that Assumption \ref{assumption:noise} holds and $d\gtrsim\log^{5}n$, 
\begin{subequations}\label{eq:conditions-thm-pca-simple} 
\begin{equation}
\frac{\omega_{\max}^{2}}{p\sigma_{r}^{\star2}}\sqrt{\frac{d}{n}}\lesssim\frac{1}{\log^{7/2}\left(n+d\right)},\qquad\frac{\omega_{\max}}{\sigma_{r}^{\star}}\sqrt{\frac{d}{np}}\lesssim\frac{1}{\log^{3}\left(n+d\right)},
\end{equation}
\begin{equation}
ndp^{2}\gtrsim\log^{9}\left(n+d\right),\qquad np\gtrsim\log^{7}\left(n+d\right).
\end{equation}
\end{subequations} Suppose, in addition, that the number of iterations
exceeds 
\begin{equation}
t_{0}\gtrsim\log\left[\left(\frac{\log^{2}\left(n+d\right)}{\sqrt{nd}p}+\frac{\omega_{\max}^{2}}{p\sigma_{r}^{\star2}}\sqrt{\frac{d}{n}}\log\left(n+d\right)+\frac{\log\left(n+d\right)}{\sqrt{np}}+\frac{\omega_{\max}}{\sigma_{r}^{\star}}\sqrt{\frac{d\log\left(n+d\right)}{np}}\right)^{-1}\right].\label{eq:main-iteration}
\end{equation}
Let $\bm{R}$ be the $r\times r$ rotation matrix $\bm{R}=\mathsf{sgn}(\bm{U}^{\top}\bm{U}^{\star})$.
Then the estimate $\bm{U}$ returned by Algorithm \ref{alg:PCA-HeteroPCA}
obeys: for all $1\leq l\leq d$,
\[
\sup_{\mathcal{C}\in\mathscr{C}^{r}}\left|\mathbb{P}\left(\left[\bm{U}\bm{R}-\bm{U}^{\star}\right]_{l,\cdot}\in\mathcal{C}\right)-\mathcal{N}\big(\bm{0},\bm{\Sigma}_{U,l}^{\star}\big)\left\{ \mathcal{C}\right\} \right|=o\left(1\right),
\]
where $\mathscr{C}^{r}$ represents the set of all convex sets in
$\mathbb{R}^{r}$, and 
\begin{align}
\bm{\Sigma}_{U,l}^{\star} & \coloneqq\left(\frac{1-p}{np}\left\Vert \bm{U}_{l,\cdot}^{\star}\bm{\Sigma}^{\star}\right\Vert _{2}^{2}+\frac{\omega_{l}^{\star2}}{np}\right)\left(\bm{\Sigma}^{\star}\right)^{-2}+\frac{2\left(1-p\right)}{np}\bm{U}_{l,\cdot}^{\star\top}\bm{U}_{l,\cdot}^{\star}\nonumber \\
 & \quad+\left(\bm{\Sigma}^{\star}\right)^{-2}\bm{U}^{\star\top}\mathsf{diag}\left\{ \left[d_{l,i}^{\star}\right]_{1\leq i\leq d}\right\} \bm{U}^{\star}(\bm{\Sigma}^{\star})^{-2}\label{eq:pca-heteropca-true-covariance}
\end{align}
with 
\[
d_{l,i}^{\star}\coloneqq\frac{1}{np^{2}}\left[\omega_{l}^{\star2}+\left(1-p\right)\left\Vert \bm{U}_{l,\cdot}^{\star}\bm{\Sigma}^{\star}\right\Vert _{2}^{2}\right]\left[\omega_{i}^{\star2}+\left(1-p\right)\left\Vert \bm{U}_{i,\cdot}^{\star}\bm{\Sigma}^{\star}\right\Vert _{2}^{2}\right]+\frac{2\left(1-p\right)^{2}}{np^{2}}S_{l,i}^{\star2}.
\]

\end{theorem}

Theorem~\ref{thm:pca} asserts that each row of the estimate $\bm{U}$
returned by \textsf{HeteroPCA} is nearly unbiased and admits a nearly
tight Gaussian approximation, whose covariance matrix can be determined
via the closed-form expression (\ref{eq:pca-heteropca-true-covariance}).
Given that $\bm{U}\bm{R}$ and $\bm{U}$ represent the same subspace,
this theorem delivers a fine-grained row-wise distributional characterization
for the estimator \textsf{HeteroPCA}.

Let us briefly mention the key error decomposition behind this theorem, which might help illuminate how Gaussian approximation emerges.
Letting $\bm{E}\coloneqq n^{-1/2}(p^{-1}\bm{Y}-\bm{X})$ (which captures the randomness from both the noise and random subsampling), we can decompose
\begin{equation} \label{eq:decomposition}
\bm{U}\bm{R}-\bm{U}^{\star}=\underbrace{\left[\bm{E}\bm{X}^{\top}+\mathcal{P}_{\mathsf{off}\text{-}\mathsf{diag}}\left(\bm{E}\bm{E}^{\top}\right)\right]\bm{U}^{\star}\left(\bm{\Sigma}^{\star}\right)^{-2}}_{\eqqcolon\,\bm{Z}\,\,\text{(first- and second-order approximation)}}+\underbrace{\left[\bm{U}\bm{R}-\bm{U}^{\star}-\bm{Z}\right]}_{\eqqcolon\,\bm{\Psi}\,\,\text{(residual term)}}.
\end{equation}
Here, $\bm{Z}$ contains not only a linear mapping of $\bm{E}$ but
also a certain quadratic mapping, the latter of which is crucial when
coping with the regime $n\gg d$. As a consequence of the central limit theorem (which will be solidified in the proof),
$\bm{Z}$ admits the following Gaussian approximation 
\begin{equation}
\bm{Z}_{l,\cdot}\overset{\mathrm{d}}{\approx}\mathcal{N}\big(\bm{0},\bm{\Sigma}_{U,l}^{\star}\big),\qquad1\leq l\leq d.
	\label{eq:Gaussian-approximation-hetero-PCA-simple}
\end{equation}
At the same time, the $\ell_{2}$ norm of the residual term 
	$\bm{\Psi}_{l,\cdot}$ is well controlled and provably negligible compared to the corresponding component in $\bm{Z}_{l,\cdot}$,
thus ascertaining the tightness of the advertised Gaussian approximation. 

\begin{remark} \label{remark: kappa-omega-necessity} 
	The decomposition \eqref{eq:decomposition} also sheds light on why our current theory assumes a finite $\kappa_{\omega}$ (cf.~\eqref{eq:defn-sigma-min-sigma-max}) when conducting statistical inference (which is unnecessary for the task of estimation). 
	Consider, for example, a simple case when (i) there is no missing data ($p=1$), and (ii) for some $1\leq l \leq d$, one has $\omega_l^\star=0$ (and hence $\kappa_{\omega}=\infty$) and  $\Vert\bm{U}_{l,\cdot}^\star\Vert_2>0$. In this case, $\bm{Z}_{l,\cdot}=\bm{0}$ since $\bm{\Sigma}_{U,l}^{\star}=\bm{0}$, although $[\bm{U}\bm{R}-\bm{U}^{\star}]_{l,\cdot}$ is in general non-zero. This implies that our Gaussian approximation  --- and the inference procedure developed based on this approximation ---
	might fall short of efficacy when $\kappa_{\omega}=\infty$. 
\end{remark}

\paragraph{Construction of confidence regions for the principal subspace. }

With the above distributional theory in place, we are well-equipped
to construct fine-grained confidence regions for $\bm{U}^{\star}$,
provided that the covariance matrix $\bm{\Sigma}_{U,l}^{\star}$ can
be estimated in a faithful manner. In Algorithm \ref{alg:PCA-HeteroPCA-CR},
we propose a procedure to estimate $\bm{\Sigma}_{U,l}^{\star}$, 
which in turn allows us to build confidence regions. As before, our estimator
for $\bm{\Sigma}_{U,l}^{\star}$ can be viewed as a sort of ``plug-in''
method in accordance with the expression (\ref{eq:pca-heteropca-true-covariance}).

\begin{algorithm}[h]
\caption{Confidence regions for $\bm{U}_{l,\cdot}^{\star}$ $(1\protect\leq l\protect\leq d)$
based  on \textsf{HeteroPCA}.}

\label{alg:PCA-HeteroPCA-CR}\begin{algorithmic}

\STATE \textbf{{Input}}: output $(\bm{U},\bm{\Sigma},\bm{S})$
of Algorithm \ref{alg:PCA-HeteroPCA}, sampling rate $p$, coverage
level $1-\alpha$.

\STATE \textbf{{Compute}} estimates of the noise levels $\left\{ \omega_{l}^{\star}\right\} _{1\leq l\leq d}$
as follows 
\[
\omega_{l}^{2}\coloneqq\frac{\sum_{j=1}^{n}y_{l,j}^{2}\ind_{(l,j)\in\Omega}}{\sum_{j=1}^{n}\ind_{(l,j)\in\Omega}}-S_{l,l}\qquad\text{for all}\quad1\leq l\leq d.
\]

\STATE \textbf{{Compute} }an estimate of $\bm{\Sigma}_{U,l}^{\star}$
(cf.~(\ref{eq:pca-heteropca-true-covariance})) as follows: 
\begin{align*}
\bm{\Sigma}_{U,l} & \coloneqq\left(\frac{1-p}{np}\left\Vert \bm{U}_{l,\cdot}\bm{\Sigma}\right\Vert _{2}^{2}+\frac{\omega_{l}^{2}}{np}\right)\bm{\Sigma}^{-2}+\frac{2\left(1-p\right)}{np}\bm{U}_{l,\cdot}^{\top}\bm{U}_{l,\cdot}+\left(\bm{\Sigma}\right)^{-2}\bm{U}^{\top}\mathsf{diag}\left\{ \left[d_{l,i}\right]_{1\leq i\leq d}\right\} \bm{U}(\bm{\Sigma})^{-2},
\end{align*}
where 
\[
d_{l,i}\coloneqq\frac{1}{np^{2}}\left[\omega_{l}^{2}+\left(1-p\right)\left\Vert \bm{U}_{l,\cdot}\bm{\Sigma}\right\Vert _{2}^{2}\right]\left[\omega_{i}^{2}+\left(1-p\right)\left\Vert \bm{U}_{i,\cdot}\bm{\Sigma}\right\Vert _{2}^{2}\right]+\frac{2\left(1-p\right)^{2}}{np^{2}}S_{l,i}^{2}.
\]

\STATE \textbf{{Compute} }the $(1-\alpha)$-quantile $\tau_{1-\alpha}$
of $\chi_r^{2}$ and construct a Euclidean ball: 
\[
\mathcal{B}_{1-\alpha}\coloneqq\left\{ \bm{z}\in\mathbb{R}^{r}:\left\Vert \bm{z}\right\Vert _{2}^{2}\leq\tau_{1-\alpha}\right\} .
\]

\STATE \textbf{{Output}} the $(1-\alpha)$-confidence region 
\[
\mathsf{CR}_{U,l}^{1-\alpha}\coloneqq\bm{U}_{l,\cdot}+\big(\bm{\Sigma}_{U,l}\big)^{1/2}\mathcal{B}_{1-\alpha}=\left\{ \bm{U}_{l,\cdot}+\big(\bm{\Sigma}_{U,l}\big)^{1/2}\bm{z}:\bm{z}\in\mathcal{B}_{1-\alpha}\right\} .
\]

\end{algorithmic} 
\end{algorithm}

The following theorem confirms the validity of the proposed inference
procedure when $\kappa,\mu,r,\kappa_{\omega}\asymp1$. The more general
case will be studied in Theorem~\ref{thm:pca-cr-complete} in Appendix~\ref{sec:A-list-of-general-thms}.

\begin{theorem}\label{thm:pca-cr}Suppose that the conditions of
Theorem \ref{thm:pca} hold. Then there exists a $r\times r$ rotation
matrix $\bm{R}=\mathsf{sgn}\left(\bm{U}^{\top}\bm{U}^{\star}\right)$
such that the confidence regions $\mathsf{CR}_{U,l}^{1-\alpha}$ ($1\leq l\leq d$)
computed in Algorithm \ref{alg:PCA-HeteroPCA-CR} obey 
\[
\sup_{1\leq l\leq d}\left|\mathbb{P}\left(\bm{U}_{l,\cdot}^{\star}\bm{R}^{\top}\in\mathsf{CR}_{U,l}^{1-\alpha}\right)-\left(1-\alpha\right)\right|=o\left(1\right).
\]
\end{theorem}In words, Theorem~\ref{thm:pca-cr} uncovers that:
a valid ground-truth subspace representation is contained --- in
a row-wise reliable manner --- within the confidence regions $\mathsf{CR}_{U,l}^{1-\alpha}$
($1\leq l\leq d$) we construct. In the special case with $r=1$,
this result leads to valid entrywise confidence intervals for the
principal component.

\paragraph{Interpretations and implications. }

We now take a moment to interpret the conditions required in Theorem
\ref{thm:pca} and Theorem \ref{thm:pca-cr}, and discuss some appealing
attributes of our methods. As before, the discussion below focuses
on the scenario where $\mu,\kappa,r,\kappa_{\omega}\asymp1$ for the
sake of simplicity. 
\begin{itemize}
\item \textit{Missing data. }Both theorems accommodate the case when a large
fraction of data are missing, namely, they cover the range 
\[
p\geq\widetilde{\Omega}\Big(\frac{1}{n\land\sqrt{nd}}\Big)
\]
for both distributional characterizations and confidence region construction
using \textsf{HeteroPCA}. In particular, if $n\gg d$, then the sampling
rate $p$ only needs to exceed 
\[
p\geq\widetilde{\Omega}\Big(\frac{1}{\sqrt{nd}}\Big);
\]
this range can include some sampling rate much smaller than $1/d$
(with $d$ the ambient dimension of each sample vector), and cannot
be improved in general according to \cite[Theorem 3.4]{cai2019subspace}. 
\item \textit{Tolerable noise levels.} The noise conditon required in both
Theorem \ref{thm:pca} and Theorem \ref{thm:pca-cr} is given by 
\[
\omega_{\max}^{2}\leq\widetilde{O}\left(\Big(\frac{n}{d}\land\sqrt{\frac{n}{d}}\Big)p\sigma_{r}^{\star2}\right).
\]
Note that when $\kappa,\mu,r\asymp1$, the variance obeys 
\[
\max_{(l,j)\in\Omega}\mathsf{var}\left(x_{l,j}\right)=\max_{l\in[d]}S_{l,l}^{\star}\asymp\max_{l\in[d]}\big\|\bm{U}_{l,\cdot}^{\star}\big\|_{2}^{2}\sigma_{1}^{\star2}\asymp\frac{1}{d}\sigma_{1}^{\star2}.
\]
This implies that: when $p\geq\widetilde{\Omega}\big(1/(n\land\sqrt{nd})\big)$,
our tolerable entrywise noise level $\omega_{\max}^{2}$ is allowed
to be significantly (i.e., $\widetilde{\Omega}(np\land\sqrt{ndp^{2}})$
times) larger than the largest variance of $x_{l,j}$ for all $(l,j)\in\Omega$,
thereby accommodating a wide range of noise levels. 
\item \emph{Adaptivity to heteroskedasticity and unknown noise levels}.
Our proposed inferential procedure is fully data-driven: it is automatically
adaptive to unknown heteroskedastic noise, without requiring prior
knowledge of the noise levels. 
\end{itemize}

\paragraph{Comparison with prior estimation theory. }

While the main purpose of the current paper is to enable efficient
statistical inference for the principal subspace, our theory (see Lemmas~\ref{lemma:pca-1st-err} and \ref{lemma:pca-noise-level-est} in the appendix)
also enables improved estimation guarantees compared to prior works.
\begin{itemize}
\item Recall that the estimation algorithm \textsf{HeteroPCA} was originally
proposed and studied by \citet{zhang2018heteroskedastic}. Our results
broaden the sample size range supported by their theory. More specifically,
note that \citet[Theorem 6 and Remark 10]{zhang2018heteroskedastic}
requires the sampling rate $p$ to satisfy 
\[
ndp\gtrsim\max\left\{ d^{1/3}n^{2/3},d\right\} \mathsf{polylog}\left(n,d\right)
\]
in order to guarantee consistent estimation, while our theoretical
guarantees only require 
\[
ndp\gtrsim\max\left\{ \sqrt{nd},d\right\} \mathsf{polylog}\left(n,d\right).
\]
When $n\gg d$, the sample size requirement in \citet{zhang2018heteroskedastic}
is $(n/d)^{1/6}$ times more stringent than the one imposed in our
theory. 
\item Let us discuss the advantage of \textsf{HeteroPCA}
compared to the diagonal-deleted spectral method studied in \citet[Algorithms 1 and 3]{cai2019subspace}.
Due to diagonal deletion, there is an additional bias term (see the
last term $\mu_{\mathsf{ce}}\kappa_{\mathsf{ce}}r/d$ in Equation~(4.16)
in \citet{cai2019subspace}), which turns out to negatively affect
our capability of performing inference. In contrast, \textsf{HeteroPCA}
eliminates this bias term by means of successive refining, thus facilitating
the subsequent inference stage. 
\end{itemize}

\subsubsection{Distributional theory and inference for the covariance matrix $\bm{S}^{\star}$} \label{subsec:Inference-for-Covariance}

As it turns out, the above distributional theory for $\bm{U}^{\star}$
further hints at how to perform statistical inference for the covariance
matrix $\bm{S}^{\star}$. 
In the sequel, we shall first develop an entrywise distributional theory 
for the estimate $\bm{S}$ returned by \textsf{HeteroPCA} (see Theorem~\ref{thm:ce}), 
followed by a data-driven inference procedure to conduct entrywise confidence intervals for $\bm{S}^{\star}$ (see Algorithm~\ref{alg:CE-HeteroPCA-CI} and Theorem~\ref{thm:ce-CI}).

\paragraph{Entrywise distributional guarantees.}

We now focus attention on characterizing the distribution of the $(i,j)$-th
entry of $\bm{S}$ returned by Algorithm \ref{alg:PCA-HeteroPCA},
which in turn suggests how to construct entrywise confidence intervals
for $\bm{S}^{\star}$. Before proceeding, let us define a set of variance
parameters $\{v_{i,j}^{\star}\}_{1\leq i,j\leq d}$ which, as we shall
demonstrate momentarily, correspond to the (approximate) variance
of the entries of $\bm{S}$. 
\begin{itemize}
\item For any $1\leq i,j\leq d$ obeying $i\neq j$, we define 
\begin{align}
v_{i,j}^{\star} & \coloneqq\frac{2-p}{np}S_{i,i}^{\star}S_{j,j}^{\star}+\frac{4-3p}{np}S_{i,j}^{\star2}+\frac{1}{np}\left(\omega_{i}^{\star2}S_{j,j}^{\star}+\omega_{j}^{\star2}S_{i,i}^{\star}\right)\nonumber \\
 & \quad+\frac{1}{np^{2}}\sum_{k=1}^{d}\left\{ \left[\omega_{i}^{\star2}+\left(1-p\right)S_{i,i}^{\star}\right]\left[\omega_{k}^{\star2}+\left(1-p\right)S_{k,k}^{\star}\right]+2\left(1-p\right)^{2}S_{i,k}^{\star2}\right\} \left(\bm{U}_{k,\cdot}^{\star}\bm{U}_{j,\cdot}^{\star\top}\right)^{2}\nonumber \\
 & \quad+\frac{1}{np^{2}}\sum_{k=1}^{d}\left\{ \left[\omega_{j}^{\star2}+\left(1-p\right)S_{j,j}^{\star}\right]\left[\omega_{k}^{\star2}+\left(1-p\right)S_{k,k}^{\star}\right]+2\left(1-p\right)^{2}S_{j,k}^{\star2}\right\} \left(\bm{U}_{k,\cdot}^{\star}\bm{U}_{i,\cdot}^{\star\top}\right)^{2}.\label{eq:ce-true-variance-i-j}
\end{align}
\item For any $1\leq i\leq d$, we set 
\begin{align}
v_{i,i}^{\star} & \coloneqq\frac{12-9p}{np}S_{i,i}^{\star2}+\frac{4}{np}\omega_{i}^{\star2}S_{i,i}^{\star}\nonumber \\
 & \quad+\frac{4}{np^{2}}\sum_{k=1}^{d}\left\{ \left[\omega_{i}^{\star2}+\left(1-p\right)S_{i,i}^{\star}\right]\left[\omega_{k}^{\star2}+\left(1-p\right)S_{k,k}^{\star}\right]+2\left(1-p\right)^{2}S_{i,k}^{\star2}\right\} \left(\bm{U}_{k,\cdot}^{\star}\bm{U}_{i,\cdot}^{\star\top}\right)^{2}.\label{eq:ce-true-variance-i}
\end{align}
\end{itemize}
We are now positioned to present our distributional theory for the
scenario where $\kappa,\mu,r,\kappa_{\omega}\asymp1$, with the more
general version deferred to Theorem~\ref{thm:ce-complete} in Appendix~\ref{sec:A-list-of-general-thms}.
Here and throughout, $S_{i,j}$ (resp.~$S_{i,j}^{\star}$) represents
the $(i,j)$-th entry of the matrix $\bm{S}$ (resp.~$\bm{S}^{\star}$).

\begin{theorem}\label{thm:ce}Suppose that $p<1-\delta$ for some
arbitrary constant $0<\delta<1$ or $p=1$, and $\kappa,\mu,r,\kappa_{\omega}\asymp1$. 
Consider any $1\leq i,j\leq d$. Assume that $\bm{U}^{\star}$ is
$\mu$-incoherent and satisfies the following condition
\begin{equation}
\left\Vert \bm{U}_{i,\cdot}^{\star}\right\Vert _{2}+\left\Vert \bm{U}_{j,\cdot}^{\star}\right\Vert _{2}\gtrsim\left[\frac{\omega_{\max}}{\sigma_{r}^{\star}}\sqrt{\frac{d\log^{5}\left(n+d\right)}{np}}+\frac{\omega_{\max}^{2}}{p\sigma_{r}^{\star2}}\sqrt{\frac{d\log^{5}\left(n+d\right)}{n}}+\sqrt{\frac{\log^{7}\left(n+d\right)}{ndp^{2}}}\right]\sqrt{\frac{1}{d}}.\label{eq:ce-Ui-Uj-max-lb}
\end{equation}
In addition, suppose that Assumption \ref{assumption:noise} holds, and
\[
d\gtrsim\log^{5}n,\qquad np\gtrsim\log^{7}(n+d),\qquad ndp^{2}\gtrsim\log^{7}(n+d), 
\]
\[
\frac{\omega_{\max}}{\sigma_{r}^{\star}}\sqrt{\frac{d}{np}}\lesssim\frac{1}{\log^{3}\left(n+d\right)},\qquad\frac{\omega_{\max}^{2}}{p\sigma_{r}^{\star2}}\sqrt{\frac{d}{n}}\lesssim\frac{1}{\log^{7/2}\left(n+d\right)}.
\]
Assume that the number of iterations satisfies (\ref{eq:main-iteration}).
Then the matrix $\bm{S}$ computed by Algorithm \ref{alg:PCA-HeteroPCA}
obeys 
\[
\sup_{t\in\mathbb{R}}\left|\mathbb{P}\left(\frac{S_{i,j}-S_{i,j}^{\star}}{\sqrt{v_{i,j}^{\star}}}\leq t\right)-\Phi\left(t\right)\right|=o\left(1\right),
\]
where $\Phi(\cdot)$ denotes the CDF of the standard Gaussian distribution.
\end{theorem}

In words, the above theorem indicates that under the conditions in
Theorem \ref{thm:pca}, if the sum of the $\ell_{2}$ norm of the
rows $\bm{U}_{i,\cdot}^{\star}$ and $\bm{U}_{j,\cdot}^{\star}$ are
not exceedingly small, then the estimation error $S_{i,j}-S_{i,j}^{\star}$
of \textsf{HeteroPCA} is approximately a zero-mean Gaussian with variance
$v_{i,j}^{\star}$. 

\begin{remark}When it comes to inference for $S_{i,j}^{\star}$,
our theorems impose the following condition (cf.~\eqref{eq:ce-Ui-Uj-max-lb}):
\[
\left\Vert \bm{U}_{i,\cdot}^{\star}\right\Vert _{2}+\left\Vert \bm{U}_{j,\cdot}^{\star}\right\Vert _{2}\geq\widetilde{\Omega}\left(\frac{1}{\sqrt{ndp^{2}}}+\frac{\omega_{\max}}{\sigma_{r}^{\star}}\sqrt{\frac{d}{np}}+\frac{\omega_{\max}^{2}}{p\sigma_{r}^{\star2}}\sqrt{\frac{d}{n}}\right)\cdot\sqrt{\frac{1}{d}}\left\Vert \bm{U}^{\star}\right\Vert _{\mathrm{F}}.
\]
Note that the typical $\ell_{2}$ norm of a row of $\bm{U}^{\star}$
is $\Vert\bm{U}^{\star}\Vert_{\mathrm{F}}/\sqrt{d}$ when the energy
is uniformly spread out across all rows. This means that under our
sampling rate condition, our results allow $\Vert\bm{U}_{i,\cdot}^{\star}\Vert_{2}+\Vert\bm{U}_{j,\cdot}^{\star}\Vert_{2}$
to be much smaller than its typical size. 
As it turns out, 
a lower bound on $\Vert\bm{U}_{i,\cdot}^\star\Vert_2+\Vert\bm{U}_{j,\cdot}^\star\Vert_2$ 
might be necessary for $S_{i,j}-S_{i,j}^\star$ to be approximately Gaussian. Consider, for example, the case when $\Vert\bm{U}_{i,\cdot}\Vert_2=\Vert\bm{U}_{j,\cdot}\Vert_2=0$. It can be seen from our analysis  that 
\[
S_{i,j}-S_{i,j}^\star \approx \bm{Z}_{i,\cdot}\bm{\Sigma}^{\star2}\bm{Z}_{j,\cdot}^\top+A_{i,j}
\]
where $\bm{Z}_{i,\cdot}$, $\bm{Z}_{j,\cdot}$ and $A_{i,j}$ are all (approximately) Gaussian. 
	This means that $S_{i,j}-S_{i,j}^\star$ might not follow the (approximate) Gaussian distribution claimed in Theorem \ref{thm:ce}
	if $\Vert\bm{U}_{i,\cdot}^\star\Vert_2+\Vert\bm{U}_{j,\cdot}^\star\Vert_2$ is too small. 
\end{remark}


\paragraph{Construction of entrywise confidence intervals.}

The distributional characterization in Theorem~\ref{thm:ce} enables
valid construction of entrywise confidence intervals for $\bm{S}^{\star}$,
as long as we can obtain reliable estimate of the variance $v_{i,j}^{\star}$.
In what follows, we come up with an algorithm --- as summarized in
Algorithm \ref{alg:CE-HeteroPCA-CI} --- that attempts to estimate
$v_{i,j}^{\star}$ and build confidence intervals in a data-driven
manner, as confirmed by the following theorem for the scenario with
$\kappa,\mu,r,\kappa_{\omega}\asymp1$. The more general result is
postponed to Theorem~\ref{thm:ce-CI-complete} in Appendix~\ref{sec:A-list-of-general-thms}.

\begin{algorithm}[h]
\caption{Confidence intervals for $S_{i,j}^{\star}$ $(1\protect\leq i,j\protect\leq d)$
based on \textsf{HeteroPCA}.}

\label{alg:CE-HeteroPCA-CI}\begin{algorithmic}

\STATE \textbf{{Input}}: output $(\bm{U},\bm{\Sigma},\bm{S})$
of Algorithm \ref{alg:PCA-HeteroPCA}, sampling rate $p$, coverage
level $1-\alpha$.

\STATE \textbf{{Compute}} estimates of the noise level $\omega_{l}^{\star}$
as follows 
\[
\omega_{l}^{2}\coloneqq\frac{\sum_{j=1}^{n}y_{l,j}^{2}\ind_{(l,j)\in\Omega}}{\sum_{j=1}^{n}\ind_{(l,j)\in\Omega}}-S_{l,l}.
\]

\STATE \textbf{{Compute} }an estimate of $v_{i,j}^{\star}$ (cf.~(\ref{eq:ce-true-variance-i-j})
or (\ref{eq:ce-true-variance-i}))as follows: if $i\neq j$ then 
\begin{align*}
v_{i,j} & \coloneqq\frac{2-p}{np}S_{i,i}S_{j,j}+\frac{4-3p}{np}S_{i,j}^{2}+\frac{1}{np}\left(\omega_{i}^{2}S_{j,j}^{\star}+\omega_{j}^{2}S_{i,i}\right)\\
 & \quad+\frac{1}{np^{2}}\sum_{k=1}^{d}\left\{ \left[\omega_{i}^{2}+\left(1-p\right)S_{i,i}\right]\left[\omega_{k}^{2}+\left(1-p\right)S_{k,k}\right]+2\left(1-p\right)^{2}S_{i,k}^{2}\right\} \left(\bm{U}_{k,\cdot}\bm{U}_{j,\cdot}^{\top}\right)^{2}\\
 & \quad+\frac{1}{np^{2}}\sum_{k=1}^{d}\left\{ \left[\omega_{j}^{2}+\left(1-p\right)S_{j,j}\right]\left[\omega_{k}^{2}+\left(1-p\right)S_{k,k}\right]+2\left(1-p\right)^{2}S_{j,k}^{2}\right\} \left(\bm{U}_{k,\cdot}\bm{U}_{i,\cdot}^{\top}\right)^{2};
\end{align*}
If $i=j$ then 
\begin{align*}
v_{i,i} & \coloneqq\frac{12-9p}{np}S_{i,i}^{2}+\frac{4}{np}\omega_{i}^{2}S_{i,i}\\
 & \quad+\frac{4}{np^{2}}\sum_{k=1}^{d}\left\{ \left[\omega_{i}^{2}+\left(1-p\right)S_{i,i}\right]\left[\omega_{k}^{2}+\left(1-p\right)S_{k,k}\right]+2\left(1-p\right)^{2}S_{i,k}^{2}\right\} \left(\bm{U}_{k,\cdot}\bm{U}_{i,\cdot}^{\top}\right)^{2}.
\end{align*}

\STATE \textbf{{Output}} the $(1-\alpha)$-confidence interval
\[
\mathsf{CI}_{i,j}^{1-\alpha}\coloneqq\left[S_{i,j}\pm\Phi^{-1}\left(1-\alpha/2\right)\sqrt{v_{i,j}}\right].
\]

\end{algorithmic} 
\end{algorithm}

\begin{theorem}\label{thm:ce-CI}Suppose that the conditions in Theorem
\ref{thm:ce} hold. Assume that $ndp^{2}\gtrsim\log^{8}(n+d)$. Then
the confidence interval computed in Algorithm~\ref{alg:CE-HeteroPCA-CI}
obeys 
\[
\mathbb{P}\left(S_{i,j}^{\star}\in\mathsf{CI}_{i,j}^{1-\alpha}\right)=1-\alpha+o\left(1\right).
\]
\end{theorem}

Compared with \citet[Corollary 2]{cai2019subspace}, we can see that
when consistent estimation is possible --- namely, under the sampling
rate condition $p\geq\widetilde{\Omega}((n\land\sqrt{nd})^{-1})$
and the noise conditions $\omega_{\max}\leq\widetilde{\Omega}((\sqrt{n/d}\land\sqrt[4]{n/d})\sqrt{p}\sigma_{r}^{\star})$
--- it is plausible to construct fine-grained confidence interval
for $S_{i,j}^{\star}$, provided that the size of $\Vert\bm{U}_{i,\cdot}^{\star}\Vert_{2}+\Vert\bm{U}_{j,\cdot}^{\star}\Vert_{2}$
is not exceedingly small.

\begin{remark}
	Before concluding this subsection, 
	we note that the sampling rate $p$ might be unknown {\em a priori} in practice. 
	If this is the case, then 
	one plausible strategy is to replace $p$ in Algorithms~\ref{alg:PCA-HeteroPCA}, \ref{alg:PCA-HeteroPCA-CR} and \ref{alg:CE-HeteroPCA-CI} 
	with the following empirical estimate:  
	\[
		\widehat{p}=\frac{\sum_{l=1}^{d}\sum_{j=1}^{n}\ind\{(l,j)\in\Omega\}}{nd}.
	\] 
	In view of the standard concentration results $\widehat{p}=(1+o(1))p$,
	it is straightforward to verify that all of these inference procedure and the accompanying theory remain valid.  
	We omit the details for the sake of brevity.
\end{remark}

\subsection{A glimpse of key technical ingredients} \label{sec:technical-contribution}

Let us take a moment to highlight several technical ingredients of the current theory,
which might be applicable to other high-dimensional statistical problems beyond the analysis of \textsf{HeteroPCA}. 

\paragraph{Second-order perturbation theory for principal subspace.} 
At the core of our analysis lies a ``second-order'' perturbation theory tailored to general subspace estimation problems, 
to be presented in Section~\ref{sec:detour-subspace}. 
More concretely, we establish a second-order expansion of the subspace perturbation error (see Theorem \ref{thm:hpca_inference_general}) that makes explicit the following two parts: 
 (i) nearly tight first- and second-order terms, which can be expressed succinctly as linear and quadratic mappings of the perturbation matrix; 
 (ii) the remaining higher-order terms that are provably negligible. 
 Given that \textsf{HeteroPCA} is an iterative algorithm,  developing such a refined perturbation theory for \textsf{HeteroPCA} becomes substantially more challenging than the vanilla SVD-based approach. 
 Our refined perturbation theory allows us to tighten prior estimation theory (e.g., the Davis-Kahan $\sin\Theta$ Theorem \citep{davis1970rotation} or recent $\ell_{2,\infty}$-type perturbation bounds \citep{abbe2017entrywise,cai2019subspace,chen2020spectral}),
 the latter of which focused mainly on providing orderwise estimation error bounds.


\paragraph{Fine-grained distributional characterizations for the principal subspace $\bm{U}^\star$.}
As alluded to previously, 
we establish the distributional characterization for the principal subspace (i.e., Theorem~\ref{thm:pca}) based on a key error decomposition
\begin{equation} \label{eq:decomposition-U}
	\bm{U}\bm{R}-\bm{U}^{\star}=\underbrace{\left[\bm{E}\bm{X}^{\top}+\mathcal{P}_{\mathsf{off}\text{-}\mathsf{diag}}\left(\bm{E}\bm{E}^{\top}\right)\right]\bm{U}^{\star}\left(\bm{\Sigma}^{\star}\right)^{-2}}_{\eqqcolon\,\bm{Z}\,\,\text{(first- and second-order approximation)}}+\underbrace{\left[\bm{U}\bm{R}-\bm{U}^{\star}-\bm{Z}\right]}_{\eqqcolon\,\bm{\Psi}\,\,\text{(residual term)}},
\end{equation}
where $\bm{E}\coloneqq n^{-1/2}(p^{-1}\bm{Y}-\bm{X})$. For each $l\in[d]$, the multivariate Berry-Esseen Theorem reveals the approximate Gaussianity of $\bm{Z}_{l,\cdot}$, while at the same time, our second-order perturbation theory (cf.~Theorem \ref{thm:hpca_inference_general}) ensures that $\bm{\Psi}_{l,\cdot}$ is stochastically dominated by $\bm{Z}_{l,\cdot}$.
Additionally, rather than providing general $\ell_{2,\infty}$ bounds (as in the prior work \citet{cai2019subspace}), 
our proof relies crucially on more delicate row-dependent error control (so that the size of $\bm{\Psi}_{l,\cdot}$ is carefully bounded in accordance with the $l$-th row of $\bm{U}^\star$ and $\bm{N}$).


\paragraph{Entrywise distributional characterizations when estimating the covariance matrix $\bm{S}^\star$.} 
Moving one step further, we derive the following key error decomposition w.r.t.~the covariance matrix $\bm{S}^\star$:
\begin{equation} \label{eq:decomposition-S}
	\begin{aligned}
		\bm{S}-\bm{S}^\star=\underbrace{\bm{U}^\star \bm{\Sigma}^{\star2} \bm{Z}^\top+\bm{Z}\bm{\Sigma}^{\star2}\bm{U}^\star}_{\eqqcolon\,\bm{W}}+ \underbrace{n^{-1}\bm{X}\bm{X}^\top-\bm{S}^{\star}}_{\eqqcolon\,\bm{A}}+\underbrace{\left[\bm{S}-\bm{S}^\star-\bm{W}-\bm{A}\right]}_{\eqqcolon\,\bm{\Phi}\,\,\text{(residual term)}},
	\end{aligned}
\end{equation}
where $\bm{Z}$ is defined in \eqref{eq:decomposition-U} and approximately Gaussian.  
Here, $\bm{W}$ serves as the main component as induced by the subspace estimation error, 
$\bm{A}$ indicates the discrepancy between the empirical covariance (using clean and fully observed data) and the true covariance, 
whereas $\bm{\Phi}$ is some higher-order term that is provably negligible in a strong entrywise sense. 
This in turn allows us to pin down tight entrywise distributional characterizations  for $\bm{S}-\bm{S}^\star$.

\section{Numerical experiments\label{sec:Numerical-experiments}}

\paragraph*{Setup.}

This section conducts a series of numerical experiments to validate
our distributional and inference theory developed in Section \ref{sec:Distributional-theory-inference}.
Throughout this section, unless otherwise noted, we fix the
dimension to be $d=100$ and the number of sample vectors to be $n=2000$,
and we generate the covariance matrix as $\bm{S}^{\star}=\bm{U}^{\star}\bm{U}^{\star\top}$
with $\bm{U}^{\star}\in\mathbb{R}^{n\times r}$ being a random orthonormal
matrix following the Haar distribution over the Grassmann manifold $G_{d,r}$ \citep[Section 5.2.6]{vershynin2016high}. In each Monte Carlo trial, the observed data are produced
according to the model described in Section \ref{subsec:intro-model}.
For the purpose of introducing heteroskedasticity, we will introduce
a parameter $\omega^{\star}$ that controls the noise level: in each
independent trial, each noise level $\omega_{l}^{\star}$ ($1\leq l\leq d$)
is independently drawn from $\mathsf{Uniform}[0.1\omega^{\star},2\omega^{\star}]$;
the random noise component $\eta_{l,j}$ is then drawn from $\mathcal{N}(0,\omega_{l}^{\star2})$
independently for every $l\in[d]$ and $j\in[n]$.

\paragraph{Superiority of \textsf{HeteroPCA} to the SVD-based approach in estimation.}

\begin{figure}[t]
\centering

\begin{tabular}{cc}
\includegraphics[scale=0.35]{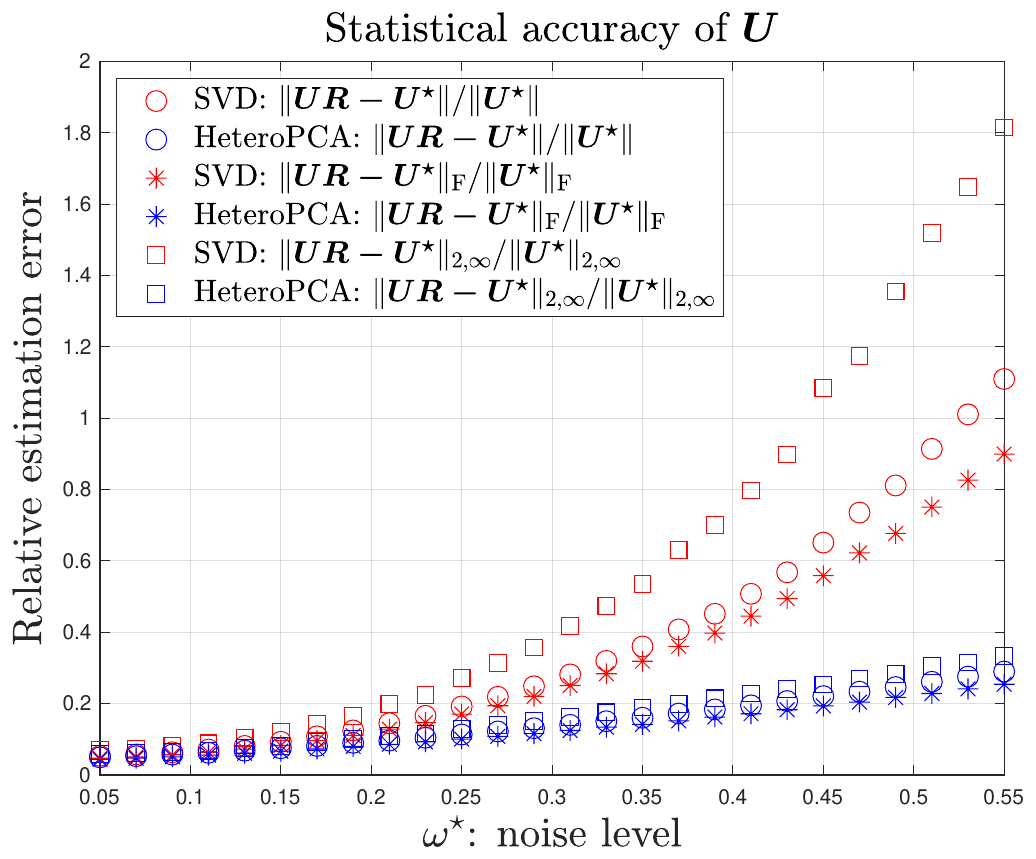}  & \includegraphics[scale=0.35]{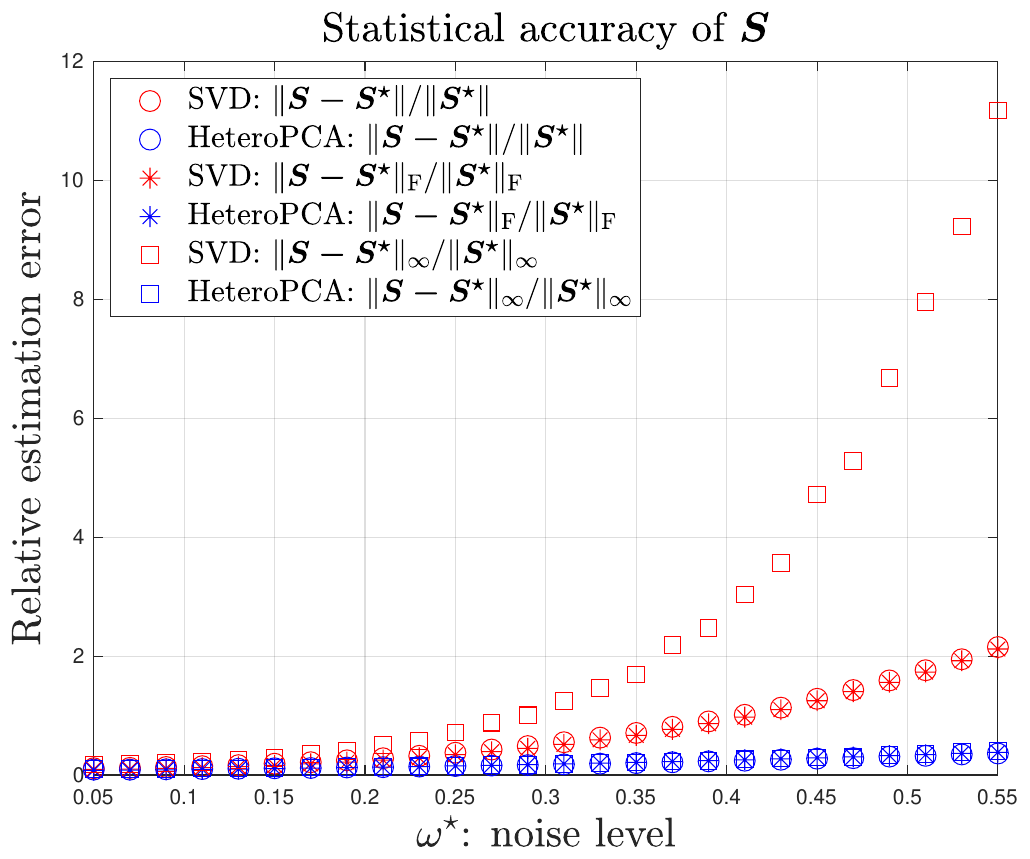}\tabularnewline
$\quad$(a)  & $\quad\quad$(b)\tabularnewline
\end{tabular}

\caption{The relative estimation error of $\bm{U}$ and $\bm{S}$ returned
by both SVD-based approach (cf.~Algorithm \ref{alg:PCA-SVD}) and
\textsf{HeteroPCA} (cf.~Algorithm \ref{alg:PCA-HeteroPCA}) over different noise level $\omega^{\star}$. (a)
Relative estimation errors of $\bm{U}\bm{R}-\bm{U}^{\star}$ measured
by $\Vert\cdot\Vert$, $\Vert\cdot\Vert_{\mathrm{F}}$ and $\Vert\cdot\Vert_{2,\infty}$
vs.~the noise level $\omega^{\star}$; (b) Relative estimation errors
of $\bm{S}-\bm{S}^{\star}$ measured by $\Vert\cdot\Vert$, $\Vert\cdot\Vert_{\mathrm{F}}$
and $\Vert\cdot\Vert_{\infty}$ vs.~the noise level $\omega^{\star}$.
The results are reported over $200$ independent trials for $r=3$
and $p=0.6$. \label{fig:comparison-omega}}
\end{figure}

\begin{figure}[t]
\centering

\begin{tabular}{cc}
\includegraphics[scale=0.35]{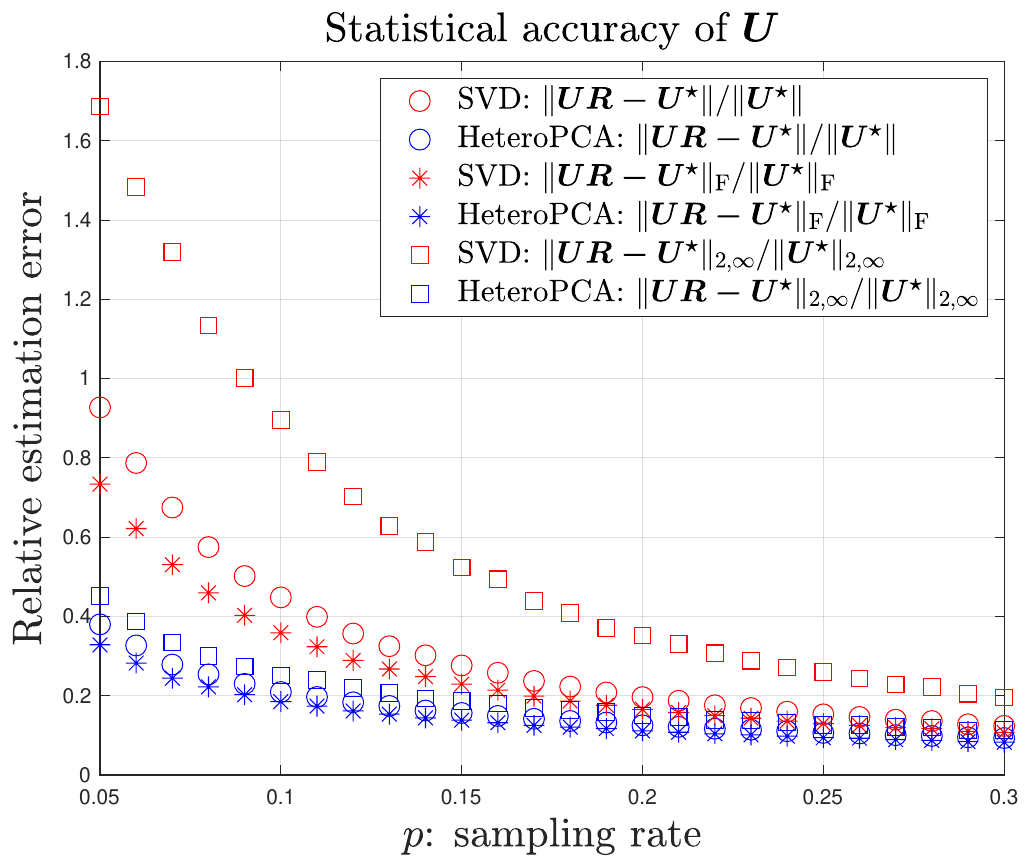}  & \includegraphics[scale=0.35]{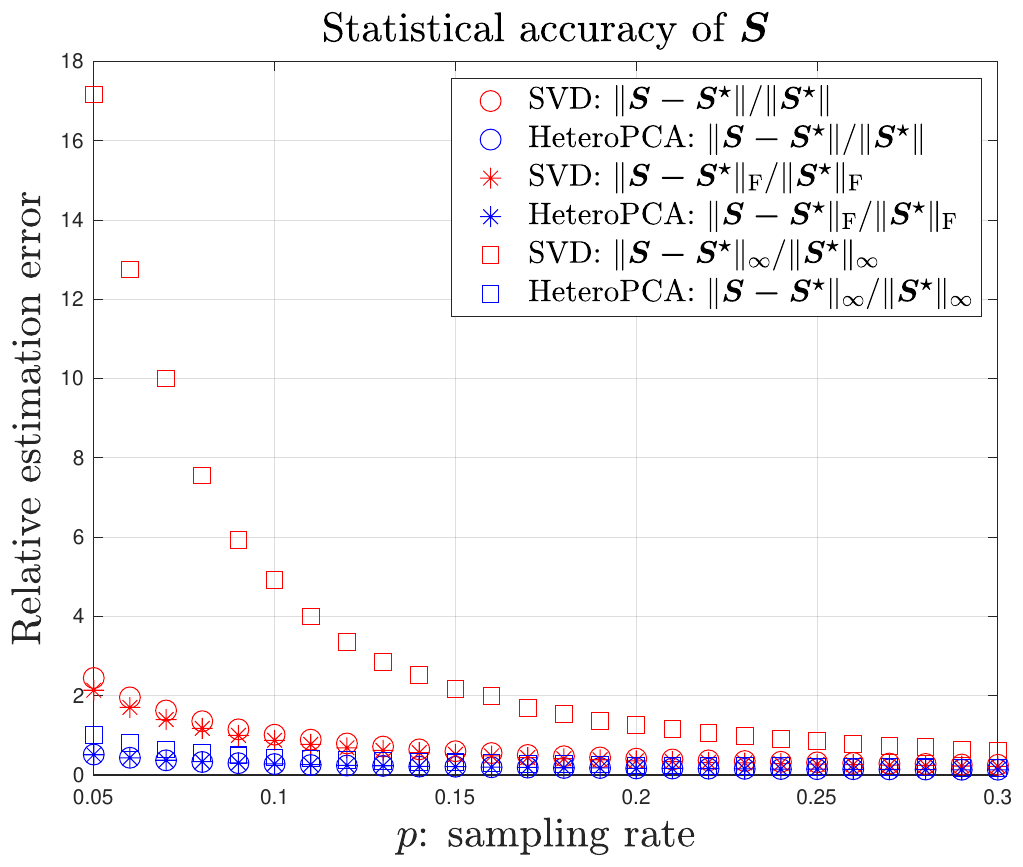}\tabularnewline
$\quad$(a)  & $\quad\quad$(b)\tabularnewline
\end{tabular}

\caption{The relative estimation error of $\bm{U}$ and $\bm{S}$ returned
by both SVD-based approach (cf.~Algorithm \ref{alg:PCA-SVD}) and
\textsf{HeteroPCA} (cf.~Algorithm \ref{alg:PCA-HeteroPCA}) across different missing probability $p$. (a)
Relative estimation errors of $\bm{U}\bm{R}-\bm{U}^{\star}$ measured
by $\Vert\cdot\Vert$, $\Vert\cdot\Vert_{\mathrm{F}}$ and $\Vert\cdot\Vert_{2,\infty}$
vs.~the missing rate $p$; (b) Relative estimation errors of $\bm{S}-\bm{S}^{\star}$
measured by $\Vert\cdot\Vert$, $\Vert\cdot\Vert_{\mathrm{F}}$ and
$\Vert\cdot\Vert_{\infty}$ vs.~the missing rate $p$. The results
are reported over $200$ independent trials for $r=3$ and $\omega^{\star}=0.05$.
\label{fig:comparison-p}}
\end{figure}

To begin with, we first compare the empirical estimation accuracy
of the SVD approach (cf.~Algorithm \ref{alg:PCA-SVD}) and \textsf{HeteroPCA}
(cf.~Algorithm \ref{alg:PCA-HeteroPCA}). Figure \ref{fig:comparison-omega}
displays the relative estimation errors --- including the ones tailored
to the principal subspace: $\Vert\bm{U}\bm{R}-\bm{U}^{\star}\Vert/\Vert\bm{U}^{\star}\Vert$,
$\Vert\bm{U}\bm{R}-\bm{U}^{\star}\Vert_{\mathrm{F}}/\Vert\bm{U}^{\star}\Vert_{2,\infty}$,
$\Vert\bm{U}\bm{R}-\bm{U}^{\star}\Vert_{2,\infty}/\Vert\bm{U}^{\star}\Vert_{2,\infty}$,
and the ones tailored to the covariance matrix: $\Vert\bm{S}-\bm{S}^{\star}\Vert/\Vert\bm{S}^{\star}\Vert$,
$\Vert\bm{S}-\bm{S}^{\star}\Vert_{\mathrm{F}}/\Vert\bm{S}^{\star}\Vert_{\mathrm{F}}$,
$\Vert\bm{S}-\bm{S}^{\star}\Vert_{\infty}/\Vert\bm{S}^{\star}\Vert_{\infty}$
--- of both algorithms as the noise level $\omega^{\star}$ varies,
with $r=3$ and $p=0.6$. Similarly, Figure \ref{fig:comparison-p}
shows the relative numerical estimation errors of both algorithms
vs.~the sampling rate $p$, with $r=3$ and $\omega^{\star}=0.05$.
As we shall see from both figures, \textsf{HeteroPCA} uniformly outperforms
the SVD-based approach in all experiments, and is able to achieve
appealing performance for a much wider range of noise levels and sampling
rates.

\begin{figure}[t]
\centering

\begin{tabular}{cc}
\includegraphics[scale=0.35]{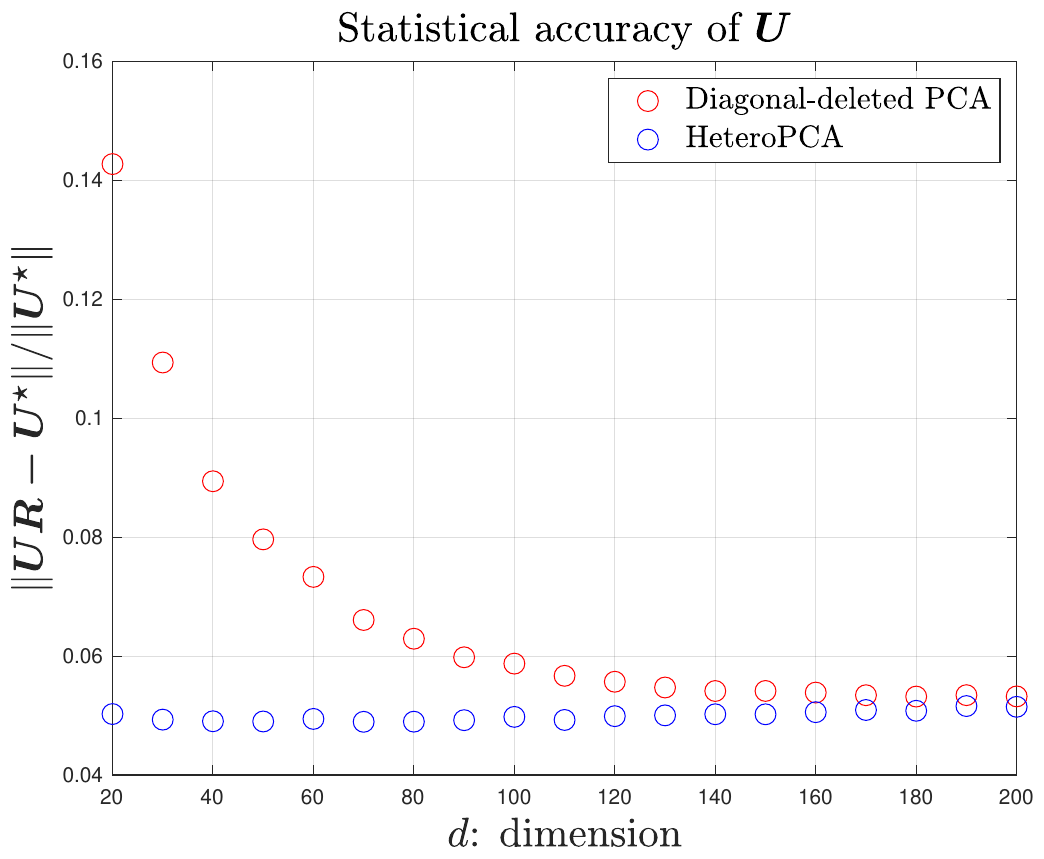} & \includegraphics[scale=0.35]{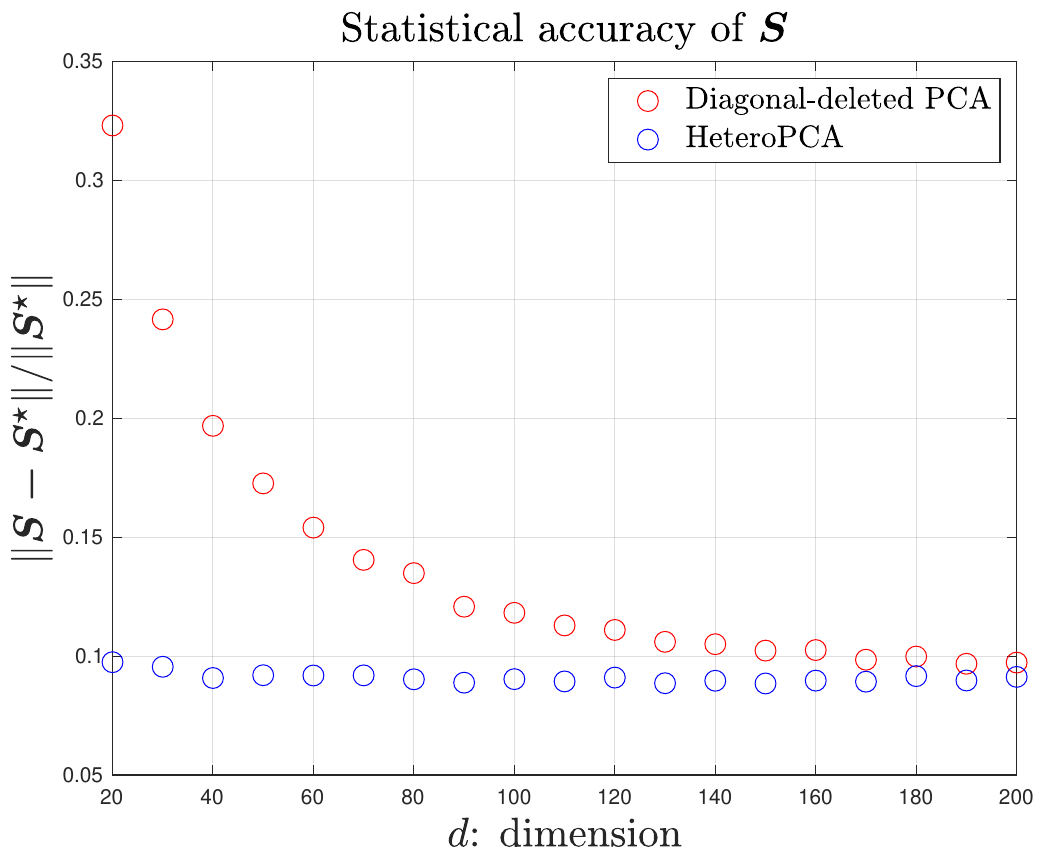}\tabularnewline
$\quad$(a) & $\quad\quad$(b)\tabularnewline
\end{tabular}

\caption{The relative estimation error of $\bm{U}$ and $\bm{S}$ returned
by both diagonal-deleted spectral method \citep{cai2019subspace}
and \textsf{HeteroPCA} (cf.~Algorithm \ref{alg:PCA-HeteroPCA}).
(a) Relative estimation error $\Vert\bm{U}\bm{R}-\bm{U}^{\star}\Vert/\Vert\bm{U}^{\star}\Vert$
vs.~dimension $d$; (b) Relative estimation error $\Vert\bm{S}-\bm{S}^{\star}\Vert/\Vert\bm{S}^{\star}\Vert$
vs.~the dimension $d$. The results are reported over $200$ independent
trials for $r=3$, $\omega^{\star}=0.05$ and $p=0.6$. \label{fig:comparison-dd}}
\end{figure}

\paragraph{Superiority of \textsf{HeteroPCA} to diagonal-deleted PCA in
estimation.}
Let us also compare the empirical estimation
accuracy of the diagonal-deleted spectral method \citep{cai2019subspace}
and \textsf{HeteroPCA} (cf.~Algorithm \ref{alg:PCA-HeteroPCA}).
Recall from Section~\ref{subsec:Inference-for-HeteroPCA} that the main difference between the
estimation error bounds of these two algorithms lies in an
additional bias term  due to the diagonal deletion operation (see the last term $\mu_{\mathsf{ce}}\kappa_{\mathsf{ce}}r/d$
in Equation~(4.16) in \citet{cai2019subspace}).
Figure~\ref{fig:comparison-dd} displays the relative estimation errors
for estimating the principal subspace $\Vert\bm{U}\bm{R}-\bm{U}^{\star}\Vert/\Vert\bm{U}^{\star}\Vert$
and for estimating the covariance matrix $\Vert\bm{S}-\bm{S}^{\star}\Vert/\Vert\bm{S}^{\star}\Vert$
as the dimension $d$ varies, with $r=3$, $\omega^{\star}=0.05$ and
$p=0.6$. As can be seen from the plots, \textsf{HeteroPCA} uniformly outperforms
the diagonal-deleted spectral method, especially when $d$ is not
too large. This numerical evidence corroborates 
the efficacy of the diagonal refinement scheme adopted in \textsf{HeteroPCA.}

\paragraph{Confidence regions for the principal subspace $\bm{U}^{\star}$.}

\begin{figure}[t]
\centering

\begin{tabular}{cc}
\includegraphics[scale=0.35]{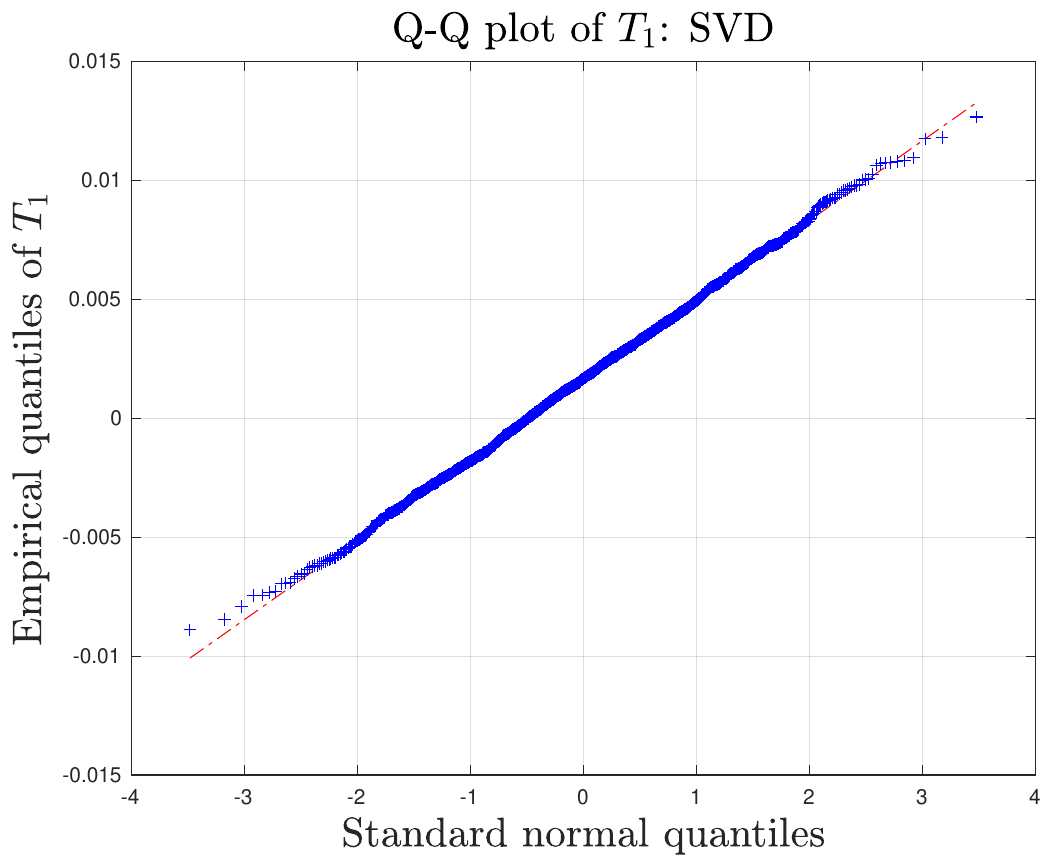}  & \includegraphics[scale=0.35]{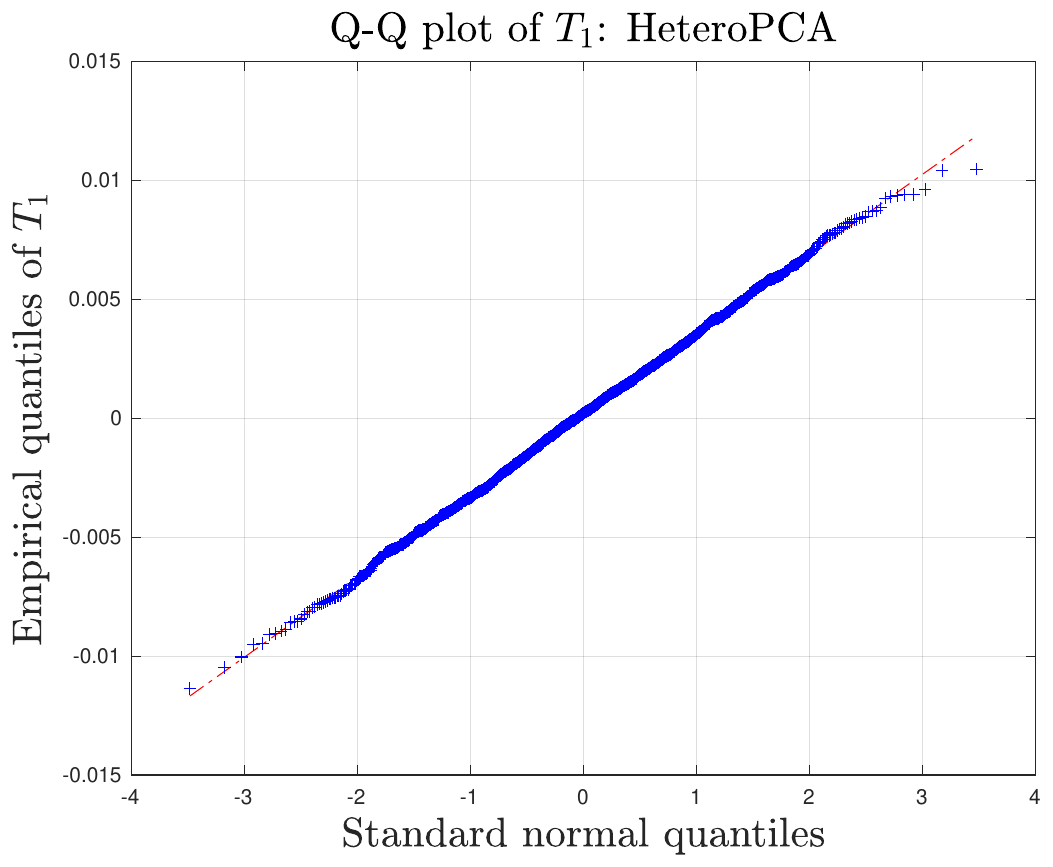}\tabularnewline
$\quad$(a)  & $\quad\quad$(b)\tabularnewline
\end{tabular}

\caption{(a) Q-Q (quantile-quantile) plot of $T_{1}$ vs.~the standard normal
distribution for the SVD-based approach; (b) Q-Q (quantile-quantile)
plot of $T_{1}$ vs.~the standard normal distribution for \textsf{HeteroPCA}.
The results are reported over $2000$ independent trials for $r=1$,
$p=0.6$ and $\omega^{\star}=0.05$. \label{fig:pca}}
\end{figure}

Next, we carry out a series of experiments to corroborate the practical
validity of the confidence regions constructed using the SVD-based
approach \citep[Algorithm 3]{yan2021inference} and \textsf{HeteroPCA}
(cf.~Algorithm \ref{alg:PCA-HeteroPCA-CR}). To this end, we define
$\widehat{\mathsf{Cov}}_{U}(i)$ to be the empirical probability that
the constructed confidence interval $\mathsf{CR}_{U,i}^{0.95}$ covers
$\bm{U}_{i,\cdot}^{\star}\mathsf{sgn}(\bm{U}^{\star\top}\bm{U})$
over 200 Monte Carlo trials, where $\bm{U}$ is the estimate returned
by either algorithm. We also let $\mathsf{Mean}(\widehat{\mathsf{Cov}}_{U})$
(resp.~$\mathsf{std}(\widehat{\mathsf{Cov}}_{U})$) be the empirical
mean (resp.~standard deviation) of $\widehat{\mathsf{Cov}}_{U}(i)$
over $i\in[d]$. Table \ref{table:pca} gathers $\mathsf{Mean}(\widehat{\mathsf{Cov}})$
and $\mathsf{std}(\widehat{\mathsf{Cov}})$ for $r=3$ and different
choices of $(p,\omega^{\star})$ for both algorithms. Encouragingly,
the empirical coverage rates are all close to $95\%$ for both methods
when $p$ is not too small and $\omega^{\star}$ is not too large.
When $p$ becomes smaller or $\omega^{\star}$ grows larger, \textsf{HeteroPCA}
is still capable of performing valid statistical inference, while
the SVD-based approach fails. This provides another empirical evidence
on the advantage and broader applicability of \textsf{HeteroPCA} compared
to the SVD-based approach. In addition, for the rank-1 case ($r=1$),
we define $T_{i}\coloneqq[\bm{U}-\mathsf{sign}(\bm{U}^{\top}\bm{U}^{\star})U^{\star}]_{i}/\sqrt{\Sigma_{U,i}}$.
Figure \ref{fig:pca} displays the Q-Q (quantile-quantile) plot of
$T_{1}\coloneqq[U-\mathsf{sign}(\bm{U}^{\top}\bm{U}^{\star})U^{\star}]_{1}/\sqrt{\Sigma_{U,1}}$
vs.~the standard Gaussian random variable over $2000$ Monte Carlo
simulations for both algorithms (when $p=0.6$ and $\omega^{\star}=0.05$);
the near-Gaussian empirical distribution of $T_{1}$ also corroborates
our distributional guarantees.

\begin{table}[b]
\caption{Empirical coverage rates of $\bm{U}^{\star}\mathsf{sgn}(\bm{U}^{\star\top}\bm{U})$
for different $(p,\omega^{\star})$'s over 200 Monte Carlo trials\label{table:pca}}

\centering

\begin{tabular}{c|c|c|c|c|c}
\hline 
 &  & \multicolumn{2}{c|}{The SVD-based Approach} & \multicolumn{2}{c}{\textsf{HeteroPCA}}\tabularnewline
\hline 
$p$  & $\omega^{\star}$  & $\mathsf{Mean}(\widehat{\mathsf{Cov}})$  & $\mathsf{Std}(\widehat{\mathsf{Cov}})$  & $\mathsf{Mean}(\widehat{\mathsf{Cov}})$  & $\mathsf{Std}(\widehat{\mathsf{Cov}})$\tabularnewline
\hline 
$0.6$  & $0.05$  & $0.9270$  & $0.0292$  & $0.9523$  & $0.0157$\tabularnewline
$0.6$  & $0.1$  & $0.8989$  & $0.0521$  & $0.9484$  & $0.0154$\tabularnewline
\hline 
$0.4$  & $0.05$  & $0.8849$  & $0.0501$  & $0.9448$  & $0.0184$\tabularnewline
$0.4$  & $0.1$  & $0.8458$  & $0.0853$  & $0.9405$  & $0.0182$\tabularnewline
\hline 
$0.2$  & $0.05$  & $0.7370$  & $0.1196$  & $0.9287$  & $0.0204$\tabularnewline
$0.2$  & $0.1$  & $0.6856$  & $0.1569$  & $0.9219$  & $0.0204$\tabularnewline
\hline 
\end{tabular}
\end{table}

\begin{table}[b]
\caption{Empirical coverage rates of $S_{i,j}^{\star}$ for different $(\omega^{\star},p)$'s
over 200 Monte Carlo trials\label{table:ce}}

\centering

\begin{tabular}{c|c|c|c|c|c}
\hline 
 \multicolumn{2}{c|}{ }  & \multicolumn{2}{c|}{The SVD-based Approach} & \multicolumn{2}{c}{\textsf{HeteroPCA}}\tabularnewline
\hline 
$p$  & $\omega^{\star}$  & $\mathsf{Mean}(\widehat{\mathsf{Cov}})$  & $\mathsf{Std}(\widehat{\mathsf{Cov}})$  & $\mathsf{Mean}(\widehat{\mathsf{Cov}})$  & $\mathsf{Std}(\widehat{\mathsf{Cov}})$\tabularnewline
\hline 
$0.6$  & $0.05$  & $0.9380$  & $0.0244$  & $0.9475$  & $0.0153$\tabularnewline
$0.6$  & $0.1$  & $0.9243$  & $0.0425$  & $0.9484$  & $0.0151$\tabularnewline
\hline 
$0.4$  & $0.05$  & $0.9200$  & $0.0509$  & $0.9485$  & $0.0156$\tabularnewline
$0.4$  & $0.1$  & $0.9027$  & $0.0713$  & $0.9490$  & $0.0153$\tabularnewline
\hline 
$0.2$  & $0.05$  & $0.8657$  & $0.1031$  & $0.9494$  & $0.0164$\tabularnewline
$0.2$  & $0.1$  & $0.8488$  & $0.1186$  & $0.9491$  & $0.0162$\tabularnewline
\hline 
\end{tabular}
\end{table}

\begin{figure}[t]
\centering

\begin{tabular}{cc}
\includegraphics[scale=0.35]{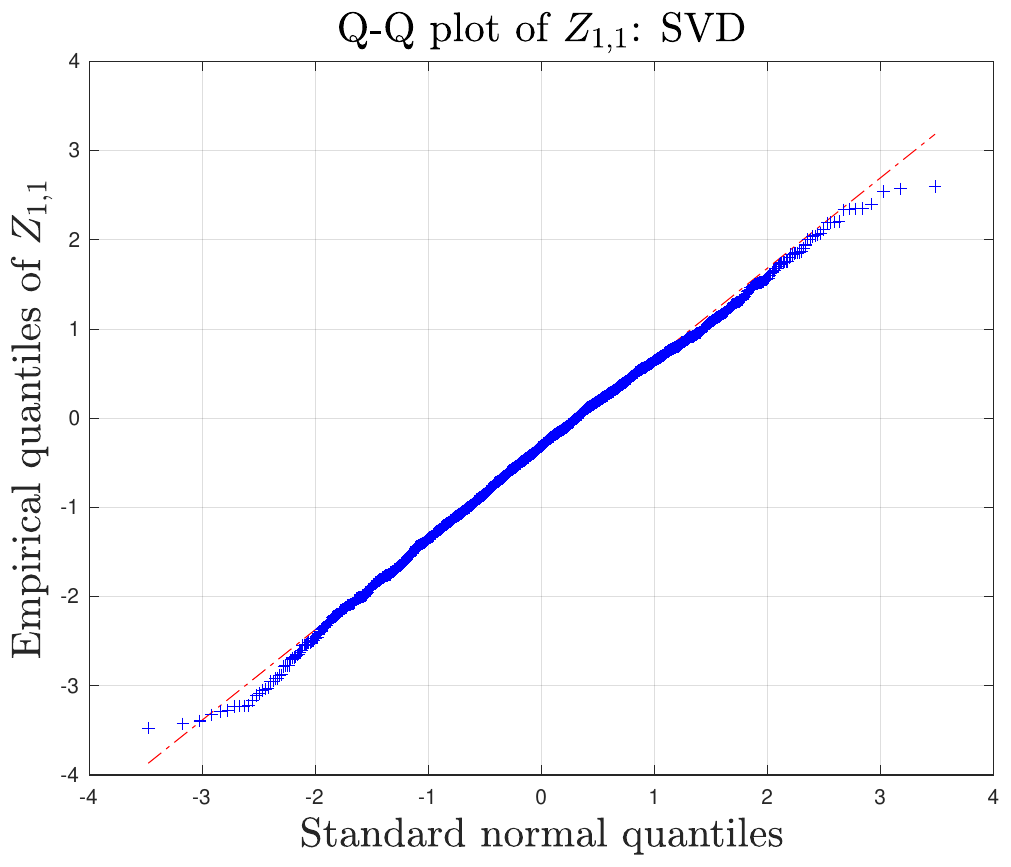}  & \includegraphics[scale=0.35]{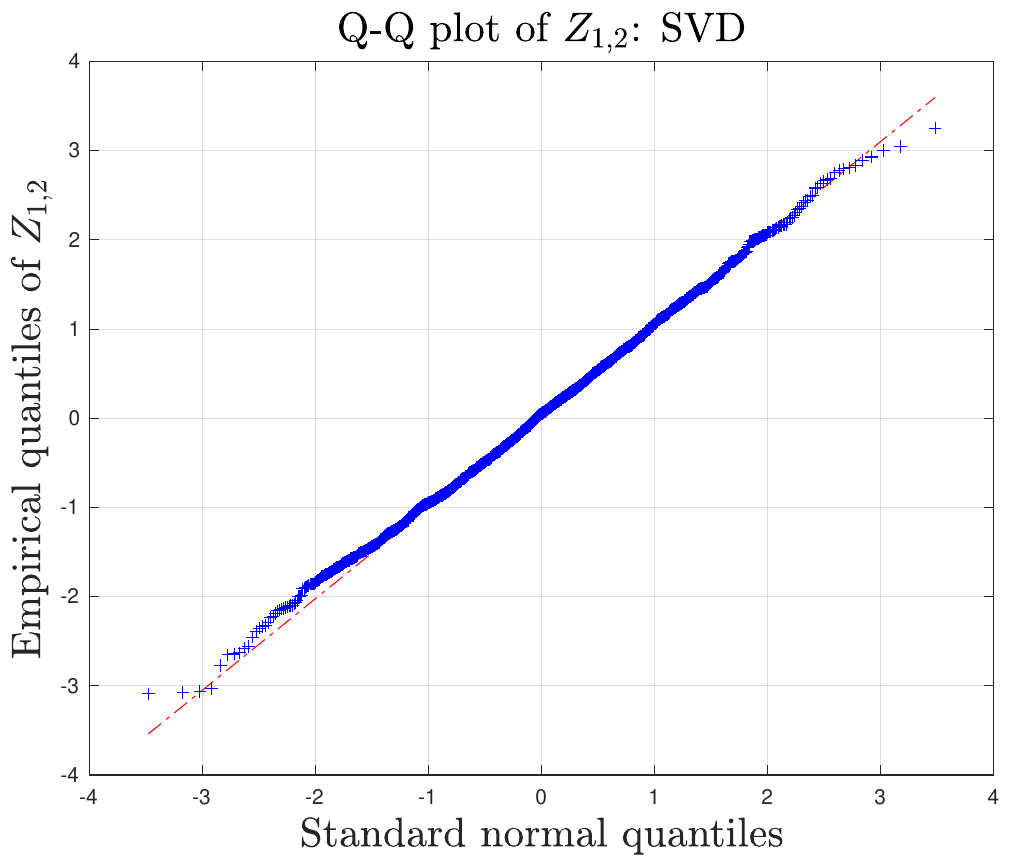}\tabularnewline
$\quad$(a)  & $\quad\quad$(b)\tabularnewline
\end{tabular}

\caption{(a) Q-Q (quantile-quantile) plot of $Z_{1,1}$ vs.~the standard normal
distribution for the SVD-based approach; (b) Q-Q (quantile-quantile)
plot of $Z_{1,2}$ vs.~a standard Gaussian distribution for the SVD-based
approach. The results are reported over $2000$ independent trials
for $r=3$, $p=0.6$, $\omega^{\star}=0.05$. \label{fig:ce-svd}}
\end{figure}

\begin{figure}[t]
\centering

\begin{tabular}{cc}
\includegraphics[scale=0.35]{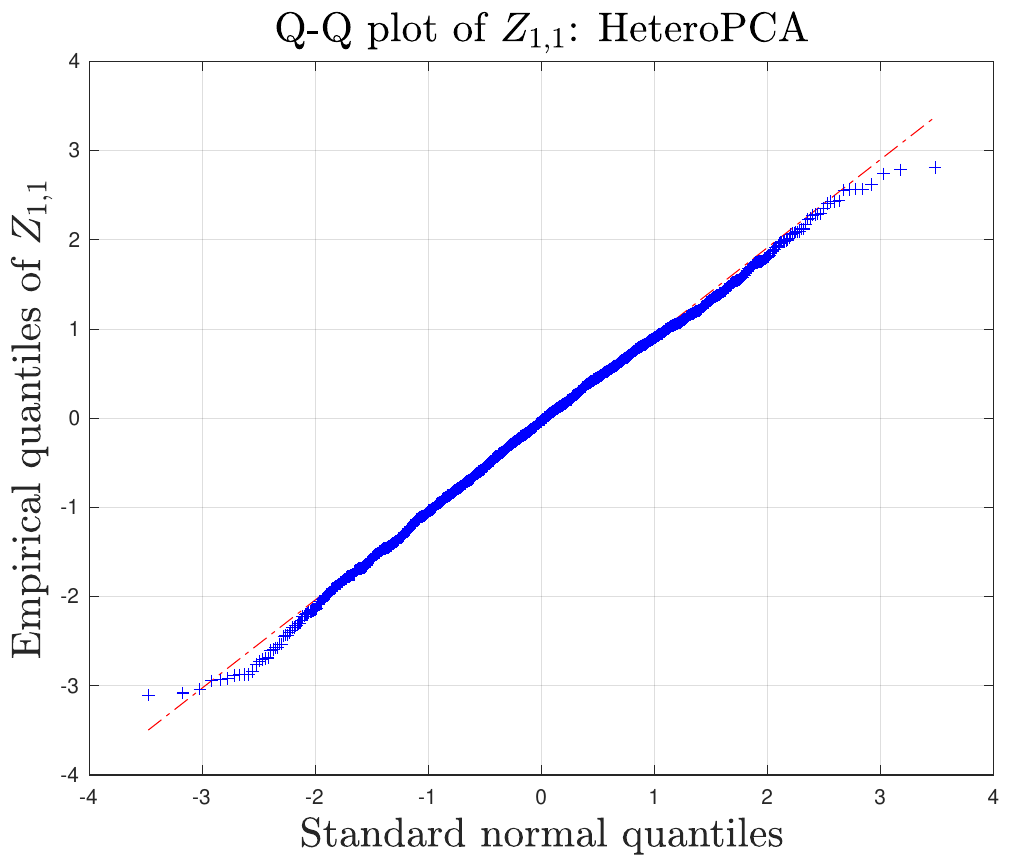}  & \includegraphics[scale=0.35]{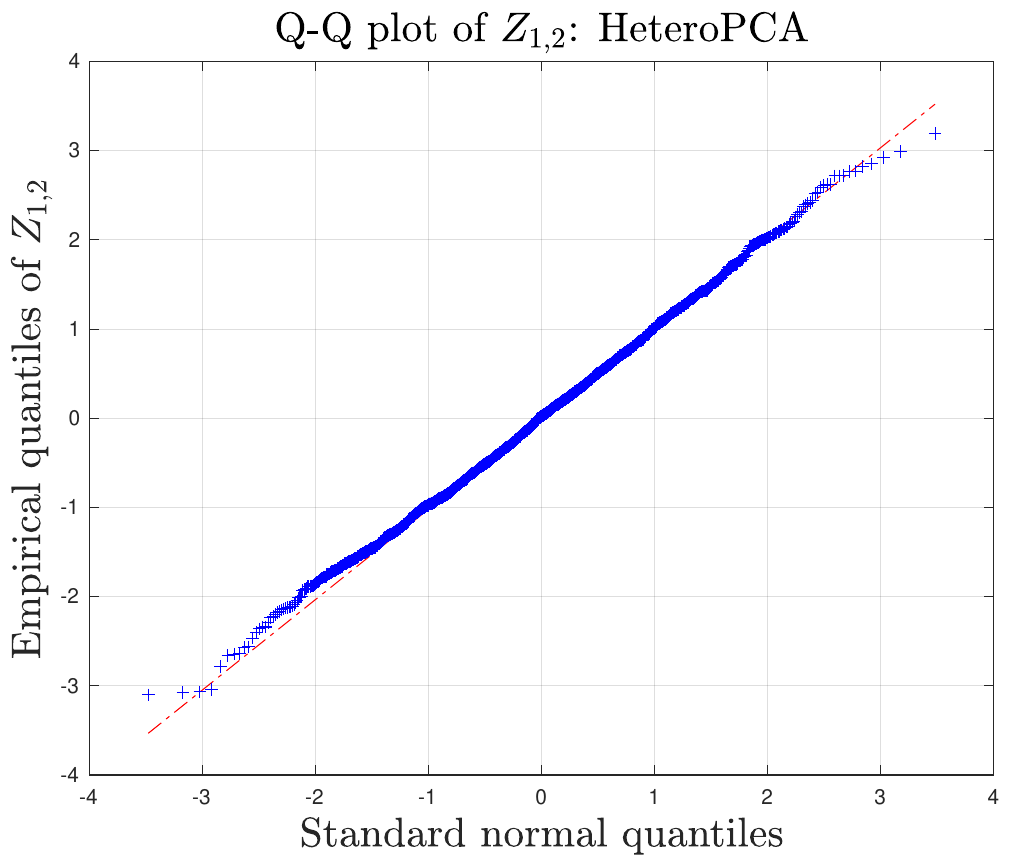}\tabularnewline
$\quad$(a)  & $\quad\quad$(b)\tabularnewline
\end{tabular}

\caption{(a) Q-Q (quantile-quantile) plot of $Z_{1,1}$ vs.~the standard normal
distribution for \textsf{HeteroPCA}; (b) Q-Q (quantile-quantile) plot
of $Z_{1,2}$ vs.~a standard Gaussian distribution for \textsf{HeteroPCA}.
The results are reported over $2000$ independent trials for $r=3$,
$p=0.6$, $\omega^{\star}=0.05$. \label{fig:ce-hpca}}
\end{figure}

\paragraph{Entrywise confidence intervals for $\bm{S}^{\star}$.}

Finally, we provide numerical evidence that confirms the validity
of the confidence interval constructed on the basis of the SVD-based
approach \citep[Algorithm 4]{yan2021inference} and \textsf{HeteroPCA}
(cf.~Algorithm \ref{alg:CE-HeteroPCA-CI}). Define $\widehat{\mathsf{Cov}}_{S}(i,j)$
to be the empirical probability that the $95\%$ confidence interval
$[S_{i,j}\pm1.96\sqrt{v_{i,j}}]$ covers $S_{i,j}^{\star}$ over 200
Monte Carlo trials, where $S_{i,j}$ is the $(i,j)$-th entry of the
estimate $\bm{S}$ returned by either algorithm. Let $\mathsf{Mean}(\widehat{\mathsf{Cov}}_{S})$
(resp.~$\mathsf{std}(\widehat{\mathsf{Cov}}_{S})$) be the empirical
mean (resp.~standard deviation) of $\widehat{\mathsf{Cov}}_{S}(i,j)$
over all $i,j\in[d]$. Table \ref{table:ce} collects $\mathsf{Mean}(\widehat{\mathsf{Cov}})$
and $\mathsf{std}(\widehat{\mathsf{Cov}})$ for $r=3$ and accounts
for different choices of $(p,\omega^{\star})$ for both algorithms.
Similar to previous experiments, \textsf{HeteroPCA} uniformly outperforms
the SVD-based approach, which again suggests that \textsf{HeteroPCA}
is the method of choice. In addition, we define $Z_{i,j}\coloneqq(S_{i,j}-S_{i,j}^{\star})/\sqrt{v_{i,j}}$.
For both algorithms, Figure \ref{fig:ce-svd} and Figure \ref{fig:ce-hpca}
depict the Q-Q (quantile-quantile) plot of $Z_{1,1}$ and $Z_{1,2}$
vs.~standard Gaussian distributions over $2000$ Monte Carlo trials
for the case with $r=3$, $p=0.6$ and $\omega^{\star}=0.05$, which
again confirm the practical validity of our distributional theory. 

\section{Other related works\label{sec:Related-works}}

Low-rank matrix denoising serves as a common model to study the effectiveness
of spectral methods \citep{chen2020spectral}, and has been the main
subject of many prior works including \citet{cai2018rate,chen2018asymmetry,ding2020high,abbe2017entrywise,bao2018singular,montanari2018adapting,lei2019unified,xia2019normal,cape2019two,agterberg2021entrywise},
among others. Several recent works began to pursue a distributional
theory for the eigenvector or singular vectors of the observed data
matrix \citep{fan2020asymptotic,cheng2020tackling,bao2018singular,xia2019normal}.
To name a few examples, \citet{bao2022statistical} studied the limiting
distribution of the inner product between an empirical singular vector
and the corresponding ground truth, assuming that the associated spectral
gap is sufficient large and that the noise components are homoskedastic;
\citet{xia2019normal} established non-asymptotic Gaussian approximation
for certain projection distance in the presence of i.i.d.~Gaussian
noise. Furthermore, the presence of missing data forms another source
of technical challenges, leading to a problem often dubbed as noisy
low-rank matrix completion \citep{CanPla10,Negahban2012restricted,chen2019noisy}.

Spectral methods have been successfully applied to tackle noisy matrix
completion \citep{KesMonSew2010,sun2016guaranteed,chen2015fast,zheng2016convergence,ma2017implicit,chen2019nonconvex,cho2017asymptotic},
which commonly serve as an effective initialization scheme for nonconvex
optimization methods \citep{chi2018nonconvex}. While statistical
inference for noisy matrix completion has been investigated recently
\citep{chen2019inference,xia2021statistical,chernozhukov2021inference},
these prior works focused on performing inference based on optimization-based
estimators. How to construct fine-grained confidence intervals based
on spectral methods remains previously out of reach for noisy matrix
completion. It is also noteworthy that the inferential procedures proposed in \citet{chen2019inference,xia2021statistical} (for noisy matrix completion)
were developed for the regime where reliable estimation of the full low-rank matrix is feasible. 
This, however, falls short of covering the most challenging regime considered herein (where one might only be able to estimate the column subspace but not the row subspace).  This crucial difference in the regimes of interest leads to substantial challenges unaddressed by these prior works.

Additionally, the recent work \citep{xia2018confidence} tackled the
confidence regions for spectral estimators tailored to the low-rank
matrix regression problem, without accommodating the noisy matrix
completion context. Most importantly, while the SVD-based vanilla
spectral method often works well for the balanced case (such that
the column dimension and the row dimension are on the same order),
sub-optimality has been well recognized when estimating the column
subspace of interest in the highly unbalanced case (so that the column
dimension far exceeds the row dimension); this issue is also present
when it comes to existing optimization-based methods like nuclear
norm minimization. As a result, all prior schemes mentioned in this
paragraph failed to tackle the highly balanced case in an statistically
efficient manner.

Turning to PCA or subspace estimation, there has been an enormous
literature dedicated to this topic; see \citet{johnstone2018pca,balzano2018streaming}
for an overview of prior development. Noteworthily, the need to handle
the diagonals of the sample covariance matrix in the presence of heteroskedastic
noise and/or missing data has been pointed out in many prior works,
e.g., \citet{loh2012high,lounici2014high,montanari2018spectral,florescu2016spectral,cai2019subspace}.
The iterative refinement scheme proposed by \citet{zhang2018heteroskedastic}
turns out to be among the most effective and adaptive schemes in handling
the diagonals. Aimed at designing fine-grained estimators for the
principal components, \citet{koltchinskii2020efficient,li2021minimax}
proposed statistically efficient de-biased estimators for linear functionals
of principal components, and moreover, the estimator proposed in \citet{koltchinskii2020efficient}
has also been shown to exhibit asymptotic normality in the presence
of i.i.d.~Gaussian noise. \citet{bloemendal2016principal} also pinned
down the asymptotic distributions of certain principal components
under a spiked covariance model. However, these papers fell short
of presenting valid and data-driven uncertainty quantification methods
for the proposed estimators, and their results operates under the
assumptions of homoskedastic noise without any missing data, a scenario
that is remarkably more restricted than ours. Under the spiked covariance
model, \citet{bao2022statistical} studied the limiting distribution
of the angle between the eigenvectors of the sample covariance matrix
and any fixed vector, under the ``balanced'' scenario where the
aspect ratio $n/d$ is a constant. In addition, recent years have
witnessed much activity in high-dimensional PCA in the face of missing
data \citep{zhang2018heteroskedastic,zhu2019high,cai2019subspace,pavez2020covariance};
these works, however, focused primarily on developing estimation guarantees,
which did not provide either distributional guarantees for the estimators
or concrete procedures that allow for confidence region construction.
Additionally, the \textsf{HeteroPCA} algorithm has been further extended by two follow-up works \citet{zhou2023deflated,zhou2023heteroskedastic} to accommodate 
the scenario with large condition numbers as well as tensor clustering in the presence of heteroskedastic noise.

From a technical viewpoint, it is worth mentioning that the $\ell_{\infty}$
and $\ell_{2,\infty}$ perturbation theory has been an active research
direction in recent years \citep{fan2018eigenvector,cape2019two,chen2018asymmetry,eldridge2018unperturbed,agterberg2021entrywise,xie2021entrywise}.
Among multiple existing technical frameworks, the leave-one-out analysis
idea --- which has been applied to a variety of statistical estimation
problems \citep{el2013robust,el2015impact,zhong2017near,chen2018gradient,chen2020bridging,cai2019nonconvex,ling2020near,cai2020uncertainty,chen2020partial,chen2020convex}
--- provides a powerful and flexible framework that enables $\ell_{\infty}$
and $\ell_{2,\infty}$ statistical guarantees for spectral methods
\citep{chen2017spectral,abbe2017entrywise,cai2019subspace}; see \citep[Chapter 4]{chen2020spectral}
for an accessible introduction of this powerful framework. Our analysis for the
\textsf{HeteroPCA} approach is influenced by the one in \citet{cai2019subspace}.
Note, however, that \citet{cai2019subspace} didn't come with any distributional
guarantees for spectral methods, which we seek to accomplish in the current version of this paper.

It is important to note that although the current version of this paper focuses primarily on the \textsf{HeteroPCA} method, a preliminary version available on arXiv \citep{yan2021inference} includes a discussion on distributional theory and inferential procedures for PCA using the SVD-based approach (cf.~Algorithm \ref{alg:PCA-SVD}). This content was subsequently omitted during the revision phase based on editiorial suggestions. Interested readers are referred to \citet{yan2021inference} for a set of inferential results developed for the SVD-based approach, in parallel to Theorems \ref{thm:pca} to \ref{thm:ce-CI} in this paper.

Finally, we note in passing that constructing confidence intervals
for sparse regression (based on, say, the Lasso estimator or other
sparsity-promoting estimator), has attracted a flurry of research
activity in the past few years \citep{celentano2020lasso,zhang2014confidence,van2014asymptotically,javanmard2014confidence,cai2017confidence,ning2017general,ren2015asymptotic}.
The methods derived therein, however, are not directly applicable
to perform statistical inference for PCA and/or other low-rank models.

\section{A detour: subspace estimation\label{sec:detour-subspace}}

We now take a detour to look at an intimately related problem, which
we shall refer to as \emph{subspace estimation} and will play a crucial
role in understanding the \textsf{HeteroPCA} approach. We will set
out to develop a fine-grained statistical theory for \textsf{HeteroPCA}
when applied to this subspace estimation setting. The resulting theory
will be invoked in Appendix~\ref{sec:Analysis:-the-approach-heteroPCA}
to analyze the PCA context. 

\subsection{Model and algorithm\label{sec:subspace-setting}}

\noindent\paragraph{Model and assumptions.}

Suppose that we are interested in a rank-$r$ matrix $\bm{M}^{\natural}\in\mathbb{R}^{n_{1}\times n_{2}}$,
whose SVD is given by 
\begin{equation}
\bm{M}^{\natural}=\sum_{i=1}^{r}\sigma_{i}^{\natural}\bm{u}_{i}^{\natural}\bm{v}_{i}^{\natural\top}=\bm{U}^{\natural}\bm{\Sigma}^{\natural}\bm{V}^{\natural\top}\in\mathbb{R}^{n_{1}\times n_{2}}\label{eq:Astar-definition-subspace}
\end{equation}
Here, $\bm{U}^{\natural}=[\bm{u}_{1}^{\natural},\ldots,\bm{u}_{r}^{\natural}]$
(resp.~$\bm{V}^{\natural}=[\bm{v}_{1}^{\natural},\ldots,\bm{v}_{r}^{\natural}]$)
consists of orthonormal columns that correspond to the left (resp.~right)
singular vectors of $\bm{M}^{\natural}$, and $\bm{\Sigma}^{\natural}=\mathsf{diag}\{\sigma_{1}^{\natural},\ldots,\sigma_{r}^{\natural}\}$
is a diagonal matrix consisting of the singular values of $\bm{M}^{\natural}$.
Without loss of generality, we assume that 
\[
n=\max\{n_{1},n_{2}\}.
\]
It is assumed that the singular values are sorted (in magnitude) in
descending order, namely, 
\begin{equation}
\sigma_{1}^{\natural}\geq\cdots\geq\sigma_{r}^{\natural}\geq0,\label{eq:singular-value-sorted-subspace}
\end{equation}
with the condition number denoted by 
\begin{equation}
\kappa^{\natural}\coloneqq\sigma_{1}^{\natural}/\sigma_{r}^{\natural}.\label{eq:condition-number-subspace}
\end{equation}
What we have observed is a noisy copy of $\bm{M}^{\natural}$, namely,
\begin{equation}
\bm{M}=\bm{M}^{\natural}+\bm{E},\label{eq:observations-M-subspace}
\end{equation}
where $\bm{E}=[E_{i,j}]_{1\leq i,j\leq n}$ stands for a noise matrix.
We focus on \emph{estimating the column subspace} represented by $\bm{U}^{\natural}$
and the singular values encapsulated in $\bm{\Sigma}^{\natural}$,
but not the row space $\bm{V}^{\natural}$. An important special scenario
one should bear in mind is the highly unbalanced case where the column
dimension $n_{2}$ far exceeds the row dimension $n_{1}$; in this
case, it is common to encounter situations where reliable estimation
of $\bm{M}^{\natural}$ and $\bm{V}^{\natural}$ is infeasible but
that of $\bm{U}^{\natural}$ shows promise. For this reason, we refer
to this setting as \emph{subspace estimation} in order to differentiate
it from matrix denoising, emphasizing that we are only interested
in column subspace estimation.

With the new aim in mind, we shall modify our incoherence and noise
assumptions accordingly. Here, we abuse the notation with the understanding
that the following set of assumptions will be used only when analyzing
the approach based on \textsf{HeteroPCA. }We shall also denote $n\coloneqq\max\{n_{1},n_{2}\}$.

\begin{assumption}[\bf Incoherence]\label{assumption:subspace-incoherence}The
rank-$r$ matrix $\bm{M}^{\natural}\in\mathbb{R}^{n_{1}\times n_{2}}$
defined in (\ref{eq:Astar-definition-subspace}) is said to be $\mu^{\natural}$-incoherent
if the following holds: 
\begin{align*}
\left\Vert \bm{U}^{\natural}\right\Vert _{2,\infty} & \leq\sqrt{\frac{\mu^{\natural}r}{n_{1}}},\ \ \left\Vert \bm{V}^{\natural}\right\Vert _{2,\infty}\leq\sqrt{\frac{\mu^{\natural}r}{n_{2}}},\ \ \text{and}\ \left\Vert \bm{M}^{\natural}\right\Vert _{\infty}\leq\sqrt{\frac{\mu^{\natural}}{n_{1}n_{2}}}\left\Vert \bm{M}^{\natural}\right\Vert _{\mathrm{F}}.
\end{align*}
\end{assumption}

\begin{assumption}[\bf Heteroskedastic random noise]\label{assumption:subspace-noise}Assume
that the $E_{i,j}$'s are independently generated, and suppose that
there exist non-negative quantities $\{\sigma_{i}\}_{i=1}^{n_{1}}$,
$\{B_{i}\}_{i=1}^{n_{1}}$, $\sigma$ and $B$ obeying 
\[
\forall(i,j)\in[n_{1}]\times[n_{2}]:\qquad\mathbb{E}\left[E_{i,j}\right]=0,\qquad\mathsf{var}\left(E_{i,j}^{2}\right)=\sigma_{i,j}^{2}\leq\sigma_{i}^{2}\leq\sigma^{2},\qquad\left|E_{i,j}\right|\leq B_{i}\leq B,
\]
where for all $i\in[n_1]$
\begin{equation}
\qquad B_i\lesssim\frac{\sigma_i\min\left\{ \sqrt{n_{2}},\sqrt[4]{n_{1}n_{2}}\right\} }{\sqrt{\log n}},\qquad\text{and}\qquad B\lesssim\frac{\sigma\min\left\{ \sqrt{n_{2}},\sqrt[4]{n_{1}n_{2}}\right\} }{\sqrt{\log n}}.\label{eq:B-UB-subspace}
\end{equation}
\end{assumption}

\paragraph{Algorithm: \textsf{HeteroPCA} for subspace estimation.}

The paradigm \textsf{HeteroPCA} can naturally be applied to tackle
the above subspace estimation task. Let us introduce the ground-truth
gram matrix as follows 
\begin{equation}
\bm{G}^{\natural}\coloneqq\bm{M}^{\natural}\bm{M}^{\natural\top}.\label{eq:Gstar-Astar-subspace}
\end{equation}
Given that $\bm{M}=\bm{M}^{\natural}+\bm{E}$ is an unbiased estimate
of $\bm{M}^{\natural}$, one might naturally attempt to estimate the
column space of $\bm{M}$ by looking at the eigenspace of the sample
Gram matrix $\bm{M}\bm{M}^{\top}$. It can be easily seen that 
\begin{equation}
\bm{E}\left[\bm{M}\bm{M}^{\top}\right]=\bm{M}^{\natural}\bm{M}^{\natural\top}+\mathsf{diag}\left\{ \left[\sum_{j=1}^{n_{2}}\sigma_{i,j}^{2}\right]_{1\leq i\leq n_{1}}\right\} ,\label{eq:E-gram-decompose}
\end{equation}
where the diagonal term on the right-hand side of (\ref{eq:E-gram-decompose})
might incur significant bias in the most challenging regime. The \textsf{HeteroPCA}
algorithm seeks to handle the diagonal part in an iterative manner,
alternating between imputing the values of the diagonal entries and
eigen-decomposition of $\bm{M}\bm{M}^{\top}$ with the diagonal replaced
by the imputed values. The procedure is summarized in Algorithm \ref{alg:HeteroPCA-subspace}.

\begin{algorithm}[h]
\caption{HeteroPCA for general subspace estimation (\textsf{HeteroPCA}).}

\label{alg:HeteroPCA-subspace}\begin{algorithmic}

\STATE \textbf{{Initialization}}: set $\bm{G}^{0}=\mathcal{P}_{\mathsf{off}\text{-}\mathsf{diag}}\left(\bm{M}\bm{M}^{\top}\right)$.  
\STATE \textbf{{Updates}}: \textbf{for }$t=0,1,\ldots,t_{0}$ \textbf{do}

\STATE \vspace{-1em}
 \begin{subequations}\label{subeq:HeteroPCA_update-subspace} 
\begin{align}
\left(\bm{U}^{t},\bm{\Lambda}^{t}\right) & =\mathsf{eigs}\left(\bm{G}^{t},r\right);\label{eq:HeteroPCA_eig-subspace}\\
\bm{G}^{t+1} & =\mathcal{P}_{\mathsf{off}\text{-}\mathsf{diag}}\left(\bm{M}\bm{M}^{\top}\right)+\mathcal{P}_{\mathsf{diag}}\left(\bm{U}^{t}\bm{\Lambda}^{t}\bm{U}^{t\top}\right).\label{eq:HeteroPCA_imputation-subspace}
\end{align}
\end{subequations} Here, $\mathsf{eigs}(\bm{G},r)$ returns $(\bm{U},\bm{\Lambda})$
where $\bm{U}\bm{\Lambda}\bm{U}^{\top}$ is the top-$r$ eigen-decomposition
of $\bm{G}$.

\STATE \textbf{{Output}}: $\bm{U}=\bm{U}^{t_{0}}$, $\bm{\Lambda}=\bm{\Lambda}^{t_{0}}$,
$\bm{\Sigma}=(\bm{\Lambda}^{t_{0}})^{1/2}$, $\bm{S}=\bm{U}^{t_{0}}\bm{\Lambda}^{t_{0}}\bm{U}^{t_{0}\top}$.

\end{algorithmic} 
\end{algorithm}

\subsection{Fine-grained statistical guarantees for \textsf{HeteroPCA} }

We now move on to present our theoretical guarantees for Algorithm~\ref{alg:HeteroPCA-subspace}.
In order to account for the potential global rotational ambiguity,
we introduce the following rotation matrix as before 
\begin{equation}
\bm{R}_{\bm{U}}\coloneqq\arg\min_{\bm{O}\in\mathcal{O}^{r\times r}}\left\Vert \bm{U}\bm{O}-\bm{U}^{\natural}\right\Vert _{\mathrm{F}}^{2},\label{eq:rotation-matrix-subspace}
\end{equation}
where we recall that $\mathcal{O}^{r\times r}$ represents the set
of $r\times r$ orthonormal matrices. It is also helpful to define
the following quantities: for all $m\in[n_{1}]$, \begin{subequations}\label{eq:error-op-subspace-general}
\begin{align}
\zeta_{\mathsf{op}} & \coloneqq\sigma^{2}\sqrt{n_{1}n_{2}}\log n+\sigma\sigma_{1}^{\natural}\sqrt{n_{1}\log n},\\
\zeta_{\mathsf{op},m} & \coloneqq\sigma\sigma_{m}\sqrt{n_{1}n_{2}}\log n+\sigma_{m}\sigma_{1}^{\natural}\sqrt{n_{1}\log n},
\end{align}
\end{subequations}Our result is as follows, with the proof postponed
to Appendix~\ref{subsec:Analysis-for-HeteroPCA}.

\begin{theorem}\label{thm:hpca_inference_general} Suppose that Assumptions
\ref{assumption:subspace-incoherence}-\ref{assumption:subspace-noise}
hold. Assume that 
\begin{equation}
n_{1}\gtrsim\kappa^{\natural4}\mu^{\natural}r+\mu^{\natural2}r\log^{2}n,\qquad n_{2}\gtrsim r\log^{4}n,\qquad\text{and}\qquad\zeta_{\mathsf{op}}\ll\frac{\sigma_{r}^{\natural2}}{\kappa^{\natural2}},\label{eq:hpca-inference-conditions}
\end{equation}
and that the algorithm is run for $t_{0}\geq\log\left(\frac{\sigma_{1}^{\star2}}{\zeta_{\mathsf{op}}}\right)$
iterations. With probability exceeding $1-O(n^{-10})$, there exist
two matrices $\bm{Z}$ and $\bm{\Psi}$ such that the estimates returned
by \textsf{HeteroPCA} obey 
\begin{align}
\bm{U}\bm{R}_{\bm{U}}-\bm{U}^{\natural} & =\bm{Z}+\bm{\Psi},\label{eq:U-decompose-subspace-est}
\end{align}
where\begin{subequations}\label{eq:U-decompose-subspace-est-ZPsi}
\begin{align}
\bm{Z} & \coloneqq\bm{E}\bm{V}^{\natural}\left(\bm{\Sigma}^{\natural}\right)^{-1}+\mathcal{P}_{\mathsf{off}\text{-}\mathsf{diag}}\left(\bm{E}\bm{E}^{\top}\right)\bm{U}^{\natural}\left(\bm{\Sigma}^{\natural}\right)^{-2},\\
\left\Vert \bm{\Psi}\right\Vert _{2,\infty} & \lesssim\kappa^{\natural2}\frac{\mu^{\natural}r}{n_{1}}\frac{\zeta_{\mathsf{op}}}{\sigma_{r}^{\natural2}}+\kappa^{\natural2}\frac{\zeta_{\mathsf{op}}^{2}}{\sigma_{r}^{\natural4}}\sqrt{\frac{\mu^{\natural}r}{n_{1}}}.
\end{align}
In fact, for each $m\in[n_{1}]$, we further have
\begin{align}
\left\Vert \bm{Z}_{m,\cdot}\right\Vert _{2} & \lesssim\frac{\zeta_{\mathsf{op},m}}{\sigma_{r}^{\natural2}}\sqrt{\frac{\mu^{\natural}r}{n_{1}}}+\left\Vert \bm{U}_{m,\cdot}^{\natural}\right\Vert _{2}\left(\kappa^{\natural2}\sqrt{\frac{\mu^{\natural}r}{n_{1}}}\frac{\zeta_{\mathsf{op}}}{\sigma_{r}^{\natural2}}+\kappa^{\natural2}\frac{\zeta_{\mathsf{op}}^{2}}{\sigma_{r}^{\natural4}}\right),\\
\left\Vert \bm{\Psi}_{m,\cdot}\right\Vert _{2} & \lesssim\kappa^{\natural2}\frac{\zeta_{\mathsf{op}}\zeta_{\mathsf{op},m}}{\sigma_{r}^{\natural4}}\sqrt{\frac{\mu^{\natural}r}{n_{1}}}+\left\Vert \bm{U}_{m,\cdot}^{\natural}\right\Vert _{2}\left(\kappa^{\natural2}\sqrt{\frac{\mu^{\natural}r}{n_{1}}}\frac{\zeta_{\mathsf{op}}}{\sigma_{r}^{\natural2}}+\kappa^{\natural2}\frac{\zeta_{\mathsf{op}}^{2}}{\sigma_{r}^{\natural4}}\right).
\end{align}
\end{subequations}\end{theorem}

\begin{remark}The interested reader might wonder why Theorem \ref{thm:hpca_inference_general} is not valid for small $n_1$ (e.g., \eqref{eq:hpca-inference-conditions} does not hold when $n_1=2$), and we provide some intuition here. Recall from \eqref{eq:E-gram-decompose} that the diagonal of the sample Gram matrix can be significantly biased, and the \textsf{HeteroPCA} algorithm uses the off-diagonal information to iteratively estimate and refine the diagonal. When $n_1$ is small, the (untrustworthy) diagonal entries account for a non-negligible fraction of all entries of the entire sample Gram matrix, and as a result, we cannot hope to debias the diagonal reliably by \textsf{HeteroPCA} using only off-diagonal observations.\end{remark}

The expressions (\ref{eq:U-decompose-subspace-est}) and (\ref{eq:U-decompose-subspace-est-ZPsi})
make apparent a key decomposition of the estimation error. As we shall
see, the term $\bm{Z}$ is often the dominant term, which captures
both the first-order and second-order approximation (w.r.t.~the noise
matrix $\bm{E}$) of the estimation error. Unless the noise level
$\sigma$ is very small, we cannot simply ignore the second-order
term $\mathcal{P}_{\mathsf{off}\text{-}\mathsf{diag}}\left(\bm{E}\bm{E}^{\top}\right)\bm{U}^{\natural}\left(\bm{\Sigma}^{\natural}\right)^{-2}$,
as it is not necessarily dominated in size by the linear mapping term
$\bm{E}\bm{V}^{\natural}\left(\bm{\Sigma}^{\natural}\right)^{-1}$.
The simple and closed-form expression of $\bm{Z}$ --- in conjunction
with the fact that $\bm{\Psi}$ is well-controlled --- plays a crucial
role when developing a non-asymptotic distributional theory.

While Theorem \ref{thm:hpca_inference_general} is established mainly
to help derive distributional characterizations for PCA, we remark
that our analysis also delivers $\ell_{2,\infty}$ statistical guarantees
in terms of estimating $\bm{U}^{\natural}$ (see Lemma \ref{lemma:hpca-estimation-two-to-infty}
in the appendix). More specifically, our analysis asserts that 
\begin{align}
\left\Vert \bm{U}\bm{R}_{\bm{U}}-\bm{U}^{\natural}\right\Vert _{2,\infty} & \lesssim\frac{\zeta_{\mathsf{op}}}{\sigma_{r}^{\natural2}}\sqrt{\frac{\mu^{\natural}r}{n_{1}}}\label{eq:hpca-2-infty-guarantee}
\end{align}
with high probability, under the conditions of Theorem \ref{thm:hpca_inference_general}.
It is perhaps helpful to compare (\ref{eq:hpca-2-infty-guarantee})
with prior $\ell_{2,\infty}$ theory concerning estimation of $\bm{U}^{\natural}$. 
\begin{itemize}
\item We first compare Theorem \ref{thm:hpca_inference_general} with the
recent work \citep[Theorem 2]{agterberg2021entrywise}, which focused
on the regime $n_{2}\gtrsim n_{1}$ and showed that 
\[
\inf_{\bm{O}\in\mathcal{O}^{r\times r}}\left\Vert \bm{U}\bm{O}-\bm{U}^{\natural}\right\Vert _{2,\infty}\lesssim\left(\frac{\sigma^{2}}{\sigma_{r}^{\natural2}}\sqrt{rn_{1}n_{2}}\log n+\kappa^{\natural}\frac{\sigma}{\sigma_{r}^{\natural}}\sqrt{rn_{1}\log n}\right)\sqrt{\frac{\mu^{\natural}r}{n_{1}}}\asymp\frac{\zeta_{\mathsf{op}}}{\sigma_{r}^{\natural2}}\sqrt{\frac{\mu^{\natural}r^{2}}{n_{1}}}
\]
under the noise condition $\sigma\sqrt{n_{2}}\ll\sigma_{r}^{\natural}/(\kappa^{\natural}\sqrt{r\log n})$
(in addition to a few other conditions omitted here). Note that when
$\kappa^{\natural},\mu^{\natural},r\asymp1$, their $\ell_{2,\infty}$
error bound resembles (\ref{eq:hpca-2-infty-guarantee}), but the
condition $\sigma\sqrt{n_{2}}\ll\sigma_{r}^{\natural}/\sqrt{\log n}$
required therein is much stronger than the noise condition $\zeta_{\mathsf{op}}\ll\sigma_{r}^{\natural2}$
--- which is equivalent to $\sigma\sqrt[4]{n_{1}n_{2}}\ll\sigma_{r}^{\natural}/\sqrt{\log n}$
when $n_{2}\gtrsim n_{1}$ --- imposed by our theory (see (\ref{eq:hpca-inference-conditions})).
It is also worth emphasizing that the theory of \citet{agterberg2021entrywise}
is capable of accommodating dependent data (i.e.~they only require
the rows of $\bm{E}$ to be independent and allow dependence within
rows), which is beyond the scope of the present paper. 
\item Compared with the $\ell_{2,\infty}$ estimation error guarantees for
the diagonal-deleted spectral method in \citet[Theorem 1]{cai2019subspace},
our bound (\ref{eq:hpca-2-infty-guarantee}) is able to get rid of
the bias term incurred by diagonal deletion (see \citet[Equation (17)]{cai2019subspace}),
thus improving upon this prior result. 
\end{itemize}

It should be noted that fine-grained perturbation results akin to Theorem \ref{thm:hpca_inference_general} were also developed for the SVD algorithm in an earlier version of this paper on arXiv, as detailed in \citet[Section 6.1]{yan2021inference}. Subsequently, \citet{yan2024entrywise} presented more refined results for cases where the entries of the noise matrix $\bm{E}$ follow a sub-Gaussian distribution, with further information available in Appendix F therein.

Before concluding this section, it is natural to ask whether Theorem~\ref{thm:hpca_inference_general} 
can be used to conduct subspace inference when every entry of $\bm{E}$ is allowed to have completely difference variance.  
To begin with, for a broad class of $\bm{E}$ with independent and heteroskedastic  components, 
we can readily apply Theorem \ref{thm:hpca_inference_general} to obtain a distributional theory for 
\textsf{HeteroPCA} when estimating $\bm{U}^{\natural}$. 
Caution needs to be exercised, however, when it comes to confidence interval construction. 
On closer inspection, 
	evaluating $\bm{Z}$ (i.e., the first- and second-order approximation of the subspace estimation error) in Theorem \ref{thm:hpca_inference_general} requires knowledge about the right singular subspace $\bm{V}^\natural$ of $\bm{M}^{\natural}$, 
	which might sometimes be difficult or even infeasible to estimate in the unbalanced regime where $n_2\gg n_1$. 
	As a result, our theory is not guaranteed to deliver useful inferential methods for such cases, unless additional information about $\bm{V}^\natural$ is available.

\section{Discussion\label{sec:Discussion}}

%
%

In this paper, we have developed a suite of statistical inference
procedures to construct confidence regions for PCA in the presence
of missing data and heterogeneous corruption, which should be easy-to-use
in practice due to their data-driven nature. 
Compared to other prior algorithms like the SVD-based approach and the diagonal-deleted spectral method, 
the solution developed
based on \textsf{HeteroPCA} enjoys a broadened applicability range
without compromising statistical efficiency. The fine-grained distributional
characterizations we have developed are non-asymptotic, which naturally
lend themselves to high-dimensional settings. 

Moving forward, there are a variety of directions that are worthy of further investigation.
\begin{itemize}
	\item \emph{Improved dependency on $\kappa$, $\mu$, $r$ and $\kappa_{\omega}$.} In our general theorems (see Theorems~\ref{thm:pca-complete}-\ref{thm:ce-CI-complete} in the appendix), we allow $\kappa$, $\mu$, $r$ and $\kappa_{\omega}$ to grow. 
		However, our theoretical results scale suboptimally with these problem parameters. 
		It remains unclear how to sharpen the dependency on these parameters, which might require developing more refined analysis techniques.

	\item \emph{Approximate low-rank structure.} Our results assume exact low-rank structure of the spiked component $\bm{S}^\star$ of the covariance matrix. In reality, there is no shortage of applications where  $\bm{S}^\star$ is at best approximately low-rank. 
How to develop trustworthy inference procedures in the presence of approximate low-rank structure?  
		Unfortunately, our current leave-one-out analysis framework relies heavily on the exact rank-$r$ structure (unless $\sigma_{r+1}^{\star}$ is extremely small); new analysis ideas are needed in order to tackle approximate low-rank structure.

	\item \emph{General missing pattern.} Uncertainty quantification in the face of heterogeneous missing patterns is another important topic of practical value. Consider, for example, the case where the entries in the same row of $\bm{X}$ are sampled with the same rate (i.e.,  the $(l,j)$-th entry of $\bm{X}$ is observed with probability $p_l$). Then by constructing the following data matrix via inverse probability weighting
		\[
			\Big[ \mathsf{diag}\big( p_1, p_2,\cdots, p_d\big) \Big]^{-1} \bm{Y}  ,
		\]
		we obtain an unbiased estimate of $\bm{X}$, and the theory developed can be readily extended to perform valid inference. Note that we can also replace $\{p_l\}$ via their empirical estimates in the inference procedures.  
		Nevertheless, in the more general case where the sampling rates are allowed to vary across all locations, it is unclear how to construct an unbiased estimate of $\bm{X}$ without knowing the per-entry sampling rates in advance; hence, our theory fails to accommodate this general scenario. 
		Extending our current results to such general sampling patterns might call for new analysis tools.

	\item \emph{Inference for individual principal components.} Moving beyond inference and uncertainty quantification for the principal subspace and the spiked covariance matrix, it is interesting to investigate how to conduct valid inference on individual principal components, particularly when the associated eigengap is vanishingly
	small \citep{li2021minimax}.
	\item \emph{Extension to unknown mean, dependent or adversarial noise.}  If the observed data are inherently biased with \emph{a priori}
	unknown means, how to properly compensate for the bias? What if the
	noise components are inter-dependent, and what if the observed data
	samples are further corrupted by a non-negligible fraction of adversarial
	outliers? 
	\item \emph{Minimax-optimal estimation and inference.} As recognized in the matrix completion literature \citep{KesMonSew2010,ma2017implicit,chen2019nonconvex}, spectral methods alone are in general unable to yield minimax-optimal statistical accuracy in the presence of missing data, 
		given that spectral methods inherently treat the missingness effect as some sort of ``noise''. 
		The same message --- namely, sub-optimality of \textsf{HeteroPCA} in the face of missing data --- carries over to the PCA setting considered herein.  We conjecture that a subsequent refinement procedure (e.g., gradient descent tailored to compute the maximum likelihood estimate) is needed in order to reach minimax optimality, and we leave this for future investigation.

	\item   \emph{Applications in financial econometrics}.  In addition to applications to the uncertainty quantification in the matrix completion problems in recommender system, the inferential procedure and analysis tools we have developed in this paper have applications in finance and econometrics. For example, our analysis and results for principal subspace are useful in testing factor structures in famous Fama-French factor models, and can also be used in sector/industry clustering using stock returns \citep{porter1998clusters}; our results on uncertainty quantification for the spiked covariance matrix could also shed light on how to better quantify the risk in portfolio optimization that takes into account on the uncertainty in the risk estimation.

\end{itemize}

\section*{Acknowledgements}

Y.~Chen is supported in part by the Alfred P. Sloan Research Fellowship, the Google Research Scholar Award, 
the AFOSR YIP award FA9550-19-1-0030, by the ONR grant N00014-19-1-2120,
by the ARO grant W911NF-20-1-0097, and by the NSF grants CCF-1907661,
IIS-2218713, IIS-2218773, DMS-2014279, CCF-2221009. 
J.~Fan is supported in part by the ONR grants N00014-19-1-2120, N00014-22-1-2340, by the NSF grants DMS-1662139, DMS-1712591, DMS-2052926, DMS-2053832, DMS-2210833, and by the NIH grant 2R01-GM072611-15.
Y.~Yan is supported in part by the Charlotte Elizabeth Procter Honorific
Fellowship from Princeton University and the Norbert Wiener Postdoctoral Fellowship from MIT.

\appendix

\section{Additional notation and organization of the appendix\label{sec:Additional-notation}}

For any matrix $\bm{U}$ with orthonormal columns, we denote by $\mathcal{P}_{\bm{U}}(\bm{M})\coloneqq\bm{U}\bm{U}^{\top}\bm{M}$
the Euclidean projection of a matrix $\bm{M}$ onto the column space
of $\bm{U}$, and let $\mathcal{P}_{\bm{U}^{\perp}}(\bm{M})=\bm{M}-\mathcal{P}_{\bm{U}}(\bm{M})$
denote the Euclidean projection of $\bm{M}$ onto the orthogonal complement
of the column space of $\bm{U}$. For any matrix $\bm{B}\in\mathbb{R}^{n_{1}\times n_{2}}$
and some $l\in[n_{1}]$, define $\mathcal{P}_{-l,\cdot}(\bm{B})$
to be the orthogonal projection of the matrix $\bm{B}$ onto the subspace
of matrix that vanishes outside the $l$-th row, namely, 
\begin{equation}
\left[\mathcal{P}_{-l,\cdot}\left(\bm{B}\right)\right]_{i,j}=\begin{cases}
B_{i,j}, & \text{if }i\neq l,\\
0, & \text{otherwise,}
\end{cases}\qquad\forall\,(i,j)\in[n_{1}]\times[n_{2}].\label{eq:P-no-l-definition}
\end{equation}
For any point $x\in\mathbb{R}^{d}$ and any non-empty convex set $\mathcal{C}\in\mathscr{C}^{d}$
satisfying $\mathcal{C}\neq\mathbb{R}^{d}$, let us define the signed
distance function as follows 
\begin{equation}
\delta_{\mathcal{C}}\left(x\right)\coloneqq\begin{cases}
-\mathsf{dist}\left(x,\mathbb{R}^{d}\setminus\mathcal{C}\right), & \text{if }x\in\mathcal{C};\\
\mathsf{dist}\left(x,\mathcal{C}\right), & \text{if }x\notin\mathcal{C}.
\end{cases}\label{eq:defn-sign-distance}
\end{equation}
Here, $\mathsf{dist}(x,\mathcal{A})$ is the Euclidean distance between
a point $x\in\mathbb{R}^{d}$ and a non-empty set $\mathcal{A}\subseteq\mathbb{R}^{d}$.
Also, for any $\varepsilon\in\mathbb{R}$, define 
\begin{equation}
\mathcal{C}^{\varepsilon}\coloneqq\left\{ x\in\mathbb{R}^{d}:\delta_{\mathcal{C}}\left(x\right)\leq\varepsilon\right\} \label{eq:defn-C-epsilon}
\end{equation}
for any non-empty convex set $\mathcal{C}\in\mathscr{C}^{d}$ satisfying
$\mathcal{C}\neq\mathbb{R}^{d}$, and define $\varnothing^{\varepsilon}=\varnothing$
and $(\mathbb{R}^{d})^{\varepsilon}=\mathbb{R}^{d}$.

\paragraph{Organization of the appendix.}
The rest of the appendix is organized as follows. 
Appendix~\ref{subsec:Analysis-for-HeteroPCA} is devoted to proving Theorem \ref{thm:hpca_inference_general}. 
Appendices~\ref{sec:Analysis:-the-approach-heteroPCA}
and \ref{appendix:hpca-auxiliary-lemmas} are dedicated to establishing
our theory for \textsf{HeteroPCA} (i.e., Theorems \ref{thm:pca-complete}-\ref{thm:ce-CI-complete}).

\section{A list of general theorems}

\label{sec:A-list-of-general-thms}

For simplicity of presentation, the theorems presented in the main
text (i.e., Section~\ref{sec:Distributional-theory-inference}) concentrate
on the scenario where $\kappa,\mu,r\asymp1$. Note, however, that
our theoretical framework can certainly allow these parameters to
grow with the problem dimension. In this section, we provide a list
of theorems accommodating more general scenarios; all subsequent proofs
are dedicated to establishing these general theorems.

The following theorems generalize Theorem~\ref{thm:pca}, Theorem~\ref{thm:pca-cr},
Theorem~\ref{thm:ce} and Theorem~\ref{thm:ce-CI}, respectively.

\begin{theorem}\label{thm:pca-complete}Suppose that $p<1-\delta$
for some arbitrary constant $0<\delta<1$ or $p=1$. In addition, suppose that
Assumption \ref{assumption:noise} holds and $n\gtrsim\kappa^{8}\mu^{2}r^{4}\kappa_{\omega}^{2}\log^{4}(n+d)$,
$d\gtrsim\kappa^{7}\mu^{3}r^{7/2}\kappa_{\omega}^{2}\log^{5}(n+d)$,
\[
\frac{\omega_{\max}^{2}}{p\sigma_{r}^{\star2}}\sqrt{\frac{d}{n}}\lesssim\frac{1}{\kappa^{4}\mu^{3/2}r^{11/4}\kappa_{\omega}\log^{7/2}\left(n+d\right)},\qquad\frac{\omega_{\max}}{\sigma_{r}^{\star}}\sqrt{\frac{d}{np}}\lesssim\frac{1}{\kappa^{3/2}\mu r^{5/4}\kappa_{\omega}^{1/2}\log^{3}\left(n+d\right)},
\]
\[
ndp^{2}\gtrsim\kappa^{9}\mu^{4}r^{13/2}\kappa_{\omega}^{2}\log^{9}\left(n+d\right),\qquad np\gtrsim\kappa^{9}\mu^{3}r^{11/2}\kappa_{\omega}^{2}\log^{7}\left(n+d\right).
\]
Then the estimate $\bm{U}$ returned by Algorithm \ref{alg:PCA-HeteroPCA}
with number of iterations satisfying (\ref{eq:main-iteration}) satisfies
\[
\sup_{1\leq l\leq d}\,\sup_{\mathcal{C}\in\mathscr{C}^{r}}\left|\mathbb{P}\left(\left[\bm{U}\mathsf{sgn}\left(\bm{U}^{\top}\bm{U}^{\star}\right)-\bm{U}^{\star}\right]_{l,\cdot}\in\mathcal{C}\right)-\mathcal{N}\left(\bm{0},\bm{\Sigma}_{U,l}^{\star}\right)\left\{ \mathcal{C}\right\} \right|=o\left(1\right).
\]
where $\mathscr{C}^{r}$ represents the the set of all convex sets
in $\mathbb{R}^{r}$, and for each $l\in[d]$,
\begin{align*}
\bm{\Sigma}_{U,l}^{\star} & \coloneqq\left(\frac{1-p}{np}\left\Vert \bm{U}_{l,\cdot}^{\star}\bm{\Sigma}^{\star}\right\Vert _{2}^{2}+\frac{\omega_{l}^{\star2}}{np}\right)\left(\bm{\Sigma}^{\star}\right)^{-2}+\frac{2\left(1-p\right)}{np}\bm{U}_{l,\cdot}^{\star\top}\bm{U}_{l,\cdot}^{\star}+\left(\bm{\Sigma}^{\star}\right)^{-2}\bm{U}^{\star\top}\mathsf{diag}\left\{ \left[d_{l,i}^{\star}\right]_{1\leq i\leq d}\right\} \bm{U}^{\star}(\bm{\Sigma}^{\star})^{-2}
\end{align*}
where 
\[
d_{l,i}^{\star}\coloneqq\frac{1}{np^{2}}\left[\omega_{l}^{\star2}+\left(1-p\right)\left\Vert \bm{U}_{l,\cdot}^{\star}\bm{\Sigma}^{\star}\right\Vert _{2}^{2}\right]\left[\omega_{i}^{\star2}+\left(1-p\right)\left\Vert \bm{U}_{i,\cdot}^{\star}\bm{\Sigma}^{\star}\right\Vert _{2}^{2}\right]+\frac{2\left(1-p\right)^{2}}{np^{2}}S_{l,i}^{\star2}.
\]
\end{theorem} 

\begin{theorem}\label{thm:pca-cr-complete} Suppose that the conditions
of Theorem \ref{thm:pca-complete} hold. Further assume that $n\gtrsim\kappa^{12}\mu^{3}r^{11/2}\kappa_{\omega}\log^{5}(n+d)$,
\[
\frac{\omega_{\max}^{2}}{p\sigma_{r}^{\star2}}\sqrt{\frac{d}{n}}\lesssim\frac{1}{\kappa^{9/2}\mu^{3/2}r^{9/4}\kappa_{\omega}^{3/2}\log^{7/2}\left(n+d\right)},\qquad\frac{\omega_{\max}}{\sigma_{r}^{\star}}\sqrt{\frac{d}{np}}\lesssim\frac{1}{\kappa^{5}\mu^{3/2}r^{9/4}\kappa_{\omega}^{3/2}\log^{3}\left(n+d\right)},
\]
\[
ndp^{2}\gtrsim\kappa^{11}\mu^{5}r^{13/2}\kappa_{\omega}^{3}\log^{9}\left(n+d\right),\qquad np\gtrsim\kappa^{11}\mu^{4}r^{11/2}\kappa_{\omega}^{3}\log^{7}\left(n+d\right).
\]
Then the confidence region $\mathsf{CR}_{U,l}^{1-\alpha}$ computed
in Algorithm \ref{alg:PCA-HeteroPCA-CR} obeys 
\[
\mathbb{P}\left(\bm{U}_{l,\cdot}^{\star}\mathsf{sgn}\left(\bm{U}^{\star\top}\bm{U}\right)\in\mathsf{CR}_{U,l}^{1-\alpha}\right)=1-\alpha+o\left(1\right).
\]
\end{theorem} 

\begin{theorem}\label{thm:ce-complete}Suppose that $p<1-\delta$
for some arbitrary constant $0<\delta<1$ or $p=1$. Consider any $1\leq i,j\leq d$.
Assume that $\bm{U}^{\star}$ is $\mu$-incoherent and satisfies the
following condition 
\[
\left\Vert \bm{U}_{i,\cdot}^{\star}\right\Vert _{2}+\left\Vert \bm{U}_{j,\cdot}^{\star}\right\Vert _{2}\gtrsim\kappa r^{1/2}\log^{5/2}\left(n+d\right)\left(\frac{\omega_{\max}}{\sigma_{r}^{\star}}\sqrt{\frac{d}{np}}+\kappa_{\omega}\frac{\omega_{\max}^{2}}{p\sigma_{r}^{\star2}}\sqrt{\frac{d}{n}}+\frac{\kappa\mu r\kappa_{\omega}\log\left(n+d\right)}{\sqrt{ndp^{2}}}\right)\sqrt{\frac{r}{d}}.
\]
In addition, suppose that Assumption \ref{assumption:noise} holds
and $n\gtrsim r\log^{4}(n+d)$, $d\gtrsim\kappa^{8}\mu^{3}r^{3}\kappa_{\omega}^{2}\log^{5}(n+d)$,
\[
\frac{\omega_{\max}^{2}}{p\sigma_{r}^{\star2}}\sqrt{\frac{d}{n}}\lesssim\frac{1}{\kappa^{3/2}\mu r\kappa_{\omega}^{1/2}\log^{7/2}\left(n+d\right)},\qquad\frac{\omega_{\max}}{\sigma_{r}^{\star}}\sqrt{\frac{d}{np}}\lesssim\frac{1}{\kappa^{4}\mu r\kappa_{\omega}^{1/2}\log^{3}\left(n+d\right)},
\]
\[
ndp^{2}\gtrsim\kappa^{10}\mu^{4}r^{4}\kappa_{\omega}^{2}\log^{7}\left(n+d\right),\qquad np\gtrsim\kappa^{10}\mu^{3}r^{3}\kappa_{\omega}^{2}\log^{7}\left(n+d\right),
\]
Then the matrix $\bm{S}$ computed by Algorithm \ref{alg:PCA-HeteroPCA}
with number of iterations satisfying (\ref{eq:main-iteration}) obeys
\[
\sup_{t\in\mathbb{R}}\left|\mathbb{P}\left(\frac{S_{i,j}-S_{i,j}^{\star}}{\sqrt{v_{i,j}^{\star}}}\leq t\right)-\Phi\left(t\right)\right|=o\left(1\right).
\]
\end{theorem} 

\begin{theorem}\label{thm:ce-CI-complete}Suppose that the conditions
of Theorem \ref{thm:ce-complete} hold. Further assume that $n\gtrsim\kappa^{9}\mu^{3}r^{4}\kappa_{\omega}^{3}\log^{4}(n+d)$,
\[
\frac{\omega_{\max}^{2}}{p\sigma_{r}^{\star2}}\sqrt{\frac{d}{n}}\lesssim\frac{1}{\kappa^{3}\mu^{3/2}r^{5/2}\kappa_{\omega}^{3/2}\log^{3}\left(n+d\right)},\qquad\frac{\omega_{\max}}{\sigma_{r}^{\star}}\sqrt{\frac{d}{np}}\lesssim\frac{1}{\kappa^{7/2}\mu^{3/2}r^{5/2}\kappa_{\omega}^{3/2}\log^{5/2}\left(n+d\right)},
\]
\[
ndp^{2}\gtrsim\kappa^{8}\mu^{5}r^{7}\kappa_{\omega}^{3}\log^{8}\left(n+d\right),\qquad np\gtrsim\kappa^{8}\mu^{4}r^{6}\kappa_{\omega}^{3}\log^{6}\left(n+d\right),
\]
and 
\[
\left\Vert \bm{U}_{i,\cdot}^{\star}\right\Vert _{2}+\left\Vert \bm{U}_{j,\cdot}^{\star}\right\Vert _{2}\gtrsim\kappa^{2}\mu^{3/2}r^{5/2}\kappa_{\omega}^{3/2}\log^{5/2}\left(n+d\right)\left[\frac{\kappa\mu r\log\left(n+d\right)}{\sqrt{nd}p}+\frac{\omega_{\max}^{2}}{p\sigma_{r}^{\star2}}\sqrt{\frac{d}{n}}+\frac{\omega_{\max}}{\sigma_{r}^{\star}}\sqrt{\frac{d}{np}}\right]\sqrt{\frac{r}{d}}.
\]
Then the confidence interval computed in Algorithm~\ref{alg:CE-HeteroPCA-CI}
obeys 
\[
\mathbb{P}\left(S_{i,j}^{\star}\in\mathsf{CI}_{i,j}^{1-\alpha}\right)=1-\alpha+o\left(1\right).
\]
\end{theorem}

\begin{remark}
	While the above theorems allow the problem parameters $\kappa, \mu, r, \kappa_{\omega}$ to grow, 
	our results remain suboptimal in terms of the dependency on these parameters. 
	For instance, our theorems require  stringent dependency on the condition number $\kappa$,
	which has been a common issue in analyzing spectral methods for low-rank matrix estimation using leave-one-out arguments
	(e.g., in order to obtain $\ell_{2,\infty}$ guarantees for matrix completion, the sample complexity requirement in the prior work by \citet{chen2019nonconvex} scales as $\kappa^{10}$). 
	It is also worth noting that inference and uncertainty quantification might require stronger conditions compared to the estimation task,
	 which further complicates matters (e.g., the multivariate Berry-Esseen theorem (cf.~Theorem~\ref{thm:berry-esseen}) employed in the current paper might already exhibit suboptimal scaling with $r$).   Improving the dependency on all these parameters is a fundamentally important direction for future investigation.
\end{remark}

\section{Analysis for \textsf{HeteroPCA} applied to subspace estimation (Theorem
\ref{thm:hpca_inference_general})\label{subsec:Analysis-for-HeteroPCA}}

This section outlines the proof that establishes our statistical guarantees
stated in Theorem \ref{thm:hpca_inference_general}. We shall begin
by isolating several useful lemmas, and then combine these lemmas
to complete the proof.

\subsection{A few key lemmas}

We now state below a couple of key lemmas that lead to improved statistical
guarantees of Algorithm \ref{alg:HeteroPCA-subspace}. From now on,
we will denote 
\[
\bm{G}=\bm{G}^{t_{0}},\qquad\bm{U}=\bm{U}^{t_{0}}\qquad\text{and}\qquad\bm{\Sigma}=(\bm{\Lambda}^{t_{0}})^{1/2}
\]
and let 
\[
\bm{H}\coloneqq\bm{U}^{\top}\bm{U}^{\natural}.
\]

To begin with, the first lemma controls the discrepancy between the
sample gram matrix $\bm{G}$ and the ground truth $\bm{G}^{\natural}$,
which improves upon prior theory developed in \citet{zhang2018heteroskedastic}.

\begin{lemma}\label{lemma:hpca-estimation}Suppose that the assumptions
of Theorem \ref{thm:hpca_inference_general} hold. Suppose that the
number of iterations exceeds $t_{0}\geq\log\left(\frac{\sigma_{1}^{\star2}}{\zeta_{\mathsf{op}}}\right)$.
Then with probability exceeding $1-O(n^{-10})$, the iterate $\bm{G}\coloneqq\bm{G}^{t_{0}}$
computed in Algorithm \ref{alg:HeteroPCA-subspace} satisfies 
\begin{align}
\left\Vert \mathcal{P}_{\mathsf{diag}}\left(\bm{G}-\bm{G}^{\natural}\right)\right\Vert  & \lesssim\kappa^{\natural2}\sqrt{\frac{\mu^{\natural}r}{n_{1}}}\zeta_{\mathsf{op}},\label{eq:hpca-estimation-diagonal}\\
\left\Vert \bm{G}-\bm{G}^{\natural}\right\Vert  & \lesssim\zeta_{\mathsf{op}},\label{eq:hpca-estimation-full}
\end{align}
where $\zeta_{\mathsf{op}}$ is defined in (\ref{eq:error-op-subspace-general}).
In addition, for each $m\in[n_{1}]$,
\[
\left|G_{m,m}-G_{m,m}^{\natural}\right|\lesssim\kappa^{\natural2}\zeta_{\mathsf{op}}\left(\left\Vert \bm{U}_{m,\cdot}^{\natural}\right\Vert _{2}+\left\Vert \bm{U}_{m,\cdot}\right\Vert _{2}\right).
\]
 \end{lemma}\begin{proof}See Appendix \ref{appendix:proof-lemma-HeteroPCA-estimation}.\end{proof}

Lemma \ref{lemma:hpca-estimation} makes clear that $\bm{G}$ converges
to $\bm{G}^{\natural}$ as the signal-to-noise ratio increases. We
pause to compare this lemma with its counterpart in \citet{zhang2018heteroskedastic}.
Specifically, \citet[Theorem 7]{zhang2018heteroskedastic} and its
proof demonstrated that $\Vert\mathcal{P}_{\mathsf{diag}}(\bm{G}-\bm{G}^{\natural})\Vert\lesssim\zeta_{\mathsf{op}}$.
In comparison our result (\ref{eq:hpca-estimation-diagonal}) strengthens
the prior estimation bound by a factor of $1/\sqrt{n_{1}}$ for the
scenario with $\kappa^{\natural},\mu^{\natural},r=O(1)$. This improvement
serves as one of the key analysis ingredients that allows us to sharpen
the statistical guarantees for \textsf{HeteroPCA}.

The above bound on the difference between $\bm{G}$ and $\bm{G}^{\natural}$
in turn allows one to develop perturbation bounds for the eigenspace
measured by the spectral norm, as stated in the following lemma. In
the meantime, this lemma also contains some basic facts regarding
$\bm{H}$ and $\bm{R}_{\bm{U}}$.

\begin{lemma}\label{lemma:hpca-basic-facts}Suppose that the assumptions
of Theorem \ref{thm:hpca_inference_general} hold, and recall the
definition of $\zeta_{\mathsf{op}}$ in (\ref{eq:error-op-subspace-general}).
Then with probability exceeding $1-O(n^{-10})$, we have 
\begin{equation}
\left\Vert \bm{U}\bm{H}-\bm{U}^{\star}\right\Vert \lesssim\frac{\zeta_{\mathsf{op}}}{\sigma_{r}^{\natural2}}\qquad\text{and}\qquad\left\Vert \bm{U}\bm{R}_{\bm{U}}-\bm{U}^{\star}\right\Vert \lesssim\frac{\zeta_{\mathsf{op}}}{\sigma_{r}^{\natural2}},\label{eq:pca-UH-U-star-spectral}
\end{equation}
\[
\left\Vert \bm{H}-\bm{R}_{\bm{U}}\right\Vert \lesssim\frac{\zeta_{\mathsf{op}}^{2}}{\sigma_{r}^{\natural4}}\qquad\text{and}\qquad\text{\ensuremath{\left\Vert \bm{H}^{\top}\bm{H}-\bm{I}_{r}\right\Vert }}\lesssim\frac{\zeta_{\mathsf{op}}^{2}}{\sigma_{r}^{\natural4}},
\]
and 
\[
\frac{1}{2}\leq\sigma_{i}\left(\bm{H}\right)\leq2,\qquad\forall1\leq i\leq r.
\]
\end{lemma}\begin{proof}See Appendix \ref{appendix:proof-hpca-H-R-proximity}.\end{proof}

While the above two lemmas focus on spectral norm metrics, the following
lemma takes one substantial step further by characterizing the difference
between $\bm{G}$ and $\bm{G}^{\natural}$ in each row, when projected
onto the subspace spanned by $\bm{U}^{\natural}$.

\begin{lemma}\label{lemma:hpca-1}Suppose that $n_{1}\gtrsim\kappa^{\natural4}\mu^{\natural}r$.
Then with probability exceeding $1-O(n^{-10})$, we have 
\begin{align*}
\left\Vert \left(\bm{G}-\bm{G}^{\natural}\right)_{m,\cdot}\bm{U}^{\natural}\right\Vert _{2} & \lesssim\zeta_{\mathsf{op},m}\sqrt{\frac{\mu^{\natural}r}{n_{1}}}+\sigma\sigma_{1}^{\natural}\sqrt{\mu^{\natural}r\log n}\left\Vert \bm{U}_{m,\cdot}^{\natural}\right\Vert _{2}+\kappa^{\natural2}\zeta_{\mathsf{op}}\left\Vert \bm{U}_{m,\cdot}^{\natural}\right\Vert _{2}^{2}\\
 & \quad+\kappa^{\natural2}\zeta_{\mathsf{op}}\left\Vert \bm{U}_{m,\cdot}\bm{H}-\bm{U}_{m,\cdot}^{\natural}\right\Vert _{2}\left\Vert \bm{U}_{m,\cdot}^{\natural}\right\Vert _{2}
\end{align*}
simultaneously for each $m\in[n_{1}]$, with $\zeta_{\mathsf{op}}$
and $\zeta_{\mathsf{op},m}$ defined in (\ref{eq:error-op-subspace-general}).
\end{lemma}\begin{proof}See Appendix \ref{sec:proof-lemma-pca-1}.\end{proof}

Next, we present a crucial technical lemma that uncovers the intertwined
relation between $\bm{U}\bm{\Sigma}^{2}\bm{H}$ and $\bm{G}\bm{U}^{\natural}$
in an $\ell_{2,\infty}$ sense. To establish the $\ell_{2,\infty}$
bound in this lemma, we invoke the powerful leave-one-out analysis
framework to decouple complicated statistical dependency.

\begin{lemma}\label{lemma:hpca-approx-1}Suppose that $n_{1}\gtrsim\kappa^{\natural4}\mu^{\natural}r$
and $\zeta_{\mathsf{op}}\ll\sigma_{r}^{\natural2}/\kappa^{\natural2}$.
Then with probability exceeding $1-O(n^{-10})$, 
\begin{align*}
\left\Vert \left(\bm{U}\bm{\Sigma}^{2}\bm{H}-\bm{G}\bm{U}^{\natural}\right)_{m,\cdot}\right\Vert _{2} & =\left\Vert \bm{G}_{m,\cdot}\left(\bm{U}\bm{H}-\bm{U}^{\natural}\right)\right\Vert _{2}\\
 & \lesssim\zeta_{\mathsf{op},m}\left(\kappa^{\natural2}\frac{\zeta_{\mathsf{op}}}{\sigma_{r}^{\natural2}}\sqrt{\frac{\mu^{\natural}r}{n_{1}}}+\left\Vert \bm{U}\bm{H}-\bm{U}^{\natural}\right\Vert _{2,\infty}\right)+\kappa^{\natural2}\frac{\zeta_{\mathsf{op}}^{2}}{\sigma_{r}^{\natural2}}\left\Vert \bm{U}_{m,\cdot}^{\natural}\right\Vert _{2}\\
 & \quad+\kappa^{\natural2}\zeta_{\mathsf{op}}\left\Vert \bm{U}_{m,\cdot}^{\natural}\right\Vert _{2}\left\Vert \bm{U}_{m,\cdot}\bm{H}-\bm{U}_{m,\cdot}^{\natural}\right\Vert _{2}+\kappa^{\natural2}\zeta_{\mathsf{op}}\left\Vert \bm{U}_{m,\cdot}\bm{H}-\bm{U}_{m,\cdot}^{\natural}\right\Vert _{2}^{2}.
\end{align*}
holds simultaneously for each $m\in[n_{1}]$. Here, $\zeta_{\mathsf{op}}$
and $\zeta_{\mathsf{op},m}$ are quantities defined in (\ref{eq:error-op-subspace-general}).\end{lemma}\begin{proof}See
Appendix \ref{sec:proof-lemma-HeteroPCA-approx-1}.\end{proof}

Furthermore, the lemma below reveals that $\bm{\Sigma}^{2}$ and $\bm{\Sigma}^{\natural2}$
remain close even after $\bm{\Sigma}^{2}$ is rotated by the rotation
matrix $\bm{R}_{\bm{U}}$.

\begin{lemma}\label{lemma:hpca-approx-2}Suppose that $n_{1}\gtrsim\mu^{\natural2}r\log^{2}n$
and $n_{2}\gtrsim r\log^{4}n$. Then with probability exceeding $1-O(n^{-10})$
we have 
\[
\left\Vert \bm{R}_{\bm{U}}^{\top}\bm{\Sigma}^{2}\bm{R}_{\bm{U}}-\bm{\Sigma}^{\natural2}\right\Vert \lesssim\kappa^{\natural2}\sqrt{\frac{\mu^{\natural}r}{n_{1}}}\zeta_{\mathsf{op}}+\kappa^{\natural2}\frac{\zeta_{\mathsf{op}}^{2}}{\sigma_{r}^{\natural2}},
\]
where $\zeta_{\mathsf{op}}$ is defined in (\ref{eq:error-op-subspace-general}).\end{lemma}\begin{proof}See
Appendix \ref{sec:proof-lemma-HeteroPCA-approx-2}.\end{proof}

With the above auxiliary results in place, we can put them together
to yield the following $\ell_{2,\infty}$ statistical guarantees.

\begin{lemma}\label{lemma:hpca-estimation-two-to-infty}Suppose that
$\zeta_{\mathsf{op}}\ll\sigma_{r}^{\natural2}/\kappa^{\natural2}$
and that 
\[
n_{1}\gtrsim\kappa^{\natural4}\mu^{\natural}r+\mu^{\natural2}r\log^{2}n,\qquad n_{2}\gtrsim r\log^{4}n.
\]
Then with probability at least $1-O(n^{-10})$, 
\begin{align*}
\left\Vert \bm{U}_{m,\cdot}\bm{R}_{\bm{U}}-\bm{U}_{m,\cdot}^{\natural}\right\Vert _{2} & \lesssim\frac{\zeta_{\mathsf{op},m}}{\sigma_{r}^{\natural2}}\sqrt{\frac{\mu^{\natural}r}{n_{1}}}+\left\Vert \bm{U}_{m,\cdot}^{\natural}\right\Vert _{2}\left(\kappa^{\natural2}\sqrt{\frac{\mu^{\natural}r}{n_{1}}}\frac{\zeta_{\mathsf{op}}}{\sigma_{r}^{\natural2}}+\kappa^{\natural2}\frac{\zeta_{\mathsf{op}}^{2}}{\sigma_{r}^{\natural4}}\right)
\end{align*}
holds simultaneously for each $m\in[n_{1}]$, and 
\begin{align*}
\left\Vert \bm{U}\bm{R}_{\bm{U}}-\bm{U}^{\natural}\right\Vert _{2,\infty} & \lesssim\frac{\zeta_{\mathsf{op}}}{\sigma_{r}^{\natural2}}\sqrt{\frac{\mu^{\natural}r}{n_{1}}},
\end{align*}
where $\zeta_{\mathsf{op}}$ and $\zeta_{\mathsf{op},m}$ are defined
in (\ref{eq:error-op-subspace-general}).\end{lemma}\begin{proof}
See Appendix \ref{appendix:proof-lemma-hpca-estimation-two-to-infty}.\end{proof}

\subsection{Proof of Theorem \ref{thm:hpca_inference_general}\label{section:proof-subspace-inference}}

Armed with the above lemmas, we are in a position to establish Theorem
\ref{thm:hpca_inference_general}. It is worth noting that, while
Lemma~\ref{lemma:hpca-estimation-two-to-infty} delivers $\ell_{2,\infty}$
perturbation bounds for the eigenspace, it falls short of revealing
the relation between the estimation error and the desired approximation
\begin{equation}
\bm{Z}=\bm{E}\bm{V}^{\natural}(\bm{\Sigma}^{\natural})^{-1}+\mathcal{P}_{\mathsf{off}\text{-}\mathsf{diag}}(\bm{E}\bm{E}^{\top})\bm{U}^{\natural}(\bm{\Sigma}^{\natural})^{-2}.\label{eq:defn-Z-subspace-proof}
\end{equation}
In order to justify the tightness of this approximation $\bm{Z}$,
we intend to establish each of the following steps: 
\begin{equation}
\bm{U}\bm{R}_{\bm{U}}\bm{\Sigma}^{\natural2}\overset{\text{Step 1}}{\approx}\bm{U}\bm{\Sigma}^{2}\bm{R}_{\bm{U}}\overset{\text{Step 2}}{\approx}\bm{U}\bm{\Sigma}^{2}\bm{H}\overset{\text{Step 3}}{\approx}\bm{G}\bm{U}^{\natural}\overset{\text{Step 4}}{\approx}\bm{U}^{\natural}\bm{\Sigma}^{\natural2}+\underset{=\,\bm{Z}(\bm{\Sigma}^{\natural})^{2}}{\underbrace{\bm{E}\bm{V}^{\natural}\bm{\Sigma}^{\natural}+\mathcal{P}_{\mathsf{off}\text{-}\mathsf{diag}}\left(\bm{E}\bm{E}^{\top}\right)\bm{U}^{\natural}}},\label{eq:step-1234-subspace}
\end{equation}
which would in turn ensure that 
\[
\bm{U}\bm{R}_{\bm{U}}\approx\bm{U}^{\natural}+\bm{Z}
\]
as advertised.

\subsubsection{Step 1: establishing the proximity of $\bm{U}\bm{R}_{\bm{U}}\bm{\Sigma}^{\natural2}$
and $\bm{U}\bm{\Sigma}^{2}\bm{R}_{\bm{U}}$. }

For each $m\in[n_{1}]$, Lemma \ref{lemma:hpca-approx-2} tells us
that 
\begin{align}
\left\Vert \bm{U}_{m,\cdot}\bm{R}_{\bm{U}}\bm{\Sigma}^{\natural2}-\bm{U}_{m,\cdot}\bm{\Sigma}^{2}\bm{R}_{\bm{U}}\right\Vert _{2} & =\left\Vert \bm{U}_{m,\cdot}\bm{R}_{\bm{U}}\left(\bm{\Sigma}^{\natural2}-\bm{R}_{\bm{U}}^{\top}\bm{\Sigma}^{2}\bm{R}_{\bm{U}}\right)\right\Vert _{2}\leq\left\Vert \bm{U}_{m,\cdot}\right\Vert _{2}\left\Vert \bm{\Sigma}^{\natural2}-\bm{R}_{\bm{U}}^{\top}\bm{\Sigma}^{2}\bm{R}_{\bm{U}}\right\Vert \nonumber \\
 & \lesssim\left(\left\Vert \bm{U}_{m,\cdot}^{\natural}\right\Vert _{2}+\frac{\zeta_{\mathsf{op},m}}{\sigma_{r}^{\natural2}}\sqrt{\frac{\mu^{\natural}r}{n_{1}}}\right)\left\Vert \bm{\Sigma}^{\natural2}-\bm{R}_{\bm{U}}^{\top}\bm{\Sigma}^{2}\bm{R}_{\bm{U}}\right\Vert \nonumber \\
 & \lesssim\left(\left\Vert \bm{U}_{m,\cdot}^{\natural}\right\Vert _{2}+\frac{\zeta_{\mathsf{op},m}}{\sigma_{r}^{\natural2}}\sqrt{\frac{\mu^{\natural}r}{n_{1}}}\right)\left(\kappa^{\natural2}\sqrt{\frac{\mu^{\natural}r}{n_{1}}}\zeta_{\mathsf{op}}+\kappa^{\natural2}\frac{\zeta_{\mathsf{op}}^{2}}{\sigma_{r}^{\natural2}}\right).\label{eq:hpca-inference-proof-inter-1}
\end{align}
Here, the penultimate inequality follows since --- according to Lemma
\ref{lemma:hpca-estimation-two-to-infty} --- we have 
\begin{align}
\left\Vert \bm{U}_{m,\cdot}\right\Vert _{2} & =\left\Vert \bm{U}_{m,\cdot}\bm{R}_{\bm{U}}\right\Vert _{2}\leq\left\Vert \bm{U}_{m,\cdot}^{\natural}\right\Vert _{2}+\left\Vert \bm{U}_{m,\cdot}\bm{R}_{\bm{U}}-\bm{U}_{m,\cdot}^{\natural}\right\Vert _{2}\nonumber \\
 & \lesssim\left\Vert \bm{U}_{m,\cdot}^{\natural}\right\Vert _{2}+\frac{\zeta_{\mathsf{op},m}}{\sigma_{r}^{\natural2}}\sqrt{\frac{\mu^{\natural}r}{n_{1}}}+\left\Vert \bm{U}_{m,\cdot}^{\natural}\right\Vert _{2}\left(\kappa^{\natural2}\sqrt{\frac{\mu^{\natural}r}{n_{1}}}\frac{\zeta_{\mathsf{op}}}{\sigma_{r}^{\natural2}}+\kappa^{\natural2}\frac{\zeta_{\mathsf{op}}^{2}}{\sigma_{r}^{\natural4}}\right)\nonumber \\
 & \lesssim\left\Vert \bm{U}_{m,\cdot}^{\natural}\right\Vert _{2}+\frac{\zeta_{\mathsf{op},m}}{\sigma_{r}^{\natural2}}\sqrt{\frac{\mu^{\natural}r}{n_{1}}}\label{eq:hpca-U-two-to-infty}
\end{align}
under the condition $\zeta_{\mathsf{op}}/\sigma_{r}^{\natural2}\lesssim1/\kappa^{\natural2}$
and $n_{1}\gtrsim\mu^{\natural}r$, while the last relation arises
from Lemma \ref{lemma:hpca-approx-2}.

\subsubsection{Step 2: replacing $\bm{U}\bm{\Sigma}^{2}\bm{R}_{\bm{U}}$ with $\bm{U}\bm{\Sigma}^{2}\bm{H}$.}

Given that $\bm{R}_{\bm{U}}$ and $\bm{H}$ are fairly close, one
can expect that replacing $\bm{R}_{\bm{U}}$ with $\bm{H}$ in $\bm{U}\bm{\Sigma}^{2}\bm{R}_{\bm{U}}$
does not change the matrix by much. To formalize this, we invoke Lemma
\ref{lemma:hpca-basic-facts} to reach 
\begin{align}
\left\Vert \bm{U}_{m,\cdot}\bm{\Sigma}^{2}\bm{H}-\bm{U}_{m,\cdot}\bm{\Sigma}^{2}\bm{R}_{\bm{U}}\right\Vert _{2} & \leq\left\Vert \bm{U}_{m,\cdot}\right\Vert _{2}\left\Vert \bm{\Sigma}^{2}\right\Vert \left\Vert \bm{H}-\bm{R}_{\bm{U}}\right\Vert \nonumber \\
 & \lesssim\left(\left\Vert \bm{U}_{m,\cdot}^{\natural}\right\Vert _{2}+\frac{\zeta_{\mathsf{op},m}}{\sigma_{r}^{\natural2}}\sqrt{\frac{\mu^{\natural}r}{n_{1}}}\right)\sigma_{1}^{\natural2}\frac{\zeta_{\mathsf{op}}^{2}}{\sigma_{r}^{\natural4}}\nonumber \\
 & \asymp\kappa^{\natural2}\frac{\zeta_{\mathsf{op}}^{2}}{\sigma_{r}^{\natural2}}\left(\left\Vert \bm{U}_{m,\cdot}^{\natural}\right\Vert _{2}+\frac{\zeta_{\mathsf{op},m}}{\sigma_{r}^{\natural2}}\sqrt{\frac{\mu^{\natural}r}{n_{1}}}\right).\label{eq:hpca-inference-proof-inter-2}
\end{align}
Here, the penultimate relation follows from Lemma \ref{lemma:hpca-basic-facts},
(\ref{eq:hpca-U-two-to-infty}), as well as a direct consequence of
Lemma \ref{lemma:hpca-approx-2}: 
\begin{equation}
\left\Vert \bm{\Sigma}^{2}\right\Vert =\left\Vert \bm{R}_{\bm{U}}^{\top}\bm{\Sigma}^{2}\bm{R}_{\bm{U}}\right\Vert \leq\left\Vert \bm{\Sigma}^{\natural2}\right\Vert +\left\Vert \bm{R}_{\bm{U}}^{\top}\bm{\Sigma}^{2}\bm{R}_{\bm{U}}-\bm{\Sigma}^{\natural2}\right\Vert \leq\sigma_{1}^{\natural2}+\kappa^{\natural2}\sqrt{\frac{\mu^{\natural}r}{n_{1}}}\zeta_{\mathsf{op}}+\kappa^{\natural2}\frac{\zeta_{\mathsf{op}}^{2}}{\sigma_{r}^{\natural2}}\asymp\sigma_{1}^{\natural2},\label{eq:hpca-Sigma-spectral}
\end{equation}
provided that $\zeta_{\mathsf{op}}/\sigma_{r}^{\natural2}\lesssim1$
and $n_{1}\gtrsim\mu^{\natural}r$.

\subsubsection{Step 3: establishing the proximity of $\bm{U}\bm{\Sigma}^{2}\bm{H}$
and $\bm{G}\bm{U}^{\natural}$.}

It is readily seen from Lemma \ref{lemma:hpca-approx-1} that 
\begin{align}
\left\Vert \left(\bm{U}\bm{\Sigma}^{2}\bm{H}-\bm{G}\bm{U}^{\natural}\right)_{m,\cdot}\right\Vert _{2} & \lesssim\zeta_{\mathsf{op},m}\left(\kappa^{\natural2}\frac{\zeta_{\mathsf{op}}}{\sigma_{r}^{\natural2}}\sqrt{\frac{\mu^{\natural}r}{n_{1}}}+\left\Vert \bm{U}\bm{H}-\bm{U}^{\natural}\right\Vert _{2,\infty}\right)+\kappa^{\natural2}\frac{\zeta_{\mathsf{op}}^{2}}{\sigma_{r}^{\natural2}}\left\Vert \bm{U}_{m,\cdot}^{\natural}\right\Vert _{2}\nonumber \\
 & \quad+\kappa^{\natural2}\zeta_{\mathsf{op}}\left\Vert \bm{U}_{m,\cdot}^{\natural}\right\Vert _{2}\left\Vert \bm{U}_{m,\cdot}\bm{H}-\bm{U}_{m,\cdot}^{\natural}\right\Vert _{2}+\kappa^{\natural2}\zeta_{\mathsf{op}}\left\Vert \bm{U}_{m,\cdot}\bm{H}-\bm{U}_{m,\cdot}^{\natural}\right\Vert _{2}^{2}\nonumber \\
 & \lesssim\zeta_{\mathsf{op},m}\kappa^{\natural2}\frac{\zeta_{\mathsf{op}}}{\sigma_{r}^{\natural2}}\sqrt{\frac{\mu^{\natural}r}{n_{1}}}+\kappa^{\natural2}\frac{\zeta_{\mathsf{op}}^{2}}{\sigma_{r}^{\natural2}}\left\Vert \bm{U}_{m,\cdot}^{\natural}\right\Vert _{2}\nonumber \\
 & \quad+\kappa^{\natural2}\zeta_{\mathsf{op}}\left\Vert \bm{U}_{m,\cdot}^{\natural}\right\Vert _{2}\left(\frac{\zeta_{\mathsf{op},m}}{\sigma_{r}^{\natural2}}\sqrt{\frac{\mu^{\natural}r}{n_{1}}}+\left\Vert \bm{U}_{m,\cdot}^{\natural}\right\Vert _{2}\kappa^{\natural2}\sqrt{\frac{\mu^{\natural}r}{n_{1}}}\frac{\zeta_{\mathsf{op}}}{\sigma_{r}^{\natural2}}+\left\Vert \bm{U}_{m,\cdot}^{\natural}\right\Vert _{2}\kappa^{\natural2}\frac{\zeta_{\mathsf{op}}^{2}}{\sigma_{r}^{\natural4}}\right)\nonumber \\
 & \quad+\kappa^{\natural2}\zeta_{\mathsf{op}}\left(\frac{\zeta_{\mathsf{op},m}}{\sigma_{r}^{\natural2}}\sqrt{\frac{\mu^{\natural}r}{n_{1}}}+\left\Vert \bm{U}_{m,\cdot}^{\natural}\right\Vert _{2}\kappa^{\natural2}\sqrt{\frac{\mu^{\natural}r}{n_{1}}}\frac{\zeta_{\mathsf{op}}}{\sigma_{r}^{\natural2}}+\left\Vert \bm{U}_{m,\cdot}^{\natural}\right\Vert _{2}\kappa^{\natural2}\frac{\zeta_{\mathsf{op}}^{2}}{\sigma_{r}^{\natural4}}\right)^{2}\nonumber \\
 & \lesssim\zeta_{\mathsf{op},m}\kappa^{\natural2}\frac{\zeta_{\mathsf{op}}}{\sigma_{r}^{\natural2}}\sqrt{\frac{\mu^{\natural}r}{n_{1}}}+\kappa^{\natural2}\frac{\zeta_{\mathsf{op}}^{2}}{\sigma_{r}^{\natural2}}\left\Vert \bm{U}_{m,\cdot}^{\natural}\right\Vert _{2},\label{eq:hpca-inference-proof-inter-3}
\end{align}
where the penultimate relation follows from Lemma \ref{lemma:hpca-estimation-two-to-infty},
while the last relation holds provided that $n_{1}\gtrsim\kappa^{\natural4}\mu^{\natural}r$
and $\zeta_{\mathsf{op}}\lesssim\sigma_{r}^{\natural2}$.

\subsubsection{Step 4: investigating the statistical properties of $\bm{G}\bm{U}^{\natural}$.}

It then remains to decompose $\bm{G}\bm{U}^{\natural}$ as claimed
in (\ref{eq:step-1234-subspace}). To begin with, we make the observation
that 
\begin{align*}
\bm{G} & =\mathcal{P}_{\mathsf{off}\text{-}\mathsf{diag}}\left[\left(\bm{M}^{\natural}+\bm{E}\right)\left(\bm{M}^{\natural}+\bm{E}\right)^{\top}\right]+\mathcal{P}_{\mathsf{diag}}\left(\bm{G}^{\natural}\right)+\mathcal{P}_{\mathsf{diag}}\left(\bm{G}-\bm{G}^{\natural}\right)\\
 & =\bm{G}^{\natural}+\mathcal{P}_{\mathsf{off}\text{-}\mathsf{diag}}\left[\bm{E}\bm{M}^{\natural\top}+\bm{M}^{\natural}\bm{E}^{\top}+\bm{E}\bm{E}^{\top}\right]+\mathcal{P}_{\mathsf{diag}}\left(\bm{G}-\bm{G}^{\natural}\right),
\end{align*}
which together with the eigen-decomposition $\bm{G}^{\natural}=\bm{U}^{\natural}\bm{\Sigma}^{\natural2}\bm{U}^{\natural\top}$
and the definition (\ref{eq:defn-Z-subspace-proof}) of $\bm{Z}$
allows one to derive 
\begin{align*}
\bm{G}\bm{U}^{\natural}-\left(\bm{U}^{\natural}+\bm{Z}\right)\bm{\Sigma}^{\natural2} & =\bm{U}^{\natural}\bm{\Sigma}^{\natural2}+\mathcal{P}_{\mathsf{off}\text{-}\mathsf{diag}}\left[\bm{E}\bm{M}^{\natural\top}+\bm{M}^{\natural}\bm{E}^{\top}+\bm{E}\bm{E}^{\top}\right]\bm{U}^{\natural}+\mathcal{P}_{\mathsf{diag}}\left(\bm{G}-\bm{G}^{\natural}\right)\bm{U}^{\natural}\\
 & \quad-\bm{U}^{\natural}\bm{\Sigma}^{\natural2}-\left[\bm{E}\bm{M}^{\natural\top}+\mathcal{P}_{\mathsf{off}\text{-}\mathsf{diag}}\left(\bm{E}\bm{E}^{\top}\right)\right]\bm{U}^{\natural}\\
 & =\underbrace{\bm{M}^{\natural}\bm{E}^{\top}\bm{U}^{\natural}}_{\eqqcolon\,\bm{R}_{1}}-\underbrace{\mathcal{P}_{\mathsf{diag}}\left(\bm{E}\bm{M}^{\natural\top}+\bm{M}^{\natural}\bm{E}^{\top}\right)\bm{U}^{\natural}}_{\eqqcolon\,\bm{R}_{2}}+\underbrace{\mathcal{P}_{\mathsf{diag}}\left(\bm{G}-\bm{G}^{\natural}\right)\bm{U}^{\natural}}_{\eqqcolon\,\bm{R}_{3}}.
\end{align*}
This motivates us to control the terms $\bm{R}_{1}$, $\bm{R}_{2}$
and $\bm{R}_{3}$ separately. 
\begin{itemize}
\item Regarding $\bm{R}_{1}$, note that we have shown in (\ref{eq:U-top-E-V-spectral})
that 
\[
\left\Vert \bm{U}^{\natural\top}\bm{E}\bm{V}^{\natural}\right\Vert \lesssim\sigma\sqrt{r\log n}+\frac{B\mu^{\natural}r\log n}{\sqrt{n_{1}n_{2}}}\lesssim\sigma\sqrt{r\log n}+\frac{\sigma\mu^{\natural}r}{\sqrt[4]{n_{1}n_{2}}}
\]
with probability exceeding $1-O(n^{-10})$, where the last inequality
holds under Assumption \ref{assumption:subspace-noise}. Therefore,
we can derive 
\begin{align*}
\left\Vert \bm{e}_{m}^{\top}\bm{R}_{1}\right\Vert _{2} & =\left\Vert \bm{U}_{m,\cdot}^{\natural}\bm{\Sigma}^{\natural}\left(\bm{U}^{\natural\top}\bm{E}\bm{V}^{\natural}\right)^{\top}\right\Vert _{2}\leq\left\Vert \bm{U}_{m,\cdot}^{\natural}\right\Vert _{2}\left\Vert \bm{\Sigma}^{\natural}\right\Vert \left\Vert \bm{U}^{\natural\top}\bm{E}\bm{V}^{\natural}\right\Vert \\
 & \lesssim\left\Vert \bm{U}_{m,\cdot}^{\natural}\right\Vert _{2}\sigma_{1}^{\natural}\left(\sigma\sqrt{r\log n}+\frac{\sigma\mu^{\natural}r}{\sqrt[4]{n_{1}n_{2}}}\right)\lesssim\left\Vert \bm{U}_{m,\cdot}^{\natural}\right\Vert _{2}\sqrt{\frac{\mu^{\natural}r}{n_{1}}}\zeta_{\mathsf{op}},
\end{align*}
where the last relation holds provided that $n_{1}n_{2}\gtrsim\mu^{\natural2}r^{2}$. 
\item With regards to $\bm{R}_{2}$, we observe that 
\[
\left\Vert \bm{e}_{m}^{\top}\bm{R}_{2}\right\Vert _{2}=\left|\sum_{j=1}^{n_{2}}E_{m,j}M_{m,j}^{\natural}\right|\left\Vert \bm{U}_{m,\cdot}^{\natural}\right\Vert _{2}
\]
Recalling from (\ref{eq:pca-Eij-Aij-sum}) that 
\[
\max_{1\leq i\leq n_{1}}\left|\sum_{j=1}^{n_{2}}E_{i,j}M_{i,j}^{\natural}\right|\lesssim\frac{\sqrt{\mu^{\natural}r}}{n_{1}}\zeta_{\mathsf{op}}
\]
holds with probability exceeding $1-O(n^{-10})$, we can upper bound
\begin{align*}
\left\Vert \bm{e}_{m}^{\top}\bm{R}_{2}\right\Vert _{2} & \lesssim\left\Vert \bm{U}_{m,\cdot}^{\natural}\right\Vert _{2}\frac{\sqrt{\mu^{\natural}r}}{n_{1}}\zeta_{\mathsf{op}}.
\end{align*}
\item It remains to control $\Vert\bm{R}_{3}\Vert_{2,\infty}$, towards
which we can apply Lemma \ref{lemma:hpca-estimation} to reach 
\[
\left\Vert \bm{e}_{m}^{\top}\bm{R}_{3}\right\Vert _{2}=\left|G_{m,m}-G_{m,m}^{\natural}\right|\left\Vert \bm{U}_{m,\cdot}^{\natural}\right\Vert _{2}\lesssim\kappa^{\natural2}\sqrt{\frac{\mu^{\natural}r}{n_{1}}}\zeta_{\mathsf{op}}\left\Vert \bm{U}_{m,\cdot}^{\natural}\right\Vert _{2}.
\]
\end{itemize}
Taking the preceding bounds on $\Vert\bm{e}_{m}^{\top}\bm{R}_{1}\Vert_{2}$,
$\Vert\bm{e}_{m}^{\top}\bm{R}_{2}\Vert_{2}$ and $\Vert\bm{e}_{m}^{\top}\bm{R}_{3}\Vert_{2}$
collectively, we arrive at 
\begin{equation}
\left\Vert \bm{G}_{m,\cdot}\bm{U}^{\natural}-\left(\bm{U}^{\natural}+\bm{Z}\right)_{m,\cdot}\bm{\Sigma}^{\natural2}\right\Vert _{2}\leq\left\Vert \bm{e}_{m}^{\top}\bm{R}_{1}\right\Vert _{2}+\left\Vert \bm{e}_{m}^{\top}\bm{R}_{2}\right\Vert _{2}+\left\Vert \bm{e}_{m}^{\top}\bm{R}_{3}\right\Vert _{2}\lesssim\kappa^{\natural2}\sqrt{\frac{\mu^{\natural}r}{n_{1}}}\zeta_{\mathsf{op}}\left\Vert \bm{U}_{m,\cdot}^{\natural}\right\Vert _{2}.\label{eq:hpca-inference-proof-inter-4}
\end{equation}

\subsubsection{Step 5: putting everything together.}

To finish up, combine (\ref{eq:hpca-inference-proof-inter-1}), (\ref{eq:hpca-inference-proof-inter-2}),
(\ref{eq:hpca-inference-proof-inter-3}) and (\ref{eq:hpca-inference-proof-inter-4})
to conclude that 
\begin{align*}
\left\Vert \left(\bm{U}\bm{R}_{\bm{U}}-\bm{U}^{\star}-\bm{Z}\right)_{m,\cdot}\right\Vert _{2} & \leq\frac{1}{\sigma_{r}^{\natural2}}\left\Vert \bm{U}_{m,\cdot}\bm{R}_{\bm{U}}\bm{\Sigma}^{\natural2}-\left(\bm{U}^{\natural}+\bm{Z}\right)_{m,\cdot}\bm{\Sigma}^{\natural2}\right\Vert _{2}\\
 & \leq\frac{1}{\sigma_{r}^{\natural2}}\left\Vert \bm{U}_{m,\cdot}\bm{R}_{\bm{U}}\bm{\Sigma}^{\natural2}-\bm{U}_{m,\cdot}\bm{\Sigma}^{2}\bm{R}_{\bm{U}}\right\Vert _{2}+\frac{1}{\sigma_{r}^{\natural2}}\left\Vert \bm{U}_{m,\cdot}\bm{\Sigma}^{2}\bm{H}-\bm{U}_{m,\cdot}\bm{\Sigma}^{2}\bm{R}_{\bm{U}}\right\Vert _{2}\\
 & \quad+\frac{1}{\sigma_{r}^{\natural2}}\left\Vert \bm{U}_{m,\cdot}\bm{\Sigma}^{2}\bm{H}-\bm{G}_{m,\cdot}\bm{U}^{\natural}\right\Vert _{2}+\frac{1}{\sigma_{r}^{\natural2}}\left\Vert \bm{G}_{m,\cdot}\bm{U}^{\natural}-\left(\bm{U}^{\natural}+\bm{Z}\right)_{m,\cdot}\bm{\Sigma}^{\natural2}\right\Vert _{2}\\
 & \lesssim\left\Vert \bm{U}_{m,\cdot}^{\natural}\right\Vert _{2}\left(\kappa^{\natural2}\sqrt{\frac{\mu^{\natural}r}{n_{1}}}\frac{\zeta_{\mathsf{op}}}{\sigma_{r}^{\natural2}}+\kappa^{\natural2}\frac{\zeta_{\mathsf{op}}^{2}}{\sigma_{r}^{\natural4}}\right)+\kappa^{\natural2}\frac{\zeta_{\mathsf{op}}\zeta_{\mathsf{op},m}}{\sigma_{r}^{\natural4}}\sqrt{\frac{\mu^{\natural}r}{n_{1}}}
\end{align*}
as claimed, provided that $\zeta_{\mathsf{op}}\lesssim\sigma_{r}^{\natural2}$

\subsection{Proof of auxiliary lemmas}

In this section, we establish the auxiliary lemmas needed in the proof
of Theorem \ref{thm:hpca_inference_general}. At the core of our analysis
lies a leave-one-out (and leave-two-out) analysis framework that has
been previously adopted to analyze spectral methods \citep{abbe2017entrywise,chen2017spectral,cai2019subspace,chen2020spectral},
which we shall introduce below.

\subsubsection{Preparation: leave-one-out and leave-two-out auxiliary estimates}

To begin with, let us introduce the leave-one-out and leave-two-out
auxiliary sequences and the associated estimates, which play a pivotal
role in enabling fine-grained statistical analysis.

\paragraph{Leave-one-out auxiliary estimates. }

For each $m\in[n_{1}]$, define 
\begin{equation}
\bm{M}^{(m)}\coloneqq\mathcal{P}_{-m,\cdot}\left(\bm{M}\right)+\mathcal{P}_{m,\cdot}\left(\bm{M}^{\natural}\right)\qquad\text{and}\qquad\bm{G}^{(m)}\coloneqq\mathcal{P}_{\mathsf{off}\text{-}\mathsf{diag}}\left(\bm{M}^{(m)}\bm{M}^{(m)\top}\right)+\mathcal{P}_{\mathsf{diag}}\left(\bm{M}^{\natural}\bm{M}^{\natural\top}\right),\label{eq:defn-M-LOO-subspace}
\end{equation}
where we recall that $\mathcal{P}_{-m,\cdot}\left(\bm{M}\right)\in\mathbb{R}^{n_{1}\times n_{2}}$
is obtained by setting to zero all entries in the $m$-th row of $\bm{M}$,
and $\mathcal{P}_{m,\cdot}\left(\bm{M}\right)=\bm{M}-\mathcal{P}_{-m,\cdot}\left(\bm{M}\right)$.
Throughout this section, we let $\bm{U}^{(m)}\bm{\Lambda}^{(m)}\bm{U}^{(m)\top}$
be the top-$r$ eigen-decomposition of $\bm{G}^{(m)}$, and we define
\begin{equation}
\bm{H}^{(m)}\coloneqq\bm{U}^{(m)\top}\bm{U}^{\natural}.\label{eq:defn-Hm-subspace}
\end{equation}
Two features are particularly worth emphasizing: 
\begin{itemize}
\item The matrices $\bm{U}^{(m)}$, $\bm{G}^{(m)}$ and $\bm{H}^{(m)}$
are all statistically independent of the $m$-th row of the noise
matrix $\bm{E}$, given that $\bm{M}^{(m)}$ does not contain any
randomness arising in the $m$-th row of $\bm{E}$. 
\item The estimate $\bm{U}^{(m)}$ (resp.~$\bm{H}^{(m)}$ and $\bm{G}^{(m)}$)
is expected to be extremely close to the original estimate $\bm{U}$
(resp.~$\bm{H}=\bm{U}^{\top}\bm{U}^{\natural}$ and $\bm{G}$), given
that we have only dropped a small fraction of data when generating
the leave-one-out estimate. 
\end{itemize}
Informally, the above two features taken together allow one to show
the weak statistical dependency between $(\bm{U},\bm{G},\bm{H})$
and the $m$-th row of $\bm{E}$.

\paragraph{Leave-two-out auxiliary estimates. }

As it turns out, we are also in need of a collection of slightly more
complicated leave-two-out estimates to assist in our analysis. Specifically,
for each $m\in[n_{1}]$ and $l\in[n_{2}]$, we let 
\begin{equation}
\bm{M}^{(m,l)}=\mathcal{P}_{-m,-l}\left(\bm{M}\right)+\mathcal{P}_{m,l}\left(\bm{M}^{\natural}\right)\qquad\text{and}\qquad\bm{G}^{(m,l)}=\mathcal{P}_{\mathsf{off}\text{-}\mathsf{diag}}\left(\bm{M}^{(m,l)}\bm{M}^{(m,l)\top}\right)+\mathcal{P}_{\mathsf{diag}}\left(\bm{M}^{\natural}\bm{M}^{\natural\top}\right).\label{eq:defn-M-leave-two-out}
\end{equation}
Here, $\mathcal{P}_{-m,-l}\left(\bm{M}\right)\in\mathbb{R}^{n_{1}\times n_{2}}$
is obtained by zeroing out all entries in the $m$-th row and the
$l$-th column of $\bm{M}$, and $\mathcal{P}_{m,l}\left(\bm{M}\right)=\bm{M}-\mathcal{P}_{-m,-l}\left(\bm{M}\right)$.
We let $\bm{U}^{(m,l)}\bm{\Lambda}^{(m,l)}\bm{U}^{(m,l)\top}$ represent
the top-$r$ eigen-decomposition of $\bm{G}^{(m,l)}$, and define
\begin{equation}
\bm{H}^{(m,l)}\coloneqq\bm{U}^{(m,l)\top}\bm{U}^{\natural}.\label{eq:defn-H-ml-subspace}
\end{equation}
Akin to the leave-one-out counterpart, the matrices $\bm{U}^{(m,l)}$,
$\bm{G}^{(m,l)}$ and $\bm{H}^{(m,l)}$ we construct are all statistically
independent of the $m$-th row and the $l$-th column of $\bm{E}$,
and we shall also exploit the proximity of $\bm{U}^{(m,l)}$, $\bm{U}^{(m)}$
and $\bm{U}$ in the subsequent analysis.

\paragraph{Useful properties concerning these auxiliary estimates. }

In the sequel, we collect a couple of useful lemmas that are concerned
with the leave-one-out and leave-two-out estimates, several of which
are adapted from \citet{cai2019subspace}. The first lemma controls
the difference between the leave-one-out gram matrix $\bm{G}^{(m)}$
and the original gram matrix $\bm{G}$, as well as the difference
between $\bm{G}^{(m)}$ and the ground-truth gram matrix.

\begin{lemma}\label{lemma:hpca-loo-basics}Suppose that $n_{1}\gtrsim\kappa^{\natural4}\mu^{\natural}r$.
Then with probability exceeding $1-O(n^{-10})$, 
\begin{align*}
\left\Vert \bm{G}^{(m)}-\bm{G}\right\Vert  & \lesssim\sigma^{2}\sqrt{n_{1}n_{2}}\log n+\kappa^{\natural2}\sigma\sigma_{1}^{\natural}\sqrt{\mu^{\natural}r\log n}\\
\left\Vert \bm{G}^{(m)}-\bm{G}^{\natural}\right\Vert  & \lesssim\zeta_{\mathsf{op}}\asymp\sigma^{2}\sqrt{n_{1}n_{2}}\log n+\sigma\sigma_{1}^{\natural}\sqrt{n_{1}\log n}
\end{align*}
hold simultaneously for all $1\leq m\leq n_{1}$.\end{lemma}\begin{proof}See
Appendix \ref{sec:proof-lemma-HeteroPCA-loo-basics}.\end{proof}

The next lemma confirms that the leave-one-out estimate $\bm{U}^{(m)}$
and the leave-two-out estimate $\bm{U}^{(m,l)}$ are exceedingly close.

\begin{lemma}\label{lemma:hpca-lto-perturbation}Suppose that the
assumptions of Theorem \ref{thm:hpca_inference_general} hold. Then
with probability exceeding $1-O(n^{-10})$, we have 
\begin{align*}
\left\Vert \bm{U}^{(m)}\bm{U}^{(m)\top}-\bm{U}^{(m,l)}\bm{U}^{(m,l)\top}\right\Vert  & \lesssim\frac{1}{\sigma_{r}^{\natural2}}\left(B\log n+\sigma\sqrt{n_{1}\log n}\right)^{2}\left\Vert \bm{U}^{(m)}\bm{H}^{(m)}\right\Vert _{2,\infty}+\frac{\sigma^{2}}{\sigma_{r}^{\natural2}}\\
 & \quad+\frac{1}{\sigma_{r}^{\natural2}}\left(B\log n+\sigma\sqrt{n_{1}\log n}\right)\left\Vert \bm{M}^{\natural\top}\right\Vert _{2,\infty}
\end{align*}
simultaneously for all $1\leq m\leq n_{1}$ and $1\leq l\leq n_{2}$.
\end{lemma}\begin{proof}See Appendix \ref{sec:proof-lemma-HeteroPCA-lto-perturbation}.\end{proof}

The following two lemmas study the size of $\bm{U}^{(m)}\bm{H}^{(m)}-\bm{U}^{\natural}$
and $\bm{U}^{(m)}\bm{H}^{(m)}$ when projected towards several important
directions.

\begin{lemma}\label{lemma:hpca-useful-2}Suppose that the assumptions
of Theorem \ref{thm:hpca_inference_general} hold. Then with probability
exceeding $1-O(n^{-10})$, we have, for all $l\in[n_{2}]$ and all
$m\in[n_{1}]$, 
\begin{align}
 & \left\Vert \bm{e}_{l}^{\top}\left[\mathcal{P}_{-m,\cdot}\left(\bm{M}\right)\right]^{\top}\left(\bm{U}^{(m)}\bm{H}^{(m)}-\bm{U}^{\natural}\right)\right\Vert _{2}\lesssim\frac{\zeta_{\mathsf{op}}^{2}}{\sigma_{r}^{\natural4}}\left\Vert \bm{M}^{\natural\top}\right\Vert _{2,\infty}+\left(B\log n+\sigma\sqrt{n_{1}\log n}\right)\left\Vert \bm{U}^{(m)}\bm{H}^{(m)}-\bm{U}^{\natural}\right\Vert _{2,\infty}\nonumber \\
 & \quad\qquad\qquad\qquad+\left(\left\Vert \bm{M}^{\natural\top}\right\Vert _{2,\infty}+B\log n+\sigma\sqrt{n_{1}\log n}\right)\left\Vert \bm{U}^{(m)}\bm{U}^{(m)\top}-\bm{U}^{(m,l)}\bm{U}^{(m,l)\top}\right\Vert ;\label{eq:hpca-useful-2-1}
\end{align}
\begin{equation}
\left\Vert \left[\mathcal{P}_{-m,\cdot}\left(\bm{M}\right)\right]^{\top}\bm{U}^{(m)}\bm{H}^{(m)}\right\Vert _{2,\infty}\lesssim\left(B\log n+\sigma\sqrt{n_{1}\log n}\right)\left\Vert \bm{U}^{(m)}\bm{H}^{(m)}\right\Vert _{2,\infty}+\sigma_{1}^{\natural}\sqrt{\frac{\mu^{\natural}r}{n_{2}}}.\label{eq:hpca-useful-2-2}
\end{equation}
\end{lemma}\begin{proof}See Appendix \ref{subsec:proof-lemma-pca-useful-2}.\end{proof}

\begin{lemma}\label{lemma:hpca-useful-1} Suppose that the assumptions
of Theorem \ref{thm:hpca_inference_general} hold. Then with probability
exceeding $1-O(n^{-10})$, we have 
\begin{align}
 & \left\Vert \bm{E}_{m,\cdot}\left[\mathcal{P}_{-m,\cdot}\left(\bm{M}\right)\right]^{\top}\left(\bm{U}^{(m)}\bm{H}^{(m)}-\bm{U}^{\natural}\right)\right\Vert _{2}\lesssim\frac{\zeta_{\mathsf{op}}}{\sigma_{r}^{\natural2}}\sigma_{m}\sigma_{1}^{\natural}\sqrt{n_{1}\log n}\sqrt{\frac{\mu^{\natural}r}{n_{1}}}\nonumber \\
 & \qquad\qquad\qquad+\zeta_{\mathsf{op},m}\left(\left\Vert \bm{U}^{(m)}\bm{H}^{(m)}-\bm{U}^{\natural}\right\Vert _{2,\infty}+\frac{\zeta_{\mathsf{op}}}{\sigma_{r}^{\natural2}}\left\Vert \bm{U}^{(m)}\bm{H}^{(m)}\right\Vert _{2,\infty}\right);\label{eq:hpca-useful-1-1}
\end{align}
\begin{align}
\left\Vert \bm{E}_{m,\cdot}\left[\mathcal{P}_{-m,\cdot}\left(\bm{M}\right)\right]^{\top}\bm{U}^{(m)}\bm{H}^{(m)}\right\Vert _{2} & \lesssim\zeta_{\mathsf{op}}\left(\left\Vert \bm{U}^{(m)}\bm{H}^{(m)}\right\Vert _{2,\infty}+\sqrt{\frac{\mu^{\natural}r}{n_{1}}}\right)\label{eq:hpca-useful-1-2}
\end{align}
simultaneously for all $1\leq m\leq n_{1}$. \end{lemma}\begin{proof}See
Appendix \ref{subsec:proof-lemma-pca-useful-1}.\end{proof}

Finally, we justify the proximity of the leave-one-out estimate $\bm{U}^{(m)}$
and the original estimate $\bm{U}$, which turns out to be a consequence
of the preceding results.

\begin{lemma}\label{lemma:hpca-loo-perturbation}Suppose that the
assumptions of Theorem \ref{thm:hpca_inference_general} hold. Then
with probability exceeding $1-O(n^{-10})$, 
\[
\left\Vert \bm{U}^{(m)}\bm{U}^{(m)\top}-\bm{U}\bm{U}^{\top}\right\Vert \lesssim\kappa^{\natural2}\frac{\zeta_{\mathsf{op}}}{\sigma_{r}^{\natural2}}\left(\left\Vert \bm{U}\bm{H}\right\Vert _{2,\infty}+\sqrt{\frac{\mu^{\natural}r}{n_{1}}}\right)
\]
holds simultaneously for all $1\leq m\leq n_{1}$.\end{lemma}\begin{proof}See
Appendix \ref{sec:proof-lemma-HeteroPCA-loo-perturbation}.\end{proof}

\subsubsection{Proof of Lemma \ref{lemma:hpca-estimation}\label{appendix:proof-lemma-HeteroPCA-estimation}}

Recall that the initialization $\bm{G}^{0}$ is obtained by dropping
all diagonal entries of $\bm{M}\bm{M}^{\top}$, namely, 
\[
\bm{G}^{0}=\mathcal{P}_{\mathsf{off}\text{-}\mathsf{diag}}\left(\bm{M}\bm{M}^{\top}\right),
\]
and $\bm{U}^{0}$ consists of the top-$r$ eigenvectors of $\bm{G}^{0}$.
These are precisely the subjects studied in \citet{cai2019subspace}.
In light of this, we state below the following two lemmas borrowed
from \citet{cai2019subspace}, which assist in proving Lemma \ref{lemma:hpca-estimation}.

\begin{lemma}\label{lemma:HPCA-G0-spectral}Suppose that the assumptions
of Theorem \ref{thm:hpca_inference_general} hold. With probability
exceeding $1-O(n^{-10})$, 
\begin{align*}
\left\Vert \bm{G}^{0}-\mathcal{P}_{\mathsf{off}\text{-}\mathsf{diag}}\left(\bm{G}^{\natural}\right)\right\Vert  & =\left\Vert \mathcal{P}_{\mathsf{off}\text{-}\mathsf{diag}}\left(\bm{G}^{0}-\bm{G}^{\natural}\right)\right\Vert \leq C^{\mathsf{dd}}\zeta_{\mathsf{op}}\\
\left\Vert \bm{G}^{0}-\bm{G}^{\natural}\right\Vert  & \leq C^{\mathsf{dd}}\zeta_{\mathsf{op}}+\left\Vert \bm{M}^{\natural}\right\Vert _{2,\infty}^{2}
\end{align*}
hold for some universal constant $C^{\mathsf{dd}}>0$. \end{lemma}\begin{proof}See
\citet[Lemma 1]{cai2019subspace}.\end{proof}

\begin{lemma}\label{lemma:HeteroPCA-U0-results}Suppose that the
assumptions of Theorem \ref{thm:hpca_inference_general} hold. With
probability exceeding $1-O(n^{-10})$, 
\begin{align*}
\left\Vert \bm{U}^{0}\bm{R}^{0}-\bm{U}^{\natural}\right\Vert  & \leq C_{\mathsf{op}}^{\mathsf{dd}}\frac{\zeta_{\mathsf{op}}+\left\Vert \bm{M}^{\natural}\right\Vert _{2,\infty}^{2}}{\sigma_{r}^{\natural2}}\\
\left\Vert \bm{U}^{0}\bm{R}^{0}-\bm{U}^{\natural}\right\Vert _{2,\infty} & \leq C_{\infty}^{\mathsf{dd}}\frac{\kappa^{\natural2}\left(\zeta_{\mathsf{op}}+\left\Vert \bm{M}^{\natural}\right\Vert _{2,\infty}^{2}\right)}{\sigma_{r}^{\natural2}}\sqrt{\frac{\mu^{\natural}r}{n_{1}}}
\end{align*}
hold for some universal constants $C_{\mathsf{op}}^{\mathsf{dd}},C_{\infty}^{\mathsf{dd}}>0$.
Here, $\bm{R}_{0}$ is defined to be the following rotation matrix
\begin{equation}
\bm{R}^{0}\coloneqq\arg\min_{\bm{O}\in\mathcal{O}^{r\times r}}\left\Vert \bm{U}^{0}\bm{O}-\bm{U}^{\natural}\right\Vert _{\mathrm{F}}^{2}.\label{eq:defn-R0-subspace-est}
\end{equation}
\end{lemma}\begin{proof}See \citet[Theorem 1]{cai2019subspace}.\end{proof}

Armed with the above lemmas, we are ready to show by induction that:
for $1\leq s\leq t_{0}$, 
\begin{align}
\left\Vert \mathcal{P}_{\mathsf{diag}}\left(\bm{G}^{s}-\bm{G}^{\natural}\right)\right\Vert  & \leq C_{0}\kappa^{\natural2}\sqrt{\frac{\mu^{\natural}r}{n_{1}}}\left\Vert \bm{G}^{s-1}-\bm{G}^{\natural}\right\Vert \label{eq:HeteroPCA-induction}
\end{align}
for some sufficiently large constant $C_{0}>0$.

\paragraph{Base case: $s=1$.}

Let us start with the base case by making the observation that 
\begin{align*}
\left\Vert \mathcal{P}_{\mathsf{diag}}\left(\bm{G}^{1}-\bm{G}^{\natural}\right)\right\Vert  & =\left\Vert \mathcal{P}_{\mathsf{diag}}\left(\bm{U}^{0}\bm{\Lambda}^{0}\bm{U}^{0\top}-\bm{G}^{\natural}\right)\right\Vert =\left\Vert \mathcal{P}_{\mathsf{diag}}\left(\mathcal{P}_{\bm{U}^{0}}(\bm{G}^{0})-\bm{G}^{\natural}\right)\right\Vert \\
 & \leq\underbrace{\left\Vert \mathcal{P}_{\mathsf{diag}}\left[\mathcal{P}_{\bm{U}^{0}}\left(\bm{G}^{0}-\bm{G}^{\natural}\right)\right]\right\Vert }_{\eqqcolon\alpha_{1}}+\underbrace{\left\Vert \mathcal{P}_{\mathsf{diag}}\left(\mathcal{P}_{\bm{U}_{\perp}^{0}}\bm{G}^{\natural}\right)\right\Vert }_{\eqqcolon\alpha_{2}},
\end{align*}
where $\bm{U}_{\perp}^{0}\in\mathbb{R}^{n_{1}\times(n_{1}-r)}$ is
a matrix whose orthonormal columns span the orthogonal complement
of the column space of $\bm{U}^{0}$. 
\begin{itemize}
\item Regarding the first term $\alpha_{1}$, we have seen from Lemma \ref{lemma:HeteroPCA-U0-results}
that 
\begin{align}
\left\Vert \bm{U}^{0}\right\Vert _{2,\infty} & =\left\Vert \bm{U}^{0}\bm{R}^{0}\right\Vert _{2,\infty}\leq\left\Vert \bm{U}^{0}\bm{R}^{0}-\bm{U}^{\natural}\right\Vert _{2,\infty}+\left\Vert \bm{U}^{\natural}\right\Vert _{2,\infty}\leq C_{\infty}^{\mathsf{dd}}\frac{\kappa^{\natural2}\left(\zeta_{\mathsf{op}}+\left\Vert \bm{M}^{\natural}\right\Vert _{2,\infty}^{2}\right)}{\sigma_{r}^{\natural2}}\sqrt{\frac{\mu^{\natural}r}{n_{1}}}+\sqrt{\frac{\mu^{\natural}r}{n_{1}}}\nonumber \\
 & \overset{\text{(i)}}{\leq}C_{\infty}^{\mathsf{dd}}\kappa^{\natural2}\left(\frac{\zeta_{\mathsf{op}}}{\sigma_{r}^{\natural2}}+\frac{\mu^{\natural}r}{d}\right)\sqrt{\frac{\mu^{\natural}r}{n_{1}}}+\sqrt{\frac{\mu^{\natural}r}{n_{1}}}\overset{\text{(ii)}}{\leq}\sqrt{\frac{4\mu^{\natural}r}{n_{1}}},\label{eq:U0-two-to-infty}
\end{align}
where (i) holds since 
\begin{equation}
\left\Vert \bm{M}^{\natural}\right\Vert _{2,\infty}=\left\Vert \bm{U}^{\natural}\bm{\Sigma}^{\natural}\bm{V}^{\natural\top}\right\Vert _{2,\infty}\leq\left\Vert \bm{U}^{\natural}\right\Vert _{2,\infty}\left\Vert \bm{\Sigma}^{\natural}\right\Vert \left\Vert \bm{V}^{\natural}\right\Vert \leq\sqrt{\frac{\mu^{\natural}r}{n_{1}}}\sigma_{1}^{\natural},\label{eq:M-natural-two-to-infty}
\end{equation}
and (ii) holds provided that $\zeta_{\mathsf{op}}/\sigma_{r}^{\natural2}\ll1/\kappa^{\natural2}$
and $d\gg\kappa^{\natural2}\mu^{\natural}r$. This in turn allows
us to use \citet[Lemma 1]{zhang2018heteroskedastic} to reach 
\begin{equation}
\alpha_{1}\leq2\sqrt{\frac{\mu^{\natural}r}{n_{1}}}\left\Vert \bm{G}^{0}-\bm{G}^{\natural}\right\Vert .\label{eq:HeteroPCA-inter1}
\end{equation}
\item Regarding the second term $\alpha_{2}$, invoke \citet[Lemma 1]{zhang2018heteroskedastic}
again to arrive at 
\begin{equation}
\alpha_{2}=\left\Vert \mathcal{P}_{\mathsf{diag}}\left(\mathcal{P}_{\bm{U}_{\perp}^{0}}\bm{G}^{\natural}\mathcal{P}_{\bm{U}^{\natural}}\right)\right\Vert \leq\sqrt{\frac{\mu^{\natural}r}{n_{1}}}\left\Vert \mathcal{P}_{\bm{U}_{\perp}^{0}}\bm{G}^{\natural}\right\Vert ,\label{eq:HeteroPCA-inter2}
\end{equation}
which has made use of the fact that $\bm{G}^{\natural}\mathcal{P}_{\bm{U}^{\natural}}=\bm{G}^{\natural}$.
Moreover, it is seen that 
\begin{align}
\left\Vert \mathcal{P}_{\bm{U}_{\perp}^{0}}\bm{G}^{\natural}\right\Vert  & =\left\Vert \bm{U}_{\perp}^{0}\bm{U}_{\perp}^{0\top}\bm{U}^{\natural}\bm{\Sigma}^{\natural2}\bm{U}^{\natural\top}\right\Vert \leq\sigma_{1}^{\natural2}\left\Vert \bm{U}_{\perp}^{0\top}\bm{U}^{\natural}\right\Vert =\sigma_{1}^{\natural2}\left\Vert \sin\bm{\Theta}\left(\bm{U}^{0},\bm{U}^{\natural}\right)\right\Vert ,\label{eq:HeteroPCA-inter3}
\end{align}
where $\bm{\Theta}\left(\bm{U}^{0},\bm{U}^{\natural}\right)$ is a
diagonal matrix whose diagonal entries correspond to the principal
angles between $\bm{U}^{0}$ and $\bm{U}^{\natural}$, and the last
identity follows from \citet[Lemma 2.1.2]{chen2020spectral}. In view
of the Davis-Kahan $\sin\bm{\Theta}$ Theorem \citep[Theorem 2.2.1]{chen2020spectral},
we can demonstrate that 
\begin{equation}
\left\Vert \sin\bm{\Theta}\left(\bm{U}^{0},\bm{U}^{\natural}\right)\right\Vert \leq\frac{\sqrt{2}\left\Vert \bm{G}^{0}-\bm{G}^{\natural}\right\Vert }{\lambda_{r}\left(\bm{G}^{0}\right)-\lambda_{r+1}\left(\bm{G}^{\natural}\right)}=\frac{\sqrt{2}\left\Vert \bm{G}^{0}-\bm{G}^{\natural}\right\Vert }{\lambda_{r}\left(\bm{G}^{0}\right)}\leq\frac{2\sqrt{2}\left\Vert \bm{G}^{0}-\bm{G}^{\natural}\right\Vert }{\sigma_{r}^{\natural2}}.\label{eq:HeteroPCA-inter4}
\end{equation}
Here, the identity comes from the fact $\lambda_{r+1}(\bm{G}^{\natural})=0$;
the last inequality follows from a direct application of Weyl's inequality:
\begin{align*}
\lambda_{r}\left(\bm{G}^{0}\right) & \geq\lambda_{r}\left(\bm{G}^{\natural}\right)-\left\Vert \bm{G}^{0}-\bm{G}^{\natural}\right\Vert \overset{\text{(i)}}{\geq}\sigma_{r}^{\natural2}-\zeta_{\mathsf{op}}-\left\Vert \bm{M}^{\natural}\right\Vert _{2,\infty}^{2}\overset{\text{(ii)}}{\geq}\sigma_{r}^{\natural2}-\zeta_{\mathsf{op}}-\frac{\mu^{\natural}r}{n_{1}}\sigma_{1}^{\natural2}\overset{\text{(iii)}}{\geq}\frac{1}{2}\sigma_{r}^{\natural2},
\end{align*}
where (i) is a consequence of Lemma \ref{lemma:HPCA-G0-spectral},
(ii) follows from (\ref{eq:M-natural-two-to-infty}), and (iii) holds
as long as $\zeta_{\mathsf{op}}\ll\sigma_{r}^{\natural2}$ and $n_{1}\gg\kappa^{\natural2}\mu^{\natural}r$.
Combine (\ref{eq:HeteroPCA-inter2}), (\ref{eq:HeteroPCA-inter3})
and (\ref{eq:HeteroPCA-inter4}) to reach 
\[
\alpha_{2}\leq2\sqrt{2}\kappa^{\natural2}\sqrt{\frac{\mu^{\natural}r}{n_{1}}}\left\Vert \bm{G}^{0}-\bm{G}^{\natural}\right\Vert .
\]
\end{itemize}
Combine the above bounds on $\alpha_{1}$ and $\alpha_{2}$ to yield
\[
\left\Vert \mathcal{P}_{\mathsf{diag}}\left(\bm{G}^{1}-\bm{G}^{\natural}\right)\right\Vert \leq\alpha_{1}+\alpha_{2}\leq\left(2+2\sqrt{2}\kappa^{\natural2}\right)\sqrt{\frac{\mu^{\natural}r}{n_{1}}}\left\Vert \bm{G}^{0}-\bm{G}^{\natural}\right\Vert \leq C_{0}\kappa^{\natural2}\sqrt{\frac{\mu^{\natural}r}{n_{1}}}\left\Vert \bm{G}^{0}-\bm{G}^{\natural}\right\Vert ,
\]
with the proviso that the constant $C_{0}$ is sufficiently large.

\paragraph{Induction step.}

For any given $t>1$, Suppose that (\ref{eq:HeteroPCA-induction})
holds for all $s=1,2\ldots,t$, and we'd like to show that it continues
to hold for $s=t+1$. From the induction hypothesis, we know that
for any $1\leq\tau\leq t$, one has 
\begin{align}
\left\Vert \mathcal{P}_{\mathsf{diag}}\left(\bm{G}^{\tau}-\bm{G}^{\natural}\right)\right\Vert  & \leq C_{0}\kappa^{\natural2}\sqrt{\frac{\mu^{\natural}r}{n_{1}}}\left\Vert \bm{G}^{\tau-1}-\bm{G}^{\natural}\right\Vert \leq C_{0}\kappa^{\natural2}\sqrt{\frac{\mu^{\natural}r}{n_{1}}}\left(\left\Vert \mathcal{P}_{\mathsf{diag}}\left(\bm{G}^{\tau-1}-\bm{G}^{\natural}\right)\right\Vert +\left\Vert \mathcal{P}_{\mathsf{off}\text{-}\mathsf{diag}}\left(\bm{G}^{\tau-1}-\bm{G}^{\natural}\right)\right\Vert \right)\nonumber \\
 & =C_{0}\kappa^{\natural2}\sqrt{\frac{\mu^{\natural}r}{n_{1}}}\left(\left\Vert \mathcal{P}_{\mathsf{diag}}\left(\bm{G}^{\tau-1}-\bm{G}^{\natural}\right)\right\Vert +\left\Vert \mathcal{P}_{\mathsf{off}\text{-}\mathsf{diag}}\left(\bm{G}^{0}-\bm{G}^{\natural}\right)\right\Vert \right),\label{eq:HeteroPCA_inter5-1}
\end{align}
where the last line holds since, by construction, $\mathcal{P}_{\mathsf{off}\text{-}\mathsf{diag}}\left(\bm{G}^{\tau}\right)=\mathcal{P}_{\mathsf{off}\text{-}\mathsf{diag}}\left(\bm{G}^{0}\right)$
for all $\tau$. Applying the above inequality recursively gives 
\begin{align}
\left\Vert \mathcal{P}_{\mathsf{diag}}\left(\bm{G}^{t}-\bm{G}^{\natural}\right)\right\Vert  & \leq C_{0}\kappa^{\natural2}\sqrt{\frac{\mu^{\natural}r}{n_{1}}}\left(\left\Vert \mathcal{P}_{\mathsf{diag}}\left(\bm{G}^{t-1}-\bm{G}^{\natural}\right)\right\Vert +\left\Vert \mathcal{P}_{\mathsf{off}\text{-}\mathsf{diag}}\left(\bm{G}^{0}-\bm{G}^{\natural}\right)\right\Vert \right)\leq\cdots\nonumber \\
 & \leq\left(C_{0}\kappa^{\natural2}\sqrt{\frac{\mu^{\natural}r}{n_{1}}}\right)^{t}\left\Vert \mathcal{P}_{\mathsf{diag}}\left(\bm{G}^{0}-\bm{G}^{\natural}\right)\right\Vert +\sum_{i=1}^{t}\left(C_{0}\kappa^{\natural2}\sqrt{\frac{\mu^{\natural}r}{n_{1}}}\right)^{i}\left\Vert \mathcal{P}_{\mathsf{off}\text{-}\mathsf{diag}}\left(\bm{G}^{0}-\bm{G}^{\natural}\right)\right\Vert \nonumber \\
 & \leq\left(C_{0}\kappa^{\natural2}\sqrt{\frac{\mu^{\natural}r}{n_{1}}}\right)^{t}\left\Vert \bm{M}^{\natural}\right\Vert _{2,\infty}^{2}+2C_{0}\kappa^{\natural2}\sqrt{\frac{\mu^{\natural}r}{n_{1}}}\left\Vert \mathcal{P}_{\mathsf{off}\text{-}\mathsf{diag}}\left(\bm{G}^{0}-\bm{G}^{\natural}\right)\right\Vert \nonumber \\
 & \leq\left(C_{0}\kappa^{\natural2}\sqrt{\frac{\mu^{\natural}r}{n_{1}}}\right)^{t}\frac{\mu^{\natural}r}{n_{1}}\sigma_{1}^{\natural2}+2C_{0}C^{\mathsf{dd}}\kappa^{\natural2}\sqrt{\frac{\mu^{\natural}r}{n_{1}}}\zeta_{\mathsf{op}}.\label{eq:HeteroPCA_inter5}
\end{align}
Here, the penultimate line holds as long as $n_{1}\gg\kappa^{\natural4}\mu^{\natural}r$,
and the last line follows from (\ref{eq:M-natural-two-to-infty})
and Lemma \ref{lemma:HPCA-G0-spectral}. An immediate consequence
is that 
\begin{align*}
\left\Vert \mathcal{P}_{\mathsf{diag}}\left(\bm{G}^{t}\right)\right\Vert  & \leq\left\Vert \mathcal{P}_{\mathsf{diag}}\left(\bm{G}^{\natural}\right)\right\Vert +\left\Vert \mathcal{P}_{\mathsf{diag}}\left(\bm{G}^{t}-\bm{G}^{\natural}\right)\right\Vert =\left\Vert \bm{M}^{\natural}\right\Vert _{2,\infty}^{2}+\left\Vert \mathcal{P}_{\mathsf{diag}}\left(\bm{G}^{t}-\bm{G}^{\natural}\right)\right\Vert \\
 & \leq\frac{\mu^{\natural}r}{n_{1}}\sigma_{1}^{\natural2}+\left(C_{0}\kappa^{\natural2}\sqrt{\frac{\mu^{\natural}r}{n_{1}}}\right)^{t}\frac{\mu^{\natural}r}{n_{1}}\sigma_{1}^{\natural2}+2C_{0}C^{\mathsf{dd}}\kappa^{\natural2}\sqrt{\frac{\mu^{\natural}r}{n_{1}}}\zeta_{\mathsf{op}}\\
 & \leq\frac{2\mu^{\natural}r}{n_{1}}\sigma_{1}^{\natural2}+2C_{0}C^{\mathsf{dd}}\kappa^{\natural2}\sqrt{\frac{\mu^{\natural}r}{n_{1}}}\zeta_{\mathsf{op}},
\end{align*}
where the second line follows from (\ref{eq:M-natural-two-to-infty})
and (\ref{eq:HeteroPCA_inter5}), and the last relation is valid as
long as $n_{1}\gg\kappa^{\natural4}\mu^{\natural}r$. This together
with the fact $\mathcal{P}_{\mathsf{off}\text{-}\mathsf{diag}}\left(\bm{G}^{t}\right)=\mathcal{P}_{\mathsf{off}\text{-}\mathsf{diag}}\left(\bm{G}^{0}\right)$
(by construction) allows one to obtain 
\begin{equation}
\left\Vert \bm{G}^{t}-\bm{G}^{0}\right\Vert =\left\Vert \mathcal{P}_{\mathsf{diag}}\left(\bm{G}^{t}\right)\right\Vert \leq\frac{2\mu^{\natural}r}{n_{1}}\sigma_{1}^{\natural2}+2C_{0}C^{\mathsf{dd}}\kappa^{\natural2}\sqrt{\frac{\mu^{\natural}r}{n_{1}}}\zeta_{\mathsf{op}}.\label{eq:HeteroPCA_inter6}
\end{equation}
In view of Weyl's inequality, we have 
\begin{align}
\lambda_{r}\left(\bm{G}^{t}\right) & \geq\lambda_{r}\left(\bm{G}^{\natural}\right)-\left\Vert \bm{G}^{t}-\bm{G}^{\natural}\right\Vert \geq\lambda_{r}\left(\bm{G}^{\natural}\right)-\left\Vert \bm{G}^{t}-\bm{G}^{0}\right\Vert -\left\Vert \bm{G}^{\natural}-\bm{G}^{0}\right\Vert \nonumber \\
 & \overset{\text{(i)}}{\geq}\sigma_{r}^{\natural2}-\frac{2\mu^{\natural}r}{n_{1}}\sigma_{1}^{\natural2}-2C_{0}C^{\mathsf{dd}}\kappa^{\natural2}\sqrt{\frac{\mu^{\natural}r}{n_{1}}}\zeta_{\mathsf{op}}-C^{\mathsf{dd}}\zeta_{\mathsf{op}}-\left\Vert \bm{M}^{\natural}\right\Vert _{2,\infty}^{2}\nonumber \\
 & \overset{\text{(ii)}}{\geq}\sigma_{r}^{\natural2}-\frac{3\mu^{\natural}r}{n_{1}}\sigma_{1}^{\natural2}-2C_{0}C^{\mathsf{dd}}\kappa^{\natural2}\sqrt{\frac{\mu^{\natural}r}{n_{1}}}\zeta_{\mathsf{op}}-C^{\mathsf{dd}}\zeta_{\mathsf{op}}\overset{\text{(iii)}}{\geq}\frac{\sigma_{r}^{\natural2}}{2}.\label{eq:HeteroPCA-Gt-lb}
\end{align}
Here, (i) follows from (\ref{eq:HeteroPCA_inter6}) and Lemma \ref{lemma:HPCA-G0-spectral},
(ii) comes from (\ref{eq:M-natural-two-to-infty}), and (iii) is guaranteed
to hold as long as $n_{1}\gg\kappa^{\natural4}\mu^{\natural}r$ and
$\zeta_{\mathsf{op}}\ll\sigma_{r}^{\natural2}$. Then by virtue of
Davis Kahan's $\sin\bm{\Theta}$ Theorem \citep[Theorem 2.2.1]{chen2020spectral},
\begin{align*}
\left\Vert \bm{U}^{t}\bm{R}^{t}-\bm{U}^{0}\right\Vert  & \leq\frac{\left\Vert \bm{G}^{t}-\bm{G}^{0}\right\Vert }{\lambda_{r}\left(\bm{G}^{t}\right)-\lambda_{r+1}\left(\bm{G}^{\natural}\right)}\leq\frac{4\kappa^{\natural2}\mu^{\natural}r}{n_{1}}+4C_{0}C^{\mathsf{dd}}\kappa^{\natural2}\sqrt{\frac{\mu^{\natural}r}{n_{1}}}\frac{\zeta_{\mathsf{op}}}{\sigma_{r}^{\natural2}},
\end{align*}
where the last relation arises from $\lambda_{r+1}(\bm{G}^{\natural})=0$,
(\ref{eq:HeteroPCA_inter6}) and (\ref{eq:HeteroPCA-Gt-lb}). This
indicates that 
\begin{align}
\left\Vert \bm{U}^{t}\right\Vert _{2,\infty} & \leq\left\Vert \bm{U}^{0}\right\Vert _{2,\infty}+\left\Vert \bm{U}^{t}\bm{R}^{t}-\bm{U}^{0}\right\Vert _{2,\infty}\nonumber \\
 & \leq\sqrt{\frac{4\mu^{\natural}r}{n_{1}}}+\frac{4\kappa^{\natural2}\mu^{\natural}r}{n_{1}}+4C_{0}C^{\mathsf{dd}}\kappa^{\natural2}\sqrt{\frac{\mu^{\natural}r}{n_{1}}}\frac{\zeta_{\mathsf{op}}}{\sigma_{r}^{\natural2}}\nonumber \\
 & \leq3\sqrt{\frac{\mu^{\natural}r}{n_{1}}},\label{eq:Ut-2-infty-UB-hpca}
\end{align}
where the penultimate relation follows from (\ref{eq:U0-two-to-infty}),
and the last relation holds provided that $n_{1}\gg\kappa^{\natural4}\mu^{\natural}r$
and $\zeta_{\mathsf{op}}/\sigma_{r}^{\natural2}\ll1/\kappa^{\natural2}$.

To proceed, we recall that the diagonal entries of $\bm{G}^{t+1}$
are set to be the diagonal entries of $\bm{U}^{t}\bm{\Lambda}^{t}\bm{U}^{t\top}$,
thus revealing that 
\begin{align*}
\left\Vert \mathcal{P}_{\mathsf{diag}}\left(\bm{G}^{t+1}-\bm{G}^{\natural}\right)\right\Vert  & =\left\Vert \mathcal{P}_{\mathsf{diag}}\left(\bm{U}^{t}\bm{\Lambda}^{t}\bm{U}^{t\top}-\bm{G}^{\natural}\right)\right\Vert =\left\Vert \mathcal{P}_{\mathsf{diag}}\left(\mathcal{P}_{\bm{U}^{t}}\left(\bm{G}^{t}\right)-\bm{G}^{\natural}\right)\right\Vert \\
 & \leq\underbrace{\left\Vert \mathcal{P}_{\mathsf{diag}}\left[\mathcal{P}_{\bm{U}^{t}}\left(\bm{G}^{t}-\bm{G}^{\natural}\right)\right]\right\Vert }_{\eqqcolon\beta_{1}}+\underbrace{\left\Vert \mathcal{P}_{\mathsf{diag}}\left(\mathcal{P}_{\bm{U}_{\perp}^{t}}\bm{G}^{\natural}\right)\right\Vert }_{\eqqcolon\beta_{2}}.
\end{align*}
Here, $\bm{U}_{\perp}^{t}$ represents the orthogonal complement of
the subspace $\bm{U}^{t}$. Similar to how we bound $\alpha_{1}$
and $\alpha_{2}$ for the base case, we can invoke \citet[Lemma 1]{zhang2018heteroskedastic}
and (\ref{eq:Ut-2-infty-UB-hpca}) to reach 
\begin{align}
\beta_{1} & \leq3\sqrt{\frac{\mu^{\natural}r}{n_{1}}}\left\Vert \bm{G}^{t}-\bm{G}^{\natural}\right\Vert ,\\
\beta_{2} & \leq\sqrt{\frac{\mu^{\natural}r}{n_{1}}}\left\Vert \mathcal{P}_{\bm{U}_{\perp}^{t}}\bm{G}^{\natural}\right\Vert \leq\sigma_{1}^{\natural2}\left\Vert \big(\bm{U}_{\perp}^{t}\big)^{\top}\bm{U}^{\natural}\right\Vert =\sigma_{1}^{\natural2}\left\Vert \sin\bm{\Theta}\left(\bm{U}^{t},\bm{U}^{\natural}\right)\right\Vert ,\label{eq:HeteroPCA_inter7}
\end{align}
where $\bm{\Theta}\left(\bm{U}^{t},\bm{U}^{\natural}\right)$ is a
diagonal matrix whose diagonal entries are the principal angles between
$\bm{U}^{t}$ and $\bm{U}^{\natural}$. Apply the Davis Kahan $\sin\bm{\Theta}$
Theorem \citep[Theorem 2.2.1]{chen2020spectral} to obtain 
\begin{equation}
\left\Vert \sin\bm{\Theta}\left(\bm{U}^{t},\bm{U}^{\natural}\right)\right\Vert \leq\frac{\left\Vert \bm{G}^{t}-\bm{G}^{\natural}\right\Vert }{\lambda_{r}\left(\bm{G}^{t}\right)-\lambda_{r+1}\left(\bm{G}^{\natural}\right)}\leq\frac{2\left\Vert \bm{G}^{t}-\bm{G}^{\natural}\right\Vert }{\sigma_{r}^{\natural2}}.\label{eq:HeteroPCA_inter8}
\end{equation}
Here, the last inequality comes from $\lambda_{r+1}(\bm{G}^{\natural})=0$
and (\ref{eq:HeteroPCA-Gt-lb}). Combine (\ref{eq:HeteroPCA_inter7})
and (\ref{eq:HeteroPCA_inter8}) to reach 
\[
\beta_{2}\leq\sqrt{\frac{\mu^{\natural}r}{n_{1}}}\sigma_{1}^{\natural2}\frac{2\left\Vert \bm{G}^{t}-\bm{G}^{\natural}\right\Vert }{\sigma_{r}^{\natural2}}\leq2\kappa^{\natural2}\sqrt{\frac{\mu^{\natural}r}{n_{1}}}\left\Vert \bm{G}^{t}-\bm{G}^{\natural}\right\Vert .
\]
Taking together the preceding bounds on $\beta_{1}$ and $\beta_{2}$,
we arrive at 
\[
\left\Vert \mathcal{P}_{\mathsf{diag}}\left(\bm{G}^{t+1}-\bm{G}^{\natural}\right)\right\Vert \leq\beta_{1}+\beta_{2}\leq\left(3+2\kappa^{\natural2}\right)\sqrt{\frac{\mu^{\natural}r}{n_{1}}}\left\Vert \bm{G}^{t}-\bm{G}^{\natural}\right\Vert \leq C_{0}\kappa^{\natural2}\sqrt{\frac{\mu^{\natural}r}{n_{1}}}\left\Vert \bm{G}^{t}-\bm{G}^{\natural}\right\Vert ,
\]
provided that the constant $C_{0}$ is large enough.

\paragraph{Invoking the inequality (\ref{eq:HeteroPCA-induction}) to establish
the lemma. }

The above induction steps taken together establish the hypothesis
(\ref{eq:HeteroPCA-induction}), namely, for all $t\geq1$, 
\[
\left\Vert \mathcal{P}_{\mathsf{diag}}\left(\bm{G}^{t}-\bm{G}^{\natural}\right)\right\Vert \leq C_{0}\kappa^{\natural2}\sqrt{\frac{\mu^{\natural}r}{n_{1}}}\left\Vert \bm{G}^{t-1}-\bm{G}^{\natural}\right\Vert .
\]
Follow the same procedure used to derive (\ref{eq:HeteroPCA_inter5}),
we can show that 
\[
\left\Vert \mathcal{P}_{\mathsf{diag}}\left(\bm{G}^{t}-\bm{G}^{\natural}\right)\right\Vert \leq\left(C_{0}\kappa^{\natural2}\sqrt{\frac{\mu^{\natural}r}{n_{1}}}\right)^{t}\frac{\mu^{\natural}r}{n_{1}}\sigma_{1}^{\natural2}+2C_{0}C^{\mathsf{dd}}\kappa^{\natural2}\sqrt{\frac{\mu^{\natural}r}{n_{1}}}\zeta_{\mathsf{op}}.
\]
If the number of iterations $t_{0}$ satisfies 
\[
t_{0}\geq\frac{\log\left(\left(C^{\mathsf{dd}}\right)^{-1}(\kappa^{\natural})^{-2}\sqrt{\mu^{\natural}r/n_{1}}\sigma_{1}^{\natural2}/\zeta_{\mathsf{op}}\right)}{\log\left((\kappa^{\natural})^{-2}\sqrt{n_{1}/\left(\mu^{\natural}r\right)}\right)},
\]
we can guarantee that 
\[
\left\Vert \mathcal{P}_{\mathsf{diag}}\left(\bm{G}^{t_{0}}-\bm{G}^{\natural}\right)\right\Vert \leq3C_{0}C^{\mathsf{dd}}\kappa^{\natural2}\sqrt{\frac{\mu^{\natural}r}{n_{1}}}\zeta_{\mathsf{op}}
\]
as long as $n_{1}\gg\kappa^{\natural4}\mu^{\natural}r$. In addition,
note that when $n_{1}\gg\kappa^{\natural4}\mu^{\natural}r$, one has
\[
\frac{\log\left(\left(C^{\mathsf{dd}}\right)^{-1}(\kappa^{\natural})^{-2}\sqrt{\mu^{\natural}r/n_{1}}\sigma_{1}^{\natural2}/\zeta_{\mathsf{op}}\right)}{\log\left((\kappa^{\natural})^{-2}\sqrt{n_{1}/\left(\mu^{\natural}r\right)}\right)}\leq\log\left(\frac{\sigma_{1}^{\natural2}}{\zeta_{\mathsf{op}}}\right).
\]
Therefore, it suffices to take $t_{0}\geq\log\left(\frac{\sigma_{1}^{\natural2}}{\zeta_{\mathsf{op}}}\right)$
as claimed.

To finish up, we observe that 
\begin{align*}
\left\Vert \bm{G}^{t_{0}}-\bm{G}^{\natural}\right\Vert  & \leq\left\Vert \mathcal{P}_{\mathsf{diag}}\left(\bm{G}^{t_{0}}-\bm{G}^{\natural}\right)\right\Vert +\left\Vert \mathcal{P}_{\mathsf{off}\text{-}\mathsf{diag}}\left(\bm{G}^{t_{0}}-\bm{G}^{\natural}\right)\right\Vert \\
 & \overset{\text{(i)}}{=}\left\Vert \mathcal{P}_{\mathsf{diag}}\left(\bm{G}^{t_{0}}-\bm{G}^{\natural}\right)\right\Vert +\left\Vert \bm{G}^{0}-\mathcal{P}_{\mathsf{off}\text{-}\mathsf{diag}}\left(\bm{G}^{\natural}\right)\right\Vert \\
 & \overset{\text{(ii)}}{\lesssim}\kappa^{\natural2}\sqrt{\frac{\mu^{\natural}r}{n_{1}}}\zeta_{\mathsf{op}}+\zeta_{\mathsf{op}}\overset{\text{(iii)}}{\lesssim}\zeta_{\mathsf{op}}.
\end{align*}
Here, (i) follows from the construction $\bm{G}^{0}=\mathcal{P}_{\mathsf{off}\text{-}\mathsf{diag}}\left(\bm{G}^{0}\right)=\mathcal{P}_{\mathsf{off}\text{-}\mathsf{diag}}\left(\bm{G}^{t_{0}}\right)$;
(ii) follows from Lemma \ref{lemma:HPCA-G0-spectral}; and (iii) holds
provided that $n_{1}\gtrsim\kappa^{\natural4}\mu^{\natural}r$. 

Last but not least, for each $m\in[n_{1}]$ and any given $t>1$, we have
\begin{align*}
	\left|G_{m,m}^{t}-G_{m,m}^{\natural}\right| & =\left|\bm{e}_{m}^{\top}\left(\bm{G}^{t}-\bm{G}^{\natural}\right)\bm{e}_{m}\right|\leq\left|\bm{e}_{m}^{\top}\mathcal{P}_{\bm{U}^{t_{0}}}\left(\bm{G}^{t}-\bm{G}^{\natural}\right)\bm{e}_{m}\right|+\left|\bm{e}_{m}^{\top}\mathcal{P}_{\bm{U}_{\perp}^{t}}\left(\bm{G}^{t}-\bm{G}^{\natural}\right)\bm{e}_{m}\right|\\
	& =\underbrace{\left|\bm{e}_{m}^{\top}\mathcal{P}_{\bm{U}^{t}}\left(\bm{G}^{t}-\bm{G}^{\natural}\right)\bm{e}_{m}\right|}_{\eqqcolon\gamma_{1}}+\underbrace{\left|\bm{e}_{m}^{\top}\mathcal{P}_{\bm{U}_{\perp}^{t}}\bm{G}^{\natural}\bm{e}_{m}\right|}_{\eqqcolon\gamma_{2}}.
\end{align*}
The first term $\gamma_1$ can be upper bounded by
\begin{align*}
	\gamma_{1} & =\left|\bm{e}_{m}^{\top}\bm{U}^{t}\bm{U}^{t\top}\left(\bm{G}^{t}-\bm{G}^{\natural}\right)\bm{e}_{m}\right|=\left|\bm{U}_{m,\cdot}^{t}\bm{U}^{t\top}\left(\bm{G}^{t}-\bm{G}^{\natural}\right)_{\cdot,m}\right|\leq\left\Vert \bm{U}_{m,\cdot}^{t}\right\Vert _{2}\left\Vert \left(\bm{G}^{t}-\bm{G}^{\natural}\right)_{\cdot,m}\right\Vert _{2}\\
	& \leq\left\Vert \bm{U}_{m,\cdot}^{t}\right\Vert _{2}\left\Vert \left(\bm{G}^{t}-\bm{G}^{\natural}\right)\right\Vert _{2}.
\end{align*}
The second term $\gamma_2$ admits the following upper bound
\begin{align*}
	\gamma_{2} & =\left|\bm{e}_{m}^{\top}\bm{U}_{\perp}^{t}\bm{U}_{\perp}^{t\top}\bm{U}^{\natural}\bm{\Sigma}^{\natural2}\bm{U}_{m,\cdot}^{\natural\top}\right|\leq\sigma_{1}^{\natural2}\left\Vert \bm{U}_{\perp}^{t\top}\bm{U}^{\natural}\right\Vert \left\Vert \bm{U}_{m,\cdot}^{\natural}\right\Vert _{2}\leq\sigma_{1}^{\natural2}\left\Vert \sin\bm{\Theta}\left(\bm{U}^{t},\bm{U}^{\natural}\right)\right\Vert \left\Vert \bm{U}_{m,\cdot}^{\natural}\right\Vert _{2}\\
	& \leq2\kappa^{\natural2}\left\Vert \bm{G}^{t}-\bm{G}^{\natural}\right\Vert \left\Vert \bm{U}_{m,\cdot}^{\natural}\right\Vert _{2},
\end{align*}
where the last relation follows from \eqref{eq:HeteroPCA_inter8}. Taking the bounds on $\gamma_{1}$
and $\gamma_{2}$ collectively and let $t=t_{0}$ gives
\[
\left|G_{m,m}^{t_{0}}-G_{m,m}^{\natural}\right|\lesssim\kappa^{\natural2}\left\Vert \bm{G}^{t}-\bm{G}^{\natural}\right\Vert \left(\left\Vert \bm{U}_{m,\cdot}^{\natural}\right\Vert _{2}+\left\Vert \bm{U}_{m,\cdot}^{t}\right\Vert _{2}\right)\lesssim\kappa^{\natural2} \zeta_{\mathsf{op}} \left(\left\Vert \bm{U}_{m,\cdot}^{\natural}\right\Vert _{2}+\left\Vert \bm{U}_{m,\cdot}^{t}\right\Vert _{2}\right).
\]

\subsubsection{Proof of Lemma \ref{lemma:hpca-basic-facts}\label{appendix:proof-hpca-H-R-proximity}}

We first apply Weyl's inequality to demonstrate that 
\[
\lambda_{r}\left(\bm{G}\right)\geq\sigma_{r}^{\natural2}-\left\Vert \bm{G}-\bm{G}^{\natural}\right\Vert \geq\sigma_{r}^{\natural2}-\widetilde{C}\zeta_{\mathsf{op}}\geq\frac{1}{2}\sigma_{r}^{\natural2}
\]
for some constant $\widetilde{C}>0$, where the penultimate step comes
from Lemma \ref{lemma:hpca-estimation}, and the last step holds true
provided that $\zeta_{\mathsf{op}}\ll\sigma_{r}^{\natural2}$. In
view of the Davis-Kahan $\sin\bm{\Theta}$ Theorem \citep[Theorem 2.2.1]{chen2020spectral},
we obtain 
\begin{equation}
\left\Vert \bm{U}\bm{U}^{\top}-\bm{U}^{\natural}\bm{U}^{\natural\top}\right\Vert \leq\frac{\sqrt{2}\left\Vert \bm{G}-\bm{G}^{\natural}\right\Vert }{\lambda_{r}\left(\bm{G}\right)-\lambda_{r+1}\left(\bm{G}^{\natural}\right)}\leq\frac{2\sqrt{2}\left\Vert \bm{G}-\bm{G}^{\natural}\right\Vert }{\sigma_{r}^{\natural2}}\lesssim\frac{\zeta_{\mathsf{op}}}{\sigma_{r}^{\natural2}},\label{eq:hpca-general-davis-kahan}
\end{equation}
which immediately leads to the advertised bound on $\left\Vert \bm{U}\bm{H}-\bm{U}^{\natural}\right\Vert $
as follows 
\[
\left\Vert \bm{U}\bm{H}-\bm{U}^{\natural}\right\Vert =\left\Vert \bm{U}\bm{U}^{\top}\bm{U}^{\natural}-\bm{U}^{\natural}\bm{U}^{\natural\top}\bm{U}^{\natural}\right\Vert \leq\left\Vert \bm{U}\bm{U}^{\top}-\bm{U}^{\natural}\bm{U}^{\natural\top}\right\Vert \lesssim\frac{\zeta_{\mathsf{op}}}{\sigma_{r}^{\natural2}}.
\]

Next, we turn attention to $\big\|\bm{H}-\bm{R}_{\bm{U}}\big\|$.
Given that both $\bm{U}$ and $\bm{U}^{\natural}$ have orthonormal
columns, the SVD of $\bm{H}=\bm{U}^{\top}\bm{U}^{\natural}$ can be
written as 
\[
\bm{H}=\bm{X}(\cos\bm{\Theta})\bm{Y}^{\top},
\]
where $\bm{X},\bm{Y}\in\mathbb{R}^{r\times r}$ are orthonormal matrices
and $\bm{\Theta}$ is a diagonal matrix composed of the principal
angles between $\bm{U}$ and $\bm{U}^{\natural}$ (see \citet[Section 2.1]{chen2020spectral}).
It is well known that one can write $\bm{R}_{\bm{U}}=\mathsf{sgn}(\bm{H})=\bm{X}\bm{Y}^{\top}$,
and therefore, 
\begin{align*}
\left\Vert \bm{H}-\bm{R}_{\bm{U}}\right\Vert  & =\left\Vert \bm{X}\left(\cos\bm{\Theta}-\bm{I}_{r}\right)\bm{Y}^{\top}\right\Vert =\left\Vert \bm{I}_{r}-\cos\bm{\Theta}\right\Vert \\
 & =\left\Vert 2\sin^{2}\left(\bm{\Theta}/2\right)\right\Vert \lesssim\left\Vert \sin\bm{\Theta}\right\Vert ^{2}\asymp\left\Vert \bm{U}\bm{U}^{\top}-\bm{U}^{\natural}\bm{U}^{\natural\top}\right\Vert ^{2}\lesssim\frac{\zeta_{\mathsf{op}}^{2}}{\sigma_{r}^{\natural4}}.
\end{align*}
Here the penultimate relation comes from \citet[Lemma 2.1.2]{chen2020spectral}
and the last relation invokes (\ref{eq:hpca-general-davis-kahan}).
Given that $\bm{R}_{\bm{U}}$ is a square orthonormal matrix, we immediately
have 
\begin{align*}
\sigma_{\max}(\bm{H}_{\bm{U}}) & \leq\sigma_{\max}(\bm{R}_{\bm{U}})+\left\Vert \bm{H}_{\bm{U}}-\bm{R}_{\bm{U}}\right\Vert \leq1+O\Big(\frac{\sigma^{2}n}{\sigma_{r}^{\natural2}}\Big)\leq2,\\
\sigma_{r}(\bm{H}_{\bm{U}}) & \geq\sigma_{\max}(\bm{R}_{\bm{U}})-\left\Vert \bm{H}_{\bm{U}}-\bm{R}_{\bm{U}}\right\Vert \geq1-O\Big(\frac{\sigma^{2}n}{\sigma_{r}^{\natural2}}\Big)\geq1/2,
\end{align*}
provided that $\zeta_{\mathsf{op}}\ll\sigma_{r}^{\natural2}.$ In
addition, we can similarly derive 
\[
\left\Vert \bm{H}^{\top}\bm{H}-\bm{I}_{r}\right\Vert =\left\Vert \bm{Y}\left[\cos^{2}\bm{\Theta}-\bm{I}_{r}\right]\bm{Y}^{\top}\right\Vert =\left\Vert \cos^{2}\bm{\Theta}-\bm{I}_{r}\right\Vert =\left\Vert \sin^{2}\bm{\Theta}\right\Vert =\left\Vert \sin\bm{\Theta}\right\Vert ^{2}\lesssim\frac{\zeta_{\mathsf{op}}^{2}}{\sigma_{r}^{\natural4}}.
\]
An immediate consequence of the proximity between $\bm{H}$ and $\bm{R}_{\bm{U}}$
is that 
\begin{align*}
\left\Vert \bm{U}\bm{H}-\bm{U}^{\natural}\right\Vert  & \leq\left\Vert \bm{U}\bm{R}_{\bm{U}}-\bm{U}^{\natural}\right\Vert +\left\Vert \bm{U}\left(\bm{H}-\bm{R}_{\bm{U}}\right)\right\Vert \leq\left\Vert \bm{U}\bm{R}_{\bm{U}}-\bm{U}^{\natural}\right\Vert +\left\Vert \bm{H}-\bm{R}_{\bm{U}}\right\Vert \\
 & \lesssim\frac{\zeta_{\mathsf{op}}}{\sigma_{r}^{\natural2}}+\frac{\zeta_{\mathsf{op}}^{2}}{\sigma_{r}^{\natural4}}\lesssim\frac{\zeta_{\mathsf{op}}}{\sigma_{r}^{\natural2}},
\end{align*}
provided that $\zeta_{\mathsf{op}}\ll\sigma_{r}^{\natural2}$.

\subsubsection{Proof of Lemma \ref{lemma:hpca-1}\label{sec:proof-lemma-pca-1}}

We begin with the triangle inequality: 
\begin{align*}
\left\Vert \left(\bm{G}-\bm{G}^{\natural}\right)_{m,\cdot}\bm{U}^{\natural}\right\Vert _{2} & \leq\left\Vert \left[\mathcal{P}_{\mathsf{off}\text{-}\mathsf{diag}}\left(\bm{G}-\bm{G}^{\natural}\right)\right]_{m,\cdot}\bm{U}^{\natural}\right\Vert _{2}+\left\Vert \left[\mathcal{P}_{\mathsf{diag}}\left(\bm{G}-\bm{G}^{\natural}\right)\right]_{m,\cdot}\bm{U}^{\natural}\right\Vert _{2}\\
 & =\underbrace{\left\Vert \left[\mathcal{P}_{\mathsf{off}\text{-}\mathsf{diag}}\left(\bm{G}^{0}-\bm{G}^{\natural}\right)\right]_{m,\cdot}\bm{U}^{\natural}\right\Vert _{2}}_{\eqqcolon\,\alpha_{1}}+\underbrace{\left\Vert \left[\mathcal{P}_{\mathsf{diag}}\left(\bm{G}-\bm{G}^{\natural}\right)\right]_{m,\cdot}\bm{U}^{\natural}\right\Vert _{2}}_{\eqqcolon\,\alpha_{2}},
\end{align*}
where the last inequality makes use of our construction $\mathcal{P}_{\mathsf{off}\text{-}\mathsf{diag}}\left(\bm{G}\right)=\mathcal{P}_{\mathsf{off}\text{-}\mathsf{diag}}\left(\bm{G}^{0}\right)$.
Given that the properties of $\bm{G}^{0}$ (which is the diagonal-deleted
version of the sample gram matrix) have been studied in \citet{cai2019subspace},
we can readily borrow \citet[Lemma 2]{cai2019subspace} to bound 
\begin{align*}
\alpha_{1} & \lesssim\left(\sigma\sigma_{m}\sqrt{n_{1}n_{2}\log n}+\sigma\sqrt{n_{1}\log n}\left\Vert \bm{M}_{m,\cdot}^{\natural}\right\Vert _{2}+B\log n\left\Vert \bm{M}_{m,\cdot}^{\natural}\right\Vert _{\infty}\right)\left\Vert \bm{U}^{\natural}\right\Vert _{2,\infty}+\sigma_{m}\sqrt{n_{2}\log n}\left\Vert \bm{M}^{\natural\top}\right\Vert _{2,\infty}
\end{align*}
with probability exceeding $1-O(n^{-10})$. The careful reader might
remark that the above bound is slightly different from \citet[Lemma 2]{cai2019subspace}
in the sense that the bound therein contains an additional term $\Vert\bm{M}^{\natural}\Vert_{2,\infty}^{2}\sqrt{\mu^{\natural}r/n_{1}}$;
note, however, that this extra term is caused by the effect of the
diagonal part $\mathcal{P}_{\mathsf{diag}}(\bm{G}^{\natural})$ in
their analysis, which has been removed in the above term $\alpha_{1}$.
The interested reader is referred to \citet[Appendix B.3]{cai2019subspace}
for details. In view of the bounds
\begin{equation}
\left\Vert \bm{M}_{m,\cdot}^{\natural}\right\Vert _{2}=\left\Vert \bm{U}_{m,\cdot}^{\natural}\bm{\Sigma}^{\natural}\bm{V}^{\natural\top}\right\Vert _{2}\leq\left\Vert \bm{U}_{m,\cdot}^{\natural}\right\Vert _{2}\left\Vert \bm{\Sigma}^{\natural}\right\Vert \left\Vert \bm{V}^{\natural}\right\Vert \leq\sigma_{1}^{\natural}\left\Vert \bm{U}_{m,\cdot}^{\natural}\right\Vert _{2},\label{eq:M-natural-m-2}
\end{equation}
\begin{equation}
\left\Vert \bm{M}_{m,\cdot}^{\natural}\right\Vert _{\infty}=\left\Vert \bm{U}_{m,\cdot}^{\natural}\bm{\Sigma}^{\natural}\bm{V}^{\natural\top}\right\Vert _{\infty}\leq\left\Vert \bm{U}_{m,\cdot}^{\natural}\right\Vert _{2}\left\Vert \bm{\Sigma}^{\natural}\right\Vert \left\Vert \bm{V}^{\natural}\right\Vert _{2,\infty}\leq\sigma_{1}^{\natural}\left\Vert \bm{U}_{m,\cdot}^{\natural}\right\Vert _{2}\sqrt{\frac{\mu^{\natural}r}{n_{2}}},\label{eq:M-natural-m-infty}
\end{equation}
and
\begin{equation}
\left\Vert \bm{M}^{\natural\top}\right\Vert _{2,\infty}=\left\Vert \bm{V}^{\natural}\bm{\Sigma}^{\natural}\bm{U}^{\natural\top}\right\Vert _{2,\infty}\leq\left\Vert \bm{V}^{\natural}\right\Vert _{2,\infty}\left\Vert \bm{\Sigma}^{\natural}\right\Vert \left\Vert \bm{U}^{\natural}\right\Vert \leq\sigma_{1}^{\natural}\sqrt{\frac{\mu^{\natural}r}{n_{2}}},\label{eq:M-natural-T-2-infty-subspace}
\end{equation}
we further have
\begin{align*}
\alpha_{1} & \lesssim\left(\sigma\sigma_{m}\sqrt{n_{1}n_{2}\log n}+\sigma\sigma_{1}^{\natural}\sqrt{n_{1}\log n}\left\Vert \bm{U}_{m,\cdot}^{\natural}\right\Vert _{2}+\sigma\sigma_{1}^{\natural}\sqrt[4]{n_{1}\mu^{\natural}r\log n}\left\Vert \bm{U}_{m,\cdot}^{\natural}\right\Vert _{2}\right)\left\Vert \bm{U}^{\natural}\right\Vert _{2,\infty}\\
 & \quad+\sigma_{m}\sigma_{1}^{\natural}\sqrt{\mu^{\natural}r\log n}\\
 & \lesssim\zeta_{\mathsf{op},m}\sqrt{\frac{\mu^{\natural}r}{n_{1}}}+\sigma\sigma_{1}^{\natural}\sqrt{\mu^{\natural}r\log n}\left\Vert \bm{U}_{m,\cdot}^{\natural}\right\Vert _{2},
\end{align*}
provided that $n_{1}\gtrsim\mu^{\natural}r$. In addition, it is seen
that 
\begin{align*}
\alpha_{2} & =\left|G_{m,m}-G_{m,m}^{\natural}\right|\left\Vert \bm{U}_{m,\cdot}^{\natural}\right\Vert _{2}\lesssim\kappa^{\natural2}\zeta_{\mathsf{op}}\left(\left\Vert \bm{U}_{m,\cdot}^{\natural}\right\Vert _{2}+\left\Vert \bm{U}_{m,\cdot}\right\Vert _{2}\right)\left\Vert \bm{U}_{m,\cdot}^{\natural}\right\Vert _{2}\\
 & \lesssim\kappa^{\natural2}\zeta_{\mathsf{op}}\left\Vert \bm{U}_{m,\cdot}^{\natural}\right\Vert _{2}^{2}+\kappa^{\natural2}\zeta_{\mathsf{op}}\left\Vert \bm{U}_{m,\cdot}\bm{H}-\bm{U}_{m,\cdot}^{\natural}\right\Vert _{2}\left\Vert \bm{U}_{m,\cdot}^{\natural}\right\Vert _{2},
\end{align*}
where the penultimate relation invokes Lemma \ref{lemma:hpca-estimation}.
We can thus conclude that 
\begin{align*}
\left\Vert \left(\bm{G}-\bm{G}^{\natural}\right)_{m,\cdot}\bm{U}^{\natural}\right\Vert _{2} & \leq\alpha_{1}+\alpha_{2}\\
 & \lesssim\zeta_{\mathsf{op},m}\sqrt{\frac{\mu^{\natural}r}{n_{1}}}+\sigma\sigma_{1}^{\natural}\sqrt{\mu^{\natural}r\log n}\left\Vert \bm{U}_{m,\cdot}^{\natural}\right\Vert _{2}+\kappa^{\natural2}\zeta_{\mathsf{op}}\left\Vert \bm{U}_{m,\cdot}^{\natural}\right\Vert _{2}^{2}\\
 & \quad+\kappa^{\natural2}\zeta_{\mathsf{op}}\left\Vert \bm{U}_{m,\cdot}\bm{H}-\bm{U}_{m,\cdot}^{\natural}\right\Vert _{2}\left\Vert \bm{U}_{m,\cdot}^{\natural}\right\Vert _{2}.
\end{align*}

\subsubsection{Proof of Lemma \ref{lemma:hpca-approx-1}\label{sec:proof-lemma-HeteroPCA-approx-1}}

Recall that for each $1\leq m\leq n_{1}$, we employ the notation
$\mathcal{P}_{-m,\cdot}\left(\bm{M}\right)$ to represent an $n_{1}\times n_{2}$
matrix such that 
\[
\big[\mathcal{P}_{-m,\cdot}\left(\bm{M}\right)\big]_{i,\cdot}=\begin{cases}
\bm{0}, & \text{if }i=m,\\
\bm{M}_{i,\cdot},\quad & \text{if }i\neq m.
\end{cases}
\]
In other words, it is obtained by zeroing out the $m$-th row of $\bm{M}$.
Simple algebra then reveals that the $m$-th row of $\bm{G}$ can
be decomposed into 
\begin{align*}
\bm{G}_{m,\cdot} & =\bm{M}_{m,\cdot}\left[\mathcal{P}_{-m,\cdot}\left(\bm{M}\right)\right]^{\top}+G_{m,m}\bm{e}_{m}^{\top}\\
 & =\bm{M}_{m,\cdot}^{\natural}\left[\mathcal{P}_{-m,\cdot}\left(\bm{M}^{\natural}\right)\right]^{\top}+\bm{M}_{m,\cdot}^{\natural}\left[\mathcal{P}_{-m,\cdot}\left(\bm{E}\right)\right]^{\top}+\bm{E}_{m,\cdot}\left[\mathcal{P}_{-m,\cdot}\left(\bm{M}\right)\right]^{\top}+G_{m,m}\bm{e}_{m}^{\top}\\
 & =\bm{M}_{m,\cdot}^{\natural}\bm{M}^{\natural\top}+\bm{M}_{m,\cdot}^{\natural}\left[\mathcal{P}_{-m,\cdot}\left(\bm{E}\right)\right]^{\top}+\bm{E}_{m,\cdot}\left[\mathcal{P}_{-m,\cdot}\left(\bm{M}\right)\right]^{\top}+\left(G_{m,m}-G_{m,m}^{\natural}\right)\bm{e}_{m}^{\top},
\end{align*}
where we have used $\bm{G}^{\natural}=\bm{M}^{\natural}\bm{M}^{\natural\top}$.
Apply the triangle inequality once again to yield 
\begin{align*}
\left\Vert \bm{G}_{m,\cdot}\left(\bm{U}\bm{H}-\bm{U}^{\natural}\right)\right\Vert _{2} & \leq\underbrace{\left\Vert \bm{M}_{m,\cdot}^{\natural}\bm{M}^{\natural\top}\left(\bm{U}\bm{H}-\bm{U}^{\natural}\right)\right\Vert _{2}}_{\eqqcolon\alpha_{1}}+\underbrace{\left\Vert \bm{M}_{m,\cdot}^{\natural}\left[\mathcal{P}_{-m,\cdot}\left(\bm{E}\right)\right]^{\top}\left(\bm{U}\bm{H}-\bm{U}^{\natural}\right)\right\Vert _{2}}_{\eqqcolon\alpha_{2}}\\
 & \quad+\underbrace{\left\Vert \bm{E}_{m,\cdot}\left[\mathcal{P}_{-m,\cdot}\left(\bm{M}\right)\right]^{\top}\left(\bm{U}\bm{H}-\bm{U}^{\natural}\right)\right\Vert _{2}}_{\eqqcolon\alpha_{3}}+\underbrace{\left\Vert \left(G_{m,m}-G_{m,m}^{\natural}\right)\bm{e}_{m}^{\top}\left(\bm{U}\bm{H}-\bm{U}^{\natural}\right)\right\Vert _{2}}_{\eqqcolon\alpha_{4}}
\end{align*}
for each $1\leq m\leq n_{1}$. In what follows, we shall bound the
terms $\alpha_{1},\cdots,\alpha_{4}$ separately. 
\begin{itemize}
\item Let us begin with $\alpha_{1}$. Write 
\begin{equation}
\left\Vert \bm{U}^{\natural\top}\left(\bm{U}\bm{H}-\bm{U}^{\natural}\right)\right\Vert =\left\Vert \bm{U}^{\natural\top}\bm{U}\bm{U}^{\top}\bm{U}^{\natural}-\bm{I}_{r}\right\Vert =\left\Vert \bm{H}^{\top}\bm{H}-\bm{I}_{r}\right\Vert \lesssim\frac{\zeta_{\mathsf{op}}^{2}}{\sigma_{r}^{\natural4}},\label{eq:pca-approx-1-1}
\end{equation}
where the last relation follows from Lemma \ref{lemma:hpca-basic-facts}.
This allows us to upper bound $\alpha_{1}$ as follows: 
\begin{align*}
\alpha_{1} & =\left\Vert \bm{e}_{m}^{\top}\bm{M}^{\natural}\bm{M}^{\natural\top}\left(\bm{U}\bm{H}-\bm{U}^{\natural}\right)\right\Vert _{2}\leq\left\Vert \bm{e}_{m}^{\top}\bm{U}^{\natural}\bm{\Sigma}^{\natural2}\bm{U}^{\natural\top}\left(\bm{U}\bm{H}-\bm{U}^{\natural}\right)\right\Vert _{2}\\
 & \leq\left\Vert \bm{U}_{m,\cdot}^{\natural}\right\Vert _{2}\left\Vert \bm{\Sigma}^{\natural}\right\Vert ^{2}\left\Vert \bm{U}^{\natural\top}\left(\bm{U}\bm{H}-\bm{U}^{\natural}\right)\right\Vert \\
 & \lesssim\sigma_{1}^{\natural2}\frac{\zeta_{\mathsf{op}}^{2}}{\sigma_{r}^{\natural4}}\left\Vert \bm{U}_{m,\cdot}^{\natural}\right\Vert _{2}\lesssim\kappa^{\natural2}\frac{\zeta_{\mathsf{op}}^{2}}{\sigma_{r}^{\natural2}}\left\Vert \bm{U}_{m,\cdot}^{\natural}\right\Vert _{2}.
\end{align*}
\item Regarding $\alpha_{2}$, it is readily seen that 
\begin{align*}
\alpha_{2} & \leq\left\Vert \bm{M}_{m,\cdot}^{\natural}\left[\mathcal{P}_{-m,\cdot}\left(\bm{E}\right)\right]^{\top}\right\Vert _{2}\left\Vert \bm{U}\bm{H}-\bm{U}^{\natural}\right\Vert .
\end{align*}
Regarding the first term $\Vert\bm{M}_{m,\cdot}^{\natural}[\mathcal{P}_{-m,\cdot}(\bm{E})]^{\top}\Vert_{2}$,
we notice that 
\begin{align*}
\left\Vert \bm{M}_{m,\cdot}^{\natural}\left[\mathcal{P}_{-m,\cdot}\left(\bm{E}\right)\right]^{\top}\right\Vert _{2}^{2} & \leq\left\Vert \bm{M}_{m,\cdot}^{\natural}\bm{E}^{\top}\right\Vert _{2}^{2}=\sum_{i=1}^{n_{1}}\left(\sum_{j=1}^{n_{2}}M_{m,j}^{\natural}E_{i,j}\right)^{2}\\
 & \lesssim\left\Vert \bm{M}_{m,\cdot}^{\natural}\right\Vert _{2}^{2}\left(\sigma^{2}n_{1}+\sigma^{2}\log^{2}n\right)+\left\Vert \bm{M}_{m,\cdot}^{\natural}\right\Vert _{\infty}^{2}B^{2}\log^{3}n\\
 & \lesssim\left\Vert \bm{U}_{m,\cdot}^{\natural}\right\Vert _{2}^{2}\sigma_{1}^{\natural2}\left(\sigma^{2}n_{1}+\sigma^{2}\log^{2}n\right)+\frac{\mu^{\natural}r}{n_{2}}\left\Vert \bm{U}_{m,\cdot}^{\natural}\right\Vert _{2}^{2}\sigma_{1}^{\natural2}B^{2}\log^{3}n\\
 & \lesssim\sigma^{2}\sigma_{1}^{\natural2}n_{1}\left\Vert \bm{U}_{m,\cdot}^{\natural}\right\Vert _{2}^{2},
\end{align*}
where the second line follows from \citet[Lemma 14]{cai2019subspace},
the third line comes from (\ref{eq:M-natural-m-2}) and (\ref{eq:M-natural-m-infty}),
and the last line holds provided that $n_{1}\gtrsim\mu^{\natural2}r\log^{2}n$,
$n_{2}\gtrsim r\log^{2}n$ and $B\lesssim\sigma\sqrt[4]{n_{1}n_{2}}/\sqrt{\log n}$.
Therefore, we can demonstrate that 
\begin{align*}
\alpha_{2} & \lesssim\left\Vert \bm{M}_{m,\cdot}^{\natural}\left[\mathcal{P}_{-m,\cdot}\left(\bm{E}\right)\right]^{\top}\right\Vert _{2}\left\Vert \bm{U}\bm{H}-\bm{U}^{\natural}\right\Vert \overset{\text{(i)}}{\lesssim}\sigma\sigma_{1}^{\natural}\sqrt{n_{1}}\left\Vert \bm{U}_{m,\cdot}^{\natural}\right\Vert _{2}\frac{\zeta_{\mathsf{op}}}{\sigma_{r}^{\natural2}}\overset{\text{(ii)}}{\lesssim}\left\Vert \bm{U}_{m,\cdot}^{\natural}\right\Vert _{2}\frac{\zeta_{\mathsf{op}}^{2}}{\sigma_{r}^{\natural2}}.
\end{align*}
Here, (i) utilizes Lemma \ref{lemma:hpca-basic-facts} and (ii) comes
from the definition of $\zeta_{\mathsf{op}}\geq\sigma\sigma_{1}^{\natural}\sqrt{n_{1}\log n}$. 
\item When it comes to $\alpha_{3}$, our starting point is 
\begin{align*}
\alpha_{3} & \leq\underbrace{\left\Vert \bm{E}_{m,\cdot}\left[\mathcal{P}_{-m,\cdot}\left(\bm{M}\right)\right]^{\top}\left(\bm{U}^{(m)}\bm{H}^{(m)}-\bm{U}^{\natural}\right)\right\Vert _{2}}_{\eqqcolon\,\beta_{1}}+\underbrace{\left\Vert \bm{E}_{m,\cdot}\left[\mathcal{P}_{-m,\cdot}\left(\bm{M}\right)\right]^{\top}\left(\bm{U}\bm{H}-\bm{U}^{(m)}\bm{H}^{(m)}\right)\right\Vert _{2}}_{\eqqcolon\,\beta_{2}}.
\end{align*}

\begin{itemize}
\item The first term $\beta_{1}$ can be controlled using Lemma \ref{lemma:hpca-useful-1}
as follows 
\begin{align}
\beta_{1} & \lesssim\zeta_{\mathsf{op},m}\left(\left\Vert \bm{U}^{(m)}\bm{H}^{(m)}-\bm{U}^{\natural}\right\Vert _{2,\infty}+\frac{\zeta_{\mathsf{op}}}{\sigma_{r}^{\natural2}}\left\Vert \bm{U}^{(m)}\bm{H}^{(m)}\right\Vert _{2,\infty}\right)\nonumber \\
 & \quad+\frac{\zeta_{\mathsf{op}}}{\sigma_{r}^{\natural2}}\sigma_{m}\sigma_{1}^{\natural}\sqrt{n_{1}\log n}\sqrt{\frac{\mu^{\natural}r}{n_{1}}}.\label{eq:beta-1-initial-bound}
\end{align}
Additionally, it follows from Lemma \ref{lemma:hpca-loo-perturbation}
that 
\begin{align*}
\left\Vert \bm{U}^{(m)}\bm{H}^{(m)}-\bm{U}\bm{H}\right\Vert  & \leq\left\Vert \bm{U}^{(m)}\bm{U}^{(m)\top}-\bm{U}\bm{U}^{\top}\right\Vert \leq\kappa^{\natural2}\frac{\zeta_{\mathsf{op}}}{\sigma_{r}^{\natural2}}\left(\left\Vert \bm{U}\bm{H}\right\Vert _{2,\infty}+\sqrt{\frac{\mu^{\natural}r}{n_{1}}}\right)\\
 & \leq\kappa^{\natural2}\frac{\zeta_{\mathsf{op}}}{\sigma_{r}^{\natural2}}\left(\left\Vert \bm{U}\bm{H}-\bm{U}^{\natural}\right\Vert _{2,\infty}+\left\Vert \bm{U}^{\natural}\right\Vert _{2,\infty}+\sqrt{\frac{\mu^{\natural}r}{n_{1}}}\right)\\
 & \asymp\kappa^{\natural2}\frac{\zeta_{\mathsf{op}}}{\sigma_{r}^{\natural2}}\left(\left\Vert \bm{U}\bm{H}-\bm{U}^{\natural}\right\Vert _{2,\infty}+\sqrt{\frac{\mu^{\natural}r}{n_{1}}}\right),
\end{align*}
which together with the triangle inequality gives 
\begin{align}
\left\Vert \bm{U}^{(m)}\bm{H}^{(m)}-\bm{U}^{\natural}\right\Vert _{2,\infty} & \leq\left\Vert \bm{U}^{(m)}\bm{H}^{(m)}-\bm{U}\bm{H}\right\Vert +\left\Vert \bm{U}\bm{H}-\bm{U}^{\natural}\right\Vert _{2,\infty}\nonumber \\
 & \lesssim\kappa^{\natural2}\frac{\zeta_{\mathsf{op}}}{\sigma_{r}^{\natural2}}\left(\left\Vert \bm{U}\bm{H}-\bm{U}^{\natural}\right\Vert _{2,\infty}+\sqrt{\frac{\mu^{\natural}r}{n_{1}}}\right)+\left\Vert \bm{U}\bm{H}-\bm{U}^{\natural}\right\Vert _{2,\infty}\nonumber \\
 & \asymp\kappa^{\natural2}\frac{\zeta_{\mathsf{op}}}{\sigma_{r}^{\natural2}}\sqrt{\frac{\mu^{\natural}r}{n_{1}}}+\left\Vert \bm{U}\bm{H}-\bm{U}^{\natural}\right\Vert _{2,\infty}\label{eq:Um-Hm-Ustar-2-infty}
\end{align}
as long as $\zeta_{\mathsf{op}}\ll\sigma_{r}^{\natural2}/\kappa^{\natural2}$.
An immediate consequence is that 
\begin{align}
\left\Vert \bm{U}^{(m)}\bm{H}^{(m)}\right\Vert _{2,\infty} & \leq\left\Vert \bm{U}^{\natural}\right\Vert _{2,\infty}+\left\Vert \bm{U}^{(m)}\bm{H}^{(m)}-\bm{U}^{\natural}\right\Vert _{2,\infty}\nonumber \\
 & \lesssim\sqrt{\frac{\mu^{\natural}r}{n_{1}}}+\kappa^{\natural2}\frac{\zeta_{\mathsf{op}}}{\sigma_{r}^{\natural2}}\sqrt{\frac{\mu^{\natural}r}{n_{1}}}+\left\Vert \bm{U}\bm{H}-\bm{U}^{\natural}\right\Vert _{2,\infty}\nonumber \\
 & \asymp\sqrt{\frac{\mu^{\natural}r}{n_{1}}}+\left\Vert \bm{U}\bm{H}-\bm{U}^{\natural}\right\Vert _{2,\infty},\label{eq:Um-2-infty}
\end{align}
with the proviso that $\zeta_{\mathsf{op}}\ll\sigma_{r}^{\natural2}/\kappa^{\natural2}$.
Therefore we can invoke (\ref{eq:Um-Hm-Ustar-2-infty}) and (\ref{eq:Um-2-infty})
to refine (\ref{eq:beta-1-initial-bound}) as 
\[
\beta_{1}\lesssim\zeta_{\mathsf{op},m}\left(\kappa^{\natural2}\frac{\zeta_{\mathsf{op}}}{\sigma_{r}^{\natural2}}\sqrt{\frac{\mu^{\natural}r}{n_{1}}}+\left\Vert \bm{U}\bm{H}-\bm{U}^{\natural}\right\Vert _{2,\infty}\right)
\]
provided that $\zeta_{\mathsf{op}}\ll\sigma_{r}^{\natural2}$. 
\item With regards to the second term $\beta_{2}$, we see that 
\begin{align*}
\beta_{2} & =\left\Vert \bm{E}_{m,\cdot}\left[\mathcal{P}_{-m,\cdot}\left(\bm{M}\right)\right]^{\top}\left(\bm{U}\bm{U}^{\top}-\bm{U}^{(m)}\bm{U}^{(m)\top}\right)\bm{U}^{\natural}\right\Vert _{2}\\
 & \leq\left\Vert \bm{E}_{m,\cdot}\left[\mathcal{P}_{-m,\cdot}\left(\bm{M}\right)\right]^{\top}\right\Vert _{2}\left\Vert \bm{U}\bm{U}^{\top}-\bm{U}^{(m)}\bm{U}^{(m)\top}\right\Vert .
\end{align*}
It has already been proved in \citet[Appendix C.2]{cai2019subspace}
that with probability exceeding $1-O(n^{-10})$, 
\begin{align*}
\left\Vert \bm{E}_{m,\cdot}\left[\mathcal{P}_{-m,\cdot}\left(\bm{M}\right)\right]^{\top}\right\Vert _{2} & \lesssim\sigma\sigma_{m}\sqrt{n_{1}n_{2}\log n}+\sigma_{m}\sqrt{n_{2}\log n}\left\Vert \bm{M}^{\natural\top}\right\Vert _{2,\infty}\\
 & \lesssim\sigma\sigma_{m}\sqrt{n_{1}n_{2}\log n}+\sigma_{m}\sqrt{\mu^{\natural}r\log n}\sigma_{1}^{\natural}
\end{align*}
where the second relation follows from (\ref{eq:M-natural-T-2-infty-subspace}).
This together with Lemma \ref{lemma:hpca-loo-perturbation} provides
an upper bound on $\beta_{2}$ as follows: 
\begin{align*}
\beta_{2} & \leq\left\Vert \bm{E}_{m,\cdot}\left[\mathcal{P}_{-m,\cdot}\left(\bm{M}\right)\right]^{\top}\right\Vert _{2}\left\Vert \bm{U}\bm{H}-\bm{U}^{(m)}\bm{H}^{(m)}\right\Vert \\
 & \lesssim\left(\sigma\sigma_{m}\sqrt{n_{1}n_{2}\log n}+\sigma_{m}\sqrt{\mu^{\natural}r\log n}\sigma_{1}^{\natural}\right)\kappa^{\natural2}\frac{\zeta_{\mathsf{op}}}{\sigma_{r}^{\natural2}}\left(\left\Vert \bm{U}\bm{H}\right\Vert _{2,\infty}+\sqrt{\frac{\mu^{\natural}r}{n_{1}}}\right)\\
 & \lesssim\left(\sigma\sigma_{m}\sqrt{n_{1}n_{2}\log n}+\sigma_{m}\sqrt{\mu^{\natural}r\log n}\sigma_{1}^{\natural}\right)\kappa^{\natural2}\frac{\zeta_{\mathsf{op}}}{\sigma_{r}^{\natural2}}\left(\left\Vert \bm{U}\bm{H}-\bm{U}^{\natural}\right\Vert _{2,\infty}+\sqrt{\frac{\mu^{\natural}r}{n_{1}}}\right).
\end{align*}
\item Combine the preceding bounds on $\beta_{1}$ and $\beta_{2}$ to arrive
at 
\[
\alpha_{3}\leq\beta_{1}+\beta_{2}\lesssim\zeta_{\mathsf{op},m}\left(\kappa^{\natural2}\frac{\zeta_{\mathsf{op}}}{\sigma_{r}^{\natural2}}\sqrt{\frac{\mu^{\natural}r}{n_{1}}}+\left\Vert \bm{U}\bm{H}-\bm{U}^{\natural}\right\Vert _{2,\infty}\right),
\]
provided that $n_{1}\gtrsim\mu^{\natural}r$ and $\zeta_{\mathsf{op}}/\sigma_{r}^{\natural2}\lesssim1/\kappa^{\natural2}$. 
\end{itemize}
\item For $\alpha_{4}$, Lemma \ref{lemma:hpca-estimation} tells us that
\begin{align*}
\alpha_{4} & \leq\left|G_{m,m}-G_{m,m}^{\natural}\right|\left\Vert \bm{U}_{m,\cdot}\bm{H}-\bm{U}_{m,\cdot}^{\natural}\right\Vert _{2}\lesssim\kappa^{\natural2}\zeta_{\mathsf{op}}\left(\left\Vert \bm{U}_{m,\cdot}^{\natural}\right\Vert _{2}+\left\Vert \bm{U}_{m,\cdot}\right\Vert _{2}\right)\left\Vert \bm{U}_{m,\cdot}\bm{H}-\bm{U}_{m,\cdot}^{\natural}\right\Vert _{2}\\
 & \lesssim\kappa^{\natural2}\zeta_{\mathsf{op}}\left\Vert \bm{U}_{m,\cdot}^{\natural}\right\Vert _{2}\left\Vert \bm{U}_{m,\cdot}\bm{H}-\bm{U}_{m,\cdot}^{\natural}\right\Vert _{2}+\kappa^{\natural2}\zeta_{\mathsf{op}}\left\Vert \bm{U}_{m,\cdot}\bm{H}-\bm{U}_{m,\cdot}^{\natural}\right\Vert _{2}^{2}
\end{align*}
\end{itemize}
Thus far, we have developed upper bounds on $\alpha_{1},\cdots,\alpha_{4}$,
which taken collectively lead to 
\begin{align*}
\left\Vert \bm{G}_{m,\cdot}\left(\bm{U}\bm{H}-\bm{U}^{\natural}\right)\right\Vert _{2} & \leq\alpha_{1}+\alpha_{2}+\alpha_{3}+\alpha_{4}\\
 & \lesssim\zeta_{\mathsf{op},m}\left(\kappa^{\natural2}\frac{\zeta_{\mathsf{op}}}{\sigma_{r}^{\natural2}}\sqrt{\frac{\mu^{\natural}r}{n_{1}}}+\left\Vert \bm{U}\bm{H}-\bm{U}^{\natural}\right\Vert _{2,\infty}\right)+\kappa^{\natural2}\frac{\zeta_{\mathsf{op}}^{2}}{\sigma_{r}^{\natural2}}\left\Vert \bm{U}_{m,\cdot}^{\natural}\right\Vert _{2}\\
 & \quad+\kappa^{\natural2}\zeta_{\mathsf{op}}\left\Vert \bm{U}_{m,\cdot}^{\natural}\right\Vert _{2}\left\Vert \bm{U}_{m,\cdot}\bm{H}-\bm{U}_{m,\cdot}^{\natural}\right\Vert _{2}+\kappa^{\natural2}\zeta_{\mathsf{op}}\left\Vert \bm{U}_{m,\cdot}\bm{H}-\bm{U}_{m,\cdot}^{\natural}\right\Vert _{2}^{2}.
\end{align*}

\subsubsection{Proof of Lemma \ref{lemma:hpca-approx-2}\label{sec:proof-lemma-HeteroPCA-approx-2}}

To control the target quantity, we make note of the following decomposition
\begin{align*}
\left\Vert \bm{R}^{\top}\bm{\Sigma}^{2}\bm{R}-\bm{\Sigma}^{\natural2}\right\Vert  & \leq\underbrace{\left\Vert \bm{R}^{\top}\bm{\Sigma}^{2}\bm{R}-\bm{H}^{\top}\bm{\Sigma}^{2}\bm{H}\right\Vert }_{\eqqcolon\,\alpha_{1}}+\underbrace{\left\Vert \bm{H}^{\top}\bm{\Sigma}^{2}\bm{H}-\bm{U}^{\natural\top}\bm{G}\bm{U}^{\natural}\right\Vert }_{\eqqcolon\,\alpha_{2}}+\underbrace{\left\Vert \bm{U}^{\natural\top}\bm{G}\bm{U}^{\natural}-\bm{\Sigma}^{\natural2}\right\Vert }_{\eqqcolon\,\alpha_{3}}.
\end{align*}
In the sequel, we shall upper bound each of these terms separately.

\paragraph{Step 1: bounding $\alpha_{1}$.}

Lemma \ref{lemma:hpca-basic-facts} tells us that 
\begin{align*}
\alpha_{1} & \leq\left\Vert \left(\bm{H}-\bm{R}\right)^{\top}\bm{\Sigma}^{2}\bm{H}\right\Vert +\left\Vert \bm{R}^{\top}\bm{\Sigma}^{2}\left(\bm{H}-\bm{R}\right)\right\Vert \\
 & \leq\left\Vert \left(\bm{H}-\bm{R}\right)\right\Vert \left\Vert \bm{\Lambda}\right\Vert \left(\left\Vert \bm{H}\right\Vert +\left\Vert \bm{R}\right\Vert \right)\\
 & \lesssim\frac{\zeta_{\mathsf{op}}^{2}}{\sigma_{r}^{\natural4}}\sigma_{1}^{\natural2}\asymp\kappa^{\natural2}\frac{\zeta_{\mathsf{op}}^{2}}{\sigma_{r}^{\natural2}}.
\end{align*}
Here, the penultimate inequality results from an application of Weyl's
inequality: 
\[
\left\Vert \bm{\Lambda}\right\Vert \leq\left\Vert \bm{\Lambda}^{\natural}\right\Vert +\left\Vert \bm{G}-\bm{G}^{\natural}\right\Vert \lesssim\sigma_{1}^{\natural2}+\zeta_{\mathsf{op}}\asymp\sigma_{1}^{\natural2},
\]
where we have used Lemma \ref{lemma:hpca-estimation} and the assumption
$\zeta_{\mathsf{op}}\lesssim\sigma_{1}^{\natural2}$.

\paragraph{Step 2: bounding $\alpha_{2}$.}

Regarding $\alpha_{2}$, it is easily seen that 
\begin{align*}
\bm{H}^{\top}\bm{\Sigma}^{2}\bm{H}-\bm{U}^{\natural\top}\bm{G}\bm{U}^{\natural} & =\bm{U}^{\natural\top}\bm{U}\bm{\Sigma}^{2}\bm{U}^{\top}\bm{U}^{\natural}-\bm{U}^{\natural\top}\bm{G}\bm{U}^{\natural}=\bm{U}^{\natural\top}\left(\bm{U}\bm{\Sigma}^{2}\bm{U}^{\top}-\bm{G}\right)\bm{U}^{\natural}\\
 & =-\bm{U}^{\natural\top}\bm{U}_{\perp}\bm{\Lambda}_{\perp}\bm{U}_{\perp}^{\top}\bm{U}^{\natural},
\end{align*}
where we denote the full SVD of $\bm{G}$ as 
\[
\bm{G}=\left[\begin{array}{cc}
\bm{U} & \bm{U}_{\perp}\end{array}\right]\left[\begin{array}{cc}
\bm{\Lambda} & \bm{0}\\
\bm{0} & \bm{\Lambda}_{\perp}
\end{array}\right]\left[\begin{array}{c}
\bm{V}^{\top}\\
\bm{V}_{\perp}^{\top}
\end{array}\right]=\bm{U}\bm{\Sigma}^{2}\bm{V}^{\top}+\bm{U}_{\perp}\bm{\Sigma}_{\perp}^{2}\bm{V}_{\perp}^{\top}.
\]
From \citet[Lemma 2.1.2]{chen2020spectral} and (\ref{eq:hpca-general-davis-kahan}),
we have learned that 
\[
\left\Vert \bm{U}^{\natural\top}\bm{U}_{\perp}\right\Vert =\left\Vert \bm{U}\bm{U}^{\top}-\bm{U}^{\natural}\bm{U}^{\natural\top}\right\Vert \lesssim\frac{\zeta_{\mathsf{op}}}{\sigma_{r}^{\natural2}}.
\]
In view of Weyl's inequality and Lemma \ref{lemma:hpca-estimation},
we obtain 
\[
\left\Vert \bm{\Lambda}_{\perp}\right\Vert \leq\lambda_{r+1}\left(\bm{G}^{\natural}\right)+\left\Vert \bm{G}-\bm{G}^{\natural}\right\Vert \lesssim\zeta_{\mathsf{op}},
\]
thus indicating that 
\[
\alpha_{2}=\left\Vert \bm{U}^{\natural\top}\bm{U}_{\perp}\bm{\Lambda}_{\perp}\bm{U}_{\perp}^{\top}\bm{U}^{\natural}\right\Vert \leq\left\Vert \bm{U}^{\natural\top}\bm{U}_{\perp}\right\Vert ^{2}\left\Vert \bm{\Lambda}_{\perp}\right\Vert \lesssim\frac{\zeta_{\mathsf{op}}^{3}}{\sigma_{r}^{\natural4}}.
\]

\paragraph{Step 3: bounding $\alpha_{3}$.}

Let us decompose 
\begin{align*}
\bm{G} & =\mathcal{P}_{\mathsf{off}\text{-}\mathsf{diag}}\left[\left(\bm{M}^{\natural}+\bm{E}\right)\left(\bm{M}^{\natural}+\bm{E}\right)^{\top}\right]+\mathcal{P}_{\mathsf{diag}}\left(\bm{G}^{\natural}\right)+\mathcal{P}_{\mathsf{diag}}\left(\bm{G}-\bm{G}^{\natural}\right)\\
 & =\bm{G}^{\natural}+\mathcal{P}_{\mathsf{off}\text{-}\mathsf{diag}}\left(\bm{E}\bm{M}^{\natural\top}+\bm{M}^{\natural}\bm{E}^{\top}+\bm{E}\bm{E}^{\top}\right)+\mathcal{P}_{\mathsf{diag}}\left(\bm{G}-\bm{G}^{\natural}\right)\\
 & =\bm{G}^{\natural}+\bm{E}\bm{M}^{\natural\top}+\bm{M}^{\natural}\bm{E}^{\top}+\mathcal{P}_{\mathsf{off}\text{-}\mathsf{diag}}\left(\bm{E}\bm{E}^{\top}\right)+\mathcal{P}_{\mathsf{diag}}\left(\bm{G}-\bm{G}^{\natural}\right)-\mathcal{P}_{\mathsf{diag}}\left(\bm{E}\bm{M}^{\natural\top}+\bm{M}^{\natural}\bm{E}^{\top}\right),
\end{align*}
which in turn implies that 
\begin{align*}
\bm{U}^{\natural\top}\bm{G}\bm{U}^{\natural}-\bm{\Sigma}^{\natural2} & =\underbrace{\bm{U}^{\natural\top}\left(\bm{E}\bm{M}^{\natural\top}+\bm{M}^{\natural}\bm{E}^{\top}\right)\bm{U}^{\natural}}_{\eqqcolon\bm{J}_{1}}+\underbrace{\bm{U}^{\natural\top}\mathcal{P}_{\mathsf{off}\text{-}\mathsf{diag}}\left(\bm{E}\bm{E}^{\top}\right)\bm{U}^{\natural}}_{\eqqcolon\bm{J}_{2}}\\
 & \quad+\underbrace{\bm{U}^{\natural\top}\mathcal{P}_{\mathsf{diag}}\left(\bm{G}-\bm{G}^{\natural}\right)\bm{U}^{\natural}}_{\eqqcolon\bm{J}_{3}}-\underbrace{\bm{U}^{\natural\top}\mathcal{P}_{\mathsf{diag}}\left(\bm{E}\bm{M}^{\natural\top}+\bm{M}^{\natural}\bm{E}^{\top}\right)\bm{U}^{\natural}}_{\eqqcolon\bm{J}_{4}}.
\end{align*}
We shall then bound the spectral norm of $\bm{J}_{1}$, $\bm{J}_{2}$,
$\bm{J}_{3}$ and $\bm{J}_{4}$ separately.

{\bigskip{}
}\emph{Step 3.1: bounding $\Vert\bm{J}_{1}\Vert$}. Note that 
\[
\left\Vert \bm{U}^{\natural\top}\bm{E}\bm{M}^{\natural\top}\bm{U}^{\natural}\right\Vert =\left\Vert \bm{U}^{\natural\top}\bm{E}\bm{V}^{\natural}\bm{\Sigma}^{\natural}\right\Vert \leq\sigma_{1}^{\natural}\left\Vert \bm{U}^{\natural\top}\bm{E}\bm{V}^{\natural}\right\Vert .
\]
It boils down to controlling the spectral norm of $\bm{U}^{\natural\top}\bm{E}\bm{V}^{\natural}$.
This matrix admits the following decomposition 
\begin{equation}
\bm{U}^{\natural\top}\bm{E}\bm{V}^{\natural}=\sum_{i=1}^{n_{1}}\sum_{j=1}^{n_{2}}E_{i,j}\bm{U}_{i,\cdot}^{\natural\top}\bm{V}_{j,\cdot}^{\natural},\label{eq:UEV-decompose-123}
\end{equation}
which can be controlled through the matrix Bernstein inequality. To
do so, we are in need of calculating the following quantities: 
\begin{align*}
L & \coloneqq\max_{\left(i,j\right)\in\left[n_{1}\right]\times\left[n_{2}\right]}\left\Vert E_{i,j}\bm{U}_{i,\cdot}^{\natural\top}\bm{V}_{j,\cdot}^{\natural}\right\Vert \leq B\sqrt{\frac{\mu^{\natural}r}{n_{1}}}\sqrt{\frac{\mu^{\natural}r}{n_{2}}}=\frac{B\mu^{\natural}r}{\sqrt{n_{1}n_{2}}},\\
V & \coloneqq\max\left\{ \left\Vert \sum_{i=1}^{n_{1}}\sum_{j=1}^{n_{2}}\mathbb{E}\left[\left(E_{i,j}\big(\bm{U}_{i,\cdot}^{\natural}\big)^{\top}\bm{V}_{j,\cdot}^{\natural}\right)\left(E_{i,j}\big(\bm{U}_{i,\cdot}^{\natural}\big)^{\top}\bm{V}_{j,\cdot}^{\natural}\right)^{\top}\right]\right\Vert ,\right.\\
 & \qquad\qquad\qquad\left.\left\Vert \sum_{i=1}^{n_{1}}\sum_{j=1}^{n_{2}}\mathbb{E}\left[\left(E_{i,j}\big(\bm{U}_{i,\cdot}^{\natural}\big)^{\top}\bm{V}_{j,\cdot}^{\natural}\right)^{\top}\left(E_{i,j}\big(\bm{U}_{i,\cdot}^{\natural}\big)^{\top}\bm{V}_{j,\cdot}^{\natural}\right)\right]\right\Vert \right\} \\
 & \leq\max\left\{ \left\Vert \sum_{i=1}^{n_{1}}\sum_{j=1}^{n_{2}}\sigma^{2}\left\Vert \bm{V}_{j,\cdot}^{\natural}\right\Vert _{2}^{2}\big(\bm{U}_{i,\cdot}^{\natural}\big)^{\top}\bm{U}_{i,\cdot}^{\natural}\right\Vert ,\left\Vert \sum_{i=1}^{n_{1}}\sum_{j=1}^{n_{2}}\sigma^{2}\left\Vert \bm{U}_{j,\cdot}^{\natural}\right\Vert _{2}^{2}\big(\bm{V}_{i,\cdot}^{\natural}\big)^{\top}\bm{V}_{i,\cdot}^{\natural}\right\Vert \right\} \\
 & =\sigma^{2}\max\left\{ \left\Vert \bm{V}^{\natural}\right\Vert _{\mathrm{F}}^{2}\left\Vert \bm{U}^{\natural\top}\bm{U}^{\natural}\right\Vert ,\left\Vert \bm{U}^{\natural}\right\Vert _{\mathrm{F}}^{2}\left\Vert \bm{V}^{\natural\top}\bm{V}^{\natural}\right\Vert \right\} =\sigma^{2}r,
\end{align*}
where we have made use of Assumptions~\ref{assumption:subspace-incoherence}-\ref{assumption:subspace-noise}.
Applying matrix Bernstein inequality \citep[Theorem 6.1.1]{Tropp:2015:IMC:2802188.2802189}
to the decomposition (\ref{eq:UEV-decompose-123}) then tells us that,
with probability exceeding $1-O(n^{-10})$, 
\begin{equation}
\left\Vert \bm{U}^{\natural\top}\bm{E}\bm{V}^{\natural}\right\Vert \lesssim\sqrt{V\log n}+L\log n\lesssim\sigma\sqrt{r\log n}+\frac{B\mu^{\natural}r\log n}{\sqrt{n_{1}n_{2}}}\lesssim\sigma\sqrt{r\log n}+\frac{\sigma\mu^{\natural}r\sqrt{\log n}}{\sqrt[4]{n_{1}n_{2}}},\label{eq:U-top-E-V-spectral}
\end{equation}
provided that $B\lesssim\sigma\sqrt[4]{n_{1}n_{2}}/\sqrt{\log n}$.
Therefore, we obtain 
\begin{align*}
\left\Vert \bm{J}_{1}\right\Vert  & \leq2\left\Vert \bm{U}^{\natural\top}\bm{E}\bm{M}^{\natural\top}\bm{U}^{\natural}\right\Vert =2\left\Vert \bm{U}^{\natural\top}\bm{E}\bm{V}^{\natural\top}\bm{\Sigma}^{\natural}\right\Vert \leq2\sigma_{1}^{\natural}\left\Vert \bm{U}^{\natural\top}\bm{E}\bm{V}^{\natural}\right\Vert \\
 & \lesssim\sigma\sigma_{1}^{\natural}\sqrt{r\log n}+\frac{\sigma\sigma_{1}^{\natural}\mu^{\natural}r\sqrt{\log n}}{\sqrt[4]{n_{1}n_{2}}}\lesssim\sqrt{\frac{\mu^{\natural}r}{n_{1}}}\zeta_{\mathsf{op}}
\end{align*}
as long as $n_{1}n_{2}\gtrsim\mu^{\natural2}r^{2}$.

{\bigskip{}
}\emph{Step 3.2: bounding $\Vert\bm{J}_{2}\Vert$.} This matrix $\bm{J}_{2}$
can be expressed as a sum of independent and zero-mean random matrices:
\begin{align*}
\bm{J}_{2} & =\sum_{l=1}^{n_{2}}\underbrace{\bm{U}^{\natural\top}\left(\bm{E}_{\cdot,l}\bm{E}_{\cdot,l}^{\top}-\bm{D}_{l}\right)\bm{U}^{\natural}}_{\eqqcolon\bm{X}_{l}};
\end{align*}
here and below, we denote 
\begin{equation}
\bm{D}_{l}\coloneqq\mathsf{diag}\big\{ E_{1,l}^{2},\ldots,E_{n_{1},l}^{2}\big\}.\label{eq:defn-Dl-HPCA}
\end{equation}
We will invoke the truncated matrix Bernstein inequality \citep[Theorem 3.1.1]{chen2020spectral}
to bound its spectral norm. To do so, we need to calculate the following
quantities: 
\begin{itemize}
\item We first study 
\begin{align*}
v & \coloneqq\left\Vert \sum_{l=1}^{n_{2}}\mathbb{E}\left[\bm{X}_{l}\bm{X}_{l}^{\top}\right]\right\Vert \leq\sup_{\left\Vert \bm{v}\right\Vert _{2}=1}\sum_{l=1}^{n_{2}}\mathbb{E}\left[\bm{v}^{\top}\bm{U}^{\natural\top}\left(\bm{E}_{\cdot,l}\bm{E}_{\cdot,l}^{\top}-\bm{D}_{l}\right)\bm{U}^{\natural}\bm{U}^{\natural\top}\left(\bm{E}_{\cdot,l}\bm{E}_{\cdot,l}^{\top}-\bm{D}_{l}\right)\bm{U}^{\natural}\bm{v}\right]\\
 & =\sup_{\left\Vert \bm{v}\right\Vert _{2}=1}\sum_{l=1}^{n_{2}}\sum_{k=1}^{r}\mathbb{E}\left[\bm{v}^{\top}\bm{U}^{\natural\top}\left(\bm{E}_{\cdot,l}\bm{E}_{\cdot,l}^{\top}-\bm{D}_{l}\right)\bm{U}_{\cdot,k}^{\natural}\bm{U}_{\cdot,k}^{\natural\top}\left(\bm{E}_{\cdot,l}\bm{E}_{\cdot,l}^{\top}-\bm{D}_{l}\right)\bm{U}^{\natural}\bm{v}\right]\\
 & =\sup_{\left\Vert \bm{v}\right\Vert _{2}=1}\sum_{l=1}^{n_{2}}\sum_{k=1}^{r}\mathbb{E}\left[\left(\bm{v}^{\top}\bm{U}^{\natural\top}\left(\bm{E}_{\cdot,l}\bm{E}_{\cdot,l}^{\top}-\bm{D}_{l}\right)\bm{U}_{\cdot,k}^{\natural}\right)^{2}\right].
\end{align*}
For any $\bm{v}\in\mathbb{R}^{r}$ with unit norm, let $\bm{x}=\bm{U}^{\natural}\bm{v}$
and derive 
\begin{align*}
 & \sum_{l=1}^{n_{2}}\sum_{k=1}^{r}\mathbb{E}\left[\left(\bm{v}^{\top}\bm{U}^{\natural\top}\left(\bm{E}_{\cdot,l}\bm{E}_{\cdot,l}^{\top}-\bm{D}_{l}\right)\bm{U}_{\cdot,k}^{\natural}\right)^{2}\right]=\sum_{l=1}^{n_{2}}\sum_{k=1}^{r}\mathbb{E}\left[\left(\bm{x}^{\top}\left(\bm{E}_{\cdot,l}\bm{E}_{\cdot,l}^{\top}-\bm{D}_{l}\right)\bm{U}_{\cdot,k}^{\natural}\right)^{2}\right]\\
 & \quad=\sum_{l=1}^{n_{2}}\sum_{k=1}^{r}\mathsf{var}\left(\bm{x}^{\top}\bm{E}_{\cdot,l}\bm{E}_{\cdot,l}^{\top}\bm{U}_{\cdot,k}^{\natural}\right)=\sum_{l=1}^{n_{2}}\sum_{k=1}^{r}\mathsf{var}\left(\sum_{i=1}^{n_{1}}\sum_{j=1}^{n_{1}}x_{i}U_{j,k}^{\natural}E_{i,l}E_{j,l}\right)\\
 & \quad=\sum_{l=1}^{n_{2}}\sum_{k=1}^{r}\sum_{i=1}^{n_{1}}\sum_{j=1}^{n_{1}}x_{i}^{2}U_{j,k}^{\natural2}\sigma_{i,l}^{2}\sigma_{j,l}^{2}\leq\sigma^{4}n_{2}\left\Vert \bm{x}\right\Vert _{2}^{2}\left\Vert \bm{U}^{\natural}\right\Vert _{\mathrm{F}}^{2}=\sigma^{4}n_{2}r.
\end{align*}
To sum up, we have obtained 
\[
v\leq\sigma^{4}n_{2}r.
\]
\item Note that for each $l\in[n_{2}]$, one can decompose 
\begin{align*}
\left\Vert \bm{X}_{l}\right\Vert  & \leq\left\Vert \bm{U}^{\natural\top}\bm{E}_{\cdot,l}\right\Vert _{2}^{2}+\left\Vert \bm{U}^{\natural\top}\bm{D}_{l}\bm{U}^{\natural}\right\Vert =\underbrace{\left\Vert \sum_{i=1}^{n_{1}}E_{i,l}\bm{U}_{i,\cdot}^{\natural}\right\Vert _{2}^{2}}_{\eqqcolon\beta_{1}}+\underbrace{\left\Vert \sum_{i=1}^{n_{1}}E_{i,l}^{2}\bm{U}_{i,\cdot}^{\natural\top}\bm{U}_{i,\cdot}^{\natural}\right\Vert }_{\eqqcolon\beta_{2}}.
\end{align*}

\begin{itemize}
\item For $\beta_{1}$, it is straightforward to calculate 
\begin{align*}
L_{\beta_{1}} & \coloneqq\max_{i\in[n_{1}]}\left\Vert E_{i,l}\bm{U}_{i,\cdot}^{\natural}\right\Vert _{2}\leq B\left\Vert \bm{U}^{\natural}\right\Vert _{2,\infty}\leq B\sqrt{\frac{\mu^{\natural}r}{n_{1}}},\\
V_{\beta_{1}} & \coloneqq\sum_{i=1}^{n_{1}}\mathbb{E}\left[E_{i,l}^{2}\right]\left\Vert \bm{U}_{i,\cdot}^{\natural}\right\Vert _{2}^{2}\leq\sigma^{2}\left\Vert \bm{U}^{\natural}\right\Vert _{\mathrm{F}}^{2}=\sigma^{2}r.
\end{align*}
The matrix Bernstein inequality \citep[Theorem 6.1.1]{Tropp:2015:IMC:2802188.2802189}
tells us that 
\[
\mathbb{P}\left(\beta_{1}\geq t\right)\leq r\exp\left(\frac{-t/2}{V_{\beta_{1}}+L_{\beta_{1}}\sqrt{t}/3}\right).
\]
\item Regarding $\beta_{2}$, we can also calculate 
\begin{align*}
L_{\beta_{2}} & \coloneqq\max_{i\in[n_{1}]}\left\Vert E_{i,l}^{2}\bm{U}_{i,\cdot}^{\natural\top}\bm{U}_{i,\cdot}^{\natural}\right\Vert _{2}\leq B^{2}\left\Vert \bm{U}^{\natural}\right\Vert _{2,\infty}^{2}\leq B^{2}\frac{\mu^{\natural}r}{n_{1}},\\
V_{\beta_{2}} & \coloneqq\left\Vert \sum_{i=1}^{n_{1}}\mathbb{E}\left[E_{i,l}^{4}\right]\left\Vert \bm{U}_{i,\cdot}^{\natural}\right\Vert _{2}^{2}\bm{U}_{i,\cdot}^{\natural\top}\bm{U}_{i,\cdot}^{\natural}\right\Vert \leq\sigma^{2}B^{2}\left\Vert \bm{U}^{\natural}\right\Vert _{2,\infty}^{2}\leq\sigma^{2}B^{2}\frac{\mu^{\natural}r}{n_{1}}.
\end{align*}
By virtue of the matrix Bernstein inequality \citep[Theorem 6.1.1]{Tropp:2015:IMC:2802188.2802189},
we can obtain 
\[
\mathbb{P}\left(\beta_{2}\geq t\right)\leq r\exp\left(\frac{-t^{2}/2}{V_{\beta_{2}}+L_{\beta_{2}}t/3}\right).
\]
\item Combine these tail probability bounds for $\beta_{1}$ and $\beta_{2}$
to achieve 
\begin{align*}
\mathbb{P}\left(\left\Vert \bm{X}_{l}\right\Vert \geq t\right) & \leq\mathbb{P}\left(\beta_{1}\geq\frac{t}{2}\right)+\mathbb{P}\left(\beta_{2}\geq\frac{t}{2}\right)\\
 & \leq r\exp\left(\frac{-t/4}{V_{\beta_{1}}+L_{\beta_{1}}\sqrt{t}/5}\right)+r\exp\left(\frac{-t^{2}/8}{V_{\beta_{2}}+L_{\beta_{2}}t/6}\right)\\
 & \leq r\exp\left(-\min\left\{ \frac{t}{8V_{\beta_{1}}},\frac{5\sqrt{t}}{8L_{\beta_{1}}}\right\} \right)+r\exp\left(-\min\left\{ \frac{t^{2}}{16V_{\beta_{2}}},\frac{3t}{8L_{\beta_{2}}}\right\} \right).
\end{align*}
Therefore, we can set 
\[
L\coloneqq\widetilde{C}\left(V_{\beta_{1}}\log^{2}n+L_{\beta_{1}}^{2}\log^{2}n+\sqrt{V_{\beta_{2}}}\log n+L_{\beta_{2}}\log n\right)\asymp\sigma^{2}r\log^{2}n+\frac{B^{2}\mu^{\natural}r}{n_{1}}\log^{2}n
\]
for some sufficiently large constant $\widetilde{C}>0$, so as to
guarantee that 
\[
\mathbb{P}\left(\left\Vert \bm{X}_{l}\right\Vert \geq L\right)\leq q_{0}\coloneqq\frac{1}{n^{10}},\qquad\text{for all }l\in[n_{1}].
\]
\end{itemize}
\item Based on the above choice of $L$, we can see that 
\begin{align*}
q_{1} & \coloneqq\left\Vert \mathbb{E}\left[\bm{X}_{l}\ind\left\{ \left\Vert \bm{X}_{l}\right\Vert \geq L\right\} \right]\right\Vert \leq\mathbb{E}\left[\left\Vert \bm{X}_{l}\right\Vert \ind\left\{ \left\Vert \bm{X}_{l}\right\Vert \geq L\right\} \right]\\
 & =\int_{0}^{\infty}\mathbb{P}\left(\left\Vert \bm{X}_{l}\right\Vert \ind\left\{ \left\Vert \bm{X}_{l}\right\Vert \geq L\right\} >t\right)\mathrm{d}t\\
 & =\int_{0}^{L}\mathbb{P}\left(\left\Vert \bm{X}_{l}\right\Vert \geq L\right)\mathrm{d}t+\int_{L}^{\infty}\mathbb{P}\left(\left\Vert \bm{X}_{l}\right\Vert >t\right)\mathrm{d}t\\
 & \leq\frac{L}{n^{10}}+\int_{L}^{\infty}\mathbb{P}\left(\left\Vert \bm{X}_{l}\right\Vert >t\right)\mathrm{d}t,
\end{align*}
where the second inequality uses Jensen's inequality. Notice that
for $t\geq L$, we have $t\gg V_{\beta_{1}}$ and $t\gg\sqrt{V_{\beta_{2}}}$
as long as the constant $\widetilde{C}$ is sufficiently large. As
a consequence, 
\begin{align*}
\mathbb{P}\left(\left\Vert \bm{X}_{l}\right\Vert \geq t\right) & \leq r\exp\left(-\min\left\{ \frac{t}{8V_{\beta_{1}}},\frac{5\sqrt{t}}{8L_{\beta_{1}}}\right\} \right)+r\exp\left(-\min\left\{ \frac{t^{2}}{16V_{\beta_{2}}},\frac{3t}{8L_{\beta_{2}}}\right\} \right)\\
 & \leq r\exp\left(-\min\left\{ \frac{\sqrt{t}}{\sqrt{V_{\beta_{1}}}},\frac{5\sqrt{t}}{8L_{\beta_{1}}}\right\} \right)+r\exp\left(-\min\left\{ \frac{t}{\sqrt{V_{\beta_{2}}}},\frac{3t}{8L_{\beta_{2}}}\right\} \right)\\
 & \leq r\exp\left(-\frac{\sqrt{t}}{\max\left\{ \sqrt{V_{\beta_{1}}},2L_{\beta_{1}}\right\} }\right)+r\exp\left(-\frac{t}{\max\left\{ \sqrt{V_{\beta_{2}}},3L_{\beta_{2}}\right\} }\right),
\end{align*}
provided that $\widetilde{C}$ is sufficiently large. Consequently,
we can deduce that 
\begin{align*}
\int_{L}^{\infty}\mathbb{P}\left(\left\Vert \bm{X}_{l}\right\Vert \geq t\right)\mathrm{d}t & \leq r\underbrace{\int_{L}^{\infty}\exp\left(-\frac{\sqrt{t}}{\max\left\{ \sqrt{V_{\beta_{1}}},2L_{\beta_{1}}\right\} }\right)\mathrm{d}t}_{\eqqcolon I_{1}}+r\underbrace{\int_{L}^{\infty}\exp\left(-\frac{t}{\max\left\{ \sqrt{V_{\beta_{2}}},3L_{\beta_{2}}\right\} }\right)\mathrm{d}t}_{\eqqcolon I_{2}}.
\end{align*}

\begin{itemize}
\item The first integral obeys 
\begin{align*}
I_{1} & \overset{\text{(i)}}{=}2\max\left\{ \sqrt{V_{\beta_{1}}},2L_{\beta_{1}}\right\} \sqrt{L}\exp\left(-\frac{\sqrt{L}}{\max\left\{ \sqrt{V_{\beta_{1}}},2L_{\beta_{1}}\right\} }\right)\\
 & \quad+2\max\left\{ \sqrt{V_{\beta_{1}}},2L_{\beta_{1}}\right\} ^{2}\exp\left(-\frac{\sqrt{L}}{\max\left\{ \sqrt{V_{\beta_{1}}},2L_{\beta_{1}}\right\} }\right)\\
 & \leq\left(\frac{4}{\sqrt{\widetilde{C}}}+\frac{8}{\widetilde{C}}\right)L\exp\left(-\frac{\sqrt{\widetilde{C}}}{2}\log n\right)\overset{\text{(iii)}}{\leq}\frac{L}{n^{20}}.
\end{align*}
Here, (i) follows from the following formula that holds for any constants
$\alpha,\beta>0$: 
\begin{align*}
\int_{\beta}^{\infty}\exp\left(-\alpha\sqrt{x}\right)\mathrm{d}x & \overset{y=\sqrt{x}}{=}2\int_{\sqrt{\beta}}^{\infty}y\exp\left(-\alpha y\right)\mathrm{d}y=\left[-\frac{2}{\alpha}y\exp\left(-\alpha y\right)\right]\Bigg|_{\sqrt{\beta}}^{\infty}+\frac{2}{\alpha}\int_{\sqrt{\beta}}^{\infty}\exp\left(-\alpha y\right)\mathrm{d}y\\
 & =\frac{2\sqrt{\beta}}{\alpha}\exp\left(-\alpha\sqrt{\beta}\right)+\left[-\frac{2}{\alpha^{2}}\exp\left(-\alpha y\right)\right]\Bigg|_{\sqrt{\beta}}^{\infty}\\
 & =\frac{2\sqrt{\beta}}{\alpha}\exp\left(-\alpha\sqrt{\beta}\right)+\frac{2}{\alpha^{2}}\exp\left(-\alpha\sqrt{\beta}\right);
\end{align*}
(ii) comes from the definition of $L$: 
\[
\sqrt{L}\geq\sqrt{\widetilde{C}\left(V_{\beta_{1}}\log^{2}n+L_{\beta_{1}}^{2}\log^{2}n\right)}\geq\sqrt{\widetilde{C}}\max\left\{ \sqrt{V_{\beta_{1}}},L_{\beta_{1}}\right\} \log n;
\]
and (iii) holds provided that $\widetilde{C}$ is sufficiently large. 
\item The second integral satisfies 
\begin{align*}
I_{2} & =\max\left\{ \sqrt{V_{\beta_{2}}},3L_{\beta_{2}}\right\} \exp\left(-\frac{L}{\max\left\{ \sqrt{V_{\beta_{2}}},3L_{\beta_{2}}\right\} }\right)\leq\frac{3L}{\widetilde{C}}\exp\left(-\frac{\widetilde{C}}{3}\log n\right)\leq\frac{L}{n^{20}}.
\end{align*}
Here the penultimate inequality follows from 
\[
L\geq\widetilde{C}\left(\sqrt{V_{\beta_{2}}}\log n+L_{\beta_{2}}\log n\right)\geq\widetilde{C}\max\left\{ \sqrt{V_{\beta_{2}}},L_{\beta_{2}}\right\} \log n,
\]
and the last inequality holds provided that $\widetilde{C}$ is sufficiently
large. 
\end{itemize}
Therefore we conclude that 
\begin{align*}
q_{1} & \leq\frac{L}{n^{10}}+rI_{1}+rI_{2}\leq\frac{L}{n^{10}}+2r\frac{L}{n^{20}}\leq\frac{L}{n^{9}}.
\end{align*}

\item With the above quantities in mind, we are ready to use the truncated
matrix Bernstein inequality \citep[Theorem 3.1.1]{chen2020spectral}
to show that with probability exceeding $1-O(n^{-10})$, 
\begin{align*}
\left\Vert \bm{J}_{2}\right\Vert  & \lesssim\sqrt{v\log n}+L\log n+nq_{1}\lesssim\sigma^{2}\sqrt{n_{2}r\log n}+\sigma^{2}r\log^{3}n+\frac{B^{2}\mu^{\natural}r}{n_{1}}\log^{3}n\\
 & \overset{\text{(i)}}{\lesssim}\sigma^{2}\sqrt{n_{2}r\log n}+\sigma^{2}r\log^{3}n+\frac{\sigma^{2}\sqrt{n_{1}n_{2}}\mu^{\natural}r}{n_{1}}\log^{2}n\\
 & \asymp\left(\sqrt{\frac{r}{n_{1}}}+\frac{r\log^{2}n}{\sqrt{n_{1}n_{2}}}+\frac{\mu^{\natural}r\log n}{n_{1}}\right)\sigma^{2}\sqrt{n_{1}n_{2}}\log n\overset{\text{(ii)}}{\lesssim}\sqrt{\frac{r}{n_{1}}}\zeta_{\mathsf{op}}.
\end{align*}
Here, (i) makes use of the noise condition $B\lesssim\sigma\sqrt[4]{n_{1}n_{2}}/\sqrt{\log n}$,
whereas (ii) holds provided that $n_{1}\gtrsim\mu^{\natural2}r\log^{2}dn$
and $n_{2}\gtrsim r\log^{4}n$. 
\end{itemize}
{\bigskip{}
}\emph{Step 3.3: bounding $\Vert\bm{J}_{3}\Vert$.} Regarding $\bm{J}_{3}$,
it is easy to show that 
\[
\left\Vert \bm{J}_{3}\right\Vert \leq\left\Vert \mathcal{P}_{\mathsf{diag}}\left(\bm{G}-\bm{G}^{\natural}\right)\right\Vert \lesssim\kappa^{\natural2}\sqrt{\frac{\mu^{\natural}r}{n_{1}}}\zeta_{\mathsf{op}},
\]
which results from Lemma \ref{lemma:hpca-estimation}.

{\bigskip{}
}\emph{Step 3.4: bounding $\Vert\bm{J}_{4}\Vert$.} We are now left
with bounding $\bm{J}_{4}$. Note that 
\[
\left\Vert \bm{J}_{4}\right\Vert \leq\left\Vert \mathcal{P}_{\mathsf{diag}}\left(\bm{E}\bm{M}^{\natural\top}+\bm{M}^{\natural}\bm{E}^{\top}\right)\right\Vert \leq2\max_{i\in[n_{1}]}\left|\sum_{j=1}^{n_{2}}E_{i,j}M_{i,j}^{\natural}\right|.
\]
For any $i\in[n_{1}]$, it is straightforward to calculate that 
\begin{align*}
L_{i} & \coloneqq\max_{j\in[n_{2}]}\left|E_{i,j}M_{i,j}^{\natural}\right|\leq B\left\Vert \bm{M}^{\natural}\right\Vert _{\infty},\\
V_{i} & \coloneqq\mathsf{var}\left(\sum_{j=1}^{n_{2}}E_{i,j}M_{i,j}^{\natural}\right)=\sum_{j=1}^{n_{2}}\sigma_{i,j}^{2}\big|M_{i,j}^{\natural}\big|^{2}\leq\sigma^{2}\left\Vert \bm{M}^{\natural}\right\Vert _{2,\infty}^{2}.
\end{align*}
In view of the Bernstein inequality \citep[Theorem 2.8.4]{vershynin2016high},
\begin{align}
\left|\sum_{j=1}^{n_{2}}E_{i,j}M_{i,j}^{\natural}\right| & \lesssim\sqrt{V_{i}\log n}+L_{i}\log n\lesssim\sigma\left\Vert \bm{M}^{\natural}\right\Vert _{2,\infty}\sqrt{\log n}+B\left\Vert \bm{M}^{\natural}\right\Vert _{\infty}\log n\nonumber \\
 & \lesssim\sigma\sigma_{1}^{\natural}\left\Vert \bm{U}^{\natural}\right\Vert _{2,\infty}\sqrt{\log n}+\sigma\sqrt{n_{2}}\sqrt{\frac{\mu^{\natural}r}{n_{1}n_{2}}}\sigma_{1}^{\natural}\nonumber \\
 & \lesssim\sigma\sigma_{1}^{\natural}\sqrt{\frac{\mu^{\natural}r\log n}{n_{1}}}\lesssim\frac{\sqrt{\mu^{\natural}r}}{n_{1}}\zeta_{\mathsf{op}}\label{eq:pca-Eij-Aij-sum}
\end{align}
with probability exceeding $1-O(n^{-11})$. Combining the above bounds
and applying the union bound show that with probability exceeding
$1-O(n^{-10})$, 
\[
\left\Vert \bm{J}_{4}\right\Vert \lesssim\frac{\sqrt{\mu^{\natural}r}}{n_{1}}\zeta_{\mathsf{op}}.
\]

{\bigskip{}
}\emph{Step 3.5: putting all this together. }Taking the previous
bounds on the spectral norm of $\bm{J}_{1}$, $\bm{J}_{2}$, $\bm{J}_{3}$,
$\bm{J}_{4}$ collectively yields 
\begin{align*}
\alpha_{3} & \leq\left\Vert \bm{J}_{1}\right\Vert +\left\Vert \bm{J}_{2}\right\Vert +\left\Vert \bm{J}_{3}\right\Vert +\left\Vert \bm{J}_{4}\right\Vert \lesssim\kappa^{\natural2}\sqrt{\frac{\mu^{\natural}r}{n_{1}}}\zeta_{\mathsf{op}}+\frac{\mu^{\natural}r}{n_{1}}\zeta_{\mathsf{op}}\asymp\kappa^{\natural2}\sqrt{\frac{\mu^{\natural}r}{n_{1}}}\zeta_{\mathsf{op}},
\end{align*}
where the last relation holds as long as $n_{1}\gtrsim\mu^{\natural}r$.

\paragraph{Step 4: combining the bounds on $\alpha_{1}$, $\alpha_{2}$ and
$\alpha_{3}$. }

Taking the bounds on $\alpha_{1}$, $\alpha_{2}$, $\alpha_{3}$ together
leads to 
\begin{align*}
\left\Vert \bm{R}^{\top}\bm{\Sigma}^{2}\bm{R}-\bm{\Sigma}^{\natural2}\right\Vert  & \leq\alpha_{1}+\alpha_{2}+\alpha_{3}\lesssim\kappa^{\natural2}\sqrt{\frac{\mu^{\natural}r}{n_{1}}}\zeta_{\mathsf{op}}+\kappa^{\natural2}\frac{\zeta_{\mathsf{op}}^{2}}{\sigma_{r}^{\natural2}}+\frac{\zeta_{\mathsf{op}}^{3}}{\sigma_{r}^{\natural4}}\asymp\kappa^{\natural2}\sqrt{\frac{\mu^{\natural}r}{n_{1}}}\zeta_{\mathsf{op}}+\kappa^{\natural2}\frac{\zeta_{\mathsf{op}}^{2}}{\sigma_{r}^{\natural2}}
\end{align*}
with probability exceeding $1-O(n^{-10})$, where the last relation
is valid under the condition that $\zeta_{\mathsf{op}}\ll\sigma_{r}^{\natural2}$.

\subsubsection{Proof of Lemma \ref{lemma:hpca-estimation-two-to-infty}\label{appendix:proof-lemma-hpca-estimation-two-to-infty}}

To begin with, observe that $\bm{G}^{\natural}=\bm{U}^{\natural}\bm{\Sigma}^{\natural2}\bm{U}^{\natural\top}$
(and hence $\bm{G}^{\natural}\bm{U}^{\natural}(\bm{\Sigma}^{\natural})^{-2}=\bm{U}^{\natural}$),
which allows us to decompose 
\begin{align*}
\left\Vert \bm{U}_{m,\cdot}\bm{H}-\bm{U}_{m,\cdot}^{\natural}\right\Vert _{2} & =\left\Vert \left(\bm{U}\bm{H}-\bm{G}\bm{U}^{\natural}\left(\bm{\Sigma}^{\natural}\right)^{-2}+\bm{G}\bm{U}^{\natural}\left(\bm{\Sigma}^{\natural}\right)^{-2}-\bm{G}^{\natural}\bm{U}^{\natural}\left(\bm{\Sigma}^{\natural}\right)^{-2}\right)_{m,\cdot}\right\Vert _{2}\\
 & \leq\underbrace{\left\Vert \left(\bm{U}\bm{H}-\bm{G}\bm{U}^{\natural}\left(\bm{\Sigma}^{\natural}\right)^{-2}\right)_{m,\cdot}\right\Vert _{2}}_{\eqqcolon\,\alpha_{1}}+\underbrace{\left\Vert \left(\bm{G}-\bm{G}^{\natural}\right)_{m,\cdot}\bm{U}^{\natural}\left(\bm{\Sigma}^{\natural}\right)^{-2}\right\Vert _{2}}_{\eqqcolon\,\alpha_{2}}.
\end{align*}
We then proceed to bound the terms $\alpha_{1}$ the $\alpha{}_{2}$. 
\begin{itemize}
\item Regarding $\alpha_{1}$, we have the following decomposition
\begin{align*}
\alpha_{1} & \leq\left\Vert \left(\bm{U}\bm{H}\bm{\Sigma}^{\natural2}-\bm{G}\bm{U}^{\natural}\right)_{m,\cdot}\left(\bm{\Sigma}^{\natural}\right)^{-2}\right\Vert _{2}\leq\frac{1}{\sigma_{r}^{\natural2}}\left\Vert \left(\bm{U}\bm{H}\bm{\Sigma}^{\natural2}-\bm{G}\bm{U}^{\natural}\right)_{m,\cdot}\right\Vert _{2}\\
 & \leq\underbrace{\frac{1}{\sigma_{r}^{\natural2}}\left\Vert \bm{U}_{m,\cdot}\left(\bm{H}\bm{\Sigma}^{\natural2}-\bm{\Sigma}^{2}\bm{H}\right)\right\Vert _{2}}_{\eqqcolon\beta_{1}}+\underbrace{\frac{1}{\sigma_{r}^{\natural2}}\left\Vert \left(\bm{U}\bm{\Sigma}^{2}\bm{H}-\bm{G}\bm{U}^{\natural}\right)_{m,\cdot}\right\Vert _{2}}_{\eqqcolon\beta_{2}}.
\end{align*}
We first bound $\beta_{1}$, where we can use the triangle inequality
to achieve
\begin{align*}
\beta_{1} & \leq\frac{1}{\sigma_{r}^{\natural2}}\left\Vert \bm{U}_{m,\cdot}\left(\bm{R}_{\bm{U}}\bm{\Sigma}^{\natural2}-\bm{\Sigma}^{2}\bm{R}_{\bm{U}}\right)\right\Vert _{2}+\frac{1}{\sigma_{r}^{\natural2}}\left\Vert \bm{U}_{m,\cdot}\left(\bm{H}-\bm{R}_{\bm{U}}\right)\bm{\Sigma}^{\natural2}\right\Vert _{2}+\frac{1}{\sigma_{r}^{\natural2}}\left\Vert \bm{U}_{m,\cdot}\bm{\Sigma}^{2}\left(\bm{H}-\bm{R}_{\bm{U}}\right)\right\Vert _{2}\\
 & \leq\frac{1}{\sigma_{r}^{\natural2}}\left\Vert \bm{U}_{m,\cdot}\bm{R}_{\bm{U}}\left(\bm{\Sigma}^{\natural2}-\bm{R}_{\bm{U}}^{\top}\bm{\Sigma}^{2}\bm{R}_{\bm{U}}\right)\right\Vert _{2}+\frac{1}{\sigma_{r}^{\natural2}}\left\Vert \bm{U}_{m,\cdot}\right\Vert _{2}\left\Vert \bm{H}-\bm{R}\right\Vert \left(\left\Vert \bm{\Sigma}^{\natural2}\right\Vert +\left\Vert \bm{\Sigma}^{2}\right\Vert \right)\\
 & \lesssim\frac{1}{\sigma_{r}^{\natural2}}\left\Vert \bm{U}_{m,\cdot}\right\Vert _{2}\left(\left\Vert \bm{\Sigma}^{\natural2}-\bm{R}_{\bm{U}}^{\top}\bm{\Sigma}^{2}\bm{R}_{\bm{U}}\right\Vert +\left\Vert \bm{\Sigma}^{\natural2}\right\Vert \left\Vert \bm{H}-\bm{R}\right\Vert +\left\Vert \bm{\Sigma}^{\natural2}-\bm{R}_{\bm{U}}^{\top}\bm{\Sigma}^{2}\bm{R}_{\bm{U}}\right\Vert \left\Vert \bm{H}-\bm{R}\right\Vert \right)\\
 & \lesssim\frac{1}{\sigma_{r}^{\natural2}}\left\Vert \bm{U}_{m,\cdot}\bm{H}\right\Vert _{2}\left(\kappa^{\natural2}\sqrt{\frac{\mu^{\natural}r}{n_{1}}}\zeta_{\mathsf{op}}+\kappa^{\natural2}\frac{\zeta_{\mathsf{op}}^{2}}{\sigma_{r}^{\natural2}}\right)\\
 & \lesssim\left(\left\Vert \bm{U}_{m,\cdot}^{\natural}\right\Vert _{2}+\left\Vert \bm{U}_{m,\cdot}\bm{H}-\bm{U}_{m,\cdot}^{\natural}\right\Vert _{2}\right)\left(\kappa^{\natural2}\sqrt{\frac{\mu^{\natural}r}{n_{1}}}\frac{\zeta_{\mathsf{op}}}{\sigma_{r}^{\natural2}}+\kappa^{\natural2}\frac{\zeta_{\mathsf{op}}^{2}}{\sigma_{r}^{\natural4}}\right).
\end{align*}
Here the penultimate line follows from Lemma \ref{lemma:hpca-basic-facts},
Lemma \ref{lemma:hpca-approx-2}, and the assumption that $\zeta_{\mathsf{op}}\lesssim\sigma_{r}^{\natural2}$.
Additionally, Lemma \ref{lemma:hpca-approx-1} tells us that 
\begin{align*}
\beta_{2} & \lesssim\frac{\zeta_{\mathsf{op},m}}{\sigma_{r}^{\natural2}}\left(\kappa^{\natural2}\frac{\zeta_{\mathsf{op}}}{\sigma_{r}^{\natural2}}\sqrt{\frac{\mu^{\natural}r}{n_{1}}}+\left\Vert \bm{U}\bm{H}-\bm{U}^{\natural}\right\Vert _{2,\infty}\right)+\kappa^{\natural2}\frac{\zeta_{\mathsf{op}}^{2}}{\sigma_{r}^{\natural4}}\left\Vert \bm{U}_{m,\cdot}^{\natural}\right\Vert _{2}\\
 & \quad+\kappa^{\natural2}\frac{\zeta_{\mathsf{op}}}{\sigma_{r}^{\natural2}}\left\Vert \bm{U}_{m,\cdot}^{\natural}\right\Vert _{2}\left\Vert \bm{U}_{m,\cdot}\bm{H}-\bm{U}_{m,\cdot}^{\natural}\right\Vert _{2}+\kappa^{\natural2}\frac{\zeta_{\mathsf{op}}}{\sigma_{r}^{\natural2}}\left\Vert \bm{U}_{m,\cdot}\bm{H}-\bm{U}_{m,\cdot}^{\natural}\right\Vert _{2}^{2}.
\end{align*}
Therefore, for all $m\in[n_{1}]$ we have
\begin{align*}
\alpha_{1} & \lesssim\beta_{1}+\beta_{2}\\
 & \lesssim\left(\left\Vert \bm{U}_{m,\cdot}^{\natural}\right\Vert _{2}+\left\Vert \bm{U}_{m,\cdot}\bm{H}-\bm{U}_{m,\cdot}^{\natural}\right\Vert _{2}\right)\left(\kappa^{\natural2}\sqrt{\frac{\mu^{\natural}r}{n_{1}}}\frac{\zeta_{\mathsf{op}}}{\sigma_{r}^{\natural2}}+\kappa^{\natural2}\frac{\zeta_{\mathsf{op}}^{2}}{\sigma_{r}^{\natural4}}\right)\\
 & \quad+\frac{\zeta_{\mathsf{op},m}}{\sigma_{r}^{\natural2}}\left(\kappa^{\natural2}\frac{\zeta_{\mathsf{op}}}{\sigma_{r}^{\natural2}}\sqrt{\frac{\mu^{\natural}r}{n_{1}}}+\left\Vert \bm{U}\bm{H}-\bm{U}^{\natural}\right\Vert _{2,\infty}\right)\\
 & \quad+\kappa^{\natural2}\frac{\zeta_{\mathsf{op}}}{\sigma_{r}^{\natural2}}\left\Vert \bm{U}_{m,\cdot}^{\natural}\right\Vert _{2}\left\Vert \bm{U}_{m,\cdot}\bm{H}-\bm{U}_{m,\cdot}^{\natural}\right\Vert _{2}+\kappa^{\natural2}\frac{\zeta_{\mathsf{op}}}{\sigma_{r}^{\natural2}}\left\Vert \bm{U}_{m,\cdot}\bm{H}-\bm{U}_{m,\cdot}^{\natural}\right\Vert _{2}^{2}.
\end{align*}
\item In view of Lemma \ref{lemma:hpca-1}, the second term $\alpha_{2}$
can be bounded by 
\begin{align*}
\alpha_{2} & \leq\frac{1}{\sigma_{r}^{\natural2}}\left\Vert \left(\bm{G}-\bm{G}^{\natural}\right)_{m,\cdot}\bm{U}^{\natural}\right\Vert _{2}\\
 & \lesssim\frac{\zeta_{\mathsf{op},m}}{\sigma_{r}^{\natural2}}\sqrt{\frac{\mu^{\natural}r}{n_{1}}}+\frac{1}{\sigma_{r}^{\natural2}}\sigma\sigma_{1}^{\natural}\sqrt{\mu^{\natural}r\log n}\left\Vert \bm{U}_{m,\cdot}^{\natural}\right\Vert _{2}+\kappa^{\natural2}\frac{\zeta_{\mathsf{op}}}{\sigma_{r}^{\natural2}}\left\Vert \bm{U}_{m,\cdot}^{\natural}\right\Vert _{2}^{2}\\
 & \quad+\kappa^{\natural2}\frac{\zeta_{\mathsf{op}}}{\sigma_{r}^{\natural2}}\left\Vert \bm{U}_{m,\cdot}\bm{H}-\bm{U}_{m,\cdot}^{\natural}\right\Vert _{2}\left\Vert \bm{U}_{m,\cdot}^{\natural}\right\Vert _{2}.
\end{align*}
\end{itemize}
The preceding bounds taken together allow us to conclude that 
\begin{align}
\left\Vert \bm{U}_{m,\cdot}\bm{H}-\bm{U}_{m,\cdot}^{\natural}\right\Vert _{2} & \leq\alpha_{1}+\alpha_{2}\nonumber \\
 & \lesssim\frac{\zeta_{\mathsf{op},m}}{\sigma_{r}^{\natural2}}\left(\sqrt{\frac{\mu^{\natural}r}{n_{1}}}+\left\Vert \bm{U}\bm{H}-\bm{U}^{\natural}\right\Vert _{2,\infty}\right)+\left\Vert \bm{U}_{m,\cdot}^{\natural}\right\Vert _{2}\left(\kappa^{\natural2}\sqrt{\frac{\mu^{\natural}r}{n_{1}}}\frac{\zeta_{\mathsf{op}}}{\sigma_{r}^{\natural2}}+\kappa^{\natural2}\frac{\zeta_{\mathsf{op}}^{2}}{\sigma_{r}^{\natural4}}\right)\nonumber \\
 & \quad+\left\Vert \bm{U}_{m,\cdot}\bm{H}-\bm{U}_{m,\cdot}^{\natural}\right\Vert _{2}\left(\kappa^{\natural2}\sqrt{\frac{\mu^{\natural}r}{n_{1}}}\frac{\zeta_{\mathsf{op}}}{\sigma_{r}^{\natural2}}+\kappa^{\natural2}\frac{\zeta_{\mathsf{op}}^{2}}{\sigma_{r}^{\natural4}}\right)+\kappa^{\natural2}\frac{\zeta_{\mathsf{op}}}{\sigma_{r}^{\natural2}}\left\Vert \bm{U}_{m,\cdot}\bm{H}-\bm{U}_{m,\cdot}^{\natural}\right\Vert _{2}^{2},\label{eq:UH-U-star-2-infty-m}
\end{align}
provided that $\zeta_{\mathsf{op}}\lesssim\sigma_{r}^{\natural2}/\kappa^{\natural2}$and
$n_{1}\gtrsim\mu^{\natural}r$. By taking supremum over $m\in[n_{1}]$,
we have
\begin{align*}
\left\Vert \bm{U}\bm{H}-\bm{U}^{\natural}\right\Vert _{2,\infty} & \lesssim\frac{\zeta_{\mathsf{op}}}{\sigma_{r}^{\natural2}}\left(\sqrt{\frac{\mu^{\natural}r}{n_{1}}}+\left\Vert \bm{U}\bm{H}-\bm{U}^{\natural}\right\Vert _{2,\infty}\right)+\left\Vert \bm{U}^{\natural}\right\Vert _{2,\infty}\left(\kappa^{\natural2}\sqrt{\frac{\mu^{\natural}r}{n_{1}}}\frac{\zeta_{\mathsf{op}}}{\sigma_{r}^{\natural2}}+\kappa^{\natural2}\frac{\zeta_{\mathsf{op}}^{2}}{\sigma_{r}^{\natural4}}\right)\\
 & \quad+\left\Vert \bm{U}\bm{H}-\bm{U}^{\natural}\right\Vert _{2,\infty}\left(\kappa^{\natural2}\sqrt{\frac{\mu^{\natural}r}{n_{1}}}\frac{\zeta_{\mathsf{op}}}{\sigma_{r}^{\natural2}}+\kappa^{\natural2}\frac{\zeta_{\mathsf{op}}^{2}}{\sigma_{r}^{\natural4}}\right)+\kappa^{\natural2}\frac{\zeta_{\mathsf{op}}}{\sigma_{r}^{\natural2}}\left\Vert \bm{U}\bm{H}-\bm{U}^{\natural}\right\Vert _{2,\infty}^{2}\\
 & \lesssim\frac{\zeta_{\mathsf{op}}}{\sigma_{r}^{\natural2}}\sqrt{\frac{\mu^{\natural}r}{n_{1}}}+\frac{\zeta_{\mathsf{op}}}{\sigma_{r}^{\natural2}}\left\Vert \bm{U}\bm{H}-\bm{U}^{\natural}\right\Vert _{2,\infty}+\kappa^{\natural2}\frac{\zeta_{\mathsf{op}}}{\sigma_{r}^{\natural2}}\left\Vert \bm{U}\bm{H}-\bm{U}^{\natural}\right\Vert _{2,\infty}^{2},
\end{align*}
provided that $\zeta_{\mathsf{op}}\lesssim\sigma_{r}^{\natural2}/\kappa^{\natural2}$
and $n_{1}\gtrsim\kappa^{\natural4}\mu^{\natural}r$. Rearrange terms
to show that
\begin{equation}
\left\Vert \bm{U}\bm{H}-\bm{U}^{\natural}\right\Vert _{2,\infty}\lesssim\frac{\zeta_{\mathsf{op}}}{\sigma_{r}^{\natural2}}\sqrt{\frac{\mu^{\natural}r}{n_{1}}},\label{eq:UH-U-star-2-infty}
\end{equation}
provided that $\zeta_{\mathsf{op}}\ll\sigma_{r}^{\natural2}/\kappa^{\natural2}$.
Then we can use (\ref{eq:UH-U-star-2-infty}) to refine the bound
(\ref{eq:UH-U-star-2-infty-m}) as
\begin{align*}
\left\Vert \bm{U}_{m,\cdot}\bm{H}-\bm{U}_{m,\cdot}^{\natural}\right\Vert _{2} & \lesssim\frac{\zeta_{\mathsf{op},m}}{\sigma_{r}^{\natural2}}\sqrt{\frac{\mu^{\natural}r}{n_{1}}}+\left\Vert \bm{U}_{m,\cdot}^{\natural}\right\Vert _{2}\left(\kappa^{\natural2}\sqrt{\frac{\mu^{\natural}r}{n_{1}}}\frac{\zeta_{\mathsf{op}}}{\sigma_{r}^{\natural2}}+\kappa^{\natural2}\frac{\zeta_{\mathsf{op}}^{2}}{\sigma_{r}^{\natural4}}\right)\\
 & \quad+\left\Vert \bm{U}_{m,\cdot}\bm{H}-\bm{U}_{m,\cdot}^{\natural}\right\Vert _{2}\left(\kappa^{\natural2}\sqrt{\frac{\mu^{\natural}r}{n_{1}}}\frac{\zeta_{\mathsf{op}}}{\sigma_{r}^{\natural2}}+\kappa^{\natural2}\frac{\zeta_{\mathsf{op}}^{2}}{\sigma_{r}^{\natural4}}\right)+\kappa^{\natural2}\frac{\zeta_{\mathsf{op}}}{\sigma_{r}^{\natural2}}\left\Vert \bm{U}_{m,\cdot}\bm{H}-\bm{U}_{m,\cdot}^{\natural}\right\Vert _{2}^{2},
\end{align*}
provided that $\zeta_{\mathsf{op}}\lesssim\sigma_{r}^{\natural2}/\kappa^{\natural2}$.
Again, we can rearrange terms to achieve
\begin{align}
\left\Vert \bm{U}_{m,\cdot}\bm{H}-\bm{U}_{m,\cdot}^{\natural}\right\Vert _{2} & \lesssim\frac{\zeta_{\mathsf{op},m}}{\sigma_{r}^{\natural2}}\sqrt{\frac{\mu^{\natural}r}{n_{1}}}+\left\Vert \bm{U}_{m,\cdot}^{\natural}\right\Vert _{2}\left(\kappa^{\natural2}\sqrt{\frac{\mu^{\natural}r}{n_{1}}}\frac{\zeta_{\mathsf{op}}}{\sigma_{r}^{\natural2}}+\kappa^{\natural2}\frac{\zeta_{\mathsf{op}}^{2}}{\sigma_{r}^{\natural4}}\right),\label{eq:UH-U-natural-m}
\end{align}
provided that $\zeta_{\mathsf{op}}\ll\sigma_{r}^{\natural2}/\kappa^{\natural2}$
and $n_{1}\gtrsim\mu^{\natural}r$. Combine (\ref{eq:UH-U-star-2-infty})
and Lemma \ref{lemma:hpca-basic-facts} gives 
\begin{align*}
\left\Vert \bm{U}\bm{R}_{\bm{U}}-\bm{U}^{\natural}\right\Vert _{2,\infty} & \leq\left\Vert \bm{U}\bm{H}-\bm{U}^{\natural}\right\Vert _{2,\infty}+\left\Vert \bm{U}\left(\bm{H}-\bm{R}_{\bm{U}}\right)\right\Vert _{2,\infty}\lesssim\left\Vert \bm{U}\bm{H}-\bm{U}^{\natural}\right\Vert _{2,\infty}+\left\Vert \bm{U}\right\Vert _{2,\infty}\left\Vert \bm{H}-\bm{R}_{\bm{U}}\right\Vert \\
 & \lesssim\left\Vert \bm{U}\bm{H}-\bm{U}^{\natural}\right\Vert _{2,\infty}+\left(\left\Vert \bm{U}^{\natural}\right\Vert _{2,\infty}+\left\Vert \bm{U}\bm{R}_{\bm{U}}-\bm{U}^{\natural}\right\Vert _{2,\infty}\right)\left\Vert \bm{H}-\bm{R}_{\bm{U}}\right\Vert \\
 & \lesssim\kappa^{\natural2}\frac{\zeta_{\mathsf{op}}}{\sigma_{r}^{\natural2}}\sqrt{\frac{\mu^{\natural}r}{n_{1}}}+\frac{\zeta_{\mathsf{op}}^{2}}{\sigma_{r}^{\natural4}}\sqrt{\frac{\mu^{\natural}r}{n_{1}}}+\frac{\zeta_{\mathsf{op}}^{2}}{\sigma_{r}^{\natural4}}\left\Vert \bm{U}\bm{R}_{\bm{U}}-\bm{U}^{\natural}\right\Vert _{2,\infty}.
\end{align*}
Once again, when $\zeta_{\mathsf{op}}\ll\sigma_{r}^{\natural2}$,
one can rearrange terms to yield 
\[
\left\Vert \bm{U}\bm{R}_{\bm{U}}-\bm{U}^{\natural}\right\Vert _{2,\infty}\lesssim\frac{\zeta_{\mathsf{op}}}{\sigma_{r}^{\natural2}}\sqrt{\frac{\mu^{\natural}r}{n_{1}}}.
\]
Similarly, we can combine (\ref{eq:UH-U-natural-m}) and Lemma \ref{lemma:hpca-basic-facts}
to achieve
\begin{align*}
\left\Vert \bm{U}_{m,\cdot}\bm{R}_{\bm{U}}-\bm{U}_{m,\cdot}^{\natural}\right\Vert _{2} & \lesssim\frac{\zeta_{\mathsf{op},m}}{\sigma_{r}^{\natural2}}\sqrt{\frac{\mu^{\natural}r}{n_{1}}}+\left\Vert \bm{U}_{m,\cdot}^{\natural}\right\Vert _{2}\left(\kappa^{\natural2}\sqrt{\frac{\mu^{\natural}r}{n_{1}}}\frac{\zeta_{\mathsf{op}}}{\sigma_{r}^{\natural2}}+\kappa^{\natural2}\frac{\zeta_{\mathsf{op}}^{2}}{\sigma_{r}^{\natural4}}\right).
\end{align*}

\subsubsection{Proof of Lemma \ref{lemma:hpca-loo-basics}\label{sec:proof-lemma-HeteroPCA-loo-basics}}

We start with $\Vert\bm{G}^{(m)}-\bm{G}\Vert$, for which the triangle
inequality yields 
\begin{equation}
\left\Vert \bm{G}^{(m)}-\bm{G}\right\Vert \leq\left\Vert \mathcal{P}_{\mathsf{off}\text{-}\mathsf{diag}}\left(\bm{G}^{(m)}-\bm{G}\right)\right\Vert +\left\Vert \mathcal{P}_{\mathsf{diag}}\left(\bm{G}^{(m)}-\bm{G}\right)\right\Vert .\label{eq:Gm-G-decompose-LOO}
\end{equation}

\begin{itemize}
\item Regarding the first term on the right-hand side of (\ref{eq:Gm-G-decompose-LOO}),
it is observed that $\mathcal{P}_{\mathsf{off}\text{-}\mathsf{diag}}(\bm{G}^{(m)}-\bm{G})$
in the current paper is the same as the matrix $\bm{G}^{(m)}-\bm{G}$
in \citet{cai2019subspace} (due to the diagonal deletion strategy
employed therein). One can then apply \citet[Lemma 6]{cai2019subspace}
to show that 
\begin{align*}
\left\Vert \mathcal{P}_{\mathsf{off}\text{-}\mathsf{diag}}\left(\bm{G}^{(m)}-\bm{G}\right)\right\Vert  & \lesssim\sigma\sqrt{n_{2}}\left(\sigma\sqrt{n_{1}}+\left\Vert \bm{M}^{\natural\top}\right\Vert _{2,\infty}\right)\sqrt{\log n}\\
 & \lesssim\sigma^{2}\sqrt{n_{1}n_{2}\log n}+\sigma\sqrt{n_{2}}\sqrt{\frac{\mu^{\natural}r}{n_{2}}}\sigma_{1}^{\natural}\sqrt{\log n}\\
 & \asymp\sigma^{2}\sqrt{n_{1}n_{2}\log n}+\sigma\sigma_{1}^{\natural}\sqrt{\mu^{\natural}r\log n}
\end{align*}
with probability exceeding $1-O(n^{-11}$), where the second inequality
relies on (\ref{eq:M-natural-T-2-infty-subspace}). Here, we have
replaced $\sigma_{\mathsf{col}}$ (resp.~$\sigma_{\mathsf{row}}$)
in \citet[Lemma 6]{cai2019subspace} with $\sigma\sqrt{n_{1}}$ (resp.~$\sigma\sqrt{n_{2}}$)
under our setting. 
\item Recalling that the diagonal of $\bm{G}^{(m)}$ coincides with the
true diagonal of $\bm{G}^{\natural}$, we can invoke Lemma \ref{lemma:hpca-estimation}
to bound the second term on the right-hand side of (\ref{eq:Gm-G-decompose-LOO})
as follows 
\[
\left\Vert \mathcal{P}_{\mathsf{diag}}\left(\bm{G}^{(m)}-\bm{G}\right)\right\Vert =\left\Vert \mathcal{P}_{\mathsf{diag}}\left(\bm{G}^{\natural}-\bm{G}\right)\right\Vert \lesssim\kappa^{\natural2}\sqrt{\frac{\mu^{\natural}r}{n_{1}}}\zeta_{\mathsf{op}}.
\]
\end{itemize}
Combining the above two bounds and invoke the union bound lead to
\begin{align*}
\left\Vert \bm{G}^{(m)}-\bm{G}\right\Vert  & \leq\left\Vert \mathcal{P}_{\mathsf{off}\text{-}\mathsf{diag}}\left(\bm{G}^{(m)}-\bm{G}\right)\right\Vert +\left\Vert \mathcal{P}_{\mathsf{diag}}\left(\bm{G}^{(m)}-\bm{G}\right)\right\Vert \\
 & \lesssim\left\{ \sigma^{2}\sqrt{n_{1}n_{2}\log n}+\sigma\sigma_{1}^{\natural}\sqrt{\mu^{\natural}r\log n}\right\} +\kappa^{\natural2}\sqrt{\frac{\mu^{\natural}r}{n_{1}}}\zeta_{\mathsf{op}}\\
 & \lesssim\sigma^{2}\sqrt{n_{1}n_{2}\log n}+\sigma\sigma_{1}^{\natural}\sqrt{\mu^{\natural}r\log n}+\kappa^{\natural2}\sqrt{\frac{\mu^{\natural}r}{n_{1}}}\left(\sigma^{2}\sqrt{n_{1}n_{2}}\log n+\sigma\sigma_{1}^{\natural}\sqrt{n_{1}\log n}\right)\\
 & \lesssim\sigma^{2}\sqrt{n_{1}n_{2}}\log n+\kappa^{\natural2}\sigma\sigma_{1}^{\natural}\sqrt{\mu^{\natural}r\log n}
\end{align*}
simultaneously for all $1\leq m\leq n_{1}$, as long as $n_{1}\gtrsim\kappa^{\natural4}\mu^{\natural}r$.
Here, the penultimate inequality relies on the definition of $\zeta_{\mathsf{op}}$.

To finish up, taking the above inequality and Lemma \ref{lemma:hpca-estimation}
collectively yields 
\[
\left\Vert \bm{G}^{(m)}-\bm{G}^{\natural}\right\Vert \leq\left\Vert \bm{G}^{(m)}-\bm{G}\right\Vert +\left\Vert \bm{G}-\bm{G}^{\natural}\right\Vert \lesssim\sigma^{2}\sqrt{n_{1}n_{2}\log n}+\kappa^{\natural2}\sigma\sigma_{1}^{\natural}\sqrt{\mu^{\natural}r\log n}+\zeta_{\mathsf{op}}\asymp\zeta_{\mathsf{op}},
\]
with the proviso that $n_{1}\gtrsim\kappa^{\natural4}\mu^{\natural}r$.

\subsubsection{Proof of Lemma \ref{lemma:hpca-lto-perturbation}\label{sec:proof-lemma-HeteroPCA-lto-perturbation}}

The proof of Lemma \ref{lemma:hpca-lto-perturbation} can be directly
adapted from the proof of \citet[Lemma 9]{cai2019subspace} (see \citet[Appendix C.5]{cai2019subspace}).
Two observations are crucial: 
\begin{itemize}
\item The off-diagonal part of $\bm{G}^{(m)}$ (resp.~$\bm{G}^{(m,l)}$)
in this paper is the same as that of $\bm{G}^{(m)}$ (resp.~$\bm{G}^{(m,l)}$)
defined in \citet{cai2019subspace}. 
\item Regarding the diagonal, this paper imputes the diagonals of both $\bm{G}^{(m)}$
and $\bm{G}^{(m,l)}$ with the diagonal of the ground truth $\bm{G}^{\natural}$,
while the diagonal entries of both $\bm{G}^{(m)}$ and $\bm{G}^{(m,l)}$
in \cite{cai2019subspace} are all zeros. 
\end{itemize}
Therefore, in both the current paper and \cite{cai2019subspace},
we end up dealing with the same matrix $\bm{G}^{(m)}-\bm{G}^{(m,l)}$,
and consequently, the bound on $\Vert\bm{G}^{(m)}-\bm{G}^{(m,l)}\Vert$
established in \citet[Appendix C.5.1]{cai2019subspace} remains valid
when it comes to our setting.

In addition, we can easily check that the bound on $\Vert(\bm{G}^{(m)}-\bm{G}^{(m,l)})\bm{U}^{(m,l)}\Vert$
derived in \citet[Appendix C.5.2]{cai2019subspace} also holds under
our setting; this is simply because the analysis in \citet[Appendix C.5.2]{cai2019subspace}
remains valid if $\bm{G}$ and $\bm{G}^{(l)}$ have the same deterministic
diagonal. By replacing $\sigma_{\mathsf{col}}$ with $\sigma\sqrt{n_{1}}$
in \citet[Lemma 9]{cai2019subspace}, we arrive at the result claimed
in Lemma \ref{lemma:hpca-lto-perturbation}.

\subsubsection{Proof of Lemma \ref{lemma:hpca-useful-2}\label{subsec:proof-lemma-pca-useful-2}}

For the sake of brevity, we shall only focus on proving (\ref{eq:hpca-useful-2-1}).
The proof of (\ref{eq:hpca-useful-2-2}) is similar to --- and in
fact, simpler than --- the proof of (\ref{eq:hpca-useful-2-1}),
and can also be directly adapted from the proof of \citet[Lemma 7]{cai2019subspace}.

For notational simplicity, we denote $\bm{B}\coloneqq\mathcal{P}_{-m,\cdot}(\bm{M})$.
For any $l\in[n_{2}]$, we can write 
\begin{align*}
 & \left\Vert \bm{e}_{l}^{\top}\left[\mathcal{P}_{-m,\cdot}\left(\bm{M}\right)\right]^{\top}\left(\bm{U}^{(m)}\bm{H}^{(m)}-\bm{U}^{\natural}\right)\right\Vert _{2}=\left\Vert \bm{B}_{\cdot,l}^{\top}\left(\bm{U}^{(m)}\bm{H}^{(m)}-\bm{U}^{\natural}\right)\right\Vert _{2}\\
 & \quad\leq\underbrace{\left\Vert \mathbb{E}\left(\bm{B}_{\cdot,l}\right)^{\top}\left(\bm{U}^{(m)}\bm{H}^{(m)}-\bm{U}^{\natural}\right)\right\Vert _{2}}_{\eqqcolon\alpha_{1}}+\underbrace{\left\Vert \left[\bm{B}_{\cdot,l}-\mathbb{E}\left(\bm{B}_{\cdot,l}\right)\right]^{\top}\left(\bm{U}^{(m)}\bm{H}^{(m)}-\bm{U}^{(m,l)}\bm{H}^{(m,l)}\right)\right\Vert _{2}}_{\eqqcolon\alpha_{2}}\\
 & \quad\quad+\underbrace{\left\Vert \left[\bm{B}_{\cdot,l}-\mathbb{E}\left(\bm{B}_{\cdot,l}\right)\right]^{\top}\left(\bm{U}^{(m,l)}\bm{H}^{(m,l)}-\bm{U}^{\natural}\right)\right\Vert _{2}}_{\eqqcolon\alpha_{3}}.
\end{align*}
Therefore, we seek to bound $\alpha_{1}$, $\alpha_{2}$ and $\alpha_{3}$
separately. 
\begin{itemize}
\item Let us begin with the quantity $\alpha_{1}$. It is straightforward
to see that 
\begin{align}
\alpha_{1} & =\left\Vert \bm{M}_{\cdot,l}^{\natural\top}\left(\bm{U}^{(m)}\bm{H}^{(m)}-\bm{U}^{\natural}\right)\right\Vert _{2}=\left\Vert \bm{M}_{\cdot,l}^{\natural\top}\bm{U}^{\natural}\bm{U}^{\natural\top}\left(\bm{U}^{(m)}\bm{H}^{(m)}-\bm{U}^{\natural}\right)\right\Vert _{2}\nonumber \\
 & \leq\left\Vert \bm{M}^{\natural\top}\right\Vert _{2,\infty}\left\Vert \bm{U}^{\natural\top}\left(\bm{U}^{(m)}\bm{H}^{(m)}-\bm{U}^{\natural}\right)\right\Vert \nonumber \\
 & =\left\Vert \bm{M}^{\natural\top}\right\Vert _{2,\infty}\left\Vert \bm{H}^{(m)\top}\bm{H}^{(m)}-\bm{I}_{r}\right\Vert ,\label{eq:alpha1-bound-1256}
\end{align}
where the last identity arises from the following relation 
\[
\left\Vert \bm{U}^{\natural\top}\left(\bm{U}^{(m)}\bm{H}^{(m)}-\bm{U}^{\natural}\right)\right\Vert =\left\Vert \bm{U}^{\natural\top}\bm{U}^{(m)}\bm{U}^{(m)\top}\bm{U}^{\natural}-\bm{I}_{r}\right\Vert =\left\Vert \bm{H}^{(m)\top}\bm{H}^{(m)}-\bm{I}_{r}\right\Vert .
\]
Let us write the SVD of $\bm{H}^{(m)}=\bm{U}^{(m)\top}\bm{U}^{\natural}$
as $\bm{X}(\cos\bm{\Theta})\bm{Y}^{\top}$, where $\bm{X},\bm{Y}\in\mathbb{R}^{r\times r}$
are square orthonormal matrices and $\bm{\Theta}$ is a diagonal matrix
composed of the principal angles between $\bm{U}^{(m)}$ and $\bm{U}^{\natural}$.
This allows us to deduce that 
\begin{align*}
\left\Vert \bm{H}^{(m)\top}\bm{H}^{(m)}-\bm{I}_{r}\right\Vert  & =\left\Vert \bm{Y}\left(\bm{I}_{r}-\cos^{2}\bm{\Theta}\right)\bm{Y}^{\top}\right\Vert =\left\Vert \bm{I}_{r}-\cos^{2}\bm{\Theta}\right\Vert =\left\Vert \sin^{2}\bm{\Theta}\right\Vert =\left\Vert \sin\bm{\Theta}\right\Vert ^{2}.
\end{align*}
In view of Davis-Kahan's $\sin\bm{\Theta}$ Theorem \citep[Theorem 2.2.1]{chen2020spectral},
we have 
\[
\left\Vert \sin\bm{\Theta}\right\Vert \leq\frac{\left\Vert \bm{G}^{(m)}-\bm{G}^{\natural}\right\Vert }{\lambda_{r}\left(\bm{G}^{(m)}\right)-\lambda_{r+1}\left(\bm{G}^{\natural}\right)}\lesssim\frac{\zeta_{\mathsf{op}}}{\sigma_{r}^{\natural2}},
\]
where the last inequality follows from Lemma \ref{lemma:hpca-loo-basics}
and an application of Weyl's inequality: 
\begin{align*}
\lambda_{r}\left(\bm{G}^{(m)}\right) & \geq\lambda_{r}\left(\bm{G}^{\natural}\right)-\left\Vert \bm{G}^{(m)}-\bm{G}^{\natural}\right\Vert \overset{\text{(i)}}{\geq}\sigma_{r}^{\natural2}-\widetilde{C}\zeta_{\mathsf{op}}\overset{\text{(ii)}}{\geq}\frac{1}{2}\sigma_{r}^{\natural2},
\end{align*}
with $\widetilde{C}>0$ representing some absolute constant. Here,
(i) results from Lemma \ref{lemma:hpca-loo-basics}, while (ii) holds
provided that $\zeta_{\mathsf{op}}\ll\sigma_{r}^{\natural2}$.. Therefore,
we arrive at 
\[
\left\Vert \bm{H}^{(m)\top}\bm{H}^{(m)}-\bm{I}_{r}\right\Vert \leq\left\Vert \sin\bm{\Theta}\right\Vert ^{2}\lesssim\frac{\zeta_{\mathsf{op}}^{2}}{\sigma_{r}^{\natural4}}.
\]
Substitution into (\ref{eq:alpha1-bound-1256}) yields 
\[
\alpha_{1}\lesssim\frac{\zeta_{\mathsf{op}}^{2}}{\sigma_{r}^{\natural4}}\left\Vert \bm{M}^{\natural\top}\right\Vert _{2,\infty}.
\]
\item Regarding $\alpha_{2}$, it is observed that 
\begin{align*}
\alpha_{2} & \leq\left(\left\Vert \bm{B}_{\cdot,l}\right\Vert _{2}+\left\Vert \mathbb{E}\left(\bm{B}_{\cdot,l}\right)\right\Vert _{2}\right)\left\Vert \bm{U}^{(m)}\bm{H}^{(m)}-\bm{U}^{(m,l)}\bm{H}^{(m,l)}\right\Vert \\
 & \leq\left(\left\Vert \bm{M}_{\cdot,l}\right\Vert _{2}+\left\Vert \bm{M}_{\cdot,l}^{\natural}\right\Vert _{2}\right)\left\Vert \left(\bm{U}^{(m)}\bm{U}^{(m)\top}-\bm{U}^{(m,l)}\bm{U}^{(m,l)\top}\right)\bm{U}^{\natural}\right\Vert \\
 & \leq\left(\left\Vert \bm{M}^{\natural\top}\right\Vert _{2,\infty}+B\sqrt{\log n}+\sigma\sqrt{n_{1}}\right)\left\Vert \bm{U}^{(m)}\bm{U}^{(m)\top}-\bm{U}^{(m,l)}\bm{U}^{(m,l)\top}\right\Vert .
\end{align*}
Here, the last inequality arises from \citet[Lemma 12]{cai2019subspace}. 
\item We are now left with the quantity $\alpha_{3}$, which can be expressed
as 
\[
\alpha_{3}=\left\Vert \sum_{i\in[n_{1}]\setminus\{m\}}E_{i,l}\left(\bm{U}^{(m,l)}\bm{H}^{(m,l)}-\bm{U}^{\natural}\right)_{i,\cdot}\right\Vert _{2}.
\]
Conditional on $\bm{U}^{(m,l)}\bm{H}^{(m,l)}-\bm{U}^{\natural}$,
this term can be viewed as the spectral norm of a sum of independent
mean-zero random vectors (where the randomness comes from $\{E_{i,l}\}_{i\in[n_{1}]\setminus\{m\}}$).
To control this term, we first calculate 
\begin{align*}
L & \coloneqq\max_{i\in[n_{1}]\setminus\{m\}}\left\Vert E_{i,l}\left(\bm{U}^{(m,l)}\bm{H}^{(m,l)}-\bm{U}^{\natural}\right)_{i,\cdot}\right\Vert \leq B\left\Vert \bm{U}^{(m,l)}\bm{H}^{(m,l)}-\bm{U}^{\natural}\right\Vert _{2,\infty},\\
V & \coloneqq\sum_{i\in[n_{1}]\setminus\{m\}}\mathbb{E}\left[E_{i,l}^{2}\right]\left\Vert \left(\bm{U}^{(m,l)}\bm{H}^{(m,l)}-\bm{U}^{\natural}\right)_{i,\cdot}\right\Vert _{2}^{2}\leq\sigma^{2}n_{1}\left\Vert \bm{U}^{(m,l)}\bm{H}^{(m,l)}-\bm{U}^{\natural}\right\Vert _{2,\infty}^{2}.
\end{align*}
In view of the matrix Bernstein inequality \citep[Theorem 6.1.1]{Tropp:2015:IMC:2802188.2802189},
\begin{align*}
\alpha_{3} & =\left\Vert \sum_{i\in[n_{1}]\setminus\{m\}}E_{i,l}\left(\bm{U}^{(m,l)}\bm{H}^{(m,l)}-\bm{U}^{\natural}\right)_{i,\cdot}\right\Vert \lesssim\sqrt{V\log n}+L\log n\\
 & \lesssim\left(\sigma\sqrt{n_{1}\log n}+B\log n\right)\left\Vert \bm{U}^{(m,l)}\bm{H}^{(m,l)}-\bm{U}^{\natural}\right\Vert _{2,\infty}
\end{align*}
with probability exceeding $1-O(n^{-10})$ . In addition, the triangle
inequality gives 
\begin{align*}
\left\Vert \bm{U}^{(m,l)}\bm{H}^{(m,l)}-\bm{U}^{\natural}\right\Vert _{2,\infty} & \leq\left\Vert \bm{U}^{(m)}\bm{H}^{(m)}-\bm{U}^{\natural}\right\Vert _{2,\infty}+\left\Vert \bm{U}^{(m,l)}\bm{H}^{(m,l)}-\bm{U}^{(m)}\bm{H}^{(m)}\right\Vert _{2,\infty}\\
 & \leq\left\Vert \bm{U}^{(m)}\bm{H}^{(m)}-\bm{U}^{\natural}\right\Vert _{2,\infty}+\left\Vert \bm{U}^{(m,l)}\bm{H}^{(m,l)}-\bm{U}^{(m)}\bm{H}^{(m)}\right\Vert \\
 & \leq\left\Vert \bm{U}^{(m)}\bm{H}^{(m)}-\bm{U}^{\natural}\right\Vert _{2,\infty}+\left\Vert \bm{U}^{(m)}\bm{U}^{(m)\top}-\bm{U}^{(m,l)}\bm{U}^{(m,l)\top}\right\Vert ,
\end{align*}
and therefore, 
\[
\alpha_{3}\lesssim\left(\sigma\sqrt{n_{1}\log n}+B\log n\right)\left(\left\Vert \bm{U}^{(m)}\bm{H}^{(m)}-\bm{U}^{\natural}\right\Vert _{2,\infty}+\left\Vert \bm{U}^{(m)}\bm{U}^{(m)\top}-\bm{U}^{(m,l)}\bm{U}^{(m,l)\top}\right\Vert \right).
\]
\end{itemize}
Combine the preceding bounds on $\alpha_{1}$, $\alpha_{2}$ and $\alpha_{3}$
to arrive at 
\begin{align*}
 & \left\Vert \bm{e}_{l}^{\top}\left[\mathcal{P}_{-m,\cdot}\left(\bm{M}\right)\right]^{\top}\left(\bm{U}^{(m)}\bm{H}^{(m)}-\bm{U}^{\natural}\right)\right\Vert _{2}\leq\alpha_{1}+\alpha_{2}+\alpha_{3}\\
 & \quad\lesssim\frac{\zeta_{\mathsf{op}}^{2}}{\sigma_{r}^{\natural4}}\left\Vert \bm{M}^{\natural\top}\right\Vert _{2,\infty}+\left(\sigma\sqrt{n_{1}\log n}+B\log n\right)\left\Vert \bm{U}^{(m)}\bm{H}^{(m)}-\bm{U}^{\natural}\right\Vert _{2,\infty}\\
 & \quad\quad+\left(\left\Vert \bm{M}^{\natural\top}\right\Vert _{2,\infty}+B\log n+\sigma\sqrt{n_{1}\log n}\right)\left\Vert \bm{U}^{(m)}\bm{U}^{(m)\top}-\bm{U}^{(m,l)}\bm{U}^{(m,l)\top}\right\Vert 
\end{align*}
as claimed.

\subsubsection{Proof of Lemma \ref{lemma:hpca-useful-1}\label{subsec:proof-lemma-pca-useful-1}}

In this subsection, for the sake of brevity, we shall only focus on
proving (\ref{eq:hpca-useful-1-1}). The proof of (\ref{eq:hpca-useful-1-2})
is similar to that of (\ref{eq:hpca-useful-1-1}), and can also be
easily adapted from the proof of \citet[Lemma 7]{cai2019subspace}.

For notational simplicity, we denote $\bm{B}\coloneqq\mathcal{P}_{-m,\cdot}(\bm{M})$,
allowing us to express 
\[
\bm{E}_{m,\cdot}\left[\mathcal{P}_{-m,\cdot}\left(\bm{M}\right)\right]^{\top}\left(\bm{U}^{(m)}\bm{H}^{(m)}-\bm{U}^{\natural}\right)=\sum_{j=1}^{n_{2}}E_{m,j}\left[\bm{B}^{\top}\left(\bm{U}^{(m)}\bm{H}^{(m)}-\bm{U}^{\natural}\right)\right]_{j,\cdot}.
\]
Conditional on $\bm{B}$ and $\bm{U}^{(m)}\bm{H}^{(m)}-\bm{U}^{\natural}$,
the above term can be viewed as a sum of independent zero-mean random
vectors, where the randomness comes from $\{E_{m,j}\}_{j\in[n_{2}]}$.
We can calculate 
\begin{align*}
L & \coloneqq\max_{j\in[n_{2}]}\left\Vert E_{m,j}\left[\bm{B}^{\top}\left(\bm{U}^{(m)}\bm{H}^{(m)}-\bm{U}^{\natural}\right)\right]_{j,\cdot}\right\Vert \leq B_{m}\left\Vert \bm{B}^{\top}\left(\bm{U}^{(m)}\bm{H}^{(m)}-\bm{U}^{\natural}\right)\right\Vert _{2,\infty},\\
V & \coloneqq\sum_{j\in[n_{2}]}\mathbb{E}\left(E_{m,j}^{2}\right)\left\Vert \bm{B}^{\top}\left(\bm{U}^{(m)}\bm{H}^{(m)}-\bm{U}^{\natural}\right)_{j,\cdot}\right\Vert _{2}^{2}\leq\sigma_{m}^{2}\left\Vert \bm{B}^{\top}\left(\bm{U}^{(m)}\bm{H}^{(m)}-\bm{U}^{\natural}\right)\right\Vert _{\mathrm{F}}^{2}.
\end{align*}
In view of the matrix Bernstein inequality \citep[Theorem 6.1.1]{Tropp:2015:IMC:2802188.2802189},
with probability exceeding $1-O(n^{-11})$ 
\begin{align*}
 & \left\Vert \sum_{j=1}^{n_{2}}E_{m,j}\left[\bm{B}^{\top}\left(\bm{U}^{(m)}\bm{H}^{(m)}-\bm{U}^{\natural}\right)\right]_{j,\cdot}\right\Vert _{2}\lesssim\sqrt{V\log n}+L\log n\\
 & \quad\lesssim\underbrace{\sigma_{m}\left\Vert \bm{B}^{\top}\left(\bm{U}^{(m)}\bm{H}^{(m)}-\bm{U}^{\natural}\right)\right\Vert _{\mathrm{F}}\sqrt{\log n}}_{\eqqcolon\alpha_{1}}+\underbrace{B_{m}\left\Vert \bm{B}^{\top}\left(\bm{U}^{(m)}\bm{H}^{(m)}-\bm{U}^{\natural}\right)\right\Vert _{2,\infty}\log n}_{\eqqcolon\alpha_{2}}.
\end{align*}

For the first term $\alpha_{1}$, with probability exceeding $1-O(n^{-11})$
we have 
\begin{align*}
\left\Vert \bm{B}^{\top}\left(\bm{U}^{(m)}\bm{H}^{(m)}-\bm{U}^{\natural}\right)\right\Vert _{\mathrm{F}} & \leq\left\Vert \bm{B}\right\Vert \left\Vert \bm{U}^{(m)}\bm{H}^{(m)}-\bm{U}^{\natural}\right\Vert _{\mathrm{F}}\leq\left\Vert \bm{M}\right\Vert \left\Vert \bm{U}^{(m)}\bm{H}^{(m)}-\bm{U}^{\natural}\right\Vert _{\mathrm{F}}\\
 & \leq\left(\left\Vert \bm{M}^{\natural}\right\Vert +\left\Vert \bm{E}\right\Vert \right)\left\Vert \bm{U}^{(m)}\bm{H}^{(m)}-\bm{U}^{\natural}\right\Vert _{\mathrm{F}}\\
 & \lesssim\left(\sigma_{1}^{\natural}+\sigma\sqrt{n}\right)\left\Vert \bm{U}^{(m)}\bm{H}^{(m)}-\bm{U}^{\natural}\right\Vert _{\mathrm{F}}\\
 & \asymp\left(\sigma_{1}^{\natural}+\sigma\sqrt{n_{2}}\right)\left\Vert \bm{U}^{(m)}\bm{H}^{(m)}-\bm{U}^{\natural}\right\Vert _{\mathrm{F}}.
\end{align*}
Here, the penultimate inequality uses $\Vert\bm{E}\Vert\lesssim\sigma\sqrt{n}$
with probability exceeding $1-O(n^{-11})$, which follows from standard
matrix tail bounds (e.g., \citet[Theorem 3.1.4]{chen2020spectral});
and the last relation holds since $n=\max\{n_{1},n_{2}\}$ and $\sigma\sqrt{n_{1}}\lesssim\zeta_{\mathsf{op}}/\sigma_{1}^{\natural}\ll\sigma_{1}^{\natural}$.
As a result, we reach
\begin{align*}
\alpha_{1} & \lesssim\sigma_{m}\left(\sigma_{1}^{\natural}+\sigma\sqrt{n_{2}}\right)\left\Vert \bm{U}^{(m)}\bm{H}^{(m)}-\bm{U}^{\natural}\right\Vert _{\mathrm{F}}\sqrt{\log n}\\
 & \lesssim\left(\sigma_{m}\sigma_{1}^{\natural}\sqrt{n_{1}\log n}+\sigma\sigma_{m}\sqrt{n_{1}n_{2}\log n}\right)\left\Vert \bm{U}^{(m)}\bm{H}^{(m)}-\bm{U}^{\natural}\right\Vert _{2,\infty}.
\end{align*}

Regarding the second term $\alpha_{2}$, we know from (\ref{eq:hpca-useful-2-1})
in Lemma \ref{lemma:hpca-useful-2} that with probability exceeding
$1-O(n^{-10})$, 
\begin{align*}
\alpha_{2} & \leq\underbrace{\left(B_{m}\log n\right)\frac{\zeta_{\mathsf{op}}^{2}}{\sigma_{r}^{\natural4}}\left\Vert \bm{M}^{\natural\top}\right\Vert _{2,\infty}}_{\eqqcolon\beta_{1}}+\underbrace{\left(B_{m}\log n\right)\left(\sigma\sqrt{n_{1}\log n}+B\log n\right)\left\Vert \bm{U}^{(m)}\bm{H}^{(m)}-\bm{U}^{\natural}\right\Vert _{2,\infty}}_{\eqqcolon\beta_{2}}\\
 & \quad+\underbrace{\left(B_{m}\log n\right)\left(\left\Vert \bm{M}^{\natural\top}\right\Vert _{2,\infty}+B\log n+\sigma\sqrt{n_{1}\log n}\right)\left\Vert \bm{U}^{(m)}\bm{U}^{(m)\top}-\bm{U}^{(m,l)}\bm{U}^{(m,l)\top}\right\Vert }_{\eqqcolon\beta_{3}}
\end{align*}
holds for all $m\in[n_{1}]$. In what follows, we shall bound $\beta_{1}$,
$\beta_{2}$ and $\beta_{3}$ respectively. 
\begin{itemize}
\item Regarding $\beta_{1}$, we first observe that 
\begin{align}
\left(B_{m}\log n\right)\left\Vert \bm{M}^{\natural\top}\right\Vert _{2,\infty} & \leq\left(B_{m}\log n\right)\sqrt{\frac{\mu^{\natural}r}{n_{2}}}\sigma_{1}^{\natural}\lesssim\sigma_{m}\sqrt{n_{2}\log n}\sqrt{\frac{\mu^{\natural}r}{n_{2}}}\sigma_{1}^{\natural}\lesssim\sigma_{m}\sigma_{1}^{\natural}\sqrt{\mu^{\natural}r\log n},\label{eq:pca-useful-1-inter-1}
\end{align}
where we have used the noise condition $B_{m}\lesssim\sigma_{m}\sqrt{n_{2}/\log n}$.
Therefore 
\[
\beta_{1}\lesssim\sigma_{m}\sigma_{1}^{\natural}\sqrt{n_{1}\log n}\frac{\zeta_{\mathsf{op}}^{2}}{\sigma_{r}^{\natural4}}\sqrt{\frac{\mu^{\natural}r}{n_{1}}}.
\]
\item When it comes to $\beta_{2}$, we notice that 
\begin{align}
\left(B_{m}\log n\right)\left(\sigma\sqrt{n_{1}\log n}+B\log n\right) & \lesssim\sigma_{m}\sigma\sqrt{n_{1}n_{2}}\log n.\label{eq:pca-useful-1-inter-3}
\end{align}
Here, we have used the noise assumption $B\lesssim\sigma\sqrt[4]{n_{1}n_{2}}/\sqrt{\log n}$
and $B_{m}\lesssim\sigma_{m}\min\{\sqrt[4]{n_{1}n_{2}},\sqrt{n_{2}}\}/\sqrt{\log n}$.
This in turn leads us to 
\[
\beta_{2}\leq\sigma_{m}\sigma\sqrt{n_{1}n_{2}}\log n\left\Vert \bm{U}^{(m)}\bm{H}^{(m)}-\bm{U}^{\natural}\right\Vert _{2,\infty}.
\]
\item We are left with the term $\beta_{3}$. From Lemma \ref{lemma:hpca-lto-perturbation},
we see that with probability exceeding $1-O(n^{-10})$, 
\begin{align*}
 & \left(B_{m}\log n\right)\left\Vert \bm{U}^{(m)}\bm{U}^{(m)\top}-\bm{U}^{(m,l)}\bm{U}^{(m,l)\top}\right\Vert \\
 & \quad\lesssim\frac{B_{m}\log n}{\sigma_{r}^{\natural2}}\left[\left(B\log n+\sigma\sqrt{n_{1}\log n}\right)^{2}\left\Vert \bm{U}^{(m)}\bm{H}^{(m)}\right\Vert _{2,\infty}+\sigma^{2}\right]+\frac{B_{m}\log n}{\sigma_{r}^{\natural2}}\left(B\log n+\sigma\sqrt{n_{1}\log n}\right)\left\Vert \bm{M}^{\natural\top}\right\Vert _{2,\infty}\\
 & \quad\lesssim\frac{B_{m}\log n}{\sigma_{r}^{\natural2}}\left[\zeta_{\mathsf{op}}\left\Vert \bm{U}^{(m)}\bm{H}^{(m)}\right\Vert _{2,\infty}+\sigma^{2}\right]+\left(B\log n+\sigma\sqrt{n_{1}\log n}\right)\frac{\sigma_{m}\sigma_{1}^{\natural}}{\sigma_{r}^{\natural2}}\sqrt{\mu^{\natural}r\log n}
\end{align*}
simultaneously for all $m$ and $l$, where the penultimate inequality
follows from (\ref{eq:pca-useful-1-inter-1}) and 
\begin{equation}
\left(\sigma\sqrt{n_{1}\log n}+B\log n\right)^{2}\asymp\sigma^{2}n_{1}\log n+B^{2}\log^{2}d\lesssim\sigma^{2}n_{1}\log n+\sigma^{2}\sqrt{n_{1}n_{2}}\log n\lesssim\zeta_{\mathsf{op}}.\label{eq:pca-useful-1-inter-2}
\end{equation}
Therefore, with probability exceeding $1-O(n^{-10})$, we have 
\begin{align*}
\beta_{3} & \lesssim\left(\left\Vert \bm{M}^{\natural\top}\right\Vert _{2,\infty}+B\log n+\sigma\sqrt{n_{1}\log n}\right)\frac{B_{m}\log n}{\sigma_{r}^{\natural2}}\left[\zeta_{\mathsf{op}}\left\Vert \bm{U}^{(m)}\bm{H}^{(m)}\right\Vert _{2,\infty}+\sigma^{2}\right]\\
 & \quad+\left(\left\Vert \bm{M}^{\natural\top}\right\Vert _{2,\infty}+B\log n+\sigma\sqrt{n_{1}\log n}\right)\left(B\log n+\sigma\sqrt{n_{1}\log n}\right)\frac{\sigma_{m}\sigma_{1}^{\natural}}{\sigma_{r}^{\natural2}}\sqrt{\mu^{\natural}r\log n}\\
 & \lesssim\left(\left\Vert \bm{M}^{\natural\top}\right\Vert _{2,\infty}\frac{B_{m}\log n}{\sigma_{r}^{\natural2}}+\frac{B_{m}\log n\left(B\log n+\sigma\sqrt{n_{1}\log n}\right)}{\sigma_{r}^{\natural2}}\right)\left[\zeta_{\mathsf{op}}\left\Vert \bm{U}^{(m)}\bm{H}^{(m)}\right\Vert _{2,\infty}+\sigma^{2}\right]\\
 & \quad+\left\Vert \bm{M}^{\natural\top}\right\Vert _{2,\infty}\left(B\log n+\sigma\sqrt{n_{1}\log n}\right)\frac{\sigma_{m}\sigma_{1}^{\natural}}{\sigma_{r}^{\natural2}}\sqrt{\mu^{\natural}r\log n}+\left(B\log n+\sigma\sqrt{n_{1}\log n}\right)^{2}\frac{\sigma_{m}\sigma_{1}^{\natural}}{\sigma_{r}^{\natural2}}\sqrt{\mu^{\natural}r\log n}\\
 & \overset{\text{(i)}}{\lesssim}\left(\sigma_{m}\sigma_{1}^{\natural}\sqrt{\mu^{\natural}r\log n}+\sigma_{m}\sigma\sqrt{n_{1}n_{2}}\log n\right)\left(\frac{\zeta_{\mathsf{op}}}{\sigma_{r}^{\natural2}}\left\Vert \bm{U}^{(m)}\bm{H}^{(m)}\right\Vert _{2,\infty}+\frac{\sigma^{2}}{\sigma_{r}^{\natural2}}\right)\\
 & \quad+\left\Vert \bm{M}^{\natural\top}\right\Vert _{2,\infty}\left(B\log n+\sigma\sqrt{n_{1}\log n}\right)\frac{\sigma_{m}\sigma_{1}^{\natural}}{\sigma_{r}^{\natural2}}\sqrt{\mu^{\natural}r\log n}+\frac{\zeta_{\mathsf{op}}}{\sigma_{r}^{\natural2}}\sigma_{m}\sigma_{1}^{\natural}\sqrt{n_{1}\log n}\sqrt{\frac{\mu^{\natural}r}{n_{1}}}\\
 & \overset{\text{(ii)}}{\lesssim}\left(\sigma_{m}\sigma_{1}^{\natural}\sqrt{\mu^{\natural}r\log n}+\sigma_{m}\sigma\sqrt{n_{1}n_{2}}\log n\right)\left(\frac{\zeta_{\mathsf{op}}}{\sigma_{r}^{\natural2}}\left\Vert \bm{U}^{(m)}\bm{H}^{(m)}\right\Vert _{2,\infty}+\frac{\sigma^{2}}{\sigma_{r}^{\natural2}}\right)\\
 & \quad+\left(\frac{\mu^{\natural}r}{n_{1}}+\frac{\mu^{\natural}r}{\sqrt{n_{1}n_{2}}}\right)\frac{\sigma\sigma_{1}^{\natural}\sqrt{n_{1}\log n}\sigma_{m}\sigma_{1}^{\natural}\sqrt{n_{1}\log n}}{\sigma_{r}^{\natural2}}+\frac{\zeta_{\mathsf{op}}}{\sigma_{r}^{\natural2}}\sigma_{m}\sigma_{1}^{\natural}\sqrt{n_{1}\log n}\sqrt{\frac{\mu^{\natural}r}{n_{1}}}\\
 & \overset{\text{(iii)}}{\lesssim}\left(\sigma_{m}\sigma_{1}^{\natural}\sqrt{\mu^{\natural}r\log n}+\sigma_{m}\sigma\sqrt{n_{1}n_{2}}\log n\right)\frac{\zeta_{\mathsf{op}}}{\sigma_{r}^{\natural2}}\left\Vert \bm{U}^{(m)}\bm{H}^{(m)}\right\Vert _{2,\infty}+\frac{\zeta_{\mathsf{op}}}{\sigma_{r}^{\natural2}}\sigma_{m}\sigma_{1}^{\natural}\sqrt{n_{1}\log n}\sqrt{\frac{\mu^{\natural}r}{n_{1}}}
\end{align*}
simultaneously for all $m\in[n_{1}]$, where (i) uses (\ref{eq:pca-useful-1-inter-1}),
(\ref{eq:pca-useful-1-inter-3}) and (\ref{eq:pca-useful-1-inter-2});
(ii) comes from 
\begin{align*}
\left\Vert \bm{M}^{\natural\top}\right\Vert _{2,\infty}\left(B\log n+\sigma\sqrt{n_{1}\log n}\right)\frac{\sigma_{m}\sigma_{1}^{\natural}}{\sigma_{r}^{\natural2}}\sqrt{\mu^{\natural}r\log n} & \lesssim\sqrt{\frac{\mu^{\natural}r}{n_{2}}}\sigma_{1}^{\natural}\sigma\sqrt{n\log n}\frac{\sigma_{m}\sigma_{1}^{\natural}}{\sigma_{r}^{\natural2}}\sqrt{\mu^{\natural}r\log n}\\
 & \lesssim\left(\frac{\mu^{\natural}r}{n_{1}}+\frac{\mu^{\natural}r}{\sqrt{n_{1}n_{2}}}\right)\frac{\sigma\sigma_{1}^{\natural}\sqrt{n_{1}\log n}\sigma_{m}\sigma_{1}^{\natural}\sqrt{n_{1}\log n}}{\sigma_{r}^{\natural2}}
\end{align*}
where we use $B\lesssim\sigma\sqrt{n_{2}/\log n}$; and (iii) holds
provided that $\zeta_{\mathsf{op}}\ll\sigma_{r}^{\natural2}$, $n_{1}\gtrsim\mu^{\natural}r$
and $n_{2}\gtrsim\mu^{\natural}r$. 
\end{itemize}
Putting the above pieces together, we conclude that 
\begin{align*}
 & \left\Vert \bm{E}_{m,\cdot}\left[\mathcal{P}_{-m,\cdot}\left(\bm{M}\right)\right]^{\top}\left(\bm{U}^{(m)}\bm{H}^{(m)}-\bm{U}^{\natural}\right)\right\Vert _{2}\leq\alpha_{1}+\alpha_{2}\leq\alpha_{1}+\beta_{1}+\beta_{2}+\beta_{3}\\
 & \quad\lesssim\left(\sigma_{m}\sigma_{1}^{\natural}\sqrt{n_{1}\log n}+\sigma_{m}\sigma\sqrt{n_{1}n_{2}}\log n\right)\left(\left\Vert \bm{U}^{(m)}\bm{H}^{(m)}-\bm{U}^{\natural}\right\Vert _{2,\infty}+\frac{\zeta_{\mathsf{op}}}{\sigma_{r}^{\natural2}}\left\Vert \bm{U}^{(m)}\bm{H}^{(m)}\right\Vert _{2,\infty}\right)\\
 & \quad+\frac{\zeta_{\mathsf{op}}}{\sigma_{r}^{\natural2}}\sigma_{m}\sigma_{1}^{\natural}\sqrt{n_{1}\log n}\sqrt{\frac{\mu^{\natural}r}{n_{1}}}
\end{align*}
simultaneously for all $m\in[n_{1}]$, provided that $\zeta_{\mathsf{op}}\lesssim\sigma_{r}^{\natural2}$.

\subsubsection{Proof of Lemma \ref{lemma:hpca-loo-perturbation}\label{sec:proof-lemma-HeteroPCA-loo-perturbation}}

In view of the Davis-Kahan $\sin\bm{\Theta}$ Theorem \citep[Theorem 2.2.1]{chen2020spectral},
we have 
\begin{align*}
\left\Vert \bm{U}^{(m)}\bm{U}^{(m)\top}-\bm{U}\bm{U}^{\top}\right\Vert  & \leq\frac{\left\Vert \left(\bm{G}-\bm{G}^{(m)}\right)\bm{U}^{(m)}\right\Vert }{\lambda_{r}\left(\bm{G}^{(m)}\right)-\lambda_{r+1}\left(\bm{G}\right)}\leq\frac{2\left\Vert \left(\bm{G}-\bm{G}^{(m)}\right)\bm{U}^{(m)}\right\Vert }{\sigma_{r}^{\natural2}}\\
 & \leq\underbrace{\frac{2\left\Vert \mathcal{P}_{\mathsf{diag}}\left(\bm{G}-\bm{G}^{(m)}\right)\bm{U}^{(m)}\right\Vert }{\sigma_{r}^{\natural2}}}_{\eqqcolon\alpha_{1}}+\underbrace{\frac{2\left\Vert \mathcal{P}_{\mathsf{off}\text{-}\mathsf{diag}}\left(\bm{G}-\bm{G}^{(m)}\right)\bm{U}^{(m)}\right\Vert }{\sigma_{r}^{\natural2}}}_{\eqqcolon\alpha_{2}},
\end{align*}
where the penultimate inequality follows from Weyl's inequality: 
\begin{align*}
\lambda_{r}\left(\bm{G}^{(m)}\right) & \geq\lambda_{r}\left(\bm{G}^{\natural}\right)-\left\Vert \bm{G}^{(m)}-\bm{G}^{\natural}\right\Vert \overset{\text{(i)}}{\geq}\sigma_{r}^{\natural2}-\widetilde{C}\zeta_{\mathsf{op}}\overset{\text{(ii)}}{\geq}\frac{3}{4}\sigma_{r}^{\natural2},\\
\lambda_{r+1}\left(\bm{G}\right) & \leq\lambda_{r+1}\left(\bm{G}^{\natural}\right)+\left\Vert \bm{G}-\bm{G}^{\natural}\right\Vert \overset{\text{(iii)}}{\leq}\widetilde{C}\zeta_{\mathsf{op}}\overset{\text{(iv)}}{\leq}\frac{1}{4}\sigma_{r}^{\natural2},
\end{align*}
with $\widetilde{C}>0$ some absolute constant. Here, (i) comes from
Lemma \ref{lemma:hpca-loo-basics}; (iii) comes from Lemma \ref{lemma:hpca-estimation}
and the fact that $\lambda_{r+1}(\bm{G}^{\natural})=0$; (ii) and
(iv) are valid provided that $\zeta_{\mathsf{op}}\ll\sigma_{r}^{\natural2}$.

\paragraph{Bounding $\alpha_{1}$.}

Recalling that $\mathcal{P}_{\mathsf{diag}}(\bm{G}^{(m)})=\mathcal{P}_{\mathsf{diag}}(\bm{G}^{\natural})$,
we obtain 
\[
\alpha_{1}=\frac{2\left\Vert \mathcal{P}_{\mathsf{diag}}\left(\bm{G}-\bm{G}^{\natural}\right)\right\Vert \left\Vert \bm{U}^{(m)}\right\Vert }{\sigma_{r}^{\natural2}}\lesssim\kappa^{\natural2}\sqrt{\frac{\mu^{\natural}r}{n_{1}}}\frac{\zeta_{\mathsf{op}}}{\sigma_{r}^{\natural2}},
\]
where the last inequality follows from Lemma \ref{lemma:hpca-estimation}.

\paragraph{Bounding $\alpha_{2}$.}

Observe that the symmetric matrix $\mathcal{P}_{\mathsf{off}\text{-}\mathsf{diag}}(\bm{G}-\bm{G}^{(m)})$
is supported on the $m$-th row and the $m$-th column, and 
\begin{equation}
\left(\bm{G}-\bm{G}^{(m)}\right)_{m,i}=\left(\bm{G}-\bm{G}^{(m)}\right)_{i,m}=\bm{E}_{m,\cdot}\bm{M}_{i,\cdot}^{\top},\qquad\forall i\in[n_{1}].\label{eq:G-Gm-entries}
\end{equation}
Therefore, we can derive 
\begin{align*}
\left\Vert \mathcal{P}_{\mathsf{off}\text{-}\mathsf{diag}}\left(\bm{G}-\bm{G}^{(m)}\right)\bm{U}^{(m)}\right\Vert  & \leq\left\Vert \mathcal{P}_{\mathsf{off}\text{-}\mathsf{diag}}\left(\bm{G}-\bm{G}^{(m)}\right)\bm{U}^{(m)}\right\Vert _{\mathrm{F}}\overset{\text{(i)}}{\leq}2\left\Vert \mathcal{P}_{\mathsf{off}\text{-}\mathsf{diag}}\left(\bm{G}-\bm{G}^{(m)}\right)\bm{U}^{(m)}\bm{H}^{(m)}\right\Vert _{\mathrm{F}}\\
 & \leq2\left\Vert \mathcal{P}_{m,\cdot}\left(\bm{G}-\bm{G}^{(m)}\right)\bm{U}^{(m)}\bm{H}^{(m)}\right\Vert _{\mathrm{F}}+2\left\Vert \mathcal{P}_{\cdot,m}\left(\bm{G}-\bm{G}^{(m)}\right)\bm{U}^{(m)}\bm{H}^{(m)}\right\Vert _{\mathrm{F}}\\
 & \overset{\text{(ii)}}{=}2\underbrace{\left\Vert \bm{E}_{m,\cdot}\left[\mathcal{P}_{-m,\cdot}\left(\bm{M}\right)\right]^{\top}\bm{U}^{(m)}\bm{H}^{(m)}\right\Vert _{\mathrm{F}}}_{\eqqcolon\alpha_{2,1}}+2\underbrace{\left\Vert \left(\bm{G}-\bm{G}^{(m)}\right)_{m,\cdot}\right\Vert _{2}\left\Vert \bm{U}_{m,\cdot}^{(m)}\bm{H}^{(m)}\right\Vert _{2}}_{\eqqcolon\alpha_{2,2}}.
\end{align*}
Here (i) follows from Lemma \ref{lemma:hpca-basic-facts}, and (ii)
follows from (\ref{eq:G-Gm-entries}). 
\begin{itemize}
\item Regarding $\alpha_{1}$, we can invoke (\ref{eq:hpca-useful-1-2})
in Lemma \ref{lemma:hpca-useful-1} to achieve that with probability
exceeding $1-O(n^{-11})$, 
\[
\alpha_{2,1}\lesssim\zeta_{\mathsf{op}}\left(\left\Vert \bm{U}^{(m)}\bm{H}^{(m)}\right\Vert _{2,\infty}+\sqrt{\frac{\mu^{\natural}r}{n_{1}}}\right).
\]
\item With regards to $\alpha_{2}$, we can invoke Lemma \ref{lemma:hpca-loo-basics}
to show that with probability exceeding $1-O(n^{-11})$, 
\[
\alpha_{2,2}\leq\left\Vert \bm{G}-\bm{G}^{(m)}\right\Vert _{2}\left\Vert \bm{U}^{(m)}\bm{H}^{(m)}\right\Vert _{2,\infty}\leq\left(\sigma^{2}\sqrt{n_{1}n_{2}}\log n+\kappa^{\natural2}\sigma\sigma_{1}^{\natural}\sqrt{\mu^{\natural}r\log n}\right)\left\Vert \bm{U}^{(m)}\bm{H}^{(m)}\right\Vert _{2,\infty}.
\]
\end{itemize}
Taking the above bounds on $\alpha_{2,1}$ and $\alpha_{2,2}$ collectively
yields 
\begin{align*}
 & \left\Vert \mathcal{P}_{\mathsf{off}\text{-}\mathsf{diag}}\left(\bm{G}-\bm{G}^{(m)}\right)\bm{U}^{(m)}\right\Vert \lesssim\alpha_{2,1}+\alpha_{2,2}\\
 & \qquad\lesssim\zeta_{\mathsf{op}}\left(\left\Vert \bm{U}^{(m)}\bm{H}^{(m)}\right\Vert _{2,\infty}+\sqrt{\frac{\mu^{\natural}r}{n_{1}}}\right)+\left(\sigma^{2}\sqrt{n_{1}n_{2}}\log n+\kappa^{\natural2}\sigma\sigma_{1}^{\natural}\sqrt{\mu^{\natural}r\log n}\right)\left\Vert \bm{U}^{(m)}\bm{H}^{(m)}\right\Vert _{2,\infty}\\
 & \qquad\lesssim\zeta_{\mathsf{op}}\left(\left\Vert \bm{U}^{(m)}\bm{H}^{(m)}\right\Vert _{2,\infty}+\sqrt{\frac{\mu^{\natural}r}{n_{1}}}\right),
\end{align*}
provided that $n_{1}\gtrsim\kappa^{\natural4}\mu^{\natural}r$. Therefore,
we arrive at 
\[
\alpha_{2}\lesssim\frac{\zeta_{\mathsf{op}}}{\sigma_{r}^{\natural2}}\left(\left\Vert \bm{U}^{(m)}\bm{H}^{(m)}\right\Vert _{2,\infty}+\sqrt{\frac{\mu^{\natural}r}{n_{1}}}\right).
\]

\paragraph{Combining bounds on $\alpha_{1}$ and $\alpha_{2}$.}

The preceding bounds on $\alpha_{1}$ and $\alpha_{2}$ combined allow
one to derive 
\begin{align*}
\left\Vert \bm{U}^{(m)}\bm{U}^{(m)\top}-\bm{U}\bm{U}^{\top}\right\Vert  & \leq\alpha_{1}+\alpha_{2}\lesssim\frac{\zeta_{\mathsf{op}}}{\sigma_{r}^{\natural2}}\left(\left\Vert \bm{U}^{(m)}\bm{H}^{(m)}\right\Vert _{2,\infty}+\sqrt{\frac{\mu^{\natural}r}{n_{1}}}\right)+\kappa^{\natural2}\sqrt{\frac{\mu^{\natural}r}{n_{1}}}\frac{\zeta_{\mathsf{op}}}{\sigma_{r}^{\natural2}}\\
 & \lesssim\frac{\zeta_{\mathsf{op}}}{\sigma_{r}^{\natural2}}\left\Vert \bm{U}^{(m)}\bm{H}^{(m)}\right\Vert _{2,\infty}+\kappa^{\natural2}\sqrt{\frac{\mu^{\natural}r}{n_{1}}}\frac{\zeta_{\mathsf{op}}}{\sigma_{r}^{\natural2}}
\end{align*}
with probability at least $1-O(n^{-11})$. This together with the
union bound concludes the proof. 

\section{Analysis for PCA: the approach based on \textsf{HeteroPCA\label{sec:Analysis:-the-approach-heteroPCA}}}

In this section, we establish our theoretical guarantees for the approach
based on \textsf{HeteroPCA ---} the ones presented in Section~\ref{subsec:Inference-for-HeteroPCA}.
. We shall establish our inference guarantees by connecting the PCA
model with the subspace estimation problem studied in Section~\ref{sec:detour-subspace}.

Before embarking on the analysis, we claim that: without loss of generality,
we can work with the assumption that 
\begin{equation}
\max_{j\in[n]}\left|\eta_{l,j}\right|\leq C_{\mathsf{noise}}\omega_{l}^{\star}\sqrt{\log\left(n+d\right)}\qquad\text{for all }l\in[d],\label{eq:noise-bounded}
\end{equation}
for some absolute constant $C_{\mathsf{noise}}>0$, in addition to
the noise condition we have already imposed in Section \ref{subsec:intro-model}.
To see why this is valid, we note that by applying Lemma \ref{lemma:truncate}
(with $\delta=C_{\delta}(n+d)^{-100}$ for some sufficiently small
constant $C_{\delta}>0$) to $\left\{ \eta_{l,j}\right\} $, we can
produce a set of auxiliary random variables $\left\{ \widetilde{\eta}_{l,j}\right\} $
such that 
\begin{itemize}
\item the $\widetilde{\eta}_{l,j}$'s are independent sub-Gaussian random
variables satisfying 
\[
\mathbb{E}[\widetilde{\eta}_{l,j}]=0,\qquad\widetilde{\omega}_{l}^{\star2}\coloneqq\mathbb{E}[\widetilde{\eta}_{l,j}^{2}]=\left[1+O\left(\left(n+d\right)^{-50}\right)\right]\omega_{l}^{\star2},\qquad\left\Vert \widetilde{\eta}_{l,j}\right\Vert _{\psi_{2}}\lesssim\omega_{l}^{\star},\qquad\left|\widetilde{\eta}_{l,j}\right|\lesssim\omega_{l}^{\star}\sqrt{\log\left(n+d\right)};
\]
\item these auxiliary variables satisfy 
\begin{equation} \label{eq:truncated-noise-identical}
\mathbb{P}\left(\eta_{l,j}=\widetilde{\eta}_{l,j}\text{ for all }l\in[d]\text{ and }j\in[n]\right)\geq1-O\left(\left(n+d\right)^{-98}\right).
\end{equation}
\end{itemize}
The above properties suggest that: if we replace the noise matrix
$\bm{N}=[\eta_{l,j}]_{l\in[d],j\in[n]}$ with $\widetilde{\bm{N}}=[\widetilde{\eta}_{l,j}]_{l\in[d],j\in[n]}$,
then the observations remain unchanged (i.e.~$\mathcal{P}_{\Omega}(\bm{X}+\bm{N})=\mathcal{P}_{\Omega}(\bm{X}+\widetilde{\bm{N}})$)
with high probability, and consequently, the resulting estimates are
also unchanged. In light of this, our analysis proceeds with the following
steps. 
\begin{itemize}
\item We shall start by proving Theorems \ref{thm:pca-complete}-\ref{thm:ce-CI-complete}
under the additional assumption (\ref{eq:noise-bounded}); if this
can be accomplished, then the results are clearly valid if we replace
$\bm{N}$ (resp.~$\{\omega_{l}^{\star}\}_{l\in[d]}$) with $\widetilde{\bm{N}}$
(resp.~$\{\widetilde{\omega}_{l}^{\star}\}_{l=1}^{d}$). 
\item Then we can invoke the proximity of $\widetilde{\omega}_l^\star$ and $\omega_l^\star$ as well as \eqref{eq:truncated-noise-identical} to establish
these theorems without assuming (\ref{eq:noise-bounded}).
\end{itemize}

\subsection{Connection between PCA and subspace estimation\label{subsec:Connection-between-PCA-subspace}}

In order to invoke our theoretical guarantees for subspace estimation
(i.e., Theorem~\ref{thm:hpca_inference_general}) to assist in understanding
PCA, we need to establish an explicit connection between these two
models. This forms the main content of this subsection.

\paragraph{Specification of $\bm{M}^{\natural}$, $\bm{M}$ and $\bm{E}$. }

Recall that in our PCA model in Section~\ref{subsec:intro-model},
we assume the covariance matrix $\bm{S}^{\star}$ admits the eigen-decomposition
$\bm{U}^{\star}\bm{\Sigma}^{\star2}\bm{U}^{\star\top}$. As a result,
the matrix $\bm{X}$ (cf.~\eqref{eq:definition-X-matrix}) can be
equivalently expressed as 
\[
\bm{X}=\bm{U}^{\star}\bm{\Sigma}^{\star}\left[\bm{f}_{1},\ldots,\bm{f}_{n}\right]=\bm{U}^{\star}\bm{\Sigma}^{\star}\bm{F},\qquad\text{where}\qquad\bm{f}_{i}\overset{\text{i.i.d.}}{\sim}\mathcal{N}\left(\bm{0},\bm{I}_{r}\right).
\]
Recalling that our observation matrix is $\bm{Y}=\mathcal{P}_{\Omega}(\bm{X}+\bm{N})$,
we can specify the matrices $\bm{M}^{\natural}$, $\bm{M}$ and $\bm{E}$
as defined in Appendix \ref{sec:subspace-setting} as follows: 
\begin{equation}
\bm{M}\coloneqq\frac{1}{\sqrt{n}p}\bm{Y},\qquad\bm{M}^{\natural}\coloneqq\mathbb{E}\left[\bm{M}\big|\bm{F}\right]=\frac{1}{\sqrt{n}}\bm{U}^{\star}\bm{\Sigma}^{\star}\bm{F}\qquad\text{and}\qquad\bm{E}\coloneqq\bm{M}-\bm{M}^{\natural}.\label{eq:defn-M-Mnatural-E-denoising-PCA}
\end{equation}
As usual, we shall let the SVD of $\bm{M}^{\natural}$ be $\bm{M}^{\natural}=\bm{U}^{\natural}\bm{\Sigma}^{\natural}\bm{V}^{\natural\top}$,
and use $\kappa^{\natural},\mu^{\natural},\sigma_{r}^{\natural},\sigma_{1}^{\natural}$
to denote the conditional number, the incoherence parameter, the minimum
and maximum singular value of $\bm{M}^{\natural}$, respectively.
Here and below, we shall focus on the randomness of noise and missing
data while treating $\bm{F}$ as given (even though it is generated
randomly).

We shall also specify several useful relations between $(\bm{U}^{\natural},\bm{V}^{\natural},\bm{R}_{\bm{U}})$
and $(\bm{U}^{\star},\bm{V}^{\star},\bm{R})$, where the rotation
matrix $\bm{R}$ is defined as $\bm{R}=\arg\min_{\bm{O}\in\mathcal{O}^{r\times r}}\Vert\bm{U}\bm{O}-\bm{U}^{\star}\Vert_{\mathrm{F}}$.
First, observe that 
\begin{equation}
\bm{V}^{\natural}=\bm{M}^{\natural\top}\bm{U}^{\natural}\left(\bm{\Sigma}^{\natural}\right)^{-1}=\frac{1}{\sqrt{n}}\bm{F}^{\top}\underset{\eqqcolon\,\bm{J}}{\underbrace{\bm{\Sigma}^{\star}\bm{U}^{\star\top}\bm{U}^{\natural}\left(\bm{\Sigma}^{\natural}\right)^{-1}}}.\label{eq:defn-J-denoising}
\end{equation}
Given that $\bm{U}^{\star}$ and $\bm{U}^{\natural}$ represent the
same column space, there exists $\bm{Q}\in\mathcal{O}^{r\times r}$
such that 
\begin{equation}
\bm{U}^{\natural}=\bm{U}^{\star}\bm{Q},\label{eq:defn-Q-denoising-PCA}
\end{equation}
thus leading to the expression 
\begin{equation}
\bm{J}=\bm{\Sigma}^{\star}\bm{Q}(\bm{\Sigma}^{\natural})^{-1}.\label{eq:defn-J-denoising-PCA}
\end{equation}
The definition of $\bm{R}$ together with $\bm{U}^{\natural}=\bm{U}^{\star}\bm{Q}$
also allows one to derive 
\begin{equation}
\bm{R}_{\bm{U}}=\arg\min_{\bm{O}\in\mathcal{O}^{r\times r}}\Vert\bm{U}\bm{O}-\bm{U}^{\natural}\Vert_{\mathrm{F}}=\arg\min_{\bm{O}\in\mathcal{O}^{r\times r}}\Vert\bm{U}\bm{O}\bm{Q}^{\top}-\bm{U}^{\star}\Vert_{\mathrm{F}}=\bm{R}\bm{Q}.\label{eq:connection-RU-R-denoisingPCA}
\end{equation}

\paragraph{Statistical properties of $\bm{E}=[E_{i,j}]$. }

We now describe the statistical properties of the perturbation matrix
$\bm{E}$. To begin with, it is readily seen from the definition of
$\bm{E}$ that $\mathbb{E}[\bm{E}\mid\bm{F}]=\bm{0}$. Moreover, from
the definition of $\bm{M}$ and $\bm{M}^{\natural}$, it is observed
that 
\begin{equation}
\bm{E}=\bm{M}-\bm{M}^{\natural\top}=\frac{1}{\sqrt{n}}\left[\frac{1}{p}\mathcal{P}_{\Omega}\left(\bm{U}^{\star}\bm{\Sigma}^{\star}\bm{F}+\bm{N}\right)-\bm{U}^{\star}\bm{\Sigma}^{\star}\bm{F}\right]=\frac{1}{\sqrt{n}}\left[\frac{1}{p}\mathcal{P}_{\Omega}\left(\bm{U}^{\star}\bm{\Sigma}^{\star}\bm{F}\right)-\bm{U}^{\star}\bm{\Sigma}^{\star}\bm{F}\right]+\frac{1}{\sqrt{n}}\frac{1}{p}\mathcal{P}_{\Omega}\left(\bm{N}\right),\label{eq:E-expression-denoising}
\end{equation}
which is clearly a zero-mean matrix conditional on $\bm{F}$. In addition,
for location $(i,j)$, we have 
\[
M_{i,j}^{\natural}=\frac{1}{\sqrt{n}}\left(\bm{U}^{\star}\bm{\Sigma}^{\star}\bm{F}\right)_{i,j}=\frac{1}{\sqrt{n}}\left(\bm{U}^{\star}\bm{\Sigma}^{\star}\bm{F}^{\star}\right)_{i,j}=\frac{1}{\sqrt{n}}\bm{U}_{i,\cdot}^{\star}\bm{\Sigma}^{\star}\bm{f}_{j},
\]
and hence the variance of $E_{i,j}$ can be calculated as 
\begin{align}
\sigma_{i,j}^{2} & \coloneqq\mathsf{var}\left(E_{i,j}|\bm{F}\right)=\mathbb{E}\left[E_{i,j}^{2}|\bm{F}\right]=\frac{1-p}{np}\left(\bm{U}^{\star}\bm{\Sigma}^{\star}\bm{F}\right)_{i,j}^{2}+\frac{\omega_{i}^{\star2}}{np}=\frac{1-p}{np}\left(\bm{U}_{i,\cdot}^{\star}\bm{\Sigma}^{\star}\bm{f}_{j}\right)^{2}+\frac{\omega_{i}^{\star2}}{np}.\label{eq:sigma-ij-square-denoising-PCA}
\end{align}

\paragraph{A good event $\mathcal{E}_{\mathsf{good}}$. }

Finally, the lemma below defines a high-probability event $\mathcal{E}_{\mathsf{good}}$
under which the random quantities defined above enjoy appealing properties.
The proof can be found in Appendix~\ref{appendix:good-event}.

\begin{lemma}\label{lemma:useful-property-good-event-hpca}There
is an event $\mathcal{E}_{\mathsf{good}}$ with $\mathbb{P}(\mathcal{E}_{\mathsf{good}})\geq1-O\big((n+d)^{-10}\big)$,
on which the following properties hold. 
\begin{itemize}
\item $\mathcal{E}_{\mathsf{good}}$ is $\sigma(\bm{F})$-measurable, where
$\sigma(\bm{F})$ is the $\sigma$-algebra generated by $\bm{F}$. 
\item If $n\gg\kappa^{2}(r+\log(n+d))$, then one has 
\begin{equation}
\left\Vert \frac{1}{n}\bm{F}\bm{F}^{\top}-\bm{I}_{r}\right\Vert \lesssim\sqrt{\frac{r+\log\left(n+d\right)}{n}},\label{eq:good-event-gaussian-concentration-hpca}
\end{equation}
\begin{equation}
\left\Vert \bm{\Sigma}^{\natural}-\bm{\Sigma}^{\star}\right\Vert \lesssim\kappa\sqrt{\frac{r+\log\left(n+d\right)}{n}}\sigma_{r}^{\star},\label{eq:good-event-Sigma-diff-hpca}
\end{equation}
\begin{equation}
\left\Vert \bm{\Sigma}^{\natural2}-\bm{\Sigma}^{\star2}\right\Vert \lesssim\sqrt{\frac{r+\log\left(n+d\right)}{n}}\sigma_{1}^{\star2},\label{eq:good-event-Sigma-2-diff-hpca}
\end{equation}
\begin{equation}
\sigma_{r}^{\natural}\asymp\sigma_{r}^{\star}\qquad\text{and}\qquad\sigma_{1}^{\natural}\asymp\sigma_{1}^{\star}.\label{eq:good-event-sigma-least-largest-hpca}
\end{equation}
\item The conditional number $\kappa^{\natural}$ and the incoherence parameter
$\mu^{\natural}$ of $\bm{M}^{\natural}$ obey 
\begin{align}
\kappa^{\natural} & \asymp\sqrt{\kappa},\label{eq:good-event-kappa-hpca}\\
\mu^{\natural} & \lesssim\kappa\mu\log\left(n+d\right).\label{eq:good-event-mu-hpca}
\end{align}
In addition, we have the following $\ell_{2,\infty}$ norm bound for
$\bm{U}^{\natural}$ and $\bm{V}^{\natural}$: 
\begin{equation}
\left\Vert \bm{U}^{\natural}\right\Vert _{2,\infty}\leq\sqrt{\frac{\mu r}{d}}\qquad\text{and}\qquad\left\Vert \bm{V}^{\natural}\right\Vert _{2,\infty}\lesssim\sqrt{\frac{r\log\left(n+d\right)}{n}}.\label{eq:good-event-V-natural-2-infty-hpca}
\end{equation}
\item The noise levels $\{\sigma_{i,j}\}$ are upper bounded by 
\begin{equation}
\sigma^{2}\coloneqq\max_{i\in[d],j\in[n]}\sigma_{i,j}^{2}\lesssim\frac{\mu r\log\left(n+d\right)}{ndp}\sigma_{1}^{\star2}+\frac{\omega_{\max}^{2}}{np}\eqqcolon\sigma_{\mathsf{ub}}^{2},\label{eq:good-event-sigma-hpca}
\end{equation}
\begin{align}
\max_{i\in[d],j\in[n]}\left|E_{i,j}\right| & \lesssim\frac{1}{p}\sqrt{\frac{\mu r\log\left(n+d\right)}{nd}}\sigma_{1}^{\star}+\frac{\omega_{\max}}{p}\sqrt{\frac{\log\left(n+d\right)}{n}}\eqqcolon B.\label{eq:good-event-B-hpca}
\end{align}
In addition, for each $i\in[d]$, 
\[
\sigma_{i}^{2}\coloneqq\max_{j\in[n]}\sigma_{i,j}^{2}\lesssim\frac{\log\left(n+d\right)}{np}\left\Vert \bm{U}_{i,\cdot}^{\star}\bm{\Sigma}^{\star}\right\Vert _{2}^{2}+\frac{\omega_{i}^{\star2}}{np}\coloneqq\sigma_{\mathsf{ub},i}^{2}
\]
and
\[
\max_{j\in[n]}\left|E_{i,j}\right|\lesssim\frac{1}{p}\sqrt{\frac{\log\left(n+d\right)}{n}}\left\Vert \bm{U}_{i,\cdot}^{\star}\bm{\Sigma}^{\star}\right\Vert _{2}+\frac{\omega_{l}^{\star}}{p}\sqrt{\frac{\log\left(n+d\right)}{n}}\eqqcolon B_{i}.
\]
\item The matrices $\bm{J}$ and $\bm{Q}$ are close in the sense that 
\begin{equation}
\left\Vert \bm{Q}-\bm{J}\right\Vert \lesssim\frac{1}{\sigma_{r}^{\star}}\left\Vert \bm{Q}\bm{\Sigma}^{\natural}-\bm{\Sigma}^{\star}\bm{Q}\right\Vert \lesssim\kappa\sqrt{\frac{r+\log\left(n+d\right)}{n}}.\label{eq:good-event-J-Q-hpca}
\end{equation}
\item For each $i\in[d]$, one has 
\begin{equation}
\max_{j\in[n]}\left|\bm{U}_{i,\cdot}^{\star}\bm{\Sigma}^{\star}\bm{f}_{j}\right|\lesssim\left\Vert \bm{U}_{i,\cdot}^{\star}\bm{\Sigma}^{\star}\right\Vert _{2}\sqrt{\log\left(n+d\right)};\label{eq:good-event-Ui-Sigma-f-hpca}
\end{equation}
for each $i,l\in[d]$, we have\begin{subequations} 
\begin{equation}
\biggl\Vert\frac{1}{n}\sum_{j=1}^{n}\left(\bm{U}_{l,\cdot}^{\star}\bm{\Sigma}^{\star}\bm{f}_{j}\right)^{2}\bm{f}_{j}\bm{f}_{j}^{\top}-\left\Vert \bm{U}_{l,\cdot}^{\star}\bm{\Sigma}^{\star}\right\Vert _{2}^{2}\bm{I}_{r}-2\bm{\Sigma}^{\star}\bm{U}_{l,\cdot}^{\star\top}\bm{U}_{l,\cdot}^{\star}\bm{\Sigma}^{\star}\biggr\Vert\lesssim\sqrt{\frac{r\log^{3}\left(n+d\right)}{n}}\left\Vert \bm{U}_{l,\cdot}^{\star}\bm{\Sigma}^{\star}\right\Vert _{2}^{2},\label{eq:good-event-concentration-1-hpca}
\end{equation}
\begin{equation}
\left|\frac{1}{n}\sum_{j=1}^{n}\left(\bm{U}_{l,\cdot}^{\star}\bm{\Sigma}^{\star}\bm{f}_{j}\right)^{2}\left(\bm{U}_{i,\cdot}^{\star}\bm{\Sigma}^{\star}\bm{f}_{j}\right)^{2}-\left(S_{l,l}^{\star}S_{i,i}^{\star}+2S_{i,l}^{\star2}\right)\right|\lesssim\sqrt{\frac{\log\left(n+d\right)^{3}}{n}}\left\Vert \bm{U}_{l,\cdot}^{\star}\bm{\Sigma}^{\star}\right\Vert _{2}^{2}\left\Vert \bm{U}_{i,\cdot}^{\star}\bm{\Sigma}^{\star}\right\Vert _{2}^{2}\label{eq:good-event-concentration-2-hpca}
\end{equation}
\begin{equation}
\text{and}\qquad\left|\frac{1}{n}\sum_{j=1}^{n}\left(\bm{U}_{i,\cdot}^{\star}\bm{\Sigma}^{\star}\bm{f}_{j}\right)^{2}-S_{i,i}^{\star}\right|\lesssim\sqrt{\frac{\log\left(n+d\right)}{n}}\left\Vert \bm{U}_{i,\cdot}^{\star}\bm{\Sigma}^{\star}\right\Vert _{2}^{2}.\label{eq:good-event-concentration-3-hpca}
\end{equation}
\end{subequations} 
\end{itemize}
\end{lemma}\begin{proof}See Appendix \ref{appendix:good-event}.\end{proof}

\subsection{Distributional characterization for principal subspace (Proof of
Theorem \ref{thm:pca-complete})}

With the above connection between the subspace estimation model and
PCA in place, we can readily move on to invoke Theorem \ref{thm:hpca_inference_general}
to establish our distributional characterization of $\bm{U}\bm{R}-\bm{U}^{\star}$
for \textsf{HeteroPCA} (as stated in Theorem \ref{thm:pca-complete}).

\paragraph{Step 1: first- and second-order approximation and the tightness. }

As a starting point, Theorem \ref{thm:hpca_inference_general} taken
together with the explicit connection between the subspace estimation
model and PCA allows one to approximate $\bm{U}\bm{R}-\bm{U}^{\star}$
in a concise form, which is a crucial first step that enables our
subsequent development of the distributional theory. The proof can
be found in Appendix \ref{appendix:proof-pca-2nd-error}.

\begin{lemma} \label{lemma:pca-2nd-error}Assume that $d\gtrsim\kappa^{3}\mu^{2}r\log^{4}(n+d)$,
$n\gtrsim r\log^{4}(n+d)$, $ndp^{2}\gg\kappa^{4}\mu^{2}r^{2}\log^{4}(n+d)$,
$np\gg\kappa^{4}\mu r\log^{2}\left(n+d\right)$, 
\[
\frac{\omega_{\max}^{2}}{p\sigma_{r}^{\star2}}\sqrt{\frac{d}{n}}\ll\frac{1}{\kappa\log\left(n+d\right)}\qquad\text{and}\qquad\frac{\omega_{\max}}{\sigma_{r}^{\star}}\sqrt{\frac{d}{np}}\ll\frac{1}{\sqrt{\kappa^{3}\log\left(n+d\right)}}.
\]
Then for each $l\in[d]$, we have
\[
\mathbb{P}\left(\left\Vert \left(\bm{U}\bm{R}-\bm{U}^{\star}-\bm{Z}\right)_{l,\cdot}\right\Vert _{2}\ind_{\mathcal{E}_{\mathsf{good}}}\lesssim\zeta_{\mathsf{2nd},l}\,\big|\,\bm{F}\right)\geq1-O\left(\left(n+d\right)^{-10}\right)
\]
almost surely, where 
\begin{equation}
\bm{Z}\coloneqq\left[\bm{E}\bm{M}^{\natural\top}+\mathcal{P}_{\mathsf{off}\text{-}\mathsf{diag}}\left(\bm{E}\bm{E}^{\top}\right)\right]\bm{U}^{\natural}\left(\bm{\Sigma}^{\natural}\right)^{-2}\bm{Q}^{\top}\label{eq:defn-Z-distribution-hpca}
\end{equation}
and\begin{subequations}\label{eq:zeta-1st-2nd-UB-hpca} 
\begin{align}
\zeta_{\mathsf{2nd},l} & \coloneqq\left\Vert \bm{U}_{l,\cdot}^{\star}\right\Vert _{2}\left(\sqrt{\frac{\kappa^{3}\mu r\log\left(n+d\right)}{d}}\frac{\zeta_{\mathsf{1st}}}{\sigma_{r}^{\star2}}+\kappa\frac{\zeta_{\mathsf{1st}}^{2}}{\sigma_{r}^{\star4}}\right)+\frac{\zeta_{\mathsf{1st}}\zeta_{\mathsf{1st},l}}{\sigma_{r}^{\star4}}\sqrt{\frac{\kappa^{3}\mu r\log\left(n+d\right)}{d}},\label{eq:zeta-2nd-UB-hpca}\\
\zeta_{\mathsf{1st}} & \coloneqq\frac{\mu r\log^{2}\left(n+d\right)}{\sqrt{nd}p}\sigma_{1}^{\star2}+\frac{\omega_{\max}^{2}}{p}\sqrt{\frac{d}{n}}\log\left(n+d\right)+\sigma_{1}^{\star2}\sqrt{\frac{\mu r}{np}}\log\left(n+d\right)+\sigma_{1}^{\star}\omega_{\max}\sqrt{\frac{d\log\left(n+d\right)}{np}},\label{eq:zeta-1st-UB-hpca}\\
\zeta_{\mathsf{1st},l} & \coloneqq\sqrt{\frac{\mu r\log^{4}\left(n+d\right)}{np^{2}}}\sigma_{1}^{\star}\left\Vert \bm{U}_{l,\cdot}^{\star}\bm{\Sigma}^{\star}\right\Vert _{2}+\frac{\omega_{l}^{\star}\omega_{\max}}{p}\sqrt{\frac{d}{n}}\log\left(n+d\right)+\sigma_{1}^{\star}\sqrt{\frac{d}{np}}\left\Vert \bm{U}_{l,\cdot}^{\star}\bm{\Sigma}^{\star}\right\Vert _{2}\log\left(n+d\right)\nonumber \\
 & \quad+\sigma_{1}^{\star}\omega_{l}^{\star}\sqrt{\frac{d\log\left(n+d\right)}{np}}+\frac{\omega_{\max}}{p}\sqrt{\frac{d}{n}}\log^{3/2}\left(n+d\right)\left\Vert \bm{U}_{l,\cdot}^{\star}\bm{\Sigma}^{\star}\right\Vert _{2}+\frac{\omega_{l}^{\star}}{p}\sqrt{\frac{\mu r\log^{3}\left(n+d\right)}{n}}\sigma_{1}^{\star}.\label{eq:zeta-1st-l-UB-hpca}
\end{align}
\end{subequations}\end{lemma}

In particular, we expect each $\zeta_{\mathsf{2nd},l}$ to be negligible,
so that the approximation $\bm{U}\bm{R}-\bm{U}^{\star}\approx\bm{Z}$
is nearly tight. It is worth noting that the approximation $\bm{Z}$
consists of both linear and second-order effects of the perturbation
matrix $\bm{E}$.

\paragraph{Step 2: computing the covariance of the first- and second-order approximation.}

In order to pin down the distribution of $\bm{Z}$, an important step
lies in characterizing its covariance. To be precise, observe that
the $l$-th row of $\bm{Z}$ ($1\leq l\leq d$) defined in \eqref{eq:defn-Z-distribution-hpca}
satisfies 
\begin{align}
\bm{Z}_{l,\cdot} & =\bm{e}_{l}^{\top}\left[\bm{E}\bm{M}^{\natural\top}+\mathcal{P}_{\mathsf{off}\text{-}\mathsf{diag}}\left(\bm{E}\bm{E}^{\top}\right)\right]\bm{U}^{\natural}(\bm{\Sigma}^{\natural})^{-2}\bm{Q}^{\top}\nonumber \\
 & =\bm{E}_{l,\cdot}\bm{V}^{\natural}(\bm{\Sigma}^{\natural})^{-1}+\bm{E}_{l,\cdot}\left[\mathcal{P}_{-l,\cdot}\left(\bm{E}\right)\right]^{\top}\bm{U}^{\natural}(\bm{\Sigma}^{\natural})^{-2}\bm{Q}^{\top}\nonumber \\
 & =\sum_{j=1}^{n}E_{l,j}\left[\bm{V}_{j,\cdot}^{\natural}(\bm{\Sigma}^{\natural})^{-1}+\left[\mathcal{P}_{-l,\cdot}\left(\bm{E}_{\cdot,j}\right)\right]^{\top}\bm{U}^{\natural}(\bm{\Sigma}^{\natural})^{-2}\right]\bm{Q}^{\top},\label{eq:Zl-decomposition-E-iid-hpca}
\end{align}
where we recall that $\mathcal{P}_{-l,\cdot}\left(\bm{E}\right)$
is obtained by zeroing out the $l$-th row of $\bm{E}$. Conditional
on $\bm{F}$, the covariance matrix of this zero-mean random vector
$\bm{Z}_{l,\cdot}$ can thus be calculated as 
\begin{align}
\widetilde{\bm{\Sigma}}_{l} & \coloneqq\bm{Q}(\bm{\Sigma}^{\natural})^{-1}\bm{V}^{\natural\top}\mathsf{diag}\left\{ \sigma_{l,1}^{2},\ldots,\sigma_{l,n}^{2}\right\} \bm{V}^{\natural}(\bm{\Sigma}^{\natural})^{-1}\bm{Q}^{\top}\nonumber \\
 & \qquad+\bm{Q}(\bm{\Sigma}^{\natural})^{-2}\bm{U}^{\natural\top}\mathsf{diag}\Biggl\{\sum_{j:j\neq l}\sigma_{l,j}^{2}\sigma_{1,j}^{2},\ldots,\sum_{j:j\neq l}\sigma_{l,j}^{2}\sigma_{d,j}^{2}\Biggr\}\bm{U}^{\natural}(\bm{\Sigma}^{\natural})^{-2}\bm{Q}^{\top}.\label{eq:Sigma-tilde-distribution-hpca}
\end{align}

Given that the above expression of $\widetilde{\bm{\Sigma}}_{l}$
contains components like $\bm{\Sigma}^{\natural}$ and $\bm{V}^{\natural}$
(which are introduced in order to use the subspace estimation model),
it is natural to see whether one can express $\widetilde{\bm{\Sigma}}_{l}$
directly in terms of the corresponding quantities introduced for the
PCA model. To do so, we first single out a deterministic matrix as
follows 
\begin{equation}
\bm{\Sigma}_{U,l}^{\star}\coloneqq\left(\frac{1-p}{np}\left\Vert \bm{U}_{l,\cdot}^{\star}\bm{\Sigma}^{\star}\right\Vert _{2}^{2}+\frac{\omega_{l}^{\star2}}{np}\right)\left(\bm{\Sigma}^{\star}\right)^{-2}+\frac{2\left(1-p\right)}{np}\bm{U}_{l,\cdot}^{\star\top}\bm{U}_{l,\cdot}^{\star}+\left(\bm{\Sigma}^{\star}\right)^{-2}\bm{U}^{\star\top}\mathsf{diag}\left\{ \left[d_{l,i}^{\star}\right]_{i=1}^{d}\right\} \bm{U}^{\star}\left(\bm{\Sigma}^{\star}\right)^{-2},\label{eq:defn-Sigma-Ul-star-hpca}
\end{equation}
where we define 
\begin{equation}
d_{l,i}^{\star}\coloneqq\frac{1}{np^{2}}\left[\omega_{l}^{\star2}+\left(1-p\right)\left\Vert \bm{U}_{l,\cdot}^{\star}\bm{\Sigma}^{\star}\right\Vert _{2}^{2}\right]\left[\omega_{i}^{\star2}+\left(1-p\right)\left\Vert \bm{U}_{i,\cdot}^{\star}\bm{\Sigma}^{\star}\right\Vert _{2}^{2}\right]+\frac{2\left(1-p\right)^{2}}{np^{2}}S_{i,l}^{\star2}.\label{eq:defn-dli-star-hpca}
\end{equation}
In view of the following lemma, $\bm{\Sigma}_{U,l}^{\star}$ approximates
$\widetilde{\bm{\Sigma}}_{l}$ in a reasonably well fashion.

\begin{lemma}\label{lemma:pca-covariance-concentration}Suppose that
$n\gg\kappa^{8}\mu^{2}r^{3}\kappa_{\omega}^{2}\log^{3}(n+d)$. On
the event $\mathcal{E}_{\mathsf{good}}$ (cf.~Lemma~\ref{lemma:useful-property-good-event-hpca})
, we have 
\begin{align*}
\left\Vert \widetilde{\bm{\Sigma}}_{l}-\bm{\Sigma}_{U,l}^{\star}\right\Vert  & \lesssim\sqrt{\frac{\kappa^{8}\mu^{2}r^{3}\kappa_{\omega}^{2}\log^{3}\left(n+d\right)}{n}}\lambda_{\min}\left(\bm{\Sigma}_{U,l}^{\star}\right),\\
\max\left\{ \lambda_{\max}\left(\widetilde{\bm{\Sigma}}_{l}\right),\lambda_{\max}\left(\bm{\Sigma}_{U,l}^{\star}\right)\right\}  & \lesssim\frac{1-p}{np\sigma_{r}^{\star2}}\left\Vert \bm{U}_{l,\cdot}^{\star}\bm{\Sigma}^{\star}\right\Vert _{2}^{2}+\frac{\omega_{l}^{\star2}}{np\sigma_{r}^{\star2}}+\frac{\kappa\mu r\left(1-p\right)^{2}}{ndp^{2}\sigma_{r}^{\star2}}\left\Vert \bm{U}_{l,\cdot}^{\star}\bm{\Sigma}^{\star}\right\Vert _{2}^{2}+\frac{\kappa\mu r\left(1-p\right)}{ndp^{2}\sigma_{r}^{\star2}}\omega_{l}^{\star2}\\
 & \quad+\frac{1-p}{np^{2}\sigma_{r}^{\star4}}\omega_{\max}^{2}\left\Vert \bm{U}_{l,\cdot}^{\star}\bm{\Sigma}^{\star}\right\Vert _{2}^{2}+\frac{\omega_{l}^{\star2}\omega_{\max}^{2}}{np^{2}\sigma_{r}^{\star4}},\\
\min\left\{ \lambda_{\min}\left(\widetilde{\bm{\Sigma}}_{l}\right),\lambda_{\min}\left(\bm{\Sigma}_{U,l}^{\star}\right)\right\}  & \gtrsim\frac{1-p}{np\sigma_{1}^{\star2}}\left\Vert \bm{U}_{l,\cdot}^{\star}\bm{\Sigma}^{\star}\right\Vert _{2}^{2}+\frac{\omega_{l}^{\star2}}{np\sigma_{1}^{\star2}}+\frac{\left(1-p\right)^{2}}{ndp^{2}\kappa\sigma_{1}^{\star2}}\left\Vert \bm{U}_{l,\cdot}^{\star}\bm{\Sigma}^{\star}\right\Vert _{2}^{2}+\frac{1-p}{ndp^{2}\kappa\sigma_{1}^{\star2}}\omega_{l}^{\star2}\\
 & \quad+\frac{1-p}{np^{2}\sigma_{1}^{\star4}}\omega_{\min}^{2}\left\Vert \bm{U}_{l,\cdot}^{\star}\bm{\Sigma}^{\star}\right\Vert _{2}^{2}+\frac{\omega_{l}^{\star2}\omega_{\min}^{2}}{np^{2}\sigma_{1}^{\star4}}.
\end{align*}
In addition, the condition number of $\widetilde{\bm{\Sigma}}_{l}$
is bounded above by $O(\kappa^{3}\mu r\kappa_{\omega})$.\end{lemma}\begin{proof}See
Appendix \ref{appendix:proof-pca-covariance-concentration}. \end{proof}

\paragraph{Step 3: establishing distributional guarantees for $(\bm{U}\bm{R}-\bm{U}^{\star})_{l,\cdot}$. }

By virtue of the decomposition \eqref{eq:Zl-decomposition-E-iid-hpca},
each row $\bm{Z}_{l,\cdot}$ can be viewed as the sum of a collection
of independent zero-mean random variables/vectors. This suggests that
$\bm{Z}_{l,\cdot}$ might be well approximated by certain multivariate
Gaussian distributions. If so, then the zero-mean nature of $\bm{Z}_{l,\cdot}$
in conjunction with the above-mentioned covariance formula allows
us to pin down the distribution of each row of $\bm{U}\bm{R}-\bm{U}^{\star}$
approximately.

\begin{lemma}[\bf Gaussian approximation of $\bm{Z}_{l,\cdot}$]\label{lemma:pca-normal-approximation}
Suppose that the following conditions hold: 
\[
n\gtrsim\kappa^{8}\mu^{2}r^{4}\kappa_{\omega}^{2}\log^{4}(n+d),\qquad d\gtrsim\kappa^{7}\mu^{3}r^{7/2}\kappa_{\omega}^{2}\log^{5}(n+d),
\]
\[
np\gtrsim\kappa^{9}\mu^{3}r^{11/2}\kappa_{\omega}^{2}\log^{7}\left(n+d\right),\qquad ndp^{2}\gtrsim\kappa^{9}\mu^{4}r^{13/2}\kappa_{\omega}^{2}\log^{9}\left(n+d\right),
\]
and
\[
\frac{\omega_{\max}}{\sigma_{r}^{\star}}\sqrt{\frac{d}{np}}\lesssim\frac{1}{\kappa^{3/2}\mu r^{5/4}\kappa_{\omega}^{1/2}\log^{3}\left(n+d\right)},\qquad\frac{\omega_{\max}^{2}}{p\sigma_{r}^{\star2}}\sqrt{\frac{d}{n}}\lesssim\frac{1}{\kappa^{4}\mu^{3/2}r^{11/4}\kappa_{\omega}\log^{7/2}\left(n+d\right)}.
\]
Then it is guaranteed that 
\[
\sup_{\mathcal{C}\in\mathscr{C}^{r}}\left|\mathbb{P}\left(\left(\bm{U}\bm{R}-\bm{U}^{\star}\right)_{l,\cdot}\in\mathcal{C}\right)-\mathcal{N}\left(\bm{0},\bm{\Sigma}_{U,l}^{\star}\right)\left\{ \mathcal{C}\right\} \right|\lesssim\frac{1}{\sqrt{\log\left(n+d\right)}}=o\left(1\right),
\]
where $\mathscr{C}^{r}$ denotes the set of all convex sets in $\mathbb{R}^{r}$.\end{lemma}\begin{proof}See
Appendix \ref{appendix:proof-pca-normal-approximation}. \end{proof}

Once Lemma \ref{lemma:pca-normal-approximation} is established, we
have solidified the advertised Gaussian approximation of each row
of $\bm{U}\bm{R}-\bm{U}^{\star}$, thus concluding the proof of Theorem
\ref{thm:pca-complete}. 

\subsection{Validity of confidence regions (Proof of Theorem \ref{thm:pca-cr-complete})\label{appendix:proof-pca-cr}}

Our distributional characterization of $(\bm{U}\bm{R}-\bm{U}^{\star})_{l,\cdot}$
hints at the possibility of constructing valid confidence region for
$\bm{U}^{\star}$, provided that the covariance matrix $\bm{\Sigma}_{U,l}^{\star}$
(cf.~\eqref{eq:defn-Sigma-Ul-star-hpca}) can be reliably estimated.
Similar to the SVD-based approach, we attempt to estimate $\bm{\Sigma}_{U,l}^{\star}$
by means of the following plug-in estimator: 
\begin{align}
\bm{\Sigma}_{U,l} & \coloneqq\left(\frac{1-p}{np}\left\Vert \bm{U}_{l,\cdot}\bm{\Sigma}\right\Vert _{2}^{2}+\frac{\omega_{l}^{2}}{np}\right)\bm{\Sigma}^{-2}+\frac{2\left(1-p\right)}{np}\bm{U}_{l,\cdot}^{\top}\bm{U}_{l,\cdot}+\left(\bm{\Sigma}\right)^{-2}\bm{U}^{\top}\mathsf{diag}\left\{ \left[d_{l,i}\right]_{1\leq i\leq d}\right\} \bm{U}(\bm{\Sigma})^{-2},\label{eq:defn-plug-in-Sigma-Ul-hpca}
\end{align}
where for each $i\in[d]$, we define 
\begin{align}
d_{l,i} & \coloneqq\frac{1}{np^{2}}\left[\omega_{l}^{2}+\left(1-p\right)\left\Vert \bm{U}_{l,\cdot}\bm{\Sigma}\right\Vert _{2}^{2}\right]\left[\omega_{i}^{2}+\left(1-p\right)\left\Vert \bm{U}_{i,\cdot}\bm{\Sigma}\right\Vert _{2}^{2}\right]+\frac{2\left(1-p\right)^{2}}{np^{2}}S_{l,i}^{2},\label{eq:defn-d-li-square-CR-hpca}\\
\omega_{i}^{2} & \coloneqq\frac{\sum_{j=1}^{n}y_{i,j}^{2}\ind_{(i,j)\in\Omega}}{\sum_{j=1}^{n}\ind_{(i,j)\in\Omega}}-S_{i,i}.\label{eq:defn-omega-i-square-CR-hpca}
\end{align}

\paragraph{Step 1: fine-grained estimation guarantees for $\bm{U}$ and $\bm{\Sigma}$.}

To begin with, we need to show that the components in the plug-in
estimator are all reliable estimates of their deterministic counterpart
(after proper rotation).

\begin{lemma}\label{lemma:pca-1st-err}Recall the definition of $\zeta_{\mathsf{1st}}$
and $\zeta_{\mathsf{2nd},l}$ in \eqref{eq:zeta-1st-2nd-UB-hpca}.
Assume that $\zeta_{\mathsf{1st}}/\sigma_{r}^{\star2}\lesssim1/\sqrt{\kappa^{3}\mu}$,
$n\gg r+\log(n+d)$ and $d\gtrsim\kappa^{4}\mu^{2}r\log(n+d)$. Let
\begin{equation}
\theta\coloneqq\sqrt{\frac{\kappa r\log^{2}\left(n+d\right)}{np}}\left(1+\frac{\sigma_{\mathsf{ub}}}{\sigma_{r}^{\star}}\sqrt{n}\right)\label{eq:theta-definition}
\end{equation}
Then with probability exceeding $1-O((n+d)^{-10})$, we have\begin{subequations}\label{eq:pca-U-Sigma-R-bound-together}
\begin{align}
\left\Vert \left(\bm{U}\bm{R}-\bm{U}^{\star}\right)_{l,\cdot}\right\Vert _{2} & \lesssim\frac{\theta}{\sqrt{\kappa}\sigma_{r}^{\star}}\left(\left\Vert \bm{U}_{l,\cdot}^{\star}\bm{\Sigma}^{\star}\right\Vert _{2}+\omega_{l}^{\star}\right)+\zeta_{\mathsf{2nd},l},\label{eq:pca-U-l-R-error}\\
\left\Vert \left(\bm{U}\bm{\Sigma}\bm{R}-\bm{U}^{\star}\bm{\Sigma}^{\star}\right)_{l,\cdot}\right\Vert _{2} & \lesssim\theta\left(\left\Vert \bm{U}_{l,\cdot}^{\star}\bm{\Sigma}^{\star}\right\Vert _{2}+\omega_{l}^{\star}\right)+\left\Vert \bm{U}_{l,\cdot}^{\star}\right\Vert _{2}\sqrt{\frac{\kappa^{2}\left(r+\log\left(n+d\right)\right)}{n}}\sigma_{1}^{\star}+\zeta_{\mathsf{2nd},l}\sigma_{1}^{\star},\label{eq:pca-U-l-Sigma-R-error}\\
\left\Vert \bm{R}\left(\bm{\Sigma}^{\star}\right)^{-2}\bm{R}^{\top}-\bm{\Sigma}^{-2}\right\Vert  & \lesssim\sqrt{\frac{\kappa^{3}\mu r\log\left(n+d\right)}{d}}\frac{\zeta_{\mathsf{1st}}}{\sigma_{r}^{\star4}}+\kappa\frac{\zeta_{\mathsf{1st}}^{2}}{\sigma_{r}^{\star6}}+\sqrt{\frac{\kappa^{3}\left(r+\log\left(n+d\right)\right)}{n}}\frac{1}{\sigma_{r}^{\star2}},\label{eq:pca-Sigma-2-error}\\
\left\Vert \bm{U}\bm{\Sigma}^{-2}\bm{R}-\bm{U}^{\star}\left(\bm{\Sigma}^{\star}\right)^{-2}\right\Vert  & \lesssim\frac{\zeta_{\mathsf{1st}}}{\sigma_{r}^{\star4}}+\sqrt{\frac{\kappa^{3}\left(r+\log\left(n+d\right)\right)}{n}}\frac{1}{\sigma_{r}^{\star2}},\label{eq:pca-U-Sigma-n2-R-error}\\
\left\Vert \bm{U}\bm{R}-\bm{U}^{\star}\right\Vert _{2,\infty} & \lesssim\frac{\zeta_{\mathsf{1st}}}{\sigma_{r}^{\star2}}\sqrt{\frac{r\log\left(n+d\right)}{d}}.\label{eq:pca-two-to-infty-error}
\end{align}
\end{subequations}\end{lemma}\begin{proof} See Appendix \ref{appendix:proof-pca-1st-err}.\end{proof}

\paragraph{Step 2: faithfulness of the plug-in estimator.}

With Lemma \ref{lemma:pca-1st-err} in hand, we move forward to show
that $\bm{S}$ and $\{\omega_{i}^{2}\}_{i=1}^{d}$ are reliable estimators
of the covariance matrix $\bm{S}^{\star}$ and the noise levels $\{\omega_{i}^{\star2}\}_{i=1}^{d}$,
respectively.

\begin{lemma}\label{lemma:pca-noise-level-est}Suppose that the conditions
of Lemma \ref{lemma:pca-1st-err} hold. In addition, assume that $n\gtrsim\kappa^{3}r\log(n+d)$,
$\zeta_{\mathsf{1st}}/\sigma_{r}^{\star2}\lesssim1/\sqrt{\log(n+d)}$.
Then with probability exceeding $1-O((n+d)^{-10})$, we have 
\begin{equation}
\left\Vert \bm{S}-\bm{S}^{\star}\right\Vert _{\infty}\lesssim\left(\frac{\zeta_{\mathsf{1st}}}{\sigma_{r}^{\star}}\sqrt{\frac{\kappa r\log\left(n+d\right)}{d}}+\sqrt{\frac{\kappa^{2}\mu r^{2}\log\left(n+d\right)}{nd}}\sigma_{1}^{\star}\right)\sqrt{\frac{\mu r}{d}}\sigma_{1}^{\star}\label{eq:pca-S-infty-err}
\end{equation}
and, for each $i,j\in[d]$, 
\begin{align}
\left|S_{i,j}-S_{i,j}^{\star}\right| & \lesssim\left(\theta+\sqrt{\frac{\kappa^{3}r\log\left(n+d\right)}{n}}\right)\left\Vert \bm{U}_{i,\cdot}^{\star}\bm{\Sigma}^{\star}\right\Vert _{2}\left\Vert \bm{U}_{j,\cdot}^{\star}\bm{\Sigma}^{\star}\right\Vert _{2}+\theta\left(\omega_{i}^{\star}\left\Vert \bm{U}_{j,\cdot}^{\star}\bm{\Sigma}^{\star}\right\Vert _{2}+\omega_{j}^{\star}\left\Vert \bm{U}_{i,\cdot}^{\star}\bm{\Sigma}^{\star}\right\Vert _{2}\right)\nonumber \\
 & \quad+\sigma_{1}^{\star}\left(\zeta_{\mathsf{2nd},i}\left\Vert \bm{U}_{j,\cdot}^{\star}\bm{\Sigma}^{\star}\right\Vert _{2}+\zeta_{\mathsf{2nd},j}\left\Vert \bm{U}_{i,\cdot}^{\star}\bm{\Sigma}^{\star}\right\Vert _{2}\right)+\theta^{2}\omega_{i}^{\star}\omega_{j}^{\star}+\zeta_{\mathsf{2nd},i}\zeta_{\mathsf{2nd},j}\sigma_{1}^{\star2}.\label{eq:pca-S-entrywise-err}
\end{align}
Here, the quantity $\theta$ is defined in (\ref{eq:theta-definition}).
In addition, with probability exceeding $1-O((n+d)^{-10})$, we have
\begin{equation}
\left|\omega_{i}^{2}-\omega_{i}^{\star2}\right|\lesssim\sqrt{\frac{\log^{2}\left(n+d\right)}{np}}\omega_{i}^{\star2}+\zeta_{\mathsf{1st}}\frac{\sqrt{\kappa^{2}\mu r^{2}\log\left(n+d\right)}}{d}+\sqrt{\frac{\kappa^{2}\mu^{2}r^{3}\log\left(n+d\right)}{nd^{2}}}\sigma_{1}^{\star2}\label{eq:pca-noise-est-all}
\end{equation}
for all $i\in[d]$, and 
\begin{align}
\left|\omega_{l}^{2}-\omega_{l}^{\star2}\right| & \lesssim\left(\sqrt{\frac{\log^{2}\left(n+d\right)}{np}}+\theta^{2}\right)\omega_{l}^{\star2}+\left(\theta+\sqrt{\frac{\kappa^{3}r\log\left(n+d\right)}{n}}\right)\left\Vert \bm{U}_{l,\cdot}^{\star}\bm{\Sigma}^{\star}\right\Vert _{2}^{2}\nonumber \\
 & \quad+\left(\theta\omega_{l}^{\star}+\zeta_{\mathsf{2nd},l}\sigma_{1}^{\star}\right)\left\Vert \bm{U}_{l,\cdot}^{\star}\bm{\Sigma}^{\star}\right\Vert _{2}+\zeta_{\mathsf{2nd},l}^{2}\sigma_{1}^{\star2}.\label{eq:pca-noise-est-l}
\end{align}
\end{lemma}\begin{proof}See Appendix \ref{appendix:proof-pca-noise-level-est}.\end{proof}

The above two lemmas taken together allow us to demonstrate that our
plug-in estimator $\bm{\Sigma}_{U,l}$ is a faithful estimate of $\bm{\Sigma}_{U,l}^{\star}$,
as stated below.

\begin{lemma}\label{lemma:pca-covariance-estimation} Suppose that
the conditions of Lemma \ref{lemma:pca-normal-approximation}, Lemma
\ref{lemma:pca-1st-err} and Lemma \ref{lemma:pca-covariance-estimation}
hold. Consider any $\delta\in(0,1)$, and we further suppose that
$n\gtrsim\delta^{-2}\kappa^{9}\mu^{2}r^{3}\log(n+d)$, $d\gtrsim\kappa^{3}\mu r\log(n+d)$,
\[
ndp^{2}\gtrsim\delta^{-2}\kappa^{8}\mu^{4}r^{4}\kappa_{\omega}^{2}\log^{5}\left(n+d\right),\qquad np\gtrsim\delta^{-2}\kappa^{8}\mu^{3}r^{3}\kappa_{\omega}^{2}\log^{3}\left(n+d\right),
\]
\[
\frac{\omega_{\max}^{2}}{p\sigma_{r}^{\star2}}\sqrt{\frac{d}{n}}\lesssim\frac{\delta}{\kappa^{3}\mu r\kappa_{\omega}\log^{3/2}\left(n+d\right)}\qquad\text{and}\qquad\frac{\omega_{\max}}{\sigma_{r}^{\star}}\sqrt{\frac{d}{np}}\lesssim\frac{\delta}{\kappa^{7/2}\mu r\kappa_{\omega}\log\left(n+d\right)}.
\]
Then with probability exceeding $1-O((n+d)^{-10})$, we have 
\[
\left\Vert \bm{\Sigma}_{U,l}-\bm{R}\bm{\Sigma}_{U,l}^{\star}\bm{R}^{\top}\right\Vert \lesssim\delta\lambda_{\min}\left(\bm{\Sigma}_{U,l}^{\star}\right).
\]
\end{lemma}\begin{proof}See Appendix \ref{appendix:proof-pca-covariance-estimation}.
\end{proof}

\paragraph{Step 3: validity of the constructed confidence regions.}

Finally, we can combine the Gaussian approximation of $\bm{U}_{l,\cdot}\bm{R}-\bm{U}^{\star}$
established in Theorem \ref{thm:pca-complete} and the estimation
guarantee of the covariance matrix established in Lemma \ref{lemma:pca-covariance-estimation}
to justify the validity of the constructed confidence region $\mathsf{CR}_{U,l}^{1-\alpha}$
(cf.~Algorithm \ref{alg:PCA-HeteroPCA-CR}).

\begin{lemma}\label{lemma:pca-cr-validity}Suppose that the conditions
of Theorem \ref{thm:pca-complete} hold. Suppose that $n\gtrsim\kappa^{12}\mu^{3}r^{11/2}\kappa_{\omega}\log^{5}(n+d)$,
$d\gtrsim\kappa^{3}\mu r\log(n+d)$, 
\[
ndp^{2}\gtrsim\kappa^{11}\mu^{5}r^{13/2}\kappa_{\omega}^{3}\log^{9}\left(n+d\right),\qquad np\gtrsim\kappa^{11}\mu^{4}r^{11/2}\kappa_{\omega}^{3}\log^{7}\left(n+d\right),
\]
and 
\[
\frac{\omega_{\max}^{2}}{p\sigma_{r}^{\star2}}\sqrt{\frac{d}{n}}\lesssim\frac{1}{\kappa^{9/2}\mu^{3/2}r^{9/4}\kappa_{\omega}^{3/2}\log^{7/2}\left(n+d\right)},\qquad\frac{\omega_{\max}}{\sigma_{r}^{\star}}\sqrt{\frac{d}{np}}\lesssim\frac{1}{\kappa^{5}\mu^{3/2}r^{9/4}\kappa_{\omega}^{3/2}\log^{3}\left(n+d\right)}.
\]
Then it holds that 
\[
\mathbb{P}\left(\bm{U}_{l,\cdot}^{\star}\bm{R}^{\top}\in\mathsf{CR}_{U,l}^{1-\alpha}\right)=1-\alpha+O\bigg(\frac{1}{\sqrt{\log\left(n+d\right)}}\bigg)=1-\alpha+o\left(1\right).
\]
\end{lemma}\begin{proof}See Appendix~\ref{appendix:proof-lemma-pca-cr-validity}.\end{proof}With
this lemma in place, we have concluded the proof of Theorem \ref{thm:pca-cr-complete}.

\subsection{Entrywise distributional characterization for $\bm{S}^{\star}$ (Proof
of Theorem \ref{thm:ce-complete})}

The preceding distributional guarantees for principal subspace in
turn allow one to perform inference on the covariance matrix $\bm{S}^{\star}$
of the noiseless data vectors $\{\bm{x}_{j}\}_{1\leq j\leq n}$. In
this subsection, we demonstrate how to establish the advertised distributional
characterization for the matrix $\bm{S}$ returned by Algorithm~\ref{alg:CE-HeteroPCA-CI}.
Towards this end, we begin with the following decomposition 
\[
\bm{S}-\bm{S}^{\star}=\underbrace{\bm{S}-\bm{M}^{\natural}\bm{M}^{\natural\top}}_{\eqqcolon\,\bm{W}}+\underbrace{\bm{M}^{\natural}\bm{M}^{\natural\top}-\bm{S}^{\star}}_{\eqqcolon\,\bm{A}},
\]
and we shall write $\bm{W}=[W_{i,j}]_{1\leq i,j\leq d}$ and $\bm{A}=[A_{i,j}]$
from now on. In what follows, our proof consists of the following
main steps: 
\begin{enumerate}
\item Show that conditional on $\bm{F}$, each entry $W_{i,j}$ is approximately
a zero-mean Gaussian, whose variance concentrates around some deterministic
quantity $\widetilde{v}_{i,j}$. 
\item Show that each entry $A_{i,j}$ is approximately Gaussian with mean
zero and variance $\overline{v}_{i,j}$. 
\item Utilize the (near) independence of $W_{i,j}$ and $A_{i,j}$ to demonstrate
that $S_{i,j}-S_{i,j}^{\star}$ is approximately Gaussian with mean
zero and variance $\widetilde{v}_{i,j}+\overline{v}_{i,j}$. 
\end{enumerate}

\paragraph{Step 1: first- and second-order approximation of $\bm{W}$.}

To begin with, Lemma \ref{lemma:pca-2nd-error} allows one to derive
a reasonably accurate first- and second-order approximation for $\bm{W}=\bm{S}-\bm{M}^{\natural}\bm{M}^{\natural\top}$.
This is formalized in the following lemma, providing an explicit form
of this approximation and the goodness of the approximation.

\begin{lemma} \label{lemma:ce-2nd-error}Suppose that the assumptions
of Lemma \ref{lemma:pca-2nd-error} hold. Then one can write 
\[
\bm{W}=\bm{S}-\bm{M}^{\natural}\bm{M}^{\natural\top}=\bm{X}+\bm{\Phi},
\]
where 
\begin{equation}
\bm{X}\coloneqq\bm{E}\bm{M}^{\natural\top}+\bm{M}^{\natural}\bm{E}^{\top}+\mathcal{P}_{\mathsf{off}\text{-}\mathsf{diag}}\left(\bm{E}\bm{E}^{\top}\right)\bm{U}^{\natural}\bm{U}^{\natural\top}+\bm{U}^{\natural}\bm{U}^{\natural\top}\mathcal{P}_{\mathsf{off}\text{-}\mathsf{diag}}\left(\bm{E}\bm{E}^{\top}\right)\label{eq:defn-X-lemma-ce-2nd-error}
\end{equation}
and the residual matrix $\bm{\Phi}$ satisfies: conditional on $\bm{F}$
and on the $\sigma(\bm{F})$-measurable event $\mathcal{E}_{\mathsf{good}}$
(see Lemma~\ref{lemma:useful-property-good-event-hpca}), 
\begin{align}
\left|\Phi_{i,j}\right| & \lesssim\theta^{2}\left(\left\Vert \bm{U}_{i,\cdot}^{\star}\bm{\Sigma}^{\star}\right\Vert _{2}+\omega_{i}^{\star}\right)\left(\left\Vert \bm{U}_{j,\cdot}^{\star}\bm{\Sigma}^{\star}\right\Vert _{2}+\omega_{j}^{\star}\right)+\sigma_{1}^{\star2}\zeta_{\mathsf{2nd},i}\zeta_{\mathsf{2nd},j}\nonumber \\
 & \quad+\sigma_{1}^{\star2}\zeta_{\mathsf{2nd},i}\left(\left\Vert \bm{U}_{j,\cdot}^{\star}\right\Vert _{2}+\theta\frac{\omega_{j}^{\star}}{\sigma_{1}^{\star}}\right)+\zeta_{\mathsf{2nd},j}\sigma_{1}^{\star2}\left(\left\Vert \bm{U}_{i,\cdot}^{\star}\right\Vert _{2}+\theta\frac{\omega_{i}^{\star}}{\sigma_{1}^{\star}}\right)\eqqcolon\zeta_{i,j}\label{eq:zeta-i-j-def}
\end{align}
holds for any $i,j\in[d]$ with probability exceeding $1-O((n+d)^{-10})$.
Here, $\zeta_{\mathsf{2nd},i}$ and $\zeta_{\mathsf{2nd},j}$ are
defined in Lemma \ref{lemma:pca-2nd-error}, and $\theta$ is defined
in (\ref{eq:theta-definition}).\end{lemma}\begin{proof}See Appendix
\ref{appendix:proof-ce-2nd-error}.\end{proof}

\paragraph{Step 2: computing the entrywise variance of our approximation.}

We can check that the $(i,j)$-th entry of the matrix $\bm{X}$ (cf.~\eqref{eq:defn-X-lemma-ce-2nd-error})
is given by 
\begin{align}
X_{i,j} & =\left[\bm{E}\bm{M}^{\natural\top}+\bm{M}^{\natural}\bm{E}^{\top}+\mathcal{P}_{\mathsf{off}\text{-}\mathsf{diag}}\left(\bm{E}\bm{E}^{\top}\right)\bm{U}^{\natural}\bm{U}^{\natural\top}+\bm{U}^{\natural}\bm{U}^{\natural\top}\mathcal{P}_{\mathsf{off}\text{-}\mathsf{diag}}\left(\bm{E}\bm{E}^{\top}\right)\right]_{i,j}\nonumber \\
 & =\sum_{l=1}^{n}\left\{ M_{j,l}^{\natural}E_{i,l}+M_{i,l}^{\natural}E_{j,l}+E_{i,l}\left[\mathcal{P}_{-i,\cdot}\left(\bm{E}_{\cdot,l}\right)\right]^{\top}\bm{U}^{\star}\big(\bm{U}_{j,\cdot}^{\star}\big)^{\top}+E_{j,l}\left[\mathcal{P}_{-j,\cdot}\left(\bm{E}_{\cdot,l}\right)\right]^{\top}\bm{U}^{\star}\big(\bm{U}_{i,\cdot}^{\star}\big)^{\top}\right\} \nonumber \\
 & =\sum_{l=1}^{n}\bigg[M_{j,l}^{\natural}E_{i,l}+M_{i,l}^{\natural}E_{j,l}+\sum_{k:k\neq i}E_{i,l}E_{k,l}\left(\bm{U}_{k,\cdot}^{\star}\big(\bm{U}_{j,\cdot}^{\star}\big)^{\top}\right)+\sum_{k:k\neq j}E_{j,l}E_{k,l}\left(\bm{U}_{k,\cdot}^{\star}\big(\bm{U}_{i,\cdot}^{\star}\big)^{\top}\right)\bigg].\label{eq:defn-Xij-hpca}
\end{align}
Here, the penultimate step uses the fact that $\bm{U}^{\natural}=\bm{U}^{\star}\bm{Q}$
for some orthonormal matrix $\bm{Q}$ (see \eqref{eq:defn-Q-denoising-PCA}).
This allows one to calculate the variance of $X_{i,j}$ conditional
on $\bm{F}$: when $i\neq j$, 
\begin{align}
\mathsf{var}\left(X_{i,j}|\bm{F}\right) & =\sum_{l=1}^{n}M_{j,l}^{\natural2}\sigma_{i,l}^{2}+\sum_{l=1}^{n}M_{i,l}^{\natural2}\sigma_{j,l}^{2}+\sum_{l=1}^{n}\sum_{k:k\neq i}\sigma_{i,l}^{2}\sigma_{k,l}^{2}\left(\bm{U}_{k,\cdot}^{\star}\bm{U}_{j,\cdot}^{\star\top}\right)^{2}+\sum_{l=1}^{n}\sum_{k:k\neq j}\sigma_{j,l}^{2}\sigma_{k,l}^{2}\left(\bm{U}_{k,\cdot}^{\star}\bm{U}_{i,\cdot}^{\star\top}\right)^{2};\label{eq:var-Xij-F-hpca}
\end{align}
when $i=j$, we have 
\begin{align*}
\mathsf{var}\left(X_{i,i}|\bm{F}\right) & =4\sum_{l=1}^{n}M_{i,l}^{\natural2}\sigma_{i,l}^{2}+4\sum_{l=1}^{n}\sum_{k:k\neq i}\sigma_{i,l}^{2}\sigma_{k,l}^{2}\left(\bm{U}_{k,\cdot}^{\natural}\bm{U}_{i,\cdot}^{\natural\top}\right)^{2}\\
 & =4\sum_{l=1}^{n}M_{i,l}^{\natural2}\sigma_{i,l}^{2}+4\sum_{l=1}^{n}\sum_{k:k\neq i}\sigma_{i,l}^{2}\sigma_{k,l}^{2}\left(\bm{U}_{k,\cdot}^{\star}\bm{U}_{i,\cdot}^{\star\top}\right)^{2}.
\end{align*}
The next lemma states that $\mathsf{var}(X_{i,j}|\bm{F})$ concentrates
around some deterministic quantity $\widetilde{v}_{i,j}$ defined
as follows: \begin{subequations}\label{eq:defn-tilde-v-ijii} 
\begin{align}
\widetilde{v}_{i,j} & \coloneqq\frac{2\left(1-p\right)}{np}\left(S_{i,i}^{\star}S_{j,j}^{\star}+2S_{i,j}^{\star2}\right)+\frac{1}{np}\left(\omega_{i}^{\star2}S_{j,j}^{\star}+\omega_{j}^{\star2}S_{i,i}^{\star}\right)\nonumber \\
 & \quad+\frac{1}{np^{2}}\sum_{k=1}^{d}\left\{ \left[\omega_{i}^{\star2}+\left(1-p\right)S_{i,i}^{\star}\right]\left[\omega_{k}^{\star2}+\left(1-p\right)S_{k,k}^{\star}\right]+2\left(1-p\right)^{2}S_{i,k}^{\star2}\right\} \left(\bm{U}_{k,\cdot}^{\star}\bm{U}_{j,\cdot}^{\star\top}\right)^{2}\nonumber \\
 & \quad+\frac{1}{np^{2}}\sum_{k=1}^{d}\left\{ \left[\omega_{j}^{\star2}+\left(1-p\right)S_{j,j}^{\star}\right]\left[\omega_{k}^{\star2}+\left(1-p\right)S_{k,k}^{\star}\right]+2\left(1-p\right)^{2}S_{j,k}^{\star2}\right\} \left(\bm{U}_{k,\cdot}^{\star}\bm{U}_{i,\cdot}^{\star\top}\right)^{2}\label{eq:defn-tilde-vij}
\end{align}
for any $i\neq j$, and 
\begin{align}
\widetilde{v}_{i,i} & \coloneqq\frac{12\left(1-p\right)}{np}S_{i,i}^{\star2}+\frac{4}{np}\omega_{i}^{\star2}S_{i,i}^{\star}\nonumber \\
 & \quad+\frac{4}{np^{2}}\sum_{k=1}^{d}\left\{ \left[\omega_{i}^{\star2}+\left(1-p\right)S_{i,i}^{\star}\right]\left[\omega_{k}^{\star2}+\left(1-p\right)S_{k,k}^{\star}\right]+2\left(1-p\right)^{2}S_{i,k}^{\star2}\right\} \left(\bm{U}_{k,\cdot}^{\star}\bm{U}_{i,\cdot}^{\star\top}\right)^{2}\label{eq:defn-tilde-vii}
\end{align}
for any $i\in[d]$. \end{subequations}\begin{lemma}\label{lemma:ce-variance-concentration}Suppose
that $n\gg\log^{3}(n+d)$, and recall the definition of $\widetilde{v}_{i,j}$
in \eqref{eq:defn-tilde-v-ijii}. On the event $\mathcal{E}_{\mathsf{good}}$
(cf.~Lemma~\ref{lemma:useful-property-good-event-hpca}), we have
\begin{equation}
\mathsf{var}\left(X_{i,j}|\bm{F}\right)=\widetilde{v}_{i,j}+O\left(\sqrt{\frac{\log^{3}\left(n+d\right)}{n}}+\frac{\kappa\mu^{2}r^{2}+\kappa_{\omega}\mu r}{d}\right)\widetilde{v}_{i,j}\label{eq:Var-Xij-concentrate-vij}
\end{equation}
for any $i,j\in[d]$. In addition, for any $i,j\in[d]$ it holds that
\begin{equation}
\widetilde{v}_{i,j}\gtrsim\frac{1}{ndp^{2}\kappa\land np}\left\Vert \bm{U}_{i,\cdot}^{\star}\bm{\Sigma}^{\star}\right\Vert _{2}^{2}\left\Vert \bm{U}_{j,\cdot}^{\star}\bm{\Sigma}^{\star}\right\Vert _{2}^{2}+\left(\frac{\sigma_{r}^{\star2}}{ndp^{2}\land np}+\frac{\omega_{\min}^{2}}{np^{2}}\right)\left(\omega_{j}^{\star2}\left\Vert \bm{U}_{i,\cdot}^{\star}\right\Vert _{2}^{2}+\omega_{i}^{\star2}\left\Vert \bm{U}_{j,\cdot}^{\star}\right\Vert _{2}^{2}\right).\label{eq:vij-tilde-lower-bound-hpca}
\end{equation}
\end{lemma}\begin{proof}See Appendix \ref{appendix:proof-lemma-ce-variance-concentration}.\end{proof}

\paragraph{Step 3: establishing approximate Gaussianity of $W_{i,j}$.}

We are now ready to invoke the Berry-Esseen Theorem to show that $W_{i,j}$
is approximately Gaussian with mean zero and variance $\widetilde{v}_{i,j}$,
as stated in the next lemma.

\begin{lemma}\label{lemma:ce-normal-approx-1} Suppose that $d\gtrsim\kappa^{8}\mu^{3}r^{3}\kappa_{\omega}^{2}\log^{5}(n+d)$,
\begin{align*}
ndp^{2}\gtrsim\kappa^{10}\mu^{4}r^{4}\kappa_{\omega}^{2}\log^{9}\left(n+d\right), & \qquad np\gtrsim\kappa^{10}\mu^{3}r^{3}\kappa_{\omega}^{2}\log^{7}\left(n+d\right),\\
\frac{\omega_{\max}^{2}}{p\sigma_{r}^{\star2}}\sqrt{\frac{d}{n}}\lesssim\frac{1}{\sqrt{\kappa^{7}\mu^{2}r^{2}\kappa_{\omega}\log^{7}\left(n+d\right)}}, & \qquad\frac{\omega_{\max}}{\sigma_{r}^{\star}}\sqrt{\frac{d}{np}}\lesssim\frac{1}{\sqrt{\kappa^{8}\mu^{2}r^{2}\kappa_{\omega}\log^{6}\left(n+d\right)}},
\end{align*}
and 
\[
\left\Vert \bm{U}_{i,\cdot}^{\star}\right\Vert _{2}+\left\Vert \bm{U}_{j,\cdot}^{\star}\right\Vert _{2}\gtrsim\kappa r\log^{5/2}\left(n+d\right)\left[\frac{\omega_{\max}}{\sigma_{r}^{\star}}\sqrt{\frac{d}{np}}+\kappa_{\omega}\frac{\omega_{\max}^{2}}{p\sigma_{r}^{\star2}}\sqrt{\frac{d}{n}}+\frac{\kappa\mu r\kappa_{\omega}\log\left(n+d\right)}{\sqrt{ndp^{2}}}\right]\sqrt{\frac{1}{d}}.
\]
Then it holds that 
\[
\sup_{t\in\mathbb{R}}\left|\mathbb{P}\left(\big(\widetilde{v}_{i,j}\big)^{-1/2}W_{i,j}\leq t\,\big|\,\bm{F}\right)-\Phi\left(t\right)\right|\lesssim\frac{1}{\sqrt{\log\left(n+d\right)}}=o\left(1\right).
\]
\end{lemma}\begin{proof}See Appendix \ref{appendix:proof-ce-normal-approx-1}.\end{proof}

An important observation that one should bear in mind is that: conditional
on $\bm{F}$, the distribution of $W_{i,j}$ is approximately $\mathcal{N}(0,\widetilde{v}_{i,j})$,
where $\widetilde{v}_{i,j}$ does not depend on $\bm{F}$. This suggests
that $W_{i,j}$ is nearly independent of the $\sigma$-algebra $\sigma(\bm{F})$.

\paragraph{Step 4: establishing approximate Gaussianity of $A_{i,j}$.}

We now move on to the matrix $\bm{A}$. It follows from \eqref{eq:defn-M-Mnatural-E-denoising-PCA}
that 
\begin{align}
A_{i,j}=\left(\bm{M}^{\natural}\bm{M}^{\natural\top}-\bm{S}^{\star}\right)_{i,j} & =\frac{1}{n}\sum_{l=1}^{n}\left(\bm{U}_{i,\cdot}^{\star}\bm{\Sigma}^{\star}\bm{f}_{l}\right)\left(\bm{U}_{j,\cdot}^{\star}\bm{\Sigma}^{\star}\bm{f}_{l}\right)-\bm{U}_{i,\cdot}^{\star}\bm{\Sigma}^{\star2}\bm{U}_{j,\cdot}^{\star\top}.\label{eq:expression-Aij}
\end{align}
In view of the independence of $\{\bm{f}_{l}\}_{1\leq l\leq n}$,
the variance of $A_{i,j}$ can be calculated as follows 
\begin{align}
\overline{v}_{i,j} & \coloneqq\mathsf{var}\left(A_{i,j}\right)=\frac{1}{n}\mathsf{var}\left[\left(\bm{U}_{i,\cdot}^{\star}\bm{\Sigma}^{\star}\bm{f}_{1}\right)\left(\bm{U}_{j,\cdot}^{\star}\bm{\Sigma}^{\star}\bm{f}_{1}\right)\right]\nonumber \\
 & =\frac{1}{n}\left[\left\Vert \bm{U}_{i,\cdot}^{\star}\bm{\Sigma}^{\star}\right\Vert _{2}^{2}\left\Vert \bm{U}_{j,\cdot}^{\star}\bm{\Sigma}^{\star}\right\Vert _{2}^{2}+2\left(\bm{U}_{i,\cdot}^{\star}\bm{\Sigma}^{\star2}\bm{U}_{j,\cdot}^{\star\top}\right)^{2}\right]-\frac{1}{n}\left(\bm{U}_{i,\cdot}^{\star}\bm{\Sigma}^{\star2}\bm{U}_{j,\cdot}^{\star\top}\right)^{2}\nonumber \\
 & =\frac{1}{n}\left[\left\Vert \bm{U}_{i,\cdot}^{\star}\bm{\Sigma}^{\star}\right\Vert _{2}^{2}\left\Vert \bm{U}_{j,\cdot}^{\star}\bm{\Sigma}^{\star}\right\Vert _{2}^{2}+\left(\bm{U}_{i,\cdot}^{\star}\bm{\Sigma}^{\star2}\bm{U}_{j,\cdot}^{\star\top}\right)^{2}\right]=\frac{1}{n}\left(S_{i,i}^{\star}S_{j,j}^{\star}+S_{i,j}^{\star2}\right).\label{eq:defn-bar-vij}
\end{align}
With the assistance of the Berry-Esseen Theorem, the lemma below demonstrates
that $A_{i,j}$ is approximately Gaussian.

\begin{lemma}\label{lemma:ce-normal-approx-2}It holds that 
\[
\sup_{z\in\mathbb{R}}\left|\mathbb{P}\left(\big(\overline{v}_{i,j}\big)^{-1/2}A_{i,j}\leq z\right)-\Phi\left(z\right)\right|\lesssim\frac{1}{\sqrt{n}}.
\]
\end{lemma}\begin{proof}See Appendix \ref{appendix:ce-normal-approx-2}.\end{proof}

Note that $A_{i,j}$ is $\sigma(\bm{F})$-measurable. This fact taken
collectively with the near independence between $W_{i,j}$ and $\sigma(\bm{F})$
implies that $W_{i,j}$ and $A_{i,j}$ are nearly statistically independent.

\paragraph{Step 5: distributional characterization of $S_{i,j}-S_{i,j}^{\star}$.}

The approximate Gaussianity of $W_{i,j}$ and $A_{i,j}$, as well
as the near independence between them, leads to the conjecture that
$S_{i,j}-S_{i,j}^{\star}=W_{i,j}+A_{i,j}$ is approximately distributed
as $\mathcal{N}(0,v_{i,j}^{\star})$ with 
\[
v_{ij}^{\star}\coloneqq\widetilde{v}_{i,j}+\overline{v}_{i,j}.
\]
The following lemma rigorizes this conjecture, which in turn concludes
the proof of Theorem \ref{thm:ce-complete}.

\begin{lemma}\label{lemma:ce-normal-approx-all}Suppose that the
assumptions of Lemma \ref{lemma:ce-normal-approx-1} hold, and suppose
that $n\gtrsim\log(n+d)$. Then we have 
\[
\sup_{t\in\mathbb{R}}\left|\mathbb{P}\left(\left(S_{i,j}-S_{i,j}^{\star}\right)/\sqrt{v_{i,j}^{\star}}\leq t\right)-\Phi\left(t\right)\right|\lesssim\frac{1}{\sqrt{\log\left(n+d\right)}}=o\left(1\right).
\]
\end{lemma}\begin{proof}See Appendix \ref{appendix:ce-normal-approx-all}.\end{proof}

\subsection{Validity of confidence intervals (Proof of Theorem \ref{thm:ce-CI-complete})}

Armed with the above entrywise distributional theory for $\bm{S}-\bm{S}^{\star}$,
we hope to construct valid confidence interval for each $S_{i,j}^{\star}$
based on the estimate $\bm{S}$ returned by \textsf{HeteroPCA}. This
requires a faithful estimate of the variance $v_{i,j}^{\star}$. Towards
this, we define the following plug-in estimator: if $i\neq j$, let
\begin{subequations}\label{eq:defn-vij-vii-square} 
\begin{align}
v_{i,j} & =\frac{2-p}{np}S_{i,i}S_{j,j}+\frac{4-3p}{np}S_{i,j}^{2}+\frac{1}{np}\left(\omega_{i}^{2}S_{j,j}+\omega_{j}^{2}S_{i,i}\right)\nonumber \\
 & \quad+\frac{1}{np^{2}}\sum_{k=1}^{d}\left\{ \left[\omega_{i}^{2}+\left(1-p\right)S_{i,i}\right]\left[\omega_{k}^{2}+\left(1-p\right)S_{k,k}\right]+2\left(1-p\right)^{2}S_{i,k}^{2}\right\} \left(\bm{U}_{k,\cdot}\bm{U}_{j,\cdot}^{\top}\right)^{2}\nonumber \\
 & \quad+\frac{1}{np^{2}}\sum_{k=1}^{d}\left\{ \left[\omega_{j}^{2}+\left(1-p\right)S_{j,j}\right]\left[\omega_{k}^{2}+\left(1-p\right)S_{k,k}\right]+2\left(1-p\right)^{2}S_{j,k}^{2}\right\} \left(\bm{U}_{k,\cdot}\bm{U}_{i,\cdot}^{\top}\right)^{2};\label{eq:defn-vij-square}
\end{align}
otherwise, let 
\begin{align}
v_{i,i} & =\frac{12-9p}{np}S_{i,i}^{2}+\frac{4}{np}\omega_{i}^{2}S_{i,i}\nonumber \\
 & \quad+\frac{4}{np^{2}}\sum_{k=1}^{d}\left\{ \left[\omega_{i}^{2}+\left(1-p\right)S_{i,i}\right]\left[\omega_{k}^{2}+\left(1-p\right)S_{k,k}\right]+2\left(1-p\right)^{2}S_{i,k}^{2}\right\} \left(\bm{U}_{k,\cdot}\bm{U}_{i,\cdot}^{\top}\right)^{2}.\label{eq:defn-vii-square}
\end{align}
\end{subequations}

\paragraph{Step 1: faithfulness of the plug-in estimator.}

With the fine-grained estimation guarantees in Lemma \ref{lemma:pca-1st-err}
and Lemma \ref{lemma:pca-noise-level-est} in place, we can demonstrate
that $v_{i,j}$ is a reliable estimate of $v_{i,j}^{\star}$, as formally
stated below.

\begin{lemma}\label{lemma:ce-var-est}Suppose that the conditions
of Lemma \ref{lemma:ce-normal-approx-1} hold. For any $\delta\in(0,1)$,
we further assume that $n\gtrsim\delta^{-2}\kappa^{7}\mu^{2}r^{3}\kappa_{\omega}^{2}\log(n+d)$,
$d\gtrsim\kappa^{2}\mu\log(n+d)$, 
\begin{align*}
ndp^{2}\gtrsim\delta^{-2}\kappa^{6}\mu^{4}r^{6}\kappa_{\omega}^{2}\log^{5}\left(n+d\right), & \qquad np\gtrsim\delta^{-2}\kappa^{6}\mu^{3}r^{5}\kappa_{\omega}^{2}\log^{3}\left(n+d\right),\\
\frac{\omega_{\max}^{2}}{p\sigma_{r}^{\star2}}\sqrt{\frac{d}{n}}\lesssim\frac{\delta}{\kappa^{2}\mu r^{2}\kappa_{\omega}\log^{3/2}\left(n+d\right)}, & \qquad\frac{\omega_{\max}}{\sigma_{r}^{\star}}\sqrt{\frac{d}{np}}\lesssim\frac{\delta}{\kappa^{5/2}\mu r^{2}\kappa_{\omega}\log\left(n+d\right)},
\end{align*}
and 
\[
\left\Vert \bm{U}_{i,\cdot}^{\star}\right\Vert _{2}+\left\Vert \bm{U}_{j,\cdot}^{\star}\right\Vert _{2}\gtrsim\delta^{-1}\kappa\mu r^{2}\kappa_{\omega}\log\left(n+d\right)\left[\frac{\kappa\mu r\log\left(n+d\right)}{\sqrt{nd}p}+\frac{\omega_{\max}^{2}}{p\sigma_{r}^{\star2}}\sqrt{\frac{d}{n}}+\frac{\omega_{\max}}{\sigma_{r}^{\star}}\sqrt{\frac{d}{np}}\right]\sqrt{\frac{r}{d}}.
\]
Then with probability exceeding $1-O((n+d)^{-10})$, the quantity
$v_{i,j}$ defined in \eqref{eq:defn-vij-vii-square} obeys 
\begin{align*}
\left|v_{i,j}^{\star}-v_{i,j}\right| & \lesssim\delta v_{i,j}^{\star}.
\end{align*}
\end{lemma}\begin{proof}See Appendix \ref{appendix:proof-lemma-ce-var-est}.\end{proof}

\paragraph{Step 2: validity of the constructed confidence intervals.}

Armed with the Gaussian approximation in Lemma \ref{lemma:ce-normal-approx-all}
as well as the faithfulness of $v_{i,j}$ as an estimate of $v_{i,j}^{\star}$
in the previous lemma, we show in the next lemma that the confidence
interval constructed in Algorithm \ref{alg:CE-HeteroPCA-CI} is valid
and nearly accurate.

\begin{lemma}\label{lemma:ce-CI-validity}Suppose that the conditions
of Lemma \ref{lemma:ce-normal-approx-1} hold. Further suppose that
$n\gtrsim\kappa^{9}\mu^{3}r^{4}\kappa_{\omega}^{3}\log^{4}(n+d)$,
$d\gtrsim\kappa^{2}\mu\log(n+d)$, 
\begin{align*}
ndp^{2}\gtrsim\kappa^{8}\mu^{5}r^{7}\kappa_{\omega}^{3}\log^{8}\left(n+d\right), & \qquad np\gtrsim\kappa^{8}\mu^{4}r^{6}\kappa_{\omega}^{3}\log^{6}\left(n+d\right),\\
\frac{\omega_{\max}^{2}}{p\sigma_{r}^{\star2}}\sqrt{\frac{d}{n}}\lesssim\frac{1}{\kappa^{3}\mu^{3/2}r^{5/2}\kappa_{\omega}^{3/2}\log^{3}\left(n+d\right)}, & \qquad\frac{\omega_{\max}}{\sigma_{r}^{\star}}\sqrt{\frac{d}{np}}\lesssim\frac{1}{\kappa^{7/2}\mu^{3/2}r^{5/2}\kappa_{\omega}^{3/2}\log^{5/2}\left(n+d\right)},
\end{align*}
and 
\[
\left\Vert \bm{U}_{i,\cdot}^{\star}\right\Vert _{2}+\left\Vert \bm{U}_{j,\cdot}^{\star}\right\Vert _{2}\gtrsim\kappa^{2}\mu^{3/2}r^{5/2}\kappa_{\omega}^{3/2}\log^{5/2}\left(n+d\right)\left[\frac{\kappa\mu r\log\left(n+d\right)}{\sqrt{nd}p}+\frac{\omega_{\max}^{2}}{p\sigma_{r}^{\star2}}\sqrt{\frac{d}{n}}+\frac{\omega_{\max}}{\sigma_{r}^{\star}}\sqrt{\frac{d}{np}}\right]\sqrt{\frac{r}{d}}.
\]
Then the confidence region $\mathsf{CI}_{i,j}^{1-\alpha}$ returned
from Algorithm~\ref{alg:CE-HeteroPCA-CI} satisfies 
\[
\mathbb{P}\left(S_{i,j}^{\star}\in\mathsf{CI}_{i,j}^{1-\alpha}\right)=1-\alpha+O\left(\frac{1}{\sqrt{\log\left(n+d\right)}}\right)=1-\alpha+o\left(1\right).
\]
\end{lemma}\begin{proof}See Appendix \ref{appendix:proof-ce-CI-validity}.\end{proof} 

\section{Auxiliary lemmas: the approach based on \textsf{HeteroPCA \label{appendix:hpca-auxiliary-lemmas}}}

\subsection{Proof of Lemma \ref{lemma:useful-property-good-event-hpca} \label{appendix:good-event}}

In this section we will establish a couple of useful properties that
occur with high probability. 
\begin{enumerate}
	\item In view of standard results on Gaussian random matrices, we have 
	\begin{equation}
		\left\Vert \frac{1}{n}\bm{F}\bm{F}^{\top}-\bm{I}_{r}\right\Vert \lesssim\sqrt{\frac{r}{n}}+\sqrt{\frac{\log\left(n+d\right)}{n}}+\frac{r}{n}+\frac{\log\left(n+d\right)}{n}\lesssim\sqrt{\frac{r+\log\left(n+d\right)}{n}}\label{eq:good-event-gaussian-concentration}
	\end{equation}
	with probability exceeding $1-O((n+d)^{-100})$, provided that $n\gtrsim r+\log(n+d)$.
	As an immediate consequence, it is seen from the definition \eqref{eq:defn-M-Mnatural-E-denoising-PCA}
	that 
	\begin{align}
		\left\Vert \bm{M}^{\natural}\bm{M}^{\natural\top}-\bm{S}^{\star}\right\Vert  & =\left\Vert \bm{U}^{\star}\bm{\Sigma}^{\star}\left(\frac{1}{n}\bm{F}\bm{F}^{\top}-\bm{I}_{r}\right)\bm{\Sigma}^{\star}\bm{U}^{\star\top}\right\Vert \leq\left\Vert \bm{U}^{\star}\right\Vert ^{2}\left\Vert \bm{\Sigma}^{\star}\right\Vert ^{2}\left\Vert \frac{1}{n}\bm{F}\bm{F}^{\top}-\bm{I}_{r}\right\Vert \nonumber \\
		& \lesssim\sqrt{\frac{r+\log\left(n+d\right)}{n}}\sigma_{1}^{\star2}.\label{eq:good-event-AA-S}
	\end{align}
	Weyl's inequality then tells us that 
	\begin{equation}
		\left\Vert \bm{\Sigma}^{\natural2}-\bm{\Sigma}^{\star2}\right\Vert \leq\left\Vert \bm{M}^{\natural}\bm{M}^{\natural\top}-\bm{S}^{\star}\right\Vert \lesssim\sqrt{\frac{r+\log\left(n+d\right)}{n}}\sigma_{1}^{\star2},\label{eq:good-event-Sigma-2-diff}
	\end{equation}
	thus indicating that 
	\begin{equation}
		\sigma_{r}^{\natural}\asymp\sigma_{r}^{\star}\qquad\text{and}\qquad\sigma_{1}^{\natural}\asymp\sigma_{1}^{\star}\label{eq:good-event-singular-value-same-order}
	\end{equation}
	as long as $n\gg\frac{\sigma_{1}^{\star4}}{\sigma_{r}^{\star4}}\big(r+\log(n+d)\big)=\kappa^{2}\big(r+\log(n+d)\big)$
	(see the definition of $\kappa$ in \eqref{eq:defn-condition-number}).
	Consequently, 
	\begin{equation}
		\left\Vert \bm{\Sigma}^{\natural}-\bm{\Sigma}^{\star}\right\Vert \leq\frac{\left\Vert \bm{\Sigma}^{\natural2}-\bm{\Sigma}^{\star2}\right\Vert }{\sigma_{r}^{\star}}\lesssim\sqrt{\frac{r+\log\left(n+d\right)}{n}}\frac{\sigma_{1}^{\star2}}{\sigma_{r}^{\star}}\asymp\kappa\sqrt{\frac{\big(r+\log\left(n+d\right)\big)}{n}}\sigma_{r}^{\star}.\label{eq:good-event-Sigma-2-diff-1}
	\end{equation}
	\item Given the calculation \eqref{eq:sigma-ij-square-denoising-PCA}, let
	us define 
	\begin{align*}
		\sigma^{2} & \coloneqq\max_{i\in[d],j\in[n]}\sigma_{i,j}^{2}\asymp\frac{1-p}{np}\max_{i\in[d],j\in[n]}\left(\bm{U}_{i,\cdot}^{\star}\bm{\Sigma}^{\star}\bm{f}_{j}\right)^{2}+\frac{\omega_{i}^{\star2}}{np},
	\end{align*}
	and
	\[
	\sigma_{i}^{2}\coloneqq\max_{j\in[n]}\sigma_{i,j}^{2}\asymp\frac{1-p}{np}\max_{j\in[n]}\left(\bm{U}_{i,\cdot}^{\star}\bm{\Sigma}^{\star}\bm{f}_{j}\right)^{2}+\frac{\omega_{i}^{\star2}}{np}
	\]
	for each $i\in[d]$. Note that for each $i\in[d]$ and $j\in[n]$,
	one has $\bm{U}_{i,\cdot}^{\star}\bm{\Sigma}^{\star}\bm{f}_{j}\sim\mathcal{N}(0,\Vert\bm{U}_{i,\cdot}^{\star}\bm{\Sigma}^{\star}\Vert_{2}^{2})$,
	thus revealing that 
	\begin{equation}
		\max_{j\in[n]}\left|\bm{U}_{i,\cdot}^{\star}\bm{\Sigma}^{\star}\bm{f}_{j}\right|\lesssim\left\Vert \bm{U}_{i,\cdot}^{\star}\bm{\Sigma}^{\star}\right\Vert _{2}\sqrt{\log\left(n+d\right)}\label{eq:good-event-Ui-Sigma-f}
	\end{equation}
	with probability exceeding $1-O((n+d)^{-100})$. As a result, taking
	the union bound gives 
	\begin{align*}
		\sigma^{2} & \lesssim\frac{\left\Vert \bm{U}^{\star}\bm{\Sigma}^{\star}\right\Vert _{2,\infty}^{2}\log\left(n+d\right)+\omega_{\max}^{2}}{np}\lesssim\frac{\mu r\log\left(n+d\right)}{ndp}\sigma_{1}^{\star2}+\frac{\omega_{\max}^{2}}{np}\eqqcolon\sigma_{\mathsf{ub}}^{2}
	\end{align*}
	and
	\[
	\sigma_{i}^{2}\lesssim\frac{\log\left(n+d\right)}{np}\left\Vert \bm{U}_{i,\cdot}^{\star}\bm{\Sigma}^{\star}\right\Vert _{2}^{2}+\frac{\omega_{i}^{\star2}}{np}\coloneqq\sigma_{\mathsf{ub},i}^{2}.
	\]
	\item In addition, it follows from the expression \eqref{eq:E-expression-denoising}
	and the property \eqref{eq:good-event-Ui-Sigma-f} that 
	\begin{align*}
		\max_{i\in[d],j\in[n]}\left|E_{i,j}\right| & \leq\frac{1}{\sqrt{n}p}\max_{i\in[d],j\in[n]}\left|\bm{U}_{i,\cdot}^{\star}\bm{\Sigma}^{\star}\bm{f}_{j}\right|+\frac{1}{\sqrt{n}p}\max_{i\in[d],j\in[n]}\left|N_{i,j}\right|\\
		& \lesssim\frac{1}{\sqrt{n}p}\left\Vert \bm{U}^{\star}\bm{\Sigma}^{\star}\right\Vert _{2,\infty}\sqrt{\log\left(n+d\right)}+\frac{\omega_{\max}\sqrt{\log\left(n+d\right)}}{\sqrt{n}p}\\
		& \lesssim\frac{1}{p}\sqrt{\frac{\mu r\log\left(n+d\right)}{nd}}\sigma_{1}^{\star}+\frac{\omega_{\max}}{p}\sqrt{\frac{\log\left(n+d\right)}{n}}
	\end{align*}
	occurs with probability exceeding $1-O((n+d)^{-100})$. Therefore,
	we shall take 
	\[
	B\coloneqq\frac{1}{p}\sqrt{\frac{\mu r\log\left(n+d\right)}{nd}}\sigma_{1}^{\star}+\frac{\omega_{\max}}{p}\sqrt{\frac{\log\left(n+d\right)}{n}}\asymp\sigma_{\mathsf{ub}}\sqrt{\frac{\log\left(n+d\right)}{p}}
	\]
	as an upper bound for $\{\vert E_{i,j}\vert:i\in[d],j\in[n]\}$. Similarly,
	for each $i\in[d]$, we can take
	\[
	B_{i}\coloneqq\frac{1}{p}\sqrt{\frac{\log\left(n+d\right)}{n}}\left\Vert \bm{U}_{i,\cdot}^{\star}\bm{\Sigma}^{\star}\right\Vert _{2}+\frac{\omega_{l}^{\star}}{p}\sqrt{\frac{\log\left(n+d\right)}{n}}
	\]
	as an upper bound for $\{\vert E_{i,j}\vert:j\in[n]\}$. 
	\item Recall that the top-$r$ eigen-decomposition of $\bm{M}^{\natural}\bm{M}^{\natural\top}$
	and $\bm{S}^{\star}$ are denoted by $\bm{U}^{\natural}\bm{\Sigma}^{\natural2}\bm{U}^{\natural\top}$
	and $\bm{U}^{\star}\bm{\Sigma}^{\star}\bm{U}^{\star\top}$, respectively,
	and that $\bm{Q}$ is a rotation matrix such that $\bm{U}^{\natural}=\bm{U}^{\star}\bm{Q}$.
	Therefore, the matrix $\bm{J}$ defined in \eqref{eq:defn-J-denoising-PCA}
	obeys 
	\begin{align*}
		\left\Vert \bm{Q}-\bm{J}\right\Vert  & =\left\Vert \bm{Q}-\bm{\Sigma}^{\star}\bm{Q}\left(\bm{\Sigma}^{\natural}\right)^{-1}\right\Vert \leq\frac{1}{\sigma_{r}^{\natural}}\left\Vert \bm{Q}\bm{\Sigma}^{\natural}-\bm{\Sigma}^{\star}\bm{Q}\right\Vert =\frac{1}{\sigma_{r}^{\natural}}\left\Vert \bm{Q}\bm{\Sigma}^{\natural}\bm{Q}^{\top}-\bm{\Sigma}^{\star}\right\Vert \\
		& =\frac{1}{\sigma_{r}^{\natural}}\left\Vert \bm{U}^{\star}\left(\bm{Q}\bm{\Sigma}^{\natural}\bm{Q}^{\top}-\bm{\Sigma}^{\star}\right)\bm{U}^{\star\top}\right\Vert =\frac{1}{\sigma_{r}^{\natural}}\left\Vert \bm{U}^{\natural}\bm{\Sigma}^{\natural}\bm{U}^{\natural\top}-\bm{U}^{\star}\bm{\Sigma}^{\star}\bm{U}^{\star\top}\right\Vert .
	\end{align*}
	Invoke the perturbation bound for matrix square roots \citep[Lemma 2.1]{MR1176461}
	to derive 
	\begin{align*}
		\left\Vert \bm{U}^{\natural}\bm{\Sigma}^{\natural}\bm{U}^{\natural\top}-\bm{U}^{\star}\bm{\Sigma}^{\star}\bm{U}^{\star\top}\right\Vert  & \lesssim\frac{1}{\sigma_{r}^{\natural}+\sigma_{r}^{\star}}\left\Vert \bm{U}^{\natural}\bm{\Sigma}^{\natural2}\bm{U}^{\natural\top}-\bm{U}^{\star}\bm{\Sigma}^{\star2}\bm{U}^{\star\top}\right\Vert \\
		& \asymp\frac{1}{\sigma_{r}^{\star}}\left\Vert \bm{M}^{\natural}\bm{M}^{\natural\top}-\bm{S}^{\star}\right\Vert ,
	\end{align*}
	where the last step follows from (\ref{eq:good-event-singular-value-same-order}).
	To summarize, with probability exceeding $1-O((n+d)^{-100})$ one
	has 
	\begin{equation}
		\left\Vert \bm{Q}-\bm{J}\right\Vert \lesssim\frac{1}{\sigma_{r}^{\star}}\left\Vert \bm{Q}\bm{\Sigma}^{\natural}-\bm{\Sigma}^{\star}\bm{Q}\right\Vert \lesssim\frac{1}{\sigma_{r}^{\star2}}\left\Vert \bm{M}^{\natural}\bm{M}^{\natural\top}-\bm{S}^{\star}\right\Vert \lesssim\kappa\sqrt{\frac{r+\log\left(n+d\right)}{n}},\label{eq:good-event-J-Q}
	\end{equation}
	where we have made used of the properties (\ref{eq:good-event-singular-value-same-order})
	and (\ref{eq:good-event-AA-S}). 
	\item In view of (\ref{eq:good-event-singular-value-same-order}), we know
	that with probability exceeding $1-O((n+d)^{-100})$, the conditional
	number of $\bm{M}^{\natural}$ satisfies 
	\begin{equation}
		\kappa^{\natural}\asymp\sqrt{\kappa}.\label{eq:good-event-kappa}
	\end{equation}
	Recalling that $\bm{U}^{\natural}=\bm{U}^{\star}\bm{Q}$, we can see
	from the incoherence assumption that 
	\begin{equation}
		\left\Vert \bm{U}^{\natural}\right\Vert _{2,\infty}=\left\Vert \bm{U}^{\star}\bm{Q}\right\Vert _{2,\infty}=\left\Vert \bm{U}^{\star}\right\Vert _{2,\infty}\leq\sqrt{\frac{\mu r}{d}}.\label{eq:Unatural-UB-1234}
	\end{equation}
	In addition, it is readily seen from \eqref{eq:defn-J-denoising}
	that 
	\begin{equation}
		\left\Vert \bm{V}^{\natural}\right\Vert _{2,\infty}=\left\Vert \frac{1}{\sqrt{n}}\bm{F}^{\top}\bm{J}\right\Vert _{2,\infty}\overset{\text{(i)}}{\lesssim}\sqrt{\frac{1}{n}}\left\Vert \bm{F}^{\top}\right\Vert _{2,\infty}\overset{\text{(ii)}}{\lesssim}\sqrt{\frac{r\log\left(n+d\right)}{n}}\label{eq:good-event-V-natural-2-infty}
	\end{equation}
	with probability exceeding $1-O((n+d)^{-100})$. Here, (i) holds since,
	according to (\ref{eq:good-event-J-Q}), 
	\begin{equation}
		\left\Vert \bm{J}\right\Vert \leq\left\Vert \bm{Q}\right\Vert +\left\Vert \bm{J}-\bm{Q}\right\Vert \leq1+\kappa\sqrt{\frac{r+\log\left(n+d\right)}{n}}\leq2\label{eq:good-event-J-bounded}
	\end{equation}
	holds as long as $n\gg\kappa^{2}(r+\log(n+d))$; and (ii) follows
	from the standard Gaussian concentration inequality. Combine \eqref{eq:Unatural-UB-1234}
	, \eqref{eq:good-event-V-natural-2-infty} and
	\begin{align*}
		\left\Vert \bm{M}^{\natural}\right\Vert _{\infty} & =\max_{i,j}\left|\bm{U}_{i,\cdot}^{\natural}\bm{\Sigma}^{\natural}\bm{V}_{j,\cdot}^{\natural\top}\right|\leq\sigma_{1}^{\natural}\left\Vert \bm{U}^{\natural}\right\Vert _{2,\infty}\left\Vert \bm{V}^{\natural}\right\Vert _{2,\infty}\overset{\text{(i)}}{\lesssim}\sqrt{\frac{\mu\log\left(n+d\right)}{nd}}\sigma_{1}^{\natural}\sqrt{r}\\
		& \lesssim\sqrt{\frac{\mu\log\left(n+d\right)}{nd}}\kappa^{\natural}\left\Vert \bm{M}^{\natural}\right\Vert _{\mathrm{F}}\overset{\text{(ii)}}{\asymp}\sqrt{\frac{\kappa\mu\log\left(n+d\right)}{nd}}\left\Vert \bm{M}^{\natural}\right\Vert _{\mathrm{F}}.
	\end{align*}
	to reach 
	\[
	\mu^{\natural}\lesssim\kappa\mu\log\left(n+d\right)
	\]
	with probability exceeding $1-O((n+d)^{-100})$. Here, (i) follows
	from (\ref{eq:good-event-V-natural-2-infty}), whereas (ii) arises
	from (\ref{eq:good-event-kappa}).
	\item Invoke \citet[Lemma 14]{chen2018gradient} to show that, for all $l\in[d]$,
	\begin{equation}
		\biggl\Vert\frac{1}{n}\sum_{j=1}^{n}\left(\bm{U}_{l,\cdot}^{\star}\bm{\Sigma}^{\star}\bm{f}_{j}\right)^{2}\bm{f}_{j}\bm{f}_{j}^{\top}-\left\Vert \bm{U}_{l,\cdot}^{\star}\bm{\Sigma}^{\star}\right\Vert _{2}^{2}\bm{I}_{r}-2\bm{\Sigma}^{\star}\bm{U}_{l,\cdot}^{\star\top}\bm{U}_{l,\cdot}^{\star}\bm{\Sigma}^{\star}\biggr\Vert\lesssim\sqrt{\frac{r\log^{3}\left(n+d\right)}{n}}\left\Vert \bm{U}_{l,\cdot}^{\star}\bm{\Sigma}^{\star}\right\Vert _{2}^{2}\label{eq:good-event-concentration-1}
	\end{equation}
	holds with probability exceeding $1-O((n+d)^{-100})$, provided that
	$n\gg r\log^{3}(n+d)$. In addition, we make note of the following
	lemmas.\begin{lemma}\label{lemma:gaussian-concentration-1}Assume
		that $n\gg\log(n+d)$. For any fixed vector $\bm{u}\in\mathbb{R}^{r}$,
		with probability exceeding $1-O((n+d)^{-100})$ we have 
		\[
		\left|\frac{1}{n}\sum_{j=1}^{n}\left(\bm{u}^{\top}\bm{f}_{j}\right)^{2}-\left\Vert \bm{u}\right\Vert _{2}^{2}\right|\lesssim\sqrt{\frac{\log\left(n+d\right)}{n}}\left\Vert \bm{u}\right\Vert _{2}^{2}.
		\]
	\end{lemma}\begin{lemma} \label{lemma:gaussian-concentration-2}Assume
		that $n\gg\log(n+d)$. For any fixed unit vectors $\bm{u},\bm{v}\in\mathbb{R}^{r}$,
		with probability exceeding $1-O((n+d)^{-100})$ we have 
		\[
		\left|\frac{1}{n}\sum_{j=1}^{n}\left(\bm{u}^{\top}\bm{f}_{j}\right)^{2}\left(\bm{v}^{\top}\bm{f}_{j}\right)^{2}-\left\Vert \bm{u}\right\Vert _{2}^{2}\left\Vert \bm{v}\right\Vert _{2}^{2}-2\left(\bm{u}^{\top}\bm{v}\right)^{2}\right|\lesssim\sqrt{\frac{\log\left(n+d\right)^{3}}{n}}\left\Vert \bm{u}\right\Vert _{2}^{2}\left\Vert \bm{v}\right\Vert _{2}^{2}.
		\]
	\end{lemma}In view of Lemma \ref{lemma:gaussian-concentration-1}
	and Lemma \ref{lemma:gaussian-concentration-2}, we know that for
	each $i,l\in[d]$\begin{subequations} 
		\begin{equation}
			\left|\frac{1}{n}\sum_{j=1}^{n}\left(\bm{U}_{l,\cdot}^{\star}\bm{\Sigma}^{\star}\bm{f}_{j}\right)^{2}\left(\bm{U}_{i,\cdot}^{\star}\bm{\Sigma}^{\star}\bm{f}_{j}\right)^{2}-\left[\left\Vert \bm{U}_{l,\cdot}^{\star}\bm{\Sigma}^{\star}\right\Vert _{2}^{2}\left\Vert \bm{U}_{i,\cdot}^{\star}\bm{\Sigma}^{\star}\right\Vert _{2}^{2}+2S_{i,l}^{\star2}\right]\right|\lesssim\sqrt{\frac{\log\left(n+d\right)^{3}}{n}}\left\Vert \bm{U}_{l,\cdot}^{\star}\bm{\Sigma}^{\star}\right\Vert _{2}^{2}\left\Vert \bm{U}_{i,\cdot}^{\star}\bm{\Sigma}^{\star}\right\Vert _{2}^{2}\label{eq:good-event-concentration-2}
		\end{equation}
		\begin{equation}
			\text{and}\qquad\left|\frac{1}{n}\sum_{j=1}^{n}\left(\bm{U}_{i,\cdot}^{\star}\bm{\Sigma}^{\star}\bm{f}_{j}\right)^{2}-\left\Vert \bm{U}_{i,\cdot}^{\star}\bm{\Sigma}^{\star}\right\Vert _{2}^{2}\right|\lesssim\sqrt{\frac{\log\left(n+d\right)}{n}}\left\Vert \bm{U}_{i,\cdot}^{\star}\bm{\Sigma}^{\star}\right\Vert _{2}^{2}\label{eq:good-event-concentration-3}
		\end{equation}
	\end{subequations}hold with probability exceeding $1-O((n+d)^{-100})$. 
	\item With the above properties in place, we can now formally define the
	``good'' event $\mathcal{E}_{\mathsf{good}}$ as follows 
	\[
	\mathcal{E}_{\mathsf{good}}\coloneqq\left\{ \text{All equations from (\ref{eq:good-event-gaussian-concentration}) to (\ref{eq:good-event-concentration-3}) hold}\right\} .
	\]
	It is immediately seen from the above analysis that 
	\[
	\mathbb{P}\left(\mathcal{E}_{\mathsf{good}}\right)\geq1-O\left(\left(n+d\right)^{-100}\right).
	\]
	By construction, it is self-evident that the event $\mathcal{E}_{\mathsf{good}}$
	is $\sigma(\bm{F})$-measurable. 
\end{enumerate}
\begin{proof}[Proof of Lemma \ref{lemma:gaussian-concentration-1}]Recognizing
	that $\bm{u}^{\top}\bm{f}_{j}\overset{\text{i.i.d.}}{\sim}\mathcal{N}(0,\Vert\bm{u}\Vert_{2}^{2})$,
	we have 
	\[
	\frac{1}{n}\sum_{j=1}^{n}\left(\bm{u}^{\top}\bm{f}_{j}\right)^{2}/\left\Vert \bm{u}\right\Vert _{2}^{2}\sim\chi_{n}^{2},
	\]
	where $\chi_{n}^{2}$ denotes the chi-square distribution with $n$
	degrees of freedom. One can then apply the tail bound for chi-square
	random variables (see, e.g., \citet[Example 2.5]{wainwright2019high})
	to establish the desired result. \end{proof}

\begin{proof}[Proof of Lemma \ref{lemma:gaussian-concentration-2}]
	Given that $\bm{u}^{\top}\bm{f}_{j}\sim\mathcal{N}(0,\Vert\bm{u}\Vert_{2}^{2})$
	for all $j\in[n]$, we have, with probability exceeding $1-O((n+d)^{-100})$,
	that 
	\[
	\max_{j\in\left[n\right]}\left|\bm{u}^{\top}\bm{f}_{j}\right|\leq C_{1}\left\Vert \bm{u}\right\Vert _{2}\sqrt{\log\left(n+d\right)}
	\]
	for some sufficiently large constant $C_{1}>0$. For each $j\in[n]$,
	let $X_{j}=(\bm{u}^{\top}\bm{f}_{j})^{2}$ and $Y_{j}=(\bm{v}^{\top}\bm{f}_{j})^{2}$
	and define the event 
	\[
	\mathcal{A}_{j}\coloneqq\left\{ \left|\bm{u}^{\top}\bm{f}_{j}\right|\leq C_{1}\left\Vert \bm{u}\right\Vert _{2}\sqrt{\log\left(n+d\right)}\right\} .
	\]
	Then with probability exceeding $1-O((n+d)^{-100})$, it holds that
	\[
	\frac{1}{n}\sum_{j=1}^{n}X_{j}Y_{j}=\frac{1}{n}\sum_{j=1}^{n}X_{j}Y_{j}\ind_{\mathcal{A}_{j}},
	\]
	which motivates us to decompose 
	\begin{align*}
		\left|\frac{1}{n}\sum_{j=1}^{n}\left(X_{j}Y_{j}-\mathbb{E}\left[X_{j}Y_{j}\right]\right)\right| & \leq\underbrace{\left|\frac{1}{n}\sum_{j=1}^{n}\left\{ X_{j}Y_{j}\ind_{\mathcal{A}_{j}}-\mathbb{E}\left[X_{j}Y_{j}\ind_{\mathcal{A}_{j}}\right]\right\} \right|}_{\eqqcolon\alpha_{1}}+\underbrace{\left|\mathbb{E}\left[X_{1}Y_{1}\ind_{\mathcal{A}_{1}^{\mathsf{c}}}\right]\right|}_{\eqqcolon\alpha_{2}}.
	\end{align*}
	
	\begin{itemize}
		\item Let us first bound $\alpha_{2}$. It is straightforward to derive
		that 
		\begin{align*}
			\alpha_{2} & =\mathbb{E}\left[\left(\bm{u}^{\top}\bm{f}_{1}\right)^{2}\left(\bm{v}^{\top}\bm{f}_{1}\right)^{2}\ind_{\left|\bm{u}^{\top}\bm{f}_{1}\right|>C_{1}\left\Vert \bm{u}\right\Vert _{2}\sqrt{\log\left(n+d\right)}}\right]\\
			& \overset{\text{(i)}}{\leq}\left(\mathbb{E}\left[\left(\bm{u}^{\top}\bm{f}_{1}\right)^{6}\right]\right)^{\frac{1}{3}}\left(\mathbb{E}\left[\left(\bm{v}^{\top}\bm{f}_{1}\right)^{6}\right]\right)^{\frac{1}{3}}\left[\mathbb{P}\left(\left|\bm{u}^{\top}\bm{f}_{1}\right|>C_{1}\left\Vert \bm{u}\right\Vert _{2}\sqrt{\log\left(n+d\right)}\right)\right]^{1/3}\\
			& \lesssim\left\Vert \bm{u}\right\Vert _{2}^{2}\left\Vert \bm{v}\right\Vert _{2}^{2}\exp\left(-\frac{C_{1}^{2}\log\left(n+d\right)}{6}\right)\overset{\text{(ii)}}{\lesssim}\left\Vert \bm{u}\right\Vert _{2}^{2}\left\Vert \bm{v}\right\Vert _{2}^{2}\left(n+d\right)^{-100},
		\end{align*}
		where (i) comes from H{ö}lder's inequality, and (ii) holds for $C_{1}$
		large enough. 
		\item Next, we shall apply the Bernstein inequality \citep[Theorem 2.8.2]{vershynin2016high}
		to bound $\alpha_{1}$. Note that for each $j\in[n]$ 
		\begin{align}
			\left\Vert X_{j}Y_{j}\ind_{\mathcal{A}_{j}}-\mathbb{E}\left[X_{j}Y_{j}\ind_{\mathcal{A}_{j}}\right]\right\Vert _{\psi_{1}} & \leq\left\Vert X_{j}Y_{j}\ind_{\mathcal{A}_{j}}\right\Vert _{\psi_{1}}+\left\Vert \mathbb{E}\left[X_{j}Y_{j}\ind_{\mathcal{A}_{j}}\right]\right\Vert _{\psi_{1}},\label{eq:XY-decompose-123}
		\end{align}
		where $\|\cdot\|_{\psi_{1}}$ denotes the sub-exponential norm \citep{Vershynin2012}.
		The first term of \eqref{eq:XY-decompose-123} obeys 
		\begin{align*}
			\left\Vert X_{j}Y_{j}\ind_{\mathcal{A}_{j}}\right\Vert _{\psi_{1}} & \leq C_{1}^{2}\left\Vert \bm{u}\right\Vert _{2}^{2}\log\left(n+d\right)\left\Vert Y_{j}\right\Vert _{\psi_{1}}\lesssim C_{1}^{2}\left\Vert \bm{u}\right\Vert _{2}^{2}\log\left(n+d\right)\left\Vert \bm{v}^{\top}\bm{f}_{j}\right\Vert _{\psi_{2}}^{2}\\
			& \lesssim C_{1}^{2}\log\left(n+d\right)\left\Vert \bm{u}\right\Vert _{2}^{2}\left\Vert \bm{v}\right\Vert _{2}^{2},
		\end{align*}
		where $\|\cdot\|_{\psi_{2}}$ denotes the sub-Gaussian norm \citep{Vershynin2012}.
		Turning to the second term of \eqref{eq:XY-decompose-123}, we have
		\begin{align*}
			\left\Vert \mathbb{E}\left[X_{j}Y_{j}\ind_{\mathcal{A}_{j}}\right]\right\Vert _{\psi_{1}} & \lesssim\mathbb{E}\left[X_{j}Y_{j}\ind_{\mathcal{A}_{j}}\right]\leq\mathbb{E}\left[\left(\bm{u}^{\top}\bm{f}_{j}\right)^{2}\left(\bm{v}^{\top}\bm{f}_{j}\right)^{2}\right]\\
			& =\left\Vert \bm{u}\right\Vert _{2}^{2}\left\Vert \bm{v}\right\Vert _{2}^{2}+2\left(\bm{u}^{\top}\bm{v}\right)^{2}\leq3\left\Vert \bm{u}\right\Vert _{2}^{2}\left\Vert \bm{v}\right\Vert _{2}^{2},
		\end{align*}
		where the last step arises from Cauchy-Schwarz. The above results
		taken collectively give 
		\[
		\left\Vert X_{j}Y_{j}\ind_{\mathcal{A}_{j}}-\mathbb{E}\left[X_{j}Y_{j}\ind_{\mathcal{A}_{j}}\right]\right\Vert _{\psi_{1}}\lesssim C_{1}^{2}\log\left(n+d\right)\left\Vert \bm{u}\right\Vert _{2}^{2}\left\Vert \bm{v}\right\Vert _{2}^{2}.
		\]
		Applying the Bernstein inequality \citep[Theorem 2.8.2]{vershynin2016high}
		then yields 
		\[
		\alpha_{1}\lesssim\sqrt{\frac{\log\left(n+d\right)^{3}}{n}}\left\Vert \bm{u}\right\Vert _{2}^{2}\left\Vert \bm{v}\right\Vert _{2}^{2}+\frac{\log^{2}\left(n+d\right)}{n}\left\Vert \bm{u}\right\Vert _{2}^{2}\left\Vert \bm{v}\right\Vert _{2}^{2}\lesssim\sqrt{\frac{\log\left(n+d\right)^{3}}{n}}\left\Vert \bm{u}\right\Vert _{2}^{2}\left\Vert \bm{v}\right\Vert _{2}^{2}
		\]
		with probability exceeding $1-O((n+d)^{-100})$, provided that $n\gg\log(n+d)$. 
	\end{itemize}
	Combine the preceding bounds on $\alpha_{1}$ and $\alpha_{2}$ to
	achieve 
	\begin{align*}
		\left|\frac{1}{n}\sum_{j=1}^{n}\left[X_{j}Y_{j}-\mathbb{E}\left[X_{j}Y_{j}\right]\right]\right| & \leq\alpha_{1}+\alpha_{2}\lesssim\sqrt{\frac{\log\left(n+d\right)^{3}}{n}}\left\Vert \bm{u}\right\Vert _{2}^{2}\left\Vert \bm{v}\right\Vert _{2}^{2}+\left(n+d\right)^{-10}\left\Vert \bm{u}\right\Vert _{2}^{2}\left\Vert \bm{v}\right\Vert _{2}^{2}\\
		& \lesssim\sqrt{\frac{\log\left(n+d\right)^{3}}{n}}\left\Vert \bm{u}\right\Vert _{2}^{2}\left\Vert \bm{v}\right\Vert _{2}^{2}
	\end{align*}
	with probability exceeding $1-O((n+d)^{-100})$. It is straightforward
	to verify that 
	\[
	\mathbb{E}\left[X_{j}Y_{j}\right]=\left\Vert \bm{u}\right\Vert _{2}^{2}\left\Vert \bm{v}\right\Vert _{2}^{2}+2\left(\bm{u}^{\top}\bm{v}\right)^{2}
	\]
	for each $j\in[n]$, thus concluding the proof. \end{proof}

\subsection{Auxiliary lemmas for Theorem \ref{thm:pca-complete}}

\subsubsection{Proof of Lemma \ref{lemma:pca-2nd-error}\label{appendix:proof-pca-2nd-error}}

To begin with, we remind the reader that $\bm{R}\in\mathcal{O}^{r\times r}$
represents the rotation matrix that best aligns $\bm{U}$ and $\bm{U}^{\star}$,
and the rotation matrix $\bm{Q}$ is chosen to satisfy $\bm{U}^{\star}\bm{Q}=\bm{U}^{\natural}$.
In addition, we have also shown in \eqref{eq:connection-RU-R-denoisingPCA}
that $\bm{R}\bm{Q}\in\mathcal{O}^{r\times r}$ is the rotation matrix
that best aligns $\bm{U}$ and $\bm{U}^{\natural}$. Suppose for the
moment that the assumptions of Theorem \ref{thm:hpca_inference_general}
are satisfied (which we shall verify shortly). Then conditional on
$\bm{F}$, invoking Theorem \ref{thm:hpca_inference_general} leads
to 
\begin{equation}
\bm{U}\bm{R}\bm{Q}-\bm{U}^{\natural}=\left[\bm{E}\bm{M}^{\natural\top}+\mathcal{P}_{\mathsf{off}\text{-}\mathsf{diag}}\left(\bm{E}\bm{E}^{\top}\right)\right]\bm{U}^{\natural}\left(\bm{\Sigma}^{\natural}\right)^{-2}+\bm{\Psi},\label{eq:decompose-URQ-1234}
\end{equation}
where the residual matrix $\bm{\Psi}$ can be controlled with high
probability as follows: for each $l\in[d]$ 
\[
\mathbb{P}\left(\left\Vert \bm{\Psi}_{l,\cdot}\right\Vert _{2}\lesssim\zeta_{\mathsf{2nd},l}\left(\bm{F}\right)\,\Big|\,\bm{F}\right)\geq1-O\left(\left(n+d\right)^{-10}\right).
\]
Here, the quantity $\zeta_{\mathsf{2nd}}\left(\bm{F}\right)$ is defined
as 
\begin{align*}
\zeta_{\mathsf{2nd},l}\left(\bm{F}\right) & \coloneqq\left\Vert \bm{U}_{l,\cdot}^{\natural}\right\Vert _{2}\left(\kappa^{\natural2}\sqrt{\frac{\mu^{\natural}r}{d}}\frac{\zeta_{\mathsf{1st}}\left(\bm{F}\right)}{\sigma_{r}^{\natural2}}+\kappa^{\natural2}\frac{\zeta_{\mathsf{1st}}^{2}\left(\bm{F}\right)}{\sigma_{r}^{\natural4}}\right)+\kappa^{\natural2}\frac{\zeta_{\mathsf{1st}}\left(\bm{F}\right)\zeta_{\mathsf{1st},l}\left(\bm{F}\right)}{\sigma_{r}^{\natural4}}\sqrt{\frac{\mu^{\natural}r}{d}}
\end{align*}
with the quantities $\zeta_{\mathsf{1st}}\left(\bm{F}\right)$ and
$\zeta_{\mathsf{1st},l}\left(\bm{F}\right)$ given by 
\begin{align*}
\zeta_{\mathsf{1st}}\left(\bm{F}\right) & \coloneqq\sigma^{2}\sqrt{nd}\log\left(n+d\right)+\sigma\sigma_{1}^{\natural}\sqrt{d\log\left(n+d\right)},\\
\zeta_{\mathsf{1st},l}\left(\bm{F}\right) & \coloneqq\sigma\sigma_{l}\sqrt{nd}\log\left(n+d\right)+\sigma_{l}\sigma_{1}^{\natural}\sqrt{d\log\left(n+d\right)}.
\end{align*}
Given that $\bm{U}^{\star}=\bm{U}^{\natural}\bm{Q}^{\top}$, the above
decomposition \eqref{eq:decompose-URQ-1234} can alternatively be
written as 
\[
\bm{U}\bm{R}-\bm{U}^{\star}=\underbrace{\left[\bm{E}\bm{M}^{\natural\top}+\mathcal{P}_{\mathsf{off}\text{-}\mathsf{diag}}\left(\bm{E}\bm{E}^{\top}\right)\right]\bm{U}^{\natural}\left(\bm{\Sigma}^{\natural}\right)^{-2}\bm{Q}^{\top}}_{\eqqcolon\bm{Z}}+\bm{\Psi}\bm{Q}^{\top}.
\]
When the event $\mathcal{E}_{\mathsf{good}}$ occurs, we can see from
(\ref{eq:good-event-sigma-least-largest-hpca}), (\ref{eq:good-event-kappa-hpca}),
(\ref{eq:good-event-mu-hpca}) and (\ref{eq:good-event-sigma-hpca})
that 
\begin{align*}
\zeta_{\mathsf{1st},l}\left(\bm{F}\right) & \lesssim\sigma_{\mathsf{ub}}\sigma_{\mathsf{ub},l}\sqrt{nd}\log\left(n+d\right)+\sigma_{\mathsf{ub},l}\sigma_{1}^{\star}\sqrt{d\log\left(n+d\right)}\\
 & \asymp\sqrt{\frac{\mu r\log^{4}\left(n+d\right)}{np^{2}}}\sigma_{1}^{\star}\left\Vert \bm{U}_{l,\cdot}^{\star}\bm{\Sigma}^{\star}\right\Vert _{2}+\frac{\omega_{l}^{\star}\omega_{\max}}{p}\sqrt{\frac{d}{n}}\log\left(n+d\right)+\sigma_{1}^{\star}\sqrt{\frac{d}{np}}\left\Vert \bm{U}_{l,\cdot}^{\star}\bm{\Sigma}^{\star}\right\Vert _{2}\log\left(n+d\right)\\
 & \quad+\sigma_{1}^{\star}\omega_{l}^{\star}\sqrt{\frac{d\log\left(n+d\right)}{np}}+\frac{\omega_{\max}}{p}\sqrt{\frac{d}{n}}\log^{3/2}\left(n+d\right)\left\Vert \bm{U}_{l,\cdot}^{\star}\bm{\Sigma}^{\star}\right\Vert _{2}+\frac{\omega_{l}^{\star}}{p}\sqrt{\frac{\mu r\log^{3}\left(n+d\right)}{n}}\sigma_{1}^{\star}\\
 & \eqqcolon\zeta_{\mathsf{1st},l},
\end{align*}
\begin{align*}
\zeta_{\mathsf{1st}}\left(\bm{F}\right) & \lesssim\sigma_{\mathsf{ub}}^{2}\sqrt{nd}\log\left(n+d\right)+\sigma_{\mathsf{ub}}\sigma_{1}^{\natural}\sqrt{d\log\left(n+d\right)}\\
 & \asymp\frac{\mu r\log^{2}\left(n+d\right)}{\sqrt{nd}p}\sigma_{1}^{\star2}+\frac{\omega_{\max}^{2}}{p}\sqrt{\frac{d}{n}}\log\left(n+d\right)+\sigma_{1}^{\star2}\sqrt{\frac{\mu r}{np}}\log\left(n+d\right)+\sigma_{1}^{\star}\omega_{\max}\sqrt{\frac{d\log\left(n+d\right)}{np}}\\
 & \eqqcolon\zeta_{\mathsf{1st}},
\end{align*}
and 
\begin{align*}
\zeta_{\mathsf{2nd},l}\left(\bm{F}\right) & \lesssim\left\Vert \bm{U}_{l,\cdot}^{\star}\right\Vert _{2}\left(\sqrt{\frac{\kappa^{3}\mu r\log\left(n+d\right)}{d}}\frac{\zeta_{\mathsf{1st}}}{\sigma_{r}^{\star2}}+\kappa\frac{\zeta_{\mathsf{1st}}^{2}}{\sigma_{r}^{\star4}}\right)+\frac{\zeta_{\mathsf{1st}}\zeta_{\mathsf{1st},l}}{\sigma_{r}^{\star4}}\sqrt{\frac{\kappa^{3}\mu r\log\left(n+d\right)}{d}}\\
 & \eqqcolon\zeta_{\mathsf{2nd},l}.
\end{align*}
These bounds taken together imply that 
\[
\mathbb{P}\left(\left\Vert \bm{\Psi}_{l,\cdot}\right\Vert _{2}\ind_{\mathcal{E}_{\mathsf{good}}}\lesssim\zeta_{\mathsf{2nd},l}\,\big|\,\bm{F}\right)\geq1-O\left(\left(n+d\right)^{-10}\right).
\]
Additionally, in view of the facts that $\mathcal{E}_{\mathsf{good}}$
is $\sigma(\bm{F})$-measurable and 
\[
\left\Vert \bm{\Psi}_{l,\cdot}\right\Vert _{2}=\left\Vert \bm{\Psi}_{l,\cdot}\bm{Q}^{\top}\right\Vert _{2}=\left\Vert \left(\bm{U}\bm{R}-\bm{U}^{\star}-\bm{Z}\right)_{l,\cdot}\right\Vert _{2,\infty},
\]
one can readily demonstrate that 
\[
\mathbb{P}\left(\left\Vert \left(\bm{U}\bm{R}-\bm{U}^{\star}-\bm{Z}\right)_{l,\cdot}\right\Vert _{2}\lesssim\zeta_{\mathsf{2nd},l}\,\big|\,\bm{F}\right)\geq1-O\left(\left(n+d\right)^{-10}\right)
\]
on the high-probability event $\mathcal{E}_{\mathsf{good}}$.

It remains to verify the assumptions of Theorem \ref{thm:hpca_inference_general},
which requires 
\begin{equation}
d\gtrsim\kappa^{\natural4}\mu^{\natural}r+\mu^{\natural2}r\log^{2}\left(n+d\right),\quad n\gtrsim r\log^{4}\left(n+d\right),\quad B\lesssim\frac{\sigma_{\mathsf{ub}}\min\left\{ \sqrt{n_{2}},\sqrt[4]{n_{1}n_{2}}\right\} }{\sqrt{\log n}},\quad\zeta_{\mathsf{1st}}\ll\frac{\sigma_{r}^{\natural2}}{\kappa^{\natural2}}\label{eq:condition-check-distribution-hpca}
\end{equation}
whenever $\mathcal{E}_{\mathsf{good}}$ occurs. According to (\ref{eq:good-event-sigma-least-largest-hpca}),
(\ref{eq:good-event-kappa-hpca}), (\ref{eq:good-event-mu-hpca}),
(\ref{eq:good-event-B-hpca}) and the definition of $\zeta_{\mathsf{2nd}}$,
the conditions in \eqref{eq:condition-check-distribution-hpca} are
guaranteed to hold as long as 
\[
d\gtrsim\kappa^{3}\mu^{2}r\log^{4}\left(n+d\right),\qquad n\gtrsim r\log^{4}\left(n+d\right),\qquad\text{and}\qquad\zeta_{\mathsf{1st}}\ll\frac{\sigma_{r}^{\star2}}{\kappa}.
\]
Finally, we would like to take a closer inspection on the condition
$\zeta_{\mathsf{1st}}\ll\sigma_{r}^{\star2}/\kappa$; in fact, we
seek to derive sufficient conditions which guarantee that $\zeta_{\mathsf{1st}}/\sigma_{r}^{\star2}\ll\delta$
for any given $\delta>0$, which will also be useful in subsequent
analysis. It is straightforward to check that $\zeta_{\mathsf{1st}}/\sigma_{r}^{\star2}\ll\delta$
is equivalent to\begin{subequations}\label{subeq:pca-1st-equivalent}
\begin{align}
ndp^{2}\gg\delta^{-2}\kappa^{2}\mu^{2}r^{2}\log^{4}\left(n+d\right), & \qquad np\gg\delta^{-2}\kappa^{2}\mu r\log^{2}\left(n+d\right),\\
\text{and}\qquad\frac{\omega_{\max}^{2}}{p\sigma_{r}^{\star2}}\sqrt{\frac{d}{n}}\ll\frac{\delta}{\log\left(n+d\right)}, & \qquad\frac{\omega_{\max}}{\sigma_{r}^{\star}}\sqrt{\frac{d}{np}}\ll\frac{\delta}{\sqrt{\kappa\log\left(n+d\right)}}.
\end{align}
\end{subequations}By taking $\delta\coloneqq1/\kappa$, we see that
$\zeta_{\mathsf{1st}}\ll\sigma_{r}^{\star2}/\kappa$ is guaranteed
as long as the following conditions hold: 
\[
ndp^{2}\gg\kappa^{4}\mu^{2}r^{2}\log^{4}\left(n+d\right),\qquad np\gg\kappa^{4}\mu r\log^{2}\left(n+d\right)
\]
\[
\text{and}\qquad\frac{\omega_{\max}^{2}}{p\sigma_{r}^{\star2}}\sqrt{\frac{d}{n}}\ll\frac{1}{\kappa\log\left(n+d\right)},\qquad\frac{\omega_{\max}}{\sigma_{r}^{\star}}\sqrt{\frac{d}{np}}\ll\frac{1}{\sqrt{\kappa^{3}\log\left(n+d\right)}}.
\]
This concludes the proof.

\subsubsection{Proof of Lemma \ref{lemma:pca-covariance-concentration}\label{appendix:proof-pca-covariance-concentration}}

Let us write $\widetilde{\bm{\Sigma}}_{l}$ as the superposition of
two components: 
\begin{align*}
\widetilde{\bm{\Sigma}}_{l} & =\underbrace{\bm{Q}(\bm{\Sigma}^{\natural})^{-1}\bm{V}^{\natural\top}\mathsf{diag}\left\{ \sigma_{l,1}^{2},\ldots,\sigma_{l,n}^{2}\right\} \bm{V}^{\natural}(\bm{\Sigma}^{\natural})^{-1}\bm{Q}^{\top}}_{\eqqcolon\widetilde{\bm{\Sigma}}_{l,1}}\\
 & \quad+\underbrace{\bm{Q}(\bm{\Sigma}^{\natural})^{-2}\bm{U}^{\natural\top}\mathsf{diag}\Biggl\{\sum_{j:j\neq l}\sigma_{l,j}^{2}\sigma_{1,j}^{2},\ldots,\sum_{j:j\neq l}\sigma_{l,j}^{2}\sigma_{d,j}^{2}\Biggr\}\bm{U}^{\natural}(\bm{\Sigma}^{\natural})^{-2}\bm{Q}^{\top}}_{\eqqcolon\widetilde{\bm{\Sigma}}_{l,2}}.
\end{align*}
We shall control $\widetilde{\bm{\Sigma}}_{l,1}$ and $\widetilde{\bm{\Sigma}}_{l,2}$
separately. Throughout this subsection we assume that $\mathcal{E}_{\mathsf{good}}$
happens.

\paragraph{Step 1: identifying a good approximation of $\widetilde{\bm{\Sigma}}_{l,1}$.}

We start with the concentration of the first component $\widetilde{\bm{\Sigma}}_{l,1}$.
The matrices in the middle part satisfy 
\begin{align*}
\bm{V}^{\natural\top}\mathsf{diag}\left\{ \sigma_{l,1}^{2},\ldots,\sigma_{l,n}^{2}\right\} \bm{V}^{\natural} & =\frac{1}{n}\bm{J}^{\top}\bm{F}\mathsf{diag}\left\{ \sigma_{l,1}^{2},\ldots,\sigma_{l,n}^{2}\right\} \bm{F}^{\top}\bm{J}=\bm{J}^{\top}\biggl(\frac{1}{n}\sum_{j=1}^{n}\sigma_{l,j}^{2}\bm{f}_{j}\bm{f}_{j}^{\top}\biggr)\bm{J},
\end{align*}
where we have made use of the identity \eqref{eq:defn-J-denoising}.
To control the term within the parentheses of the above identity,
we can make use of the variance calculation in \eqref{eq:sigma-ij-square-denoising-PCA}
to obtain 
\[
\frac{1}{n}\sum_{j=1}^{n}\sigma_{l,j}^{2}\bm{f}_{j}\bm{f}_{j}^{\top}=\frac{1-p}{n^{2}p}\sum_{j=1}^{n}\left(\bm{U}_{l,\cdot}^{\star}\bm{\Sigma}^{\star}\bm{f}_{j}\right)^{2}\bm{f}_{j}\bm{f}_{j}^{\top}+\frac{\omega_{l}^{\star2}}{n^{2}p}\sum_{j=1}^{n}\bm{f}_{j}\bm{f}_{j}^{\top}.
\]
In addition, it is seen from (\ref{eq:good-event-concentration-1-hpca})
and (\ref{eq:good-event-gaussian-concentration-hpca}) that 
\[
\left\Vert \frac{1}{n}\sum_{j=1}^{n}\left(\bm{U}_{l,\cdot}^{\star}\bm{\Sigma}^{\star}\bm{f}_{j}\right)^{2}\bm{f}_{j}\bm{f}_{j}^{\top}-\left\Vert \bm{U}_{l,\cdot}^{\star}\bm{\Sigma}^{\star}\right\Vert _{2}^{2}\bm{I}_{r}-2\bm{\Sigma}^{\star}\bm{U}_{l,\cdot}^{\star\top}\bm{U}_{l,\cdot}^{\star}\bm{\Sigma}^{\star}\right\Vert \lesssim\sqrt{\frac{r\log^{3}\left(n+d\right)}{n}}\left\Vert \bm{U}_{l,\cdot}^{\star}\bm{\Sigma}^{\star}\right\Vert _{2}^{2}
\]
\[
\text{and}\qquad\left\Vert \frac{1}{n}\sum_{j=1}^{n}\bm{f}_{i}\bm{f}_{i}^{\top}-\bm{I}_{r}\right\Vert =\left\Vert \frac{1}{n}\bm{F}\bm{F}^{\top}-\bm{I}_{r}\right\Vert \lesssim\sqrt{\frac{r+\log\left(n+d\right)}{n}}.
\]
These bounds taken together allow us to express 
\[
\frac{1}{n}\sum_{j=1}^{n}\sigma_{l,j}^{2}\bm{f}_{j}\bm{f}_{j}^{\top}=\frac{1-p}{np}\left(\left\Vert \bm{U}_{l,\cdot}^{\star}\bm{\Sigma}^{\star}\right\Vert _{2}^{2}\bm{I}_{r}+2\bm{\Sigma}^{\star}\bm{U}_{l,\cdot}^{\star\top}\bm{U}_{l,\cdot}^{\star}\bm{\Sigma}^{\star}\right)+\frac{\omega_{l}^{\star2}}{np}\bm{I}_{r}+\bm{R}_{1,1}
\]
for some residual matrix $\bm{R}_{1,1}$ satisfying 
\begin{align*}
\left\Vert \bm{R}_{1,1}\right\Vert  & \lesssim\frac{1-p}{np}\sqrt{\frac{r\log^{3}\left(n+d\right)}{n}}\left\Vert \bm{U}_{l,\cdot}^{\star}\bm{\Sigma}^{\star}\right\Vert _{2}^{2}+\frac{\omega_{l}^{\star2}}{np}\sqrt{\frac{r+\log\left(n+d\right)}{n}}.
\end{align*}
Putting the above pieces together, we arrive at 
\begin{align*}
\widetilde{\bm{\Sigma}}_{l,1} & =\bm{Q}(\bm{\Sigma}^{\natural})^{-1}\bm{J}^{\top}\biggl(\frac{1}{n}\sum_{j=1}^{n}\sigma_{l,j}^{2}\bm{f}_{j}^{\star}\bm{f}_{j}^{\star\top}\biggr)\bm{J}(\bm{\Sigma}^{\natural})^{-1}\bm{Q}^{\top}\\
 & =\underbrace{\left(\frac{1-p}{np}\left\Vert \bm{U}_{l,\cdot}^{\star}\bm{\Sigma}^{\star}\right\Vert _{2}^{2}+\frac{\omega_{l}^{\star2}}{np}\right)\bm{Q}(\bm{\Sigma}^{\natural})^{-1}\bm{J}^{\top}\bm{J}(\bm{\Sigma}^{\natural})^{-1}\bm{Q}^{\top}}_{\eqqcolon\,\bm{\Sigma}_{1,1}}\\
 & \quad+\underbrace{\frac{2\left(1-p\right)}{np}\bm{Q}(\bm{\Sigma}^{\natural})^{-1}\bm{J}^{\top}\bm{\Sigma}^{\star}\bm{U}_{l,\cdot}^{\star\top}\bm{U}_{l,\cdot}^{\star}\bm{\Sigma}^{\star}\bm{J}(\bm{\Sigma}^{\natural})^{-1}\bm{Q}^{\top}}_{\eqqcolon\,\bm{\Sigma}_{1,2}}+\bm{R}_{1,2}
\end{align*}
for some residual matrix 
\[
\bm{R}_{1,2}=\bm{Q}(\bm{\Sigma}^{\natural})^{-1}\bm{J}^{\top}\bm{R}_{1,1}\bm{J}(\bm{\Sigma}^{\natural})^{-1}\bm{Q}^{\top}.
\]
This motivates us to look at $\bm{\Sigma}_{1,1}$ and $\bm{\Sigma}_{1,2}$
separately. 
\begin{itemize}
\item Regarding the matrix $\bm{\Sigma}_{1,1}$, we make the observation
that 
\begin{align}
\left\Vert (\bm{\Sigma}^{\natural})^{-1}\bm{J}^{\top}\bm{J}(\bm{\Sigma}^{\natural})^{-1}-(\bm{\Sigma}^{\natural})^{-2}\right\Vert  & \overset{\text{(i)}}{=}\left\Vert (\bm{\Sigma}^{\natural})^{-1}\left(\bm{J}^{\top}\bm{J}-\bm{Q}^{\top}\bm{Q}\right)(\bm{\Sigma}^{\natural})^{-1}\right\Vert \nonumber \\
 & \leq\left\Vert (\bm{\Sigma}^{\natural})^{-1}\left(\bm{J}-\bm{Q}\right)^{\top}\bm{J}(\bm{\Sigma}^{\natural})^{-1}\right\Vert +\left\Vert (\bm{\Sigma}^{\natural})^{-1}\bm{Q}^{\top}\left(\bm{J}-\bm{Q}\right)(\bm{\Sigma}^{\natural})^{-1}\right\Vert \nonumber \\
 & \overset{\mathrm{(ii)}}{\lesssim}\frac{1}{\sigma_{r}^{\star2}}\left\Vert \bm{J}-\bm{Q}\right\Vert \left(\left\Vert \bm{Q}\right\Vert +\left\Vert \bm{J}\right\Vert \right)\overset{\text{(iii)}}{\lesssim}\frac{\kappa}{\sigma_{r}^{\star2}}\sqrt{\frac{r+\log\left(n+d\right)}{n}}.\label{eq:pca-cov-concentration-inter-1}
\end{align}
Here, (i) comes from the fact that $\bm{Q}$ is a orthonormal matrix,
(ii) follows from the property \eqref{eq:good-event-sigma-least-largest-hpca},
whereas (iii) utilizes the property (\ref{eq:good-event-J-Q-hpca})
and its direct application 
\begin{equation}
\left\Vert \bm{J}\right\Vert \leq\left\Vert \bm{Q}\right\Vert +\left\Vert \bm{J}-\bm{Q}\right\Vert \leq1+\kappa\sqrt{\frac{r+\log\left(n+d\right)}{n}}\leq2,\label{eq:J-bounded}
\end{equation}
provided that $n\gg\kappa^{2}(r+\log(n+d))$. In addition, from (\ref{eq:good-event-J-Q-hpca})
and (\ref{eq:good-event-sigma-least-largest-hpca}) we know that 
\begin{align}
\left\Vert \bm{Q}(\bm{\Sigma}^{\natural})^{-1}-\left(\bm{\Sigma}^{\star}\right)^{-1}\bm{Q}\right\Vert  & \leq\left\Vert \left(\bm{\Sigma}^{\star}\right)^{-1}\left(\bm{\Sigma}^{\star}\bm{Q}-\bm{Q}\bm{\Sigma}^{\natural}\right)(\bm{\Sigma}^{\natural})^{-1}\right\Vert \lesssim\frac{1}{\sigma_{r}^{\star2}}\left\Vert \bm{\Sigma}^{\star}\bm{Q}-\bm{Q}\bm{\Sigma}^{\natural}\right\Vert \nonumber \\
 & \lesssim\frac{\kappa}{\sigma_{r}^{\star}}\sqrt{\frac{r+\log\left(n+d\right)}{n}}.\label{eq:Q-sigma-inv-interchange}
\end{align}
This immediately leads to 
\begin{align}
\left\Vert \bm{Q}(\bm{\Sigma}^{\natural})^{-2}\bm{Q}^{\top}-\left(\bm{\Sigma}^{\star}\right)^{-2}\right\Vert  & =\left\Vert \bm{Q}(\bm{\Sigma}^{\natural})^{-1}(\bm{\Sigma}^{\natural})^{-1}\bm{Q}^{\top}-\left(\bm{\Sigma}^{\star}\right)^{-1}\bm{Q}\bm{Q}^{\top}\left(\bm{\Sigma}^{\star}\right)^{-1}\right\Vert \nonumber \\
 & \leq\left\Vert \left[\bm{Q}(\bm{\Sigma}^{\natural})^{-1}-\left(\bm{\Sigma}^{\star}\right)^{-1}\bm{Q}\right](\bm{\Sigma}^{\natural})^{-1}\bm{Q}^{\top}\right\Vert +\left\Vert \left(\bm{\Sigma}^{\star}\right)^{-1}\bm{Q}\left[\bm{Q}(\bm{\Sigma}^{\natural})^{-1}-\left(\bm{\Sigma}^{\star}\right)^{-1}\bm{Q}\right]^{\top}\right\Vert \nonumber \\
 & \leq\left(\left\Vert (\bm{\Sigma}^{\natural})^{-1}\right\Vert +\left\Vert (\bm{\Sigma}^{\star})^{-1}\right\Vert \right)\left\Vert \bm{Q}\right\Vert \left\Vert \bm{Q}(\bm{\Sigma}^{\natural})^{-1}-\left(\bm{\Sigma}^{\star}\right)^{-1}\bm{Q}\right\Vert \nonumber \\
 & \lesssim\frac{\kappa}{\sigma_{r}^{\star2}}\sqrt{\frac{r+\log\left(n+d\right)}{n}},\label{eq:pca-cov-concentration-inter-2}
\end{align}
where we have again used (\ref{eq:good-event-sigma-least-largest-hpca}).
Taking (\ref{eq:pca-cov-concentration-inter-1}) and (\ref{eq:pca-cov-concentration-inter-2})
collectively gives 
\begin{align*}
 & \left\Vert \bm{Q}(\bm{\Sigma}^{\natural})^{-1}\bm{J}^{\top}\bm{J}(\bm{\Sigma}^{\natural})^{-1}\bm{Q}^{\top}-\left(\bm{\Sigma}^{\star}\right)^{-2}\right\Vert \\
 & \quad\leq\left\Vert \bm{Q}\left[(\bm{\Sigma}^{\natural})^{-1}\bm{J}^{\top}\bm{J}(\bm{\Sigma}^{\natural})^{-1}-(\bm{\Sigma}^{\natural})^{-2}\right]\bm{Q}^{\top}\right\Vert +\left\Vert \bm{Q}(\bm{\Sigma}^{\natural})^{-2}\bm{Q}^{\top}-\left(\bm{\Sigma}^{\star}\right)^{-2}\right\Vert \\
 & \quad\leq\left\Vert (\bm{\Sigma}^{\natural})^{-1}\bm{J}^{\top}\bm{J}(\bm{\Sigma}^{\natural})^{-1}-(\bm{\Sigma}^{\natural})^{-2}\right\Vert +\left\Vert \bm{Q}(\bm{\Sigma}^{\natural})^{-2}\bm{Q}^{\top}-\left(\bm{\Sigma}^{\star}\right)^{-2}\right\Vert \\
 & \quad\lesssim\frac{\kappa}{\sigma_{r}^{\star2}}\sqrt{\frac{r+\log\left(n+d\right)}{n}}.
\end{align*}
Substitution into the definition of $\bm{\Sigma}_{1,1}$ allows us
to conclude that 
\begin{equation}
\left\Vert \bm{\Sigma}_{1,1}-\left(\frac{1-p}{np}\left\Vert \bm{U}_{l,\cdot}^{\star}\bm{\Sigma}^{\star}\right\Vert _{2}^{2}+\frac{\sigma_{l}^{2}}{np}\right)\left(\bm{\Sigma}^{\star}\right)^{-2}\right\Vert \lesssim\left(\frac{1-p}{np}\left\Vert \bm{U}_{l,\cdot}^{\star}\bm{\Sigma}^{\star}\right\Vert _{2}^{2}+\frac{\omega_{l}^{\star2}}{np}\right)\frac{\kappa}{\sigma_{r}^{\star2}}\sqrt{\frac{r+\log\left(n+d\right)}{n}}.\label{eq:Sigma-l1-bound-123}
\end{equation}
\item When it comes to the remaining term $\bm{\Sigma}_{1,2}$, we first
notice that 
\begin{align*}
 & \left\Vert \bm{U}_{l,\cdot}^{\star}\bm{\Sigma}^{\star}\bm{J}(\bm{\Sigma}^{\natural})^{-1}-\bm{U}_{l,\cdot}^{\natural}\right\Vert _{2}\leq\left\Vert \bm{U}_{l,\cdot}^{\star}\bm{\Sigma}^{\star}\bm{Q}(\bm{\Sigma}^{\natural})^{-1}-\bm{U}_{l,\cdot}^{\natural}\right\Vert _{2}+\left\Vert \bm{U}_{l,\cdot}^{\star}\bm{\Sigma}^{\star}\left(\bm{J}-\bm{Q}\right)(\bm{\Sigma}^{\natural})^{-1}\right\Vert _{2}\\
 & \quad\leq\left\Vert \bm{U}_{l,\cdot}^{\star}\left(\bm{\Sigma}^{\star}\bm{Q}-\bm{Q}\bm{\Sigma}^{\natural}\right)(\bm{\Sigma}^{\natural})^{-1}\right\Vert _{2}+\left\Vert \bm{U}_{l,\cdot}^{\star}\bm{Q}-\bm{U}_{l,\cdot}^{\natural}\right\Vert _{2}+\left\Vert \bm{U}_{l,\cdot}^{\star}\bm{\Sigma}^{\star}\left(\bm{J}-\bm{Q}\right)(\bm{\Sigma}^{\natural})^{-1}\right\Vert _{2}\\
 & \quad\lesssim\frac{1}{\sigma_{r}^{\star}}\left\Vert \bm{U}_{l,\cdot}^{\star}\right\Vert _{2}\left\Vert \bm{\Sigma}^{\star}\bm{Q}-\bm{Q}\bm{\Sigma}^{\natural}\right\Vert +\frac{\sigma_{1}^{\star}}{\sigma_{r}^{\star}}\left\Vert \bm{U}_{l,\cdot}^{\star}\right\Vert _{2}\left\Vert \bm{J}-\bm{Q}\right\Vert \\
 & \quad\lesssim\kappa^{3/2}\sqrt{\frac{r+\log\left(n+d\right)}{n}}\left\Vert \bm{U}_{l,\cdot}^{\star}\right\Vert _{2},
\end{align*}
where the penultimate line uses (\ref{eq:good-event-sigma-least-largest-hpca})
and fact that $\bm{U}^{\star}\bm{Q}=\bm{U}^{\natural}$, and the last
line relies on the property (\ref{eq:good-event-J-Q-hpca}). An immediate
consequence is that 
\begin{align*}
\left\Vert \bm{U}_{l,\cdot}^{\star}\bm{\Sigma}^{\star}\bm{J}(\bm{\Sigma}^{\natural})^{-1}\right\Vert _{2} & \leq\left\Vert \bm{U}_{l,\cdot}^{\natural}\right\Vert _{2}+\left\Vert \bm{U}_{l,\cdot}^{\star}\bm{\Sigma}^{\star}\bm{J}(\bm{\Sigma}^{\natural})^{-1}-\bm{U}_{l,\cdot}^{\natural}\right\Vert _{2}\lesssim\left\Vert \bm{U}_{l,\cdot}^{\star}\right\Vert _{2},
\end{align*}
provided that $n\gg\kappa^{3}r+\kappa^{3}\log(n+d)$. This combined
with the fact $\bm{U}^{\star}\bm{Q}=\bm{U}^{\natural}$ immediately
yields 
\begin{align*}
 & \left\Vert \bm{Q}(\bm{\Sigma}^{\natural})^{-1}\bm{J}^{\top}\bm{\Sigma}^{\star}\bm{U}_{l,\cdot}^{\star\top}\bm{U}_{l,\cdot}^{\star}\bm{\Sigma}^{\star}\bm{J}(\bm{\Sigma}^{\natural})^{-1}\bm{Q}^{\top}-\bm{U}_{l,\cdot}^{\star\top}\bm{U}_{l,\cdot}^{\star}\right\Vert \\
 & \quad=\left\Vert (\bm{\Sigma}^{\natural})^{-1}\bm{J}^{\top}\bm{\Sigma}^{\star}\bm{U}_{l,\cdot}^{\star\top}\bm{U}_{l,\cdot}^{\star}\bm{\Sigma}^{\star}\bm{J}(\bm{\Sigma}^{\natural})^{-1}-\bm{U}_{l,\cdot}^{\natural\top}\bm{U}_{l,\cdot}^{\natural}\right\Vert \\
 & \quad\leq\left\Vert (\bm{\Sigma}^{\natural})^{-1}\bm{J}^{\top}\bm{\Sigma}^{\star}\bm{U}_{l,\cdot}^{\star\top}\left(\bm{U}_{l,\cdot}^{\star}\bm{\Sigma}^{\star}\bm{J}(\bm{\Sigma}^{\natural})^{-1}-\bm{U}_{l,\cdot}^{\natural}\right)\right\Vert +\left\Vert \left(\bm{U}_{l,\cdot}^{\star}\bm{\Sigma}^{\star}\bm{J}(\bm{\Sigma}^{\natural})^{-1}-\bm{U}_{l,\cdot}^{\natural}\right)^{\top}\bm{U}_{l,\cdot}^{\natural}\right\Vert \\
 & \quad\leq\left\Vert \bm{U}_{l,\cdot}^{\star}\bm{\Sigma}^{\star}\bm{J}(\bm{\Sigma}^{\natural})^{-1}-\bm{U}_{l,\cdot}^{\natural}\right\Vert _{2}\left(\left\Vert \bm{U}_{l,\cdot}^{\star}\bm{\Sigma}^{\star}\bm{J}(\bm{\Sigma}^{\natural})^{-1}\right\Vert _{2}+\left\Vert \bm{U}_{l,\cdot}^{\natural}\right\Vert _{2}\right)\\
 & \quad\lesssim\kappa^{3/2}\sqrt{\frac{r+\log\left(n+d\right)}{n}}\left\Vert \bm{U}_{l,\cdot}^{\star}\right\Vert _{2}^{2},
\end{align*}
and as a result, 
\begin{equation}
\left\Vert \bm{\Sigma}_{1,2}-\frac{2\left(1-p\right)}{np}\bm{U}_{l,\cdot}^{\star\top}\bm{U}_{l,\cdot}^{\star}\right\Vert \lesssim\frac{1-p}{np}\kappa^{3/2}\sqrt{\frac{r+\log\left(n+d\right)}{n}}\left\Vert \bm{U}_{l,\cdot}^{\star}\right\Vert _{2}^{2}.\label{eq:Sigma-l2-bound-123}
\end{equation}
\end{itemize}
Combining the above bounds \eqref{eq:Sigma-l1-bound-123} and \eqref{eq:Sigma-l2-bound-123},
we can demonstrate that 
\[
\widetilde{\bm{\Sigma}}_{l,1}=\bm{\Sigma}_{1,1}+\bm{\Sigma}_{1,2}+\bm{R}_{1,2}=\underbrace{\left(\frac{1-p}{np}\left\Vert \bm{U}_{l,\cdot}^{\star}\bm{\Sigma}^{\star}\right\Vert _{2}^{2}+\frac{\omega_{l}^{\star2}}{np}\right)\left(\bm{\Sigma}^{\star}\right)^{-2}+\frac{2\left(1-p\right)}{np}\bm{U}_{l,\cdot}^{\star\top}\bm{U}_{l,\cdot}^{\star}}_{\eqqcolon\,\bm{\Sigma}_{l,1}^{\star}}+\bm{R}_{1}
\]
holds for some residual matrix $\bm{R}_{l}$ satisfying 
\begin{align*}
\left\Vert \bm{R}_{1}\right\Vert  & \leq\left\Vert \bm{R}_{1,2}\right\Vert +\left\Vert \bm{\Sigma}_{1,1}-\left(\frac{1-p}{np}\left\Vert \bm{U}_{l,\cdot}^{\star}\bm{\Sigma}^{\star}\right\Vert _{2}^{2}+\frac{\omega_{l}^{\star2}}{np}\right)\left(\bm{\Sigma}^{\star}\right)^{-2}\right\Vert +\left\Vert \bm{\Sigma}_{1,2}-\frac{2\left(1-p\right)}{np}\bm{U}_{l,\cdot}^{\star\top}\bm{U}_{l,\cdot}^{\star}\right\Vert \\
 & \overset{\text{(i)}}{\lesssim}\frac{1}{\sigma_{r}^{\star2}}\left(\frac{1-p}{np}\sqrt{\frac{r\log^{3}\left(n+d\right)}{n}}\left\Vert \bm{U}_{l,\cdot}^{\star}\bm{\Sigma}^{\star}\right\Vert _{2}^{2}+\frac{\omega_{l}^{\star2}}{np}\sqrt{\frac{r+\log\left(n+d\right)}{n}}\right)\\
 & \quad+\left(\frac{1-p}{np}\left\Vert \bm{U}_{l,\cdot}^{\star}\bm{\Sigma}^{\star}\right\Vert _{2}^{2}+\frac{\omega_{l}^{\star2}}{np}\right)\frac{\kappa}{\sigma_{r}^{\star2}}\sqrt{\frac{r+\log\left(n+d\right)}{n}}+\frac{1-p}{np}\kappa^{3/2}\sqrt{\frac{r+\log\left(n+d\right)}{n}}\left\Vert \bm{U}_{l,\cdot}^{\star}\right\Vert _{2}^{2}\\
 & \lesssim\frac{1-p}{np\sigma_{r}^{\star2}}\left(\sqrt{\frac{r\log^{3}\left(n+d\right)}{n}}+\kappa^{3/2}\sqrt{\frac{r+\log\left(n+d\right)}{n}}\right)\left\Vert \bm{U}_{l,\cdot}^{\star}\bm{\Sigma}^{\star}\right\Vert _{2}^{2}+\frac{\omega_{l}^{\star2}}{np\sigma_{r}^{\star2}}\kappa\sqrt{\frac{r+\log\left(n+d\right)}{n}}.
\end{align*}
Here, (i) has made use of (\ref{eq:good-event-sigma-least-largest-hpca})
and the property (\ref{eq:J-bounded}).

\paragraph{Step 2: controlling the spectrum of $\widetilde{\bm{\Sigma}}_{l,1}$. }

We first observe that 
\begin{align*}
\left(\frac{1-p}{np}\left\Vert \bm{U}_{l,\cdot}^{\star}\bm{\Sigma}^{\star}\right\Vert _{2}^{2}+\frac{\omega_{l}^{\star2}}{np}\right)\left(\bm{\Sigma}^{\star}\right)^{-2}+\frac{2\left(1-p\right)}{np}\bm{U}_{l,\cdot}^{\star\top}\bm{U}_{l,\cdot}^{\star}\succeq\bm{\Sigma}_{l,1}^{\star} & \succeq\left(\frac{1-p}{np}\left\Vert \bm{U}_{l,\cdot}^{\star}\bm{\Sigma}^{\star}\right\Vert _{2}^{2}+\frac{\omega_{l}^{\star2}}{np}\right)\left(\bm{\Sigma}^{\star}\right)^{-2},
\end{align*}
and as a result, \begin{subequations}\label{eq:lambda-max-min-Sigma-UL}
\begin{align}
\lambda_{\max}\left(\bm{\Sigma}_{l,1}^{\star}\right) & \leq\frac{3\left(1-p\right)}{np\sigma_{r}^{\star2}}\left\Vert \bm{U}_{l,\cdot}^{\star}\bm{\Sigma}^{\star}\right\Vert _{2}^{2}+\frac{\omega_{l}^{\star2}}{np\sigma_{r}^{\star2}},\\
\lambda_{\min}\left(\bm{\Sigma}_{l,1}^{\star}\right) & \geq\frac{1-p}{np\sigma_{1}^{\star2}}\left\Vert \bm{U}_{l,\cdot}^{\star}\bm{\Sigma}^{\star}\right\Vert _{2}^{2}+\frac{\omega_{l}^{\star2}}{np\sigma_{1}^{\star2}}.
\end{align}
\end{subequations} This immediately implies that the condition number
of $\bm{\Sigma}_{l,1}^{\star}$ is at most $3\kappa$. In addition,
it follows from the preceding bounds that 
\[
\frac{\left\Vert \bm{R}_{1}\right\Vert }{\lambda_{\min}\left(\bm{\Sigma}_{U,l}^{\star}\right)}\lesssim\sqrt{\frac{\kappa^{2}r\log^{3}\left(n+d\right)+\kappa^{5}r+\kappa^{5}\log\left(n+d\right)}{n}}\asymp\sqrt{\frac{\kappa^{5}r\log^{3}\left(n+d\right)}{n}}.
\]
This means that $\Vert\bm{R}_{l}\Vert\ll\lambda_{\min}(\bm{\Sigma}_{l,1}^{\star})$
holds as long as $n\gg\kappa^{5}r\log^{3}(n+d)$, and as a consequence,
\[
\lambda_{\min}\big(\widetilde{\bm{\Sigma}}_{l,1}\big)\in\left[\lambda_{\min}\left(\bm{\Sigma}_{l,1}^{\star}\right)-\left\Vert \bm{R}_{1}\right\Vert ,\lambda_{\min}\left(\bm{\Sigma}_{l,1}^{\star}\right)+\left\Vert \bm{R}_{1}\right\Vert \right]\quad\Longrightarrow\quad\lambda_{\min}\big(\widetilde{\bm{\Sigma}}_{l}\big)\asymp\lambda_{\min}\left(\bm{\Sigma}_{U,l}^{\star}\right).
\]
Therefore, one can conclude that 
\[
\left\Vert \widetilde{\bm{\Sigma}}_{l,1}-\bm{\Sigma}_{l,1}^{\star}\right\Vert =\left\Vert \bm{R}_{1}\right\Vert \lesssim\sqrt{\frac{\kappa^{5}r\log^{3}\left(n+d\right)}{n}}\lambda_{\min}\left(\bm{\Sigma}_{U,l}^{\star}\right).
\]

\paragraph{Step 3: identifying a good approximation of $\widetilde{\bm{\Sigma}}_{l,2}$. }

We now turn attention to approximating $\widetilde{\bm{\Sigma}}_{l,2}$,
and it suffices to study $\sum_{j\neq l}^{n}\sigma_{l,j}^{2}\sigma_{i,j}^{2}$
for each $i\in[d]$. In view of the expression \eqref{eq:sigma-ij-square-denoising-PCA},
we have

\begin{align}
\sum_{j:j\neq l}\sigma_{l,j}^{2}\sigma_{i,j}^{2} & =\sum_{j:j\neq l}\left[\frac{1-p}{np}\left(\bm{U}_{l,\cdot}^{\star}\bm{\Sigma}^{\star}\bm{f}_{j}\right)^{2}+\frac{\omega_{l}^{\star2}}{np}\right]\left[\frac{1-p}{np}\left(\bm{U}_{i,\cdot}^{\star}\bm{\Sigma}^{\star}\bm{f}_{j}\right)^{2}+\frac{\omega_{i}^{\star2}}{np}\right]\nonumber \\
 & =\left(\frac{1-p}{np}\right)^{2}\underbrace{\sum_{j:j\neq l}\left(\bm{U}_{l,\cdot}^{\star}\bm{\Sigma}^{\star}\bm{f}_{j}\right)^{2}\left(\bm{U}_{i,\cdot}^{\star}\bm{\Sigma}^{\star}\bm{f}_{j}\right)^{2}}_{\eqqcolon\alpha_{1}}+\frac{\left(1-p\right)\omega_{i}^{\star2}}{n^{2}p^{2}}\underbrace{\sum_{j:j\neq l}\left(\bm{U}_{l,\cdot}^{\star}\bm{\Sigma}^{\star}\bm{f}_{j}\right)^{2}}_{\eqqcolon\alpha_{2}}\nonumber \\
 & \quad+\frac{\left(1-p\right)\omega_{l}^{\star2}}{n^{2}p^{2}}\underbrace{\sum_{j:j\neq l}\left(\bm{U}_{i,\cdot}^{\star}\bm{\Sigma}^{\star}\bm{f}_{j}\right)^{2}}_{\eqqcolon\alpha_{3}}+\frac{\left(n-1\right)\omega_{l}^{\star2}\omega_{i}^{\star2}}{n^{2}p^{2}}.\label{eq:sum-sigma-square-1356}
\end{align}
To proceed, we can see from (\ref{eq:good-event-concentration-2-hpca}),
(\ref{eq:good-event-concentration-3-hpca}) and (\ref{eq:good-event-Ui-Sigma-f-hpca})
that 
\begin{align*}
\left|\frac{1}{n}\alpha_{1}-\left[\left\Vert \bm{U}_{l,\cdot}^{\star}\bm{\Sigma}^{\star}\right\Vert _{2}^{2}\left\Vert \bm{U}_{i,\cdot}^{\star}\bm{\Sigma}^{\star}\right\Vert _{2}^{2}+2S_{l,i}^{\star2}\right]\right| & \lesssim\sqrt{\frac{\log^{3}\left(n+d\right)}{n}}\left\Vert \bm{U}_{l,\cdot}^{\star}\bm{\Sigma}^{\star}\right\Vert _{2}^{2}\left\Vert \bm{U}_{i,\cdot}^{\star}\bm{\Sigma}^{\star}\right\Vert _{2}^{2},\\
\left|\frac{1}{n}\alpha_{2}-\left\Vert \bm{U}_{l,\cdot}^{\star}\bm{\Sigma}^{\star}\right\Vert _{2}^{2}\right| & \lesssim\sqrt{\frac{\log\left(n+d\right)}{n}}\left\Vert \bm{U}_{l,\cdot}^{\star}\bm{\Sigma}^{\star}\right\Vert _{2}^{2},\\
\left|\frac{1}{n}\alpha_{3}-\left\Vert \bm{U}_{i,\cdot}^{\star}\bm{\Sigma}^{\star}\right\Vert _{2}^{2}\right| & \lesssim\sqrt{\frac{\log\left(n+d\right)}{n}}\left\Vert \bm{U}_{i,\cdot}^{\star}\bm{\Sigma}^{\star}\right\Vert _{2}^{2}.
\end{align*}
Substitution into \eqref{eq:sum-sigma-square-1356} yields 
\begin{equation}
\sum_{j:j\neq l}\sigma_{l,j}^{2}\sigma_{i,j}^{2}=\underbrace{\frac{\left(1-p\right)^{2}}{np^{2}}\left(\left\Vert \bm{U}_{l,\cdot}^{\star}\bm{\Sigma}^{\star}\right\Vert _{2}^{2}\left\Vert \bm{U}_{i,\cdot}^{\star}\bm{\Sigma}^{\star}\right\Vert _{2}^{2}+2S_{l,i}^{\star2}\right)+\frac{1-p}{np^{2}}\left(\omega_{i}^{\star2}\left\Vert \bm{U}_{l,\cdot}^{\star}\bm{\Sigma}^{\star}\right\Vert _{2}^{2}+\omega_{l}^{\star2}\left\Vert \bm{U}_{i,\cdot}^{\star}\bm{\Sigma}^{\star}\right\Vert _{2}^{2}\right)+\frac{\omega_{l}^{\star2}\omega_{i}^{\star2}}{np^{2}}}_{\eqqcolon\,d_{l,i}}+r_{i}\label{eq:defn-d-li-156}
\end{equation}
for some residual term $r_{i}$ satisfying 
\begin{align}
\left|r_{i}\right| & \lesssim\sqrt{\frac{\log^{3}\left(n+d\right)}{n}}d_{l,i}.\label{eq:ri-UB-135}
\end{align}
As a consequence, the above calculations together with the definition
of $\widetilde{\bm{\Sigma}}_{l,2}$ allow one to write 
\begin{equation}
\widetilde{\bm{\Sigma}}_{l,2}=\bm{Q}(\bm{\Sigma}^{\natural})^{-2}\bm{U}^{\natural\top}\mathsf{diag}\left\{ d_{l,i}\right\} _{i=1}^{d}\bm{U}^{\natural}(\bm{\Sigma}^{\natural})^{-2}\bm{Q}^{\top}+\bm{R}_{2}\label{eq:Sigma-l2-decomposition-R2-355}
\end{equation}
for some residual matrix $\bm{R}_{2}$ satisfying 
\begin{align*}
\left\Vert \bm{R}_{2}\right\Vert  & \leq\left\Vert \bm{Q}(\bm{\Sigma}^{\natural})^{-2}\bm{U}^{\natural\top}\mathsf{diag}\left\{ |r_{1}|,\ldots,|r_{d}|\right\} \bm{U}^{\natural}(\bm{\Sigma}^{\natural})^{-2}\bm{Q}^{\top}\right\Vert \\
 & \lesssim\frac{1}{\big(\sigma_{r}(\bm{\Sigma}^{\natural})\big)^{4}}\max_{1\leq i\leq d}\left|r_{i}\right|\lesssim\sqrt{\frac{\log^{3}\left(n+d\right)}{n}}\max_{1\leq i\leq d}\frac{1}{\sigma_{r}^{\star4}}d_{l,i},
\end{align*}
where the last inequality makes use of \eqref{eq:good-event-sigma-least-largest-hpca}
and \eqref{eq:ri-UB-135}.

Further, it turns out that the term $\bm{Q}(\bm{\Sigma}^{\natural})^{-2}\bm{U}^{\natural\top}$
used in \eqref{eq:Sigma-l2-decomposition-R2-355} can be well approximated
by $\left(\bm{\Sigma}^{\star}\right)^{-2}\bm{Q}$. To see this, we
note that 
\begin{align*}
\left\Vert \bm{Q}(\bm{\Sigma}^{\natural})^{-2}-\left(\bm{\Sigma}^{\star}\right)^{-2}\bm{Q}\right\Vert  & \leq\left\Vert \left[\bm{Q}(\bm{\Sigma}^{\natural})^{-1}-\left(\bm{\Sigma}^{\star}\right)^{-1}\bm{Q}\right](\bm{\Sigma}^{\natural})^{-1}\right\Vert +\left\Vert \left(\bm{\Sigma}^{\star}\right)^{-1}\left[\bm{Q}(\bm{\Sigma}^{\natural})^{-1}-\left(\bm{\Sigma}^{\star}\right)^{-1}\bm{Q}\right]\right\Vert \\
 & \overset{\text{(i)}}{\asymp}\frac{1}{\sigma_{r}^{\star}}\left\Vert \bm{Q}(\bm{\Sigma}^{\natural})^{-1}-\left(\bm{\Sigma}^{\star}\right)^{-1}\bm{Q}\right\Vert \overset{\text{(ii)}}{\lesssim}\frac{\kappa}{\sigma_{r}^{\star2}}\sqrt{\frac{r+\log\left(n+d\right)}{n}},
\end{align*}
where (i) comes from (\ref{eq:good-event-sigma-least-largest-hpca}),
and (ii) results from (\ref{eq:Q-sigma-inv-interchange}). This combined
with the identity $\bm{U}^{\star}=\bm{U}^{\natural}\bm{Q}^{\top}$
gives 
\begin{align*}
\left\Vert \bm{Q}(\bm{\Sigma}^{\natural})^{-2}\bm{U}^{\natural\top}-\left(\bm{\Sigma}^{\star}\right)^{-2}\bm{U}^{\star\top}\right\Vert  & =\left\Vert \left[\bm{Q}(\bm{\Sigma}^{\natural})^{-2}-\left(\bm{\Sigma}^{\star}\right)^{-2}\bm{Q}\right]\bm{U}^{\natural\top}\right\Vert \\
 & \leq\left\Vert \bm{Q}(\bm{\Sigma}^{\natural})^{-2}-\left(\bm{\Sigma}^{\star}\right)^{-2}\bm{Q}\right\Vert \\
 & \lesssim\frac{\kappa}{\sigma_{r}^{\star2}}\sqrt{\frac{r+\log\left(n+d\right)}{n}}.
\end{align*}
Therefore, substituting this into \eqref{eq:Sigma-l2-decomposition-R2-355}
reveals that 
\begin{align*}
 & \left\Vert \widetilde{\bm{\Sigma}}_{l,2}-\left(\bm{\Sigma}^{\star}\right)^{-2}\bm{U}^{\star\top}\mathsf{diag}\left\{ d_{l,i}\right\} _{i=1}^{d}\bm{U}^{\star}\left(\bm{\Sigma}^{\star}\right)^{-2}\right\Vert \\
 & \quad\leq\left\Vert \bm{Q}(\bm{\Sigma}^{\natural})^{-2}\bm{U}^{\natural\top}\mathsf{diag}\left\{ d_{l,i}\right\} _{i=1}^{d}\bm{U}^{\natural}(\bm{\Sigma}^{\natural})^{-2}\bm{Q}^{\top}-\left(\bm{\Sigma}^{\star}\right)^{-2}\bm{U}^{\star\top}\mathsf{diag}\left\{ d_{l,i}\right\} _{i=1}^{d}\bm{U}^{\star}\left(\bm{\Sigma}^{\star}\right)^{-2}\right\Vert +\|\bm{R}_{2}\|\\
 & \quad\leq\left\Vert \left(\bm{Q}(\bm{\Sigma}^{\natural})^{-2}\bm{U}^{\natural\top}-\left(\bm{\Sigma}^{\star}\right)^{-2}\bm{U}^{\star\top}\right)\mathsf{diag}\left\{ d_{l,i}\right\} _{i=1}^{d}\bm{U}^{\natural}(\bm{\Sigma}^{\natural})^{-2}\bm{Q}^{\top}\right\Vert \\
 & \quad\quad+\left\Vert \left(\bm{\Sigma}^{\star}\right)^{-2}\bm{U}^{\star\top}\mathsf{diag}\left\{ d_{l,i}\right\} _{i=1}^{d}\left(\bm{Q}(\bm{\Sigma}^{\natural})^{-2}\bm{U}^{\natural\top}-\left(\bm{\Sigma}^{\star}\right)^{-2}\bm{U}^{\star\top}\right)^{\top}\right\Vert +\|\bm{R}_{2}\|\\
 & \quad\lesssim\frac{\kappa}{\sigma_{r}^{\star4}}\sqrt{\frac{r+\log\left(n+d\right)}{n}}\max_{1\leq i\leq d}d_{l,i}+\|\bm{R}_{2}\|,
\end{align*}
where the last relation relies on (\ref{eq:good-event-sigma-least-largest-hpca}).
These bounds taken collectively allow one to express 
\begin{equation}
\widetilde{\bm{\Sigma}}_{l,2}=\bm{\Sigma}_{l,2}^{\star}+\bm{R}_{3},\label{eq:Sigma-l2-R3-expression}
\end{equation}
where 
\begin{equation}
\bm{\Sigma}_{l,2}^{\star}\coloneqq\left(\bm{\Sigma}^{\star}\right)^{-2}\bm{U}^{\star\top}\mathsf{diag}\left\{ d_{l,i}\right\} _{i=1}^{d}\bm{U}^{\star}\left(\bm{\Sigma}^{\star}\right)^{-2}\label{eq:defn-Sigma-l2-156}
\end{equation}
and $\bm{R}_{3}$ is some residual matrix satisfying 
\begin{align*}
\left\Vert \bm{R}_{3}\right\Vert  & \lesssim\left\Vert \bm{R}_{2}\right\Vert +\frac{\kappa}{\sigma_{r}^{\star4}}\sqrt{\frac{r+\log\left(n+d\right)}{n}}\max_{1\leq i\leq d}d_{l,i}\lesssim\sqrt{\frac{\kappa^{2}r\log^{3}\left(n+d\right)}{n}}\frac{1}{\sigma_{r}^{\star4}}\max_{1\leq i\leq d}d_{l,i}.
\end{align*}

\paragraph{Step 4: controlling the spectrum of $\widetilde{\bm{\Sigma}}_{l,2}$. }

Next, we will investigate the spectrum of $\widetilde{\bm{\Sigma}}_{l,2}$.
It is straightforward to show that the matrix $\bm{\Sigma}_{l,2}^{\star}$
defined in \eqref{eq:defn-Sigma-l2-156} obeys 
\begin{align}
\left\Vert \bm{\Sigma}_{l,2}^{\star}\right\Vert  & \leq\frac{1}{\sigma_{r}^{\star4}}\max_{1\leq i\leq d}\left|d_{l,i}\right|\leq\frac{3\kappa\mu r\left(1-p\right)^{2}}{ndp^{2}\sigma_{r}^{\star2}}\left\Vert \bm{U}_{l,\cdot}^{\star}\bm{\Sigma}^{\star}\right\Vert _{2}^{2}+\frac{\kappa\mu r\left(1-p\right)}{ndp^{2}\sigma_{r}^{\star2}}\omega_{l}^{\star2}+\frac{1-p}{np^{2}\sigma_{r}^{\star4}}\omega_{\max}^{2}\left\Vert \bm{U}_{l,\cdot}^{\star}\bm{\Sigma}^{\star}\right\Vert _{2}^{2}+\frac{\omega_{l}^{\star2}\omega_{\max}^{\star2}}{np^{2}\sigma_{r}^{\star4}}.\label{eq:Sigma-l2-UB-123}
\end{align}
In addition, we claim that the following inequality is valid.\begin{claim}\label{claim:matrix-min-eigenvalue}It
holds that 
\begin{equation}
\bm{U}^{\star\top}\mathsf{diag}\left\{ \left\Vert \bm{U}_{i,\cdot}^{\star}\bm{\Sigma}^{\star}\right\Vert _{2}^{2}\right\} _{i=1}^{d}\bm{U}^{\star}\succeq\frac{\sigma_{r}^{\star2}}{4d}\bm{I}_{r}.\label{eq:matrix-min-eigenvalue}
\end{equation}
\end{claim}With this claim in place, we can combine it with the definition
\eqref{eq:defn-d-li-156} of $d_{l,i}$ to bound 
\begin{align*}
\bm{\Sigma}_{l,2}^{\star} & =\left(\bm{\Sigma}^{\star}\right)^{-2}\bm{U}^{\star\top}\mathsf{diag}\left\{ d_{l,i}\right\} _{i=1}^{d}\bm{U}^{\star}\left(\bm{\Sigma}^{\star}\right)^{-2}\\
 & \succeq\frac{\left(1-p\right)^{2}}{np^{2}}\left\Vert \bm{U}_{l,\cdot}^{\star}\bm{\Sigma}^{\star}\right\Vert _{2}^{2}\left(\bm{\Sigma}^{\star}\right)^{-2}\bm{U}^{\star\top}\mathsf{diag}\left\{ \left\Vert \bm{U}_{i,\cdot}^{\star}\bm{\Sigma}^{\star}\right\Vert _{2}^{2}\right\} _{i=1}^{d}\bm{U}^{\star}\left(\bm{\Sigma}^{\star}\right)^{-2}\\
 & \quad+\frac{1-p}{np^{2}}\omega_{l}^{\star2}\left(\bm{\Sigma}^{\star}\right)^{-2}\bm{U}^{\star\top}\mathsf{diag}\left\{ \left\Vert \bm{U}_{i,\cdot}^{\star}\bm{\Sigma}^{\star}\right\Vert _{2}^{2}\right\} _{i=1}^{d}\bm{U}^{\star}\left(\bm{\Sigma}^{\star}\right)^{-2}+\frac{1-p}{np^{2}}\omega_{\min}^{2}\left\Vert \bm{U}_{l,\cdot}^{\star}\bm{\Sigma}^{\star}\right\Vert _{2}^{2}\left(\bm{\Sigma}^{\star}\right)^{-4}\\
 & \quad+\frac{1}{np^{2}}\omega_{l}^{\star2}\left(\bm{\Sigma}^{\star}\right)^{-2}\bm{U}^{\star\top}\mathsf{diag}\left\{ \omega_{i}^{\star2}\right\} _{i=1}^{d}\bm{U}^{\star}\left(\bm{\Sigma}^{\star}\right)^{-2}\\
 & \succeq\left[\frac{\left(1-p\right)^{2}}{4ndp^{2}}\sigma_{r}^{\star2}\left\Vert \bm{U}_{l,\cdot}^{\star}\bm{\Sigma}^{\star}\right\Vert _{2}^{2}+\frac{1-p}{4ndp^{2}}\sigma_{r}^{\star2}\omega_{l}^{\star2}+\frac{1-p}{np^{2}}\omega_{\min}^{2}\left\Vert \bm{U}_{l,\cdot}^{\star}\bm{\Sigma}^{\star}\right\Vert _{2}^{2}+\frac{\omega_{l}^{\star2}\omega_{\min}^{2}}{np^{2}}\right]\left(\bm{\Sigma}^{\star}\right)^{-4},
\end{align*}
which in turn leads to 
\begin{align}
\lambda_{\min}\left(\bm{\Sigma}_{l,2}^{\star}\right) & \geq\frac{1}{\sigma_{1}^{\star4}}\left[\frac{\left(1-p\right)^{2}}{4ndp^{2}}\sigma_{r}^{\star2}\left\Vert \bm{U}_{l,\cdot}^{\star}\bm{\Sigma}^{\star}\right\Vert _{2}^{2}+\frac{1-p}{4ndp^{2}}\sigma_{r}^{\star2}\omega_{l}^{\star2}+\frac{1-p}{np^{2}}\omega_{\min}^{2}\left\Vert \bm{U}_{l,\cdot}^{\star}\bm{\Sigma}^{\star}\right\Vert _{2}^{2}+\frac{\omega_{l}^{\star2}\omega_{\min}^{2}}{np^{2}}\right]\nonumber \\
 & =\frac{\left(1-p\right)^{2}}{4ndp^{2}\kappa\sigma_{1}^{\star2}}\left\Vert \bm{U}_{l,\cdot}^{\star}\bm{\Sigma}^{\star}\right\Vert _{2}^{2}+\frac{1-p}{4ndp^{2}\kappa\sigma_{1}^{\star2}}\omega_{l}^{\star2}+\frac{1-p}{np^{2}\sigma_{1}^{\star4}}\omega_{\min}^{2}\left\Vert \bm{U}_{l,\cdot}^{\star}\bm{\Sigma}^{\star}\right\Vert _{2}^{2}+\frac{\omega_{l}^{\star2}\omega_{\min}^{2}}{np^{2}\sigma_{1}^{\star4}}.\label{eq:Sigma-l2-LB-123}
\end{align}
The above bounds \eqref{eq:Sigma-l2-UB-123} and \eqref{eq:Sigma-l2-LB-123}
imply that the condition number of $\bm{\Sigma}_{l,2}^{\star}$ ---
denoted by $\kappa(\bm{\Sigma}_{l,2}^{\star})$ --- is upper bounded
by 
\[
\kappa(\bm{\Sigma}_{l,2}^{\star})=\frac{\left\Vert \bm{\Sigma}_{l,2}^{\star}\right\Vert }{\lambda_{\min}\left(\bm{\Sigma}_{l,2}^{\star}\right)}\lesssim\kappa^{3}\mu r\kappa_{\omega}.
\]
Moreover, we can also obtain 
\begin{align*}
\left\Vert \bm{R}_{3}\right\Vert  & \lesssim\sqrt{\frac{\kappa^{2}r\log^{3}\left(n+d\right)}{n}}\frac{1}{\sigma_{r}^{\star4}}\max_{1\leq i\leq d}\left|d_{l,i}\right|\\
 & \lesssim\sqrt{\frac{\kappa^{2}r\log^{3}\left(n+d\right)}{n}}\kappa^{3}\mu r\kappa_{\omega}\lambda_{\min}\left(\bm{\Sigma}_{l,2}^{\star}\right)\\
 & \asymp\sqrt{\frac{\kappa^{8}\mu^{2}r^{3}\kappa_{\omega}^{2}\log^{3}\left(n+d\right)}{n}}\lambda_{\min}\left(\bm{\Sigma}_{l,2}^{\star}\right).
\end{align*}

\paragraph{Step 5: putting everything together.}

Combining the above results, we are allowed to express 
\begin{align*}
\widetilde{\bm{\Sigma}}_{l} & =\bm{\Sigma}_{U,l}^{\star}+\bm{R}_{l}
\end{align*}
where $\bm{\Sigma}_{U,l}^{\star}=\bm{\Sigma}_{l,1}^{\star}+\bm{\Sigma}_{l,2}^{\star}$
and the residual matrix $\bm{R}_{l}$ satisfies 
\begin{align*}
\left\Vert \bm{R}_{l}\right\Vert  & \leq\left\Vert \bm{R}_{1}\right\Vert +\left\Vert \bm{R}_{3}\right\Vert \\
 & \lesssim\sqrt{\frac{\kappa^{5}r\log^{3}\left(n+d\right)}{n}}\lambda_{\min}\left(\bm{\Sigma}_{l,1}^{\star}\right)+\sqrt{\frac{\kappa^{8}\mu^{2}r^{3}\kappa_{\omega}^{2}\log^{3}\left(n+d\right)}{n}}\lambda_{\min}\left(\bm{\Sigma}_{l,2}^{\star}\right)\\
 & \lesssim\sqrt{\frac{\kappa^{8}\mu^{2}r^{3}\kappa_{\omega}^{2}\log^{3}\left(n+d\right)}{n}}\lambda_{\min}\left(\bm{\Sigma}_{U,l}^{\star}\right).
\end{align*}
An immediate consequence is that the conditional number of $\bm{\Sigma}_{U,l}^{\star}$
is upper bounded by the maximum of the conditional numbers of $\bm{\Sigma}_{l,1}^{\star}$
and $\bm{\Sigma}_{l,2}^{\star}$, which is at most $O(\kappa^{3}\mu r\kappa_{\omega})$.
Consequently, when $n\gg\kappa^{8}\mu^{2}r^{3}\kappa_{\omega}^{2}\log^{3}(n+d)$,
the conditional number of $\widetilde{\bm{\Sigma}}_{l}$ is also upper
bounded by $O(\kappa^{3}\mu r\kappa_{\omega})$. This finishes the
proof, as long as Claim~\ref{claim:matrix-min-eigenvalue} can be
justified.

\begin{proof}[Proof of Claim \ref{claim:matrix-min-eigenvalue}]For
any $c\in(0,1)$, we define an index set 
\[
\mathcal{I}_{c}\coloneqq\{i\in[d]:\Vert\bm{U}_{i,\cdot}^{\star}\Vert_{2}\geq\sqrt{c/d}\}.
\]
Then for any $\bm{v}\in\mathbb{R}^{d}$, one has 
\begin{align*}
\bm{v}^{\top}\bm{U}^{\star\top}\mathsf{diag}\left\{ \left\Vert \bm{U}_{i,\cdot}^{\star}\bm{\Sigma}^{\star}\right\Vert _{2}^{2}\right\} _{i=1}^{d}\bm{U}^{\star}\bm{v} & =\sum_{i=1}^{d}\left(\bm{v}^{\top}\bm{U}_{i,\cdot}^{\star}\right)^{2}\left\Vert \bm{U}_{i,\cdot}^{\star}\bm{\Sigma}^{\star}\right\Vert _{2}^{2}\geq\sigma_{r}^{\star2}\sum_{i=1}^{d}\left(\bm{v}^{\top}\bm{U}_{i,\cdot}^{\star}\right)^{2}\left\Vert \bm{U}_{i,\cdot}^{\star}\right\Vert _{2}^{2}\\
 & \geq\sigma_{r}^{\star2}\sum_{i\in\mathcal{I}_{c}}\left(\bm{v}^{\top}\bm{U}_{i,\cdot}^{\star}\right)^{2}\left\Vert \bm{U}_{i,\cdot}^{\star}\right\Vert _{2}^{2}\geq\frac{c\sigma_{r}^{\star2}}{d}\sum_{i\in\mathcal{I}_{c}}\left(\bm{v}^{\top}\bm{U}_{i,\cdot}^{\star}\right)^{2}\\
 & =\frac{c\sigma_{r}^{\star2}}{d}\left[\left\Vert \bm{v}\right\Vert _{2}^{2}-\sum_{i\in[d]\setminus\mathcal{I}_{c}}\left(\bm{v}^{\top}\bm{U}_{i,\cdot}^{\star}\right)^{2}\right]\\
 & \geq\frac{c\sigma_{r}^{\star2}}{d}\left[\left\Vert \bm{v}\right\Vert _{2}^{2}-\sum_{i\in[d]\setminus\mathcal{I}_{c}}\left\Vert \bm{v}\right\Vert _{2}^{2}\left\Vert \bm{U}_{i,\cdot}^{\star}\right\Vert _{2}^{2}\right]\\
 & \geq\frac{c\sigma_{r}^{\star2}}{d}\left[\left\Vert \bm{v}\right\Vert _{2}^{2}-\sum_{i\in[d]\setminus\mathcal{I}_{c}}\left\Vert \bm{v}\right\Vert _{2}^{2}\frac{c}{d}\right]\\
 & \geq\frac{c\sigma_{r}^{\star2}}{d}\left[\left\Vert \bm{v}\right\Vert _{2}^{2}-c\left\Vert \bm{v}\right\Vert _{2}^{2}\right]\geq\frac{c\left(1-c\right)\sigma_{r}^{\star2}}{d}\left\Vert \bm{v}\right\Vert _{2}^{2}.
\end{align*}
Since the above inequality holds for arbitrary $c\in(0,1)$ and any
$\bm{v}\in\mathbb{R}^{d}$, taking $c=1/2$ leads to 
\[
\inf_{\bm{v}\in\mathbb{R}^{d}}\bm{v}^{\top}\bm{U}^{\star\top}\mathsf{diag}\left\{ \left\Vert \bm{U}_{i,\cdot}^{\star}\bm{\Sigma}^{\star}\right\Vert _{2}^{2}\right\} _{i=1}^{d}\bm{U}^{\star}\bm{v}\geq\frac{\sigma_{r}^{\star2}}{4d}\left\Vert \bm{v}\right\Vert _{2}^{2}.
\]
Therefore, we can conclude that 
\[
\bm{U}^{\star\top}\mathsf{diag}\left\{ \left\Vert \bm{U}_{i,\cdot}^{\star}\bm{\Sigma}^{\star}\right\Vert _{2}^{2}\right\} _{i=1}^{d}\bm{U}^{\star}\succeq\frac{\sigma_{r}^{\star2}}{4d}\bm{I}_{r}.
\]
\end{proof}

\subsubsection{Proof of Lemma \ref{lemma:pca-normal-approximation}\label{appendix:proof-pca-normal-approximation}}

\paragraph{Step 1: Gaussian approximation of $\bm{Z}_{l,\cdot}$ using the Berry-Esseen
Theorem. }

It is first seen that 
\[
\bm{Z}_{l,\cdot}=\sum_{j=1}^{n}\bm{Y}_{j},\qquad\text{where}\quad\bm{Y}_{j}=E_{l,j}\left[\bm{V}_{j,\cdot}^{\natural}(\bm{\Sigma}^{\natural})^{-1}+\left[\mathcal{P}_{-l,\cdot}\left(\bm{E}_{\cdot,j}\right)\right]^{\top}\bm{U}^{\natural}(\bm{\Sigma}^{\natural})^{-2}\right]\bm{Q}^{\top}.
\]
Invoking the Berry-Esseen theorem (see Theorem \ref{thm:berry-esseen-multivariate})
then yields 
\begin{equation}
\sup_{\mathcal{C}\in\mathscr{C}^{r}}\left|\mathbb{P}\left(\bm{Z}_{l,\cdot}\in\mathcal{C}\big|\bm{F}\right)-\mathbb{P}\left(\mathcal{N}\big(\bm{0},\widetilde{\bm{\Sigma}}_{l}\big)\in\mathcal{C}\big|\bm{F}\right)\right|\lesssim r^{1/4}\gamma\left(\bm{F}\right),\label{eq:sup-Z-tilde-Sigma-456}
\end{equation}
where $\widetilde{\bm{\Sigma}}_{l}$ is the covariance matrix of $\bm{Z}_{l,\cdot}$
that has been calculated in \eqref{eq:Sigma-tilde-distribution-hpca},
and $\gamma\left(\bm{F}\right)$ is defined to be 
\[
\gamma\left(\bm{F}\right)\coloneqq\sum_{i=1}^{n}\mathbb{E}\left[\left\Vert \bm{Y}_{j}\widetilde{\bm{\Sigma}}_{l}^{-1/2}\right\Vert _{2}^{3}\Big|\bm{F}\right].
\]
In the sequel, let us develop an upper bound on $\gamma\left(\bm{F}\right)$.
Towards this, we find it helpful to define the quantity $B_{l}$ such
that 
\[
\max_{j\in[n]}\left|E_{l,j}\right|\leq B_{l}.
\]

\begin{itemize}
\item Towards this, we first provide a (conditional) high probability bound
for each $\Vert\bm{Y}_{j}\Vert$. In view of \citet[Lemma 12]{cai2019subspace},
with probability exceeding $1-O((n+d)^{-101})$ we have 
\begin{align*}
\left\Vert \left[\mathcal{P}_{-l,\cdot}\left(\bm{E}_{\cdot,j}\right)\right]^{\top}\bm{U}^{\natural}\right\Vert _{2} & =\Bigg\|\sum_{i:i\neq l}E_{i,j}\bm{U}_{i,\cdot}^{\natural}\Bigg\|_{2}\lesssim\left(B\log\left(n+d\right)+\sigma\sqrt{d\log\left(n+d\right)}\right)\left\Vert \bm{U}^{\natural}\right\Vert _{2,\infty}\\
 & \asymp\left(B\log\left(n+d\right)+\sigma\sqrt{d\log\left(n+d\right)}\right)\left\Vert \bm{U}^{\star}\right\Vert _{2,\infty}
\end{align*}
holds for all $j\in[n]$, where the last relation makes use of the
fact that $\bm{U}^{\natural}=\bm{U}^{\star}\bm{Q}$ for some orthonormal
matrix $\bm{Q}$ (see \eqref{eq:defn-Q-denoising-PCA}). Consequently,
with probability exceeding $1-O((n+d)^{-101})$, 
\begin{align*}
\left\Vert \bm{Y}_{j}\right\Vert _{2} & \leq\frac{1}{\sigma_{r}^{\natural}}B_{l}\left\Vert \bm{V}^{\natural}\right\Vert _{2,\infty}+\frac{1}{\sigma_{r}^{\natural2}}B_{l}\left(B\log\left(n+d\right)+\sigma\sqrt{d\log\left(n+d\right)}\right)\left\Vert \bm{U}^{\star}\right\Vert _{2,\infty}
\end{align*}
holds for each $j\in[n]$. If we define 
\[
C_{\mathsf{prob}}\coloneqq\widetilde{C}_{1}B_{l}\left[\frac{1}{\sigma_{r}^{\natural}}\left\Vert \bm{V}^{\natural}\right\Vert _{2,\infty}+\frac{1}{\sigma_{r}^{\natural2}}\left(B\log\left(n+d\right)+\sigma\sqrt{d\log\left(n+d\right)}\right)\left\Vert \bm{U}^{\star}\right\Vert _{2,\infty}\right]
\]
for some sufficiently large constant $\widetilde{C}_{1}>0$, then
the union bound over $1\leq j\leq n$ guarantees that 
\begin{equation}
\mathbb{P}\left(\max_{1\leq j\leq n}\left\Vert \bm{Y}_{j}\right\Vert _{2}\leq C_{\mathsf{prob}}\big|\bm{F}\right)\geq1-O\left(\left(n+d\right)^{-100}\right).\label{eq:Y-j-probabilistic}
\end{equation}
\item In addition, we also know that $\Vert\bm{Y}_{j}\Vert_{2}$ admits
a trivial deterministic upper bound as follows 
\begin{align*}
\left\Vert \bm{Y}_{j}\right\Vert _{2} & \lesssim\frac{1}{\sigma_{r}^{\natural}}B_{l}\left\Vert \bm{V}^{\natural}\right\Vert _{2,\infty}+\frac{1}{\sigma_{r}^{\natural2}}B_{l}\left\Vert \bm{E}_{\cdot,j}\right\Vert _{2}\lesssim\frac{1}{\sigma_{r}^{\natural}}B_{l}\left\Vert \bm{V}^{\natural}\right\Vert _{2,\infty}+\frac{1}{\sigma_{r}^{\natural2}}B_{l}B\sqrt{d}.
\end{align*}
This means that by defining the quantity 
\[
C_{\mathsf{det}}\coloneqq\widetilde{C}_{2}B_{l}\left(\frac{1}{\sigma_{r}^{\natural}}\left\Vert \bm{V}^{\natural}\right\Vert _{2,\infty}+\frac{1}{\sigma_{r}^{\natural2}}B\sqrt{d}\right)
\]
for some sufficiently large constant $\widetilde{C}_{2}>0$, we can
ensure that 
\begin{equation}
\max_{1\leq j\leq n}\left\Vert \bm{Y}_{j}\right\Vert _{2}\leq C_{\mathsf{det}}.\label{eq:Y-j-deterministic}
\end{equation}
\end{itemize}
The above arguments taken together allow one to bound $\gamma\left(\bm{F}\right)$
as follows 
\begin{align*}
\gamma\left(\bm{F}\right) & \leq\lambda_{\min}^{-3/2}\big(\widetilde{\bm{\Sigma}}_{l}\big)\sum_{j=1}^{n}\mathbb{E}\left[\left\Vert \bm{Y}_{j}\right\Vert _{2}^{3}\big|\bm{F}\right]\\
 & =\lambda_{\min}^{-3/2}\big(\widetilde{\bm{\Sigma}}_{l}\big)\sum_{j=1}^{n}\mathbb{E}\left[\left\Vert \bm{Y}_{j}\right\Vert _{2}^{3}\ind_{\left\Vert \bm{Y}_{j}\right\Vert _{2}\leq C_{\mathsf{prob}}}\big|\bm{F}\right]+\lambda_{\min}^{-3/2}\big(\widetilde{\bm{\Sigma}}_{l}\big)\sum_{j=1}^{n}\mathbb{E}\left[\left\Vert \bm{Y}_{j}\right\Vert _{2}^{3}\ind_{\left\Vert \bm{Y}_{j}\right\Vert _{2}>C_{\mathsf{prob}}}\big|\bm{F}\right]\\
 & \overset{\text{(i)}}{\lesssim}\lambda_{\min}^{-3/2}\big(\widetilde{\bm{\Sigma}}_{l}\big)C_{\mathsf{prob}}\sum_{j=1}^{n}\mathbb{E}\left[\left\Vert \bm{Y}_{j}\right\Vert _{2}^{2}\big|\bm{F}\right]+\lambda_{\min}^{-3/2}\big(\widetilde{\bm{\Sigma}}_{l}\big)\sum_{j=1}^{n}C_{\mathsf{det}}^{3}\mathbb{P}\left(\max_{1\leq j\leq n}\left\Vert \bm{Y}_{j}\right\Vert _{2}\leq C_{\mathsf{prob}}\big|\bm{F}\right)\\
 & \overset{\text{(ii)}}{\lesssim}\lambda_{\min}^{-3/2}\big(\widetilde{\bm{\Sigma}}_{l}\big)\mathsf{tr}\big(\widetilde{\bm{\Sigma}}_{l}\big)C_{\mathsf{prob}}+\lambda_{\min}^{-3/2}\big(\widetilde{\bm{\Sigma}}_{l}\big)C_{\mathsf{det}}^{3}\left(n+d\right)^{-99}.
\end{align*}
Here (i) relies on (\ref{eq:Y-j-deterministic}); (ii) follows from
(\ref{eq:Y-j-probabilistic}) as well as 
\begin{align}
\sum_{j=1}^{n}\mathbb{E}\left[\left\Vert \bm{Y}_{j}\right\Vert _{2}^{2}\big|\bm{F}\right] & =\sum_{j=1}^{n}\mathbb{E}\left[\bm{Y}_{j}\bm{Y}_{j}^{\top}\big|\bm{F}\right]=\sum_{j=1}^{n}\mathbb{E}\left[\mathsf{tr}\left(\bm{Y}_{j}^{\top}\bm{Y}_{j}\right)\big|\bm{F}\right]=\mathsf{tr}\left[\sum_{j=1}^{n}\mathbb{E}\left(\bm{Y}_{j}^{\top}\bm{Y}_{j}\big|\bm{F}\right)\right]\nonumber \\
 & =\mathsf{tr}\left(\mathbb{E}\left[\bm{Z}_{l,\cdot}^{\top}\bm{Z}_{l,\cdot}\big|\bm{F}\right]\right)=\mathsf{tr}\big(\widetilde{\bm{\Sigma}}_{l}\big).\label{eq:mse-trace-equation}
\end{align}
Substitution into \eqref{eq:sup-Z-tilde-Sigma-456} then yields 
\begin{align}
 & \sup_{\mathcal{C}\in\mathscr{C}^{r}}\left|\mathbb{P}\left(\bm{Z}_{l,\cdot}\in\mathcal{C}\big|\bm{F}\right)-\mathbb{P}\left(\mathcal{N}\big(\bm{0},\widetilde{\bm{\Sigma}}_{l}\big)\in\mathcal{C}\big|\bm{F}\right)\right|\nonumber \\
 & \qquad\lesssim\underbrace{r^{1/4}\lambda_{\min}^{-3/2}\big(\widetilde{\bm{\Sigma}}_{l}\big)\mathsf{tr}\big(\widetilde{\bm{\Sigma}}_{l}\big)C_{\mathsf{prob}}}_{\eqqcolon\alpha}+\underbrace{r^{1/4}\lambda_{\min}^{-3/2}\big(\widetilde{\bm{\Sigma}}_{l}\big)C_{\mathsf{det}}^{3}\left(n+d\right)^{-99}}_{\eqqcolon\beta}.\label{eq:defn-alpha-beta-hpca-7892}
\end{align}

When the high-probability event $\mathcal{E}_{\mathsf{good}}$ occurs,
we know from Lemma \ref{lemma:pca-covariance-concentration} that
the condition number of $\widetilde{\bm{\Sigma}}_{l}$ is bounded
by $O(\kappa^{3}\mu r\kappa_{\omega})$. This implies that 
\[
\frac{\mathsf{tr}\big(\widetilde{\bm{\Sigma}}_{l}\big)}{\lambda_{\min}^{3/2}\big(\widetilde{\bm{\Sigma}}_{l}\big)}\leq\frac{r\big\|\widetilde{\bm{\Sigma}}_{l}\big\|}{\lambda_{\min}^{3/2}\big(\widetilde{\bm{\Sigma}}_{l}\big)}\lesssim\frac{\kappa^{3}\mu r^{2}\kappa_{\omega}}{\lambda_{\min}^{1/2}\big(\widetilde{\bm{\Sigma}}_{l}\big)}.
\]
and as a consequence, 
\begin{align*}
\alpha & \lesssim\frac{\kappa^{3}\mu r^{9/4}\kappa_{\omega}}{\lambda_{\min}^{1/2}\big(\widetilde{\bm{\Sigma}}_{l}\big)}C_{\mathsf{prob}}\\
 & \asymp\underbrace{\frac{\kappa^{3}\mu r^{9/4}\kappa_{\omega}}{\lambda_{\min}^{1/2}\big(\widetilde{\bm{\Sigma}}_{l}\big)}\frac{B_{l}}{\sigma_{r}^{\natural}}\left\Vert \bm{V}^{\natural}\right\Vert _{2,\infty}}_{\eqqcolon\alpha_{1}}+\underbrace{\frac{\kappa^{3}\mu r^{9/4}\kappa_{\omega}\log\left(n+d\right)}{\lambda_{\min}^{1/2}\big(\widetilde{\bm{\Sigma}}_{l}\big)}\frac{B_{l}B}{\sigma_{r}^{\natural2}}\left\Vert \bm{U}^{\star}\right\Vert _{2,\infty}}_{\eqqcolon\alpha_{2}}\\
 & \quad+\underbrace{\frac{\kappa^{3}\mu r^{9/4}\kappa_{\omega}\log^{3/2}\left(n+d\right)}{\lambda_{\min}^{1/2}\big(\widetilde{\bm{\Sigma}}_{l}\big)}\frac{\sigma B_{l}}{\sigma_{r}^{\natural2}}\sqrt{d}\left\Vert \bm{U}^{\star}\right\Vert _{2,\infty}}_{\eqqcolon\alpha_{3}}.
\end{align*}
All of these terms involve $\lambda_{\min}^{1/2}\big(\widetilde{\bm{\Sigma}}_{l}\big)$,
which can be lower bounded using Lemma \ref{lemma:pca-covariance-concentration}
as follows 
\begin{align}
\lambda_{\min}^{1/2}\big(\widetilde{\bm{\Sigma}}_{l}\big) & \gtrsim\frac{1}{\sqrt{np}\sigma_{1}^{\star}}\left\Vert \bm{U}_{l,\cdot}^{\star}\bm{\Sigma}^{\star}\right\Vert _{2}+\frac{\omega_{l}^{\star}}{\sqrt{np}\sigma_{1}^{\star}}+\frac{1}{\sqrt{ndp^{2}\kappa}\sigma_{1}^{\star}}\left\Vert \bm{U}_{l,\cdot}^{\star}\bm{\Sigma}^{\star}\right\Vert _{2}+\frac{1}{\sqrt{ndp^{2}\kappa}\sigma_{1}^{\star}}\omega_{l}^{\star}\nonumber \\
 & \quad+\frac{1}{\sqrt{np^{2}}\sigma_{1}^{\star2}}\omega_{\min}\left\Vert \bm{U}_{l,\cdot}^{\star}\bm{\Sigma}^{\star}\right\Vert _{2}+\frac{\omega_{l}^{\star}\omega_{\min}}{\sqrt{np^{2}}\sigma_{1}^{\star2}}.\label{eq:lambda-min-half-LB-hpca}
\end{align}
Moreover, it is seen from (\ref{eq:good-event-B-hpca}) that 
\begin{equation}
B_{l}\lesssim\frac{1}{p}\sqrt{\frac{\log\left(n+d\right)}{n}}\left\Vert \bm{U}_{l,\cdot}^{\star}\bm{\Sigma}^{\star}\right\Vert _{2}+\frac{\omega_{l}^{\star}}{p}\sqrt{\frac{\log\left(n+d\right)}{n}}\asymp\frac{1}{p}\sqrt{\frac{\log\left(n+d\right)}{n}}\left(\left\Vert \bm{U}_{l,\cdot}^{\star}\bm{\Sigma}^{\star}\right\Vert _{2}+\omega_{l}^{\star}\right).\label{eq:B-l-bound}
\end{equation}
In the sequel, we shall bound $\alpha_{1}$, $\alpha_{2}$, and $\alpha_{3}$
separately. 
\begin{itemize}
\item We start with bounding $\alpha_{1}$, where we have 
\begin{align*}
\alpha_{1} & \overset{\text{(i)}}{\lesssim}\frac{\kappa^{3}\mu r^{9/4}\kappa_{\omega}}{\lambda_{\min}^{1/2}\big(\widetilde{\bm{\Sigma}}_{l}\big)}\frac{1}{\sigma_{r}^{\star}}\frac{1}{p}\sqrt{\frac{\log\left(n+d\right)}{n}}\left(\left\Vert \bm{U}_{l,\cdot}^{\star}\bm{\Sigma}^{\star}\right\Vert _{2}+\omega_{l}^{\star}\right)\sqrt{\frac{r\log\left(n+d\right)}{n}}\\
 & \lesssim\frac{1}{\lambda_{\min}^{1/2}\big(\widetilde{\bm{\Sigma}}_{l}\big)}\frac{\kappa^{3}\mu r^{11/4}\kappa_{\omega}}{np\sigma_{r}^{\star}}\log\left(n+d\right)\left(\left\Vert \bm{U}_{l,\cdot}^{\star}\bm{\Sigma}^{\star}\right\Vert _{2}+\omega_{l}^{\star}\right)\\
 & \overset{\text{(ii)}}{\lesssim}\frac{\kappa^{3}\mu r^{11/4}\kappa_{\omega}}{np\sigma_{r}^{\star}}\log\left(n+d\right)\left(\left\Vert \bm{U}_{l,\cdot}^{\star}\bm{\Sigma}^{\star}\right\Vert _{2}+\omega_{l}^{\star}\right)\left[\frac{1}{\sqrt{np}\sigma_{1}^{\star}}\left(\left\Vert \bm{U}_{l,\cdot}^{\star}\bm{\Sigma}^{\star}\right\Vert _{2}+\omega_{l}^{\star}\right)\right]^{-1}\\
 & \asymp\frac{\kappa^{7/2}\mu r^{11/4}\kappa_{\omega}}{\sqrt{np}}\log\left(n+d\right)\overset{\text{(iii)}}{\lesssim}\frac{1}{\sqrt{\log\left(n+d\right)}}.
\end{align*}
Here (i) follows from (\ref{eq:good-event-sigma-least-largest-hpca}),
(\ref{eq:good-event-V-natural-2-infty-hpca}), (\ref{eq:good-event-B-hpca})
and (\ref{eq:B-l-bound}); (ii) makes use of \eqref{eq:lambda-min-half-LB-hpca};
and (iii) is valid with the proviso that $np\gtrsim\kappa^{7}\mu^{2}r^{11/2}\kappa_{\omega}^{2}\log^{3}(n+d)$. 
\item Regarding $\alpha_{2}$, we know that 
\begin{align*}
\alpha_{2} & \overset{\text{(i)}}{\lesssim}\frac{\kappa^{3}\mu r^{9/4}\kappa_{\omega}\log^{2}\left(n+d\right)}{\lambda_{\min}^{1/2}\big(\widetilde{\bm{\Sigma}}_{l}\big)}\frac{1}{\sigma_{r}^{\star2}}\left(\frac{1}{p^{2}}\sqrt{\frac{\mu r}{n^{2}d}}\sigma_{1}^{\star}+\frac{\omega_{\max}}{np^{2}}\right)\left(\left\Vert \bm{U}_{l,\cdot}^{\star}\bm{\Sigma}^{\star}\right\Vert _{2}+\omega_{l}^{\star}\right)\left\Vert \bm{U}^{\star}\right\Vert _{2,\infty}\\
 & \lesssim\frac{1}{\lambda_{\min}^{1/2}\big(\widetilde{\bm{\Sigma}}_{l}\big)}\frac{1}{\sigma_{r}^{\star2}}\left(\frac{\kappa^{3}\mu^{2}r^{13/4}\kappa_{\omega}}{ndp^{2}}\sigma_{1}^{\star}+\frac{\kappa^{3}\mu^{3/2}r^{11/4}\kappa_{\omega}}{n\sqrt{d}p^{2}}\omega_{\max}\right)\log^{2}\left(n+d\right)\left(\left\Vert \bm{U}_{l,\cdot}^{\star}\bm{\Sigma}^{\star}\right\Vert _{2}+\omega_{l}^{\star}\right)\\
 & \overset{\text{(ii)}}{\lesssim}\frac{1}{\sigma_{r}^{\star2}}\left(\frac{\kappa^{3}\mu^{2}r^{13/4}\kappa_{\omega}}{ndp^{2}}\sigma_{1}^{\star}+\frac{\kappa^{3}\mu^{3/2}r^{11/4}\kappa_{\omega}}{n\sqrt{d}p^{2}}\omega_{\max}\right)\log^{2}\left(n+d\right)\\
 & \qquad\cdot\left(\left\Vert \bm{U}_{l,\cdot}^{\star}\bm{\Sigma}^{\star}\right\Vert _{2}+\omega_{l}^{\star}\right)\left[\frac{1}{\sqrt{ndp^{2}\kappa}\sigma_{1}^{\star}}\left(\left\Vert \bm{U}_{l,\cdot}^{\star}\bm{\Sigma}^{\star}\right\Vert _{2}+\omega_{l}^{\star}\right)\right]^{-1}\\
 & \lesssim\left(\frac{\kappa^{9/2}\mu^{2}r^{13/4}\kappa_{\omega}}{\sqrt{ndp^{2}}}+\frac{\kappa^{4}\mu^{3/2}r^{11/4}\kappa_{\omega}}{\sqrt{np^{2}}}\frac{\omega_{\max}}{\sigma_{r}^{\star}}\right)\log^{2}\left(n+d\right)\overset{\text{(iii)}}{\lesssim}\frac{1}{\sqrt{\log\left(n+d\right)}}.
\end{align*}
Here, (i) follows from (\ref{eq:good-event-sigma-least-largest-hpca}),
(\ref{eq:good-event-B-hpca}) and (\ref{eq:B-l-bound}); (ii) makes
use of \eqref{eq:lambda-min-half-LB-hpca}; and (iii) holds provided
that $ndp^{2}\gtrsim\kappa^{9}\mu^{4}r^{13/2}\kappa_{\omega}^{2}\log^{5}(n+d)$
and 
\[
\frac{\omega_{\max}}{\sigma_{r}^{\star}\sqrt{np^{2}}}\lesssim\frac{\omega_{\max}^{2}}{p\sigma_{r}^{\star2}}\sqrt{\frac{d}{n}}+\frac{1}{\sqrt{nd}p}\lesssim\frac{1}{\kappa^{4}\mu^{3/2}r^{11/4}\kappa_{\omega}\log^{5/2}\left(n+d\right)},
\]
which can be guaranteed by $ndp^{2}\gtrsim\kappa^{8}\mu^{3}r^{11/2}\kappa_{\omega}^{2}\log^{5}(n+d)$
and 
\[
\frac{\omega_{\max}^{2}}{p\sigma_{r}^{\star2}}\sqrt{\frac{d}{n}}\lesssim\frac{1}{\kappa^{4}\mu^{3/2}r^{11/4}\kappa_{\omega}\log^{5/2}\left(n+d\right)}.
\]
\item When it comes to $\alpha_{3}$, we have 
\begin{align*}
\alpha_{3} & \overset{\text{(i)}}{\lesssim}\frac{\kappa^{3}\mu^{3/2}r^{11/4}\kappa_{\omega}\log^{2}\left(n+d\right)}{\lambda_{\min}^{1/2}\big(\widetilde{\bm{\Sigma}}_{l}\big)}\frac{1}{\sigma_{r}^{\star2}}\left(\sqrt{\frac{\mu r\log\left(n+d\right)}{ndp}}\sigma_{1}^{\star}+\frac{\omega_{\max}}{\sqrt{np}}\right)\sqrt{\frac{1}{np^{2}}}\left(\left\Vert \bm{U}_{l,\cdot}^{\star}\bm{\Sigma}^{\star}\right\Vert _{2}+\omega_{l}^{\star}\right)\\
 & \lesssim\frac{\kappa^{3}\mu^{3/2}r^{11/4}\kappa_{\omega}\log^{2}\left(n+d\right)}{\lambda_{\min}^{1/2}\big(\widetilde{\bm{\Sigma}}_{l}\big)}\frac{1}{\sigma_{r}^{\star2}}\left(\sqrt{\frac{\mu r\log\left(n+d\right)}{n^{2}dp^{3}}}\sigma_{1}^{\star}+\frac{\omega_{\max}}{\sqrt{n^{2}p^{3}}}\right)\left(\left\Vert \bm{U}_{l,\cdot}^{\star}\bm{\Sigma}^{\star}\right\Vert _{2}+\omega_{l}^{\star}\right)\\
 & \overset{\text{(ii)}}{\lesssim}\frac{1}{\sigma_{r}^{\star2}}\left(\frac{\kappa^{3}\mu^{2}r^{13/4}\kappa_{\omega}\log^{5/2}\left(n+d\right)}{\sqrt{n^{2}dp^{3}}}\sigma_{1}^{\star}+\omega_{\max}\frac{\kappa^{3}\mu^{3/2}r^{11/4}\kappa_{\omega}\log^{2}\left(n+d\right)}{\sqrt{n^{2}p^{3}}}\right)\\
 & \qquad\cdot\left(\left\Vert \bm{U}_{l,\cdot}^{\star}\bm{\Sigma}^{\star}\right\Vert _{2}+\omega_{l}^{\star}\right)\left[\frac{1}{\sqrt{np}\sigma_{1}^{\star}}\left(\left\Vert \bm{U}_{l,\cdot}^{\star}\bm{\Sigma}^{\star}\right\Vert _{2}+\omega_{l}^{\star}\right)\right]^{-1}\\
 & \lesssim\frac{\kappa^{4}\mu^{2}r^{13/4}\log^{5/2}\left(n+d\right)}{\sqrt{ndp^{2}}}+\frac{\omega_{\max}}{\sigma_{r}^{\star}}\frac{\kappa^{7/2}\mu^{3/2}r^{11/4}\kappa_{\omega}\log^{2}\left(n+d\right)}{\sqrt{np^{2}}}\overset{\text{(iii)}}{\lesssim}\frac{1}{\sqrt{\log\left(n+d\right)}}.
\end{align*}
Here, (i) follows from (\ref{eq:good-event-sigma-hpca}), (\ref{eq:good-event-B-hpca})
and (\ref{eq:B-l-bound}); (ii) makes use of \eqref{eq:lambda-min-half-LB-hpca};
and (iii) holds provided that $ndp^{2}\gtrsim\kappa^{8}\mu^{4}r^{13/2}\kappa_{\omega}^{2}\log^{6}(n+d)$
and 
\[
\frac{\omega_{\max}}{\sigma_{r}^{\star}\sqrt{np^{2}}}\lesssim\frac{\omega_{\max}^{2}}{p\sigma_{r}^{\star2}}\sqrt{\frac{d}{n}}+\frac{1}{\sqrt{nd}p}\lesssim\frac{1}{\kappa^{7/2}\mu^{3/2}r^{11/4}\kappa_{\omega}\log^{5/2}\left(n+d\right)},
\]
which can be guaranteed by $ndp^{2}\gtrsim\kappa^{7}\mu^{3}r^{11/2}\kappa_{\omega}^{2}\log^{5}(n+d)$
and 
\[
\frac{\omega_{\max}^{2}}{p\sigma_{r}^{\star2}}\sqrt{\frac{d}{n}}\lesssim\frac{1}{\kappa^{7/2}\mu^{3/2}r^{11/4}\kappa_{\omega}\log^{5/2}\left(n+d\right)}.
\]
\end{itemize}
In addition, we have learned from (\ref{eq:lambda-min-half-LB-hpca})
and (\ref{eq:B-l-bound}) that 
\begin{equation}
B_{l}\lesssim\frac{1}{p}\sqrt{\frac{\log\left(n+d\right)}{n}}\sqrt{np}\sigma_{1}^{\star}\lambda_{\min}^{1/2}\big(\widetilde{\bm{\Sigma}}_{l}\big)\asymp\sqrt{\frac{\log\left(n+d\right)}{p}}\sigma_{1}^{\star}\lambda_{\min}^{1/2}\big(\widetilde{\bm{\Sigma}}_{l}\big).\label{eq:B-l-bound-2}
\end{equation}
Therefore the term $\beta$ defined in \eqref{eq:defn-alpha-beta-hpca-7892}
can be bounded by 
\begin{align*}
\beta & \overset{\text{(i)}}{\lesssim}r^{1/4}\lambda_{\min}^{-3/2}\big(\widetilde{\bm{\Sigma}}_{l}\big)B_{l}^{3}\left(\frac{1}{\sigma_{r}^{\star}}\sqrt{\frac{r\log\left(n+d\right)}{n}}+\frac{1}{\sigma_{r}^{\star2}}B\sqrt{d}\right)^{3}\left(n+d\right)^{-99}\\
 & \overset{\text{(ii)}}{\lesssim}\kappa^{3/2}r^{1/4}\frac{\log^{3/2}\left(n+d\right)}{p^{3/2}}\left(1+\frac{1}{\sigma_{r}^{\star}}B\sqrt{d}\right)^{3}\left(n+d\right)^{-99}\\
 & \overset{\text{(iii)}}{\lesssim}\kappa^{3/2}r^{1/4}\frac{\log^{3/2}\left(n+d\right)}{p^{3/2}}d^{3/2}\left(n+d\right)^{-99}\overset{\text{(iv)}}{\lesssim}\left(n+d\right)^{-50}.
\end{align*}
Here (i) utilizes (\ref{eq:good-event-sigma-least-largest-hpca})
and (\ref{eq:good-event-V-natural-2-infty-hpca}); (ii) follows from
(\ref{eq:B-l-bound-2}); (iii) holds as long as $B\lesssim\sigma_{r}^{\star}$,
which can be guaranteed by $ndp^{2}\gtrsim\kappa\mu r\log(n+d)$ and
\begin{equation}
\frac{\omega_{\max}}{p\sigma_{r}^{\star}}\sqrt{\frac{1}{n}}\lesssim\frac{\omega_{\max}^{2}}{p\sigma_{r}^{\star2}}\sqrt{\frac{d}{n}}+\frac{1}{\sqrt{nd}p}\lesssim\frac{1}{\sqrt{\log\left(n+d\right)}};\label{eq:am-gm-useful}
\end{equation}
and (iv) is valid provided that $np\gtrsim\kappa r\log(n+d)$ .

Therefore, the preceding bounds allow one to conclude that 
\begin{align}
\sup_{\mathcal{C}\in\mathscr{C}^{r}}\left|\mathbb{P}\left(\bm{Z}_{l,\cdot}\in\mathcal{C}\big|\bm{F}\right)-\mathbb{P}\left(\mathcal{N}\left(\bm{0},\widetilde{\bm{\Sigma}}_{l}\right)\in\mathcal{C}\big|\bm{F}\right)\right| & \lesssim\alpha_{1}+\alpha_{2}+\alpha_{3}+\beta\lesssim\frac{1}{\sqrt{\log\left(n+d\right)}},\label{eq:proof-pca-normal-approximation-inter-1}
\end{align}
provided that $np\gtrsim\kappa^{7}\mu^{2}r^{11/2}\kappa_{\omega}^{2}\log^{3}(n+d)$,
$ndp^{2}\gtrsim\kappa^{9}\mu^{4}r^{13/2}\kappa_{\omega}^{2}\log^{6}(n+d)$
and 
\[
\frac{\omega_{\max}^{2}}{p\sigma_{r}^{\star2}}\sqrt{\frac{d}{n}}\lesssim\frac{1}{\kappa^{4}\mu^{3/2}r^{11/4}\kappa_{\omega}\log^{5/2}\left(n+d\right)}.
\]

\paragraph{Step 2: bounding TV distance between Gaussian distributions. }

In view of Lemma \ref{lemma:pca-covariance-concentration}, we know
that 
\begin{equation}
\left\Vert \widetilde{\bm{\Sigma}}_{l}-\bm{\Sigma}_{U,l}^{\star}\right\Vert \lesssim\sqrt{\frac{\kappa^{8}\mu^{2}r^{3}\kappa_{\omega}^{2}\log^{3}\left(n+d\right)}{n}}\lambda_{\min}\left(\bm{\Sigma}_{U,l}^{\star}\right)\label{eq:Sigma-l-Sigma-l-star-dist-pca}
\end{equation}
holds on the event $\mathcal{E}_{\mathsf{good}}$. In addition, $\bm{\Sigma}_{U,l}^{\star}$
is assumed to be non-singular. As a result, conditional on $\bm{F}$
and the event $\mathcal{E}_{\mathsf{good}}$, one can invoke Theorem
\ref{thm:gaussian-TV-distance} and \eqref{eq:Sigma-l-Sigma-l-star-dist-pca}
to arrive at 
\begin{align*}
 & \sup_{\mathcal{C}\in\mathscr{C}^{r}}\left|\mathbb{P}\left(\mathcal{N}\big(\bm{0},\widetilde{\bm{\Sigma}}_{l}\big)\in\mathcal{C}\mid\bm{F}\right)-\mathbb{P}\left(\mathcal{N}\left(\bm{0},\bm{\Sigma}_{U,l}^{\star}\right)\in\mathcal{C}\right)\right|\leq\mathsf{TV}\left(\mathcal{N}\big(\bm{0},\widetilde{\bm{\Sigma}}_{l}\big),\mathcal{N}\left(\bm{0},\bm{\Sigma}_{U,l}^{\star}\right)\right)\\
 & \quad\asymp\left\Vert \left(\bm{\Sigma}_{U,l}^{\star}\right)^{-1/2}\widetilde{\bm{\Sigma}}_{l}\left(\bm{\Sigma}_{U,l}^{\star}\right)^{-1/2}-\bm{I}_{d}\right\Vert _{\mathrm{F}}=\left\Vert \left(\bm{\Sigma}_{U,l}^{\star}\right)^{-1/2}\big(\widetilde{\bm{\Sigma}}_{l}-\bm{\Sigma}_{U,l}^{\star}\big)\left(\bm{\Sigma}_{U,l}^{\star}\right)^{-1/2}\right\Vert _{\mathrm{F}}\\
 & \quad\lesssim\left\Vert \left(\bm{\Sigma}_{U,l}^{\star}\right)^{-1/2}\right\Vert \left\Vert \widetilde{\bm{\Sigma}}_{l}-\bm{\Sigma}_{U,l}^{\star}\right\Vert _{\mathrm{F}}\left\Vert \left(\bm{\Sigma}_{U,l}^{\star}\right)^{-1/2}\right\Vert \lesssim\sqrt{r}\left\Vert \left(\bm{\Sigma}_{U,l}^{\star}\right)^{-1}\right\Vert \left\Vert \widetilde{\bm{\Sigma}}_{l}-\bm{\Sigma}_{U,l}^{\star}\right\Vert \\
 & \quad\lesssim\frac{1}{\lambda_{\min}\left(\bm{\Sigma}_{U,l}^{\star}\right)}\cdot\sqrt{\frac{\kappa^{8}\mu^{2}r^{4}\kappa_{\omega}^{2}\log^{3}\left(n+d\right)}{n}}\lambda_{\min}\left(\bm{\Sigma}_{U,l}^{\star}\right)\\
 & \quad\lesssim\sqrt{\frac{\kappa^{8}\mu^{2}r^{4}\kappa_{\omega}^{2}\log^{3}\left(n+d\right)}{n}}.
\end{align*}
where $\mathsf{TV}(\cdot,\cdot)$ represents the total-variation distance
between two distributions \citep{tsybakov2009introduction}. The above
inequality taken together with (\ref{eq:proof-pca-normal-approximation-inter-1})
yields 
\begin{equation}
\sup_{\mathcal{C}\in\mathscr{C}^{r}}\left|\mathbb{P}\left(\bm{Z}_{l,\cdot}\in\mathcal{C}\big|\bm{F}\right)-\mathbb{P}\left(\mathcal{N}\left(\bm{0},\bm{\Sigma}_{l}^{\star}\right)\in\mathcal{C}\right)\right|\lesssim\frac{1}{\sqrt{\log\left(n+d\right)}}\label{eq:proof-pca-normal-approximation-inter-2}
\end{equation}
on the event $\mathcal{E}_{\mathsf{good}}$, provided that $n\gtrsim\kappa^{8}\mu^{2}r^{4}\kappa_{\omega}^{2}\log^{4}(n+d)$.

\paragraph{Step 3: accounting for higher-order errors. }

Let us assume that the high-probability event $\mathcal{E}_{\mathsf{good}}$
happens. Lemma \ref{lemma:pca-2nd-error} tells us that 
\[
\mathbb{P}\left(\left\Vert \left(\bm{U}\bm{R}-\bm{U}^{\star}-\bm{Z}\right)_{l,\cdot}\right\Vert _{2}\lesssim\zeta_{\mathsf{2nd},l}\,\big|\,\bm{F}\right)\geq1-O\left(\left(n+d\right)^{-10}\right),
\]
which immediately gives 
\begin{equation}
\mathbb{P}\left(\left\Vert \left(\bm{\Sigma}_{U,l}^{\star}\right)^{-1/2}\left(\bm{U}\bm{R}-\bm{U}^{\star}-\bm{Z}\right)\right\Vert _{2,\infty}\leq\zeta\,\big|\,\bm{F}\right)\geq1-O\left(\left(n+d\right)^{-10}\right).\label{eq:UB-15278}
\end{equation}
Here, the quantity $\zeta$ is defined as 
\begin{equation}
\zeta\coloneqq c_{\zeta}\zeta_{\mathsf{2nd},l}\big(\lambda_{\min}\big(\bm{\Sigma}_{U,l}^{\star}\big)\big)^{-\frac{1}{2}}\label{eq:defn-zeta-denoising}
\end{equation}
for some sufficiently large constant $c_{\zeta}>0$.

For any convex set $\mathcal{C}\in\mathscr{C}^{r}$ and any $\varepsilon$,
recalling the definition of $\mathcal{C}^{\varepsilon}$ in \eqref{eq:defn-C-epsilon}
in Appendix \ref{sec:Additional-notation}, we have 
\begin{align}
 & \mathbb{P}\left(\left(\bm{\Sigma}_{U,l}^{\star}\right)^{-1/2}\bm{Z}_{l,\cdot}\in\mathcal{C}^{-\zeta}\,\big|\,\bm{F}\right)=\mathbb{P}\left(\left(\bm{\Sigma}_{U,l}^{\star}\right)^{-1/2}\bm{Z}_{l,\cdot}\in\mathcal{C}^{-\zeta},\left\Vert \left(\bm{\Sigma}_{U,l}^{\star}\right)^{-1/2}\left(\bm{U}\bm{R}-\bm{U}^{\star}-\bm{Z}\right)\right\Vert _{2,\infty}\leq\zeta\,\big|\,\bm{F}\right)\nonumber \\
 & \quad\qquad\qquad\qquad\qquad\qquad+\mathbb{P}\left(\left(\bm{\Sigma}_{U,l}^{\star}\right)^{-1/2}\bm{Z}_{l,\cdot}\in\mathcal{C}^{-\zeta},\left\Vert \left(\bm{\Sigma}_{U,l}^{\star}\right)^{-1/2}\left(\bm{U}\bm{R}-\bm{U}^{\star}-\bm{Z}\right)\right\Vert _{2,\infty}>\zeta\,\big|\,\bm{F}\right)\nonumber \\
 & \quad\leq\mathbb{P}\left(\left(\bm{\Sigma}_{U,l}^{\star}\right)^{-1/2}\left(\bm{U}\bm{R}-\bm{U}^{\star}\right)_{l,\cdot}\in\mathcal{C}\,\big|\,\bm{F}\right)+\mathbb{P}\left(\left\Vert \left(\bm{\Sigma}_{U,l}^{\star}\right)^{-1/2}\left(\bm{U}\bm{R}-\bm{U}^{\star}-\bm{Z}\right)\right\Vert _{2,\infty}>\zeta\,\big|\,\bm{F}\right)\nonumber \\
 & \quad\leq\mathbb{P}\left(\left(\bm{\Sigma}_{U,l}^{\star}\right)^{-1/2}\left(\bm{U}\bm{R}-\bm{U}^{\star}\right)_{l,\cdot}\in\mathcal{C}\,\big|\,\bm{F}\right)+O\left(\left(n+d\right)^{-10}\right),\label{eq:proof-pca-normal-approximation-inter-3}
\end{align}
where the first inequality follows from the definition of $\mathcal{C}^{-\zeta}$,
and the last inequality makes use of \eqref{eq:UB-15278}. Similarly,
\begin{align}
\mathbb{P}\left(\left(\bm{\Sigma}_{U,l}^{\star}\right)^{-1/2}\left(\bm{U}\bm{R}-\bm{U}^{\star}\right)_{l,\cdot}\in\mathcal{C}\,\big|\,\bm{F}\right) & \leq\mathbb{P}\left(\bm{\Sigma}_{l}^{-1/2}\bm{Z}_{l,\cdot}\in\mathcal{C}^{\zeta}\,\big|\,\bm{F}\right)+O\left(\left(n+d\right)^{-10}\right).\label{eq:proof-pca-normal-approximation-inter-4}
\end{align}
In addition, for any set $\mathcal{X}\subseteq\mathbb{R}^{r}$ and
any matrix $\bm{A}\in\mathbb{R}^{r\times r}$, let us denote by $\bm{A}\mathcal{X}$
be the set $\{\bm{A}\bm{x}:\bm{x}\in\mathcal{X}\}$. It is easily
seen that when $\bm{A}$ is non-singular, $\bm{A}\mathcal{C}\subseteq\mathscr{C}^{r}$
holds if and only if $\mathcal{C}\in\mathscr{C}^{r}$. We can then
deduce that 
\begin{align*}
\mathbb{P}\left(\left(\bm{\Sigma}_{U,l}^{\star}\right)^{-1/2}\bm{Z}_{l,\cdot}\in\mathcal{C}^{\zeta}\,\big|\,\bm{F}\right) & =\mathbb{P}\left(\bm{Z}_{l,\cdot}\in\left(\bm{\Sigma}_{U,l}^{\star}\right)^{1/2}\mathcal{C}^{\zeta}\,\big|\,\bm{F}\right)\\
 & \overset{\text{(i)}}{\leq}\mathbb{P}\left(\mathcal{N}\left(\bm{0},\bm{\Sigma}_{U,l}^{\star}\right)\in\left(\bm{\Sigma}_{U,l}^{\star}\right)^{1/2}\mathcal{C}^{\zeta}\right)+O\left(\frac{1}{\sqrt{\log\left(n+d\right)}}\right)\\
 & =\mathbb{P}\left(\mathcal{N}\left(\bm{0},\bm{I}_{r}\right)\in\mathcal{C}^{\zeta}\right)+O\left(\frac{1}{\sqrt{\log\left(n+d\right)}}\right)\\
 & \overset{\text{(ii)}}{\leq}\mathbb{P}\left(\mathcal{N}\left(\bm{0},\bm{I}_{r}\right)\in\mathcal{C}\right)+\zeta\left(0.59r^{1/4}+0.21\right)+\left(\frac{1}{\sqrt{\log\left(n+d\right)}}\right)\\
 & \overset{\text{(iii)}}{\leq}\mathbb{P}\left(\mathcal{N}\left(\bm{0},\bm{I}_{r}\right)\in\mathcal{C}\right)+O\left(\frac{1}{\sqrt{\log\left(n+d\right)}}\right),
\end{align*}
where (i) uses (\ref{eq:proof-pca-normal-approximation-inter-2}),
(ii) is a consequence of Theorem \ref{thm:gaussian-perimeter}, and
(iii) holds provided that $\zeta\lesssim1/(r^{1/4}\log^{1/2}(n+d))$.
Similarly we can show that 
\[
\mathbb{P}\left(\left(\bm{\Sigma}_{U,l}^{\star}\right)^{-1/2}\bm{Z}_{l,\cdot}\in\mathcal{C}^{-\zeta}\,\big|\,\bm{F}\right)\geq\mathbb{P}\left(\mathcal{N}\left(\bm{0},\bm{I}_{r}\right)\in\mathcal{C}\right)-O\left(\frac{1}{\sqrt{\log\left(n+d\right)}}\right).
\]
Combine the above two inequalities with (\ref{eq:proof-pca-normal-approximation-inter-3})
and (\ref{eq:proof-pca-normal-approximation-inter-4}) to achieve
\[
\left|\mathbb{P}\left(\left(\bm{\Sigma}_{U,l}^{\star}\right)^{-1/2}\left(\bm{U}\bm{R}-\bm{U}^{\star}\right)_{l,\cdot}\in\mathcal{C}\,\big|\,\bm{F}\right)-\mathbb{P}\left(\mathcal{N}\left(\bm{0},\bm{I}_{r}\right)\in\mathcal{C}\right)\right|\lesssim\frac{1}{\sqrt{\log\left(n+d\right)}}.
\]
It is worth noting that this inequality holds for all $\mathcal{C}\in\mathscr{C}^{r}$.
As a result, on the event $\mathcal{E}_{\mathsf{good}}$ we can obtain
\begin{align}
 & \sup_{\mathcal{C}\in\mathscr{C}^{r}}\left|\mathbb{P}\left(\left(\bm{U}\bm{R}-\bm{U}^{\star}\right)_{l,\cdot}\in\mathcal{C}\,\big|\,\bm{F}\right)-\mathbb{P}\left(\mathcal{N}\left(\bm{0},\bm{\Sigma}_{l}\right)\in\mathcal{C}\right)\right|\nonumber \\
 & \quad=\sup_{\mathcal{C}\in\mathscr{C}^{r}}\left|\mathbb{P}\left(\left(\bm{U}\bm{R}-\bm{U}^{\star}\right)_{l,\cdot}\in\bm{\Sigma}_{l}^{1/2}\mathcal{C}\right)-\mathbb{P}\left(\mathcal{N}\left(\bm{0},\bm{\Sigma}_{l}\right)\in\bm{\Sigma}_{l}^{1/2}\mathcal{C}\right)\right|\nonumber \\
 & \quad=\sup_{\mathcal{C}\in\mathscr{C}^{r}}\left|\mathbb{P}\left(\bm{\Sigma}_{l}^{-1/2}\left(\bm{U}\bm{R}-\bm{U}^{\star}\right)_{l,\cdot}\in\mathcal{C}\right)-\mathbb{P}\left(\mathcal{N}\left(\bm{0},\bm{I}_{r}\right)\in\mathcal{C}\right)\right|\nonumber \\
 & \quad\lesssim\frac{1}{\sqrt{\log\left(n+d\right)}},\label{eq:proof-pca-normal-approximation-inter-5}
\end{align}
where the first identity makes use of the fact that $\mathcal{C}\to(\bm{\Sigma}_{U,l}^{\star})^{1/2}\mathcal{C}$
is a one-to one mapping from $\mathscr{C}^{r}$ to $\mathscr{C}^{r}$
(since $\bm{\Sigma}_{U,l}^{\star}$ has full rank).

It remains to verify the conditions required to guarantee $\zeta\lesssim1/(r^{1/4}\log^{1/2}(n+d))$.
More generally, we shall check that under what conditions we can guarantee
\begin{align*}
\zeta_{\mathsf{2nd},l} & \lesssim\delta\lambda_{\min}^{1/2}\left(\bm{\Sigma}_{U,l}^{\star}\right)
\end{align*}
for some $\delta>0$. Recall from Lemma \ref{lemma:pca-2nd-error}
that 
\begin{align*}
\zeta_{\mathsf{2nd},l} & \lesssim\underbrace{\left\Vert \bm{U}_{l,\cdot}^{\star}\right\Vert _{2}\sqrt{\frac{\kappa^{3}\mu r\log\left(n+d\right)}{d}}\frac{\zeta_{\mathsf{1st}}}{\sigma_{r}^{\star2}}}_{\eqqcolon\gamma_{1}}+\underbrace{\left\Vert \bm{U}_{l,\cdot}^{\star}\right\Vert _{2}\kappa\frac{\zeta_{\mathsf{1st}}^{2}}{\sigma_{r}^{\star4}}}_{\eqqcolon\gamma_{2}}+\underbrace{\frac{\zeta_{\mathsf{1st}}\zeta_{\mathsf{1st},l}}{\sigma_{r}^{\star4}}\sqrt{\frac{\kappa^{3}\mu r\log\left(n+d\right)}{d}}}_{\eqqcolon\gamma_{3}},
\end{align*}
and from Lemma \ref{lemma:pca-covariance-concentration} that 
\begin{align*}
\lambda_{\min}^{1/2}\left(\bm{\Sigma}_{U,l}^{\star}\right) & \gtrsim\frac{1}{\sqrt{np}\sigma_{1}^{\star}}\left\Vert \bm{U}_{l,\cdot}^{\star}\bm{\Sigma}^{\star}\right\Vert _{2}+\frac{\omega_{l}^{\star}}{\sqrt{np}\sigma_{1}^{\star}}+\frac{1}{\sqrt{ndp^{2}\kappa}\sigma_{1}^{\star}}\left\Vert \bm{U}_{l,\cdot}^{\star}\bm{\Sigma}^{\star}\right\Vert _{2}+\frac{1}{\sqrt{ndp^{2}\kappa}\sigma_{1}^{\star}}\omega_{l}^{\star}\\
 & \quad+\frac{1}{\sqrt{np^{2}}\sigma_{1}^{\star2}}\omega_{\min}\left\Vert \bm{U}_{l,\cdot}^{\star}\bm{\Sigma}^{\star}\right\Vert _{2}+\frac{\omega_{l}^{\star}\omega_{\min}}{\sqrt{np^{2}}\sigma_{1}^{\star2}}.
\end{align*}

\begin{itemize}
\item Regarding $\gamma_{1}$, we can derive 
\begin{align*}
\gamma_{1} & \asymp\underbrace{\sqrt{\frac{\kappa^{5}\mu^{3}r^{3}\log^{5}\left(n+d\right)}{nd^{2}p^{2}}}\left\Vert \bm{U}_{l,\cdot}^{\star}\right\Vert _{2}}_{\eqqcolon\gamma_{1,1}}+\underbrace{\frac{\omega_{\max}^{2}}{\sigma_{r}^{\star2}p}\sqrt{\frac{\kappa^{3}\mu r\log^{3}\left(n+d\right)}{n}}\left\Vert \bm{U}_{l,\cdot}^{\star}\right\Vert _{2}}_{\eqqcolon\gamma_{1,2}}\\
 & \quad+\underbrace{\left\Vert \bm{U}_{l,\cdot}^{\star}\right\Vert _{2}\sqrt{\frac{\kappa^{5}\mu^{2}r^{2}\log^{3}\left(n+d\right)}{ndp}}}_{\eqqcolon\gamma_{1,3}}+\underbrace{\left\Vert \bm{U}_{l,\cdot}^{\star}\right\Vert _{2}\frac{\omega_{\max}}{\sigma_{r}^{\star}}\sqrt{\frac{\kappa^{3}\mu r\log^{2}\left(n+d\right)}{np}}}_{\eqqcolon\gamma_{1,4}}\\
 & \lesssim\delta\lambda_{\min}^{1/2}\left(\bm{\Sigma}_{U,l}^{\star}\right),
\end{align*}
where the last line holds since 
\begin{align*}
\gamma_{1,1} & \lesssim\delta\frac{1}{\sqrt{ndp^{2}\kappa}\sigma_{1}^{\star}}\left\Vert \bm{U}_{l,\cdot}^{\star}\bm{\Sigma}^{\star}\right\Vert _{2},\\
\gamma_{1,2} & \lesssim\delta\frac{\omega_{l}^{\star}\omega_{\min}}{\sqrt{np^{2}}\sigma_{1}^{\star2}},\\
\gamma_{1,3} & \lesssim\delta\frac{1}{\sqrt{np}\sigma_{1}^{\star}}\left\Vert \bm{U}_{l,\cdot}^{\star}\bm{\Sigma}^{\star}\right\Vert _{2},\\
\gamma_{1,4} & \lesssim\delta\frac{\omega_{\min}}{\sqrt{np}\sigma_{1}^{\star}},
\end{align*}
provided that $d\gtrsim\delta^{-2}\kappa^{7}\mu^{3}r^{3}\kappa_{\omega}^{2}\log^{4}(n+d)$.
\item When it comes to $\gamma_{2}$, we observe that 
\begin{align*}
\gamma_{2} & \asymp\underbrace{\frac{\kappa^{3}\mu^{2}r^{2}\log^{4}\left(n+d\right)}{ndp^{2}}\left\Vert \bm{U}_{l,\cdot}^{\star}\right\Vert _{2}}_{\eqqcolon\gamma_{2,1}}+\underbrace{\kappa\frac{\omega_{\max}^{4}}{p^{2}\sigma_{r}^{\star4}}\frac{d}{n}\log^{2}\left(n+d\right)\left\Vert \bm{U}_{l,\cdot}^{\star}\right\Vert _{2}}_{\eqqcolon\gamma_{2,2}}\\
 & \quad+\underbrace{\frac{\kappa^{3}\mu r\log^{2}\left(n+d\right)}{np}\left\Vert \bm{U}_{l,\cdot}^{\star}\right\Vert _{2}}_{\eqqcolon\gamma_{2,3}}+\underbrace{\frac{\omega_{\max}^{2}}{\sigma_{r}^{\star2}}\frac{d\kappa^{2}\log\left(n+d\right)}{np}\left\Vert \bm{U}_{l,\cdot}^{\star}\right\Vert _{2}}_{\eqqcolon\gamma_{2,4}}\\
 & \lesssim\delta\lambda_{\min}^{1/2}\left(\bm{\Sigma}_{U,l}^{\star}\right),
\end{align*}
where the last line holds since 
\begin{align*}
\gamma_{2,1} & \lesssim\delta\frac{1}{\sqrt{ndp^{2}\kappa}\sigma_{1}^{\star}}\left\Vert \bm{U}_{l,\cdot}^{\star}\bm{\Sigma}^{\star}\right\Vert _{2},\\
\gamma_{2,2} & \lesssim\delta\frac{\omega_{\min}\omega_{l}^{\star}}{\sqrt{np^{2}}\sigma_{1}^{\star2}},\\
\gamma_{2,3} & \lesssim\delta\frac{1}{\sqrt{np}\sigma_{1}^{\star}}\left\Vert \bm{U}_{l,\cdot}^{\star}\bm{\Sigma}^{\star}\right\Vert _{2},\\
\gamma_{2,4} & \lesssim\delta\frac{\omega_{\min}}{\sqrt{np}\sigma_{1}^{\star}},
\end{align*}
provided that $ndp^{2}\gtrsim\delta^{-2}\kappa^{8}\mu^{4}r^{4}\log^{8}(n+d)$,
$np\gtrsim\delta^{-2}\kappa^{7}\mu^{2}r^{2}\log^{4}(n+d)$, 
\[
\frac{\omega_{\max}}{\sigma_{r}^{\star}}\sqrt{\frac{d}{np}}\lesssim\frac{\delta}{\kappa^{5/2}\mu^{1/2}r^{1/2}\kappa_{\omega}^{1/2}\log\left(n+d\right)},\qquad\frac{\omega_{\max}^{2}}{p\sigma_{r}^{\star2}}\sqrt{\frac{d}{n}}\lesssim\frac{\delta}{\kappa^{2}\mu^{1/2}r^{1/2}\kappa_{\omega}\log^{2}\left(n+d\right)}.
\]
\item We are left with $\gamma_{3}$, which can be bounded by
\begin{align*}
\gamma_{3} & \asymp\underbrace{\sqrt{\frac{\kappa^{5}\mu^{2}r^{2}\log^{4}\left(n+d\right)}{ndp^{2}}}\left\Vert \bm{U}_{l,\cdot}^{\star}\right\Vert _{2}\frac{\zeta_{\mathsf{1st}}}{\sigma_{r}^{\star2}}}_{\eqqcolon\gamma_{3,1}}+\underbrace{\frac{\omega_{l}^{\star}\omega_{\max}}{p\sigma_{r}^{\star2}}\sqrt{\frac{\kappa^{3}\mu r\log^{3}\left(n+d\right)}{n}}\frac{\zeta_{\mathsf{1st}}}{\sigma_{r}^{\star2}}}_{\eqqcolon\gamma_{3,2}}+\underbrace{\sqrt{\frac{\kappa^{5}\mu r\log^{3}\left(n+d\right)}{np}}\left\Vert \bm{U}_{l,\cdot}^{\star}\right\Vert _{2}\frac{\zeta_{\mathsf{1st}}}{\sigma_{r}^{\star2}}}_{\eqqcolon\gamma_{3,3}}\\
 & \quad+\underbrace{\frac{\omega_{l}^{\star}}{\sigma_{r}^{\star}}\sqrt{\frac{\kappa^{4}\mu r\log^{2}\left(n+d\right)}{np}}\frac{\zeta_{\mathsf{1st}}}{\sigma_{r}^{\star2}}}_{\eqqcolon\gamma_{3,4}}+\underbrace{\frac{\omega_{\max}}{p\sigma_{r}^{\star}}\sqrt{\frac{\kappa^{4}\mu r\log^{4}\left(n+d\right)}{n}}\left\Vert \bm{U}_{l,\cdot}^{\star}\right\Vert _{2}\frac{\zeta_{\mathsf{1st}}}{\sigma_{r}^{\star2}}}_{\eqqcolon\gamma_{3,5}}+\underbrace{\frac{\omega_{l}^{\star}}{p\sigma_{r}^{\star}}\sqrt{\frac{\kappa^{4}\mu^{2}r^{2}\log^{4}\left(n+d\right)}{nd}}\frac{\zeta_{\mathsf{1st}}}{\sigma_{r}^{\star2}}}_{\eqqcolon\gamma_{3,6}}.
\end{align*}
Similar to $\gamma_{1}$, we can show that $\gamma_{3,1}\lesssim\delta\lambda_{\min}^{1/2}(\bm{\Sigma}_{U,l}^{\star})$,
$\gamma_{3,3}\lesssim\delta\lambda_{\min}^{1/2}(\bm{\Sigma}_{U,l}^{\star})$
and $\gamma_{3,5}\lesssim\delta\lambda_{\min}^{1/2}(\bm{\Sigma}_{U,l}^{\star})$,
provided that $ndp^{2}\gtrsim\delta^{-2}\kappa^{9}\mu^{4}r^{4}\kappa_{\omega}^{2}\log^{7}(n+d)$,
$np\gtrsim\delta^{-2}\kappa^{9}\mu^{3}r^{3}\kappa_{\omega}^{2}\log^{6}(n+d)$,
and
\[
\frac{\omega_{\max}}{p\sigma_{r}^{\star}}\sqrt{\frac{1}{n}}\lesssim\delta\sqrt{\frac{1}{\kappa^{8}\mu^{3}r^{3}\kappa_{\omega}^{2}\log^{7}(n+d)}}.
\]
The last condition can be guaranteed by $ndp^{2}\gtrsim\delta^{-2}\kappa^{9}\mu^{4}r^{4}\kappa_{\omega}^{2}\log^{7}(n+d)$
and 
\[
\frac{\omega_{\max}^{2}}{p\sigma_{r}^{\star2}}\sqrt{\frac{d}{n}}\lesssim1.
\]
We are left with bounding $\gamma_{3,2}$, $\gamma_{3,4}$ and $\gamma_{3,6}$,
where we have
\begin{align*}
\gamma_{3,2} & \lesssim\delta\frac{\omega_{l}^{\star}\omega_{\min}}{\sqrt{np^{2}}\sigma_{1}^{\star2}}\lesssim\delta\lambda_{\min}^{1/2}\left(\bm{\Sigma}_{U,l}^{\star}\right),\\
\gamma_{3,4} & \lesssim\delta\frac{\omega_{l}^{\star}}{\sqrt{np}\sigma_{1}^{\star}}\lesssim\delta\lambda_{\min}^{1/2}\left(\bm{\Sigma}_{U,l}^{\star}\right),\\
\gamma_{3,6} & \lesssim\delta\frac{1}{\sqrt{ndp^{2}\kappa}\sigma_{1}^{\star}}\omega_{l}^{\star}\lesssim\delta\lambda_{\min}^{1/2}\left(\bm{\Sigma}_{U,l}^{\star}\right),
\end{align*}
provided that
\[
\frac{\zeta_{\mathsf{1st}}}{\sigma_{r}^{\star2}}\lesssim\frac{\delta}{\sqrt{\kappa^{6}\mu^{2}r^{2}\kappa_{\omega}\log^{4}\left(n+d\right)}}.
\]
In view of (\ref{subeq:pca-1st-equivalent}), the above condition
is equivalent to
\begin{align*}
ndp^{2}\gtrsim\delta^{-2}\kappa^{8}\mu^{4}r^{4}\kappa_{\omega}\log^{8}\left(n+d\right), & \qquad np\gtrsim\delta^{-2}\kappa^{8}\mu^{3}r^{3}\kappa_{\omega}\log^{6}\left(n+d\right),\\
\text{and}\qquad\frac{\omega_{\max}^{2}}{p\sigma_{r}^{\star2}}\sqrt{\frac{d}{n}}\lesssim\frac{\delta}{\sqrt{\kappa^{6}\mu^{2}r^{2}\kappa_{\omega}\log^{6}\left(n+d\right)}}, & \qquad\frac{\omega_{\max}}{\sigma_{r}^{\star}}\sqrt{\frac{d}{np}}\lesssim\frac{\delta}{\sqrt{\kappa^{7}\mu^{2}r^{2}\kappa_{\omega}\log^{5}\left(n+d\right)}}.
\end{align*}
Therefore we have
\[
\gamma_{3}\lesssim\delta\lambda_{\min}^{1/2}\left(\bm{\Sigma}_{U,l}^{\star}\right)
\]
under the conditions
\begin{align*}
ndp^{2}\gtrsim\delta^{-2}\kappa^{9}\mu^{4}r^{4}\kappa_{\omega}^{2}\log^{8}\left(n+d\right), & \qquad np\gtrsim\delta^{-2}\kappa^{9}\mu^{3}r^{3}\kappa_{\omega}^{2}\log^{6}\left(n+d\right),\\
\text{and}\qquad\frac{\omega_{\max}^{2}}{p\sigma_{r}^{\star2}}\sqrt{\frac{d}{n}}\lesssim\frac{\delta}{\sqrt{\kappa^{6}\mu^{2}r^{2}\kappa_{\omega}\log^{6}\left(n+d\right)}}, & \qquad\frac{\omega_{\max}}{\sigma_{r}^{\star}}\sqrt{\frac{d}{np}}\lesssim\frac{\delta}{\sqrt{\kappa^{7}\mu^{2}r^{2}\kappa_{\omega}\log^{5}\left(n+d\right)}}.
\end{align*}
\end{itemize}
Therefore, $\zeta_{\mathsf{2nd},l}\lesssim\delta\lambda_{\min}^{1/2}\left(\bm{\Sigma}_{U,l}^{\star}\right)$
is guaranteed to hold as long as $d\gtrsim\delta^{-2}\kappa^{7}\mu^{3}r^{3}\kappa_{\omega}^{2}\log^{4}(n+d)$,
\begin{align*}
ndp^{2}\gtrsim\delta^{-2}\kappa^{9}\mu^{4}r^{4}\kappa_{\omega}^{2}\log^{8}\left(n+d\right), & \qquad np\gtrsim\delta^{-2}\kappa^{9}\mu^{3}r^{3}\kappa_{\omega}^{2}\log^{6}\left(n+d\right),\\
\text{and}\qquad\frac{\omega_{\max}^{2}}{p\sigma_{r}^{\star2}}\sqrt{\frac{d}{n}}\lesssim\frac{\delta}{\sqrt{\kappa^{6}\mu^{2}r^{2}\kappa_{\omega}\log^{6}\left(n+d\right)}}, & \qquad\frac{\omega_{\max}}{\sigma_{r}^{\star}}\sqrt{\frac{d}{np}}\lesssim\frac{\delta}{\sqrt{\kappa^{7}\mu^{2}r^{2}\kappa_{\omega}\log^{5}\left(n+d\right)}}.
\end{align*}
By taking $\delta=1/(r^{1/4}\log^{1/2}(n+d))$, we see that $\zeta\lesssim1/(r^{1/4}\log^{1/2}(n+d))$
holds provided that $d\gtrsim\kappa^{7}\mu^{3}r^{7/2}\kappa_{\omega}^{2}\log^{5}(n+d)$,
\begin{align*}
ndp^{2}\gtrsim\kappa^{9}\mu^{4}r^{9/2}\kappa_{\omega}^{2}\log^{9}\left(n+d\right), & \qquad np\gtrsim\kappa^{9}\mu^{3}r^{7/2}\kappa_{\omega}^{2}\log^{7}\left(n+d\right),\\
\text{and}\qquad\frac{\omega_{\max}^{2}}{p\sigma_{r}^{\star2}}\sqrt{\frac{d}{n}}\lesssim\frac{1}{\sqrt{\kappa^{6}\mu^{2}r^{5/2}\kappa_{\omega}\log^{7}\left(n+d\right)}}, & \qquad\frac{\omega_{\max}}{\sigma_{r}^{\star}}\sqrt{\frac{d}{np}}\lesssim\frac{1}{\sqrt{\kappa^{7}\mu^{2}r^{5/2}\kappa_{\omega}\log^{6}\left(n+d\right)}}.
\end{align*}

\paragraph{Step 4: distributional characterization of $\bm{Z}_{l,\cdot}$.}

For any convex set $\mathcal{C}\in\mathscr{C}^{r}$, it holds that
\begin{align*}
 & \left|\mathbb{P}\left(\left(\bm{U}\bm{R}-\bm{U}^{\star}\right)_{l,\cdot}\in\mathcal{C}\right)-\mathbb{P}\left(\mathcal{N}\left(\bm{0},\bm{\Sigma}_{U,l}^{\star}\right)\in\mathcal{C}\right)\right|\\
 & \quad=\left|\mathbb{E}\left[\mathbb{P}\left(\left(\bm{U}\bm{R}-\bm{U}^{\star}\right)_{l,\cdot}\in\mathcal{C}\big|\bm{F}\right)-\mathbb{P}\left(\mathcal{N}\left(\bm{0},\bm{\Sigma}_{U,l}^{\star}\right)\in\mathcal{C}\right)\right]\right|\\
 & \quad\leq\left|\mathbb{E}\left[\left[\mathbb{P}\left(\left(\bm{U}\bm{R}-\bm{U}^{\star}\right)_{l,\cdot}\in\mathcal{C}\big|\bm{F}\right)-\mathbb{P}\left(\mathcal{N}\left(\bm{0},\bm{\Sigma}_{U,l}^{\star}\right)\in\mathcal{C}\right)\right]\ind_{\mathcal{E}_{\mathsf{good}}}\right]\right|\\
 & \quad\quad+\left|\mathbb{E}\left[\left[\mathbb{P}\left(\left(\bm{U}\bm{R}-\bm{U}^{\star}\right)_{l,\cdot}\in\mathcal{C}\big|\bm{F}\right)-\mathbb{P}\left(\mathcal{N}\left(\bm{0},\bm{\Sigma}_{U,l}^{\star}\right)\in\mathcal{C}\right)\right]\ind_{\mathcal{E}_{\mathsf{good}}^{\mathsf{c}}}\right]\right|\\
 & \quad\leq\frac{1}{\sqrt{\log\left(n+d\right)}}+2\mathbb{P}\left(\mathcal{E}_{\mathsf{good}}^{\mathsf{c}}\right)\\
 & \quad\lesssim\frac{1}{\sqrt{\log\left(n+d\right)}}+\frac{1}{\left(n+d\right)^{100}}\lesssim\frac{1}{\sqrt{\log\left(n+d\right)}},
\end{align*}
where the penultimate line relies on (\ref{eq:proof-pca-normal-approximation-inter-5}).
This allows one to conclude that 
\[
\sup_{\mathcal{C}\in\mathscr{C}^{r}}\left|\mathbb{P}\left(\bm{Z}_{l,\cdot}\in\mathcal{C}\right)-\mathbb{P}\left(\mathcal{N}\left(\bm{0},\bm{\Sigma}_{U,l}^{\star}\right)\in\mathcal{C}\right)\right|\lesssim\frac{1}{\sqrt{\log\left(n+d\right)}}=o\left(1\right).
\]

\subsection{Auxiliary lemmas for Theorem \ref{thm:pca-cr-complete}}

\subsubsection{Proof of Lemma \ref{lemma:pca-1st-err} \label{appendix:proof-pca-1st-err}}

This lemma is concerned with several different quantities, which we
seek to control separately.

\paragraph{Bounding $\Vert(\bm{U}\bm{R}-\bm{U}^{\star})_{l,\cdot}\Vert_{2}$
and $\Vert\bm{U}\bm{R}-\bm{U}^{\star}\Vert_{2,\infty}$.}

In view of Lemma~\ref{lemma:pca-2nd-error}, we can begin with the
decomposition 
\begin{align*}
\left\Vert \bm{U}_{l,\cdot}\bm{R}-\bm{U}_{l,\cdot}^{\star}\right\Vert _{2} & \leq\left\Vert \bm{Z}_{l,\cdot}\right\Vert _{2}+\zeta_{\mathsf{2nd},l}=\left\Vert \left[\bm{E}\bm{M}^{\natural\top}+\mathcal{P}_{\mathsf{off}\text{-}\mathsf{diag}}\left(\bm{E}\bm{E}^{\top}\right)\right]_{l,\cdot}\bm{U}^{\natural}\left(\bm{\Sigma}^{\natural}\right)^{-2}\bm{Q}^{\top}\right\Vert _{2}+\zeta_{\mathsf{2nd},l}\\
 & \leq\underbrace{\left\Vert \bm{E}_{l,\cdot}\bm{V}^{\natural}(\bm{\Sigma}^{\natural})^{-1}\bm{Q}^{\top}\right\Vert _{2}}_{\eqqcolon\alpha_{1}}+\underbrace{\left\Vert \bm{E}_{l,\cdot}\left[\mathcal{P}_{-l,\cdot}\left(\bm{E}\right)\right]^{\top}\bm{U}^{\natural}\left(\bm{\Sigma}^{\natural}\right)^{-2}\bm{Q}^{\top}\right\Vert _{2}}_{\eqqcolon\alpha_{2}}+\zeta_{\mathsf{2nd},l},
\end{align*}
where $\bm{Z}$ and $\zeta_{\mathsf{2nd},l}$ are defined in \eqref{eq:defn-Z-distribution-hpca}. 

We first bound $\alpha_{1}$. One can decompose
\[
\bm{E}_{l,\cdot}\bm{V}^{\natural}(\bm{\Sigma}^{\natural})^{-1}\bm{Q}^{\top}=\sum_{j=1}^{n}E_{l,j}\bm{V}_{j,\cdot}^{\natural}(\bm{\Sigma}^{\natural})^{-1}\bm{Q}^{\top},
\]
which is an independent sum of random vectors (where the randomness
comes from $\{E_{l,j}\}_{1\leq j\leq n}$) conditional on $\bm{F}$.
By carrying out the following calculation (see \eqref{eq:E-expression-denoising}
and \eqref{eq:sigma-ij-square-denoising-PCA}) 
\begin{align*}
L & \coloneqq\max_{1\leq j\leq n}\left\Vert E_{l,j}\bm{V}_{j,\cdot}^{\natural}(\bm{\Sigma}^{\natural})^{-1}\bm{Q}^{\top}\right\Vert _{2}\leq\left\{ \max_{1\leq j\leq n}|E_{l,j}|\right\} \left\Vert \bm{V}^{\natural}(\bm{\Sigma}^{\natural})^{-1}\right\Vert _{2,\infty}\\
 & \lesssim\frac{1}{\sqrt{n}p}\left(\max_{1\leq j\leq n}\left|\bm{U}_{l,\cdot}^{\star}\bm{\Sigma}^{\star}\bm{f}_{j}\right|+\omega_{l}^{\star}\sqrt{\log\left(n+d\right)}\right)\left\Vert \bm{V}^{\natural}(\bm{\Sigma}^{\natural})^{-1}\right\Vert _{2,\infty},\\
V & \coloneqq\sum_{j=1}^{n}\mathbb{E}\left[E_{l,j}^{2}\left\Vert \bm{V}_{j,\cdot}^{\natural}(\bm{\Sigma}^{\natural})^{-1}\bm{Q}^{\top}\right\Vert _{2}^{2}\right]\leq\sum_{j=1}^{n}\mathbb{E}\left[E_{l,j}^{2}\right]\left\Vert \bm{V}_{j,\cdot}^{\natural}(\bm{\Sigma}^{\natural})^{-1}\right\Vert _{2}^{2}\leq\max_{j}\mathbb{E}\left[E_{l,j}^{2}\right]\left\Vert \bm{V}^{\natural}(\bm{\Sigma}^{\natural})^{-1}\right\Vert _{\mathrm{F}}^{2}\\
 & \lesssim\frac{1}{np}\left[\max_{1\leq j\leq n}\left(\bm{U}_{l,\cdot}^{\star}\bm{\Sigma}^{\star}\bm{f}_{j}\right)^{2}+\omega_{l}^{\star2}\right]\left\Vert \bm{V}^{\natural}(\bm{\Sigma}^{\natural})^{-1}\right\Vert _{\mathrm{F}}^{2},
\end{align*}
we can invoke the Bernstein inequality \citep[Corollary 3.1.3]{chen2020spectral}
to demonstrate that 
\begin{align*}
\left\Vert \bm{Z}_{l,\cdot}\right\Vert _{2} & \lesssim\sqrt{V\log\left(n+d\right)}+L\log\left(n+d\right)\\
 & \lesssim\sqrt{\frac{\log\left(n+d\right)}{np}}\left[\max_{1\leq j\leq n}\left|\bm{U}_{l,\cdot}^{\star}\bm{\Sigma}^{\star}\bm{f}_{j}\right|+\omega_{l}^{\star}\right]\left\Vert \bm{V}^{\natural}(\bm{\Sigma}^{\natural})^{-1}\right\Vert _{\mathrm{F}}\\
 & \quad+\frac{\log\left(n+d\right)}{\sqrt{n}p}\left(\max_{1\leq j\leq n}\left|\bm{U}_{l,\cdot}^{\star}\bm{\Sigma}^{\star}\bm{f}_{j}\right|+\omega_{l}^{\star}\sqrt{\log\left(n+d\right)}\right)\left\Vert \bm{V}^{\natural}(\bm{\Sigma}^{\natural})^{-1}\right\Vert _{2,\infty}
\end{align*}
with probability at least $1-O((n+d)^{-10})$. On the event $\mathcal{E}_{\mathsf{good}}$,
we can further derive 
\begin{align*}
\left\Vert \bm{Z}_{l,\cdot}\right\Vert _{2} & \overset{\text{(i)}}{\lesssim}\frac{1}{\sigma_{r}^{\star}}\sqrt{\frac{r\log\left(n+d\right)}{np}}\left[\left\Vert \bm{U}_{l,\cdot}^{\star}\bm{\Sigma}^{\star}\right\Vert _{2}\sqrt{\log\left(n+d\right)}+\omega_{l}^{\star}\right]\\
 & \quad+\frac{1}{\sigma_{r}^{\star}}\frac{\log\left(n+d\right)}{\sqrt{n}p}\left(\left\Vert \bm{U}_{l,\cdot}^{\star}\bm{\Sigma}^{\star}\right\Vert _{2}\sqrt{\log\left(n+d\right)}+\omega_{l}^{\star}\sqrt{\log\left(n+d\right)}\right)\sqrt{\frac{\log\left(n+d\right)}{n}}\\
 & \asymp\frac{1}{\sigma_{r}^{\star}}\left\Vert \bm{U}_{l,\cdot}^{\star}\bm{\Sigma}^{\star}\right\Vert _{2}\left(\sqrt{\frac{r\log^{2}\left(n+d\right)}{np}}+\frac{\log^{3/2}\left(n+d\right)}{np}\right)+\frac{\omega_{l}^{\star}}{\sigma_{r}^{\star}}\left(\sqrt{\frac{r\log\left(n+d\right)}{np}}+\frac{\log^{2}\left(n+d\right)}{np}\right)\\
 & \overset{\text{(ii)}}{\lesssim}\frac{1}{\sigma_{r}^{\star}}\left\Vert \bm{U}_{l,\cdot}^{\star}\bm{\Sigma}^{\star}\right\Vert _{2}\sqrt{\frac{r\log^{2}\left(n+d\right)}{np}}+\frac{\omega_{l}^{\star}}{\sigma_{r}^{\star}}\sqrt{\frac{r\log\left(n+d\right)}{np}}
\end{align*}
where (i) uses (\ref{eq:good-event-sigma-least-largest-hpca}), (\ref{eq:good-event-Ui-Sigma-f-hpca})
as well as (\ref{eq:good-event-V-natural-2-infty-hpca}), and (ii)
holds true as long as $np\gtrsim\log^{3}(n+d)$.

\[
\alpha_{1}\lesssim\frac{1}{\sigma_{r}^{\star}}\left(\left\Vert \bm{U}_{l,\cdot}^{\star}\bm{\Sigma}^{\star}\right\Vert _{2}+\omega_{l}^{\star}\right)\sqrt{\frac{r\log^{2}\left(n+d\right)}{np}}.
\]

Next we bound $\alpha_{2}$. In order to do so, we know from (\ref{eq:good-event-sigma-least-largest-hpca})
that on the event $\mathcal{E}_{\mathsf{good}}$, 
\begin{equation}
\alpha_{2}\ind_{\mathcal{E}_{\mathsf{good}}}\lesssim\frac{1}{\sigma_{r}^{\star2}}\left\Vert \bm{E}_{l,\cdot}\left[\mathcal{P}_{-l,\cdot}\left(\bm{E}\right)\right]^{\top}\bm{U}^{\natural}\right\Vert _{2}\asymp\frac{1}{\sigma_{r}^{\star2}}\left\Vert \bm{E}_{l,\cdot}\left[\mathcal{P}_{-l,\cdot}\left(\bm{E}\right)\right]^{\top}\bm{U}^{\star}\right\Vert _{2},\label{eq:alpha2-UB-125}
\end{equation}
where the last relation arises from the fact $\bm{U}^{\natural}=\bm{U}^{\star}\bm{Q}$
for some orthonormal matrix $\bm{Q}$ (see \eqref{eq:defn-Q-denoising-PCA}).
Therefore, it suffices to bound 
\[
\left\Vert \bm{E}_{l,\cdot}\left[\mathcal{P}_{-l,\cdot}\left(\bm{E}\right)\right]^{\top}\bm{U}^{\star}\right\Vert _{2}=\bigg\|\sum_{j=1}^{n}E_{l,j}\bm{C}_{j,\cdot}\bm{U}^{\star}\bigg\|_{2}=\bigg\|\sum_{j=1}^{n}\bm{X}_{j}\bigg\|_{2},
\]
where $\bm{C}=[\mathcal{P}_{-l,\cdot}(\bm{E})]^{\top}$ and $\bm{X}_{j}=E_{l,j}\bm{C}_{j,\cdot}\bm{U}^{\star}$
for $j\in[n]$. Recognizing that $\sum_{j=1}^{n}\bm{X}_{j}$ is a
sum of independent random vectors (conditional on $\bm{F}$), we can
employ the truncated Bernstein inequality (see, e.g., \citet[Theorem 3.1.1]{chen2020spectral})
to bound the above quantity. From now on, we will always assume occurrence
of $\mathcal{E}_{\mathsf{good}}$ when bounding $\Vert(\bm{U}\bm{R}-\bm{U}^{\star}\Vert_{2}$
(recall that $\mathcal{E}_{\mathsf{good}}$ is $\sigma(\bm{F})$-measurable). 
\begin{itemize}
\item It is first observed that 
\begin{align*}
\max_{1\leq j\leq n}\left\Vert \bm{X}_{j}\right\Vert _{2} & \leq\max_{1\leq j\leq n}\left|E_{l,j}\right|\left\Vert \bm{C}_{j,\cdot}\bm{U}^{\star}\right\Vert _{2}.
\end{align*}
Recall that for each $j\in[n]$, 
\[
E_{l,j}=\frac{1}{\sqrt{n}p}\left[\left(\delta_{l,j}-1\right)\bm{U}_{l,\cdot}^{\star}\bm{\Sigma}^{\star}\bm{f}_{j}+N_{l,j}\right]
\]
with $\delta_{l,j}=\ind_{(l,j)\in\Omega}$. It is seen from (\ref{eq:good-event-Ui-Sigma-f-hpca})
that 
\[
\left|E_{l,j}\right|\leq\widetilde{C}_{1}\sqrt{\frac{1}{np^{2}}}\left(\left\Vert \bm{U}_{l,\cdot}^{\star}\bm{\Sigma}^{\star}\right\Vert _{2}+\omega_{l}^{\star}\right)\sqrt{\log\left(n+d\right)}\eqqcolon L_{E}
\]
for some sufficiently large constant $\widetilde{C}_{1}>0$. In addition,
for each $j\in[n]$, we can write $\bm{C}_{j,\cdot}\bm{U}^{\star}=\sum_{i:i\neq l}E_{i,j}\bm{U}_{i,\cdot}^{\star}$.
It is then straightforward to compute 
\begin{align*}
L_{1} & \coloneqq\max_{i:i\neq l}\left\Vert E_{i,j}\bm{U}_{i,\cdot}^{\star}\right\Vert _{2}\lesssim\sqrt{\frac{\log\left(n+d\right)}{np^{2}}}\left(\left\Vert \bm{U}^{\star}\bm{\Sigma}^{\star}\right\Vert _{2,\infty}+\omega_{\max}\right)\sqrt{\frac{\mu r}{d}}\\
 & \asymp\sqrt{\frac{\mu r\log\left(n+d\right)}{ndp^{2}}}\left(\left\Vert \bm{U}^{\star}\bm{\Sigma}^{\star}\right\Vert _{2,\infty}+\omega_{\max}\right)\lesssim\sqrt{\frac{\mu r\log\left(n+d\right)}{dp}}\sigma_{\mathsf{ub}},\\
V_{1} & \coloneqq\sum_{i:i\neq l}\mathbb{E}\left[E_{i,j}^{2}\left\Vert \bm{U}_{i,\cdot}^{\star}\right\Vert _{2}^{2}\right]\leq\sum_{i:i\neq l}\sigma_{i,j}^{2}\left\Vert \bm{U}_{i,\cdot}^{\star}\right\Vert _{2}^{2}\leq\sigma_{\mathsf{ub}}^{2}\left\Vert \bm{U}^{\star}\right\Vert _{\mathrm{F}}^{2}=r\sigma_{\mathsf{ub}}^{2},
\end{align*}
where $\sigma_{\mathsf{ub}}$ is defined in \eqref{eq:good-event-sigma-hpca}.
Apply the matrix Bernstein inequality \citep[Theorem 6.1.1]{Tropp:2015:IMC:2802188.2802189}
to achieve 
\[
\mathbb{P}\left(\left\Vert \bm{C}_{j,\cdot}\bm{U}^{\star}\right\Vert _{2}\geq t\,\big|\,\bm{F}\right)\leq\left(d+1\right)\exp\left(\frac{-t^{2}/2}{V_{1}+L_{1}t/3}\right).
\]
If we take 
\[
L_{C}\coloneqq\widetilde{C}_{2}\left[\sqrt{V_{1}\log\left(n+d\right)}+L_{1}\log\left(n+d\right)\right].
\]
for some sufficiently large constant $\widetilde{C}_{2}>0$, then
the above inequality tells us that for each $j\in[n]$ 
\[
\mathbb{P}\left(\left\Vert \bm{C}_{j,\cdot}\bm{U}^{\star}\right\Vert _{2}\geq L_{C}\,\big|\,\bm{F}\right)\leq\left(n+d\right)^{-20}.
\]
Consequently, by setting 
\begin{align*}
L & \coloneqq L_{E}L_{C}\asymp\sqrt{\frac{\log\left(n+d\right)}{np^{2}}}\left(\left\Vert \bm{U}_{l,\cdot}^{\star}\bm{\Sigma}^{\star}\right\Vert _{2}+\omega_{l}^{\star}\right)\left[\sqrt{V_{1}\log\left(n+d\right)}+L_{1}\log\left(n+d\right)\right]\\
 & \lesssim\sqrt{\frac{r\log^{2}\left(n+d\right)}{np^{2}}}\sigma_{\mathsf{ub}}\left(\left\Vert \bm{U}_{l,\cdot}^{\star}\bm{\Sigma}^{\star}\right\Vert _{2}+\omega_{l}^{\star}\right)\left[1+\sqrt{\frac{\mu}{dp}}\right],
\end{align*}
we can further derive 
\begin{align*}
\mathbb{P}\left(\left\Vert \bm{X}_{j}\right\Vert \geq L\,\big|\,\bm{F}\right) & \leq\mathbb{P}\left(\left\Vert \bm{C}_{j,\cdot}\bm{U}^{\star}\right\Vert _{2}\geq L_{C}\,\big|\,\bm{F}\right)\leq\left(n+d\right)^{-20}\coloneqq q_{0}.
\end{align*}
\item Next, we can develop the following upper bound 
\begin{align*}
q_{1} & \coloneqq\left\Vert \mathbb{E}\left[\bm{X}_{j}\ind_{\left\Vert \bm{X}_{j}\right\Vert \geq L}\,\big|\,\bm{F}\right]\right\Vert \overset{\text{(i)}}{\leq}\mathbb{E}\left[\left\Vert \bm{X}_{j}\right\Vert \ind_{\left\Vert \bm{X}_{j}\right\Vert \geq L}\,\big|\,\bm{F}\right]\\
 & \overset{\text{(ii)}}{\leq}L_{E}\mathbb{E}\left[\left\Vert \bm{C}_{j,\cdot}\bm{U}^{\star}\right\Vert _{2}\ind_{\left\Vert \bm{C}_{j,\cdot}\bm{U}^{\star}\right\Vert _{2}\geq L_{C}}\,\big|\,\bm{F}\right]\\
 & =L_{E}\int_{0}^{\infty}\mathbb{P}\left(\left\Vert \bm{C}_{j,\cdot}\bm{U}^{\star}\right\Vert _{2}\ind_{\left\Vert \bm{C}_{j,\cdot}\bm{U}^{\star}\right\Vert _{2}\geq L_{C}}\geq t\,\big|\,\bm{F}\right)\mathrm{d}t\\
 & =L_{E}\int_{0}^{L_{C}}\mathbb{P}\left(\left\Vert \bm{C}_{j,\cdot}\bm{U}^{\star}\right\Vert _{2}\geq L_{C}\,\big|\,\bm{F}\right)\mathrm{d}t+L_{E}\int_{L_{C}}^{\infty}\mathbb{P}\left(\left\Vert \bm{C}_{j,\cdot}\bm{U}^{\star}\right\Vert _{2}\geq t\,\big|\,\bm{F}\right)\mathrm{d}t\\
 & \leq L_{E}L_{C}\left(n+d\right)^{-20}+L_{E}\int_{L_{C}}^{\infty}\mathbb{P}\left(\left\Vert \bm{C}_{j,\cdot}\bm{U}^{\star}\right\Vert _{2}\geq t\,\big|\,\bm{F}\right)\mathrm{d}t.
\end{align*}
Here, (i) makes use of Jensen's inequality, while (ii) holds since
the conditions $\Vert\bm{X}_{j}\Vert\geq L$ and $\vert E_{l,j}\vert\leq L_{E}$
taken together imply that $\Vert\bm{C}_{j,\cdot}\bm{U}^{\star}\Vert_{2}\geq L_{C}$.
Note that for any $t\geq L_{C}$, we have $t\gg\sqrt{V_{1}\log(n+d)}$
and $t\gg L_{1}\log(n+d)$ as long as $\widetilde{C}_{2}$ is sufficiently
large. Therefore, we arrive at 
\[
\mathbb{P}\left(\left\Vert \bm{C}_{j,\cdot}\bm{U}^{\star}\right\Vert _{2}\geq t\,\big|\,\bm{F}\right)\leq\left(d+1\right)\exp\left(-t/\max\left\{ 4\sqrt{V_{1}/\log\left(n+d\right)},4L_{1}/3\right\} \right),
\]
which immediately gives 
\begin{align*}
 & \int_{L_{C}}^{\infty}\mathbb{P}\left(\left\Vert \bm{C}_{j,\cdot}\bm{U}^{\star}\right\Vert _{2}\geq t\,\big|\,\bm{F}\right)\mathrm{d}t\leq\left(d+1\right)\int_{L_{C}}^{\infty}\exp\left(-t/\max\left\{ 4\sqrt{V_{1}/\log\left(n+d\right)},4L_{1}/3\right\} \right)\mathrm{d}t\\
 & \quad\leq\left(d+1\right)\max\left\{ 4\sqrt{V_{1}/\log\left(n+d\right)},4L_{1}/3\right\} \exp\left(-L_{C}/\max\left\{ 4\sqrt{V_{1}/\log\left(n+d\right)},4L_{1}/3\right\} \right)\\
 & \quad\leq4\left(d+1\right)\widetilde{C}_{2}^{-1}L_{C}\exp\left(-4\widetilde{C}_{2}\log\left(n+d\right)\right)\leq L_{C}\left(n+d\right)^{-20},
\end{align*}
provided that $\widetilde{C}_{2}$ is sufficiently large. As a consequence,
we reach 
\[
q_{1}\leq2L_{E}L_{C}\left(n+d\right)^{-20}\leq L\left(n+d\right)^{-19}.
\]
\item Finally, let us calculate the variance statistics as follows 
\begin{align*}
v & \coloneqq\sum_{j=1}^{n}\mathbb{E}\left[\left\Vert \bm{X}_{j}\right\Vert _{2}^{2}\,\big|\,\bm{F}\right]=\sum_{j=1}^{n}\mathbb{E}\left[\left\Vert E_{l,j}\bm{C}_{j,\cdot}\bm{U}^{\star}\right\Vert _{2}^{2}\,\big|\,\bm{F}\right]=\sum_{j=1}^{n}\sigma_{l,j}^{2}\mathbb{E}\left[\left\Vert \bm{C}_{j,\cdot}\bm{U}^{\star}\right\Vert _{2}^{2}\,\big|\,\bm{F}\right]\\
 & \leq\max_{1\leq j\leq n}\sigma_{l,j}^{2}\mathbb{E}\left[\left\Vert \bm{C}\bm{U}^{\star}\right\Vert _{\mathrm{F}}^{2}\,\big|\,\bm{F}\right]=\max_{1\leq j\leq n}\sigma_{l,j}^{2}\mathbb{E}\left[\left\Vert \left[\mathcal{P}_{-l,\cdot}\left(\bm{E}\right)\right]^{\top}\bm{U}^{\star}\right\Vert _{\mathrm{F}}^{2}\,\big|\,\bm{F}\right]\\
 & \leq\frac{\left\Vert \bm{U}_{l,\cdot}^{\star}\bm{\Sigma}^{\star}\right\Vert _{2}^{2}\log\left(n+d\right)+\omega_{l}^{\star2}}{np}\mathbb{E}\left[\left\Vert \left[\mathcal{P}_{-l,\cdot}\left(\bm{E}\right)\right]^{\top}\bm{U}^{\star}\right\Vert _{\mathrm{F}}^{2}\,\big|\,\bm{F}\right],
\end{align*}
where the last inequality results from the following relation 
\[
\max_{1\leq j\leq n}\sigma_{l,j}^{2}\asymp\frac{1-p}{np}\max_{1\leq j\leq n}\left(\bm{U}_{l,\cdot}^{\star}\bm{\Sigma}^{\star}\bm{f}_{j}\right)^{2}+\frac{\omega_{l}^{\star2}}{np}\lesssim\frac{\big\|\bm{U}_{l,\cdot}^{\star}\bm{\Sigma}^{\star}\big\|_{2}^{2}\log\left(n+d\right)+\omega_{l}^{\star2}}{np}.
\]
Notice that 
\begin{align*}
\mathbb{E}\left[\left\Vert \left[\mathcal{P}_{-l,\cdot}\left(\bm{E}\right)\right]^{\top}\bm{U}^{\star}\right\Vert _{\mathrm{F}}^{2}\,\big|\,\bm{F}\right] & =\mathbb{E}\left[\mathsf{tr}\left(\bm{U}^{\star\top}\mathcal{P}_{-l,\cdot}\left(\bm{E}\right)\left[\mathcal{P}_{-l,\cdot}\left(\bm{E}\right)\right]^{\top}\bm{U}^{\star}\right)\,\big|\,\bm{F}\right]\\
 & =\mathbb{E}\left[\mathsf{tr}\left(\mathcal{P}_{-l,\cdot}\left(\bm{E}\right)\left[\mathcal{P}_{-l,\cdot}\left(\bm{E}\right)\right]^{\top}\bm{U}^{\star}\bm{U}^{\star\top}\right)\,\big|\,\bm{F}\right]\\
 & =\mathsf{tr}\left\{ \bm{D}\bm{U}^{\star}\bm{U}^{\star\top}\right\} ,
\end{align*}
where $\bm{D}=\mathbb{E}\big[\mathcal{P}_{-l,\cdot}(\bm{E})[\mathcal{P}_{-l,\cdot}(\bm{E})]^{\top}|\bm{F}\big]\in\mathbb{R}^{d\times d}$
is a diagonal matrix with the $i$-th diagonal entry given by 
\[
D_{i,i}=\begin{cases}
\sum_{j=1}^{n}\sigma_{i,j}^{2} & \text{if }i\neq l,\\
0 & \text{if }i=l.
\end{cases}
\]
As a result, we have demonstrated that 
\begin{align*}
\mathbb{E}\left[\left\Vert \left[\mathcal{P}_{-l,\cdot}\left(\bm{E}\right)\right]^{\top}\bm{U}^{\star}\right\Vert _{\mathrm{F}}^{2}\,\big|\,\bm{F}\right] & =\mathsf{tr}\left\{ \bm{D}\bm{U}^{\star}\bm{U}^{\star\top}\right\} =\sum_{i=1}^{d}D_{i,i}\left\Vert \bm{U}_{i,\cdot}^{\star}\right\Vert _{2}^{2}\\
 & \leq\sum_{i=1}^{d}\left(\sum_{j=1}^{n}\sigma_{i,j}^{2}\right)\left\Vert \bm{U}_{i,\cdot}^{\star}\right\Vert _{2}^{2}\leq n\sigma_{\mathsf{ub}}^{2}\left\Vert \bm{U}^{\star}\right\Vert _{\mathrm{F}}^{2}=nr\sigma_{\mathsf{ub}}^{2}.
\end{align*}
Taking the above inequalities collectively yields 
\[
v\lesssim\frac{r\sigma_{\mathsf{ub}}^{2}}{p}\left(\left\Vert \bm{U}_{l,\cdot}^{\star}\bm{\Sigma}^{\star}\right\Vert _{2}^{2}\log\left(n+d\right)+\omega_{l}^{\star2}\right).
\]
\end{itemize}
Equipped with the above quantities, we are ready to invoke the truncated
matrix Bernstein inequality \citep[Theorem 3.1.1]{chen2020spectral}
to show that 
\begin{align*}
\bigg\|\sum_{j=1}^{n}\bm{X}_{j}\bigg\|_{2} & \lesssim\sqrt{v\log\left(n+d\right)}+L\log\left(n+d\right)+nq_{1}\\
 & \lesssim\sigma_{\mathsf{ub}}\left(\left\Vert \bm{U}_{l,\cdot}^{\star}\bm{\Sigma}^{\star}\right\Vert _{2}\sqrt{\log\left(n+d\right)}+\omega_{l}^{\star}\right)\sqrt{\frac{r\log\left(n+d\right)}{p}}+\sqrt{\frac{r\log^{4}\left(n+d\right)}{np^{2}}}\sigma_{\mathsf{ub}}\left(\left\Vert \bm{U}_{l,\cdot}^{\star}\bm{\Sigma}^{\star}\right\Vert _{2}+\omega_{l}^{\star}\right)\left[1+\sqrt{\frac{\mu}{dp}}\right]\\
 & \lesssim\sigma_{\mathsf{ub}}\left(\left\Vert \bm{U}_{l,\cdot}^{\star}\bm{\Sigma}^{\star}\right\Vert _{2}+\omega_{l}^{\star}\right)\sqrt{\frac{r\log^{2}\left(n+d\right)}{p}}
\end{align*}
with probability exceeding $1-O((n+d)^{-10})$, where the last line
holds as long as $np\gtrsim\log^{2}(n+d)$ and $ndp^{2}\gtrsim\mu\log^{2}(n+d)$.
Substitution into \eqref{eq:alpha2-UB-125} yields 
\[
\alpha_{2}\ind_{\mathcal{E}_{\mathsf{good}}}\lesssim\frac{1}{\sigma_{r}^{\star2}}\sigma_{\mathsf{ub}}\left(\left\Vert \bm{U}_{l,\cdot}^{\star}\bm{\Sigma}^{\star}\right\Vert _{2}+\omega_{l}^{\star}\right)\sqrt{\frac{r\log^{2}\left(n+d\right)}{p}}.
\]

Putting the above bounds together allows one to conclude that 
\begin{align}
\left\Vert \bm{U}_{l,\cdot}\bm{R}-\bm{U}_{l,\cdot}^{\star}\right\Vert _{2} & \leq\alpha_{1}+\alpha_{2}+\zeta_{\mathsf{2nd},l}\nonumber \\
 & \lesssim\frac{1}{\sigma_{r}^{\star}}\left(\left\Vert \bm{U}_{l,\cdot}^{\star}\bm{\Sigma}^{\star}\right\Vert _{2}+\omega_{l}^{\star}\right)\sqrt{\frac{r\log^{2}\left(n+d\right)}{np}}+\frac{\sigma_{\mathsf{ub}}}{\sigma_{r}^{\star2}}\left(\left\Vert \bm{U}_{l,\cdot}^{\star}\bm{\Sigma}^{\star}\right\Vert _{2}+\omega_{l}^{\star}\right)\sqrt{\frac{r\log^{2}\left(n+d\right)}{p}}+\zeta_{\mathsf{2nd},l}\nonumber \\
 & \lesssim\frac{1}{\sigma_{r}^{\star}}\sqrt{\frac{r\log^{2}\left(n+d\right)}{np}}\left(\left\Vert \bm{U}_{l,\cdot}^{\star}\bm{\Sigma}^{\star}\right\Vert _{2}+\omega_{l}^{\star}\right)\left(1+\frac{\sigma_{\mathsf{ub}}}{\sigma_{r}^{\star}}\sqrt{n}\right)+\zeta_{\mathsf{2nd},l}\nonumber \\
 & \asymp\frac{\theta}{\sqrt{\kappa}\sigma_{r}^{\star}}\left(\left\Vert \bm{U}_{l,\cdot}^{\star}\bm{\Sigma}^{\star}\right\Vert _{2}+\omega_{l}^{\star}\right)+\zeta_{\mathsf{2nd},l},\label{eq:pca-1st-inter-0}
\end{align}
where we define 
\[
\theta\coloneqq\sqrt{\frac{\kappa r\log^{2}\left(n+d\right)}{np}}\left(1+\frac{\sigma_{\mathsf{ub}}}{\sigma_{r}^{\star}}\sqrt{n}\right)
\]
for notational simplicity. By taking the supremum over $l\in[d]$,
we further arrive at 
\begin{align}
\left\Vert \bm{U}\bm{R}-\bm{U}^{\star}\right\Vert _{2,\infty} & =\max_{l\in[d]}\left\Vert \bm{U}_{l,\cdot}\bm{R}-\bm{U}_{l,\cdot}^{\star}\right\Vert _{2}\nonumber \\
 & \lesssim\frac{1}{\sigma_{r}^{\star}}\sqrt{\frac{r\log^{2}\left(n+d\right)}{np}}\left(\left\Vert \bm{U}^{\star}\bm{\Sigma}^{\star}\right\Vert _{2,\infty}+\omega_{\max}\right)\left[1+\frac{\sigma_{\mathsf{ub}}}{\sigma_{r}^{\star}}\sqrt{n}\right]+\max_{l\in[d]}\zeta_{\mathsf{2nd},l}\nonumber \\
 & \lesssim\frac{\sigma_{\mathsf{ub}}}{\sigma_{r}^{\star}}\left(1+\frac{\sigma_{\mathsf{ub}}}{\sigma_{r}^{\star}}\sqrt{n}\right)\sqrt{r\log^{2}\left(n+d\right)}+\frac{\kappa^{3/2}\mu r\log^{1/2}\left(n+d\right)}{d}\frac{\zeta_{\mathsf{1st}}}{\sigma_{r}^{\star2}}+\frac{\zeta_{\mathsf{1st}}^{2}}{\sigma_{r}^{\star4}}\sqrt{\frac{\kappa^{3}\mu r\log\left(n+d\right)}{d}}\nonumber \\
 & \asymp\frac{\zeta_{\mathsf{1st}}}{\sigma_{r}^{\star2}}\sqrt{\frac{r\log\left(n+d\right)}{d}}+\frac{\kappa^{3/2}\mu r\log^{1/2}\left(n+d\right)}{d}\frac{\zeta_{\mathsf{1st}}}{\sigma_{r}^{\star2}}+\frac{\zeta_{\mathsf{1st}}^{2}}{\sigma_{r}^{\star4}}\sqrt{\frac{\kappa^{3}\mu r\log\left(n+d\right)}{d}}\nonumber \\
 & \asymp\frac{\zeta_{\mathsf{1st}}}{\sigma_{r}^{\star2}}\sqrt{\frac{r\log\left(n+d\right)}{d}},\label{eq:pca-1st-inter-6}
\end{align}
where the last relation holds provided that $d\gtrsim\kappa^{3}\mu^{2}r\log(n+d)$
and $\zeta_{\mathsf{1st}}/\sigma_{r}^{\star2}\lesssim1/\sqrt{\kappa^{3}\mu}$.

\paragraph{Bounding $\Vert\bm{R}^{\top}\bm{\Sigma}^{-2}\bm{R}-\left(\bm{\Sigma}^{\star}\right)^{-2}\Vert$. }

Conditional on $\bm{F}$, we learn from Lemma \ref{lemma:hpca-approx-2}
that with probability exceeding $1-O((n+d)^{-10})$ 
\begin{align}
\left\Vert \bm{Q}^{\top}\bm{R}^{\top}\bm{\Sigma}^{2}\bm{R}\bm{Q}-\bm{\Sigma}^{\natural2}\right\Vert  & \overset{\text{(i)}}{=}\left\Vert \bm{R}_{\bm{U}}^{\top}\bm{\Sigma}^{2}\bm{R}_{\bm{U}}-\bm{\Sigma}^{\natural2}\right\Vert \overset{\text{(ii)}}{\lesssim}\kappa^{\natural2}\sqrt{\frac{\mu^{\natural}r}{d}}\zeta_{\mathsf{1st}}+\kappa^{\natural2}\frac{\zeta_{\mathsf{1st}}^{2}}{\sigma_{r}^{\natural2}}\nonumber \\
 & \overset{\text{(iii)}}{\lesssim}\sqrt{\frac{\kappa^{3}\mu r\log\left(n+d\right)}{d}}\zeta_{\mathsf{1st}}+\kappa\frac{\zeta_{\mathsf{1st}}^{2}}{\sigma_{r}^{\star2}}.\label{eq:pca-1st-inter-1}
\end{align}
Here, (i) follows from the fact that $\bm{R}_{\bm{U}}=\bm{R}\bm{Q}$;
(ii) follows from Lemma \ref{lemma:hpca-approx-2}; and (iii) utilizes
(\ref{eq:good-event-sigma-least-largest-hpca}), (\ref{eq:good-event-mu-hpca})
and (\ref{eq:good-event-kappa-hpca}). We can then readily derive
\begin{align}
\sigma_{1}^{2} & \overset{\text{(i)}}{\leq}\sigma_{1}^{\natural2}+\left\Vert \bm{Q}^{\top}\bm{R}^{\top}\bm{\Sigma}^{2}\bm{R}\bm{Q}-\bm{\Sigma}^{\natural2}\right\Vert \overset{\text{(ii)}}{\leq}\sigma_{1}^{\natural2}+\sqrt{\frac{\kappa^{3}\mu r\log\left(n+d\right)}{d}}\zeta_{\mathsf{1st}}+\kappa\frac{\zeta_{\mathsf{1st}}^{2}}{\sigma_{r}^{\star2}}\overset{\text{(iii)}}{\leq}2\sigma_{1}^{\natural2}\overset{\text{(iv)}}{\leq}4\sigma_{1}^{\star2},\label{eq:pca-sigma-lower}\\
\sigma_{r}^{2} & \overset{\text{(v)}}{\geq}\sigma_{r}^{\natural2}-\left\Vert \bm{Q}^{\top}\bm{R}^{\top}\bm{\Sigma}^{2}\bm{R}\bm{Q}-\bm{\Sigma}^{\natural2}\right\Vert \overset{\text{(vi)}}{\geq}\sigma_{r}^{\natural2}-\sqrt{\frac{\kappa^{3}\mu r\log\left(n+d\right)}{d}}\zeta_{\mathsf{1st}}-\kappa\frac{\zeta_{\mathsf{1st}}^{2}}{\sigma_{r}^{\star2}}\overset{\text{(vii)}}{\geq}\frac{1}{2}\sigma_{r}^{\natural2}\overset{\text{(viii)}}{\geq}\frac{1}{4}\sigma_{r}^{\star2}.\label{eq:pca-sigma-upper}
\end{align}
Here, (i) and (v) follow from Weyl's inequality; (ii) and (vi) rely
on (\ref{eq:pca-1st-inter-1}); (iii) and (vii) make use of (\ref{eq:good-event-sigma-least-largest-hpca})
and hold true provided that $\zeta_{\mathsf{1st}}/\sigma_{r}^{\star2}\ll1/\sqrt{\kappa}$
and $d\gtrsim\kappa^{2}\mu r\log(n+d)$; (iv) and (viii) utilize (\ref{eq:good-event-Sigma-2-diff})
and are valid as long as $n\gg r+\log(n+d)$. On the event $\mathcal{E}_{\mathsf{good}}$,
it is seen from (\ref{eq:good-event-J-Q-hpca}) that 
\[
\left\Vert \bm{Q}\bm{\Sigma}^{\natural}\bm{Q}^{\top}-\bm{\Sigma}^{\star}\right\Vert =\left\Vert \bm{Q}\bm{\Sigma}^{\natural}-\bm{\Sigma}^{\star}\bm{Q}\right\Vert \lesssim\sqrt{\frac{\kappa\left(r+\log\left(n+d\right)\right)}{n}}\sigma_{1}^{\star},
\]
and as a result, 
\begin{align}
\left\Vert \bm{Q}\bm{\Sigma}^{\natural2}\bm{Q}^{\top}-\bm{\Sigma}^{\star2}\right\Vert  & =\left\Vert \left(\bm{Q}\bm{\Sigma}^{\natural}\bm{Q}^{\top}-\bm{\Sigma}^{\star}\right)\bm{Q}\bm{\Sigma}^{\natural}\bm{Q}^{\top}+\bm{\Sigma}^{\star}\left(\bm{Q}\bm{\Sigma}^{\natural}\bm{Q}^{\top}-\bm{\Sigma}^{\star}\right)\right\Vert \nonumber \\
 & \leq\left\Vert \bm{Q}\bm{\Sigma}^{\natural}\bm{Q}^{\top}-\bm{\Sigma}^{\star}\right\Vert \left(\sigma_{1}^{\natural}+\sigma_{1}^{\star}\right)\nonumber \\
 & \lesssim\sqrt{\frac{\kappa\left(r+\log\left(n+d\right)\right)}{n}}\sigma_{1}^{\star2},\label{eq:pca-1st-inter-2}
\end{align}
where the last line relies on (\ref{eq:good-event-sigma-least-largest-hpca}).
Taking (\ref{eq:pca-1st-inter-1}) and (\ref{eq:pca-1st-inter-2})
together yields 
\begin{align*}
\left\Vert \bm{R}^{\top}\bm{\Sigma}^{2}\bm{R}-\bm{\Sigma}^{\star2}\right\Vert  & \leq\left\Vert \bm{R}^{\top}\bm{\Sigma}^{2}\bm{R}-\bm{Q}\bm{\Sigma}^{\natural2}\bm{Q}^{\top}\right\Vert +\left\Vert \bm{Q}\bm{\Sigma}^{\natural2}\bm{Q}^{\top}-\bm{\Sigma}^{\star2}\right\Vert \\
 & =\left\Vert \bm{Q}^{\top}\bm{R}^{\top}\bm{\Sigma}^{2}\bm{R}\bm{Q}-\bm{\Sigma}^{\natural2}\right\Vert +\left\Vert \bm{Q}\bm{\Sigma}^{\natural2}\bm{Q}^{\top}-\bm{\Sigma}^{\star2}\right\Vert \\
 & \lesssim\sqrt{\frac{\kappa^{3}\mu r\log\left(n+d\right)}{d}}\zeta_{\mathsf{1st}}+\kappa\frac{\zeta_{\mathsf{1st}}^{2}}{\sigma_{r}^{\star2}}+\sqrt{\frac{\kappa\left(r+\log\left(n+d\right)\right)}{n}}\sigma_{1}^{\star2}.
\end{align*}
In view of the perturbation bound for matrix square roots \citep[Lemma 2.2]{MR1176461},
we obtain 
\begin{align}
\left\Vert \bm{R}^{\top}\bm{\Sigma}\bm{R}-\bm{\Sigma}^{\star}\right\Vert  & \leq\frac{1}{\sigma_{r}^{\star}+\sigma_{r}}\left\Vert \bm{R}^{\top}\bm{\Sigma}^{2}\bm{R}-\bm{\Sigma}^{\star2}\right\Vert \leq\frac{1}{\sigma_{r}^{\star}}\left\Vert \bm{R}^{\top}\bm{\Sigma}^{2}\bm{R}-\bm{\Sigma}^{\star2}\right\Vert \nonumber \\
 & \lesssim\sqrt{\frac{\kappa^{3}\mu r\log\left(n+d\right)}{d}}\frac{\zeta_{\mathsf{1st}}}{\sigma_{r}^{\star}}+\kappa\frac{\zeta_{\mathsf{1st}}^{2}}{\sigma_{r}^{\star3}}+\sqrt{\frac{\kappa^{2}\left(r+\log\left(n+d\right)\right)}{n}}\sigma_{1}^{\star}.\label{eq:pca-1st-inter-3}
\end{align}
In addition, we also have
\begin{align}
\left\Vert \bm{R}^{\top}\bm{\Sigma}^{-2}\bm{R}-\left(\bm{\Sigma}^{\star}\right)^{-2}\right\Vert  & =\left\Vert \bm{R}^{\top}\bm{\Sigma}^{-2}\bm{R}\left(\bm{R}^{\top}\bm{\Sigma}^{2}\bm{R}-\bm{\Sigma}^{\star2}\right)\left(\bm{\Sigma}^{\star}\right)^{-2}\right\Vert \nonumber \\
 & \leq\left\Vert \bm{R}^{\top}\bm{\Sigma}^{-2}\bm{R}\right\Vert \left\Vert \bm{R}^{\top}\bm{\Sigma}^{2}\bm{R}-\bm{\Sigma}^{\star2}\right\Vert \left\Vert \left(\bm{\Sigma}^{\star}\right)^{-2}\right\Vert \nonumber \\
 & \lesssim\frac{1}{\sigma_{r}^{\star4}}\left\Vert \bm{R}^{\top}\bm{\Sigma}^{2}\bm{R}-\bm{\Sigma}^{\star2}\right\Vert \nonumber \\
 & \lesssim\sqrt{\frac{\kappa^{3}\mu r\log\left(n+d\right)}{d}}\frac{\zeta_{\mathsf{1st}}}{\sigma_{r}^{\star4}}+\kappa\frac{\zeta_{\mathsf{1st}}^{2}}{\sigma_{r}^{\star6}}+\sqrt{\frac{\kappa^{3}\left(r+\log\left(n+d\right)\right)}{n}}\frac{1}{\sigma_{r}^{\star2}}.\label{eq:pca-1st-inter-5}
\end{align}
Here, the penultimate line comes from (\ref{eq:pca-sigma-lower}).

\paragraph{Bounding $\Vert(\bm{U}\bm{\Sigma}\bm{R}-\bm{U}^{\star}\bm{\Sigma}^{\star})_{l,\cdot}\Vert_{2}$.}

We start by bounding 
\begin{align}
\left\Vert \bm{U}_{l\cdot}\bm{\Sigma}\bm{R}-\bm{U}_{l,\cdot}^{\star}\bm{\Sigma}^{\star}\right\Vert _{2} & \leq\left\Vert \left(\bm{U}_{l\cdot}\bm{R}-\bm{U}_{l,\cdot}^{\star}\right)\bm{R}^{\top}\bm{\Sigma}\bm{R}+\bm{U}_{l,\cdot}^{\star}\left(\bm{R}^{\top}\bm{\Sigma}\bm{R}-\bm{\Sigma}^{\star}\right)\right\Vert _{2}\nonumber \\
 & \overset{\text{(i)}}{\leq}2\left\Vert \bm{U}_{l\cdot}\bm{R}-\bm{U}_{l,\cdot}^{\star}\right\Vert _{2}\sigma_{1}^{\star}+\left\Vert \bm{U}_{l,\cdot}^{\star}\right\Vert _{2}\left\Vert \bm{R}^{\top}\bm{\Sigma}\bm{R}-\bm{\Sigma}^{\star}\right\Vert \nonumber \\
 & \overset{\text{(ii)}}{\lesssim}\theta\left(\left\Vert \bm{U}_{l,\cdot}^{\star}\bm{\Sigma}^{\star}\right\Vert _{2}+\omega_{l}^{\star}\right)+\zeta_{\mathsf{2nd},l}\sigma_{1}^{\star}\nonumber \\
 & \quad+\left\Vert \bm{U}_{l,\cdot}^{\star}\right\Vert _{2}\left(\sqrt{\frac{\kappa^{3}\mu r\log\left(n+d\right)}{d}}\frac{\zeta_{\mathsf{1st}}}{\sigma_{r}^{\star}}+\kappa\frac{\zeta_{\mathsf{1st}}^{2}}{\sigma_{r}^{\star3}}+\sqrt{\frac{\kappa^{2}\left(r+\log\left(n+d\right)\right)}{n}}\sigma_{1}^{\star}\right)\nonumber \\
 & \lesssim\theta\left(\left\Vert \bm{U}_{l,\cdot}^{\star}\bm{\Sigma}^{\star}\right\Vert _{2}+\omega_{l}^{\star}\right)+\zeta_{\mathsf{2nd},l}\sigma_{1}^{\star}+\zeta_{\mathsf{2nd},l}\sigma_{r}^{\star}+\left\Vert \bm{U}_{l,\cdot}^{\star}\right\Vert _{2}\sqrt{\frac{\kappa^{2}\left(r+\log\left(n+d\right)\right)}{n}}\sigma_{1}^{\star}\nonumber \\
 & \asymp\theta\left(\left\Vert \bm{U}_{l,\cdot}^{\star}\bm{\Sigma}^{\star}\right\Vert _{2}+\omega_{l}^{\star}\right)+\left\Vert \bm{U}_{l,\cdot}^{\star}\right\Vert _{2}\sqrt{\frac{\kappa^{2}\left(r+\log\left(n+d\right)\right)}{n}}\sigma_{1}^{\star}+\zeta_{\mathsf{2nd},l}\sigma_{1}^{\star}.\label{eq:pca-1st-inter-7}
\end{align}
Here, (i) arises from (\ref{eq:pca-sigma-upper}), whereas (ii) follows
from (\ref{eq:pca-1st-inter-0}) and (\ref{eq:pca-1st-inter-3}).
In addition, we observe that 
\begin{align}
\left\Vert \bm{U}\bm{\Sigma}\bm{R}-\bm{U}^{\star}\bm{\Sigma}^{\star}\right\Vert _{2,\infty} & \leq\left\Vert \left(\bm{U}\bm{R}-\bm{U}^{\star}\right)\bm{R}^{\top}\bm{\Sigma}\bm{R}+\bm{U}^{\star}\left(\bm{R}^{\top}\bm{\Sigma}\bm{R}-\bm{\Sigma}^{\star}\right)\right\Vert _{2,\infty}\nonumber \\
 & \overset{\text{(i)}}{\leq}\left\Vert \bm{U}\bm{R}-\bm{U}^{\star}\right\Vert _{2,\infty}\sigma_{1}^{\star}+\left\Vert \bm{U}^{\star}\right\Vert _{2,\infty}\left\Vert \bm{R}^{\top}\bm{\Sigma}\bm{R}-\bm{\Sigma}^{\star}\right\Vert \nonumber \\
 & \overset{\text{(ii)}}{\lesssim}\frac{\zeta_{\mathsf{1st}}}{\sigma_{r}^{\star}}\sqrt{\frac{\kappa r\log\left(n+d\right)}{d}}+\sqrt{\frac{\mu r}{d}}\left(\sqrt{\frac{\kappa^{3}\mu r\log\left(n+d\right)}{d}}\frac{\zeta_{\mathsf{1st}}}{\sigma_{r}^{\star}}+\kappa\frac{\zeta_{\mathsf{1st}}^{2}}{\sigma_{r}^{\star3}}+\sqrt{\frac{\kappa^{2}\left(r+\log\left(n+d\right)\right)}{n}}\sigma_{1}^{\star}\right)\nonumber \\
 & \overset{\text{(iii)}}{\lesssim}\frac{\zeta_{\mathsf{1st}}}{\sigma_{r}^{\star}}\sqrt{\frac{\kappa r\log\left(n+d\right)}{d}}+\sqrt{\frac{\kappa^{2}\mu r^{2}\log\left(n+d\right)}{nd}}\sigma_{1}^{\star}.\label{eq:pca-1st-inter-8}
\end{align}
Here, (i) results from (\ref{eq:pca-sigma-upper}); (ii) follows from
(\ref{eq:pca-1st-inter-6}) and (\ref{eq:pca-1st-inter-3}); (iii)
holds true provided that $d\gtrsim\kappa^{2}\mu^{2}r$ and $\zeta_{\mathsf{1st}}/\sigma_{r}^{\star2}\lesssim1/\sqrt{\kappa\mu}$.

\paragraph{Bounding $\Vert\bm{U}\bm{\Sigma}^{-2}\bm{R}-\bm{U}^{\star}(\bm{\Sigma}^{\star})^{-2}\Vert$.}

It is first observed that, on the event $\mathcal{E}_{\mathsf{good}}$,
\begin{equation}
\left\Vert \bm{U}\bm{R}-\bm{U}^{\star}\right\Vert \overset{\text{(i)}}{=}\left\Vert \bm{U}\bm{R}_{\bm{U}}-\bm{U}^{\natural}\right\Vert \overset{\text{(ii)}}{\lesssim}\frac{\zeta_{\mathsf{1st}}}{\sigma_{r}^{\natural2}}\overset{\text{(iii)}}{\lesssim}\frac{\zeta_{\mathsf{1st}}}{\sigma_{r}^{\star2}}.\label{eq:pca-1st-inter-4}
\end{equation}
Here, (i) comes from the facts that $\bm{R}_{\bm{U}}=\bm{R}\bm{Q}$
and $\bm{U}^{\natural}=\bm{U}^{\star}\bm{Q}$ for some orthonormal
matrix $\bm{Q}$ (see \eqref{eq:defn-Q-denoising-PCA}); (ii) follows
from Lemma \ref{lemma:hpca-basic-facts}; and (iii) is a consequence
of (\ref{eq:good-event-sigma-least-largest-hpca}). This immediately
gives 
\begin{align*}
\left\Vert \bm{U}\bm{\Sigma}^{-2}\bm{R}-\bm{U}^{\star}\left(\bm{\Sigma}^{\star}\right)^{-2}\right\Vert  & =\left\Vert \left(\bm{U}\bm{R}-\bm{U}^{\star}\right)\bm{R}^{\top}\bm{\Sigma}^{-2}\bm{R}+\bm{U}^{\star}\left[\bm{R}^{\top}\bm{\Sigma}^{-2}\bm{R}-\left(\bm{\Sigma}^{\star}\right)^{-2}\right]\right\Vert \\
 & \leq\frac{1}{\sigma_{r}^{2}}\left\Vert \bm{U}\bm{R}-\bm{U}^{\star}\right\Vert +\left\Vert \bm{U}^{\star}\right\Vert \left\Vert \bm{R}^{\top}\bm{\Sigma}^{-2}\bm{R}-\left(\bm{\Sigma}^{\star}\right)^{-2}\right\Vert \\
 & \overset{\text{(i)}}{\lesssim}\frac{\zeta_{\mathsf{1st}}}{\sigma_{r}^{\star4}}+\sqrt{\frac{\kappa^{3}\mu r\log\left(n+d\right)}{d}}\frac{\zeta_{\mathsf{1st}}}{\sigma_{r}^{\star4}}+\kappa\frac{\zeta_{\mathsf{1st}}^{2}}{\sigma_{r}^{\star6}}+\sqrt{\frac{\kappa^{3}\left(r+\log\left(n+d\right)\right)}{n}}\frac{1}{\sigma_{r}^{\star2}}\\
 & \overset{\text{(ii)}}{\asymp}\frac{\zeta_{\mathsf{1st}}}{\sigma_{r}^{\star4}}+\sqrt{\frac{\kappa^{3}\left(r+\log\left(n+d\right)\right)}{n}}\frac{1}{\sigma_{r}^{\star2}},
\end{align*}
where (i) follows from (\ref{eq:pca-1st-inter-5}) and (\ref{eq:pca-1st-inter-4}),
and (ii) holds true as long as $d\gtrsim\kappa^{3}\mu r\log(n+d)$
and $\zeta_{\mathsf{1st}}/\sigma_{r}^{\star2}\lesssim1/\sqrt{\kappa}$.

\subsubsection{Proof of Lemma \ref{lemma:pca-noise-level-est}\label{appendix:proof-pca-noise-level-est}}

Before proceeding, let us make some useful observations: for all $l\in[d]$,
\begin{align}
\left\Vert \bm{U}_{l,\cdot}\bm{\Sigma}\right\Vert _{2} & =\left\Vert \bm{U}_{l,\cdot}\bm{\Sigma}\bm{R}\right\Vert _{2}\leq\left\Vert \bm{U}_{l\cdot}\bm{\Sigma}\bm{R}-\bm{U}_{l,\cdot}^{\star}\bm{\Sigma}^{\star}\right\Vert _{2}+\left\Vert \bm{U}_{l,\cdot}^{\star}\bm{\Sigma}^{\star}\right\Vert _{2}\nonumber \\
 & \overset{\text{(i)}}{\lesssim}\theta\left(\left\Vert \bm{U}_{l,\cdot}^{\star}\bm{\Sigma}^{\star}\right\Vert _{2}+\omega_{l}^{\star}\right)+\left\Vert \bm{U}_{l,\cdot}^{\star}\right\Vert _{2}\sqrt{\frac{\kappa^{2}\left(r+\log\left(n+d\right)\right)}{n}}\sigma_{1}^{\star}+\zeta_{\mathsf{2nd},l}\sigma_{1}^{\star}+\left\Vert \bm{U}_{l,\cdot}^{\star}\bm{\Sigma}^{\star}\right\Vert _{2}\nonumber \\
 & \lesssim\theta\left(\left\Vert \bm{U}_{l,\cdot}^{\star}\bm{\Sigma}^{\star}\right\Vert _{2}+\omega_{l}^{\star}\right)+\frac{1}{\sigma_{r}^{\star}}\left\Vert \bm{U}_{l,\cdot}^{\star}\bm{\Sigma}^{\star}\right\Vert _{2}\sqrt{\frac{\kappa^{2}\left(r+\log\left(n+d\right)\right)}{n}}\sigma_{1}^{\star}+\zeta_{\mathsf{2nd},l}\sigma_{1}^{\star}+\left\Vert \bm{U}_{l,\cdot}^{\star}\bm{\Sigma}^{\star}\right\Vert _{2}\nonumber \\
 & \overset{\text{(ii)}}{\asymp}\left\Vert \bm{U}_{l,\cdot}^{\star}\bm{\Sigma}^{\star}\right\Vert _{2}+\theta\omega_{l}^{\star}+\zeta_{\mathsf{2nd},l}\sigma_{1}^{\star},\label{eq:pca-U-l-Sigma-bound}
\end{align}
where $\theta$ is defined in \eqref{eq:theta-definition}. Here,
(i) relies on (\ref{eq:pca-U-l-Sigma-R-error}), while (ii) holds
true as long as $n\gtrsim\kappa^{3}r\log(n+d)$ and $\theta\ll1$.
Additionally, note that 
\begin{align}
\theta & \asymp\sqrt{\frac{\kappa r\log^{2}\left(n+d\right)}{np}}\left(1+\sqrt{\frac{\kappa\mu r\log\left(n+d\right)}{dp}}+\frac{\omega_{\max}}{\sigma_{r}^{\star}\sqrt{p}}\right)\nonumber \\
 & \asymp\sqrt{\frac{\kappa r\log^{2}\left(n+d\right)}{np}}+\sqrt{\frac{\kappa^{2}\mu r^{2}\log^{3}\left(n+d\right)}{ndp^{2}}}+\frac{\omega_{\max}}{\sigma_{r}^{\star}}\sqrt{\frac{\kappa r\log^{2}\left(n+d\right)}{np^{2}}}.\label{eq:theta-exact}
\end{align}
In view of the following relation (which makes use of the AM-GM inequality)
\begin{equation}
\frac{\zeta_{\mathsf{1st}}}{\sigma_{r}^{\star2}}\geq\frac{\kappa\mu r\log^{2}\left(n+d\right)}{\sqrt{nd}p}+\frac{\omega_{\max}^{2}}{p\sigma_{r}^{\star2}}\sqrt{\frac{d}{n}}\log\left(n+d\right)\geq\frac{\omega_{\max}}{\sigma_{r}^{\star}}\sqrt{\frac{\kappa\mu r\log^{3}\left(n+d\right)}{np^{2}}},\label{eq:pca-1st-err-useful}
\end{equation}
we can see that 
\begin{equation}
\theta\lesssim\frac{1}{\sqrt{\mu}}\cdot\frac{\zeta_{\mathsf{1st}}}{\sigma_{r}^{\star2}}.\label{eq:theta-zeta-1st}
\end{equation}
This means that $\theta\ll1$ can be guaranteed by the condition $\zeta_{\mathsf{1st}}\ll\sigma_{r}^{\star2}$.
We are now positioned to embark on the proof.

\paragraph{Step 1: bounding $\vert S_{i,j}-S_{i,j}^{\star}\vert$.}

We first develop an entrywise upper bound on $\bm{S}-\bm{S}^{\star}$
as follows 
\begin{align}
\left\Vert \bm{S}-\bm{S}^{\star}\right\Vert _{\infty} & =\left\Vert \left(\bm{U}\bm{\Sigma}\bm{R}\right)\left(\bm{U}\bm{\Sigma}\bm{R}\right)^{\top}-\bm{U}^{\star}\bm{\Sigma}^{\star}\left(\bm{U}^{\star}\bm{\Sigma}^{\star}\right)^{\top}\right\Vert _{\infty}\nonumber \\
 & =\left\Vert \left(\bm{U}\bm{\Sigma}\bm{R}-\bm{U}^{\star}\bm{\Sigma}^{\star}\right)\left(\bm{U}\bm{\Sigma}\bm{R}\right)^{\top}+\bm{U}^{\star}\bm{\Sigma}^{\star}\left(\bm{U}\bm{\Sigma}\bm{R}-\bm{U}^{\star}\bm{\Sigma}^{\star}\right)^{\top}\right\Vert _{\infty}\nonumber \\
 & \leq\left\Vert \bm{U}\bm{\Sigma}\bm{R}-\bm{U}^{\star}\bm{\Sigma}^{\star}\right\Vert _{2,\infty}\left\Vert \bm{U}\bm{\Sigma}\right\Vert _{2,\infty}+\left\Vert \bm{U}\bm{\Sigma}\bm{R}-\bm{U}^{\star}\bm{\Sigma}^{\star}\right\Vert _{2,\infty}\left\Vert \bm{U}^{\star}\bm{\Sigma}^{\star}\right\Vert _{2,\infty}\nonumber \\
 & \lesssim\left\Vert \bm{U}\bm{\Sigma}\bm{R}-\bm{U}^{\star}\bm{\Sigma}^{\star}\right\Vert _{2,\infty}\left\Vert \bm{U}^{\star}\bm{\Sigma}^{\star}\right\Vert _{2,\infty}+\left\Vert \bm{U}\bm{\Sigma}\bm{R}-\bm{U}^{\star}\bm{\Sigma}^{\star}\right\Vert _{2,\infty}^{2}\nonumber \\
 & \lesssim\left(\frac{\zeta_{\mathsf{1st}}}{\sigma_{r}^{\star}}\sqrt{\frac{\kappa r\log\left(n+d\right)}{d}}+\sqrt{\frac{\kappa^{2}\mu r^{2}\log\left(n+d\right)}{nd}}\sigma_{1}^{\star}\right)\sqrt{\frac{\mu r}{d}}\sigma_{1}^{\star}\nonumber \\
 & \quad+\left(\frac{\zeta_{\mathsf{1st}}}{\sigma_{r}^{\star}}\sqrt{\frac{\kappa r\log\left(n+d\right)}{d}}+\sqrt{\frac{\kappa^{2}\mu r^{2}\log\left(n+d\right)}{nd}}\sigma_{1}^{\star}\right)^{2}\nonumber \\
 & \asymp\left(\frac{\zeta_{\mathsf{1st}}}{\sigma_{r}^{\star}}\sqrt{\frac{\kappa r\log\left(n+d\right)}{d}}+\sqrt{\frac{\kappa^{2}\mu r^{2}\log\left(n+d\right)}{nd}}\sigma_{1}^{\star}\right)\sqrt{\frac{\mu r}{d}}\sigma_{1}^{\star}.\label{eq:pca-noise-inter-1}
\end{align}
Here the penultimate relation follows from (\ref{eq:pca-1st-inter-8}),
and the last relation holds provided that 
\[
\frac{\zeta_{\mathsf{1st}}}{\sigma_{r}^{\star}}\sqrt{\frac{\kappa r\log\left(n+d\right)}{d}}+\sqrt{\frac{\kappa^{2}\mu r^{2}\log\left(n+d\right)}{nd}}\sigma_{1}^{\star}\lesssim\sqrt{\frac{\mu r}{d}}\sigma_{1}^{\star},
\]
which can be guaranteed by $\zeta_{\mathsf{1st}}/\sigma_{r}^{\star2}\lesssim1/\sqrt{\log(n+d)}$
and $n\gtrsim\kappa^{2}r\log(n+d)$. Focusing on the $(i,j)$-th entry,
we can obtain (without loss of generality, assume $\omega_{i}^{\star}\leq\omega_{j}^{\star}$)
\begin{align}
\left|S_{i,j}-S_{i,j}^{\star}\right| & =\left|\left(\bm{U}_{i,\cdot}\bm{\Sigma}\bm{R}\right)\left(\bm{U}_{j,\cdot}\bm{\Sigma}\bm{R}\right)^{\top}-\bm{U}_{i,\cdot}^{\star}\bm{\Sigma}^{\star}\left(\bm{U}_{j,\cdot}^{\star}\bm{\Sigma}^{\star}\right)^{\top}\right|\nonumber \\
 & \leq\left\Vert \bm{U}_{i\cdot}\bm{\Sigma}\bm{R}-\bm{U}_{i,\cdot}^{\star}\bm{\Sigma}^{\star}\right\Vert _{2}\left\Vert \bm{U}_{j,\cdot}\bm{\Sigma}\right\Vert _{2}+\left\Vert \bm{U}_{j\cdot}\bm{\Sigma}\bm{R}-\bm{U}_{j,\cdot}^{\star}\bm{\Sigma}^{\star}\right\Vert _{2}\left\Vert \bm{U}_{i,\cdot}^{\star}\bm{\Sigma}^{\star}\right\Vert _{2}\nonumber \\
 & \overset{\text{(i)}}{\lesssim}\left\Vert \bm{U}_{i\cdot}\bm{\Sigma}\bm{R}-\bm{U}_{i,\cdot}^{\star}\bm{\Sigma}^{\star}\right\Vert _{2}\left\Vert \bm{U}_{j,\cdot}^{\star}\bm{\Sigma}^{\star}\right\Vert _{2}+\left\Vert \bm{U}_{i\cdot}\bm{\Sigma}\bm{R}-\bm{U}_{i,\cdot}^{\star}\bm{\Sigma}^{\star}\right\Vert _{2}\left(\theta\omega_{j}^{\star}+\zeta_{\mathsf{2nd},j}\sigma_{1}^{\star}\right)\nonumber \\
 & \quad+\left\Vert \bm{U}_{j\cdot}\bm{\Sigma}\bm{R}-\bm{U}_{j,\cdot}^{\star}\bm{\Sigma}^{\star}\right\Vert _{2}\left\Vert \bm{U}_{i,\cdot}^{\star}\bm{\Sigma}^{\star}\right\Vert _{2}\nonumber \\
 & \overset{\text{(ii)}}{\lesssim}\left[\theta\left(\left\Vert \bm{U}_{i,\cdot}^{\star}\bm{\Sigma}^{\star}\right\Vert _{2}+\omega_{i}^{\star}\right)+\left\Vert \bm{U}_{i,\cdot}^{\star}\right\Vert _{2}\sqrt{\frac{\kappa^{2}\left(r+\log\left(n+d\right)\right)}{n}}\sigma_{1}^{\star}+\zeta_{\mathsf{2nd},i}\sigma_{1}^{\star}\right]\left\Vert \bm{U}_{j,\cdot}^{\star}\bm{\Sigma}^{\star}\right\Vert _{2}\nonumber \\
 & \quad+\left[\theta\left(\left\Vert \bm{U}_{j,\cdot}^{\star}\bm{\Sigma}^{\star}\right\Vert _{2}+\omega_{j}^{\star}\right)+\left\Vert \bm{U}_{j,\cdot}^{\star}\right\Vert _{2}\sqrt{\frac{\kappa^{2}\left(r+\log\left(n+d\right)\right)}{n}}\sigma_{1}^{\star}+\zeta_{\mathsf{2nd},j}\sigma_{1}^{\star}\right]\left\Vert \bm{U}_{i,\cdot}^{\star}\bm{\Sigma}^{\star}\right\Vert _{2}\nonumber \\
 & \quad+\left[\theta\left(\left\Vert \bm{U}_{i,\cdot}^{\star}\bm{\Sigma}^{\star}\right\Vert _{2}+\omega_{i}^{\star}\right)+\left\Vert \bm{U}_{i,\cdot}^{\star}\right\Vert _{2}\sqrt{\frac{\kappa^{2}\left(r+\log\left(n+d\right)\right)}{n}}\sigma_{1}^{\star}+\zeta_{\mathsf{2nd},i}\sigma_{1}^{\star}\right]\left(\theta\omega_{j}^{\star}+\zeta_{\mathsf{2nd},j}\sigma_{1}^{\star}\right)\nonumber \\
 & \overset{\text{(iii)}}{\lesssim}\left(\theta+\sqrt{\frac{\kappa^{3}r\log\left(n+d\right)}{n}}\right)\left\Vert \bm{U}_{i,\cdot}^{\star}\bm{\Sigma}^{\star}\right\Vert _{2}\left\Vert \bm{U}_{j,\cdot}^{\star}\bm{\Sigma}^{\star}\right\Vert _{2}+\theta\left(\omega_{i}^{\star}\left\Vert \bm{U}_{j,\cdot}^{\star}\bm{\Sigma}^{\star}\right\Vert _{2}+\omega_{j}^{\star}\left\Vert \bm{U}_{i,\cdot}^{\star}\bm{\Sigma}^{\star}\right\Vert _{2}\right)\nonumber \\
 & \quad+\sigma_{1}^{\star}\left(\zeta_{\mathsf{2nd},i}\left\Vert \bm{U}_{j,\cdot}^{\star}\bm{\Sigma}^{\star}\right\Vert _{2}+\zeta_{\mathsf{2nd},j}\left\Vert \bm{U}_{i,\cdot}^{\star}\bm{\Sigma}^{\star}\right\Vert _{2}\right)+\theta^{2}\omega_{i}^{\star}\omega_{j}^{\star}+\zeta_{\mathsf{2nd},i}\zeta_{\mathsf{2nd},j}\sigma_{1}^{\star2}.\label{eq:pca-noise-inter-2}
\end{align}
Here, (i) follows from (\ref{eq:pca-U-l-Sigma-bound}), (ii) arises
from (\ref{eq:pca-U-l-Sigma-R-error}) and the fact that $\left\Vert \bm{U}_{j,\cdot}^{\star}\right\Vert _{2}\leq\frac{1}{\sigma_{r}^{\star}}\left\Vert \bm{U}_{j,\cdot}^{\star}\bm{\Sigma}^{\star}\right\Vert _{2}$,
while (iii) uses the AM-GM inequality $\zeta_{\mathsf{2nd}}\theta\sigma_{1}^{\star}\omega_{i}^{\star}\lesssim\theta^{2}\omega_{i}^{\star}\omega_{j}^{\star}+\zeta_{\mathsf{2nd}}^{2}\sigma_{1}^{\star2}$,
and holds true as long as $\theta\ll1$ (which is guaranteed when
$\zeta_{\mathsf{1st}}\ll\sigma_{r}^{\star2}$) and $n\gtrsim\kappa^{3}r\log(n+d)$.

\paragraph{Step 2: bounding $\vert\omega_{i}^{2}-\omega_{i}^{\star2}\vert$. }

With the above bound on $\big|S_{i,i}-S_{i,i}^{\star}\big|$ in place,
we can move on to control $\left|\omega_{i}^{2}-\omega_{i}^{\star2}\right|$.
First of all, invoke Chernoff's inequality (see \cite[Exercise 2.3.5]{vershynin2016high})
to show that 
\[
\mathbb{P}\left(\sum_{j=1}^{n}\ind_{(i,j)\in\Omega}<\frac{1}{2}np\right)\leq\exp\left(-\frac{c}{4}np\right)\leq\left(n+d\right)^{-10}
\]
as long as $np\gg\log(n+d)$. In what follows, we shall define the
following event 
\[
\mathcal{E}_{i}\coloneqq\left\{ \sum_{j=1}^{n}\ind_{(l,j)\in\Omega}\geq np/2\right\} ,
\]
which occurs with probability at least $1-O((n+d)^{-10})$.

Recall the definition of $y_{i,j}$ in \eqref{eq:observed-data}.
Conditional on 
\[
\Omega_{i,\cdot}=\{(i,j):(i,j)\in\Omega\},
\]
one can verify that $\{y_{i,j}:(i,j)\in\Omega_{l,\cdot}\}$ are independent
sub-Gaussian random variables with sub-Gaussian norm bounded above
by
\[
\left\Vert y_{i,j}\right\Vert _{\psi_{2}}=\left\Vert x_{i,j}\right\Vert _{\psi_{2}}+\left\Vert \eta_{i,j}\right\Vert _{\psi_{2}}\lesssim\sqrt{S_{i,i}^{\star}}+\omega_{i}^{\star}
\]
for each $j\in[n]$. As a result, we know that: conditional on $\Omega_{i,\cdot}$,
$\{y_{i,j}^{2}:(i,j)\in\Omega_{i,\cdot}\}$ are independent sub-exponential
random variables obeying
\[
K\coloneqq\max_{1\leq j\leq n}\left\Vert y_{i,j}^{2}\right\Vert _{\psi_{1}}\leq\max_{1\leq j\leq n}\left\Vert y_{i,j}\right\Vert _{\psi_{2}}^{2}\lesssim S_{i,i}^{\star}+\omega_{i}^{\star2}.
\]
In addition, it is easily observed that $\mathbb{E}[y_{i,j}^{2}]=S_{i,i}^{\star}+\omega_{i}^{\star2}$.
Then we can apply the Bernstein inequality (see \cite[Theorem 2.8.1]{vershynin2016high})
to demonstrate that: for any $t>0$, 
\[
\mathbb{P}\left(\left|\frac{\sum_{j=1}^{n}y_{i,j}^{2}\ind_{(i,j)\in\Omega}}{\sum_{j=1}^{n}\ind_{(i,j)\in\Omega}}-\omega_{i}^{\star2}-S_{i,i}^{\star}\right|\geq t\,\Biggl|\,\Omega_{i,\cdot}\right)\leq2\exp\left[-c\left(\sum_{j=1}^{n}\ind_{(i,j)\in\Omega}\right)\min\left\{ \frac{t^{2}}{K^{2}},\frac{t}{K}\right\} \right]
\]
for some universal constant $c>0$. Therefore, when $\mathcal{E}_{i}$
happens, we have
\[
\mathbb{P}\left(\left|\frac{\sum_{j=1}^{n}y_{i,j}^{2}\ind_{(i,j)\in\Omega}}{\sum_{j=1}^{n}\ind_{(i,j)\in\Omega}}-\omega_{i}^{\star2}-S_{i,i}^{\star}\right|\geq t\,\Biggl|\,\Omega_{i,\cdot}\right)\ind_{\mathcal{E}_{i}}\leq2\exp\left[-\frac{c}{2}np\min\left\{ \frac{t^{2}}{K^{2}},\frac{t}{K}\right\} \right].
\]
Take expectation to further achieve that: for any $t>0$, 
\begin{align*}
 & \mathbb{P}\left(\left|\frac{\sum_{j=1}^{n}y_{i,j}^{2}\ind_{(i,j)\in\Omega}}{\sum_{j=1}^{n}\ind_{(i,j)\in\Omega}}-\omega_{i}^{\star2}-S_{i,i}^{\star}\right|\geq t\right)=\mathbb{E}\left[\mathbb{P}\left(\left|\frac{\sum_{j=1}^{n}y_{i,j}^{2}\ind_{(i,j)\in\Omega}}{\sum_{j=1}^{n}\ind_{(i,j)\in\Omega}}-\omega_{i}^{\star2}-S_{i,i}^{\star}\right|\geq t\,\Biggl|\,\Omega_{i,\cdot}\right)\right]\\
 & \quad=\mathbb{E}\left[\mathbb{P}\left(\left|\frac{\sum_{j=1}^{n}y_{i,j}^{2}\ind_{(i,j)\in\Omega}}{\sum_{j=1}^{n}\ind_{(i,j)\in\Omega}}-\omega_{i}^{\star2}-S_{i,i}^{\star}\right|\geq t\,\Biggl|\,\Omega_{i,\cdot}\right)\ind_{\mathcal{E}_{i}}\right]\\
 & \quad\quad+\mathbb{E}\left[\mathbb{P}\left(\left|\frac{\sum_{j=1}^{n}y_{i,j}^{2}\ind_{(i,j)\in\Omega}}{\sum_{j=1}^{n}\ind_{(i,j)\in\Omega}}-\omega_{i}^{\star2}-S_{i,i}^{\star}\right|\geq t\,\Biggl|\,\Omega_{i,\cdot}\right)\ind_{\mathcal{E}_{i}^{\mathsf{c}}}\right]\\
 & \quad\leq2\exp\left[-\frac{c}{2}np\min\left\{ \frac{t^{2}}{K^{2}},\frac{t}{K}\right\} \right]+\mathbb{P}\left(\mathcal{E}_{i}^{\mathsf{c}}\right)\lesssim2\exp\left[-\frac{c}{2}np\min\left\{ \frac{t^{2}}{K^{2}},\frac{t}{K}\right\} \right]+\left(n+d\right)^{-10}.
\end{align*}
By taking 
\[
t=\widetilde{C}K\sqrt{\frac{\log^{2}\left(n+d\right)}{np}}\asymp\widetilde{C}\left(\omega_{i}^{\star2}+S_{i,i}^{\star}\right)\sqrt{\frac{\log^{2}\left(n+d\right)}{np}}
\]
 for some sufficiently large constant $\widetilde{C}>0$, we see that
with probability exceeding $1-O((n+d)^{-10})$ 
\[
\left|\frac{\sum_{j=1}^{n}y_{i,j}^{2}\ind_{(i,j)\in\Omega}}{\sum_{j=1}^{n}\ind_{(i,j)\in\Omega}}-\omega_{i}^{\star2}-S_{i,i}^{\star}\right|\lesssim\sqrt{\frac{\log^{2}\left(n+d\right)}{np}}\left(\omega_{i}^{\star2}+S_{i,i}^{\star}\right).
\]
A direct consequence is that
\begin{align*}
\left|\omega_{i}^{2}-\omega_{i}^{\star2}\right| & \leq\left|\frac{\sum_{j=1}^{n}y_{i,j}^{2}\ind_{(i,j)\in\Omega}}{\sum_{j=1}^{n}\ind_{(i,j)\in\Omega}}-\omega_{i}^{\star2}-S_{i,i}^{\star}\right|+\left|S_{i,i}-S_{i,i}^{\star}\right|\\
 & \lesssim\sqrt{\frac{\log^{2}\left(n+d\right)}{np}}\left(\omega_{i}^{\star2}+S_{i,i}^{\star}\right)+\left|S_{i,i}-S_{i,i}^{\star}\right|.
\end{align*}
In view of (\ref{eq:pca-noise-inter-1}) and the fact that $\|\bm{S}^{\star}\|_{\infty}=\|\bm{U}^{\star}\|_{2,\infty}^{2}\|\bm{\Sigma}^{\star2}\|\leq\frac{\mu r}{d}\sigma_{1}^{\star2}$,
the above inequality implies that, for all $i\in[d]$, 
\begin{align*}
\left|\omega_{i}^{2}-\omega_{i}^{\star2}\right| & \lesssim\sqrt{\frac{\log^{2}\left(n+d\right)}{np}}\left(\omega_{i}^{\star2}+\frac{\mu r}{d}\sigma_{1}^{\star2}\right)+\left(\frac{\zeta_{\mathsf{1st}}}{\sigma_{r}^{\star}}\sqrt{\frac{\kappa r\log\left(n+d\right)}{d}}+\sqrt{\frac{\kappa^{2}\mu r^{2}\log\left(n+d\right)}{nd}}\sigma_{1}^{\star}\right)\sqrt{\frac{\mu r}{d}}\sigma_{1}^{\star}\\
 & \asymp\sqrt{\frac{\log^{2}\left(n+d\right)}{np}}\omega_{i}^{\star2}+\sqrt{\frac{\log^{2}\left(n+d\right)}{np}}\frac{\mu r}{d}\sigma_{1}^{\star2}+\zeta_{\mathsf{1st}}\frac{\sqrt{\kappa^{2}\mu r^{2}\log\left(n+d\right)}}{d}+\sqrt{\frac{\kappa^{2}\mu^{2}r^{3}\log\left(n+d\right)}{nd^{2}}}\sigma_{1}^{\star2}\\
 & \asymp\sqrt{\frac{\log^{2}\left(n+d\right)}{np}}\omega_{i}^{\star2}+\zeta_{\mathsf{1st}}\frac{\sqrt{\kappa^{2}\mu r^{2}\log\left(n+d\right)}}{d}+\sqrt{\frac{\kappa^{2}\mu^{2}r^{3}\log\left(n+d\right)}{nd^{2}}}\sigma_{1}^{\star2},
\end{align*}
where the last line makes use of the definition of $\zeta_{\mathsf{1st}}$
in \eqref{eq:zeta-1st-UB-hpca}. Additionally, in view of (\ref{eq:pca-noise-inter-2}),
we can derive for $i=l$ that 
\begin{align*}
\left|\omega_{l}^{2}-\omega_{l}^{\star2}\right| & \lesssim\sqrt{\frac{\log^{2}\left(n+d\right)}{np}}\left(\omega_{l}^{\star2}+S_{l,l}^{\star}\right)+\left(\theta+\sqrt{\frac{\kappa^{3}r\log\left(n+d\right)}{n}}\right)\left\Vert \bm{U}_{l,\cdot}^{\star}\bm{\Sigma}^{\star}\right\Vert _{2}^{2}\\
 & \quad+\left(\theta\omega_{l}^{\star}+\zeta_{\mathsf{2nd},l}\sigma_{1}^{\star}\right)\left\Vert \bm{U}_{l,\cdot}^{\star}\bm{\Sigma}^{\star}\right\Vert _{2}+\theta^{2}\omega_{l}^{\star2}+\zeta_{\mathsf{2nd},l}^{2}\sigma_{1}^{\star2}\\
 & \asymp\left(\sqrt{\frac{\log^{2}\left(n+d\right)}{np}}+\theta^{2}\right)\omega_{l}^{\star2}+\left(\theta+\sqrt{\frac{\kappa^{3}r\log\left(n+d\right)}{n}}\right)\left\Vert \bm{U}_{l,\cdot}^{\star}\bm{\Sigma}^{\star}\right\Vert _{2}^{2}\\
 & \quad+\left(\theta\omega_{l}^{\star}+\zeta_{\mathsf{2nd},l}\sigma_{1}^{\star}\right)\left\Vert \bm{U}_{l,\cdot}^{\star}\bm{\Sigma}^{\star}\right\Vert _{2}+\zeta_{\mathsf{2nd},l}^{2}\sigma_{1}^{\star2}.
\end{align*}
Here, the last relation invokes the following condition 
\[
\sqrt{\frac{\log^{2}\left(n+d\right)}{np}}S_{l,l}^{\star}\lesssim\theta\left\Vert \bm{U}_{l,\cdot}^{\star}\bm{\Sigma}^{\star}\right\Vert _{2}^{2},
\]
which is a direct consequence of (\ref{eq:theta-exact}).

\subsubsection{Proof of Lemma \ref{lemma:pca-covariance-estimation}\label{appendix:proof-pca-covariance-estimation}}

Before proceeding, let us recall a few facts as follows. To begin
with, from the definition \eqref{eq:defn-Sigma-Ul-star-hpca} of $\bm{\Sigma}_{U,l}^{\star}$,
we can decompose 
\begin{align*}
\bm{R}\bm{\Sigma}_{U,l}^{\star}\bm{R}^{\top} & =\underbrace{\frac{1-p}{np}\left\Vert \bm{U}_{l,\cdot}^{\star}\bm{\Sigma}^{\star}\right\Vert _{2}^{2}\bm{R}\left(\bm{\Sigma}^{\star}\right)^{-2}\bm{R}^{\top}}_{\eqqcolon\bm{A}_{1}}+\underbrace{\frac{\omega_{l}^{\star2}}{np}\bm{R}\left(\bm{\Sigma}^{\star}\right)^{-2}\bm{R}^{\top}}_{\eqqcolon\bm{A}_{2}}+\underbrace{\frac{2\left(1-p\right)}{np}\bm{R}\bm{U}_{l,\cdot}^{\star\top}\bm{U}_{l,\cdot}^{\star}\bm{R}^{\top}}_{\eqqcolon\bm{A}_{3}}\\
 & \quad+\underbrace{\bm{R}\left(\bm{\Sigma}^{\star}\right)^{-2}\bm{U}^{\star\top}\mathsf{diag}\left\{ \left[d_{l,i}^{\star}\right]_{i=1}^{d}\right\} \bm{U}^{\star}\left(\bm{\Sigma}^{\star}\right)^{-2}\bm{R}^{\top}}_{\eqqcolon\bm{A}_{4}},
\end{align*}
where we recall that 
\[
d_{l,i}^{\star}\coloneqq\frac{1}{np^{2}}\left[\omega_{l}^{\star2}+\left(1-p\right)\left\Vert \bm{U}_{l,\cdot}^{\star}\bm{\Sigma}^{\star}\right\Vert _{2}^{2}\right]\left[\omega_{i}^{\star2}+\left(1-p\right)\left\Vert \bm{U}_{i,\cdot}^{\star}\bm{\Sigma}^{\star}\right\Vert _{2}^{2}\right]+\frac{2\left(1-p\right)^{2}}{np^{2}}S_{i,l}^{\star2}.
\]
Regarding $\bm{\Sigma}_{U,l}$, we remind the reader of its definition
in \eqref{eq:defn-plug-in-Sigma-Ul-hpca} as follows 
\begin{align*}
\bm{\Sigma}_{U,l} & =\underbrace{\frac{1-p}{np}\left\Vert \bm{U}_{l,\cdot}\bm{\Sigma}\right\Vert _{2}^{2}\bm{\Sigma}^{-2}}_{\eqqcolon\bm{B}_{1}}+\underbrace{\frac{\omega_{l}^{2}}{np}\bm{\Sigma}^{-2}}_{\eqqcolon\bm{B}_{2}}+\underbrace{\frac{2\left(1-p\right)}{np}\bm{U}_{l,\cdot}^{\top}\bm{U}_{l,\cdot}}_{\eqqcolon\bm{B}_{3}}+\underbrace{\bm{\Sigma}^{-2}\bm{U}^{\top}\mathsf{diag}\left\{ \left[d_{l,i}\right]_{1\leq i\leq d}\right\} \bm{U}\bm{\Sigma}^{-2}}_{\eqqcolon\bm{B}_{4}},
\end{align*}
where we define 
\[
d_{l,i}\coloneqq\frac{1}{np^{2}}\left[\omega_{l}^{2}+\left(1-p\right)\left\Vert \bm{U}_{l,\cdot}\bm{\Sigma}\right\Vert _{2}^{2}\right]\left[\omega_{i}^{2}+\left(1-p\right)\left\Vert \bm{U}_{i,\cdot}\bm{\Sigma}\right\Vert _{2}^{2}\right]+\frac{2\left(1-p\right)^{2}}{np^{2}}S_{l,i}^{2}.
\]
In addition, we also recall from Lemma \ref{lemma:pca-covariance-concentration}
that 
\begin{align*}
\lambda_{\min}\left(\bm{\Sigma}_{U,l}^{\star}\right) & \gtrsim\frac{1}{np\sigma_{1}^{\star2}}\left\Vert \bm{U}_{l,\cdot}^{\star}\bm{\Sigma}^{\star}\right\Vert _{2}^{2}+\frac{\omega_{l}^{\star2}}{np\sigma_{1}^{\star2}}+\frac{1}{ndp^{2}\kappa\sigma_{1}^{\star2}}\left\Vert \bm{U}_{l,\cdot}^{\star}\bm{\Sigma}^{\star}\right\Vert _{2}^{2}+\frac{1}{ndp^{2}\kappa\sigma_{1}^{\star2}}\omega_{l}^{\star2}\\
 & \quad+\frac{1}{np^{2}\sigma_{1}^{\star4}}\omega_{\min}^{2}\left\Vert \bm{U}_{l,\cdot}^{\star}\bm{\Sigma}^{\star}\right\Vert _{2}^{2}+\frac{\omega_{l}^{\star2}\omega_{\min}^{2}}{np^{2}\sigma_{1}^{\star4}}.
\end{align*}
We are now ready to present the proof, with the focus of bounding
$\|\bm{A}_{i}-\bm{B}_{i}\|$ for each $1\leq i\leq4$ as well as $\big|d_{l,i}^{\star}-d_{l,i}\big|$.

\paragraph{Step 1: controlling $\Vert\bm{A}_{1}-\bm{B}_{1}\Vert$.}

It follows from the triangle inequality and Lemma \ref{lemma:pca-1st-err}
that 
\begin{align*}
\left|\left\Vert \bm{U}_{l,\cdot}^{\star}\bm{\Sigma}^{\star}\right\Vert _{2}-\left\Vert \bm{U}_{l,\cdot}\bm{\Sigma}\right\Vert _{2}\right| & =\left|\left\Vert \bm{U}_{l,\cdot}^{\star}\bm{\Sigma}^{\star}\right\Vert _{2}-\left\Vert \bm{U}_{l,\cdot}\bm{\Sigma}\bm{R}\right\Vert _{2}\right|\leq\left\Vert \bm{U}_{l\cdot}\bm{\Sigma}\bm{R}-\bm{U}_{l,\cdot}^{\star}\bm{\Sigma}^{\star}\right\Vert _{2}\\
 & \lesssim\theta\left(\left\Vert \bm{U}_{l,\cdot}^{\star}\bm{\Sigma}^{\star}\right\Vert _{2}+\omega_{l}^{\star}\right)+\left\Vert \bm{U}_{l,\cdot}^{\star}\right\Vert _{2}\sqrt{\frac{\kappa^{2}\left(r+\log\left(n+d\right)\right)}{n}}\sigma_{1}^{\star}+\zeta_{\mathsf{2nd},l}\sigma_{1}^{\star}.
\end{align*}
Here, the last line follows from (\ref{eq:pca-U-l-Sigma-R-error}).
This immediately gives 
\begin{align}
 & \left|\left\Vert \bm{U}_{l,\cdot}^{\star}\bm{\Sigma}^{\star}\right\Vert _{2}^{2}-\left\Vert \bm{U}_{l,\cdot}\bm{\Sigma}\right\Vert _{2}^{2}\right|\lesssim\left|\left\Vert \bm{U}_{l,\cdot}^{\star}\bm{\Sigma}^{\star}\right\Vert _{2}-\left\Vert \bm{U}_{l,\cdot}\bm{\Sigma}\right\Vert _{2}\right|\left|\left\Vert \bm{U}_{l,\cdot}^{\star}\bm{\Sigma}^{\star}\right\Vert _{2}+\left\Vert \bm{U}_{l,\cdot}\bm{\Sigma}\right\Vert _{2}\right|\nonumber \\
 & \qquad\asymp\left|\left\Vert \bm{U}_{l,\cdot}^{\star}\bm{\Sigma}^{\star}\right\Vert _{2}-\left\Vert \bm{U}_{l,\cdot}\bm{\Sigma}\right\Vert _{2}\right|\left\Vert \bm{U}_{l,\cdot}^{\star}\bm{\Sigma}^{\star}\right\Vert _{2}+\left|\left\Vert \bm{U}_{l,\cdot}^{\star}\bm{\Sigma}^{\star}\right\Vert _{2}-\left\Vert \bm{U}_{l,\cdot}\bm{\Sigma}\right\Vert _{2}\right|^{2}\nonumber \\
 & \qquad\lesssim\theta\left(\left\Vert \bm{U}_{l,\cdot}^{\star}\bm{\Sigma}^{\star}\right\Vert _{2}+\omega_{l}^{\star}\right)\left\Vert \bm{U}_{l,\cdot}^{\star}\bm{\Sigma}^{\star}\right\Vert _{2}+\left\Vert \bm{U}_{l,\cdot}^{\star}\bm{\Sigma}^{\star}\right\Vert _{2}\left\Vert \bm{U}_{l,\cdot}^{\star}\right\Vert _{2}\sqrt{\frac{\kappa^{2}\left(r+\log\left(n+d\right)\right)}{n}}\sigma_{1}^{\star}\nonumber \\
 & \qquad\quad+\zeta_{\mathsf{2nd},l}\sigma_{1}^{\star}\left\Vert \bm{U}_{l,\cdot}^{\star}\bm{\Sigma}^{\star}\right\Vert _{2}+\theta^{2}\left(\left\Vert \bm{U}_{l,\cdot}^{\star}\bm{\Sigma}^{\star}\right\Vert _{2}^{2}+\omega_{l}^{\star2}\right)+\left\Vert \bm{U}_{l,\cdot}^{\star}\right\Vert _{2}^{2}\frac{\kappa^{2}\left(r+\log\left(n+d\right)\right)}{n}\sigma_{1}^{\star2}+\zeta_{\mathsf{2nd},l}^{2}\sigma_{1}^{\star2}\nonumber \\
 & \qquad\lesssim\left\Vert \bm{U}_{l,\cdot}^{\star}\bm{\Sigma}^{\star}\right\Vert _{2}^{2}\left(\theta+\sqrt{\frac{\kappa^{3}\left(r+\log\left(n+d\right)\right)}{n}}\right)+\theta\omega_{l}^{\star}\left\Vert \bm{U}_{l,\cdot}^{\star}\bm{\Sigma}^{\star}\right\Vert _{2}\nonumber \\
 & \quad\qquad+\zeta_{\mathsf{2nd},l}\sigma_{1}^{\star}\left\Vert \bm{U}_{l,\cdot}^{\star}\bm{\Sigma}^{\star}\right\Vert _{2}+\theta^{2}\omega_{l}^{\star2}+\zeta_{\mathsf{2nd},l}^{2}\sigma_{1}^{\star2}.\label{eq:pca-cov-est-inter-1}
\end{align}
Here the last relation holds provided that $\theta\ll1$ and $n\gtrsim\kappa^{3}r\log(n+d)$.
In addition, from (\ref{eq:pca-Sigma-2-error}) in Lemma \ref{lemma:pca-1st-err},
we know that 
\begin{equation}
\left\Vert \bm{R}\left(\bm{\Sigma}^{\star}\right)^{-2}\bm{R}^{\top}-\bm{\Sigma}^{-2}\right\Vert \lesssim\sqrt{\frac{\kappa^{3}\mu r\log\left(n+d\right)}{d}}\frac{\zeta_{\mathsf{1st}}}{\sigma_{r}^{\star4}}+\kappa\frac{\zeta_{\mathsf{1st}}^{2}}{\sigma_{r}^{\star6}}+\sqrt{\frac{\kappa^{3}\left(r+\log\left(n+d\right)\right)}{n}}\frac{1}{\sigma_{r}^{\star2}}.\label{eq:pca-cov-est-inter-2}
\end{equation}
These two inequalities help us derive 
\begin{align*}
\left\Vert \bm{A}_{1}-\bm{B}_{1}\right\Vert  & \leq\frac{1-p}{np}\left|\left\Vert \bm{U}_{l,\cdot}^{\star}\bm{\Sigma}^{\star}\right\Vert _{2}^{2}-\left\Vert \bm{U}_{l,\cdot}\bm{\Sigma}\right\Vert _{2}^{2}\right|\left\Vert \bm{R}\bm{\Sigma}^{-2}\bm{R}^{\top}\right\Vert +\frac{1-p}{np}\left\Vert \bm{U}_{l,\cdot}^{\star}\bm{\Sigma}^{\star}\right\Vert _{2}^{2}\left\Vert \bm{R}\left(\bm{\Sigma}^{\star}\right)^{-2}\bm{R}^{\top}-\bm{\Sigma}^{-2}\right\Vert \\
 & \overset{\text{(i)}}{\lesssim}\underbrace{\frac{1}{np\sigma_{r}^{\star2}}\left\Vert \bm{U}_{l,\cdot}^{\star}\bm{\Sigma}^{\star}\right\Vert _{2}^{2}\left(\theta+\sqrt{\frac{\kappa^{3}\left(r+\log\left(n+d\right)\right)}{n}}\right)}_{\eqqcolon\alpha_{1,1}}+\underbrace{\frac{1}{np\sigma_{r}^{\star2}}\theta\omega_{l}^{\star}\left\Vert \bm{U}_{l,\cdot}^{\star}\bm{\Sigma}^{\star}\right\Vert _{2}}_{\eqqcolon\alpha_{1,2}}\\
 & \quad+\frac{1}{np\sigma_{r}^{\star2}}\zeta_{\mathsf{2nd},l}\sigma_{1}^{\star}\left\Vert \bm{U}_{l,\cdot}^{\star}\bm{\Sigma}^{\star}\right\Vert _{2}+\underbrace{\frac{1}{np\sigma_{r}^{\star2}}\theta^{2}\omega_{l}^{\star2}}_{\eqqcolon\alpha_{1,3}}+\frac{1}{np\sigma_{r}^{\star2}}\zeta_{\mathsf{2nd},l}^{2}\sigma_{1}^{\star2}\\
 & \quad+\underbrace{\frac{1}{np}\left\Vert \bm{U}_{l,\cdot}^{\star}\bm{\Sigma}^{\star}\right\Vert _{2}^{2}\sqrt{\frac{\kappa^{3}\mu r\log\left(n+d\right)}{d}}\frac{\zeta_{\mathsf{1st}}}{\sigma_{r}^{\star4}}}_{\eqqcolon\alpha_{1,4}}+\underbrace{\frac{1}{np}\left\Vert \bm{U}_{l,\cdot}^{\star}\bm{\Sigma}^{\star}\right\Vert _{2}^{2}\kappa\frac{\zeta_{\mathsf{1st}}^{2}}{\sigma_{r}^{\star6}}}_{\eqqcolon\alpha_{1,5}}\\
 & \overset{\text{(ii)}}{\lesssim}\delta\lambda_{\min}\left(\bm{\Sigma}_{U,l}^{\star}\right)+\frac{\sqrt{\kappa}}{np\sigma_{r}^{\star}}\zeta_{\mathsf{2nd},l}\left\Vert \bm{U}_{l,\cdot}^{\star}\bm{\Sigma}^{\star}\right\Vert _{2}+\frac{\kappa}{np}\zeta_{\mathsf{2nd},l}^{2}.
\end{align*}
Here, (i) results from (\ref{eq:pca-sigma-lower}), (\ref{eq:pca-cov-est-inter-1})
and (\ref{eq:pca-cov-est-inter-2}); (ii) holds true due to the following
inequalities: 
\begin{align*}
\alpha_{1,1}+\alpha_{1,4}+\alpha_{1,5} & \lesssim\delta\frac{1}{np\sigma_{1}^{\star2}}\left\Vert \bm{U}_{l,\cdot}^{\star}\bm{\Sigma}^{\star}\right\Vert _{2}^{2}\lesssim\delta\lambda_{\min}\left(\bm{\Sigma}_{U,l}^{\star}\right),\\
\alpha_{1,2} & \lesssim\delta\frac{1}{np\sigma_{1}^{\star2}}\omega_{l}^{\star}\left\Vert \bm{U}_{l,\cdot}^{\star}\bm{\Sigma}^{\star}\right\Vert _{2}\lesssim\delta\frac{1}{np\sigma_{1}^{\star2}}\left\Vert \bm{U}_{l,\cdot}^{\star}\bm{\Sigma}^{\star}\right\Vert _{2}^{2}+\delta\frac{\omega_{l}^{\star2}}{np\sigma_{1}^{\star2}}\lesssim\delta\lambda_{\min}\left(\bm{\Sigma}_{U,l}^{\star}\right),\\
\alpha_{1,3} & \lesssim\delta\frac{\omega_{l}^{\star2}}{np\sigma_{1}^{\star2}}\lesssim\delta\lambda_{\min}\left(\bm{\Sigma}_{U,l}^{\star}\right),
\end{align*}
provided that $\theta\lesssim\delta/\kappa$, $n\gtrsim\delta^{-2}\kappa^{5}r\log(n+d)$,
$\zeta_{\mathsf{1st}}/\sigma_{r}^{\star2}\lesssim1/\kappa$ and $d\gtrsim\kappa^{3}\mu r\log(n+d)$.

\paragraph{Step 2: controlling $\Vert\bm{A}_{2}-\bm{B}_{2}\Vert$.}

Recall from (\ref{eq:pca-noise-est-l}) that 
\begin{align}
\left|\omega_{l}^{2}-\omega_{l}^{\star2}\right| & \lesssim\left(\sqrt{\frac{\log^{2}\left(n+d\right)}{np}}+\theta^{2}\right)\omega_{l}^{\star2}+\left(\theta+\sqrt{\frac{\kappa^{3}r\log\left(n+d\right)}{n}}\right)\left\Vert \bm{U}_{l,\cdot}^{\star}\bm{\Sigma}^{\star}\right\Vert _{2}^{2}\nonumber \\
 & \quad+\left(\theta\omega_{l}^{\star}+\zeta_{\mathsf{2nd},l}\sigma_{1}^{\star}\right)\left\Vert \bm{U}_{l,\cdot}^{\star}\bm{\Sigma}^{\star}\right\Vert _{2}+\zeta_{\mathsf{2nd},l}^{2}\sigma_{1}^{\star2}.\label{eq:pca-cov-est-inter-3}
\end{align}
Armed with this inequality, we can derive 
\begin{align*}
\left\Vert \bm{A}_{2}-\bm{B}_{2}\right\Vert  & =\left\Vert \frac{\omega_{l}^{\star2}}{np}\left[\bm{R}\left(\bm{\Sigma}^{\star}\right)^{-2}\bm{R}^{\top}-\bm{\Sigma}^{-2}\right]+\frac{\omega_{l}^{\star2}-\omega_{l}^{2}}{np}\bm{\Sigma}^{-2}\right\Vert \\
 & \overset{\text{(i)}}{\lesssim}\frac{\omega_{l}^{\star2}}{np}\left\Vert \bm{R}\left(\bm{\Sigma}^{\star}\right)^{-2}\bm{R}^{\top}-\bm{\Sigma}^{-2}\right\Vert +\frac{\left|\omega_{l}^{\star2}-\omega_{l}^{2}\right|}{np\sigma_{r}^{\star2}}\\
 & \overset{\text{(ii)}}{\lesssim}\underbrace{\frac{\omega_{l}^{\star2}}{np}\left(\sqrt{\frac{\kappa^{3}\mu r\log\left(n+d\right)}{d}}\frac{\zeta_{\mathsf{1st}}}{\sigma_{r}^{\star4}}+\kappa\frac{\zeta_{\mathsf{1st}}^{2}}{\sigma_{r}^{\star6}}\right)}_{\eqqcolon\alpha_{2,1}}+\underbrace{\frac{\omega_{l}^{\star2}}{np\sigma_{r}^{\star2}}\sqrt{\frac{\kappa^{3}\left(r+\log\left(n+d\right)\right)}{n}}}_{\eqqcolon\alpha_{2,2}}\\
 & \quad+\underbrace{\frac{1}{np\sigma_{r}^{\star2}}\left(\sqrt{\frac{\log^{2}\left(n+d\right)}{np}}+\theta^{2}\right)\omega_{l}^{\star2}}_{\eqqcolon\alpha_{2,3}}+\underbrace{\frac{1}{np\sigma_{r}^{\star2}}\left(\theta+\sqrt{\frac{\kappa^{3}r\log\left(n+d\right)}{n}}\right)\left\Vert \bm{U}_{l,\cdot}^{\star}\bm{\Sigma}^{\star}\right\Vert _{2}^{2}}_{\eqqcolon\alpha_{2,4}}\\
 & \quad+\underbrace{\frac{1}{np\sigma_{r}^{\star2}}\theta\omega_{l}^{\star}\left\Vert \bm{U}_{l,\cdot}^{\star}\bm{\Sigma}^{\star}\right\Vert _{2}}_{\eqqcolon\alpha_{2,5}}+\frac{1}{np\sigma_{r}^{\star2}}\zeta_{\mathsf{2nd},l}\sigma_{1}^{\star}\left\Vert \bm{U}_{l,\cdot}^{\star}\bm{\Sigma}^{\star}\right\Vert _{2}+\frac{1}{np\sigma_{r}^{\star2}}\zeta_{\mathsf{2nd},l}^{2}\sigma_{1}^{\star2}\\
 & \overset{\text{(iii)}}{\lesssim}\delta\lambda_{\min}\left(\bm{\Sigma}_{U,l}^{\star}\right)+\frac{\sqrt{\kappa}}{np\sigma_{r}^{\star}}\zeta_{\mathsf{2nd}}\left\Vert \bm{U}_{l,\cdot}^{\star}\bm{\Sigma}^{\star}\right\Vert _{2}+\frac{\kappa}{np}\zeta_{\mathsf{2nd}}^{2}.
\end{align*}
Here, (i) arises from (\ref{eq:pca-sigma-lower}); (ii) utilizes (\ref{eq:pca-cov-est-inter-2})
and (\ref{eq:pca-cov-est-inter-3}); (iii) follows from the inequalities
below: 
\begin{align*}
\alpha_{2,1}+\alpha_{2,2}+\alpha_{2,3} & \lesssim\delta\frac{\omega_{l}^{\star2}}{np\sigma_{1}^{\star2}}\lesssim\delta\lambda_{\min}\left(\bm{\Sigma}_{U,l}^{\star}\right),\\
\alpha_{2,4} & \lesssim\delta\frac{1}{np\sigma_{1}^{\star2}}\left\Vert \bm{U}_{l,\cdot}^{\star}\bm{\Sigma}^{\star}\right\Vert _{2}^{2}\lesssim\delta\lambda_{\min}\left(\bm{\Sigma}_{U,l}^{\star}\right),\\
\alpha_{2,5} & \lesssim\delta\frac{1}{np\sigma_{1}^{\star2}}\omega_{l}^{\star}\left\Vert \bm{U}_{l,\cdot}^{\star}\bm{\Sigma}^{\star}\right\Vert _{2}\lesssim\delta\frac{1}{np\sigma_{1}^{\star2}}\left\Vert \bm{U}_{l,\cdot}^{\star}\bm{\Sigma}^{\star}\right\Vert _{2}^{2}+\delta\frac{\omega_{l}^{\star2}}{np\sigma_{1}^{\star2}}\lesssim\delta\lambda_{\min}\left(\bm{\Sigma}_{U,l}^{\star}\right),
\end{align*}
provided that $\zeta_{\mathsf{1st}}/\sigma_{r}^{\star2}\lesssim1/\kappa$,
$d\gtrsim\kappa^{3}\mu r\log(n+d)$, $n\gtrsim\delta^{-2}\kappa^{5}r\log(n+d)$,
$np\gtrsim\delta^{-2}\log^{2}(n+d)$ and $\theta\lesssim\delta/\kappa$.

\paragraph{Step 3: controlling $\Vert\bm{A}_{3}-\bm{B}_{3}\Vert$.}

We have learned from (\ref{eq:pca-U-l-R-error}) in Lemma \ref{lemma:pca-1st-err}
we know that 
\begin{align*}
\left\Vert \bm{U}_{l,\cdot}\bm{R}-\bm{U}_{l,\cdot}^{\star}\right\Vert _{2} & \lesssim\frac{\theta}{\sqrt{\kappa}\sigma_{r}^{\star}}\left(\left\Vert \bm{U}_{l,\cdot}^{\star}\bm{\Sigma}^{\star}\right\Vert _{2}+\omega_{l}^{\star}\right)+\zeta_{\mathsf{2nd},l}.
\end{align*}
This immediately gives 
\begin{align}
\left\Vert \bm{U}_{l,\cdot}\right\Vert _{2} & =\left\Vert \bm{U}_{l,\cdot}\bm{R}\right\Vert _{2}\leq\left\Vert \bm{U}_{l,\cdot}^{\star}\right\Vert _{2}+\left\Vert \bm{U}_{l,\cdot}\bm{R}-\bm{U}_{l,\cdot}^{\star}\right\Vert _{2}\lesssim\left\Vert \bm{U}_{l,\cdot}^{\star}\right\Vert _{2}+\frac{\theta}{\sqrt{\kappa}\sigma_{r}^{\star}}\left(\left\Vert \bm{U}_{l,\cdot}^{\star}\bm{\Sigma}^{\star}\right\Vert _{2}+\omega_{l}^{\star}\right)+\zeta_{\mathsf{2nd}}\nonumber \\
 & \leq\left\Vert \bm{U}_{l,\cdot}^{\star}\right\Vert _{2}+\frac{\theta}{\sqrt{\kappa}\sigma_{r}^{\star}}\left(\left\Vert \bm{U}_{l,\cdot}^{\star}\right\Vert _{2}\sigma_{1}^{\star}+\omega_{l}^{\star}\right)+\zeta_{\mathsf{2nd}}\nonumber \\
 & \lesssim\left\Vert \bm{U}_{l,\cdot}^{\star}\right\Vert _{2}+\frac{\theta}{\sqrt{\kappa}\sigma_{r}^{\star}}\omega_{l}^{\star}+\zeta_{\mathsf{2nd},l},\label{eq:pca-cov-est-inter-4}
\end{align}
where the last line holds provided that $\theta\ll1$. As a consequence,
we arrive at 
\begin{align*}
\left\Vert \bm{A}_{3}-\bm{B}_{3}\right\Vert  & =\frac{2\left(1-p\right)}{np}\left\Vert \bm{R}\bm{U}_{l,\cdot}^{\star\top}\bm{U}_{l,\cdot}^{\star}\bm{R}^{\top}-\bm{U}_{l,\cdot}^{\top}\bm{U}_{l,\cdot}\right\Vert =\frac{2\left(1-p\right)}{np}\left\Vert \bm{U}_{l,\cdot}^{\star\top}\bm{U}_{l,\cdot}^{\star}-\bm{R}^{\top}\bm{U}_{l,\cdot}^{\top}\bm{U}_{l,\cdot}\bm{R}\right\Vert \\
 & \leq\frac{2\left(1-p\right)}{np}\left\Vert \left(\bm{U}^{\star}-\bm{U}\bm{R}\right)_{l,\cdot}^{\top}\bm{U}_{l,\cdot}^{\star}\right\Vert +\frac{2\left(1-p\right)}{np}\left\Vert \left(\bm{U}_{l}\bm{R}\right)^{\top}\left(\bm{U}_{l,\cdot}^{\star}-\bm{U}_{l,\cdot}\bm{R}\right)\right\Vert \\
 & \lesssim\frac{1}{np}\left\Vert \bm{U}_{l,\cdot}^{\star}-\bm{U}_{l,\cdot}\bm{R}\right\Vert _{2}\left(\left\Vert \bm{U}_{l,\cdot}^{\star}\right\Vert _{2}+\left\Vert \bm{U}_{l,\cdot}\right\Vert _{2}\right)\\
 & \overset{\text{(i)}}{\lesssim}\frac{1}{np}\left[\frac{\theta}{\sqrt{\kappa}\sigma_{r}^{\star}}\left(\left\Vert \bm{U}_{l,\cdot}^{\star}\bm{\Sigma}^{\star}\right\Vert _{2}+\omega_{l}^{\star}\right)+\zeta_{\mathsf{2nd},l}\right]\left(\left\Vert \bm{U}_{l,\cdot}^{\star}\right\Vert _{2}+\frac{\theta}{\sqrt{\kappa}\sigma_{r}^{\star}}\omega_{l}^{\star}+\zeta_{\mathsf{2nd},l}\right)\\
 & \overset{\text{(ii)}}{\lesssim}\underbrace{\frac{1}{np}\frac{\theta}{\sqrt{\kappa}\sigma_{r}^{\star2}}\left\Vert \bm{U}_{l,\cdot}^{\star}\bm{\Sigma}^{\star}\right\Vert _{2}^{2}}_{\eqqcolon\alpha_{3,1}}+\underbrace{\frac{1}{np}\frac{\theta}{\sqrt{\kappa}\sigma_{r}^{\star2}}\omega_{l}^{\star}\left\Vert \bm{U}_{l,\cdot}^{\star}\bm{\Sigma}^{\star}\right\Vert _{2}}_{\eqqcolon\alpha_{3,2}}+\frac{1}{np}\zeta_{\mathsf{2nd},l}\left\Vert \bm{U}_{l,\cdot}^{\star}\right\Vert _{2}+\underbrace{\frac{1}{np}\frac{\theta^{2}}{\kappa\sigma_{r}^{\star2}}\omega_{l}^{\star2}}_{\eqqcolon\alpha_{3,3}}\\
 & \quad+\frac{1}{np}\zeta_{\mathsf{2nd},l}\frac{\theta}{\sqrt{\kappa}\sigma_{r}^{\star}}\omega_{l}^{\star}+\frac{1}{np}\zeta_{\mathsf{2nd},l}^{2}\\
 & \overset{\text{(iii)}}{\lesssim}\delta\lambda_{\min}\left(\bm{\Sigma}_{U,l}^{\star}\right)+\frac{1}{np}\zeta_{\mathsf{2nd},l}\left\Vert \bm{U}_{l,\cdot}^{\star}\right\Vert _{2}+\frac{1}{np}\zeta_{\mathsf{2nd},l}\frac{\theta}{\sqrt{\kappa}\sigma_{r}^{\star}}\omega_{l}^{\star}+\frac{1}{np}\zeta_{\mathsf{2nd},l}^{2}.
\end{align*}
Here, (i) follows from (\ref{eq:pca-U-l-R-error}) and (\ref{eq:pca-cov-est-inter-4});
(ii) holds provided that $\theta\ll1$; and (iii) is valid due to
the following facts 
\begin{align*}
\alpha_{3,1} & \lesssim\delta\frac{1}{np\sigma_{1}^{\star2}}\left\Vert \bm{U}_{l,\cdot}^{\star}\bm{\Sigma}^{\star}\right\Vert _{2}^{2}\lesssim\delta\lambda_{\min}\left(\bm{\Sigma}_{U,l}^{\star}\right),\\
\alpha_{3,2} & \lesssim\delta\frac{1}{np\sigma_{1}^{\star2}}\omega_{l}^{\star}\left\Vert \bm{U}_{l,\cdot}^{\star}\bm{\Sigma}^{\star}\right\Vert _{2}\lesssim\delta\frac{1}{np\sigma_{1}^{\star2}}\left\Vert \bm{U}_{l,\cdot}^{\star}\bm{\Sigma}^{\star}\right\Vert _{2}^{2}+\delta\frac{\omega_{l}^{\star2}}{np\sigma_{1}^{\star2}}\lesssim\delta\lambda_{\min}\left(\bm{\Sigma}_{U,l}^{\star}\right),\\
\alpha_{3,3} & \lesssim\delta\frac{\omega_{l}^{\star2}}{np\sigma_{1}^{\star2}}\lesssim\delta\lambda_{\min}\left(\bm{\Sigma}_{U,l}^{\star}\right),
\end{align*}
provided that $\theta\lesssim\delta/\sqrt{\kappa}$.

\paragraph{Step 4: bounding $\Vert\bm{A}_{4}-\bm{B}_{4}\Vert$.}

To begin with, we recall from (\ref{eq:pca-U-Sigma-n2-R-error}) in
Lemma \ref{lemma:pca-1st-err} that 
\begin{equation}
\left\Vert \bm{U}^{\star}\left(\bm{\Sigma}^{\star}\right)^{-2}\bm{R}^{\top}-\bm{U}\bm{\Sigma}^{-2}\right\Vert \lesssim\frac{\zeta_{\mathsf{1st}}}{\sigma_{r}^{\star4}}+\sqrt{\frac{\kappa^{3}\left(r+\log\left(n+d\right)\right)}{n}}\frac{1}{\sigma_{r}^{\star2}}.\label{eq:pca-cov-est-inter-5}
\end{equation}
This allows one to upper bound 
\begin{align}
\left\Vert \bm{A}_{4}-\bm{B}_{4}\right\Vert  & \leq\left\Vert \left(\bm{U}^{\star}\left(\bm{\Sigma}^{\star}\right)^{-2}\bm{R}^{\top}-\bm{U}\bm{\Sigma}^{-2}\right)^{\top}\mathsf{diag}\left\{ \left[d_{l,i}^{\star}\right]_{i=1}^{d}\right\} \bm{U}^{\star}\left(\bm{\Sigma}^{\star}\right)^{-2}\bm{R}^{\top}\right\Vert \nonumber \\
 & \quad+\left\Vert \bm{U}\bm{\Sigma}^{-2}\mathsf{diag}\left\{ \left[d_{l,i}^{\star}-d_{l,i}\right]_{i=1}^{d}\right\} \bm{U}^{\star}\left(\bm{\Sigma}^{\star}\right)^{-2}\bm{R}^{\top}\right\Vert \nonumber \\
 & \quad+\left\Vert \bm{U}\bm{\Sigma}^{-2}\mathsf{diag}\left\{ \left[d_{l,i}\right]_{i=1}^{d}\right\} \left(\bm{U}^{\star}\left(\bm{\Sigma}^{\star}\right)^{-2}\bm{R}^{\top}-\bm{U}\bm{\Sigma}^{-2}\right)\right\Vert \nonumber \\
 & \overset{\text{(i)}}{\lesssim}\frac{1}{\sigma_{r}^{\star2}}\left\Vert \bm{U}^{\star}\left(\bm{\Sigma}^{\star}\right)^{-2}\bm{R}^{\top}-\bm{U}\bm{\Sigma}^{-2}\right\Vert \max_{1\leq i\leq d}\left(d_{l,i}^{\star}+d_{l,i}\right)+\frac{1}{\sigma_{r}^{\star4}}\max_{1\leq i\leq d}\left|d_{l,i}^{\star}-d_{l,i}\right|\nonumber \\
 & \overset{\text{(ii)}}{\lesssim}\left(\frac{\zeta_{\mathsf{1st}}}{\sigma_{r}^{\star6}}+\frac{1}{\sigma_{r}^{\star4}}\sqrt{\frac{\kappa^{3}\left(r+\log\left(n+d\right)\right)}{n}}\right)\max_{1\leq i\leq d}\left\{ d_{l,i}^{\star}+\left|d_{l,i}^{\star}-d_{l,i}\right|\right\} +\frac{1}{\sigma_{r}^{\star4}}\max_{1\leq i\leq d}\left|d_{l,i}^{\star}-d_{l,i}\right|\nonumber \\
 & \overset{\text{(iii)}}{\lesssim}\underbrace{\left(\frac{\zeta_{\mathsf{1st}}}{\sigma_{r}^{\star6}}+\frac{1}{\sigma_{r}^{\star4}}\sqrt{\frac{\kappa^{3}\left(r+\log\left(n+d\right)\right)}{n}}\right)\max_{1\leq i\leq d}d_{l,i}^{\star}}_{\eqqcolon\alpha_{4}}+\underbrace{\frac{1}{\sigma_{r}^{\star4}}\max_{1\leq i\leq d}\left|d_{l,i}^{\star}-d_{l,i}\right|}_{\eqqcolon\beta}.\label{eq:A4-B4-d-bound}
\end{align}
Here, (i) relies on (\ref{eq:pca-sigma-lower}); (ii) comes from (\ref{eq:pca-cov-est-inter-5});
(iii) holds true provided that $\zeta_{\mathsf{1st}}\lesssim\sigma_{r}^{\star2}$
and $n\gtrsim\kappa^{3}r\log(n+d)$. Note that, for each $i\in[d]$,
\begin{align}
d_{l,i}^{\star} & =\frac{1}{np^{2}}\left[\omega_{l}^{\star2}+\left(1-p\right)\left\Vert \bm{U}_{l,\cdot}^{\star}\bm{\Sigma}^{\star}\right\Vert _{2}^{2}\right]\left[\omega_{i}^{\star2}+\left(1-p\right)\left\Vert \bm{U}_{i,\cdot}^{\star}\bm{\Sigma}^{\star}\right\Vert _{2}^{2}\right]+\frac{2\left(1-p\right)^{2}}{np^{2}}\left(\bm{U}_{l,\cdot}^{\star}\bm{\Sigma}^{\star2}\bm{U}_{i,\cdot}^{\star\top}\right)^{2}\nonumber \\
 & \lesssim\frac{1}{np^{2}}\left[\omega_{l}^{\star2}+\left(1-p\right)\left\Vert \bm{U}_{l,\cdot}^{\star}\bm{\Sigma}^{\star}\right\Vert _{2}^{2}\right]\left[\omega_{\max}^{2}+\frac{\mu r}{d}\sigma_{1}^{\star2}\right]+\frac{\mu r}{ndp^{2}}\sigma_{1}^{\star2}\left\Vert \bm{U}_{l,\cdot}^{\star}\bm{\Sigma}^{\star}\right\Vert _{2}^{2}\nonumber \\
 & \lesssim\frac{1}{np^{2}}\omega_{l}^{\star2}\omega_{\max}^{2}+\frac{\mu r}{ndp^{2}}\sigma_{1}^{\star2}\omega_{l}^{\star2}+\frac{1}{np^{2}}\omega_{\max}^{2}\left\Vert \bm{U}_{l,\cdot}^{\star}\bm{\Sigma}^{\star}\right\Vert _{2}^{2}+\frac{\mu r}{ndp^{2}}\sigma_{1}^{\star2}\left\Vert \bm{U}_{l,\cdot}^{\star}\bm{\Sigma}^{\star}\right\Vert _{2}^{2},\label{eq:pca-cov-est-inter-6}
\end{align}
which in turn results in 
\begin{align*}
\alpha_{4} & \overset{\text{(i)}}{\lesssim}\underbrace{\left(\frac{\zeta_{\mathsf{1st}}}{\sigma_{r}^{\star6}}+\frac{1}{\sigma_{r}^{\star4}}\sqrt{\frac{\kappa^{3}r\log\left(n+d\right)}{n}}\right)\frac{\omega_{l}^{\star2}\omega_{\max}^{2}}{np^{2}}}_{\eqqcolon\alpha_{4,1}}+\underbrace{\frac{\kappa\mu r}{ndp^{2}}\left(\frac{\zeta_{\mathsf{1st}}}{\sigma_{r}^{\star4}}+\frac{1}{\sigma_{r}^{\star2}}\sqrt{\frac{\kappa^{3}r\log\left(n+d\right)}{n}}\right)\omega_{l}^{\star2}}_{\eqqcolon\alpha_{4,2}}\\
 & \quad+\underbrace{\frac{1}{np^{2}}\left(\frac{\zeta_{\mathsf{1st}}}{\sigma_{r}^{\star6}}+\frac{1}{\sigma_{r}^{\star4}}\sqrt{\frac{\kappa^{3}r\log\left(n+d\right)}{n}}\right)\omega_{\max}^{2}\left\Vert \bm{U}_{l,\cdot}^{\star}\bm{\Sigma}^{\star}\right\Vert _{2}^{2}}_{\eqqcolon\alpha_{4,3}}\\
 & \quad+\underbrace{\frac{\kappa\mu r}{ndp^{2}}\left(\frac{\zeta_{\mathsf{1st}}}{\sigma_{r}^{\star4}}+\frac{1}{\sigma_{r}^{\star2}}\sqrt{\frac{\kappa^{3}r\log\left(n+d\right)}{n}}\right)\left\Vert \bm{U}_{l,\cdot}^{\star}\bm{\Sigma}^{\star}\right\Vert _{2}^{2}}_{\eqqcolon\alpha_{4,4}}\\
 & \overset{\text{(ii)}}{\lesssim}\delta\lambda_{\min}\left(\bm{\Sigma}_{U,l}^{\star}\right).
\end{align*}
Here, (i) follows from (\ref{eq:pca-cov-est-inter-6}), while (ii)
holds since 
\begin{align*}
\alpha_{4,1} & \lesssim\delta\frac{\omega_{l}^{\star2}\omega_{\min}^{2}}{np^{2}\sigma_{1}^{\star4}}\lesssim\delta\lambda_{\min}\left(\bm{\Sigma}_{U,l}^{\star}\right),\\
\alpha_{4,2} & \lesssim\delta\frac{1}{ndp^{2}\kappa\sigma_{1}^{\star2}}\omega_{l}^{\star2}\lesssim\delta\lambda_{\min}\left(\bm{\Sigma}_{U,l}^{\star}\right),\\
\alpha_{4,3} & \lesssim\delta\frac{1}{np^{2}\sigma_{1}^{\star4}}\omega_{\min}^{2}\left\Vert \bm{U}_{l,\cdot}^{\star}\bm{\Sigma}^{\star}\right\Vert _{2}^{2}\lesssim\delta\lambda_{\min}\left(\bm{\Sigma}_{U,l}^{\star}\right),\\
\alpha_{4,4} & \lesssim\delta\frac{1}{ndp^{2}\kappa\sigma_{1}^{\star2}}\left\Vert \bm{U}_{l,\cdot}^{\star}\bm{\Sigma}^{\star}\right\Vert _{2}^{2}\lesssim\delta\lambda_{\min}\left(\bm{\Sigma}_{U,l}^{\star}\right),
\end{align*}
provided that $\zeta_{\mathsf{1st}}/\sigma_{r}^{\star2}\lesssim\delta/(\kappa^{3}\mu r\kappa_{\omega})$,
$n\gtrsim\delta^{-2}\kappa^{9}r\log(n+d)$ and $n\gtrsim\kappa^{7}r\kappa_{\omega}^{2}\log(n+d)$.
Note that we still need to bound the term $\beta$ in \eqref{eq:A4-B4-d-bound},
which we leave to the next step.

\paragraph{Step 5: bounding $\vert d_{l,i}^{\star}-d_{l,i}\vert$. }

For each $i\in[d]$, we can decompose 
\begin{align*}
\frac{1}{\sigma_{r}^{\star4}}\left|d_{l,i}^{\star}-d_{l,i}\right| & \leq\frac{1}{np^{2}\sigma_{r}^{\star4}}\left|\left[\omega_{l}^{\star2}+\left(1-p\right)S_{l,l}^{\star}\right]\left[\omega_{i}^{\star2}+\left(1-p\right)S_{i,i}^{\star}\right]-\left[\omega_{l}^{2}+\left(1-p\right)S_{l,l}\right]\left[\omega_{i}^{2}+\left(1-p\right)S_{i,i}\right]\right|\\
 & \quad\quad+\frac{2\left(1-p\right)^{2}}{np^{2}\sigma_{r}^{\star4}}\left|S_{i,l}^{\star2}-S_{i,l}^{2}\right|\\
 & \leq\underbrace{\frac{1}{np^{2}\sigma_{r}^{\star4}}\left|\omega_{l}^{\star2}+\left(1-p\right)S_{l,l}^{\star}\right|\left|\omega_{i}^{\star2}+\left(1-p\right)S_{i,i}^{\star}-\omega_{i}^{2}-\left(1-p\right)S_{i,i}\right|}_{\eqqcolon\beta_{1}}\\
 & \quad+\underbrace{\frac{1}{np^{2}\sigma_{r}^{\star4}}\left|\omega_{l}^{\star2}+\left(1-p\right)S_{l,l}^{\star}-\omega_{l}^{2}-\left(1-p\right)S_{l,l}\right|\left|\omega_{i}^{\star2}+\left(1-p\right)S_{i,i}^{\star}\right|}_{\eqqcolon\beta_{2}}\\
 & \quad+\underbrace{\frac{1}{np^{2}\sigma_{r}^{\star4}}\left|\omega_{l}^{\star2}+\left(1-p\right)S_{l,l}^{\star}-\omega_{l}^{2}-\left(1-p\right)S_{l,l}\right|\left|\omega_{i}^{2}+\left(1-p\right)S_{i,i}-\omega_{i}^{\star2}-\left(1-p\right)S_{i,i}^{\star}\right|}_{\eqqcolon\beta_{3}}\\
 & \quad+\underbrace{\frac{2\left(1-p\right)^{2}}{np^{2}\sigma_{r}^{\star4}}\left|S_{i,l}^{\star2}-S_{i,l}^{2}\right|}_{\eqqcolon\beta_{4}}.
\end{align*}
Denote $\Delta_{i}\coloneqq\omega_{i}^{\star2}+(1-p)S_{i,i}^{\star}-\omega_{i}^{2}-(1-p)S_{i,i}$.
We know from (\ref{eq:pca-noise-est-all}) and (\ref{eq:pca-S-infty-err})
in Lemma \ref{lemma:pca-noise-level-est} that for each $i\in[d]$,
\begin{align}
\left|\Delta_{i}\right| & \leq\left|\omega_{i}^{\star2}-\omega_{i}^{2}\right|+\left\Vert \bm{S}-\bm{S}^{\star}\right\Vert _{\infty}\nonumber \\
 & \lesssim\sqrt{\frac{\log^{2}\left(n+d\right)}{np}}\omega_{\max}^{\star2}+\zeta_{\mathsf{1st}}\frac{\sqrt{\kappa^{2}\mu r^{2}\log\left(n+d\right)}}{d}+\sqrt{\frac{\kappa^{2}\mu^{2}r^{3}\log\left(n+d\right)}{nd^{2}}}\sigma_{1}^{\star2};\label{eq:pca-cov-est-inter-7}
\end{align}
and we have also learned from (\ref{eq:pca-noise-est-l}) and (\ref{eq:pca-S-entrywise-err})
that 
\begin{align}
\left|\Delta_{l}\right| & \leq\left|\omega_{l}^{\star2}-\omega_{l}^{2}\right|+\left|S_{l,l}-S_{l,l}^{\star}\right|\nonumber \\
 & \lesssim\left(\sqrt{\frac{\log^{2}\left(n+d\right)}{np}}+\theta^{2}\right)\omega_{l}^{\star2}+\left(\theta+\sqrt{\frac{\kappa^{3}r\log\left(n+d\right)}{n}}\right)\left\Vert \bm{U}_{l,\cdot}^{\star}\bm{\Sigma}^{\star}\right\Vert _{2}^{2}\nonumber \\
 & \quad+\theta\omega_{l}^{\star}\left\Vert \bm{U}_{l,\cdot}^{\star}\bm{\Sigma}^{\star}\right\Vert _{2}+\zeta_{\mathsf{2nd},l}\sigma_{1}^{\star}\left\Vert \bm{U}_{l,\cdot}^{\star}\bm{\Sigma}^{\star}\right\Vert _{2}+\zeta_{\mathsf{2nd},l}^{2}\sigma_{1}^{\star2}.\label{eq:pca-cov-est-inter-8}
\end{align}
In addition, it is straightforward to verify that for each $i\in[d]$,
\begin{equation}
\left|\omega_{i}^{\star2}+\left(1-p\right)S_{i,i}^{\star}\right|\leq\omega_{\max}^{2}+\frac{\mu r}{d}\sigma_{1}^{\star2}\label{eq:pca-cov-est-inter-9}
\end{equation}
and 
\begin{equation}
\left|\omega_{l}^{\star2}+\left(1-p\right)S_{l,l}^{\star}\right|\leq\omega_{l}^{\star2}+\left\Vert \bm{U}_{l,\cdot}^{\star}\bm{\Sigma}^{\star}\right\Vert _{2}^{2}.\label{eq:pca-cov-est-inter-10}
\end{equation}

\begin{itemize}
\item Regarding $\beta_{1}$, we can derive 
\begin{align*}
\beta_{1} & =\frac{1}{np^{2}\sigma_{r}^{\star4}}\left|\omega_{l}^{\star2}+\left(1-p\right)S_{l,l}^{\star}\right|\left|\Delta_{i}\right|\\
 & \overset{\text{(i)}}{\lesssim}\frac{1}{np^{2}\sigma_{r}^{\star4}}\left(\omega_{l}^{\star2}+\left\Vert \bm{U}_{l,\cdot}^{\star}\bm{\Sigma}^{\star}\right\Vert _{2}^{2}\right)\left(\sqrt{\frac{\log^{2}\left(n+d\right)}{np}}\omega_{\max}^{\star2}+\zeta_{\mathsf{1st}}\frac{\sqrt{\kappa^{2}\mu r^{2}\log\left(n+d\right)}}{d}+\sqrt{\frac{\kappa^{2}\mu^{2}r^{3}\log\left(n+d\right)}{nd^{2}}}\sigma_{1}^{\star2}\right)\\
 & \asymp\underbrace{\frac{\omega_{l}^{\star2}\omega_{\max}^{\star2}}{np^{2}\sigma_{r}^{\star4}}\sqrt{\frac{\log^{2}\left(n+d\right)}{np}}}_{\eqqcolon\beta_{1,1}}+\underbrace{\frac{\omega_{l}^{\star2}}{np^{2}\sigma_{r}^{\star4}}\zeta_{\mathsf{1st}}\frac{\sqrt{\kappa^{2}\mu r^{2}\log\left(n+d\right)}}{d}}_{\eqqcolon\beta_{1,2}}+\underbrace{\frac{\omega_{l}^{\star2}}{np^{2}\sigma_{r}^{\star2}}\sqrt{\frac{\kappa^{4}\mu^{2}r^{3}\log\left(n+d\right)}{nd^{2}}}}_{\eqqcolon\beta_{1,3}}\\
 & \quad+\underbrace{\frac{\omega_{\max}^{\star2}}{np^{2}\sigma_{r}^{\star4}}\sqrt{\frac{\log^{2}\left(n+d\right)}{np}}\left\Vert \bm{U}_{l,\cdot}^{\star}\bm{\Sigma}^{\star}\right\Vert _{2}^{2}}_{\eqqcolon\beta_{1,4}}+\underbrace{\frac{1}{np^{2}\sigma_{r}^{\star4}}\zeta_{\mathsf{1st}}\frac{\sqrt{\kappa^{2}\mu r^{2}\log\left(n+d\right)}}{d}\left\Vert \bm{U}_{l,\cdot}^{\star}\bm{\Sigma}^{\star}\right\Vert _{2}^{2}}_{\eqqcolon\beta_{1,5}}\\
 & \quad+\underbrace{\frac{1}{np^{2}\sigma_{r}^{\star2}}\sqrt{\frac{\kappa^{4}\mu^{2}r^{3}\log\left(n+d\right)}{nd^{2}}}\left\Vert \bm{U}_{l,\cdot}^{\star}\bm{\Sigma}^{\star}\right\Vert _{2}^{2}}_{\eqqcolon\beta_{1,6}}\overset{\text{(ii)}}{\lesssim}\delta\lambda_{\min}\left(\bm{\Sigma}_{U,l}^{\star}\right).
\end{align*}
Here, (i) follows from (\ref{eq:pca-cov-est-inter-7}) and (\ref{eq:pca-cov-est-inter-10});
(ii) holds since 
\begin{align*}
\beta_{1,1} & \lesssim\delta\frac{\omega_{l}^{\star2}\omega_{\min}^{2}}{np^{2}\sigma_{1}^{\star4}}\lesssim\delta\lambda_{\min}\left(\bm{\Sigma}_{U,l}^{\star}\right),\\
\beta_{1,2}+\beta_{1,3}+\beta_{1,4} & \lesssim\delta\frac{1}{ndp^{2}\kappa\sigma_{1}^{\star2}}\omega_{l}^{\star2}\lesssim\delta\lambda_{\min}\left(\bm{\Sigma}_{U,l}^{\star}\right),\\
\beta_{1,5}+\beta_{1,6} & \lesssim\delta\frac{1}{ndp^{2}\kappa\sigma_{1}^{\star2}}\left\Vert \bm{U}_{l,\cdot}^{\star}\bm{\Sigma}^{\star}\right\Vert _{2}^{2}\lesssim\delta\lambda_{\min}\left(\bm{\Sigma}_{U,l}^{\star}\right),
\end{align*}
provided that $np\gtrsim\delta^{-2}\kappa^{6}\mu^{2}r^{2}\kappa_{\omega}^{2}\log^{2}(n+d)$,
$\zeta_{\mathsf{1st}}/\sigma_{r}^{\star2}\lesssim\delta/\sqrt{\kappa^{6}\mu r^{2}\log(n+d)}$,
$n\gtrsim\delta^{-2}\kappa^{8}\mu^{2}r^{3}\log(n+d)$. 
\item When it comes to $\beta_{2}$, we can see that 
\begin{align*}
\beta_{2} & =\frac{1}{np^{2}\sigma_{r}^{\star4}}\left|\omega_{i}^{\star2}+\left(1-p\right)S_{i,i}^{\star}\right|\left|\Delta_{l}\right|\\
 & \overset{\text{(i)}}{\lesssim}\underbrace{\frac{\omega_{l}^{\star2}\omega_{\max}^{2}}{np^{2}\sigma_{r}^{\star4}}\left(\sqrt{\frac{\log^{2}\left(n+d\right)}{np}}+\theta^{2}\right)}_{\eqqcolon\beta_{2,1}}+\underbrace{\frac{\kappa\mu r}{ndp^{2}\sigma_{r}^{\star2}}\omega_{l}^{\star2}\left(\sqrt{\frac{\log^{2}\left(n+d\right)}{np}}+\theta^{2}\right)}_{\eqqcolon\beta_{2,2}}\\
 & \quad+\underbrace{\frac{1}{np^{2}\sigma_{r}^{\star4}}\omega_{\max}^{2}\left(\theta+\sqrt{\frac{\kappa^{3}r\log\left(n+d\right)}{n}}\right)\left\Vert \bm{U}_{l,\cdot}^{\star}\bm{\Sigma}^{\star}\right\Vert _{2}^{2}}_{\eqqcolon\beta_{2,3}}+\underbrace{\frac{\kappa\mu r}{ndp^{2}\sigma_{r}^{\star2}}\left(\theta+\sqrt{\frac{\kappa^{3}r\log\left(n+d\right)}{n}}\right)\left\Vert \bm{U}_{l,\cdot}^{\star}\bm{\Sigma}^{\star}\right\Vert _{2}^{2}}_{\eqqcolon\beta_{2,4}}\\
 & \quad+\underbrace{\frac{1}{np^{2}\sigma_{r}^{\star4}}\omega_{\max}^{2}\omega_{l}^{\star}\theta\left\Vert \bm{U}_{l,\cdot}^{\star}\bm{\Sigma}^{\star}\right\Vert _{2}}_{\beta_{2,5}}+\frac{\sqrt{\kappa}}{np^{2}\sigma_{r}^{\star3}}\omega_{\max}^{2}\zeta_{\mathsf{2nd},l}\left\Vert \bm{U}_{l,\cdot}^{\star}\bm{\Sigma}^{\star}\right\Vert _{2}+\frac{\kappa}{np^{2}\sigma_{r}^{\star2}}\omega_{\max}^{2}\zeta_{\mathsf{2nd},l}^{2}\\
 & \quad+\underbrace{\frac{\kappa\mu r}{ndp^{2}\sigma_{r}^{\star2}}\theta\omega_{l}^{\star}\left\Vert \bm{U}_{l,\cdot}^{\star}\bm{\Sigma}^{\star}\right\Vert _{2}}_{\eqqcolon\beta_{2,6}}+\frac{\kappa^{3/2}\mu r}{ndp^{2}\sigma_{r}^{\star}}\zeta_{\mathsf{2nd},l}\left\Vert \bm{U}_{l,\cdot}^{\star}\bm{\Sigma}^{\star}\right\Vert _{2}+\frac{\kappa^{2}\mu r}{ndp^{2}}\zeta_{\mathsf{2nd},l}^{2}\\
 & \overset{\text{(ii)}}{\lesssim}\delta\lambda_{\min}\left(\bm{\Sigma}_{U,l}^{\star}\right)+\frac{\sqrt{\kappa}}{np^{2}\sigma_{r}^{\star3}}\omega_{\max}^{2}\zeta_{\mathsf{2nd},l}\left\Vert \bm{U}_{l,\cdot}^{\star}\bm{\Sigma}^{\star}\right\Vert _{2}+\frac{\kappa}{np^{2}\sigma_{r}^{\star2}}\omega_{\max}^{2}\zeta_{\mathsf{2nd},l}^{2}+\frac{\kappa^{3/2}\mu r}{ndp^{2}\sigma_{r}^{\star}}\zeta_{\mathsf{2nd},l}\left\Vert \bm{U}_{l,\cdot}^{\star}\bm{\Sigma}^{\star}\right\Vert _{2}+\frac{\kappa^{2}\mu r}{ndp^{2}}\zeta_{\mathsf{2nd},l}^{2}.
\end{align*}
Here, (i) follows from (\ref{eq:pca-cov-est-inter-8}) and (\ref{eq:pca-cov-est-inter-9});
(ii) holds since 
\begin{align*}
\beta_{2,1} & \lesssim\delta\frac{\omega_{l}^{\star2}\omega_{\min}^{2}}{np^{2}\sigma_{1}^{\star4}}\lesssim\delta\lambda_{\min}\left(\bm{\Sigma}_{U,l}^{\star}\right),\\
\beta_{2,2} & \lesssim\delta\frac{1}{ndp^{2}\kappa\sigma_{1}^{\star2}}\omega_{l}^{\star2}\lesssim\delta\lambda_{\min}\left(\bm{\Sigma}_{U,l}^{\star}\right),\\
\beta_{2,3} & \lesssim\delta\frac{1}{np^{2}\sigma_{1}^{\star4}}\omega_{\min}^{2}\left\Vert \bm{U}_{l,\cdot}^{\star}\bm{\Sigma}^{\star}\right\Vert _{2}^{2}\lesssim\delta\lambda_{\min}\left(\bm{\Sigma}_{U,l}^{\star}\right),\\
\beta_{2,4} & \lesssim\delta\frac{1}{ndp^{2}\kappa\sigma_{1}^{\star2}}\left\Vert \bm{U}_{l,\cdot}^{\star}\bm{\Sigma}^{\star}\right\Vert _{2}^{2}\lesssim\delta\lambda_{\min}\left(\bm{\Sigma}_{U,l}^{\star}\right),,\\
\beta_{2,5} & \lesssim\delta\frac{1}{np^{2}\sigma_{1}^{\star4}}\omega_{\min}^{2}\omega_{l}^{\star}\left\Vert \bm{U}_{l,\cdot}^{\star}\bm{\Sigma}^{\star}\right\Vert _{2}\lesssim\delta\frac{1}{np^{2}\sigma_{1}^{\star4}}\omega_{\min}^{2}\left\Vert \bm{U}_{l,\cdot}^{\star}\bm{\Sigma}^{\star}\right\Vert _{2}^{2}+\delta\frac{\omega_{l}^{\star2}\omega_{\min}^{2}}{np^{2}\sigma_{1}^{\star4}}\lesssim\delta\lambda_{\min}\left(\bm{\Sigma}_{U,l}^{\star}\right),\\
\beta_{2,6} & \lesssim\delta\frac{1}{ndp^{2}\kappa\sigma_{1}^{\star2}}\left\Vert \bm{U}_{l,\cdot}^{\star}\bm{\Sigma}^{\star}\right\Vert _{2}\omega_{l}^{\star}\lesssim\delta\frac{1}{ndp^{2}\kappa\sigma_{1}^{\star2}}\left\Vert \bm{U}_{l,\cdot}^{\star}\bm{\Sigma}^{\star}\right\Vert _{2}^{2}+\delta\frac{1}{ndp^{2}\kappa\sigma_{1}^{\star2}}\omega_{l}^{\star2}\lesssim\delta\lambda_{\min}\left(\bm{\Sigma}_{U,l}^{\star}\right),
\end{align*}
provided that $np\gtrsim\delta^{-2}\kappa^{6}\mu^{2}r^{2}\log^{2}(n+d)$,
$\theta\lesssim\delta/(\kappa^{3}\mu r\kappa_{\omega})$, $n\gtrsim\delta^{-2}\kappa^{9}\mu^{2}r^{3}\log(n+d)$,
$n\gtrsim\delta^{-2}\kappa^{7}r\kappa_{\omega}^{2}\log(n+d)$. 
\item With regards to $\beta_{3}$, we notice that 
\[
\left|\Delta_{i}\right|\lesssim\sqrt{\frac{\log^{2}\left(n+d\right)}{np}}\omega_{\max}^{\star2}+\zeta_{\mathsf{1st}}\frac{\sqrt{\kappa^{2}\mu r^{2}\log\left(n+d\right)}}{d}+\sqrt{\frac{\kappa^{2}\mu^{2}r^{3}\log\left(n+d\right)}{nd^{2}}}\sigma_{1}^{\star2}\lesssim\omega_{\max}^{2}+\frac{\mu r}{d}\sigma_{1}^{\star2}
\]
provided that $np\gtrsim\log^{2}(n+d)$, $n\gtrsim\kappa^{2}r\log(n+d)$
and $\zeta_{\mathsf{1st}}/\sigma_{r}^{\star2}\lesssim1/\sqrt{\log(n+d)}$.
As a consequence, we can upper bound 
\[
\beta_{3}=\frac{1}{np^{2}\sigma_{r}^{\star4}}\left|\Delta_{i}\right|\left|\Delta_{l}\right|\lesssim\frac{1}{np^{2}\sigma_{r}^{\star4}}\left|\omega_{i}^{\star2}+\left(1-p\right)S_{i,i}^{\star}\right|\left|\Delta_{l}\right|.
\]
This immediately suggests that $\beta_{3}$ satisfies the same upper
bound we derive for $\beta_{2}$. 
\item Regarding $\beta_{4}$, it is seen that 
\begin{align*}
\beta_{4} & \lesssim\frac{1}{np^{2}\sigma_{r}^{\star4}}\left|S_{i,l}^{\star2}-S_{i,l}^{2}\right|\lesssim\frac{1}{np^{2}\sigma_{r}^{\star4}}\left|S_{i,l}^{\star}-S_{i,l}\right|\left|S_{i,l}^{\star}+S_{i,l}\right|\\
 & \lesssim\underbrace{\frac{1}{np^{2}\sigma_{r}^{\star4}}\left|S_{i,l}^{\star}-S_{i,l}\right|S_{i,l}^{\star}}_{\eqqcolon\beta_{4,1}}+\underbrace{\frac{1}{np^{2}\sigma_{r}^{\star4}}\left|S_{i,l}^{\star}-S_{i,l}\right|^{2}}_{\eqqcolon\beta_{4,2}}.
\end{align*}
Recall from (\ref{eq:pca-S-entrywise-err}) in Lemma \ref{lemma:pca-noise-level-est}
that 
\begin{align}
\left|S_{i,l}-S_{i,l}^{\star}\right| & \lesssim\left(\theta+\sqrt{\frac{\kappa^{3}r\log\left(n+d\right)}{n}}\right)\left\Vert \bm{U}_{i,\cdot}^{\star}\bm{\Sigma}^{\star}\right\Vert _{2}\left\Vert \bm{U}_{l,\cdot}^{\star}\bm{\Sigma}^{\star}\right\Vert _{2}+\theta\left(\omega_{i}^{\star}\left\Vert \bm{U}_{l,\cdot}^{\star}\bm{\Sigma}^{\star}\right\Vert _{2}+\omega_{l}^{\star}\left\Vert \bm{U}_{i,\cdot}^{\star}\bm{\Sigma}^{\star}\right\Vert _{2}\right)\nonumber \\
 & \quad+\sigma_{1}^{\star}\left(\zeta_{\mathsf{2nd},i}\left\Vert \bm{U}_{l,\cdot}^{\star}\bm{\Sigma}^{\star}\right\Vert _{2}+\zeta_{\mathsf{2nd},l}\left\Vert \bm{U}_{i,\cdot}^{\star}\bm{\Sigma}^{\star}\right\Vert _{2}\right)+\theta^{2}\omega_{\max}\omega_{l}^{\star}+\zeta_{\mathsf{2nd},i}\zeta_{\mathsf{2nd},l}\sigma_{1}^{\star2}\nonumber \\
 & \lesssim\left(\theta+\sqrt{\frac{\kappa^{3}r\log\left(n+d\right)}{n}}\right)\sqrt{\frac{\mu r}{d}}\sigma_{1}^{\star}\left\Vert \bm{U}_{l,\cdot}^{\star}\bm{\Sigma}^{\star}\right\Vert _{2}+\theta^{2}\omega_{\max}\omega_{l}^{\star}+\theta\omega_{l}^{\star}\sqrt{\frac{\mu r}{d}}\sigma_{1}^{\star}\nonumber \\
 & \quad+\theta\omega_{\max}\left\Vert \bm{U}_{l,\cdot}^{\star}\bm{\Sigma}^{\star}\right\Vert _{2}+\sigma_{1}^{\star}\left\Vert \bm{U}_{l,\cdot}^{\star}\bm{\Sigma}^{\star}\right\Vert _{2}\max_{i\in[d]}\zeta_{\mathsf{2nd},i}+\zeta_{\mathsf{2nd},l}\sqrt{\frac{\mu r}{d}}\sigma_{1}^{\star2}.\label{eq:pca-cov-est-inter-11}
\end{align}
provided that $\max_{i\in[d]}\zeta_{\mathsf{2nd},i}\sqrt{d}\leq\sqrt{\mu r}$.
The first term $\beta_{4,1}$ can be upper bounded by 
\begin{align*}
\beta_{4,1} & \lesssim\frac{1}{np^{2}\sigma_{r}^{\star4}}\left|S_{i,l}^{\star}-S_{i,l}\right|\left\Vert \bm{U}_{i,\cdot}^{\star}\bm{\Sigma}^{\star}\right\Vert _{2}\left\Vert \bm{U}_{l,\cdot}^{\star}\bm{\Sigma}^{\star}\right\Vert _{2}\lesssim\frac{1}{np^{2}\sigma_{r}^{\star4}}\left|S_{i,l}^{\star}-S_{i,l}\right|\sqrt{\frac{\mu r}{d}}\sigma_{1}^{\star}\left\Vert \bm{U}_{l,\cdot}^{\star}\bm{\Sigma}^{\star}\right\Vert _{2}\\
 & \overset{\text{(i)}}{\lesssim}\underbrace{\frac{\kappa\mu r}{ndp^{2}\sigma_{r}^{\star2}}\left(\theta+\sqrt{\frac{\kappa^{3}r\log\left(n+d\right)}{n}}\right)\left\Vert \bm{U}_{l,\cdot}^{\star}\bm{\Sigma}^{\star}\right\Vert _{2}^{2}}_{\eqqcolon\beta_{4,1,1}}+\underbrace{\frac{1}{np^{2}\sigma_{r}^{\star4}}\theta^{2}\omega_{\max}\omega_{l}^{\star}\sqrt{\frac{\mu r}{d}}\sigma_{1}^{\star}\left\Vert \bm{U}_{l,\cdot}^{\star}\bm{\Sigma}^{\star}\right\Vert _{2}}_{\eqqcolon\beta_{4,1,2}}\\
 & \quad+\underbrace{\frac{\kappa\mu r}{ndp^{2}\sigma_{r}^{\star2}}\theta\omega_{l}^{\star}\left\Vert \bm{U}_{l,\cdot}^{\star}\bm{\Sigma}^{\star}\right\Vert _{2}}_{\eqqcolon\beta_{4,1,3}}+\underbrace{\frac{\theta\omega_{\max}}{np^{2}\sigma_{r}^{\star3}}\sqrt{\frac{\kappa\mu r}{d}}\left\Vert \bm{U}_{l,\cdot}^{\star}\bm{\Sigma}^{\star}\right\Vert _{2}^{2}}_{\eqqcolon\beta_{4,1,4}}\\
 & \quad+\underbrace{\frac{1}{np^{2}\sigma_{r}^{\star2}}\sqrt{\frac{\kappa^{2}\mu r}{d}}\left\Vert \bm{U}_{l,\cdot}^{\star}\bm{\Sigma}^{\star}\right\Vert _{2}^{2}\max_{i\in[d]}\zeta_{\mathsf{2nd},i}}_{\eqqcolon\beta_{4,1,5}}+\frac{\kappa^{3/2}\mu r}{ndp^{2}\sigma_{r}^{\star}}\zeta_{\mathsf{2nd},l}\left\Vert \bm{U}_{l,\cdot}^{\star}\bm{\Sigma}^{\star}\right\Vert _{2}\\
 & \overset{\text{(ii)}}{\lesssim}\delta\lambda_{\min}\left(\bm{\Sigma}_{U,l}^{\star}\right)+\frac{\kappa^{3/2}\mu r}{ndp^{2}\sigma_{r}^{\star}}\zeta_{\mathsf{2nd},l}\left\Vert \bm{U}_{l,\cdot}^{\star}\bm{\Sigma}^{\star}\right\Vert _{2}.
\end{align*}
Here, (i) follows from (\ref{eq:pca-cov-est-inter-11}), and (ii)
holds since 
\begin{align*}
\beta_{4,1,1}+\beta_{4,1,5} & \lesssim\delta\frac{1}{ndp^{2}\kappa\sigma_{1}^{\star2}}\left\Vert \bm{U}_{l,\cdot}^{\star}\bm{\Sigma}^{\star}\right\Vert _{2}^{2}\lesssim\delta\lambda_{\min}\left(\bm{\Sigma}_{U,l}^{\star}\right),\\
\beta_{4,1,2} & \lesssim\delta\frac{1}{ndp^{2}\kappa\sigma_{1}^{\star2}}\omega_{l}^{\star2}\lesssim\delta\lambda_{\min}\left(\bm{\Sigma}_{U,l}^{\star}\right),\\
\beta_{4,1,3}+\beta_{4,1,4} & \lesssim\frac{\delta}{ndp^{2}\kappa\sigma_{1}^{\star2}}\left\Vert \bm{U}_{l,\cdot}^{\star}\bm{\Sigma}^{\star}\right\Vert _{2}\omega_{l}^{\star}\lesssim\frac{\delta}{ndp^{2}\kappa\sigma_{1}^{\star2}}\left\Vert \bm{U}_{l,\cdot}^{\star}\bm{\Sigma}^{\star}\right\Vert _{2}^{2}+\frac{\delta}{ndp^{2}\kappa\sigma_{1}^{\star2}}\omega_{l}^{\star2}\lesssim\delta\lambda_{\min}\left(\bm{\Sigma}_{U,l}^{\star}\right),
\end{align*}
provided that $\theta\lesssim\delta/(\kappa^{3}\mu r\sqrt{\kappa_{\omega}})$,
$n\gtrsim\delta^{-2}\kappa^{9}\mu^{2}r^{3}\log(n+d)$ and $\max_{i\in[d]}\zeta_{\mathsf{2nd},i}\sqrt{d}\lesssim\delta/\sqrt{\kappa^{6}\mu r}$.
The second term $\beta_{4,2}$ can be controlled as follows 
\begin{align*}
\beta_{4,2} & \overset{\text{(i)}}{\lesssim}\underbrace{\frac{\kappa\mu r}{ndp^{2}\sigma_{r}^{\star2}}\left(\theta+\sqrt{\frac{\kappa^{3}r\log\left(n+d\right)}{n}}\right)^{2}\left\Vert \bm{U}_{l,\cdot}^{\star}\bm{\Sigma}^{\star}\right\Vert _{2}^{2}}_{\eqqcolon\beta_{4,2,1}}+\underbrace{\frac{1}{np^{2}\sigma_{r}^{\star4}}\theta^{4}\omega_{\max}^{2}\omega_{l}^{\star2}}_{\eqqcolon\beta_{4,2,2}}+\underbrace{\frac{\kappa\mu r}{ndp^{2}\sigma_{r}^{\star2}}\theta^{2}\omega_{l}^{\star2}}_{\eqqcolon\beta_{4,2,3}}\\
 & \quad+\underbrace{\frac{1}{np^{2}\sigma_{r}^{\star4}}\theta^{2}\omega_{\max}^{2}\left\Vert \bm{U}_{l,\cdot}^{\star}\bm{\Sigma}^{\star}\right\Vert _{2}^{2}}_{\eqqcolon\beta_{4,2,4}}+\underbrace{\frac{\kappa}{np^{2}\sigma_{r}^{\star2}}\left\Vert \bm{U}_{l,\cdot}^{\star}\bm{\Sigma}^{\star}\right\Vert _{2}^{2}\max_{i\in[d]}\zeta_{\mathsf{2nd},i}^{2}}_{\eqqcolon\beta_{4,2,5}}+\frac{\kappa^{2}\mu r}{ndp^{2}}\zeta_{\mathsf{2nd},l}^{2}\\
 & \overset{\text{(ii)}}{\lesssim}\delta\lambda_{\min}\left(\bm{\Sigma}_{U,l}^{\star}\right)+\frac{\kappa^{2}\mu r}{ndp^{2}}\zeta_{\mathsf{2nd},l}^{2}.
\end{align*}
Here, (i) follows from (\ref{eq:pca-cov-est-inter-11}), and (ii)
holds since 
\begin{align*}
\beta_{4,2,1}+\beta_{4,2,5} & \lesssim\delta\frac{1}{ndp^{2}\kappa\sigma_{1}^{\star2}}\left\Vert \bm{U}_{l,\cdot}^{\star}\bm{\Sigma}^{\star}\right\Vert _{2}^{2}\lesssim\delta\lambda_{\min}\left(\bm{\Sigma}_{U,l}^{\star}\right),\\
\beta_{4,2,2} & \lesssim\delta\frac{\omega_{l}^{\star2}\omega_{\min}^{2}}{np^{2}\sigma_{1}^{\star4}}\lesssim\delta\lambda_{\min}\left(\bm{\Sigma}_{U,l}^{\star}\right),\\
\beta_{4,2,3} & \lesssim\delta\frac{1}{ndp^{2}\kappa\sigma_{1}^{\star2}}\omega_{l}^{\star2}\lesssim\delta\frac{1}{ndp^{2}\kappa\sigma_{1}^{\star2}}\left\Vert \bm{U}_{l,\cdot}^{\star}\bm{\Sigma}^{\star}\right\Vert _{2}^{2}+\delta\frac{1}{4ndp^{2}\kappa\sigma_{1}^{\star2}}\omega_{l}^{\star2}\lesssim\delta\lambda_{\min}\left(\bm{\Sigma}_{U,l}^{\star}\right),\\
\beta_{4,2,4} & \lesssim\delta\frac{1}{np^{2}\sigma_{1}^{\star4}}\omega_{\min}^{2}\left\Vert \bm{U}_{l,\cdot}^{\star}\bm{\Sigma}^{\star}\right\Vert _{2}^{2}\lesssim\delta\lambda_{\min}\left(\bm{\Sigma}_{U,l}^{\star}\right),
\end{align*}
provided that $\theta\lesssim\sqrt{\delta/(\kappa^{3}\mu r\kappa_{\omega})}$,
$n\gtrsim\delta^{-1}\kappa^{6}\mu r^{2}\log(n+d)$ and $\max_{i\in[d]}\zeta_{\mathsf{2nd},i}\sqrt{d}\lesssim\sqrt{\delta/\kappa^{3}}$.
Combine the bounds on $\beta_{4,1}$ and $\beta_{4,2}$ to reach 
\begin{align*}
\beta_{4} & \leq\beta_{4,1}+\beta_{4,2}\lesssim\delta\lambda_{\min}\left(\bm{\Sigma}_{U,l}^{\star}\right)+\frac{\kappa^{3/2}\mu r}{ndp^{2}\sigma_{r}^{\star}}\zeta_{\mathsf{2nd},l}\left\Vert \bm{U}_{l,\cdot}^{\star}\bm{\Sigma}^{\star}\right\Vert _{2}+\frac{\kappa^{2}\mu r}{ndp^{2}}\zeta_{\mathsf{2nd},l}^{2}.
\end{align*}
\end{itemize}
Taking together the bounds on $\beta_{1}$, $\beta_{2}$, $\beta_{3}$
and $\beta_{4}$ yields 
\begin{align*}
\frac{1}{\sigma_{r}^{\star4}}\left|d_{l,i}^{\star}-d_{l,i}\right| & \leq\beta_{1}+\beta_{2}+\beta_{3}+\beta_{4}\\
 & \lesssim\delta\lambda_{\min}\left(\bm{\Sigma}_{U,l}^{\star}\right)+\frac{\sqrt{\kappa}}{np^{2}\sigma_{r}^{\star3}}\omega_{\max}^{2}\zeta_{\mathsf{2nd},l}\left\Vert \bm{U}_{l,\cdot}^{\star}\bm{\Sigma}^{\star}\right\Vert _{2}+\frac{\kappa}{np^{2}\sigma_{r}^{\star2}}\omega_{\max}^{2}\zeta_{\mathsf{2nd},l}^{2}\\
 & \quad+\frac{\kappa^{3/2}\mu r}{ndp^{2}\sigma_{r}^{\star}}\zeta_{\mathsf{2nd},l}\left\Vert \bm{U}_{l,\cdot}^{\star}\bm{\Sigma}^{\star}\right\Vert _{2}+\frac{\kappa^{2}\mu r}{ndp^{2}}\zeta_{\mathsf{2nd},l}^{2},
\end{align*}
with the proviso that $np\gtrsim\delta^{-2}\kappa^{6}\mu^{2}r^{2}\kappa_{\omega}^{2}\log^{2}(n+d)$,
$\zeta_{\mathsf{1st}}/\sigma_{r}^{\star2}\lesssim\delta/\sqrt{\kappa^{6}\mu r^{2}\log(n+d)}$,
$\theta\lesssim\delta/(\kappa^{3}\mu r\kappa_{\omega})$, $\max_{i\in[d]}\zeta_{\mathsf{2nd},i}\sqrt{d}\lesssim\delta/\sqrt{\kappa^{6}\mu r}$,
$n\gtrsim\delta^{-2}\kappa^{9}\mu^{2}r^{3}\log(n+d)$ and $n\gtrsim\kappa^{7}r\kappa_{\omega}^{2}\log(n+d)$.
Given that this holds for all $i\in[d]$, it follows that 
\begin{align*}
\beta & =\frac{1}{\sigma_{r}^{\star4}}\max_{1\leq i\leq d}\left|d_{l,i}^{\star}-d_{l,i}\right|\\
 & \lesssim\delta\lambda_{\min}\left(\bm{\Sigma}_{U,l}^{\star}\right)+\frac{\sqrt{\kappa}}{np^{2}\sigma_{r}^{\star3}}\omega_{\max}^{2}\zeta_{\mathsf{2nd},l}\left\Vert \bm{U}_{l,\cdot}^{\star}\bm{\Sigma}^{\star}\right\Vert _{2}+\frac{\kappa}{np^{2}\sigma_{r}^{\star2}}\omega_{\max}^{2}\zeta_{\mathsf{2nd},l}^{2}\\
 & \quad+\frac{\kappa^{3/2}\mu r}{ndp^{2}\sigma_{r}^{\star}}\zeta_{\mathsf{2nd},l}\left\Vert \bm{U}_{l,\cdot}^{\star}\bm{\Sigma}^{\star}\right\Vert _{2}+\frac{\kappa^{2}\mu r}{ndp^{2}}\zeta_{\mathsf{2nd},l}^{2}.
\end{align*}

\paragraph{Step 6: putting all pieces together.}

From the above steps, we can demonstrate that 
\begin{align*}
\left\Vert \bm{R}\bm{\Sigma}_{U,l}^{\star}\bm{R}^{\top}-\bm{\Sigma}_{U,l}\right\Vert  & \leq\sum_{i=1}^{4}\left\Vert \bm{A}_{i}-\bm{B}_{i}\right\Vert \leq\sum_{i=1}^{3}\left\Vert \bm{A}_{i}-\bm{B}_{i}\right\Vert +\alpha_{4}+\beta\\
 & \lesssim\delta\lambda_{\min}\left(\bm{\Sigma}_{U,l}^{\star}\right)+\underbrace{\frac{\sqrt{\kappa}}{np\sigma_{r}^{\star}}\zeta_{\mathsf{2nd},l}\left\Vert \bm{U}_{l,\cdot}^{\star}\bm{\Sigma}^{\star}\right\Vert _{2}}_{\eqqcolon\gamma_{1}}+\underbrace{\frac{\kappa}{np}\zeta_{\mathsf{2nd},l}^{2}}_{\eqqcolon\gamma_{2}}+\underbrace{\frac{1}{np}\zeta_{\mathsf{2nd},l}\frac{\theta}{\sqrt{\kappa}\sigma_{r}^{\star}}\omega_{l}^{\star}}_{\eqqcolon\gamma_{3}}\\
 & \quad+\underbrace{\frac{\kappa^{3/2}\mu r}{ndp^{2}\sigma_{r}^{\star}}\zeta_{\mathsf{2nd},l}\left\Vert \bm{U}_{l,\cdot}^{\star}\bm{\Sigma}^{\star}\right\Vert _{2}}_{\eqqcolon\gamma_{4}}+\underbrace{\frac{\kappa^{2}\mu r}{ndp^{2}}\zeta_{\mathsf{2nd},l}^{2}}_{\eqqcolon\gamma_{5}}\\
 & \quad+\underbrace{\frac{\sqrt{\kappa}}{np^{2}\sigma_{r}^{\star3}}\omega_{\max}^{2}\zeta_{\mathsf{2nd},l}\left\Vert \bm{U}_{l,\cdot}^{\star}\bm{\Sigma}^{\star}\right\Vert _{2}}_{\eqqcolon\gamma_{6}}+\underbrace{\frac{\kappa}{np^{2}\sigma_{r}^{\star2}}\omega_{\max}^{2}\zeta_{\mathsf{2nd},l}^{2}}_{\eqqcolon\gamma_{7}}.
\end{align*}
Note that under the assumption of Lemma \ref{lemma:pca-normal-approximation},
we have shown in Appendix \ref{appendix:proof-pca-normal-approximation}
that 
\begin{equation}
\zeta_{\mathsf{2nd},l}\lambda_{\min}^{-1/2}\left(\bm{\Sigma}_{U,l}^{\star}\right)\lesssim\frac{1}{r^{1/4}\log^{1/2}\left(n+d\right)}\ll1.\label{eq:pca-cov-est-inter-12}
\end{equation}
In addition, it follows from Lemma \ref{lemma:pca-covariance-concentration}
that 
\begin{equation}
\lambda_{\min}^{1/2}\left(\bm{\Sigma}_{U,l}^{\star}\right)\gtrsim\frac{1}{\sqrt{np}\sigma_{1}^{\star}}\left\Vert \bm{U}_{l,\cdot}^{\star}\bm{\Sigma}^{\star}\right\Vert _{2}+\frac{\omega_{l}^{\star}}{\sqrt{np}\sigma_{1}^{\star}}+\frac{1}{\sqrt{ndp^{2}\kappa}\sigma_{1}^{\star}}\left\Vert \bm{U}_{l,\cdot}^{\star}\bm{\Sigma}^{\star}\right\Vert _{2}.\label{eq:pca-cov-est-inter-13}
\end{equation}
As a consequence, we can derive the following upper bounds 
\begin{align*}
\gamma_{1} & \lesssim\frac{\sqrt{\kappa}}{np\sigma_{r}^{\star}}\lambda_{\min}^{1/2}\left(\bm{\Sigma}_{U,l}^{\star}\right)\left\Vert \bm{U}_{l,\cdot}^{\star}\bm{\Sigma}^{\star}\right\Vert _{2}\lesssim\frac{\kappa}{\sqrt{np}}\lambda_{\min}\left(\bm{\Sigma}_{U,l}^{\star}\right)\lesssim\delta\lambda_{\min}\left(\bm{\Sigma}_{U,l}^{\star}\right),\\
\gamma_{2} & \lesssim\frac{\kappa}{np}\lambda_{\min}\left(\bm{\Sigma}_{U,l}^{\star}\right)\lesssim\delta\lambda_{\min}\left(\bm{\Sigma}_{U,l}^{\star}\right),\\
\gamma_{3} & \lesssim\frac{1}{np}\lambda_{\min}^{1/2}\left(\bm{\Sigma}_{U,l}^{\star}\right)\frac{\theta}{\sqrt{\kappa}\sigma_{r}^{\star}}\omega_{l}^{\star}\lesssim\frac{\theta}{\sqrt{np}}\lambda_{\min}\left(\bm{\Sigma}_{U,l}^{\star}\right)\lesssim\delta\lambda_{\min}\left(\bm{\Sigma}_{U,l}^{\star}\right),\\
\gamma_{4} & \lesssim\frac{\kappa^{3/2}\mu r}{ndp^{2}\sigma_{r}^{\star}}\lambda_{\min}^{1/2}\left(\bm{\Sigma}_{U,l}^{\star}\right)\left\Vert \bm{U}_{l,\cdot}^{\star}\bm{\Sigma}^{\star}\right\Vert _{2}\lesssim\frac{\kappa^{5/2}\mu r}{\sqrt{ndp^{2}}}\lambda_{\min}\left(\bm{\Sigma}_{U,l}^{\star}\right)\lesssim\delta\lambda_{\min}\left(\bm{\Sigma}_{U,l}^{\star}\right),\\
\gamma_{5} & \lesssim\frac{\kappa^{2}\mu r}{ndp^{2}}\lambda_{\min}\left(\bm{\Sigma}_{U,l}^{\star}\right)\lesssim\delta\lambda_{\min}\left(\bm{\Sigma}_{U,l}^{\star}\right),
\end{align*}
provided that $\theta\lesssim1$, $np\gtrsim\delta^{-2}\kappa^{2}$
and $ndp^{2}\gtrsim\delta^{-2}\kappa^{5}\mu^{2}r^{2}$. In addition,
we also have
\begin{align*}
\gamma_{6} & \overset{\text{(i)}}{\lesssim}\frac{\kappa}{n^{2}p^{4}\sigma_{r}^{\star6}\delta}\omega_{\max}^{4}\left\Vert \bm{U}_{l,\cdot}^{\star}\bm{\Sigma}^{\star}\right\Vert _{2}^{2}+\delta\zeta_{\mathsf{2nd},l}^{2}\lesssim\frac{\kappa^{3}\kappa_{\omega}}{np^{2}\sigma_{r}^{\star2}\delta}\omega_{\max}^{2}\cdot\frac{1}{np^{2}\sigma_{1}^{\star4}}\omega_{\min}^{2}\left\Vert \bm{U}_{l,\cdot}^{\star}\bm{\Sigma}^{\star}\right\Vert _{2}^{2}+\delta\zeta_{\mathsf{2nd},l}^{2}\\
 & \lesssim\frac{\kappa^{3}\kappa_{\omega}}{\delta}\cdot\frac{1}{\sqrt{ndp^{2}}}\cdot\frac{\omega_{\max}^{2}}{p\sigma_{r}^{\star2}}\sqrt{\frac{d}{n}}\cdot\lambda_{\min}\left(\bm{\Sigma}_{U,l}^{\star}\right)+\delta\lambda_{\min}\left(\bm{\Sigma}_{U,l}^{\star}\right)\overset{\text{(ii)}}{\lesssim}\delta\lambda_{\min}\left(\bm{\Sigma}_{U,l}^{\star}\right)
\end{align*}
as well as
\[
\gamma_{7}\lesssim\frac{\kappa}{np^{2}\sigma_{r}^{\star2}}\omega_{\max}^{2}\lambda_{\min}\left(\bm{\Sigma}_{U,l}^{\star}\right)\lesssim\kappa\frac{1}{\sqrt{ndp^{2}}}\cdot\frac{\omega_{\max}^{2}}{p\sigma_{r}^{\star2}}\sqrt{\frac{d}{n}}\cdot\lambda_{\min}\left(\bm{\Sigma}_{U,l}^{\star}\right)\lesssim\delta\lambda_{\min}\left(\bm{\Sigma}_{U,l}^{\star}\right),
\]
where (i) utilizes the AM-GM inequality, while (ii) and (iii) hold
when $ndp^{2}\gtrsim\delta^{-2}\kappa^{6}\kappa_{\omega}^{2}$ and
\[
\frac{\omega_{\max}^{2}}{p\sigma_{r}^{\star2}}\sqrt{\frac{d}{n}}\lesssim\delta.
\]
These allow us to conclude that 
\[
\left\Vert \bm{R}\bm{\Sigma}_{U,l}^{\star}\bm{R}^{\top}-\bm{\Sigma}_{U,l}\right\Vert \lesssim\delta\lambda_{\min}\left(\bm{\Sigma}_{U,l}^{\star}\right),
\]
as long as the following assumptions hold: $np\gtrsim\delta^{-2}\kappa^{6}\mu^{2}r^{2}\kappa_{\omega}^{2}\log^{2}(n+d)$,
$ndp^{2}\gtrsim\delta^{-2}\kappa^{6}\kappa_{\omega}^{2}$, $n\gtrsim\delta^{-2}\kappa^{9}\mu^{2}r^{3}\log(n+d)$,
$n\gtrsim\kappa^{7}r\kappa_{\omega}^{2}\log(n+d)$, $d\gtrsim\kappa^{3}\mu r\log(n+d)$,
$\zeta_{\mathsf{1st}}/\sigma_{r}^{\star2}\lesssim\delta/(\kappa^{3}\mu r\kappa_{\omega}\sqrt{\log(n+d)})$,
$\theta\lesssim\delta/(\kappa^{3}\mu r\kappa_{\omega})$, $\max_{i\in[d]}\zeta_{\mathsf{2nd},i}\sqrt{d}\lesssim\delta/\sqrt{\kappa^{6}\mu r}$
and 
\[
\frac{\omega_{\max}^{2}}{p\sigma_{r}^{\star2}}\sqrt{\frac{d}{n}}\lesssim\delta.
\]
In what follows, we take a closer look at the last three assumptions. 
\begin{itemize}
\item In view of (\ref{subeq:pca-1st-equivalent}), it is readily seen that
$\zeta_{\mathsf{1st}}/\sigma_{r}^{\star2}\lesssim\delta/(\kappa^{3}\mu r\kappa_{\omega}\sqrt{\log(n+d)})$
can be guaranteed by 
\[
ndp^{2}\gtrsim\delta^{-2}\kappa^{8}\mu^{4}r^{4}\kappa_{\omega}^{2}\log^{5}\left(n+d\right),\qquad np\gtrsim\delta^{-2}\kappa^{8}\mu^{3}r^{3}\kappa_{\omega}^{2}\log^{3}\left(n+d\right),
\]
and 
\[
\frac{\omega_{\max}^{2}}{p\sigma_{r}^{\star2}}\sqrt{\frac{d}{n}}\lesssim\frac{\delta}{\kappa^{3}\mu r\kappa_{\omega}\log^{3/2}\left(n+d\right)},\qquad\frac{\omega_{\max}}{\sigma_{r}^{\star}}\sqrt{\frac{d}{np}}\lesssim\frac{\delta}{\kappa^{7/2}\mu r\kappa_{\omega}\log\left(n+d\right)}.
\]
\item By virtue of (\ref{eq:theta-zeta-1st}), it is straightforward to
see that $\theta\lesssim\delta/(\kappa^{3}\mu r\kappa_{\omega})$
can be guaranteed by $\zeta_{\mathsf{1st}}/\sigma_{r}^{\star2}\lesssim\delta/(\kappa^{3}\mu r\kappa_{\omega})$
--- the latter is already guaranteed by $\zeta_{\mathsf{1st}}/\sigma_{r}^{\star2}\lesssim\delta/(\kappa^{3}\mu r\kappa_{\omega}\sqrt{\log(n+d)})$. 
\item In addition, we know that $\max_{i\in[d]}\zeta_{\mathsf{2nd},i}\sqrt{d}\lesssim\delta/(\sqrt{\kappa^{6}\mu r})$
is equivalent to 
\[
\frac{\zeta_{\mathsf{1st}}}{\sigma_{r}^{\star2}}\frac{\sqrt{\kappa^{3}\mu^{2}r^{2}\log\left(n+d\right)}}{\sqrt{d}}+\frac{\zeta_{\mathsf{1st}}^{2}}{\sigma_{r}^{\star4}}\sqrt{\frac{\kappa^{3}\mu r\log\left(n+d\right)}{}}\lesssim\frac{\delta}{\sqrt{\kappa^{6}\mu r}},
\]
which can be guaranteed by $\zeta_{\mathsf{1st}}/\sigma_{r}^{\star2}\lesssim\delta/(\kappa^{3}\mu r\kappa_{\omega}\sqrt{\log(n+d)})$
as long as $d\gtrsim\kappa^{3}\mu r$. 
\end{itemize}
To summarize, we can conclude that the required assumptions for the
above results to hold are: $n\gtrsim\delta^{-2}\kappa^{9}\mu^{2}r^{3}\log(n+d)$,
$n\gtrsim\kappa^{7}r\kappa_{\omega}^{2}\log(n+d)$, $d\gtrsim\kappa^{3}\mu r\log(n+d)$,
\[
ndp^{2}\gtrsim\delta^{-2}\kappa^{8}\mu^{4}r^{4}\kappa_{\omega}^{2}\log^{5}\left(n+d\right),\qquad np\gtrsim\delta^{-2}\kappa^{8}\mu^{3}r^{3}\kappa_{\omega}^{2}\log^{3}\left(n+d\right),
\]
and 
\[
\frac{\omega_{\max}^{2}}{p\sigma_{r}^{\star2}}\sqrt{\frac{d}{n}}\lesssim\frac{\delta}{\kappa^{3}\mu r\kappa_{\omega}\log^{3/2}\left(n+d\right)},\qquad\frac{\omega_{\max}}{\sigma_{r}^{\star}}\sqrt{\frac{d}{np}}\lesssim\frac{\delta}{\kappa^{7/2}\mu r\kappa_{\omega}\log\left(n+d\right)}.
\]

\subsubsection{Proof of Lemma \ref{lemma:pca-cr-validity}\label{appendix:proof-lemma-pca-cr-validity}}

Recall the definition of the Euclidean ball $\mathcal{B}_{1-\alpha}$
in Algorithm \ref{alg:PCA-HeteroPCA-CR}. It is easily seen that 
\begin{align*}
\left(\bm{U}\bm{R}-\bm{U}^{\star}\right)_{l,\cdot}\left(\bm{\Sigma}_{U,l}^{\star}\right){}^{-1/2}\in\mathcal{B}_{1-\alpha}\quad & \Longleftrightarrow\quad\left(\bm{U}\bm{R}-\bm{U}^{\star}\right)_{l,\cdot}\left(\bm{\Sigma}_{U,l}^{\star}\right){}^{-1/2}\bm{R}\in\mathcal{B}_{1-\alpha}\\
 & \Longleftrightarrow\quad\left(\bm{U}-\bm{U}^{\star}\bm{R}^{\top}\right)_{l,\cdot}\bm{R}\left(\bm{\Sigma}_{U,l}^{\star}\right){}^{-1/2}\bm{R}^{\top}\in\mathcal{B}_{1-\alpha},
\end{align*}
where the last line comes from the rotational invariance of $\mathcal{B}_{1-\alpha}$.
From the definition of $\mathsf{CR}_{U,l}^{1-\alpha}$ in Algorithm
\ref{alg:PCA-HeteroPCA-CR}, we also know that 
\[
\bm{U}_{l,\cdot}^{\star}\bm{R}^{\top}\in\mathsf{CR}_{U,l}^{1-\alpha}\quad\Longleftrightarrow\quad\left(\bm{U}-\bm{U}^{\star}\bm{R}^{\top}\right)_{l,\cdot}\bm{\Sigma}_{U,l}^{-1/2}\in\mathcal{B}_{1-\alpha}.
\]
Let us define 
\[
\bm{\Delta}\coloneqq\left(\bm{U}-\bm{U}^{\star}\bm{R}^{\top}\right)_{l,\cdot}\bm{R}\left(\bm{\Sigma}_{U,l}^{\star}\right){}^{-1/2}\bm{R}^{\top}-\left(\bm{U}-\bm{U}^{\star}\bm{R}^{\top}\right)_{l,\cdot}\bm{\Sigma}_{U,l}^{-1/2},
\]
then it is straightforward to check that with probability exceeding
$1-O((n+d)^{-10})$ 
\begin{align}
\left\Vert \bm{\Delta}\right\Vert _{2} & \leq\left\Vert \left(\bm{U}-\bm{U}^{\star}\bm{R}^{\top}\right)_{l,\cdot}\right\Vert _{2}\left\Vert \bm{R}\left(\bm{\Sigma}_{U,l}^{\star}\right){}^{-1/2}\bm{R}^{\top}-\bm{\Sigma}_{U,l}^{-1/2}\right\Vert \nonumber \\
 & =\left\Vert \left(\bm{U}-\bm{U}^{\star}\bm{R}^{\top}\right)_{l,\cdot}\right\Vert _{2}\left\Vert \bm{R}\left(\bm{\Sigma}_{U,l}^{\star}\right){}^{-1/2}\bm{R}^{\top}\left(\bm{\Sigma}_{U,l}^{1/2}-\bm{R}\left(\bm{\Sigma}_{U,l}^{\star}\right){}^{1/2}\bm{R}^{\top}\right)\bm{\Sigma}_{U,l}^{-1/2}\right\Vert \nonumber \\
 & \leq\left\Vert \left(\bm{U}-\bm{U}^{\star}\bm{R}^{\top}\right)_{l,\cdot}\right\Vert _{2}\left\Vert \bm{R}\left(\bm{\Sigma}_{U,l}^{\star}\right){}^{-1/2}\bm{R}^{\top}\right\Vert \left\Vert \bm{R}\left(\bm{\Sigma}_{U,l}^{\star}\right){}^{1/2}\bm{R}^{\top}-\bm{\Sigma}_{U,l}^{1/2}\right\Vert \left\Vert \bm{\Sigma}_{U,l}^{-1/2}\right\Vert \nonumber \\
 & \lesssim\left\Vert \left(\bm{U}-\bm{U}^{\star}\bm{R}^{\top}\right)_{l,\cdot}\right\Vert _{2}\lambda_{\min}^{-1}\left(\bm{\Sigma}_{U,l}^{\star}\right)\left\Vert \bm{R}\left(\bm{\Sigma}_{U,l}^{\star}\right){}^{1/2}\bm{R}^{\top}-\bm{\Sigma}_{U,l}^{1/2}\right\Vert .\label{eq:proof-pca-cr-inter-2}
\end{align}
Here the last line follows from an immediate result from Lemma \ref{lemma:pca-covariance-estimation}
and Weyl's inequality: 
\begin{equation}
\lambda_{\min}\left(\bm{\Sigma}_{U,l}\right)\asymp\lambda_{\min}\left(\bm{\Sigma}_{U,l}^{\star}\right),\label{eq:proof-pca-cr-inter-1}
\end{equation}
which holds as long as $\delta\ll1$. Notice that 
\begin{align}
\left\Vert \bm{R}\left(\bm{\Sigma}_{l}^{\star}\right){}^{1/2}\bm{R}^{\top}-\bm{\Sigma}_{l}^{1/2}\right\Vert  & \overset{\text{(i)}}{\lesssim}\frac{1}{\lambda_{\min}^{1/2}\left(\bm{\Sigma}_{l}^{\star}\right)+\lambda_{\min}^{1/2}\left(\bm{\Sigma}_{l}\right)}\left\Vert \bm{R}\bm{\Sigma}_{l}^{\star}\bm{R}^{\top}-\bm{\Sigma}_{l}\right\Vert \nonumber \\
 & \overset{\text{(ii)}}{\lesssim}\lambda_{\min}^{-1/2}\left(\bm{\Sigma}_{l}^{\star}\right)\left\Vert \bm{R}\bm{\Sigma}_{l}^{\star}\bm{R}^{\top}-\bm{\Sigma}_{l}\right\Vert \nonumber \\
 & \overset{\text{(iii)}}{\lesssim}\lambda_{\min}^{1/2}\left(\bm{\Sigma}_{U,l}^{\star}\right)\delta,\label{eq:proof-pca-cr-inter-3}
\end{align}
where (i) follows from the perturbation bound of matrix square root
\citep[Lemma 2.1]{MR1176461}; (ii) arises from (\ref{eq:proof-pca-cr-inter-1});
and (iii) is a consequence of Lemma \ref{lemma:pca-covariance-estimation}.
We can combine (\ref{eq:proof-pca-cr-inter-2}) and (\ref{eq:proof-pca-cr-inter-3})
to achieve 
\begin{align*}
\left\Vert \bm{\Delta}\right\Vert _{2} & \lesssim\left\Vert \left(\bm{U}-\bm{U}^{\star}\bm{R}^{\top}\right)_{l,\cdot}\right\Vert _{2}\lambda_{\min}^{-1/2}\left(\bm{\Sigma}_{U,l}^{\star}\right)\delta\\
 & \lesssim\left[\frac{\theta}{\sqrt{\kappa}\sigma_{r}^{\star}}\left(\left\Vert \bm{U}_{l,\cdot}^{\star}\bm{\Sigma}^{\star}\right\Vert _{2}+\omega_{l}^{\star}\right)+\zeta_{\mathsf{2nd},l}\right]\lambda_{\min}^{-1/2}\left(\bm{\Sigma}_{U,l}^{\star}\right)\delta.
\end{align*}
Here, $\delta\in(0,1)$ is the (unspecified) quantity appearing in
Lemma \ref{lemma:pca-covariance-estimation} such that 
\[
\left\Vert \bm{R}\bm{\Sigma}_{U,l}^{\star}\bm{R}^{\top}-\bm{\Sigma}_{U,l}\right\Vert \lesssim\delta\lambda_{\min}\left(\bm{\Sigma}_{U,l}^{\star}\right),
\]
and the last line follows from (\ref{eq:pca-U-l-R-error}) in Lemma
\ref{lemma:pca-1st-err}. Let 
\[
\zeta\coloneqq\widetilde{C}\left[\frac{\theta}{\sqrt{\kappa}\sigma_{r}^{\star}}\left(\left\Vert \bm{U}_{l,\cdot}^{\star}\bm{\Sigma}^{\star}\right\Vert _{2}+\omega_{l}^{\star}\right)\lambda_{\min}^{-1/2}\left(\bm{\Sigma}_{U,l}^{\star}\right)\delta+\zeta_{\mathsf{2nd},l}\lambda_{\min}^{-1/2}\left(\bm{\Sigma}_{U,l}^{\star}\right)\delta\right]
\]
for some sufficiently large constant $\tilde{C}>0$ such that $\mathbb{P}(\Vert\bm{\Delta}\Vert_{2}\leq\zeta)\geq1-O((n+d)^{-10})$.
Recalling the definition \eqref{eq:defn-C-epsilon} of $\mathcal{C}^{\varepsilon}$
for any convex set $\mathcal{C}$, we have 
\begin{align}
 & \mathbb{P}\left(\bm{U}_{l,\cdot}^{\star}\bm{R}^{\top}\in\mathsf{CR}_{U,l}^{1-\alpha}\right)=\mathbb{P}\left(\left(\bm{U}-\bm{U}^{\star}\bm{R}^{\top}\right)_{l,\cdot}\bm{\Sigma}_{l}^{-1/2}\in\mathcal{B}_{1-\alpha}\right)\nonumber \\
 & \quad=\mathbb{P}\left(\left(\bm{U}-\bm{U}^{\star}\bm{R}^{\top}\right)_{l,\cdot}\bm{\Sigma}_{l}^{-1/2}\in\mathcal{B}_{1-\alpha},\left\Vert \bm{\Delta}\right\Vert _{2}\leq\zeta\right)+\mathbb{P}\left(\left(\bm{U}-\bm{U}^{\star}\bm{R}^{\top}\right)_{l,\cdot}\bm{\Sigma}_{l}^{-1/2}\in\mathcal{B}_{1-\alpha},\left\Vert \bm{\Delta}\right\Vert _{2}>\zeta\right)\nonumber \\
 & \quad\leq\mathbb{P}\left(\left(\bm{U}-\bm{U}^{\star}\bm{R}^{\top}\right)_{l,\cdot}\bm{R}\left(\bm{\Sigma}_{l}^{\star}\right){}^{-1/2}\bm{R}^{\top}\in\mathcal{B}_{1-\alpha}^{\zeta}\right)+\mathbb{P}\left(\left\Vert \bm{\Delta}\right\Vert _{2}>\zeta\right)\nonumber \\
 & \quad\overset{\text{(i)}}{=}\mathbb{P}\left(\left(\bm{U}\bm{R}-\bm{U}^{\star}\right)_{l,\cdot}\left(\bm{\Sigma}_{l}^{\star}\right){}^{-1/2}\in\mathcal{B}_{1-\alpha}^{\zeta}\right)+O\left(\left(n+d\right)^{-10}\right)\nonumber \\
 & \quad\overset{\text{(ii)}}{\leq}\mathcal{N}\left(\bm{0},\bm{I}_{r}\right)\left\{ \mathcal{B}_{1-\alpha}^{\zeta}\right\} +O\left(\log^{-1/2}\left(n+d\right)\right)\nonumber \\
 & \quad\overset{\text{(iii)}}{\leq}\mathcal{N}\left(\bm{0},\bm{I}_{r}\right)\left\{ \mathcal{B}_{1-\alpha}\right\} +\zeta\left(0.59r^{1/4}+0.21\right)+O\left(\log^{-1/2}\left(n+d\right)\right)\nonumber \\
 & \quad\overset{\text{(iv)}}{\leq}1-\alpha+O\left(\log^{-1/2}\left(n+d\right)\right).\label{eq:pca-cr-inter-4}
\end{align}
Here (i) holds since $\mathcal{B}_{1-\alpha}^{\zeta}$ is rotational
invariant; (ii) uses Lemma \ref{lemma:pca-normal-approximation};
(iii) invokes Theorem \ref{thm:gaussian-perimeter}; and (iv) makes
use of the definition of $\mathcal{B}_{1-\alpha}$ in Algorithm \ref{alg:PCA-HeteroPCA-CR}
and holds under the condition $\zeta r^{1/4}\lesssim1/\sqrt{\log(n+d)}$.
Similar to (\ref{eq:pca-cr-inter-4})we can show that 
\begin{equation}
\mathbb{P}\left(\bm{U}_{l,\cdot}^{\star}\bm{R}^{\top}\in\mathsf{CR}_{U,l}^{1-\alpha}\right)\geq1-\alpha+O\left(\log^{-1/2}\left(n+d\right)\right)\label{eq:pca-cr-inter-5}
\end{equation}
Taking (\ref{eq:pca-cr-inter-4}) and (\ref{eq:pca-cr-inter-5}) collectively
yields 
\[
\mathbb{P}\left(\bm{U}_{l,\cdot}^{\star}\bm{R}^{\top}\in\mathsf{CR}_{U,l}^{1-\alpha}\right)=1-\alpha+O\left(\log^{-1/2}\left(n+d\right)\right),
\]
provided that $\zeta r^{1/4}\lesssim1/\sqrt{\log(n+d)}$, or equivalently,
\begin{equation}
\frac{\theta}{\sqrt{\kappa}\sigma_{r}^{\star}}\left(\left\Vert \bm{U}_{l,\cdot}^{\star}\bm{\Sigma}^{\star}\right\Vert _{2}+\omega_{l}^{\star}\right)\lambda_{\min}^{-1/2}\left(\bm{\Sigma}_{U,l}^{\star}\right)\delta+\zeta_{\mathsf{2nd},l}\lambda_{\min}^{-1/2}\left(\bm{\Sigma}_{U,l}^{\star}\right)\delta\lesssim\frac{1}{r^{1/4}\sqrt{\log\left(n+d\right)}}.\label{eq:condition-theta-delta-1579}
\end{equation}

It remains to choose $\delta\in(0,1)$ to satisfy the above condition.
First, we have learned from the proof of Lemma \ref{lemma:pca-normal-approximation}
(more specifically, Step 3 in Appendix \ref{appendix:proof-pca-normal-approximation})
that 
\[
\zeta_{\mathsf{2nd},l}\lambda_{\min}^{-1/2}(\bm{\Sigma}_{l}^{\star})\lesssim1/(r^{1/4}\log^{1/2}(n+d))
\]
holds under the conditions of Lemma \ref{lemma:pca-normal-approximation}.
As a result, it is sufficient to verify that 
\begin{equation}
\frac{\theta}{\sqrt{\kappa}\sigma_{r}^{\star}}\left(\left\Vert \bm{U}_{l,\cdot}^{\star}\bm{\Sigma}^{\star}\right\Vert _{2}+\omega_{l}^{\star}\right)\lambda_{\min}^{-1/2}\left(\bm{\Sigma}_{U,l}^{\star}\right)\delta\lesssim\frac{1}{r^{1/4}\sqrt{\log\left(n+d\right)}}.\label{eq:condition-theta-delta-1579-1}
\end{equation}
To this end, recall from Lemma \ref{lemma:pca-covariance-concentration}
that 
\[
\lambda_{\min}^{1/2}\left(\bm{\Sigma}_{U,l}^{\star}\right)\gtrsim\left(\frac{1}{\sqrt{np}\sigma_{1}^{\star}}+\frac{1}{\sqrt{ndp^{2}\kappa}\sigma_{1}^{\star}}\right)\left(\left\Vert \bm{U}_{l,\cdot}^{\star}\bm{\Sigma}^{\star}\right\Vert _{2}+\omega_{l}^{\star}\right)+\frac{1}{\sqrt{np^{2}}\sigma_{1}^{\star2}}\omega_{\min}\left\Vert \bm{U}_{l,\cdot}^{\star}\bm{\Sigma}^{\star}\right\Vert _{2}+\frac{\omega_{\min}\omega_{l}^{\star}}{\sqrt{np^{2}}\sigma_{1}^{\star2}}.
\]
In view of the definition of $\theta$ (cf.~(\ref{eq:theta-definition})),
we can show that 
\begin{align*}
\frac{\theta}{\sqrt{\kappa}\sigma_{r}^{\star}}\left(\left\Vert \bm{U}_{l,\cdot}^{\star}\bm{\Sigma}^{\star}\right\Vert _{2}+\omega_{l}^{\star}\right) & \asymp\underbrace{\frac{1}{\sigma_{r}^{\star}}\sqrt{\frac{r\log^{2}\left(n+d\right)}{np}}\left(\left\Vert \bm{U}_{l,\cdot}^{\star}\bm{\Sigma}^{\star}\right\Vert _{2}+\omega_{l}^{\star}\right)}_{\eqqcolon\alpha_{1}}+\underbrace{\frac{1}{\sigma_{r}^{\star}}\sqrt{\frac{\kappa\mu r^{2}\log^{3}\left(n+d\right)}{ndp^{2}}}\left(\left\Vert \bm{U}_{l,\cdot}^{\star}\bm{\Sigma}^{\star}\right\Vert _{2}+\omega_{l}^{\star}\right)}_{\eqqcolon\alpha_{2}}\\
 & \quad+\underbrace{\frac{\omega_{\max}}{\sigma_{r}^{\star2}}\sqrt{\frac{r\log^{2}\left(n+d\right)}{np^{2}}}\left\Vert \bm{U}_{l,\cdot}^{\star}\bm{\Sigma}^{\star}\right\Vert _{2}}_{\eqqcolon\alpha_{3}}+\underbrace{\frac{\omega_{\max}\omega_{l}^{\star}}{\sigma_{r}^{\star2}}\sqrt{\frac{r\log^{2}\left(n+d\right)}{np^{2}}}}_{\eqqcolon\alpha_{4}}\\
 & \lesssim\sqrt{\kappa^{3}\mu r^{2}\kappa_{\omega}\log^{3}\left(n+d\right)}\lambda_{\min}^{1/2}\left(\bm{\Sigma}_{U,l}^{\star}\right),
\end{align*}
where the last line holds since 
\begin{align*}
\alpha_{1} & \lesssim\sqrt{\kappa r\log^{2}\left(n+d\right)}\frac{1}{\sqrt{np}\sigma_{1}^{\star}}\left(\left\Vert \bm{U}_{l,\cdot}^{\star}\bm{\Sigma}^{\star}\right\Vert _{2}+\omega_{l}^{\star}\right)\lesssim\sqrt{\kappa r\log^{2}\left(n+d\right)}\lambda_{\min}^{1/2}\left(\bm{\Sigma}_{U,l}^{\star}\right),\\
\alpha_{2} & \lesssim\sqrt{\kappa^{3}\mu r^{2}\log^{3}\left(n+d\right)}\frac{1}{\sqrt{ndp^{2}\kappa}\sigma_{1}^{\star}}\left(\left\Vert \bm{U}_{l,\cdot}^{\star}\bm{\Sigma}^{\star}\right\Vert _{2}+\omega_{l}^{\star}\right)\lesssim\sqrt{\kappa^{3}\mu r^{2}\log^{3}\left(n+d\right)}\lambda_{\min}^{1/2}\left(\bm{\Sigma}_{U,l}^{\star}\right),\\
\alpha_{3} & \lesssim\sqrt{\kappa^{2}r\kappa_{\omega}\log^{2}\left(n+d\right)}\frac{1}{\sqrt{np^{2}}\sigma_{1}^{\star2}}\omega_{\min}\left\Vert \bm{U}_{l,\cdot}^{\star}\bm{\Sigma}^{\star}\right\Vert _{2}\lesssim\sqrt{\kappa^{2}r\kappa_{\omega}\log^{2}\left(n+d\right)}\lambda_{\min}^{1/2}\left(\bm{\Sigma}_{U,l}^{\star}\right),\\
\alpha_{4} & \lesssim\sqrt{\kappa^{2}r\kappa_{\omega}\log^{2}\left(n+d\right)}\frac{\omega_{\min}\omega_{l}^{\star}}{\sqrt{np^{2}}\sigma_{1}^{\star2}}\lesssim\sqrt{\kappa^{2}r\kappa_{\omega}\log^{2}\left(n+d\right)}\lambda_{\min}^{1/2}\left(\bm{\Sigma}_{U,l}^{\star}\right).
\end{align*}
The above bound immediately gives 
\[
\frac{\theta}{\sqrt{\kappa}\sigma_{r}^{\star}}\left(\left\Vert \bm{U}_{l,\cdot}^{\star}\bm{\Sigma}^{\star}\right\Vert _{2}+\omega_{l}^{\star}\right)\lambda_{\min}^{-1/2}\left(\bm{\Sigma}_{U,l}^{\star}\right)\delta\lesssim\sqrt{\kappa^{3}\mu r^{2}\kappa_{\omega}\log^{3}\left(n+d\right)}\delta.
\]
Therefore, the desired condition \eqref{eq:condition-theta-delta-1579-1}
--- and hence $\zeta r^{1/4}\lesssim1/\sqrt{\log(n+d)}$ --- can
be guaranteed by taking 
\[
\delta=\frac{1}{\kappa^{3/2}\mu^{1/2}r^{5/4}\kappa_{\omega}^{1/2}\log^{2}\left(n+d\right)}.
\]

After $\delta$ is specified, we can easily check that under our choice
of $\delta$, the assumptions in Lemma \ref{lemma:pca-covariance-estimation}
read $n\gtrsim\kappa^{12}\mu^{3}r^{11/2}\kappa_{\omega}\log^{5}(n+d)$,
$d\gtrsim\kappa^{3}\mu r\log(n+d)$, 
\[
ndp^{2}\gtrsim\kappa^{11}\mu^{5}r^{13/2}\kappa_{\omega}^{3}\log^{9}\left(n+d\right),\qquad np\gtrsim\kappa^{11}\mu^{4}r^{11/2}\kappa_{\omega}^{3}\log^{7}\left(n+d\right),
\]
and 
\[
\frac{\omega_{\max}^{2}}{p\sigma_{r}^{\star2}}\sqrt{\frac{d}{n}}\lesssim\frac{1}{\kappa^{9/2}\mu^{3/2}r^{9/4}\kappa_{\omega}^{3/2}\log^{7/2}\left(n+d\right)},\qquad\frac{\omega_{\max}}{\sigma_{r}^{\star}}\sqrt{\frac{d}{np}}\lesssim\frac{1}{\kappa^{5}\mu^{3/2}r^{9/4}\kappa_{\omega}^{3/2}\log^{3}\left(n+d\right)}.
\]
This concludes the proof.

\subsection{Auxiliary lemmas for Theorem \ref{thm:ce-complete}}

\subsubsection{Proof of Lemma \ref{lemma:ce-2nd-error}\label{appendix:proof-ce-2nd-error}}

Throughout this section, all the probabilistic arguments are conditional
on $\bm{F}$, and we shall always assume that the $\sigma(\bm{F})$-measurable
high-probability event $\mathcal{E}_{\mathsf{good}}$ occurs. Following
the same analysis as in Appendix \ref{appendix:proof-pca-2nd-error}
(Proof of Lemma \ref{lemma:pca-2nd-error}), we can obtain 
\[
\bm{U}\bm{R}_{\bm{U}}-\bm{U}^{\natural}=\bm{Z}+\bm{\Psi}_{\bm{U}}
\]
with $\bm{R}_{\bm{U}}=\arg\min_{\bm{O}\in\mathcal{O}^{r\times r}}\|\bm{U}\bm{O}-\bm{U}^{\natural}\|_{\mathrm{F}}^{2}$,
where 
\[
\bm{Z}=\left[\bm{E}\bm{M}^{\natural\top}+\mathcal{P}_{\mathsf{off}\text{-}\mathsf{diag}}\left(\bm{E}\bm{E}^{\top}\right)\right]\bm{U}^{\natural}\left(\bm{\Sigma}^{\natural}\right)^{-2}
\]
and $\bm{\Psi}_{\bm{U}}$ is a residual matrix obeying 
\begin{equation}
\mathbb{P}\left(\left\Vert \bm{e}_{l}^{\top}\bm{\Psi}_{\bm{U}}\right\Vert _{2}\ind_{\mathcal{E}_{\mathsf{good}}}\lesssim\zeta_{\mathsf{2nd},l}\,\big|\,\bm{F}\right)\geq1-O\left(\left(n+d\right)^{-10}\right)\label{eq:ce-2nd-inter-1}
\end{equation}
for all $l\in[d]$. Here, $\zeta_{\mathsf{2nd},l}$ is a quantity
defined in Lemma \ref{lemma:pca-2nd-error}. From the proof of Lemma
\ref{lemma:pca-1st-err}, we know that 
\[
\bm{R}_{\bm{U}}^{\top}\bm{\Sigma}^{2}\bm{R}_{\bm{U}}=\bm{\Sigma}^{\natural2}+\bm{\Psi}_{\bm{\Sigma}}
\]
holds for some matrix $\bm{\Psi}_{\bm{\Sigma}}$ satisfying 
\begin{equation}
\mathbb{P}\left(\left\Vert \bm{\Psi}_{\bm{\Sigma}}\right\Vert \ind_{\mathcal{E}_{\mathsf{good}}}\lesssim\sqrt{\frac{\kappa^{3}\mu r\log\left(n+d\right)}{d}}\zeta_{\mathsf{1st}}+\kappa\frac{\zeta_{\mathsf{1st}}^{2}}{\sigma_{r}^{\star2}}\,\big|\,\bm{F}\right)\geq1-O\left(\left(n+d\right)^{-10}\right).\label{eq:ce-2nd-inter-2}
\end{equation}
Armed with the above facts, we can write 
\begin{align*}
\bm{S}-\bm{M}^{\natural}\bm{M}^{\natural\top} & =\bm{U}\bm{\Sigma}^{2}\bm{U}^{\top}-\bm{U}^{\natural}\bm{\Sigma}^{\natural2}\bm{U}^{\natural\top}\\
 & =\bm{U}\bm{R}_{\bm{U}}\bm{R}_{\bm{U}}^{\top}\bm{\Sigma}^{2}\bm{R}_{\bm{U}}\bm{R}_{\bm{U}}^{\top}\bm{U}^{\top}-\bm{U}^{\natural}\bm{\Sigma}^{\natural2}\bm{U}^{\natural\top}\\
 & =\bm{U}\bm{R}_{\bm{U}}\bm{\Sigma}^{\natural2}\bm{R}_{\bm{U}}^{\top}\bm{U}^{\top}+\bm{U}\bm{R}_{\bm{U}}\bm{\Psi}_{\bm{\Sigma}}\bm{R}_{\bm{U}}^{\top}\bm{U}^{\top}-\bm{U}^{\natural}\bm{\Sigma}^{\natural2}\bm{U}^{\natural\top}\\
 & =\left(\bm{U}^{\natural}+\bm{Z}+\bm{\Psi}_{\bm{U}}\right)\bm{\Sigma}^{\natural2}\left(\bm{U}^{\natural}+\bm{Z}+\bm{\Psi}_{\bm{U}}\right)^{\top}+\bm{U}\bm{R}_{\bm{U}}\bm{\Psi}_{\bm{\Sigma}}\bm{R}_{\bm{U}}^{\top}\bm{U}^{\top}-\bm{U}^{\natural}\bm{\Sigma}^{\natural2}\bm{U}^{\natural\top}\\
 & =\bm{X}+\bm{\Phi},
\end{align*}
where 
\begin{align*}
\bm{X} & =\bm{U}^{\natural}\bm{\Sigma}^{\natural2}\bm{Z}^{\top}+\bm{Z}\bm{\Sigma}^{\natural2}\bm{U}^{\natural}\\
 & =\bm{E}\bm{M}^{\natural\top}+\bm{M}^{\natural}\bm{E}^{\top}+\mathcal{P}_{\mathsf{off}\text{-}\mathsf{diag}}\left(\bm{E}\bm{E}^{\top}\right)\bm{U}^{\natural}\bm{U}^{\natural\top}+\bm{U}^{\natural}\bm{U}^{\natural\top}\mathcal{P}_{\mathsf{off}\text{-}\mathsf{diag}}\left(\bm{E}\bm{E}^{\top}\right)
\end{align*}
and 
\[
\bm{\Phi}=\bm{U}^{\natural}\bm{\Sigma}^{\natural2}\bm{\Psi}_{\bm{U}}^{\top}+\bm{Z}\bm{\Sigma}^{\natural2}\left(\bm{Z}+\bm{\Psi}_{\bm{U}}\right)^{\top}+\bm{\Psi}_{\bm{U}}\bm{\Sigma}^{\natural2}\left(\bm{U}\bm{R}_{\bm{U}}\right)^{\top}+\bm{U}\bm{R}_{\bm{U}}\bm{\Psi}_{\bm{\Sigma}}\left(\bm{U}\bm{R}_{\bm{U}}\right)^{\top}.
\]

It has already been shown in Lemma \ref{lemma:pca-1st-err} that for
each $l\in[d]$, 
\begin{equation}
\left\Vert \bm{Z}_{l,\cdot}\right\Vert _{2}\lesssim\frac{\theta}{\sqrt{\kappa}\sigma_{r}^{\star}}\left(\left\Vert \bm{U}_{l,\cdot}^{\star}\bm{\Sigma}^{\star}\right\Vert _{2}+\omega_{l}^{\star}\right)\label{eq:ce-2nd-inter-3}
\end{equation}
and 
\[
\left\Vert \left(\bm{U}\bm{R}-\bm{U}^{\star}\right)_{l,\cdot}\right\Vert _{2}\lesssim\frac{\theta}{\sqrt{\kappa}\sigma_{r}^{\star}}\left(\left\Vert \bm{U}_{l,\cdot}^{\star}\bm{\Sigma}^{\star}\right\Vert _{2}+\omega_{l}^{\star}\right)+\zeta_{\mathsf{2nd},l}
\]
hold with probability at least $1-O\big((n+d)^{-10}\big)$. As an
immediate consequence, 
\begin{equation}
\left\Vert \bm{U}_{l,\cdot}\right\Vert _{2}\leq\left\Vert \bm{U}_{l,\cdot}^{\star}\right\Vert _{2}+\left\Vert \left(\bm{U}\bm{R}-\bm{U}^{\star}\right)_{l,\cdot}\right\Vert _{2}\lesssim\left\Vert \bm{U}_{l,\cdot}^{\star}\right\Vert _{2}+\theta\frac{\omega_{l}^{\star}}{\sigma_{1}^{\star}}+\zeta_{\mathsf{2nd},l},\label{eq:ce-2nd-inter-4}
\end{equation}
provided that $\theta\lesssim1$. Then with probability exceeding
$1-O((n+d)^{-10})$, we have 
\begin{align*}
\left|\Phi_{i,j}\right| & \overset{\text{(i)}}{\lesssim}\sigma_{1}^{\star2}\left(\left\Vert \bm{U}_{i,\cdot}^{\star}\right\Vert _{2}\left\Vert \bm{e}_{j}^{\top}\bm{\Psi}_{\bm{U}}\right\Vert _{2}+\left\Vert \bm{Z}_{i,\cdot}\right\Vert _{2}\left\Vert \bm{Z}_{j,\cdot}\right\Vert _{2}+\left\Vert \bm{Z}_{i,\cdot}\right\Vert _{2}\left\Vert \bm{e}_{j}^{\top}\bm{\Psi}_{\bm{U}}\right\Vert _{2}+\left\Vert \bm{e}_{i}^{\top}\bm{\Psi}_{\bm{U}}\right\Vert _{2}\left\Vert \bm{U}_{j,\cdot}\right\Vert _{2}\right)\\
 & \quad+\left\Vert \bm{\Psi}_{\bm{\Sigma}}\right\Vert \left\Vert \bm{U}_{i,\cdot}\right\Vert _{2}\left\Vert \bm{U}_{j,\cdot}\right\Vert _{2}\\
 & \overset{\text{(ii)}}{\lesssim}\zeta_{\mathsf{2nd},j}\sigma_{1}^{\star2}\left(\left\Vert \bm{U}_{i,\cdot}^{\star}\right\Vert _{2}+\theta\frac{\omega_{i}^{\star}}{\sigma_{1}^{\star}}\right)+\sigma_{1}^{\star2}\zeta_{\mathsf{2nd},i}\left(\left\Vert \bm{U}_{j,\cdot}^{\star}\right\Vert _{2}+\theta\frac{\omega_{j}^{\star}}{\sigma_{1}^{\star}}+\zeta_{\mathsf{2nd},j}\right)\\
 & \quad+\theta^{2}\left(\left\Vert \bm{U}_{i,\cdot}^{\star}\bm{\Sigma}^{\star}\right\Vert _{2}+\omega_{i}^{\star}\right)\left(\left\Vert \bm{U}_{j,\cdot}^{\star}\bm{\Sigma}^{\star}\right\Vert _{2}+\omega_{j}^{\star}\right)\\
 & \quad+\left(\sqrt{\frac{\kappa^{3}\mu r\log\left(n+d\right)}{d}}\zeta_{\mathsf{1st}}+\kappa\frac{\zeta_{\mathsf{1st}}^{2}}{\sigma_{r}^{\star2}}\right)\left(\left\Vert \bm{U}_{i,\cdot}^{\star}\right\Vert _{2}+\theta\frac{\omega_{i}^{\star}}{\sigma_{1}^{\star}}+\zeta_{\mathsf{2nd},i}\right)\left(\left\Vert \bm{U}_{j,\cdot}^{\star}\right\Vert _{2}+\theta\frac{\omega_{j}^{\star}}{\sigma_{1}^{\star}}+\zeta_{\mathsf{2nd},j}\right)\\
 & \overset{\text{(iii)}}{\lesssim}\theta^{2}\left(\left\Vert \bm{U}_{i,\cdot}^{\star}\bm{\Sigma}^{\star}\right\Vert _{2}+\omega_{i}^{\star}\right)\left(\left\Vert \bm{U}_{j,\cdot}^{\star}\bm{\Sigma}^{\star}\right\Vert _{2}+\omega_{j}^{\star}\right)+\sigma_{1}^{\star2}\zeta_{\mathsf{2nd},i}\zeta_{\mathsf{2nd},j}\\
 & \quad+\sigma_{1}^{\star2}\zeta_{\mathsf{2nd},i}\left(\left\Vert \bm{U}_{j,\cdot}^{\star}\right\Vert _{2}+\theta\frac{\omega_{j}^{\star}}{\sigma_{1}^{\star}}\right)+\zeta_{\mathsf{2nd},j}\sigma_{1}^{\star2}\left(\left\Vert \bm{U}_{i,\cdot}^{\star}\right\Vert _{2}+\theta\frac{\omega_{i}^{\star}}{\sigma_{1}^{\star}}\right).
\end{align*}
for each $i,j\in[d]$. Here, (i) makes use of (\ref{eq:good-event-sigma-least-largest-hpca})
and the fact that $\bm{U}^{\natural}=\bm{U}^{\star}\bm{Q}$ for some
orthonormal matrix $\bm{Q}$ (cf.~\eqref{eq:defn-Q-denoising-PCA});
(ii) follows from (\ref{eq:ce-2nd-inter-1}), (\ref{eq:ce-2nd-inter-2}),
(\ref{eq:ce-2nd-inter-3}) and (\ref{eq:ce-2nd-inter-4}), and holds
provided that $\theta\lesssim1$; (iii) holds provided that $\zeta_{\mathsf{1st}}/\sigma_{r}^{\star2}\lesssim1$
and $d\gtrsim\kappa\mu r\log(n+d)$.

\subsubsection{Proof of Lemma \ref{lemma:ce-variance-concentration}\label{appendix:proof-lemma-ce-variance-concentration}}

In this subsection, we shall focus on establishing the claimed result
for the case when $i\neq j$; the case when $i=j$ can be proved in
a similar (in fact, easier) manner. In view of the expression \eqref{eq:var-Xij-F-hpca},
we can write 
\[
\mathsf{var}\left(X_{i,j}|\bm{F}\right)=\underbrace{\sum_{l=1}^{n}M_{j,l}^{\natural2}\sigma_{i,l}^{2}}_{\eqqcolon\alpha_{1}}+\underbrace{\sum_{l=1}^{n}M_{i,l}^{\natural2}\sigma_{j,l}^{2}}_{\eqqcolon\alpha_{2}}+\underbrace{\sum_{l=1}^{n}\sum_{k:k\neq i}\sigma_{i,l}^{2}\sigma_{k,l}^{2}\left(\bm{U}_{k,\cdot}^{\star}\bm{U}_{j,\cdot}^{\star\top}\right)^{2}}_{\eqqcolon\alpha_{3}}+\underbrace{\sum_{l=1}^{n}\sum_{k:k\neq j}\sigma_{j,l}^{2}\sigma_{k,l}^{2}\left(\bm{U}_{k,\cdot}^{\star}\bm{U}_{i,\cdot}^{\star\top}\right)^{2}}_{\eqqcolon\alpha_{4}},
\]
thus motivating us to study the behavior of $\alpha_{1}$, $\alpha_{2}$,
$\alpha_{3}$ and $\alpha_{4}$ respectively. 
\begin{itemize}
\item Let us begin with the term $\alpha_{1}$. By virtue of \eqref{eq:defn-M-Mnatural-E-denoising-PCA}
and \eqref{eq:sigma-ij-square-denoising-PCA}, we have 
\begin{align*}
\alpha_{1} & =\sum_{l=1}^{n}M_{j,l}^{\natural2}\sigma_{i,l}^{2}=\sum_{l=1}^{n}\left(\frac{1}{\sqrt{n}}\bm{U}_{j,\cdot}^{\star}\bm{\Sigma}^{\star}\bm{f}_{l}\right)^{2}\left[\frac{1-p}{np}\left(\bm{U}_{i,\cdot}^{\star}\bm{\Sigma}^{\star}\bm{f}_{l}\right)^{2}+\frac{\omega_{i}^{\star2}}{np}\right]\\
 & =\frac{1-p}{n^{2}p}\sum_{l=1}^{n}\left(\bm{U}_{j,\cdot}^{\star}\bm{\Sigma}^{\star}\bm{f}_{l}\right)^{2}\left(\bm{U}_{i,\cdot}^{\star}\bm{\Sigma}^{\star}\bm{f}_{l}\right)^{2}+\frac{\omega_{i}^{\star2}}{n^{2}p}\sum_{l=1}^{n}\left(\bm{U}_{j,\cdot}^{\star}\bm{\Sigma}^{\star}\bm{f}_{l}\right)^{2}.
\end{align*}
On the event $\mathcal{E}_{\mathsf{good}}$ (see Lemma~\ref{lemma:useful-property-good-event-hpca}),
we know from the basic facts in (\ref{eq:good-event-concentration-2-hpca})
and (\ref{eq:good-event-concentration-3-hpca}) that 
\begin{align*}
\left|\frac{1}{n}\sum_{l=1}^{n}\left(\bm{U}_{j,\cdot}^{\star}\bm{\Sigma}^{\star}\bm{f}_{l}\right)^{2}\left(\bm{U}_{i,\cdot}^{\star}\bm{\Sigma}^{\star}\bm{f}_{l}\right)^{2}-S_{i,i}^{\star}S_{j,j}^{\star}-2S_{i,j}^{\star2}\right| & \lesssim\sqrt{\frac{\log^{3}\left(n+d\right)}{n}}S_{i,i}^{\star}S_{j,j}^{\star};\\
\left|\frac{1}{n}\sum_{l=1}^{n}\left(\bm{U}_{j,\cdot}^{\star}\bm{\Sigma}^{\star}\bm{f}_{l}\right)^{2}-S_{j,j}^{\star}\right| & \lesssim\sqrt{\frac{\log\left(n+d\right)}{n}}S_{j,j}^{\star}.
\end{align*}
As a result, we can express $\alpha_{1}$ as 
\begin{align*}
\alpha_{1} & =\underbrace{\frac{1-p}{np}\left(S_{i,i}^{\star}S_{j,j}^{\star}+2S_{i,j}^{\star2}\right)+\frac{\omega_{i}^{\star2}}{np}S_{j,j}^{\star}}_{\eqqcolon\alpha_{1}^{\star}}+r_{1}
\end{align*}
for some residual term $r_{1}$ satisfying 
\begin{align*}
\left|r_{1}\right| & \lesssim\sqrt{\frac{\log^{3}\left(n+d\right)}{n}}\alpha_{1}^{\star}.
\end{align*}
\item Akin to our analysis for $\alpha_{1}$, we can also demonstrate that
\[
\alpha_{2}=\underbrace{\frac{1-p}{np}\left[S_{i,i}^{\star}S_{j,j}^{\star}+2S_{i,j}^{\star2}\right]+\frac{\omega_{j}^{\star2}}{np}S_{i,i}^{\star}}_{\eqqcolon\alpha_{2}^{\star}}+r_{2}
\]
for some residual term $r_{2}$ obeying 
\begin{align*}
\left|r_{2}\right| & \lesssim\sqrt{\frac{\log^{3}\left(n+d\right)}{n}}\alpha_{2}^{\star}.
\end{align*}
\item When it comes to $\alpha_{3}$, we make the observation that 
\begin{align*}
\alpha_{3} & =\sum_{l=1}^{n}\sum_{k:k\neq i}\left[\frac{1-p}{np}\left(\bm{U}_{i,\cdot}^{\star}\bm{\Sigma}^{\star}\bm{f}_{l}\right)^{2}+\frac{\omega_{i}^{\star2}}{np}\right]\left[\frac{1-p}{np}\left(\bm{U}_{k,\cdot}^{\star}\bm{\Sigma}^{\star}\bm{f}_{l}\right)^{2}+\frac{\omega_{k}^{\star2}}{np}\right]\left(\bm{U}_{k,\cdot}^{\star}\bm{U}_{j,\cdot}^{\star\top}\right)^{2}\\
 & =\left(\frac{1-p}{np}\right)^{2}\sum_{k:k\neq i}\left[\sum_{l=1}^{n}\left(\bm{U}_{i,\cdot}^{\star}\bm{\Sigma}^{\star}\bm{f}_{l}\right)^{2}\left(\bm{U}_{k,\cdot}^{\star}\bm{\Sigma}^{\star}\bm{f}_{l}\right)^{2}\right]\left(\bm{U}_{k,\cdot}^{\star}\bm{U}_{j,\cdot}^{\star\top}\right)^{2}\\
 & \quad\quad+\frac{1-p}{n^{2}p^{2}}\sum_{k:k\neq i}\sum_{l=1}^{n}\left(\bm{U}_{i,\cdot}^{\star}\bm{\Sigma}^{\star}\bm{f}_{l}\right)^{2}\omega_{k}^{\star2}\left(\bm{U}_{k,\cdot}^{\star}\bm{U}_{j,\cdot}^{\star\top}\right)^{2}\\
 & \quad\quad+\frac{\left(1-p\right)\omega_{i}^{\star2}}{n^{2}p^{2}}\sum_{k:k\neq i}\sum_{l=1}^{n}\left(\bm{U}_{k,\cdot}^{\star}\bm{\Sigma}^{\star}\bm{f}_{l}\right)^{2}\left(\bm{U}_{k,\cdot}^{\star}\bm{U}_{j,\cdot}^{\star\top}\right)^{2}+\frac{\omega_{i}^{\star2}}{np^{2}}\sum_{k:k\neq i}\omega_{k}^{\star2}\left(\bm{U}_{k,\cdot}^{\star}\bm{U}_{j,\cdot}^{\star\top}\right)^{2}.
\end{align*}
On the event $\mathcal{E}_{\mathsf{good}}$, it is seen from (\ref{eq:good-event-concentration-2-hpca})
and (\ref{eq:good-event-concentration-3-hpca}) that 
\begin{align*}
\left|\frac{1}{n}\sum_{l=1}^{n}\left(\bm{U}_{i,\cdot}^{\star}\bm{\Sigma}^{\star}\bm{f}_{l}\right)^{2}\left(\bm{U}_{k,\cdot}^{\star}\bm{\Sigma}^{\star}\bm{f}_{l}\right)^{2}-S_{k,k}^{\star}S_{i,i}^{\star}-2S_{i,k}^{\star2}\right| & \lesssim\sqrt{\frac{\log^{3}\left(n+d\right)}{n}}S_{k,k}^{\star}S_{i,i}^{\star}\\
\left|\frac{1}{n}\sum_{l=1}^{n}\left(\bm{U}_{i,\cdot}^{\star}\bm{\Sigma}^{\star}\bm{f}_{l}\right)^{2}-S_{i,i}^{\star}\right| & \lesssim\sqrt{\frac{\log\left(n+d\right)}{n}}S_{i,i}^{\star}\\
\left|\frac{1}{n}\sum_{l=1}^{n}\left(\bm{U}_{k,\cdot}^{\star}\bm{\Sigma}^{\star}\bm{f}_{l}\right)^{2}-S_{k,k}^{\star}\right| & \lesssim\sqrt{\frac{\log\left(n+d\right)}{n}}S_{k,k}^{\star}
\end{align*}
for each $k\in[d]\backslash\{i\}$. As a consequence, one can express
\begin{align*}
\alpha_{3} & =\frac{\omega_{i}^{\star2}}{np^{2}}\sum_{k:k\neq i}\omega_{k}^{\star2}\left(\bm{U}_{k,\cdot}^{\star}\bm{U}_{j,\cdot}^{\star\top}\right)^{2}+\frac{1-p}{np^{2}}S_{i,i}^{\star}\sum_{k:k\neq i}\omega_{k}^{\star2}\left(\bm{U}_{k,\cdot}^{\star}\bm{U}_{j,\cdot}^{\star\top}\right)^{2}\\
 & \quad+\frac{\left(1-p\right)\omega_{i}^{\star2}}{np^{2}}\sum_{k:k\neq i}S_{k,k}^{\star}\left(\bm{U}_{k,\cdot}^{\star}\bm{U}_{j,\cdot}^{\star\top}\right)^{2}+\frac{\left(1-p\right)^{2}}{np^{2}}\sum_{k:k\neq i}\left(S_{k,k}^{\star}S_{i,i}^{\star}+2S_{i,k}^{\star2}\right)\left(\bm{U}_{k,\cdot}^{\star}\bm{U}_{j,\cdot}^{\star\top}\right)^{2}+\widetilde{r}_{3}
\end{align*}
for some residual term $\widetilde{r}_{3}$ satisfying 
\begin{align*}
\left|\widetilde{r}_{3}\right| & \lesssim\sqrt{\frac{\log^{3}\left(n+d\right)}{n}}\left|\alpha_{3}-\widetilde{r}_{3}\right|\asymp\sqrt{\frac{\log^{3}\left(n+d\right)}{n}}\left|\alpha_{3}\right|,
\end{align*}
where the last relation holds provided that $n\gg\log^{3}(n+d)$.
This allows one to decompose $\alpha_{3}$ as follows 
\begin{align*}
\alpha_{3} & =\frac{1}{np^{2}}\sum_{k:k\neq i}\left\{ \left[\omega_{i}^{\star2}+\left(1-p\right)S_{i,i}^{\star}\right]\left[\omega_{k}^{\star2}+\left(1-p\right)S_{k,k}^{\star}\right]+2\left(1-p\right)^{2}S_{i,k}^{\star2}\right\} \left(\bm{U}_{k,\cdot}^{\star}\bm{U}_{j,\cdot}^{\star\top}\right)^{2}+\widetilde{r}_{3}\\
 & =\underbrace{\frac{1}{np^{2}}\sum_{k=1}^{d}\left\{ \left[\omega_{i}^{\star2}+\left(1-p\right)S_{i,i}^{\star}\right]\left[\omega_{k}^{\star2}+\left(1-p\right)S_{k,k}^{\star}\right]+2\left(1-p\right)^{2}S_{i,k}^{\star2}\right\} \left(\bm{U}_{k,\cdot}^{\star}\bm{U}_{j,\cdot}^{\star\top}\right)^{2}}_{\eqqcolon\alpha_{3}^{\star}}+r_{3},
\end{align*}
where we define 
\begin{align*}
r_{3} & =\widetilde{r}_{3}-\underbrace{\frac{1}{np^{2}}\left\{ \left[\omega_{i}^{\star2}+\left(1-p\right)S_{i,i}^{\star}\right]^{2}+2\left(1-p\right)^{2}S_{i,i}^{\star2}\right\} \left(\bm{U}_{i,\cdot}^{\star}\bm{U}_{j,\cdot}^{\star\top}\right)^{2}}_{\eqqcolon\delta}.
\end{align*}
Recalling from Claim \ref{claim:matrix-min-eigenvalue} that 
\begin{equation}
\sum_{k=1}^{d}S_{k,k}^{\star}\bm{U}_{k,\cdot}^{\star\top}\bm{U}_{k,\cdot}^{\star}=\sum_{k=1}^{d}\left\Vert \bm{U}_{k,\cdot}^{\star}\bm{\Sigma}^{\star}\right\Vert _{2}^{2}\bm{U}_{k,\cdot}^{\star\top}\bm{U}_{k,\cdot}^{\star}\succeq\frac{\sigma_{r}^{\star2}}{4d}\bm{I}_{r},\label{eq:ce-variance-concen-inter-1}
\end{equation}
we can reach 
\begin{align}
\alpha_{3}^{\star} & \geq\frac{1}{np^{2}}\sum_{k=1}^{d}\left\{ \left[\omega_{i}^{\star2}+\left(1-p\right)S_{i,i}^{\star}\right]\left[\omega_{k}^{\star2}+\left(1-p\right)S_{k,k}^{\star}\right]\right\} \left(\bm{U}_{k,\cdot}^{\star}\bm{U}_{j,\cdot}^{\star\top}\right)^{2}\nonumber \\
 & \gtrsim\frac{1}{np^{2}}\left[\omega_{i}^{\star2}+\left(1-p\right)S_{i,i}^{\star}\right]\bm{U}_{j,\cdot}^{\star}\left[\sum_{k=1}^{d}\omega_{k}^{\star2}\bm{U}_{k,\cdot}^{\star\top}\bm{U}_{k,\cdot}^{\star}+\sum_{k=1}^{d}S_{k,k}^{\star}\bm{U}_{k,\cdot}^{\star\top}\bm{U}_{k,\cdot}^{\star}\right]\bm{U}_{j,\cdot}^{\star\top}\nonumber \\
 & \gtrsim\frac{1}{np^{2}}\left[\omega_{i}^{\star2}+\left(1-p\right)S_{i,i}^{\star}\right]\left[\omega_{\min}^{2}\left\Vert \bm{U}_{j,\cdot}^{\star}\right\Vert _{2}^{2}+\frac{\sigma_{r}^{\star2}}{d}\left\Vert \bm{U}_{j,\cdot}^{\star}\right\Vert _{2}^{2}\right]\label{eq:ce-variance-concen-inter-2}\\
 & \gtrsim\frac{1}{np^{2}}\left[\omega_{i}^{\star2}+\left(1-p\right)S_{i,i}^{\star}\right]\left\Vert \bm{U}_{j,\cdot}^{\star}\right\Vert _{2}^{2}\left(\omega_{\min}^{2}+\frac{\sigma_{r}^{\star2}}{d}\right),\label{eq:ce-variance-concen-inter-3}
\end{align}
where the second line uses the assumption that $p$ is strictly bounded
away from 1 (so that $1-p\asymp1$), and the penultimate line makes
use of (\ref{eq:ce-variance-concen-inter-1}). This immediately leads
to 
\begin{align*}
\delta & \lesssim\frac{1}{np^{2}}\left[\omega_{i}^{\star2}+\left(1-p\right)S_{i,i}^{\star}\right]^{2}\left\Vert \bm{U}_{i,\cdot}^{\star}\right\Vert _{2}^{2}\left\Vert \bm{U}_{j,\cdot}^{\star}\right\Vert _{2}^{2}\\
 & \lesssim\frac{\mu r}{ndp^{2}}\left[\omega_{i}^{\star2}+\left(1-p\right)S_{i,i}^{\star}\right]\left(\omega_{\max}^{2}+\frac{\mu r}{d}\sigma_{1}^{\star2}\right)\left\Vert \bm{U}_{j,\cdot}^{\star}\right\Vert _{2}^{2}\\
 & \lesssim\frac{\kappa\mu^{2}r^{2}+\kappa_{\omega}\mu r}{d}\alpha_{3}^{\star},
\end{align*}
where the penultimate line comes from the assumption $1-p\asymp1$
and the fact that $S_{i,i}^{\star}=\|\bm{U}_{i,\cdot}^{\star}\|_{2}^{2}\|\bm{\Sigma}^{\star}\|^{2}\leq\frac{\mu r}{d}\sigma_{1}^{\star2}$,
and the last inequality results from (\ref{eq:ce-variance-concen-inter-3}).
As a consequence, we arrive at 
\begin{align*}
\left|r_{3}\right| & \leq\left|\widetilde{r}_{3}\right|+\left|\delta\right|\lesssim\sqrt{\frac{\log^{3}\left(n+d\right)}{n}}\left|\alpha_{3}\right|+\frac{\kappa\mu^{2}r^{2}+\kappa_{\omega}\mu r}{d}\alpha_{3}^{\star}\\
 & \lesssim\sqrt{\frac{\log^{3}\left(n+d\right)}{n}}\left(\alpha_{3}^{\star}+\left|r_{3}\right|\right)+\frac{\kappa\mu^{2}r^{2}+\kappa_{\omega}\mu r}{d}\alpha_{3}^{\star}\\
 & \lesssim\left(\sqrt{\frac{\log^{3}\left(n+d\right)}{n}}+\frac{\kappa\mu^{2}r^{2}+\kappa_{\omega}\mu r}{d}\right)\alpha_{3}^{\star},
\end{align*}
where the last relation holds true as long as $n\gg\log^{3}(n+d)$. 
\item Repeating our analysis for $\alpha_{3}$ allows one to show that 
\[
\alpha_{4}=\underbrace{\frac{1}{np^{2}}\sum_{k=1}^{d}\left\{ \left[\omega_{j}^{\star2}+\left(1-p\right)S_{j,j}^{\star}\right]\left[\omega_{k}^{\star2}+\left(1-p\right)S_{k,k}^{\star}\right]+2\left(1-p\right)^{2}S_{j,k}^{\star2}\right\} \left(\bm{U}_{k,\cdot}^{\star}\bm{U}_{i,\cdot}^{\star\top}\right)^{2}}_{\eqqcolon\alpha_{4}^{\star}}+r_{4},
\]
where the residual term $r_{4}$ obeys 
\begin{align*}
\left|r_{4}\right| & \lesssim\left(\sqrt{\frac{\log^{3}\left(n+d\right)}{n}}+\frac{\kappa\mu^{2}r^{2}+\kappa_{\omega}\mu r}{d}\right)\alpha_{4}^{\star}.
\end{align*}
\end{itemize}
Putting the above bounds together, we can conclude that 
\begin{align*}
\mathsf{var}\left(X_{i,j}|\bm{F}\right) & =\alpha_{1}+\alpha_{2}+\alpha_{3}+\alpha_{4}=\widetilde{v}_{i,j}+r_{i,j},
\end{align*}
where 
\begin{align*}
\widetilde{v}_{i,j} & \coloneqq\alpha_{1}^{\star}+\alpha_{2}^{\star}+\alpha_{3}^{\star}+\alpha_{4}^{\star}\\
 & =\frac{2\left(1-p\right)}{np}\left(S_{i,i}^{\star}S_{j,j}^{\star}+2S_{i,j}^{\star2}\right)+\frac{1}{np}\left(\omega_{i}^{\star2}S_{j,j}^{\star}+\omega_{j}^{\star2}S_{i,i}^{\star}\right)\\
 & \quad+\frac{1}{np^{2}}\sum_{k=1}^{d}\left\{ \left[\omega_{i}^{\star2}+\left(1-p\right)S_{i,i}^{\star}\right]\left[\omega_{k}^{\star2}+\left(1-p\right)S_{k,k}^{\star}\right]+2\left(1-p\right)^{2}S_{i,k}^{\star2}\right\} \left(\bm{U}_{k,\cdot}^{\star}\bm{U}_{j,\cdot}^{\star\top}\right)^{2}\\
 & \quad+\frac{1}{np^{2}}\sum_{k=1}^{d}\left\{ \left[\omega_{j}^{\star2}+\left(1-p\right)S_{j,j}^{\star}\right]\left[\omega_{k}^{\star2}+\left(1-p\right)S_{k,k}^{\star}\right]+2\left(1-p\right)^{2}S_{j,k}^{\star2}\right\} \left(\bm{U}_{k,\cdot}^{\star}\bm{U}_{i,\cdot}^{\star\top}\right)^{2},
\end{align*}
and the residual term $r_{i,j}$ is bounded in magnitude by 
\begin{align*}
\left|r_{i,j}\right| & \leq\left|r_{1}\right|+\left|r_{2}\right|+\left|r_{3}\right|+\left|r_{4}\right|\\
 & \lesssim\sqrt{\frac{\log^{3}\left(n+d\right)}{n}}\left(\alpha_{1}^{\star}+\alpha_{2}^{\star}\right)+\left(\sqrt{\frac{\log^{3}\left(n+d\right)}{n}}+\frac{\kappa\mu^{2}r^{2}+\kappa_{\omega}\mu r}{d}\right)\left(\alpha_{3}^{\star}+\alpha_{4}^{\star}\right)\\
 & \lesssim\left(\sqrt{\frac{\log^{3}\left(n+d\right)}{n}}+\frac{\kappa\mu^{2}r^{2}+\kappa_{\omega}\mu r}{d}\right)\widetilde{v}_{i,j}.
\end{align*}

To finish up, it remains to develop a lower bound on $\widetilde{v}_{i,j}$.
Towards this end, we first observe that 
\begin{align*}
\alpha_{1}^{\star}+\alpha_{2}^{\star} & =\frac{2\left(1-p\right)}{np}\left(S_{i,i}^{\star}S_{j,j}^{\star}+2S_{i,j}^{\star2}\right)+\frac{1}{np}\left(\omega_{i}^{\star2}S_{j,j}^{\star}+\omega_{j}^{\star2}S_{i,i}^{\star}\right)\\
 & \gtrsim\frac{1}{np}\left\Vert \bm{U}_{i,\cdot}^{\star}\bm{\Sigma}^{\star}\right\Vert _{2}^{2}\left\Vert \bm{U}_{j,\cdot}^{\star}\bm{\Sigma}^{\star}\right\Vert _{2}^{2}+\frac{\sigma_{r}^{\star2}}{np}\left(\omega_{j}^{\star2}\left\Vert \bm{U}_{i,\cdot}^{\star}\right\Vert _{2}^{2}+\omega_{i}^{\star2}\left\Vert \bm{U}_{j,\cdot}^{\star}\right\Vert _{2}^{2}\right),
\end{align*}
where we have made use of the assumption $1-p\asymp1$ and the elementary
inequality $S_{i,i}^{\star}=\left\Vert \bm{U}_{i,\cdot}^{\star}\bm{\Sigma}^{\star}\right\Vert _{2}^{2}\geq\sigma_{r}^{\star2}\left\Vert \bm{U}_{i,\cdot}^{\star}\right\Vert _{2}^{2}$.
In view of (\ref{eq:ce-variance-concen-inter-2}), the assumption
$1-p\asymp1$ as well as the bound $\left\Vert \bm{U}_{j,\cdot}^{\star}\bm{\Sigma}^{\star}\right\Vert _{2}^{2}\leq\sigma_{1}^{\star2}\left\Vert \bm{U}_{j,\cdot}^{\star}\right\Vert _{2}^{2}$,
we can further lower bound 
\begin{align*}
\alpha_{3}^{\star} & \gtrsim\frac{1}{np^{2}}\left[\omega_{i}^{\star2}+\left(1-p\right)S_{i,i}^{\star}\right]\left\Vert \bm{U}_{j,\cdot}^{\star}\right\Vert _{2}^{2}\left(\omega_{\min}^{2}+\frac{\sigma_{r}^{\star2}}{d}\right)\\
 & \geq\frac{\omega_{i}^{\star2}\omega_{\min}^{2}}{np^{2}}\left\Vert \bm{U}_{j,\cdot}^{\star}\right\Vert _{2}^{2}+\frac{\omega_{i}^{\star2}}{np^{2}}\frac{\sigma_{r}^{\star2}}{d}\left\Vert \bm{U}_{j,\cdot}^{\star}\right\Vert _{2}^{2}+\frac{\left(1-p\right)S_{i,i}^{\star}}{np^{2}}\left\Vert \bm{U}_{j,\cdot}^{\star}\right\Vert _{2}^{2}\frac{\sigma_{r}^{\star2}}{d}\\
 & \gtrsim\frac{\omega_{i}^{\star2}\omega_{\min}^{2}}{np^{2}}\left\Vert \bm{U}_{j,\cdot}^{\star}\right\Vert _{2}^{2}+\frac{\omega_{i}^{\star2}\sigma_{r}^{\star2}}{ndp^{2}}\left\Vert \bm{U}_{j,\cdot}^{\star}\right\Vert _{2}^{2}+\frac{1}{ndp^{2}\kappa}\left\Vert \bm{U}_{i,\cdot}^{\star}\bm{\Sigma}^{\star}\right\Vert _{2}^{2}\left\Vert \bm{U}_{j,\cdot}^{\star}\bm{\Sigma}^{\star}\right\Vert _{2}^{2},
\end{align*}
and similarly, 
\[
\alpha_{4}^{\star}\gtrsim\frac{\omega_{j}^{\star2}\omega_{\min}^{2}}{np^{2}}\left\Vert \bm{U}_{i,\cdot}^{\star}\right\Vert _{2}^{2}+\frac{\omega_{j}^{\star2}\sigma_{r}^{\star2}}{ndp^{2}}\left\Vert \bm{U}_{i,\cdot}^{\star}\right\Vert _{2}^{2}+\frac{1}{ndp^{2}\kappa}\left\Vert \bm{U}_{i,\cdot}^{\star}\bm{\Sigma}^{\star}\right\Vert _{2}^{2}\left\Vert \bm{U}_{j,\cdot}^{\star}\bm{\Sigma}^{\star}\right\Vert _{2}^{2}.
\]
The above bounds taken collectively yield 
\begin{align*}
\widetilde{v}_{i,j} & =\alpha_{1}^{\star}+\alpha_{2}^{\star}+\alpha_{3}^{\star}+\alpha_{4}^{\star}\\
 & \gtrsim\frac{1}{\min\left\{ ndp^{2}\kappa,np\right\} }\left\Vert \bm{U}_{i,\cdot}^{\star}\bm{\Sigma}^{\star}\right\Vert _{2}^{2}\left\Vert \bm{U}_{j,\cdot}^{\star}\bm{\Sigma}^{\star}\right\Vert _{2}^{2}+\left(\frac{\sigma_{r}^{\star2}}{\min\left\{ ndp^{2},np\right\} }+\frac{\omega_{\min}^{2}}{np^{2}}\right)\left(\omega_{j}^{\star2}\left\Vert \bm{U}_{i,\cdot}^{\star}\right\Vert _{2}^{2}+\omega_{i}^{\star2}\left\Vert \bm{U}_{j,\cdot}^{\star}\right\Vert _{2}^{2}\right)
\end{align*}
as claimed.

\subsubsection{Proof of Lemma \ref{lemma:ce-normal-approx-1}\label{appendix:proof-ce-normal-approx-1}}

As before, all the probabilistic arguments in this subsection are
conditional on $\bm{F}$, and it is assumed, unless otherwise noted,
that the $\sigma(\bm{F})$-measurable high-probability event $\mathcal{E}_{\mathsf{good}}$
occurs.

\paragraph{Step 1: Gaussian approximation of $X_{i,j}$ using the Berry-Esseen
Theorem.}

Recalling the definition \eqref{eq:defn-Xij-hpca} of $X_{i,j}$,
let us denote

\begin{align*}
X_{i,j} & =\sum_{l=1}^{n}\underbrace{\left\{ M_{j,l}^{\natural}E_{i,l}+M_{i,l}^{\natural}E_{j,l}+E_{i,l}\left[\mathcal{P}_{-i,\cdot}\left(\bm{E}_{\cdot,l}\right)\right]^{\top}\bm{U}^{\star}\bm{U}_{j,\cdot}^{\star\top}+E_{j,l}\left[\mathcal{P}_{-j,\cdot}\left(\bm{E}_{\cdot,l}\right)\right]^{\top}\bm{U}^{\star}\bm{U}_{i,\cdot}^{\star\top}\right\} }_{\eqqcolon\,Y_{l}},
\end{align*}
where the $Y_{l}$'s are statistically independent. Apply the Berry-Esseen
Theorem (cf.~Theorem \ref{thm:berry-esseen}) to reach 
\[
\sup_{t\in\mathbb{R}}\left|\mathbb{P}\left(X_{i,j}/\sqrt{\mathsf{var}\left(X_{i,j}|\bm{F}\right)}\leq t\,\big|\,\bm{F}\right)-\Phi\left(t\right)\right|\lesssim\gamma\left(\bm{F}\right),
\]
where 
\[
\gamma\left(\bm{F}\right)\coloneqq\frac{1}{\mathsf{var}^{3/2}\left(X_{i,j}|\bm{F}\right)}\sum_{l=1}^{n}\mathbb{E}\left[\left|Y_{l}\right|^{3}\big|\bm{F}\right].
\]
It thus boils down to controlling the quantity $\gamma(\bm{F})$,
which forms the main content of the remaining proof. 
\begin{itemize}
\item We first develop a high-probability bound on each $\vert Y_{l}\vert$.
For any $l\in[n]$, observe that 
\[
\left[\mathcal{P}_{-i,\cdot}\left(\bm{E}_{\cdot,l}\right)\right]^{\top}\bm{U}^{\star}\bm{U}_{j,\cdot}^{\star\top}=\sum_{k:k\neq i}E_{k,l}\left(\bm{U}^{\star}\bm{U}^{\star\top}\right)_{k,j}.
\]
It is straightforward to calculate that 
\begin{align*}
L & \coloneqq\max_{k:k\neq i}\left|E_{k,l}\left(\bm{U}^{\star}\bm{U}^{\star\top}\right)_{k,j}\right|\leq B\left\Vert \bm{U}_{j,\cdot}^{\star}\right\Vert _{2}\left\Vert \bm{U}^{\star}\right\Vert _{2,\infty}\lesssim B\sqrt{\frac{\mu r}{d}}\left\Vert \bm{U}_{j,\cdot}^{\star}\right\Vert _{2},\\
V & \coloneqq\sum_{k:k\neq i}\mathsf{var}\left(E_{k,l}\left(\bm{U}^{\star}\bm{U}^{\star\top}\right)_{k,j}\right)\leq\sum_{k=1}^{d}\sigma_{k,l}^{2}\left(\bm{U}^{\star}\bm{U}^{\star\top}\right)_{k,j}^{2}\leq\sigma_{\mathsf{ub}}^{2}\left\Vert \bm{U}^{\star}\bm{U}_{j,\cdot}^{\star\top}\right\Vert _{2}^{2}.
\end{align*}
In view of the Bernstein inequality \citep[Theorem 2.8.4]{vershynin2016high},
with probability exceeding $1-O((n+d)^{-101})$ we have 
\begin{align*}
\left|\left[\mathcal{P}_{-i,\cdot}\left(\bm{E}_{\cdot,l}\right)\right]^{\top}\bm{U}^{\star}\bm{U}_{j,\cdot}^{\star\top}\right| & \lesssim\sqrt{V\log\left(n+d\right)}+L\log\left(n+d\right)\\
 & \lesssim\sigma_{\mathsf{ub}}\left\Vert \bm{U}^{\star}\bm{U}_{j,\cdot}^{\star\top}\right\Vert _{2}\sqrt{\log\left(n+d\right)}+B\left\Vert \bm{U}_{j,\cdot}^{\star}\right\Vert _{2}\sqrt{\frac{\mu r}{d}}\log\left(n+d\right)\\
 & \lesssim\sigma_{\mathsf{ub}}\sqrt{\log\left(n+d\right)}\left\Vert \bm{U}_{j,\cdot}^{\star}\right\Vert _{2}+B\left\Vert \bm{U}_{j,\cdot}^{\star}\right\Vert _{2}\sqrt{\frac{\mu r}{d}}\log\left(n+d\right).
\end{align*}
By defining two quantities $B_{i}$ and $B_{j}$ such that 
\[
\max_{l\in[n]}\left|E_{i,l}\right|\leq B_{i}\qquad\text{and}\qquad\max_{l\in[n]}\left|E_{j,l}\right|\leq B_{j},
\]
we can immediately obtain from the preceding inequality that 
\begin{equation}
\left|E_{i,l}\left[\mathcal{P}_{-i,\cdot}\left(\bm{E}_{\cdot,l}\right)\right]^{\top}\bm{U}^{\star}\bm{U}_{j,\cdot}^{\star\top}\right|\lesssim\sigma_{\mathsf{ub}}B_{i}\sqrt{\log\left(n+d\right)}\left\Vert \bm{U}_{j,\cdot}^{\star}\right\Vert _{2}+BB_{i}\left\Vert \bm{U}_{j,\cdot}^{\star}\right\Vert _{2}\sqrt{\frac{\mu r}{d}}\log\left(n+d\right).\label{eq:ce-upper-bound-probabilistic}
\end{equation}
Similarly, we can also show that with probability exceeding $1-O((n+d)^{-101})$,
\[
\left|E_{j,l}\left[\mathcal{P}_{-j,\cdot}\left(\bm{E}_{\cdot,l}\right)\right]^{\top}\bm{U}^{\star}\bm{U}_{i,\cdot}^{\star\top}\right|\lesssim\sigma_{\mathsf{ub}}B_{j}\sqrt{\log\left(n+d\right)}\left\Vert \bm{U}_{i,\cdot}^{\star}\right\Vert _{2}+BB_{j}\left\Vert \bm{U}_{i,\cdot}\right\Vert _{2}\sqrt{\frac{\mu r}{d_{1}}}\log\left(n+d\right).
\]
In addition, it is easily seen from the inequality $|M_{i,l}^{\natural}|\leq\Vert\bm{U}_{i,\cdot}^{\natural}\Vert_{2}\|\bm{\Sigma}^{\natural}\|\Vert\bm{V}^{\natural}\Vert_{2,\infty}$
that 
\begin{align*}
\left|M_{j,l}^{\natural}E_{i,l}+M_{i,l}^{\natural}E_{j,l}\right| & \leq\sigma_{1}^{\natural}\left\Vert \bm{U}_{j,\cdot}^{\natural}\right\Vert _{2}\left\Vert \bm{V}^{\natural}\right\Vert _{2,\infty}B_{i}+\sigma_{1}^{\natural}\left\Vert \bm{U}_{i,\cdot}^{\natural}\right\Vert _{2}\left\Vert \bm{V}^{\natural}\right\Vert _{2,\infty}B_{j}\\
 & =\sigma_{1}^{\natural}\left\Vert \bm{V}^{\natural}\right\Vert _{2,\infty}\left(B_{i}\left\Vert \bm{U}_{j,\cdot}^{\star}\right\Vert _{2}+B_{j}\left\Vert \bm{U}_{i,\cdot}\right\Vert _{2}\right).
\end{align*}
Therefore we know that with probability exceeding $1-O((n+d)^{-101})$
\begin{align*}
\left|Y_{l}\right| & \lesssim\left[\sigma_{\mathsf{ub}}\sqrt{\log\left(n+d\right)}+B\sqrt{\frac{\mu r}{d}}\log\left(n+d\right)+\sigma_{1}^{\natural}\left\Vert \bm{V}^{\natural}\right\Vert _{2,\infty}\right]\left(B_{i}\left\Vert \bm{U}_{j,\cdot}^{\star}\right\Vert _{2}+B_{j}\left\Vert \bm{U}_{i,\cdot}\right\Vert _{2}\right).
\end{align*}
Let 
\[
C_{\mathsf{prob}}\coloneqq\widetilde{C}_{1}\left[\sigma_{\mathsf{ub}}\sqrt{\log\left(n+d\right)}+B\sqrt{\frac{\mu r}{d}}\log\left(n+d\right)+\sigma_{1}^{\natural}\left\Vert \bm{V}^{\natural}\right\Vert _{2,\infty}\right]\left(B_{i}\left\Vert \bm{U}_{j,\cdot}^{\star}\right\Vert _{2}+B_{j}\left\Vert \bm{U}_{i,\cdot}\right\Vert _{2}\right)
\]
for some sufficiently large constant $\tilde{C}>0$ such that with
probability exceeding $1-O((n+d)^{-101})$, 
\begin{equation}
\max_{l\in[n]}\left|Y_{l}\right|\leq C_{\mathsf{prob}}.\label{eq:Yl-Cprob}
\end{equation}
\item In addition, we are also in need of a deterministic upper bound on
each $\vert Y_{l}\vert$. Observe that 
\[
\left|\left[\mathcal{P}_{-i,\cdot}\left(\bm{E}_{\cdot,l}\right)\right]^{\top}\bm{U}^{\star}\bm{U}_{j,\cdot}^{\star\top}\right|\leq\left\Vert \mathcal{P}_{-i,\cdot}\left(\bm{E}_{\cdot,l}\right)\right\Vert _{2}\left\Vert \bm{U}_{j,\cdot}^{\star}\right\Vert _{2}\leq\left\Vert \bm{E}_{\cdot,l}\right\Vert _{2}\left\Vert \bm{U}_{j,\cdot}^{\star}\right\Vert _{2}\leq\sqrt{d}B\left\Vert \bm{U}_{j,\cdot}^{\star}\right\Vert _{2},
\]
and similarly, 
\begin{equation}
\left|\left[\mathcal{P}_{-j,\cdot}\left(\bm{E}_{\cdot,l}\right)\right]^{\top}\bm{U}^{\star}\bm{U}_{i,\cdot}^{\star\top}\right|\leq\sqrt{d}B\left\Vert \bm{U}_{i,\cdot}^{\star}\right\Vert _{2}.\label{eq:ce-upper-bound-deterministic}
\end{equation}
As a result, we can derive 
\begin{align}
\left|Y_{l}\right| & \leq\sigma_{1}^{\natural2}\left\Vert \bm{V}^{\natural}\right\Vert _{2,\infty}\left(B_{i}\left\Vert \bm{U}_{j,\cdot}^{\star}\right\Vert _{2}+B_{j}\left\Vert \bm{U}_{i,\cdot}\right\Vert _{2}\right)+B_{i}\sqrt{d}B\left\Vert \bm{U}_{j,\cdot}^{\star}\right\Vert _{2}+B_{j}\sqrt{d}B\left\Vert \bm{U}_{i,\cdot}^{\star}\right\Vert _{2}\nonumber \\
 & \leq\left(\sigma_{1}^{\natural2}\left\Vert \bm{V}^{\natural}\right\Vert _{2,\infty}+B\sqrt{d}\right)\left(B_{i}\left\Vert \bm{U}_{j,\cdot}^{\star}\right\Vert _{2}+B_{j}\left\Vert \bm{U}_{i,\cdot}\right\Vert _{2}\right)\eqqcolon C_{\mathsf{det}}\label{eq:C-det-bound}
\end{align}
\end{itemize}
With the above probabilistic and deterministic bounds in place (see
\eqref{eq:Yl-Cprob} and \eqref{eq:C-det-bound}), we can decompose
$\mathbb{E}\left[\left|Y_{l}\right|^{3}|\bm{F}\right]$ for each $l\in[n]$
as follows 
\begin{align*}
\mathbb{E}\left[\left|Y_{l}\right|^{3}|\bm{F}\right] & =\mathbb{E}\left[\left|Y_{l}\right|^{3}\ind_{\left|Y_{l}\right|\leq C_{\mathsf{prob}}}\big|\bm{F}\right]+\mathbb{E}\left[\left|Y_{l}\right|^{3}\ind_{\left|Y_{l}\right|>C_{\mathsf{prob}}}\big|\bm{F}\right]\\
 & \leq C_{\mathsf{prob}}\mathbb{E}\left[Y_{l}^{2}|\bm{F}\right]+C_{\mathsf{det}}^{3}\mathbb{P}\left(\left|Y_{l}\right|>C_{\mathsf{prob}}\right)\\
 & \lesssim C_{\mathsf{prob}}\mathbb{E}\left[Y_{l}^{2}|\bm{F}\right]+C_{\mathsf{det}}^{3}\left(n+d\right)^{-101}.
\end{align*}
Recognizing that $\sum_{l=1}^{n}\mathbb{E}[Y_{l}^{2}|\bm{F}]=\mathsf{var}(X_{i,j}|\bm{F})$,
we obtain 
\begin{equation}
\gamma\left(\bm{F}\right)\leq\frac{1}{\mathsf{var}^{3/2}\left(X_{i,j}|\bm{F}\right)}\sum_{l=1}^{n}\mathbb{E}\left[\left|Y_{l}\right|^{3}\big|\bm{F}\right]\lesssim\underbrace{\frac{C_{\mathsf{prob}}}{\mathsf{var}^{1/2}\left(X_{i,j}|\bm{F}\right)}}_{\eqqcolon\alpha}+\underbrace{\frac{C_{\mathsf{det}}^{3}\left(n+d\right)^{-100}}{\mathsf{var}^{3/2}\left(X_{i,j}|\bm{F}\right)}}_{\beta}.\label{eq:gamma-F-alpha-beta-decompose}
\end{equation}

It remains to control the terms $\mathsf{var}\left(X_{i,j}|\bm{F}\right)$,
$C_{\mathsf{prob}}$ and $C_{\mathsf{det}}$. On the event $\mathcal{E}_{\mathsf{good}}$,
we have learned from Lemma \ref{lemma:ce-variance-concentration}
that 
\begin{align}
\mathsf{var}^{1/2}\left(X_{i,j}|\bm{F}\right) & \asymp\widetilde{v}_{i,j}^{1/2}\gtrsim\frac{1}{\sqrt{ndp^{2}\kappa\wedge np}}\left\Vert \bm{U}_{i,\cdot}^{\star}\bm{\Sigma}^{\star}\right\Vert _{2}\left\Vert \bm{U}_{j,\cdot}^{\star}\bm{\Sigma}^{\star}\right\Vert _{2}+\frac{\omega_{\min}}{\sqrt{n}p}\left(\omega_{j}^{\star}\left\Vert \bm{U}_{i,\cdot}^{\star}\right\Vert _{2}+\omega_{i}^{\star}\left\Vert \bm{U}_{j,\cdot}^{\star}\right\Vert _{2}\right)\nonumber \\
 & \qquad\quad+\frac{\sigma_{r}^{\star}}{\sqrt{ndp^{2}\wedge np}}\left(\omega_{j}^{\star}\left\Vert \bm{U}_{i,\cdot}^{\star}\right\Vert _{2}+\omega_{i}^{\star}\left\Vert \bm{U}_{j,\cdot}^{\star}\right\Vert _{2}\right),\label{eq:ce-normal-approx-1-inter-1}
\end{align}
provided that $n\gg\log^{3}(n+d)$ and $d\gg\kappa\mu^{2}r^{2}+\kappa_{\omega}\mu r$.
Moreover, it is seen from (\ref{eq:good-event-B-hpca}) that 
\begin{equation}
B_{i}\lesssim\frac{1}{p}\sqrt{\frac{\log\left(n+d\right)}{n}}\left(\left\Vert \bm{U}_{i,\cdot}^{\star}\bm{\Sigma}^{\star}\right\Vert _{2}+\omega_{i}^{\star}\right)\qquad\text{and}\qquad B_{j}\lesssim\frac{1}{p}\sqrt{\frac{\log\left(n+d\right)}{n}}\left(\left\Vert \bm{U}_{j,\cdot}^{\star}\bm{\Sigma}^{\star}\right\Vert _{2}+\omega_{j}^{\star}\right),\label{eq:B-i-bound}
\end{equation}
which immediately lead to 
\[
B_{i}\left\Vert \bm{U}_{j,\cdot}^{\star}\right\Vert _{2}+B_{j}\left\Vert \bm{U}_{i,\cdot}\right\Vert _{2}\lesssim\frac{1}{p}\sqrt{\frac{\log\left(n+d\right)}{n}}\left[\frac{1}{\sigma_{r}^{\star}}\left\Vert \bm{U}_{i,\cdot}^{\star}\bm{\Sigma}^{\star}\right\Vert _{2}\left\Vert \bm{U}_{j,\cdot}^{\star}\bm{\Sigma}^{\star}\right\Vert _{2}+\left(\omega_{j}^{\star}\left\Vert \bm{U}_{i,\cdot}^{\star}\right\Vert _{2}+\omega_{i}^{\star}\left\Vert \bm{U}_{j,\cdot}^{\star}\right\Vert _{2}\right)\right].
\]
As a result, we can bound 
\begin{align*}
C_{\mathsf{prob}} & \lesssim\left[\frac{\sigma_{\mathsf{ub}}}{p}\sqrt{\frac{\log^{2}\left(n+d\right)}{n}}+B\sqrt{\frac{\mu r\log\left(n+d\right)}{ndp^{2}}}\log\left(n+d\right)+\sigma_{1}^{\star}\frac{\log\left(n+d\right)}{np}\right]\\
 & \quad\cdot\left[\frac{1}{\sigma_{r}^{\star}}\left\Vert \bm{U}_{i,\cdot}^{\star}\bm{\Sigma}^{\star}\right\Vert _{2}\left\Vert \bm{U}_{j,\cdot}^{\star}\bm{\Sigma}^{\star}\right\Vert _{2}+\left(\omega_{j}^{\star}\left\Vert \bm{U}_{i,\cdot}^{\star}\right\Vert _{2}+\omega_{i}^{\star}\left\Vert \bm{U}_{j,\cdot}^{\star}\right\Vert _{2}\right)\right]
\end{align*}
as well as 
\begin{align*}
C_{\mathsf{det}} & \lesssim\left(\sigma_{1}^{\star}\frac{\log\left(n+d\right)}{np}+\frac{B}{p}\sqrt{\frac{d\log\left(n+d\right)}{n}}\right)\\
 & \quad\quad\cdot\left[\frac{1}{\sigma_{r}^{\star}}\left\Vert \bm{U}_{i,\cdot}^{\star}\bm{\Sigma}^{\star}\right\Vert _{2}\left\Vert \bm{U}_{j,\cdot}^{\star}\bm{\Sigma}^{\star}\right\Vert _{2}+\left(\omega_{j}^{\star}\left\Vert \bm{U}_{i,\cdot}^{\star}\right\Vert _{2}+\omega_{i}^{\star}\left\Vert \bm{U}_{j,\cdot}^{\star}\right\Vert _{2}\right)\right],
\end{align*}
where we have utilized (\ref{eq:good-event-sigma-least-largest-hpca})
and (\ref{eq:good-event-V-natural-2-infty-hpca}). In what follows,
we bound the terms $\alpha$ and $\beta$ in \eqref{eq:gamma-F-alpha-beta-decompose}
separately. 
\begin{itemize}
\item Regarding $\alpha$, one first observes that 
\begin{align*}
\alpha\widetilde{v}_{i,j}^{1/2} & \lesssim\underbrace{\frac{\sigma_{\mathsf{ub}}}{p\sigma_{r}^{\star}}\sqrt{\frac{\log^{2}\left(n+d\right)}{n}}\left\Vert \bm{U}_{i,\cdot}^{\star}\bm{\Sigma}^{\star}\right\Vert _{2}\left\Vert \bm{U}_{j,\cdot}^{\star}\bm{\Sigma}^{\star}\right\Vert _{2}}_{\eqqcolon\alpha_{1}}+\underbrace{\frac{\sigma_{\mathsf{ub}}}{p}\sqrt{\frac{\log^{2}\left(n+d\right)}{n}}\left(\omega_{j}^{\star}\left\Vert \bm{U}_{i,\cdot}^{\star}\right\Vert _{2}+\omega_{i}^{\star}\left\Vert \bm{U}_{j,\cdot}^{\star}\right\Vert _{2}\right)}_{\eqqcolon\alpha_{2}}\\
 & \quad+\underbrace{B\sqrt{\frac{\mu r\log^{3}\left(n+d\right)}{ndp^{2}}}\frac{1}{\sigma_{r}^{\star}}\left\Vert \bm{U}_{i,\cdot}^{\star}\bm{\Sigma}^{\star}\right\Vert _{2}\left\Vert \bm{U}_{j,\cdot}^{\star}\bm{\Sigma}^{\star}\right\Vert _{2}}_{\eqqcolon\alpha_{3}}+\underbrace{B\sqrt{\frac{\mu r\log^{3}\left(n+d\right)}{ndp^{2}}}\left(\omega_{j}^{\star}\left\Vert \bm{U}_{i,\cdot}^{\star}\right\Vert _{2}+\omega_{i}^{\star}\left\Vert \bm{U}_{j,\cdot}^{\star}\right\Vert _{2}\right)}_{\eqqcolon\alpha_{4}}\\
 & \quad+\underbrace{\frac{\sqrt{\kappa}\log\left(n+d\right)}{np}\left\Vert \bm{U}_{i,\cdot}^{\star}\bm{\Sigma}^{\star}\right\Vert _{2}\left\Vert \bm{U}_{j,\cdot}^{\star}\bm{\Sigma}^{\star}\right\Vert _{2}}_{\eqqcolon\alpha_{5}}+\underbrace{\sigma_{1}^{\star}\frac{\log\left(n+d\right)}{np}\left(\omega_{j}^{\star}\left\Vert \bm{U}_{i,\cdot}^{\star}\right\Vert _{2}+\omega_{i}^{\star}\left\Vert \bm{U}_{j,\cdot}^{\star}\right\Vert _{2}\right)}_{\eqqcolon\alpha_{6}}\\
 & \lesssim\frac{1}{\sqrt{\log\left(n+d\right)}}\widetilde{v}_{i,j}^{1/2}.
\end{align*}
To see why the last step holds, we note that due to (\ref{eq:good-event-sigma-hpca}),
(\ref{eq:good-event-B-hpca}), the following inequalities hold: 
\begin{align*}
\alpha_{1} & \asymp\frac{1}{p\sigma_{r}^{\star}}\left(\sqrt{\frac{\mu r\log\left(n+d\right)}{ndp}}\sigma_{1}^{\star}+\frac{\omega_{\max}}{\sqrt{np}}\right)\sqrt{\frac{\log^{2}\left(n+d\right)}{n}}\left\Vert \bm{U}_{i,\cdot}^{\star}\bm{\Sigma}^{\star}\right\Vert _{2}\left\Vert \bm{U}_{j,\cdot}^{\star}\bm{\Sigma}^{\star}\right\Vert _{2}\\
 & \lesssim\frac{1}{\sqrt{\log\left(n+d\right)}}\left(\frac{1}{\sqrt{ndp^{2}\kappa}}+\frac{1}{\sqrt{np}}\right)\left\Vert \bm{U}_{i,\cdot}^{\star}\bm{\Sigma}^{\star}\right\Vert _{2}\left\Vert \bm{U}_{j,\cdot}^{\star}\bm{\Sigma}^{\star}\right\Vert _{2}\lesssim\frac{1}{\sqrt{\log\left(n+d\right)}}\widetilde{v}_{i,j}^{1/2},\\
\alpha_{2} & \asymp\frac{1}{p}\left(\sqrt{\frac{\mu r\log\left(n+d\right)}{ndp}}\sigma_{1}^{\star}+\frac{\omega_{\max}}{\sqrt{np}}\right)\sqrt{\frac{\log^{2}\left(n+d\right)}{n}}\left(\omega_{j}^{\star}\left\Vert \bm{U}_{i,\cdot}^{\star}\right\Vert _{2}+\omega_{i}^{\star}\left\Vert \bm{U}_{j,\cdot}^{\star}\right\Vert _{2}\right)\\
 & \lesssim\frac{1}{\sqrt{\log\left(n+d\right)}}\left(\frac{1}{\sqrt{ndp^{2}}}+\frac{1}{\sqrt{np}}\right)\sigma_{r}^{\star}\left(\omega_{j}^{\star}\left\Vert \bm{U}_{i,\cdot}^{\star}\right\Vert _{2}+\omega_{i}^{\star}\left\Vert \bm{U}_{j,\cdot}^{\star}\right\Vert _{2}\right)\lesssim\frac{1}{\sqrt{\log\left(n+d\right)}}\widetilde{v}_{i,j}^{1/2},\\
\alpha_{3} & \asymp\left(\frac{1}{p}\sqrt{\frac{\kappa\mu^{2}r^{2}\log^{4}\left(n+d\right)}{n^{2}d^{2}p^{2}}}+\frac{\omega_{\max}}{p\sigma_{r}^{\star}}\sqrt{\frac{\mu r\log^{4}\left(n+d\right)}{n^{2}dp^{2}}}\right)\left\Vert \bm{U}_{i,\cdot}^{\star}\bm{\Sigma}^{\star}\right\Vert _{2}\left\Vert \bm{U}_{j,\cdot}^{\star}\bm{\Sigma}^{\star}\right\Vert _{2}\\
 & \lesssim\frac{1}{\sqrt{\log\left(n+d\right)}}\frac{1}{\sqrt{ndp^{2}\kappa}}\left\Vert \bm{U}_{i,\cdot}^{\star}\bm{\Sigma}^{\star}\right\Vert _{2}\left\Vert \bm{U}_{j,\cdot}^{\star}\bm{\Sigma}^{\star}\right\Vert _{2}\lesssim\frac{1}{\sqrt{\log\left(n+d\right)}}\widetilde{v}_{i,j}^{1/2},\\
\alpha_{4} & \asymp\left(\frac{\mu r\log^{2}\left(n+d\right)}{ndp^{2}}\sigma_{1}^{\star}+\frac{\omega_{\max}}{p}\sqrt{\frac{\mu r\log^{4}\left(n+d\right)}{n^{2}dp^{2}}}\right)\left(\omega_{j}^{\star}\left\Vert \bm{U}_{i,\cdot}^{\star}\right\Vert _{2}+\omega_{i}^{\star}\left\Vert \bm{U}_{j,\cdot}^{\star}\right\Vert _{2}\right)\\
 & \lesssim\frac{1}{\sqrt{\log\left(n+d\right)}}\frac{\sigma_{r}^{\star}}{\sqrt{ndp^{2}}}\left(\omega_{j}^{\star}\left\Vert \bm{U}_{i,\cdot}^{\star}\right\Vert _{2}+\omega_{i}^{\star}\left\Vert \bm{U}_{j,\cdot}^{\star}\right\Vert _{2}\right)\lesssim\frac{1}{\sqrt{\log\left(n+d\right)}}\widetilde{v}_{i,j}^{1/2},\\
\alpha_{5} & \lesssim\frac{1}{\sqrt{\log\left(n+d\right)}}\frac{1}{\sqrt{np}}\left\Vert \bm{U}_{i,\cdot}^{\star}\bm{\Sigma}^{\star}\right\Vert _{2}\left\Vert \bm{U}_{j,\cdot}^{\star}\bm{\Sigma}^{\star}\right\Vert _{2}\lesssim\frac{1}{\sqrt{\log\left(n+d\right)}}\widetilde{v}_{i,j}^{1/2},\\
\alpha_{6} & \lesssim\frac{1}{\sqrt{\log\left(n+d\right)}}\frac{\sigma_{r}^{\star}}{\sqrt{np}}\left(\omega_{j}^{\star}\left\Vert \bm{U}_{i,\cdot}^{\star}\right\Vert _{2}+\omega_{i}^{\star}\left\Vert \bm{U}_{j,\cdot}^{\star}\right\Vert _{2}\right)\lesssim\frac{1}{\sqrt{\log\left(n+d\right)}}\widetilde{v}_{i,j}^{1/2},
\end{align*}
provided that $np\gtrsim\kappa^{2}\mu r\log^{4}(n+d)$, $ndp^{2}\gtrsim\kappa^{2}\mu^{2}r^{2}\log^{5}(n+d)$,
\begin{equation}
\omega_{\max}\sqrt{\frac{d}{np}}\lesssim\frac{1}{\sqrt{\kappa\log^{3}\left(n+d\right)}},\qquad\text{and}\qquad\frac{\omega_{\max}}{\sigma_{r}^{\star}}\sqrt{\frac{1}{np^{2}}}\lesssim\frac{1}{\sqrt{\kappa\mu r\log^{5}\left(n+d\right)}}.\label{eq:omega-max-condition-1356}
\end{equation}
In view of (\ref{eq:pca-1st-err-useful}), we know that the second
condition in \eqref{eq:omega-max-condition-1356} can be guaranteed
by $ndp^{2}\gtrsim\kappa\mu r\log^{5}(n+d)$ and 
\[
\frac{\omega_{\max}^{2}}{p\sigma_{r}^{\star2}}\sqrt{\frac{d}{n}}\lesssim\frac{1}{\sqrt{\kappa\mu r\log^{5}\left(n+d\right)}}.
\]
As a result, we conclude that 
\[
\alpha\lesssim\frac{1}{\sqrt{\log\left(n+d\right)}}
\]
holds as long as $np\gtrsim\kappa^{2}\mu r\log^{4}(n+d)$, $ndp^{2}\gtrsim\kappa^{2}\mu^{2}r^{2}\log^{5}(n+d)$,
\[
\omega_{\max}\sqrt{\frac{d}{np}}\lesssim\frac{1}{\sqrt{\kappa\log^{3}\left(n+d\right)}},\qquad\text{and}\qquad\frac{\omega_{\max}^{2}}{p\sigma_{r}^{\star2}}\sqrt{\frac{d}{n}}\lesssim\frac{1}{\sqrt{\kappa\mu r\log^{5}\left(n+d\right)}}.
\]
\item Turning to the term $\beta$, we make the observation that 
\begin{align*}
C_{\mathsf{det}}\widetilde{v}_{i,j}^{-1/2} & \overset{\text{(i)}}{\lesssim}\left(\sigma_{1}^{\star}\frac{\log\left(n+d\right)}{np}+\frac{B}{p}\sqrt{\frac{d\log\left(n+d\right)}{n}}\right)\left[\frac{1}{\sigma_{r}^{\star}}\left\Vert \bm{U}_{i,\cdot}^{\star}\bm{\Sigma}^{\star}\right\Vert _{2}\left\Vert \bm{U}_{j,\cdot}^{\star}\bm{\Sigma}^{\star}\right\Vert _{2}+\left(\omega_{j}^{\star}\left\Vert \bm{U}_{i,\cdot}^{\star}\right\Vert _{2}+\omega_{i}^{\star}\left\Vert \bm{U}_{j,\cdot}^{\star}\right\Vert _{2}\right)\right]\\
 & \quad\cdot\left\{ \frac{1}{\sqrt{np}}\left[\left\Vert \bm{U}_{i,\cdot}^{\star}\bm{\Sigma}^{\star}\right\Vert _{2}\left\Vert \bm{U}_{j,\cdot}^{\star}\bm{\Sigma}^{\star}\right\Vert _{2}+\sigma_{r}^{\star}\left(\omega_{j}^{\star}\left\Vert \bm{U}_{i,\cdot}^{\star}\right\Vert _{2}+\omega_{i}^{\star}\left\Vert \bm{U}_{j,\cdot}^{\star}\right\Vert _{2}\right)\right]\right\} ^{-1}\\
 & \lesssim\left(\sigma_{1}^{\star}\frac{\log\left(n+d\right)}{np}+\frac{B}{p}\sqrt{\frac{d\log\left(n+d\right)}{n}}\right)\frac{\sqrt{np}}{\sigma_{r}^{\star}}\overset{\text{(ii)}}{\lesssim}\left(\frac{\log\left(n+d\right)}{np}+\frac{1}{p}\sqrt{\frac{d\log\left(n+d\right)}{n}}\right)\sqrt{\kappa np}\\
 & \lesssim\frac{\sqrt{\kappa}\log\left(n+d\right)}{\sqrt{np}}+\sqrt{\frac{\kappa d\log\left(n+d\right)}{p}}\overset{\text{(iii)}}{\lesssim}\sqrt{nd}.
\end{align*}
Here (i) follows from (\ref{eq:ce-normal-approx-1-inter-1}); (ii)
holds as long as $B\lesssim\sigma_{r}^{\star}$, which can be guaranteed
by $ndp^{2}\gtrsim\kappa\mu r\log(n+d)$ and 
\[
\frac{\omega_{\max}}{p\sigma_{r}^{\star}}\sqrt{\frac{1}{n}}\lesssim\frac{\omega_{\max}^{2}}{p\sigma_{r}^{\star2}}\sqrt{\frac{d}{n}}+\frac{1}{\sqrt{nd}p}\lesssim\frac{1}{\sqrt{\log\left(n+d\right)}};
\]
and (iii) holds as long as $np\gtrsim\kappa\log(n+d)$. This immediately
results in 
\[
\beta\lesssim\frac{C_{\mathsf{det}}^{3}\left(n+d\right)^{-100}}{\widetilde{v}_{i,j}^{3/2}}\lesssim\left(nd\right)^{3/2}\left(n+d\right)^{-100}\lesssim\left(n+d\right)^{-50}.
\]
\end{itemize}
Putting the above pieces together, we have shown that 
\begin{equation}
\sup_{t\in\mathbb{R}}\left|\mathbb{P}\left(X_{i,j}/\sqrt{\mathsf{var}\left(X_{i,j}|\bm{F}\right)}\leq t\,\Big|\,\bm{F}\right)-\Phi\left(t\right)\right|\lesssim\gamma\left(\bm{F}\right)\lesssim\alpha+\beta\lesssim\frac{1}{\sqrt{\log\left(n+d\right)}},\label{eq:ce-normal-approx-1-inter-2}
\end{equation}
provided that $np\gtrsim\kappa^{2}\mu r\log^{4}(n+d)$, $ndp^{2}\gtrsim\kappa^{2}\mu^{2}r^{2}\log^{5}(n+d)$,
\[
\omega_{\max}\sqrt{\frac{d}{np}}\lesssim\frac{1}{\sqrt{\kappa\log^{3}\left(n+d\right)}},\qquad\text{and}\qquad\frac{\omega_{\max}^{2}}{p\sigma_{r}^{\star2}}\sqrt{\frac{d}{n}}\lesssim\frac{1}{\sqrt{\kappa\mu r\log^{5}\left(n+d\right)}}.
\]

\paragraph{Step 2: derandomizing the conditional variance.}

In this step, we intend to replace $\mathsf{var}(X_{i,j}|\bm{F})$
in (\ref{eq:ce-normal-approx-1-inter-2}) with $\widetilde{v}_{i,j}$.
Towards this, it is first observed that 
\begin{align*}
\mathbb{P}\left(X_{i,j}/\sqrt{\widetilde{v}_{i,j}}\leq t\,\Big|\,\bm{F}\right)-\Phi\left(t\right) & =\underbrace{\mathbb{P}\left(X_{i,j}/\sqrt{\widetilde{v}_{i,j}}\leq t\,\Big|\,\bm{F}\right)-\mathbb{P}\left(X_{i,j}/\sqrt{\mathsf{var}\left(X_{i,j}|\bm{F}\right)}\leq t\,\Big|\,\bm{F}\right)}_{\eqqcolon\beta_{1}}\\
 & \quad+\underbrace{\mathbb{P}\left(X_{i,j}/\sqrt{\mathsf{var}\left(X_{i,j}|\bm{F}\right)}\leq t\,\Big|\,\bm{F}\right)-\Phi\left(t\right)}_{\eqqcolon\beta_{2}}.
\end{align*}
Regarding $\beta_{2}$, it follows from (\ref{eq:ce-normal-approx-1-inter-2})
that when $\mathcal{E}_{\mathsf{good}}$ occurs, one has 
\begin{align*}
\left|\beta_{2}\right| & \lesssim\frac{1}{\sqrt{\log\left(n+d\right)}}.
\end{align*}
We then turn attention to bounding $\beta_{1}$. In view of Lemma
\ref{lemma:ce-variance-concentration}, we know that when $\mathcal{E}_{\mathsf{good}}$
happens, one has 
\begin{equation}
\mathsf{var}\left(X_{i,j}|\bm{F}\right)=\widetilde{v}_{i,j}+O\left(\sqrt{\frac{\log^{3}\left(n+d\right)}{n}}+\frac{\kappa\mu^{2}r^{2}}{d}\right)\widetilde{v}_{i,j}.\label{eq:ce-normal-approx-1-inter-3-1}
\end{equation}
An immediate consequence is that 
\begin{equation}
\frac{1}{2}\widetilde{v}_{i,j}\leq\mathsf{var}\left(X_{i,j}|\bm{F}\right)\leq2\widetilde{v}_{i,j},\label{eq:ce-normal-approx-1-inter-4-1}
\end{equation}
with the proviso that $n\gg\log^{3}(n+d)$ and $d\gg\kappa\mu^{2}r^{2}$.
Taking the above two equations collectively yields 
\begin{equation}
\left|\sqrt{\mathsf{var}\left(X_{i,j}|\bm{F}\right)}-\sqrt{\widetilde{v}_{i,j}}\right|=\frac{\left|\mathsf{var}\left(X_{i,j}|\bm{F}\right)-\widetilde{v}_{i,j}\right|}{\sqrt{\mathsf{var}\left(X_{i,j}|\bm{F}\right)}+\sqrt{\widetilde{v}_{i,j}}}\leq\underbrace{\tilde{C}\left(\sqrt{\frac{\log\left(n+d\right)^{3}}{n}}+\frac{\kappa\mu^{2}r^{2}}{d}\right)}_{\eqqcolon\delta}\sqrt{\widetilde{v}_{i,j}}\label{eq:ce-normal-approx-1-inter-5-1}
\end{equation}
for some sufficiently large constant $\tilde{C}>0$. Consequently,
we arrive at 
\begin{align*}
\mathbb{P}\left(X_{i,j}/\sqrt{\widetilde{v}_{i,j}}\leq t\,\Big|\,\bm{F}\right) & =\mathbb{P}\left(X_{i,j}\leq t\sqrt{\widetilde{v}_{i,j}}\,\Big|\,\bm{F}\right)\overset{\text{(i)}}{\leq}\mathbb{P}\left(X_{i,j}\leq t\sqrt{\mathsf{var}\left(X_{i,j}|\bm{F}\right)}+t\delta\sqrt{\widetilde{v}_{i,j}}\,\Big|\,\bm{F}\right)\\
 & \overset{\text{(ii)}}{\leq}\mathbb{P}\left(X_{i,j}/\sqrt{\mathsf{var}\left(X_{i,j}|\bm{F}\right)}\leq t+\sqrt{2}t\delta\,\Big|\,\bm{F}\right)\\
 & \overset{\text{(iii)}}{\leq}\Phi\left(t+\sqrt{2}t\delta\right)+O\left(\frac{1}{\sqrt{\log\left(n+d\right)}}\right)\overset{\text{(iv)}}{\leq}\Phi\left(t\right)+O\left(\frac{1}{\sqrt{\log\left(n+d\right)}}\right)\\
 & \overset{\text{(v)}}{\leq}\mathbb{P}\left(X_{i,j}/\sqrt{\mathsf{var}\left(X_{i,j}|\bm{F}\right)}\leq t\,\Big|\,\bm{F}\right)+O\left(\frac{1}{\sqrt{\log\left(n+d\right)}}\right).
\end{align*}
Here, (i) follows from (\ref{eq:ce-normal-approx-1-inter-5-1}); (ii)
is a consequence of (\ref{eq:ce-normal-approx-1-inter-4-1}); (iii)
and (v) comes from (\ref{eq:ce-normal-approx-1-inter-2}); and (iv)
arises from 
\[
\Phi\left(t+\sqrt{2}t\delta\right)-\Phi\left(t\right)=\int_{t}^{t+\sqrt{2}t\delta}\phi\left(s\right)\mathrm{d}s\leq\sqrt{2}t\phi\left(t\right)\delta\leq\sqrt{2}\cdot\sqrt{\frac{1}{2\pi e}}\delta\lesssim\frac{1}{\sqrt{\log\left(n+d\right)}}
\]
provided that $n\gtrsim\log^{4}(n+d)$ and $d\gtrsim\kappa\mu^{2}r^{2}\sqrt{\log(n+d)}$,
where we use the fact that $\sup_{t\in\mathbb{R}}t\phi(t)=\phi(1)=\sqrt{1/(2\pi e)}$.
Similarly, we can show that 
\[
\mathbb{P}\left(X_{i,j}/\sqrt{\widetilde{v}_{i,j}}\leq t\,\Big|\,\bm{F}\right)\geq\left[\mathbb{P}\left(X_{i,j}/\sqrt{\mathsf{var}\left(X_{i,j}|\bm{F}\right)}\leq t\,\Big|\,\bm{F}\right)-O\left(\frac{1}{\sqrt{\log\left(n+d\right)}}\right)\right].
\]
As a consequence, we arrive at 
\[
\left|\beta_{1}\right|\lesssim\frac{1}{\sqrt{\log\left(n+d\right)}}.
\]
Combine the preceding bounds on $\beta_{1}$ and $\beta_{2}$ to reach
\[
\left|\mathbb{P}\left(X_{i,j}/\sqrt{\widetilde{v}_{i,j}}\leq t\,\Big|\,\bm{F}\right)-\Phi\left(t\right)\right|\leq\left(\left|\beta_{1}\right|+\left|\beta_{2}\right|\right)\lesssim\frac{1}{\sqrt{\log\left(n+d\right)}}
\]
for all $t\in\mathbb{R}$. 

\paragraph{Step 3: taking higher-order errors into account. }

By following the same analysis as Step 3 in Appendix \ref{appendix:proof-pca-normal-approximation}
(proof of Lemma \ref{lemma:pca-normal-approximation}), we know that
if one can show 
\begin{equation}
\mathbb{P}\left(\widetilde{v}_{i,j}^{-1/2}\left|\Phi_{i,j}\right|\lesssim\log^{-1/2}\left(n+d\right)\,\big|\,\bm{F}\right)\geq1-O\left(\left(n+d\right)^{-10}\right),\label{eq:ce-normal-approx-inter-3}
\end{equation}
then it holds that 
\[
\sup_{t\in\mathbb{R}}\left|\mathbb{P}\left(\widetilde{v}_{i,j}^{-1/2}\left(\bm{S}-\bm{M}^{\natural}\bm{M}^{\natural\top}\right)_{i,j}\leq t\,\big|\,\bm{F}\right)-\Phi\left(t\right)\right|\lesssim\frac{1}{\sqrt{\log\left(n+d\right)}}.
\]
As a result, we shall focus on proving (\ref{eq:ce-normal-approx-inter-3})
from now on. Recall from Lemma \ref{lemma:ce-2nd-error} that with
probability exceeding $1-O((n+d)^{-10})$, 
\begin{align}
\left|\Phi_{i,j}\right| & \lesssim\zeta_{i,j}\asymp\underbrace{\theta^{2}\left[\sigma_{1}^{\star}\left(\omega_{i}^{\star}\left\Vert \bm{U}_{j,\cdot}^{\star}\right\Vert _{2}+\omega_{j}^{\star}\left\Vert \bm{U}_{i,\cdot}^{\star}\right\Vert _{2}\right)+\left\Vert \bm{U}_{i,\cdot}^{\star}\bm{\Sigma}^{\star}\right\Vert _{2}\left\Vert \bm{U}_{j,\cdot}^{\star}\bm{\Sigma}^{\star}\right\Vert _{2}\right]}_{\eqqcolon\gamma_{1}}\nonumber \\
 & \qquad\quad+\underbrace{\sigma_{1}^{\star2}\left(\left\Vert \bm{U}_{i,\cdot}^{\star}\right\Vert _{2}\zeta_{\mathsf{2nd},j}+\left\Vert \bm{U}_{j,\cdot}^{\star}\right\Vert _{2}\zeta_{\mathsf{2nd},i}\right)}_{\eqqcolon\gamma_{2}}+\underbrace{\theta^{2}\omega_{i}^{\star}\omega_{j}^{\star}}_{\eqqcolon\gamma_{3}}\nonumber \\
 & \qquad\quad+\underbrace{\sigma_{1}^{\star2}\zeta_{\mathsf{2nd},i}\zeta_{\mathsf{2nd},j}}_{\eqqcolon\gamma_{4}}+\underbrace{\theta\sigma_{1}^{\star}\left(\omega_{i}^{\star}\zeta_{\mathsf{2nd},j}+\omega_{j}^{\star}\zeta_{\mathsf{2nd},i}\right)}_{\eqqcolon\gamma_{5}},\label{eq:Phi-ij-gamma-1234}
\end{align}
and from Lemma \ref{lemma:ce-variance-concentration} that 
\begin{align*}
\widetilde{v}_{i,j}^{1/2} & \gtrsim\frac{1}{\sqrt{\min\left\{ ndp^{2}\kappa,np\right\} }}\left\Vert \bm{U}_{i,\cdot}^{\star}\bm{\Sigma}^{\star}\right\Vert _{2}\left\Vert \bm{U}_{j,\cdot}^{\star}\bm{\Sigma}^{\star}\right\Vert _{2}+\frac{\omega_{\min}}{\sqrt{n}p}\left(\omega_{j}^{\star}\left\Vert \bm{U}_{i,\cdot}^{\star}\right\Vert _{2}+\omega_{i}^{\star}\left\Vert \bm{U}_{j,\cdot}^{\star}\right\Vert _{2}\right)\\
 & \quad+\frac{\sigma_{r}^{\star}}{\sqrt{ndp^{2}\wedge np}}\left(\omega_{j}^{\star}\left\Vert \bm{U}_{i,\cdot}^{\star}\right\Vert _{2}+\omega_{i}^{\star}\left\Vert \bm{U}_{j,\cdot}^{\star}\right\Vert _{2}\right).
\end{align*}
Armed with these bounds, we seek to derive sufficient conditions that
guarantee $\zeta_{i,j}\lesssim\delta\widetilde{v}_{i,j}^{1/2}$. 
\begin{itemize}
\item Regarding the quantity $\gamma_{1}$ defined in \eqref{eq:Phi-ij-gamma-1234},
we note that 
\begin{align*}
\gamma_{1} & \lesssim\left(\frac{\kappa r\log^{2}\left(n+d\right)}{np}+\frac{\kappa^{2}\mu r^{2}\log^{3}\left(n+d\right)}{ndp^{2}}+\frac{\omega_{\max}^{2}}{\sigma_{r}^{\star2}}\frac{\kappa r\log^{2}\left(n+d\right)}{np^{2}}\right)\\
 & \quad\cdot\left[\sigma_{1}^{\star}\left(\omega_{i}^{\star}\left\Vert \bm{U}_{j,\cdot}^{\star}\right\Vert _{2}+\omega_{j}^{\star}\left\Vert \bm{U}_{i,\cdot}^{\star}\right\Vert _{2}\right)+\left\Vert \bm{U}_{i,\cdot}^{\star}\bm{\Sigma}^{\star}\right\Vert _{2}\left\Vert \bm{U}_{j,\cdot}^{\star}\bm{\Sigma}^{\star}\right\Vert _{2}\right]\\
 & \lesssim\delta\frac{1}{\sqrt{\min\left\{ ndp^{2}\kappa,np\right\} }}\left\Vert \bm{U}_{i,\cdot}^{\star}\bm{\Sigma}^{\star}\right\Vert _{2}\left\Vert \bm{U}_{j,\cdot}^{\star}\bm{\Sigma}^{\star}\right\Vert _{2}+\delta\frac{\sigma_{r}^{\star}}{\sqrt{\min\left\{ ndp^{2},np\right\} }}\left(\omega_{j}^{\star}\left\Vert \bm{U}_{i,\cdot}^{\star}\right\Vert _{2}+\omega_{i}^{\star}\left\Vert \bm{U}_{j,\cdot}^{\star}\right\Vert _{2}\right)\\
 & \lesssim\delta\widetilde{v}_{i,j}^{1/2},
\end{align*}
provided that $np\gtrsim\delta^{-2}\kappa^{3}r^{2}\log^{4}(n+d)$,
$ndp^{2}\gtrsim\delta^{-2}\kappa^{5}\mu^{2}r^{4}\log^{6}(n+d)$, and
\[
\frac{\omega_{\max}^{2}}{p\sigma_{r}^{\star2}}\sqrt{\frac{d}{n}}\lesssim\frac{\delta}{\kappa^{3/2}r\log^{2}\left(n+d\right)}.
\]
\item Regarding the term $\gamma_{2}$ defined in \eqref{eq:Phi-ij-gamma-1234},
we recall from the proof of Lemma \ref{lemma:pca-normal-approximation}
(Step 3 in Appendix \ref{appendix:proof-pca-normal-approximation})
that with probability exceeding $1-O((n+d)^{-10})$, \begin{subequations}\label{eq:ce-normal-approx-1-inter-6}
\begin{align}
\zeta_{\mathsf{2nd},i} & \lesssim\frac{\delta}{\sqrt{\kappa}}\frac{\left\Vert \bm{U}_{i,\cdot}^{\star}\bm{\Sigma}^{\star}\right\Vert _{2}+\omega_{i}^{\star}}{\sqrt{\min\left\{ ndp^{2}\kappa,np\right\} }\sigma_{1}^{\star}}+\frac{\delta}{\sqrt{\kappa}}\frac{\omega_{\min}\omega_{i}^{\star}}{\sqrt{np^{2}}\sigma_{1}^{\star2}},\\
\zeta_{\mathsf{2nd},j} & \lesssim\frac{\delta}{\sqrt{\kappa}}\frac{\left\Vert \bm{U}_{j,\cdot}^{\star}\bm{\Sigma}^{\star}\right\Vert _{2}+\omega_{j}^{\star}}{\sqrt{\min\left\{ ndp^{2}\kappa,np\right\} }\sigma_{1}^{\star}}+\frac{\delta}{\sqrt{\kappa}}\frac{\omega_{\min}\omega_{j}^{\star}}{\sqrt{np^{2}}\sigma_{1}^{\star2}},
\end{align}
\end{subequations}provided that $d\gtrsim\delta^{-2}\kappa^{8}\mu^{3}r^{3}\kappa_{\omega}^{2}\log^{4}(n+d)$,
\begin{align*}
ndp^{2}\gtrsim\delta^{-2}\kappa^{10}\mu^{4}r^{4}\kappa_{\omega}^{2}\log^{8}\left(n+d\right), & \qquad np\gtrsim\delta^{-2}\kappa^{10}\mu^{3}r^{3}\kappa_{\omega}^{2}\log^{6}\left(n+d\right),\\
\text{and}\qquad\frac{\omega_{\max}^{2}}{p\sigma_{r}^{\star2}}\sqrt{\frac{d}{n}}\lesssim\frac{\delta}{\sqrt{\kappa^{7}\mu^{2}r^{2}\kappa_{\omega}\log^{6}\left(n+d\right)}}, & \qquad\frac{\omega_{\max}}{\sigma_{r}^{\star}}\sqrt{\frac{d}{np}}\lesssim\frac{\delta}{\sqrt{\kappa^{8}\mu^{2}r^{2}\kappa_{\omega}\log^{5}\left(n+d\right)}}.
\end{align*}
Consequently, we obtain 
\begin{align*}
\gamma_{2} & \lesssim\delta\frac{\left\Vert \bm{U}_{i,\cdot}^{\star}\bm{\Sigma}^{\star}\right\Vert _{2}\left\Vert \bm{U}_{j,\cdot}^{\star}\bm{\Sigma}^{\star}\right\Vert _{2}}{\sqrt{\min\left\{ ndp^{2}\kappa,np\right\} }}+\delta\sigma_{r}^{\star}\frac{\omega_{j}^{\star}\left\Vert \bm{U}_{i,\cdot}^{\star}\right\Vert _{2}+\omega_{i}^{\star}\left\Vert \bm{U}_{j,\cdot}^{\star}\right\Vert _{2}}{\sqrt{\min\left\{ ndp^{2}\kappa,np\right\} }}\\
 & \quad+\frac{\delta}{\sqrt{\kappa}}\frac{\omega_{\min}}{\sqrt{np^{2}}}\left(\omega_{j}^{\star}\left\Vert \bm{U}_{i,\cdot}^{\star}\right\Vert _{2}+\omega_{i}^{\star}\left\Vert \bm{U}_{j,\cdot}^{\star}\right\Vert _{2}\right)\\
 & \lesssim\delta\widetilde{v}_{i,j}^{1/2}.
\end{align*}
\item We now move on to the term $\gamma_{3}$ defined in \eqref{eq:Phi-ij-gamma-1234},
which obeys 
\begin{align*}
\gamma_{3} & \asymp\underbrace{\frac{\kappa r\log^{2}\left(n+d\right)}{np}\omega_{i}^{\star}\omega_{j}^{\star}}_{\gamma_{3,1}}+\underbrace{\frac{\kappa^{2}\mu r^{2}\log^{3}\left(n+d\right)}{ndp^{2}}\omega_{i}^{\star}\omega_{j}^{\star}}_{\eqqcolon\gamma_{3,2}}+\underbrace{\frac{\omega_{\max}^{2}\omega_{i}^{\star}\omega_{j}^{\star}}{\sigma_{r}^{\star2}}\frac{\kappa r\log^{2}\left(n+d\right)}{np^{2}}}_{\eqqcolon\gamma_{3,3}}\lesssim\delta\widetilde{v}_{i,j}^{1/2}.
\end{align*}
Here, the last relation holds since 
\begin{align*}
\gamma_{3,1} & \lesssim\delta\frac{\sigma_{r}^{\star}}{\sqrt{np}}\left(\omega_{j}^{\star}\left\Vert \bm{U}_{i,\cdot}^{\star}\right\Vert _{2}+\omega_{i}^{\star}\left\Vert \bm{U}_{j,\cdot}^{\star}\right\Vert _{2}\right)\lesssim\delta\widetilde{v}_{i,j}^{1/2},\\
\gamma_{3,2}+\gamma_{3,3} & \lesssim\delta\frac{\omega_{\min}}{\sqrt{n}p}\left(\omega_{j}^{\star}\left\Vert \bm{U}_{i,\cdot}^{\star}\right\Vert _{2}+\omega_{i}^{\star}\left\Vert \bm{U}_{j,\cdot}^{\star}\right\Vert _{2}\right)\lesssim\delta\widetilde{v}_{i,j}^{1/2},
\end{align*}
with the proviso that 
\begin{equation}
\left\Vert \bm{U}_{i,\cdot}^{\star}\right\Vert _{2}+\left\Vert \bm{U}_{j,\cdot}^{\star}\right\Vert _{2}\gtrsim\delta^{-1}\kappa r\log^{2}\left(n+d\right)\left[\frac{\omega_{\max}}{\sigma_{r}^{\star}}\sqrt{\frac{d}{np}}+\kappa_{\omega}\frac{\omega_{\max}^{2}}{p\sigma_{r}^{\star2}}\sqrt{\frac{d}{n}}+\frac{\kappa\mu r\kappa_{\omega}\log\left(n+d\right)}{\sqrt{ndp^{2}}}\right]\sqrt{\frac{1}{d}}.\label{eq:ce-normal-approx-1-inter-7}
\end{equation}
\item With regards to the term $\gamma_{4}$ defined in \eqref{eq:Phi-ij-gamma-1234},
we can see from (\ref{eq:ce-normal-approx-1-inter-6}) that 
\begin{align*}
\gamma_{4} & \lesssim\delta^{2}\sigma_{r}^{\star2}\left(\frac{\left\Vert \bm{U}_{i,\cdot}^{\star}\bm{\Sigma}^{\star}\right\Vert _{2}+\omega_{i}^{\star}}{\sqrt{\min\left\{ ndp^{2}\kappa,np\right\} }\sigma_{1}^{\star}}+\frac{\omega_{\min}\omega_{i}^{\star}}{\sqrt{np^{2}}\sigma_{1}^{\star2}}\right)\left(\frac{\left\Vert \bm{U}_{j,\cdot}^{\star}\bm{\Sigma}^{\star}\right\Vert _{2}+\omega_{j}^{\star}}{\sqrt{\min\left\{ ndp^{2}\kappa,np\right\} }\sigma_{1}^{\star}}+\frac{\omega_{\min}\omega_{j}^{\star}}{\sqrt{np^{2}}\sigma_{1}^{\star2}}\right)\\
 & \lesssim\underbrace{\delta^{2}\frac{\left\Vert \bm{U}_{i,\cdot}^{\star}\bm{\Sigma}^{\star}\right\Vert _{2}\left\Vert \bm{U}_{j,\cdot}^{\star}\bm{\Sigma}^{\star}\right\Vert _{2}}{\min\left\{ ndp^{2}\kappa,np\right\} \kappa}}_{\eqqcolon\gamma_{4,1}}+\underbrace{\delta^{2}\sigma_{1}^{\star}\frac{\omega_{i}^{\star}\left\Vert \bm{U}_{j,\cdot}^{\star}\right\Vert _{2}+\omega_{j}^{\star}\left\Vert \bm{U}_{i,\cdot}^{\star}\right\Vert _{2}}{\min\left\{ ndp^{2}\kappa,np\right\} \kappa}}_{\eqqcolon\gamma_{4,2}}+\underbrace{\delta^{2}\frac{\omega_{i}^{\star}\omega_{j}^{\star}}{\min\left\{ ndp^{2}\kappa,np\right\} \kappa}}_{\eqqcolon\gamma_{4,3}}\\
 & \quad+\underbrace{\delta^{2}\frac{\omega_{\min}\left(\omega_{i}^{\star}\left\Vert \bm{U}_{j,\cdot}^{\star}\right\Vert _{2}+\omega_{j}^{\star}\left\Vert \bm{U}_{i,\cdot}^{\star}\right\Vert _{2}\right)}{\sqrt{np^{2}}\sqrt{\min\left\{ ndp^{2}\kappa,np\right\} }\kappa}}_{\eqqcolon\gamma_{4,4}}+\underbrace{\delta^{2}\frac{\omega_{\min}\omega_{i}^{\star}\omega_{j}^{\star}}{\sqrt{np^{2}}\sqrt{\min\left\{ ndp^{2}\kappa,np\right\} }\sigma_{1}^{\star}\kappa}}_{\eqqcolon\gamma_{4,5}}+\underbrace{\delta^{2}\sigma_{r}^{\star2}\frac{\omega_{\min}^{2}\omega_{i}^{\star}\omega_{j}^{\star}}{np^{2}\sigma_{1}^{\star4}}}_{\eqqcolon\gamma_{4,6}}\\
 & \lesssim\delta\widetilde{v}_{i,j}^{1/2}.
\end{align*}
Here the last inequality follows from
\begin{align*}
\gamma_{4,1} & \lesssim\frac{\delta}{\sqrt{\min\left\{ ndp^{2}\kappa,np\right\} }}\left\Vert \bm{U}_{i,\cdot}^{\star}\bm{\Sigma}^{\star}\right\Vert _{2}\left\Vert \bm{U}_{j,\cdot}^{\star}\bm{\Sigma}^{\star}\right\Vert _{2}\lesssim\delta\widetilde{v}_{i,j}^{1/2},\\
\gamma_{4,2} & \lesssim\frac{\delta\sigma_{r}^{\star}}{\sqrt{ndp^{2}\wedge np}}\left(\omega_{j}^{\star}\left\Vert \bm{U}_{i,\cdot}^{\star}\right\Vert _{2}+\omega_{i}^{\star}\left\Vert \bm{U}_{j,\cdot}^{\star}\right\Vert _{2}\right)\lesssim\delta\widetilde{v}_{i,j}^{1/2},\\
\gamma_{4,3} & \lesssim\gamma_{3}\lesssim\delta\widetilde{v}_{i,j}^{1/2},\\
\gamma_{4,4} & \lesssim\delta\frac{\omega_{\min}}{\sqrt{n}p}\left(\omega_{j}^{\star}\left\Vert \bm{U}_{i,\cdot}^{\star}\right\Vert _{2}+\omega_{i}^{\star}\left\Vert \bm{U}_{j,\cdot}^{\star}\right\Vert _{2}\right)\lesssim\delta\widetilde{v}_{i,j}^{1/2},\\
\gamma_{4,5} & \lesssim\frac{\omega_{\min}}{\sqrt{np^{2}}\sigma_{1}^{\star}}\gamma_{4,3}\lesssim\left(\frac{\omega_{\min}^{2}}{p\sigma_{1}^{\star}}\sqrt{\frac{d}{n}}+\frac{1}{\sqrt{ndp^{2}}}\right)\gamma_{4,3}\lesssim\delta\widetilde{v}_{i,j}^{1/2},\\
\gamma_{4,6} & \lesssim\delta\frac{\omega_{\min}}{\sqrt{n}p}\left(\omega_{j}^{\star}\left\Vert \bm{U}_{i,\cdot}^{\star}\right\Vert _{2}+\omega_{i}^{\star}\left\Vert \bm{U}_{j,\cdot}^{\star}\right\Vert _{2}\right)\lesssim\delta\widetilde{v}_{i,j}^{1/2},
\end{align*}
provided that $\delta\lesssim1$, $ndp^{2}\gtrsim1$, $np\gtrsim1$,
\[
\frac{\omega_{\max}^{2}}{p\sigma_{r}^{\star}}\sqrt{\frac{d}{n}}\lesssim1
\]
and
\[
\left\Vert \bm{U}_{i,\cdot}^{\star}\right\Vert _{2}+\left\Vert \bm{U}_{j,\cdot}^{\star}\right\Vert _{2}\gtrsim\frac{\omega_{\max}^{2}}{p\sigma_{r}^{\star2}}\sqrt{\frac{d}{n}}\cdot\sqrt{\frac{1}{d}},
\]
which is guaranteed by (\ref{eq:ce-normal-approx-1-inter-7}).
\item We are left with $\gamma_{5}$, where we can utilize (\ref{eq:ce-normal-approx-1-inter-6})
to achieve
\begin{align*}
\gamma_{5} & \lesssim\theta\sigma_{1}^{\star}\left(\omega_{i}^{\star}\zeta_{\mathsf{2nd},j}+\omega_{j}^{\star}\zeta_{\mathsf{2nd},i}\right)\\
 & \lesssim\underbrace{\delta\theta\sigma_{r}^{\star}\frac{\omega_{i}^{\star}\left\Vert \bm{U}_{j,\cdot}^{\star}\right\Vert _{2}+\omega_{j}^{\star}\left\Vert \bm{U}_{i,\cdot}^{\star}\right\Vert _{2}}{\sqrt{\min\left\{ ndp^{2}\kappa,np\right\} }}}_{\eqqcolon\gamma_{5,1}}+\underbrace{\delta\theta\frac{\omega_{i}^{\star}\omega_{j}^{\star}}{\sqrt{\min\left\{ ndp^{2}\kappa,np\right\} \kappa}}}_{\eqqcolon\gamma_{5,2}}+\underbrace{\delta\theta\frac{\omega_{\min}\omega_{i}^{\star}\omega_{j}^{\star}}{\sqrt{np^{2}\kappa}\sigma_{1}^{\star}}}_{\eqqcolon\gamma_{5,3}}\\
 & \lesssim\delta\widetilde{v}_{i,j}^{1/2}.
\end{align*}
Here the last relation holds since
\begin{align*}
\gamma_{5,1}+\gamma_{5,2} & \lesssim\delta\frac{\sigma_{r}^{\star}}{\sqrt{ndp^{2}\wedge np}}\left(\omega_{j}^{\star}\left\Vert \bm{U}_{i,\cdot}^{\star}\right\Vert _{2}+\omega_{i}^{\star}\left\Vert \bm{U}_{j,\cdot}^{\star}\right\Vert _{2}\right)\lesssim\delta\widetilde{v}_{i,j}^{1/2},\\
\gamma_{5,3} & \lesssim\delta\frac{\omega_{\min}}{\sqrt{n}p}\left(\omega_{j}^{\star}\left\Vert \bm{U}_{i,\cdot}^{\star}\right\Vert _{2}+\omega_{i}^{\star}\left\Vert \bm{U}_{j,\cdot}^{\star}\right\Vert _{2}\right)\lesssim\delta\widetilde{v}_{i,j}^{1/2},
\end{align*}
provided that $\theta\lesssim1$ and 
\begin{align*}
\left\Vert \bm{U}_{i,\cdot}^{\star}\right\Vert _{2}+\left\Vert \bm{U}_{j,\cdot}^{\star}\right\Vert _{2} & \gtrsim\theta\frac{\omega_{\max}}{\sigma_{1}^{\star}}\asymp\left(\frac{\omega_{\max}}{\sigma_{1}^{\star}}\sqrt{\frac{d\kappa}{np}}+\frac{\omega_{\max}}{\sigma_{1}^{\star}}\sqrt{\frac{\kappa^{2}\mu r\log\left(n+d\right)}{np^{2}}}+\frac{\omega_{\max}^{2}}{p\sigma_{r}^{\star2}}\sqrt{\frac{d}{n}}\right)\sqrt{\frac{r\log^{2}\left(n+d\right)}{d}},
\end{align*}
which can be guaranteed by 
\[
\left\Vert \bm{U}_{i,\cdot}^{\star}\right\Vert _{2}+\left\Vert \bm{U}_{j,\cdot}^{\star}\right\Vert _{2}\gtrsim\left(\frac{\omega_{\max}}{\sigma_{1}^{\star}}\sqrt{\frac{d\kappa}{np}}+\frac{\omega_{\max}^{2}}{p\sigma_{r}^{\star2}}\sqrt{\frac{d}{n}}+\frac{\kappa^{2}\mu r\log\left(n+d\right)}{\sqrt{ndp^{2}}}\right)\sqrt{\frac{r\log^{2}\left(n+d\right)}{d}}
\]
where we have used (\ref{eq:theta-exact}) and a result from the AM-GM
inequality
\[
\frac{\omega_{\max}}{\sigma_{r}^{\star}}\sqrt{\frac{\kappa^{2}\mu r\log\left(n+d\right)}{np^{2}}}\lesssim\frac{\omega_{\max}^{2}}{\sigma_{r}^{\star2}}\sqrt{\frac{d}{np^{2}}}+\frac{\kappa^{2}\mu r\log\left(n+d\right)}{\sqrt{ndp^{2}}}.
\]
 
\end{itemize}
As a consequence, we have demonstrated that with probability exceeding
$1-O((n+d)^{-10})$, 
\[
\big(\widetilde{v}_{i,j}\big)^{-1/2}\zeta_{i,j}\lesssim\big(\widetilde{v}_{i,j}\big)^{-1/2}\left(\gamma_{1}+\gamma_{2}+\gamma_{3}+\gamma_{4}+\gamma_{5}\right)\lesssim\delta
\]
holds as long as $d\gtrsim\delta^{-2}\kappa^{8}\mu^{3}r^{3}\kappa_{\omega}^{2}\log^{4}(n+d)$,
\begin{align*}
ndp^{2}\gtrsim\delta^{-2}\kappa^{10}\mu^{4}r^{4}\kappa_{\omega}^{2}\log^{8}\left(n+d\right), & \qquad np\gtrsim\delta^{-2}\kappa^{10}\mu^{3}r^{3}\kappa_{\omega}^{2}\log^{6}\left(n+d\right),\\
\frac{\omega_{\max}^{2}}{p\sigma_{r}^{\star2}}\sqrt{\frac{d}{n}}\lesssim\frac{\delta}{\sqrt{\kappa^{7}\mu^{2}r^{2}\kappa_{\omega}\log^{6}\left(n+d\right)}}, & \qquad\frac{\omega_{\max}}{\sigma_{r}^{\star}}\sqrt{\frac{d}{np}}\lesssim\frac{\delta}{\sqrt{\kappa^{8}\mu^{2}r^{2}\kappa_{\omega}\log^{5}\left(n+d\right)}},
\end{align*}
and 
\[
\left\Vert \bm{U}_{i,\cdot}^{\star}\right\Vert _{2}+\left\Vert \bm{U}_{j,\cdot}^{\star}\right\Vert _{2}\gtrsim\delta^{-1}\kappa r\log^{2}\left(n+d\right)\left[\frac{\omega_{\max}}{\sigma_{r}^{\star}}\sqrt{\frac{d}{np}}+\kappa_{\omega}\frac{\omega_{\max}^{2}}{p\sigma_{r}^{\star2}}\sqrt{\frac{d}{n}}+\frac{\kappa\mu r\kappa_{\omega}\log\left(n+d\right)}{\sqrt{ndp^{2}}}\right]\sqrt{\frac{1}{d}}.
\]
Taking $\delta\asymp\log^{-1/2}(n+d)$ in the above bounds directly
establishes the advertised result. 

\subsubsection{Proof of Lemma \ref{lemma:ce-normal-approx-2}\label{appendix:ce-normal-approx-2}}

Define 
\[
a_{l}=n^{-1}(\bm{U}_{i,\cdot}^{\star}\bm{\Sigma}^{\star}\bm{f}_{l})(\bm{U}_{j,\cdot}^{\star}\bm{\Sigma}^{\star}\bm{f}_{l})
\]
for each $l=1,\ldots,n$. In view of the expression \eqref{eq:expression-Aij},
we can write 
\[
A_{i,j}=\sum_{l=1}^{n}\left(a_{l}-\mathbb{E}\left[a_{l}\right]\right).
\]
Apply Theorem \ref{thm:berry-esseen} to show that 
\[
\sup_{z\in\mathbb{R}}\left|\mathbb{P}\left(\overline{v}_{i,j}^{-1/2}A_{i,j}\leq z\right)-\Phi\left(z\right)\right|\lesssim\gamma,
\]
where $\overline{v}_{i,j}$ is defined in \eqref{eq:defn-bar-vij}
and 
\[
\gamma\coloneqq\overline{v}_{i,j}^{-3/2}\sum_{l=1}^{n}\mathbb{E}\left[\left|a_{l}^{3}\right|\right].
\]

It remains to control the term $\gamma$. Given that $\bm{U}_{i,\cdot}^{\star}\bm{\Sigma}^{\star}\bm{f}_{l}\sim\mathcal{N}(0,\Vert\bm{U}_{i,\cdot}^{\star}\bm{\Sigma}^{\star}\Vert_{2}^{2})$
and $\bm{U}_{j,\cdot}^{\star}\bm{\Sigma}^{\star}\bm{f}_{l}\sim\mathcal{N}(0,\Vert\bm{U}_{j,\cdot}^{\star}\bm{\Sigma}^{\star}\Vert_{2}^{2})$,
it is straightforward to check that for each $l\in[n]$, 
\begin{align*}
\mathbb{E}\left[\left|a_{l}^{3}\right|\right] & \leq\frac{1}{n^{3}}\mathbb{E}\left[\left(\bm{U}_{i,\cdot}^{\star}\bm{\Sigma}^{\star}\bm{f}_{l}\right)^{6}\right]^{1/2}\mathbb{E}\left[\left(\bm{U}_{j,\cdot}^{\star}\bm{\Sigma}^{\star}\bm{f}_{l}\right)^{6}\right]^{1/2}\lesssim\frac{1}{n^{3}}\left\Vert \bm{U}_{i,\cdot}^{\star}\bm{\Sigma}^{\star}\right\Vert _{2}^{3}\left\Vert \bm{U}_{j,\cdot}^{\star}\bm{\Sigma}^{\star}\right\Vert _{2}^{3}.
\end{align*}
Recognizing that 
\[
\overline{v}_{i,j}=\frac{1}{n}\left(S_{i,i}^{\star}S_{j,j}^{\star}+S_{i,j}^{\star2}\right)\geq\frac{1}{n}S_{i,i}^{\star}S_{j,j}^{\star}=\frac{1}{n}\left\Vert \bm{U}_{i,\cdot}^{\star}\bm{\Sigma}^{\star}\right\Vert _{2}^{2}\left\Vert \bm{U}_{j,\cdot}^{\star}\bm{\Sigma}^{\star}\right\Vert _{2}^{2},
\]
we can combine the above bounds to arrive at 
\[
\gamma=\overline{v}_{i,j}^{-3/2}\sum_{l=1}^{n}\mathbb{E}\left[\left|a_{l}^{3}\right|\right]\lesssim\frac{n^{-2}\left\Vert \bm{U}_{i,\cdot}^{\star}\bm{\Sigma}^{\star}\right\Vert _{2}^{3}\left\Vert \bm{U}_{j,\cdot}^{\star}\bm{\Sigma}^{\star}\right\Vert _{2}^{3}}{n^{-3/2}\left\Vert \bm{U}_{i,\cdot}^{\star}\bm{\Sigma}^{\star}\right\Vert _{2}^{3}\left\Vert \bm{U}_{j,\cdot}^{\star}\bm{\Sigma}^{\star}\right\Vert _{2}^{3}}\lesssim\frac{1}{\sqrt{n}}
\]
as claimed.

\subsubsection{Proof of Lemma \ref{lemma:ce-normal-approx-all}\label{appendix:ce-normal-approx-all}}

For any $z\in\mathbb{R}$, we can decompose 
\begin{align*}
 & \mathbb{P}\left(\left(S_{i,j}-S_{i,j}^{\star}\right)/\sqrt{v_{i,j}^{\star}}\leq z\right)=\mathbb{P}\left(A_{i,j}+W_{i,j}\leq\sqrt{v_{i,j}^{\star}}z\right)=\mathbb{E}\left[\mathbb{E}\left[\ind_{A_{i,j}+W_{i,j}\leq\sqrt{v_{i,j}^{\star}}z}\,\big|\,\bm{F}\right]\right]\\
 & \quad=\underbrace{\mathbb{E}\left[\mathbb{E}\left[\ind_{A_{i,j}+W_{i,j}\leq\sqrt{v_{i,j}^{\star}}z}\ind_{\mathcal{E}_{\mathsf{good}}}\,\big|\,\bm{F}\right]\right]}_{\eqqcolon\alpha_{1}}+\underbrace{\mathbb{E}\left[\mathbb{E}\left[\ind_{A_{i,j}+W_{i,j}\leq\sqrt{v_{i,j}^{\star}}z}\ind_{\mathcal{E}_{\mathsf{good}}^{\mathsf{c}}}\,\big|\,\bm{F}\right]\right]}_{\eqqcolon\alpha_{2}},
\end{align*}
leaving us with two terms to control. 
\begin{itemize}
\item Regarding the first term $\alpha_{1}$, we note that $\mathcal{E}_{\mathsf{good}}$
is $\sigma(\bm{F})$-measurable, and consequently, 
\begin{align*}
\mathbb{E}\left[\ind_{A_{i,j}+W_{i,j}\leq\sqrt{v_{i,j}^{\star}}z}\ind_{\mathcal{E}_{\mathsf{good}}}\,\big|\,\bm{F}\right] & =\ind_{\mathcal{E}_{\mathsf{good}}}\mathbb{E}\left[\ind_{A_{i,j}+W_{i,j}\leq\sqrt{v_{i,j}^{\star}}z}\,\big|\,\bm{F}\right]\\
 & =\ind_{\mathcal{E}_{\mathsf{good}}}\mathbb{P}\left(W_{i,j}\leq\sqrt{v_{i,j}^{\star}}z-A_{i,j}\,\big|\,\bm{F}\right).
\end{align*}
In view of Lemma \ref{lemma:ce-normal-approx-1}, we can see that
on the event $\mathcal{E}_{\mathsf{good}}$, 
\[
\sup_{t\in\mathbb{R}}\left|\mathbb{P}\left(\big(\widetilde{v}_{i,j}\big)^{-1/2}W_{i,j}\leq t\,\big|\,\bm{F}\right)-\Phi\left(t\right)\right|\lesssim\frac{1}{\sqrt{\log\left(n+d\right)}}.
\]
Since $A_{i,j}$ is $\sigma(\bm{F})$-measurable, by choosing $t=\big(\widetilde{v}_{i,j}\big)^{-1/2}(\sqrt{v_{i,j}^{\star}}z-A_{i,j})$
we have 
\[
\left|\mathbb{P}\left(W_{i,j}\leq\sqrt{v_{i,j}^{\star}}z-A_{i,j}\,\big|\,\bm{F}\right)-\Phi\left(\frac{\sqrt{v_{i,j}^{\star}}z-A_{i,j}}{\widetilde{v}_{i,j}^{1/2}}\right)\right|\ind_{\mathcal{E}_{\mathsf{good}}}\lesssim\frac{1}{\sqrt{\log\left(n+d\right)}}.
\]
This in turn leads to 
\begin{align*}
\alpha_{1} & =\mathbb{E}\left[\Phi\left(\frac{\sqrt{v_{i,j}^{\star}}z-A_{i,j}}{\widetilde{v}_{i,j}^{1/2}}\right)\ind_{\mathcal{E}_{\mathsf{good}}}\right]+O\left(\frac{1}{\sqrt{\log\left(n+d\right)}}\right)\\
 & =\mathbb{E}\left[\Phi\left(\frac{\sqrt{v_{i,j}^{\star}}z-A_{i,j}}{\widetilde{v}_{i,j}^{1/2}}\right)\right]-\mathbb{E}\left[\Phi\left(\frac{\sqrt{v_{i,j}^{\star}}z-A_{i,j}}{\widetilde{v}_{i,j}^{1/2}}\right)\ind_{\mathcal{E}_{\mathsf{good}}^{\mathsf{c}}}\right]+O\left(\frac{1}{\sqrt{\log\left(n+d\right)}}\right)\\
 & =\mathbb{E}\left[\Phi\left(\frac{\sqrt{v_{i,j}^{\star}}z-A_{i,j}}{\widetilde{v}_{i,j}^{1/2}}\right)\right]+O\left(\frac{1}{\sqrt{\log\left(n+d\right)}}\right),
\end{align*}
where the last identity is valid since 
\[
\mathbb{E}\left[\Phi\left(\frac{\sqrt{v_{i,j}^{\star}}z-A_{i,j}}{\widetilde{v}_{i,j}^{1/2}}\right)\ind_{\mathcal{E}_{\mathsf{good}}^{\mathsf{c}}}\right]\leq\mathbb{P}\left(\mathcal{E}_{\mathsf{good}}^{\mathsf{c}}\right)\lesssim\left(n+d\right)^{-10}\lesssim\frac{1}{\sqrt{\log\left(n+d\right)}}.
\]
Let $\phi(\cdot)$ denote the probability density of $\mathcal{N}(0,1)$,
then it is readily seen that 
\begin{align*}
\mathbb{E}\left[\Phi\left(\frac{\sqrt{v_{i,j}^{\star}}z-A_{i,j}}{\widetilde{v}_{i,j}^{1/2}}\right)\right] & =\mathbb{E}\left[\int_{-\infty}^{+\infty}\phi\left(t\right)\ind_{t\leq(\widetilde{v}_{i,j})^{-1/2}\left(\sqrt{v_{i,j}^{\star}}z-A_{i,j}\right)}\mathrm{d}t\right]\\
 & \overset{\text{(i)}}{=}\int_{-\infty}^{+\infty}\mathbb{E}\left[\phi\left(t\right)\ind_{t\leq(\widetilde{v}_{i,j})^{-1/2}\left(\sqrt{v_{i,j}^{\star}}z-A_{i,j}\right)}\right]\mathrm{d}t\\
 & =\int_{-\infty}^{+\infty}\phi\left(t\right)\mathbb{P}\left(t\leq\big(\widetilde{v}_{i,j}\big)^{-1/2}\left(\sqrt{v_{i,j}^{\star}}z-A_{i,j}\right)\right)\mathrm{d}t\\
 & =\int_{-\infty}^{+\infty}\phi\left(t\right)\mathbb{P}\left(A_{i,j}\leq\sqrt{v_{i,j}^{\star}}z-t\widetilde{v}_{i,j}^{1/2}\right)\mathrm{d}t\\
 & \overset{\text{(ii)}}{=}\int_{-\infty}^{+\infty}\phi\left(t\right)\Phi\left[\big(\overline{v}_{i,j}\big)^{-1/2}\left(\sqrt{v_{i,j}^{\star}}z-t\widetilde{v}_{i,j}^{1/2}\right)\right]\mathrm{d}t+O\left(\frac{1}{\sqrt{n}}\right).
\end{align*}
Here, (i) invokes Fubini's Theorem for nonnegative measurable functions,
whereas (ii) follows from Lemma \ref{lemma:ce-normal-approx-2}. Finally,
letting $U$ and $V$ be two independent $\mathcal{N}(0,1)$ random
variables, we can readily calculate 
\begin{align*}
\int_{-\infty}^{+\infty}\phi\left(t\right)\Phi\left(\frac{z\sqrt{v_{i,j}^{\star}}-t\sqrt{\widetilde{v}_{i,j}}}{\sqrt{\overline{v}_{i,j}}}\right)\mathrm{d}t & =\mathbb{P}\left(U\leq\frac{z\sqrt{v_{i,j}^{\star}}-V\sqrt{\widetilde{v}_{i,j}}}{\sqrt{\overline{v}_{i,j}}}\right)\\
 & =\mathbb{P}\left(U\sqrt{\overline{v}_{i,j}}+V\sqrt{\widetilde{v}_{i,j}}\leq z\sqrt{v_{i,j}^{\star}}\right)\\
 & =\Phi\left(z\right),
\end{align*}
where the last relation follows from the fact that $U\sqrt{\overline{v}_{i,j}}+V\sqrt{\widetilde{v}_{i,j}}\sim\mathcal{N}(0,v_{i,j}^{\star})$.
This allows us to conclude that 
\[
\alpha_{1}=\Phi\left(z\right)+O\left(\frac{1}{\sqrt{\log\left(n+d\right)}}+\frac{1}{\sqrt{n}}\right).
\]
\item We then move on to the other term $\alpha_{2}$. Towards this, it
is straightforward to derive that 
\[
\alpha_{2}\leq\mathbb{P}\left(\mathcal{E}_{\mathsf{good}}^{\mathsf{c}}\right)\lesssim\left(n+d\right)^{-10}.
\]
\end{itemize}
Note that the above analysis holds for all $z\in\mathbb{R}$. Therefore,
for any $z\in\mathbb{R}$, taking the above calculation together yields
\[
\left|\mathbb{P}\left(\left(S_{i,j}-S_{i,j}^{\star}\right)/\sqrt{v_{i,j}^{\star}}\leq z\right)-\Phi\left(z\right)\right|\leq\left|\alpha_{1}-\Phi\left(z\right)\right|+\alpha_{2}\lesssim\frac{1}{\sqrt{\log\left(n+d\right)}}+\frac{1}{\sqrt{n}}\asymp\frac{1}{\sqrt{\log\left(n+d\right)}}
\]
as claimed, provided that $n\gtrsim\log(n+d)$.

\subsection{Auxiliary lemmas for Theorem \ref{thm:ce-CI-complete}}

\subsubsection{Proof of Lemma \ref{lemma:ce-var-est}\label{appendix:proof-lemma-ce-var-est}}

In the sequel, we shall only consider the case when $i\neq j$; the
analysis for $i=j$ is similar and in fact simpler, and hence we omit
it here for brevity. Let us denote 
\begin{align*}
v_{i,j}^{\star} & =\underbrace{\frac{2-p}{np}S_{i,i}^{\star}S_{j,j}^{\star}}_{\eqqcolon\alpha_{1}}+\underbrace{\frac{4-3p}{np}S_{i,j}^{\star2}}_{\eqqcolon\alpha_{2}}+\underbrace{\frac{1}{np}\left(\omega_{i}^{\star2}S_{j,j}^{\star}+\omega_{j}^{\star2}S_{i,i}^{\star}\right)}_{\eqqcolon\alpha_{3}}\\
 & \quad+\underbrace{\frac{1}{np^{2}}\sum_{k=1}^{d}\left[\omega_{i}^{\star2}+\left(1-p\right)S_{i,i}^{\star}\right]\left[\omega_{k}^{\star2}+\left(1-p\right)S_{k,k}^{\star}\right]\left(\bm{U}_{k,\cdot}^{\star}\bm{U}_{j,\cdot}^{\star\top}\right)^{2}}_{\eqqcolon\alpha_{4}}+\underbrace{\frac{2\left(1-p\right)^{2}}{np^{2}}\sum_{k=1}^{d}S_{i,k}^{\star2}\left(\bm{U}_{k,\cdot}^{\star}\bm{U}_{j,\cdot}^{\star\top}\right)^{2}}_{\eqqcolon\alpha_{5}}\\
 & \quad+\underbrace{\frac{1}{np^{2}}\sum_{k=1}^{d}\left[\omega_{j}^{\star2}+\left(1-p\right)S_{j,j}^{\star}\right]\left[\omega_{k}^{\star2}+\left(1-p\right)S_{k,k}^{\star}\right]\left(\bm{U}_{k,\cdot}^{\star}\bm{U}_{i,\cdot}^{\star\top}\right)^{2}}_{\eqqcolon\alpha_{6}}+\underbrace{\frac{2\left(1-p\right)^{2}}{np^{2}}\sum_{k=1}^{d}S_{j,k}^{\star2}\left(\bm{U}_{k,\cdot}^{\star}\bm{U}_{i,\cdot}^{\star\top}\right)^{2}}_{\eqqcolon\alpha_{7}}
\end{align*}
and 
\begin{align*}
v_{i,j} & =\underbrace{\frac{2-p}{np}S_{i,i}S_{j,j}}_{\eqqcolon\beta_{1}}+\underbrace{\frac{4-3p}{np}S_{i,j}^{2}}_{\eqqcolon\beta_{2}}+\underbrace{\frac{1}{np}\left(\omega_{i}^{2}S_{j,j}+\omega_{j}^{2}S_{i,i}\right)}_{\eqqcolon\beta_{3}}\\
 & \quad+\underbrace{\frac{1}{np^{2}}\sum_{k=1}^{d}\left[\omega_{i}^{2}+\left(1-p\right)S_{i,i}\right]\left[\omega_{k}^{2}+\left(1-p\right)S_{k,k}\right]\left(\bm{U}_{k,\cdot}\bm{U}_{j,\cdot}^{\top}\right)^{2}}_{\eqqcolon\beta_{4}}+\underbrace{\frac{2\left(1-p\right)^{2}}{np^{2}}\sum_{k=1}^{d}S_{i,k}^{2}\left(\bm{U}_{k,\cdot}\bm{U}_{j,\cdot}^{\top}\right)^{2}}_{\eqqcolon\beta_{5}}\\
 & \quad+\underbrace{\frac{1}{np^{2}}\sum_{k=1}^{d}\left[\omega_{j}^{2}+\left(1-p\right)S_{j,j}\right]\left[\omega_{k}^{2}+\left(1-p\right)S_{k,k}\right]\left(\bm{U}_{k,\cdot}\bm{U}_{i,\cdot}^{\top}\right)^{2}}_{\eqqcolon\beta_{6}}+\underbrace{\frac{2\left(1-p\right)^{2}}{np^{2}}\sum_{k=1}^{d}S_{j,k}^{2}\left(\bm{U}_{k,\cdot}\bm{U}_{i,\cdot}^{\top}\right)^{2}}_{\eqqcolon\beta_{7}}.
\end{align*}
It follows from Lemma \ref{lemma:ce-variance-concentration} that
\[
v_{i,j}^{\star}\gtrsim\frac{1}{\min\left\{ ndp^{2}\kappa,np\right\} }\left\Vert \bm{U}_{i,\cdot}^{\star}\bm{\Sigma}^{\star}\right\Vert _{2}^{2}\left\Vert \bm{U}_{j,\cdot}^{\star}\bm{\Sigma}^{\star}\right\Vert _{2}^{2}+\left(\frac{\sigma_{r}^{\star2}}{\min\left\{ ndp^{2},np\right\} }+\frac{\omega_{\min}^{2}}{np^{2}}\right)\left(\omega_{j}^{\star2}\left\Vert \bm{U}_{i,\cdot}^{\star}\right\Vert _{2}^{2}+\omega_{i}^{\star2}\left\Vert \bm{U}_{j,\cdot}^{\star}\right\Vert _{2}^{2}\right).
\]
In view of the bounds on $\gamma_{2}$, $\gamma_{3}$, $\gamma_{4}$
and $\gamma_{5}$ in Step 3 of Appendix \ref{appendix:proof-ce-normal-approx-1}
as well as (\ref{eq:ce-normal-approx-1-inter-6}), we know that for
any $\varepsilon\in(0,1),$ one has 
\begin{equation}
\sigma_{1}^{\star2}\left(\left\Vert \bm{U}_{i,\cdot}^{\star}\right\Vert _{2}\zeta_{\mathsf{2nd},j}+\left\Vert \bm{U}_{j,\cdot}^{\star}\right\Vert _{2}\zeta_{\mathsf{2nd},i}\right)\lesssim\varepsilon\widetilde{v}_{i,j}^{1/2}\lesssim\varepsilon\left(v_{i,j}^{\star}\right)^{1/2},\label{eq:ce-var-est-inter-1}
\end{equation}
and \begin{subequations}\label{eq:ce-var-est-inter-2}
\begin{align}
\zeta_{\mathsf{2nd},i} & \lesssim\frac{\varepsilon}{\sqrt{\kappa}}\frac{\left\Vert \bm{U}_{i,\cdot}^{\star}\bm{\Sigma}^{\star}\right\Vert _{2}+\omega_{i}^{\star}}{\sqrt{\min\left\{ ndp^{2}\kappa,np\right\} }\sigma_{1}^{\star}}+\frac{\varepsilon}{\sqrt{\kappa}}\frac{\omega_{\min}\omega_{i}^{\star}}{\sqrt{np^{2}}\sigma_{1}^{\star2}},\\
\zeta_{\mathsf{2nd},j} & \lesssim\frac{\varepsilon}{\sqrt{\kappa}}\frac{\left\Vert \bm{U}_{j,\cdot}^{\star}\bm{\Sigma}^{\star}\right\Vert _{2}+\omega_{j}^{\star}}{\sqrt{\min\left\{ ndp^{2}\kappa,np\right\} }\sigma_{1}^{\star}}+\frac{\varepsilon}{\sqrt{\kappa}}\frac{\omega_{\min}\omega_{j}^{\star}}{\sqrt{np^{2}}\sigma_{1}^{\star2}},
\end{align}
\end{subequations}and
\begin{equation}
\theta^{2}\omega_{i}^{\star}\omega_{j}^{\star}+\sigma_{1}^{\star2}\zeta_{\mathsf{2nd},i}\zeta_{\mathsf{2nd},j}+\theta\sigma_{1}^{\star}\left(\omega_{i}^{\star}\zeta_{\mathsf{2nd},j}+\omega_{j}^{\star}\zeta_{\mathsf{2nd},i}\right)\lesssim\varepsilon\widetilde{v}_{i,j}^{1/2}\lesssim\varepsilon\left(v_{i,j}^{\star}\right)^{1/2},\label{eq:ce-var-est-inter-3}
\end{equation}
provided that the following conditions hold: $d\gtrsim\varepsilon^{-2}\kappa^{8}\mu^{3}r^{3}\kappa_{\omega}^{2}\log^{4}(n+d)$,
\begin{align*}
ndp^{2}\gtrsim\varepsilon^{-2}\kappa^{10}\mu^{4}r^{4}\kappa_{\omega}^{2}\log^{8}\left(n+d\right), & \qquad np\gtrsim\varepsilon^{-2}\kappa^{10}\mu^{3}r^{3}\kappa_{\omega}^{2}\log^{6}\left(n+d\right),\\
\frac{\omega_{\max}^{2}}{p\sigma_{r}^{\star2}}\sqrt{\frac{d}{n}}\lesssim\frac{\varepsilon}{\sqrt{\kappa^{7}\mu^{2}r^{2}\kappa_{\omega}\log^{6}\left(n+d\right)}}, & \qquad\frac{\omega_{\max}}{\sigma_{r}^{\star}}\sqrt{\frac{d}{np}}\lesssim\frac{\varepsilon}{\sqrt{\kappa^{8}\mu^{2}r^{2}\kappa_{\omega}\log^{5}\left(n+d\right)}},
\end{align*}
and 
\[
\left\Vert \bm{U}_{i,\cdot}^{\star}\right\Vert _{2}+\left\Vert \bm{U}_{j,\cdot}^{\star}\right\Vert _{2}\gtrsim\varepsilon^{-1}\kappa r\log^{2}\left(n+d\right)\left[\frac{\omega_{\max}}{\sigma_{r}^{\star}}\sqrt{\frac{d}{np}}+\kappa_{\omega}\frac{\omega_{\max}^{2}}{p\sigma_{r}^{\star2}}\sqrt{\frac{d}{n}}+\frac{\kappa\mu r\kappa_{\omega}\log\left(n+d\right)}{\sqrt{ndp^{2}}}\right]\sqrt{\frac{1}{d}}.
\]
These basic facts will be very useful for us to control the difference
between $v_{i,j}^{\star}$ and $v_{i,j}$, towards which we shall
bound $\alpha_{i}-\beta_{i}$, $1\leq i\leq7$, separately.

\paragraph{Step 1: bounding $\vert\alpha_{1}-\beta_{1}\vert$.}

Recall from Lemma \ref{lemma:pca-noise-level-est} that \begin{subequations}\label{eq:ce-var-est-inter-4}
\begin{align}
\left|S_{i,i}-S_{i,i}^{\star}\right| & \lesssim\left(\theta+\sqrt{\frac{\kappa^{3}r\log\left(n+d\right)}{n}}\right)\left\Vert \bm{U}_{i,\cdot}^{\star}\bm{\Sigma}^{\star}\right\Vert _{2}^{2}+\theta^{2}\omega_{i}^{\star2}+\left(\theta\omega_{i}^{\star}+\zeta_{\mathsf{2nd},i}\sigma_{1}^{\star}\right)\left\Vert \bm{U}_{i,\cdot}^{\star}\bm{\Sigma}^{\star}\right\Vert _{2}+\zeta_{\mathsf{2nd},i}^{2}\sigma_{1}^{\star2},\\
\left|S_{j,j}-S_{j,j}^{\star}\right| & \lesssim\left(\theta+\sqrt{\frac{\kappa^{3}r\log\left(n+d\right)}{n}}\right)\left\Vert \bm{U}_{j,\cdot}^{\star}\bm{\Sigma}^{\star}\right\Vert _{2}^{2}+\theta^{2}\omega_{j}^{\star2}+\left(\theta\omega_{j}^{\star}+\zeta_{\mathsf{2nd},j}\sigma_{1}^{\star}\right)\left\Vert \bm{U}_{j,\cdot}^{\star}\bm{\Sigma}^{\star}\right\Vert _{2}+\zeta_{\mathsf{2nd},j}^{2}\sigma_{1}^{\star2}.
\end{align}
\end{subequations}We proceed with the following elementary inequality:
\begin{align*}
\left|\alpha_{1}-\beta_{1}\right| & \lesssim\frac{1}{np}\left|S_{i,i}S_{j,j}-S_{i,i}^{\star}S_{j,j}^{\star}\right|\lesssim\frac{1}{np}S_{i.i}\left|S_{j,j}-S_{j,j}^{\star}\right|+\frac{1}{np}S_{j,j}^{\star}\left|S_{i,i}-S_{i,i}^{\star}\right|\\
 & \lesssim\underbrace{\frac{1}{np}S_{i.i}^{\star}\left|S_{j,j}-S_{j,j}^{\star}\right|+\frac{1}{np}S_{j,j}^{\star}\left|S_{i,i}-S_{i,i}^{\star}\right|}_{\eqqcolon\gamma_{1,1}}+\underbrace{\frac{1}{np}\left|S_{i,i}-S_{i,i}^{\star}\right|\left|S_{j,j}-S_{j,j}^{\star}\right|}_{\eqqcolon\gamma_{1,2}}.
\end{align*}

\begin{itemize}
\item Regarding $\gamma_{1,1}$, it is seen that 
\begin{align*}
\gamma_{1,1} & \lesssim\underbrace{\frac{1}{np}\left(\theta+\sqrt{\frac{\kappa^{3}r\log\left(n+d\right)}{n}}\right)\left\Vert \bm{U}_{i,\cdot}^{\star}\bm{\Sigma}^{\star}\right\Vert _{2}^{2}\left\Vert \bm{U}_{j,\cdot}^{\star}\bm{\Sigma}^{\star}\right\Vert _{2}^{2}}_{\eqqcolon\gamma_{1,1,1}}+\underbrace{\frac{1}{np}\theta^{2}\left(\omega_{j}^{\star2}\left\Vert \bm{U}_{i,\cdot}^{\star}\bm{\Sigma}^{\star}\right\Vert _{2}^{2}+\omega_{i}^{\star2}\left\Vert \bm{U}_{j,\cdot}^{\star}\bm{\Sigma}^{\star}\right\Vert _{2}^{2}\right)}_{\eqqcolon\gamma_{1,1,2}}\\
 & \quad+\underbrace{\frac{1}{np}\theta\left(\omega_{i}^{\star}\left\Vert \bm{U}_{i,\cdot}^{\star}\bm{\Sigma}^{\star}\right\Vert _{2}\left\Vert \bm{U}_{j,\cdot}^{\star}\bm{\Sigma}^{\star}\right\Vert _{2}^{2}+\omega_{j}^{\star}\left\Vert \bm{U}_{j,\cdot}^{\star}\bm{\Sigma}^{\star}\right\Vert _{2}\left\Vert \bm{U}_{i,\cdot}^{\star}\bm{\Sigma}^{\star}\right\Vert _{2}^{2}\right)}_{\eqqcolon\gamma_{1,1,3}}\\
 & \quad+\underbrace{\frac{1}{np}\sigma_{1}^{\star2}\left(\left\Vert \bm{U}_{i,\cdot}^{\star}\bm{\Sigma}^{\star}\right\Vert _{2}^{2}\zeta_{\mathsf{2nd},j}^{2}+\left\Vert \bm{U}_{j,\cdot}^{\star}\bm{\Sigma}^{\star}\right\Vert _{2}^{2}\zeta_{\mathsf{2nd},i}^{2}\right)}_{\eqqcolon\gamma_{1,1,4}}\\
 & \quad+\underbrace{\frac{1}{np}\sigma_{1}^{\star}\left\Vert \bm{U}_{i,\cdot}^{\star}\bm{\Sigma}^{\star}\right\Vert _{2}\left\Vert \bm{U}_{j,\cdot}^{\star}\bm{\Sigma}^{\star}\right\Vert _{2}\left(\left\Vert \bm{U}_{i,\cdot}^{\star}\bm{\Sigma}^{\star}\right\Vert _{2}\zeta_{\mathsf{2nd},j}+\left\Vert \bm{U}_{j,\cdot}^{\star}\bm{\Sigma}^{\star}\right\Vert _{2}\zeta_{\mathsf{2nd},i}\right)}_{\eqqcolon\gamma_{1,1,5}}\lesssim\delta v_{i,j}^{\star},
\end{align*}
where the last relation holds true since 
\begin{align*}
\gamma_{1,1,1} & \lesssim\frac{\delta}{np}\left\Vert \bm{U}_{i,\cdot}^{\star}\bm{\Sigma}^{\star}\right\Vert _{2}^{2}\left\Vert \bm{U}_{j,\cdot}^{\star}\bm{\Sigma}^{\star}\right\Vert _{2}^{2}\lesssim\delta v_{i,j}^{\star},\\
\gamma_{1,1,2} & \lesssim\delta\frac{\sigma_{r}^{\star2}}{np}\left(\omega_{j}^{\star2}\left\Vert \bm{U}_{i,\cdot}^{\star}\right\Vert _{2}^{2}+\omega_{i}^{\star2}\left\Vert \bm{U}_{j,\cdot}^{\star}\right\Vert _{2}^{2}\right)\lesssim\delta v_{i,j}^{\star},\\
\gamma_{1,1,3} & \lesssim\frac{\delta}{\sqrt{\kappa}}\frac{1}{np}\left(\omega_{i}^{\star}\left\Vert \bm{U}_{i,\cdot}^{\star}\bm{\Sigma}^{\star}\right\Vert _{2}\left\Vert \bm{U}_{j,\cdot}^{\star}\bm{\Sigma}^{\star}\right\Vert _{2}^{2}+\omega_{j}^{\star}\left\Vert \bm{U}_{j,\cdot}^{\star}\bm{\Sigma}^{\star}\right\Vert _{2}\left\Vert \bm{U}_{i,\cdot}^{\star}\bm{\Sigma}^{\star}\right\Vert _{2}^{2}\right)\\
 & \overset{\text{(i)}}{\lesssim}\delta\frac{1}{np}\left\Vert \bm{U}_{i,\cdot}^{\star}\bm{\Sigma}^{\star}\right\Vert _{2}^{2}\left\Vert \bm{U}_{j,\cdot}^{\star}\bm{\Sigma}^{\star}\right\Vert _{2}^{2}+\delta\frac{\sigma_{r}^{\star2}}{np}\left(\omega_{j}^{\star2}\left\Vert \bm{U}_{i,\cdot}^{\star}\right\Vert _{2}^{2}+\omega_{i}^{\star2}\left\Vert \bm{U}_{j,\cdot}^{\star}\right\Vert _{2}^{2}\right)\lesssim\delta v_{i,j}^{\star},\\
\gamma_{1,1,4} & \lesssim\frac{1}{np}\sigma_{1}^{\star4}\left(\left\Vert \bm{U}_{i,\cdot}^{\star}\right\Vert _{2}^{2}\zeta_{\mathsf{2nd},j}^{2}+\left\Vert \bm{U}_{j,\cdot}^{\star}\right\Vert _{2}^{2}\zeta_{\mathsf{2nd},i}^{2}\right)\overset{\text{(ii)}}{\lesssim}\frac{\varepsilon^{2}}{np}v_{i,j}^{\star}\lesssim\delta v_{i,j}^{\star},\\
\gamma_{1,1,5} & \lesssim\frac{1}{np}\left\Vert \bm{U}_{i,\cdot}^{\star}\bm{\Sigma}^{\star}\right\Vert _{2}\left\Vert \bm{U}_{j,\cdot}^{\star}\bm{\Sigma}^{\star}\right\Vert _{2}\cdot\sigma_{1}^{\star2}\left(\left\Vert \bm{U}_{i,\cdot}^{\star}\right\Vert _{2}\zeta_{\mathsf{2nd},j}+\left\Vert \bm{U}_{j,\cdot}^{\star}\right\Vert _{2}\zeta_{\mathsf{2nd},i}\right)\overset{\text{(iii)}}{\lesssim}\frac{\varepsilon}{\sqrt{np}}v_{i,j}^{\star}\lesssim\delta v_{i,j}^{\star},
\end{align*}
provided that $\theta\lesssim\delta/\sqrt{\kappa}$, $n\gtrsim\delta^{-2}\kappa^{3}r\log(n+d)$,
$np\gtrsim\delta^{-2}$ and $\varepsilon\lesssim1$. Here, the relation
(i) invokes the AM-GM inequality, while (ii) and (iii) make use of
(\ref{eq:ce-var-est-inter-1}). 
\item Regarding $\gamma_{1,2}$, we make the observation that 
\begin{align*}
\gamma_{1,2} & \lesssim\left(\theta+\sqrt{\frac{\kappa^{3}r\log\left(n+d\right)}{n}}\right)\gamma_{1,1}+\frac{1}{np}\left[\theta^{2}\omega_{i}^{\star2}+\left(\theta\omega_{i}^{\star}+\zeta_{\mathsf{2nd},i}\sigma_{1}^{\star}\right)\left\Vert \bm{U}_{i,\cdot}^{\star}\bm{\Sigma}^{\star}\right\Vert _{2}+\zeta_{\mathsf{2nd},i}^{2}\sigma_{1}^{\star2}\right]\\
 & \qquad\qquad\qquad\qquad\qquad\qquad\qquad\cdot\left[\theta^{2}\omega_{j}^{\star2}+\left(\theta\omega_{j}^{\star}+\zeta_{\mathsf{2nd},j}\sigma_{1}^{\star}\right)\left\Vert \bm{U}_{j,\cdot}^{\star}\bm{\Sigma}^{\star}\right\Vert _{2}+\zeta_{\mathsf{2nd},j}^{2}\sigma_{1}^{\star2}\right]\\
 & \lesssim\delta v_{i,j}^{\star}+\underbrace{\frac{1}{np}\theta^{4}\omega_{i}^{\star2}\omega_{j}^{\star2}}_{\eqqcolon\gamma_{1,2,1}}+\underbrace{\frac{1}{np}\left[\theta^{2}\omega_{i}^{\star}\omega_{j}^{\star}+\theta\sigma_{1}^{\star}\left(\omega_{i}^{\star}\zeta_{\mathsf{2nd},j}+\omega_{j}^{\star}\zeta_{\mathsf{2nd},i}\right)+\sigma_{1}^{\star2}\zeta_{\mathsf{2nd},i}\zeta_{\mathsf{2nd},j}\right]\left\Vert \bm{U}_{i,\cdot}^{\star}\bm{\Sigma}^{\star}\right\Vert _{2}\left\Vert \bm{U}_{j,\cdot}^{\star}\bm{\Sigma}^{\star}\right\Vert _{2}}_{\eqqcolon\gamma_{1,2,2}}\\
 & \quad+\underbrace{\frac{1}{np}\zeta_{\mathsf{2nd},i}^{2}\zeta_{\mathsf{2nd},j}^{2}\sigma_{1}^{\star4}}_{\eqqcolon\gamma_{1,2,3}}+\underbrace{\frac{1}{np}\theta^{3}\omega_{i}^{\star}\omega_{j}^{\star}\left(\omega_{j}^{\star}\left\Vert \bm{U}_{i,\cdot}^{\star}\bm{\Sigma}^{\star}\right\Vert _{2}+\omega_{i}^{\star}\left\Vert \bm{U}_{j,\cdot}^{\star}\bm{\Sigma}^{\star}\right\Vert _{2}\right)}_{\eqqcolon\gamma_{1,2,4}}\\
 & \quad+\underbrace{\frac{1}{np}\zeta_{\mathsf{2nd},i}\zeta_{\mathsf{2nd},j}\sigma_{1}^{\star4}\left(\left\Vert \bm{U}_{i,\cdot}^{\star}\right\Vert _{2}\zeta_{\mathsf{2nd},j}+\left\Vert \bm{U}_{j,\cdot}^{\star}\right\Vert _{2}\zeta_{\mathsf{2nd},i}\right)}_{\eqqcolon\gamma_{1,2,5}}\\
 & \quad+\underbrace{\frac{1}{np}\theta^{2}\sigma_{1}^{\star}\left(\omega_{j}^{\star2}\zeta_{\mathsf{2nd},i}\left\Vert \bm{U}_{i,\cdot}^{\star}\bm{\Sigma}^{\star}\right\Vert _{2}+\omega_{i}^{\star2}\zeta_{\mathsf{2nd},j}\left\Vert \bm{U}_{j,\cdot}^{\star}\bm{\Sigma}^{\star}\right\Vert _{2}\right)}_{\eqqcolon\gamma_{1,2,6}}+\underbrace{\frac{1}{np}\theta^{2}\sigma_{1}^{\star2}\left(\omega_{i}^{\star2}\zeta_{\mathsf{2nd},j}^{2}+\omega_{j}^{\star2}\zeta_{\mathsf{2nd},i}^{2}\right)}_{\eqqcolon\gamma_{1,2,7}}\\
 & \quad+\underbrace{\frac{1}{np}\theta\sigma_{1}^{\star2}\left(\omega_{i}^{\star}\left\Vert \bm{U}_{i,\cdot}^{\star}\bm{\Sigma}^{\star}\right\Vert _{2}\zeta_{\mathsf{2nd},j}^{2}+\omega_{j}^{\star}\left\Vert \bm{U}_{j,\cdot}^{\star}\bm{\Sigma}^{\star}\right\Vert _{2}\zeta_{\mathsf{2nd},i}^{2}\right)}_{\eqqcolon\gamma_{1,2,8}}\lesssim\delta v_{i,j}^{\star}.
\end{align*}
Here, the penultimate relation follows from the previous bound on
$\gamma_{1,1}$ as well as the assumptions that $\theta\lesssim1$
and $n\gtrsim\kappa^{3}r\log(n+d)$; the last line is valid since
\begin{align*}
\gamma_{1,2,1}+\gamma_{1,2,3} & \overset{\text{(i)}}{\lesssim}\frac{\varepsilon^{2}}{np}v_{i,j}^{\star}\leq\delta v_{i,j}^{\star},\\
\gamma_{1,2,2} & \overset{\text{(ii)}}{\lesssim}\frac{1}{np}\left\Vert \bm{U}_{i,\cdot}^{\star}\bm{\Sigma}^{\star}\right\Vert _{2}\left\Vert \bm{U}_{j,\cdot}^{\star}\bm{\Sigma}^{\star}\right\Vert _{2}\varepsilon\left(v_{i,j}^{\star}\right)^{1/2}\lesssim\frac{\varepsilon}{\sqrt{np}}v_{i,j}^{\star}\lesssim\delta v_{i,j}^{\star},\\
\gamma_{1,2,4} & \lesssim\sqrt{\frac{\kappa}{np}}\theta\cdot\theta^{2}\omega_{i}^{\star}\omega_{j}^{\star}\cdot\frac{\sigma_{r}^{\star}}{\sqrt{np}}\left(\omega_{j}^{\star}\left\Vert \bm{U}_{i,\cdot}^{\star}\right\Vert _{2}+\omega_{i}^{\star}\left\Vert \bm{U}_{j,\cdot}^{\star}\right\Vert _{2}\right)\overset{\text{(iii)}}{\lesssim}\sqrt{\frac{\kappa}{np}}\theta\varepsilon v_{i,j}^{\star}\lesssim\delta v_{i,j}^{\star},\\
\gamma_{1,2,5} & \lesssim\frac{1}{np}\cdot\zeta_{\mathsf{2nd},i}\zeta_{\mathsf{2nd},j}\sigma_{1}^{\star2}\cdot\sigma_{1}^{\star2}\left(\left\Vert \bm{U}_{i,\cdot}^{\star}\right\Vert _{2}\zeta_{\mathsf{2nd},j}+\left\Vert \bm{U}_{j,\cdot}^{\star}\right\Vert _{2}\zeta_{\mathsf{2nd},i}\right)\overset{\text{(iv)}}{\lesssim}\frac{\varepsilon^{2}}{np}v_{i,j}^{\star}\lesssim\delta v_{i,j}^{\star},\\
\gamma_{1,2,6} & \lesssim\sqrt{\frac{\kappa}{np}}\theta\cdot\theta\sigma_{1}^{\star}\left(\omega_{i}^{\star}\zeta_{\mathsf{2nd},j}+\omega_{j}^{\star}\zeta_{\mathsf{2nd},i}\right)\cdot\frac{\sigma_{r}^{\star}}{\sqrt{np}}\left(\omega_{j}^{\star}\left\Vert \bm{U}_{i,\cdot}^{\star}\right\Vert _{2}+\omega_{i}^{\star}\left\Vert \bm{U}_{j,\cdot}^{\star}\right\Vert _{2}\right)\overset{\text{(v)}}{\lesssim}\sqrt{\frac{\kappa}{np}}\theta\varepsilon v_{i,j}^{\star}\lesssim\delta v_{i,j}^{\star},\\
\gamma_{1,2,7} & \lesssim\frac{1}{np}\left[\sigma_{1}^{\star}\left(\omega_{i}^{\star}\zeta_{\mathsf{2nd},j}+\omega_{j}^{\star}\zeta_{\mathsf{2nd},i}\right)\right]^{2}\overset{\text{(vi)}}{\lesssim}\frac{\varepsilon}{np}v_{i,j}^{\star}\lesssim\delta v_{i,j}^{\star},\\
\gamma_{1,2,8} & \lesssim\frac{1}{np}\cdot\theta\sigma_{1}^{\star}\left(\omega_{i}^{\star}\zeta_{\mathsf{2nd},j}+\omega_{j}^{\star}\zeta_{\mathsf{2nd},i}\right)\cdot\sigma_{1}^{\star2}\left(\left\Vert \bm{U}_{i,\cdot}^{\star}\right\Vert _{2}\zeta_{\mathsf{2nd},j}+\left\Vert \bm{U}_{j,\cdot}^{\star}\right\Vert _{2}\zeta_{\mathsf{2nd},i}\right)\overset{\text{(vii)}}{\lesssim}\frac{\varepsilon^{2}}{np}v_{i,j}^{\star}\lesssim\delta v_{i,j}^{\star},
\end{align*}
provided that $np\gtrsim\delta^{-2}$, $\varepsilon\lesssim1$ and
$\theta\lesssim\delta/\sqrt{\kappa}$. Here, the inequalities (i)-(vii)
rely on (\ref{eq:ce-var-est-inter-1}) and (\ref{eq:ce-var-est-inter-3}). 
\end{itemize}
Taking the above bounds on $\gamma_{1,1}$ and $\gamma_{1,2}$ collectively,
we arrive at 
\[
\left|\alpha_{1}-\beta_{1}\right|\lesssim\left|\gamma_{1,1}\right|+\left|\gamma_{1,2}\right|\lesssim\delta v_{i,j}^{\star},
\]
provided that $\theta\lesssim\delta/\sqrt{\kappa}$, $n\gtrsim\delta^{-2}\kappa^{3}r\log(n+d)$,
$\varepsilon\lesssim1$, and $np\gtrsim\delta^{-2}$.

\paragraph{Step 2: bounding $\vert\alpha_{2}-\beta_{2}\vert$.}

It comes from Lemma \ref{lemma:pca-noise-level-est} that 
\begin{align*}
\left|S_{i,j}-S_{i,j}^{\star}\right| & \lesssim\left(\theta+\sqrt{\frac{\kappa^{3}r\log\left(n+d\right)}{n}}\right)\left\Vert \bm{U}_{i,\cdot}^{\star}\bm{\Sigma}^{\star}\right\Vert _{2}\left\Vert \bm{U}_{j,\cdot}^{\star}\bm{\Sigma}^{\star}\right\Vert _{2}+\theta\left(\omega_{i}^{\star}\left\Vert \bm{U}_{j,\cdot}^{\star}\bm{\Sigma}^{\star}\right\Vert _{2}+\omega_{j}^{\star}\left\Vert \bm{U}_{i,\cdot}^{\star}\bm{\Sigma}^{\star}\right\Vert _{2}\right)\\
 & \quad+\sigma_{1}^{\star}\left(\zeta_{\mathsf{2nd},i}\left\Vert \bm{U}_{j,\cdot}^{\star}\bm{\Sigma}^{\star}\right\Vert _{2}+\zeta_{\mathsf{2nd},j}\left\Vert \bm{U}_{i,\cdot}^{\star}\bm{\Sigma}^{\star}\right\Vert _{2}\right)+\theta^{2}\omega_{i}^{\star}\omega_{j}^{\star}+\zeta_{\mathsf{2nd},i}\zeta_{\mathsf{2nd},j}\sigma_{1}^{\star2}.
\end{align*}
With this in place, we shall proceed to decompose $\left|\alpha_{2}-\beta_{2}\right|$
as follows: 
\begin{align*}
\left|\alpha_{2}-\beta_{2}\right| & \lesssim\frac{1}{np}\left|S_{i,j}^{\star2}-S_{i,j}^{2}\right|\lesssim\underbrace{\frac{1}{np}S_{i,j}^{\star}\left|S_{i,j}-S_{i,j}^{\star}\right|}_{\eqqcolon\gamma_{2,1}}+\underbrace{\frac{1}{np}\left|S_{i,j}-S_{i,j}^{\star}\right|^{2}}_{\eqqcolon\gamma_{2,2}}.
\end{align*}

\begin{itemize}
\item With regards to $\gamma_{2,1}$, it is observed that 
\begin{align*}
\gamma_{2,1} & \lesssim\underbrace{\frac{1}{np}\left(\theta+\sqrt{\frac{\kappa^{3}r\log\left(n+d\right)}{n}}\right)\left\Vert \bm{U}_{i,\cdot}^{\star}\bm{\Sigma}^{\star}\right\Vert _{2}^{2}\left\Vert \bm{U}_{j,\cdot}^{\star}\bm{\Sigma}^{\star}\right\Vert _{2}^{2}}_{\eqqcolon\gamma_{2,1,1}}+\underbrace{\frac{1}{np}\theta^{2}\omega_{i}^{\star}\omega_{j}^{\star}\left\Vert \bm{U}_{i,\cdot}^{\star}\bm{\Sigma}^{\star}\right\Vert _{2}\left\Vert \bm{U}_{j,\cdot}^{\star}\bm{\Sigma}^{\star}\right\Vert _{2}}_{\eqqcolon\gamma_{2,1,2}}\\
 & \quad+\underbrace{\frac{1}{np}\theta\left\Vert \bm{U}_{i,\cdot}^{\star}\bm{\Sigma}^{\star}\right\Vert _{2}\left\Vert \bm{U}_{j,\cdot}^{\star}\bm{\Sigma}^{\star}\right\Vert _{2}\left(\omega_{i}^{\star}\left\Vert \bm{U}_{j,\cdot}^{\star}\bm{\Sigma}^{\star}\right\Vert _{2}+\omega_{j}^{\star}\left\Vert \bm{U}_{i,\cdot}^{\star}\bm{\Sigma}^{\star}\right\Vert _{2}\right)}_{\eqqcolon\gamma_{2,1,3}}\\
 & \quad+\underbrace{\frac{1}{np}\sigma_{1}^{\star}\left\Vert \bm{U}_{i,\cdot}^{\star}\bm{\Sigma}^{\star}\right\Vert _{2}\left\Vert \bm{U}_{j,\cdot}^{\star}\bm{\Sigma}^{\star}\right\Vert _{2}\left(\zeta_{\mathsf{2nd},i}\left\Vert \bm{U}_{j,\cdot}^{\star}\bm{\Sigma}^{\star}\right\Vert _{2}+\zeta_{\mathsf{2nd},j}\left\Vert \bm{U}_{i,\cdot}^{\star}\bm{\Sigma}^{\star}\right\Vert _{2}\right)}_{\eqqcolon\gamma_{2,1,4}}\\
 & \quad+\underbrace{\frac{1}{np}\zeta_{\mathsf{2nd},i}\zeta_{\mathsf{2nd},j}\sigma_{1}^{\star2}\left\Vert \bm{U}_{i,\cdot}^{\star}\bm{\Sigma}^{\star}\right\Vert _{2}\left\Vert \bm{U}_{j,\cdot}^{\star}\bm{\Sigma}^{\star}\right\Vert _{2}}_{\eqqcolon\gamma_{2,1,5}}\lesssim\delta v_{i,j}^{\star}.
\end{align*}
Here, the last relation holds due to the following bounds 
\begin{align*}
\gamma_{2,1,1} & \lesssim\frac{\delta}{np}\left\Vert \bm{U}_{i,\cdot}^{\star}\bm{\Sigma}^{\star}\right\Vert _{2}^{2}\left\Vert \bm{U}_{j,\cdot}^{\star}\bm{\Sigma}^{\star}\right\Vert _{2}^{2}\lesssim\delta v_{i,j}^{\star},\\
\gamma_{2,1,2} & \lesssim\frac{1}{\sqrt{np}}\cdot\theta^{2}\omega_{i}^{\star}\omega_{j}^{\star}\cdot\frac{1}{\sqrt{np}}\left\Vert \bm{U}_{i,\cdot}^{\star}\bm{\Sigma}^{\star}\right\Vert _{2}\left\Vert \bm{U}_{j,\cdot}^{\star}\bm{\Sigma}^{\star}\right\Vert _{2}\overset{\text{(i)}}{\lesssim}\frac{\varepsilon}{\sqrt{np}}v_{i,j}^{\star}\lesssim\delta v_{i,j}^{\star},\\
\gamma_{2,1,3} & \lesssim\sqrt{\kappa}\theta\cdot\frac{1}{\sqrt{np}}\left\Vert \bm{U}_{i,\cdot}^{\star}\bm{\Sigma}^{\star}\right\Vert _{2}\left\Vert \bm{U}_{j,\cdot}^{\star}\bm{\Sigma}^{\star}\right\Vert _{2}\cdot\frac{\sigma_{r}^{\star}}{\sqrt{np}}\left(\omega_{i}^{\star}\left\Vert \bm{U}_{j,\cdot}^{\star}\bm{\Sigma}^{\star}\right\Vert _{2}+\omega_{j}^{\star}\left\Vert \bm{U}_{i,\cdot}^{\star}\bm{\Sigma}^{\star}\right\Vert _{2}\right)\\
 & \lesssim\sqrt{\kappa}\theta v_{i,j}^{\star}\lesssim\delta v_{i,j}^{\star},\\
\gamma_{2,1,4} & \lesssim\frac{1}{\sqrt{np}}\cdot\frac{1}{\sqrt{np}}\left\Vert \bm{U}_{i,\cdot}^{\star}\bm{\Sigma}^{\star}\right\Vert _{2}\left\Vert \bm{U}_{j,\cdot}^{\star}\bm{\Sigma}^{\star}\right\Vert _{2}\cdot\sigma_{1}^{\star2}\left(\left\Vert \bm{U}_{i,\cdot}^{\star}\right\Vert _{2}\zeta_{\mathsf{2nd},j}+\left\Vert \bm{U}_{j,\cdot}^{\star}\right\Vert _{2}\zeta_{\mathsf{2nd},i}\right)\\
 & \overset{\text{(ii)}}{\lesssim}\frac{\varepsilon}{\sqrt{np}}v_{i,j}^{\star}\lesssim\delta v_{i,j}^{\star},\\
\gamma_{2,1,5} & \lesssim\frac{1}{np}\left\Vert \bm{U}_{i,\cdot}^{\star}\bm{\Sigma}^{\star}\right\Vert _{2}\left\Vert \bm{U}_{j,\cdot}^{\star}\bm{\Sigma}^{\star}\right\Vert _{2}\cdot\sigma_{1}^{\star2}\zeta_{\mathsf{2nd},i}\zeta_{\mathsf{2nd},j}\overset{\text{(iii)}}{\lesssim}\frac{\varepsilon}{\sqrt{np}}v_{i,j}^{\star}\lesssim\delta v_{i,j}^{\star},
\end{align*}
with the proviso that $\theta\lesssim\delta/\sqrt{\kappa}$, $n\gtrsim\delta^{-2}\kappa^{3}r\log(n+d)$,
$\varepsilon\lesssim1$ and $np\gtrsim\delta^{-2}$. Note that (i)-(iii)
follow from (\ref{eq:ce-var-est-inter-1}) and (\ref{eq:ce-var-est-inter-3}). 
\item Regarding $\gamma_{2,2}$, we obtain 
\begin{align*}
\gamma_{2,2} & \lesssim\underbrace{\frac{1}{np}\left(\theta^{2}+\frac{\kappa^{3}r\log\left(n+d\right)}{n}\right)\left\Vert \bm{U}_{i,\cdot}^{\star}\bm{\Sigma}^{\star}\right\Vert _{2}^{2}\left\Vert \bm{U}_{j,\cdot}^{\star}\bm{\Sigma}^{\star}\right\Vert _{2}^{2}}_{\eqqcolon\gamma_{2,2,1}}+\underbrace{\frac{1}{np}\theta^{4}\omega_{i}^{\star2}\omega_{j}^{\star2}}_{\eqqcolon\gamma_{2,2,2}}+\underbrace{\frac{1}{np}\zeta_{\mathsf{2nd},i}^{2}\zeta_{\mathsf{2nd},j}^{2}\sigma_{1}^{\star4}}_{\eqqcolon\gamma_{2,2,3}}\\
 & \quad+\underbrace{\frac{1}{np}\theta^{2}\left(\omega_{i}^{\star}\left\Vert \bm{U}_{j,\cdot}^{\star}\bm{\Sigma}^{\star}\right\Vert _{2}+\omega_{j}^{\star}\left\Vert \bm{U}_{i,\cdot}^{\star}\bm{\Sigma}^{\star}\right\Vert _{2}\right)^{2}}_{\eqqcolon\gamma_{2,2,4}}+\underbrace{\frac{1}{np}\sigma_{1}^{\star2}\left(\zeta_{\mathsf{2nd},i}\left\Vert \bm{U}_{j,\cdot}^{\star}\bm{\Sigma}^{\star}\right\Vert _{2}+\zeta_{\mathsf{2nd},j}\left\Vert \bm{U}_{i,\cdot}^{\star}\bm{\Sigma}^{\star}\right\Vert _{2}\right)^{2}}_{\eqqcolon\gamma_{2,2,5}}\\
 & \lesssim\delta v_{i,j}^{\star},
\end{align*}
where the last inequality holds true since 
\begin{align*}
\gamma_{2,2,1} & \lesssim\left(\theta+\sqrt{\frac{\kappa^{3}r\log\left(n+d\right)}{n}}\right)\gamma_{2,1,1}\lesssim\delta v_{i,j}^{\star},\\
\gamma_{2,2,2}+\gamma_{2,2,3} & \overset{\text{(i)}}{\lesssim}\frac{\varepsilon^{2}}{np}v_{i,j}^{\star}\lesssim\delta v_{i,j}^{\star},\\
\gamma_{2,2,4} & \lesssim\theta^{2}\kappa\cdot\frac{\sigma_{r}^{\star2}}{np}\left(\omega_{j}^{\star2}\left\Vert \bm{U}_{i,\cdot}^{\star}\right\Vert _{2}^{2}+\omega_{i}^{\star2}\left\Vert \bm{U}_{j,\cdot}^{\star}\right\Vert _{2}^{2}\right)\lesssim\theta^{2}\kappa v_{i,j}^{\star}\lesssim\delta v_{i,j}^{\star},\\
\gamma_{2,2,5} & \lesssim\frac{1}{np}\sigma_{1}^{\star4}\left(\left\Vert \bm{U}_{i,\cdot}^{\star}\right\Vert _{2}\zeta_{\mathsf{2nd},j}+\left\Vert \bm{U}_{j,\cdot}^{\star}\right\Vert _{2}\zeta_{\mathsf{2nd},i}\right)^{2}\overset{\text{(ii)}}{\lesssim}\frac{\varepsilon^{2}}{np}v_{i,j}^{\star}\lesssim\delta v_{i,j}^{\star},
\end{align*}
provided that $\theta\lesssim\delta/\sqrt{\kappa}$, $n\gtrsim\kappa^{3}r\log(n+d)$,
$\varepsilon\lesssim1$ and $np\geq\delta^{-1}$. Here, (i) follows
from (\ref{eq:ce-var-est-inter-3}), whereas (ii) follows from (\ref{eq:ce-var-est-inter-1}). 
\end{itemize}
Taking the above bounds on $\gamma_{2,1}$ and $\gamma_{2,2}$ together
yields 
\[
\left|\alpha_{2}-\beta_{2}\right|\lesssim\left|\gamma_{2,1}\right|+\left|\gamma_{2,2}\right|\lesssim\delta v_{i,j}^{\star},
\]
as long as $\theta\lesssim\delta/\sqrt{\kappa}$, $n\gtrsim\delta^{-2}\kappa^{3}r\log(n+d)$,
$\varepsilon\lesssim1$ and $np\gtrsim\delta^{-2}$.

\paragraph{Step 3: bounding $\vert\alpha_{3}-\beta_{3}\vert$.}

For each $l\in[d]$, let us first define $\mathsf{UB}_{l}$ as follows
\[
\mathsf{UB}_{l}\coloneqq\left(\theta+\sqrt{\frac{\kappa^{3}r\log\left(n+d\right)}{n}}\right)\left\Vert \bm{U}_{l,\cdot}^{\star}\bm{\Sigma}^{\star}\right\Vert _{2}^{2}+\theta^{2}\omega_{l}^{\star2}+\left(\theta\omega_{l}^{\star}+\zeta_{\mathsf{2nd},l}\sigma_{1}^{\star}\right)\left\Vert \bm{U}_{l,\cdot}^{\star}\bm{\Sigma}^{\star}\right\Vert _{2}+\zeta_{\mathsf{2nd},l}^{2}\sigma_{1}^{\star2}.
\]
According to Lemma \ref{lemma:pca-noise-level-est}, we can obtain
\[
\left|S_{l,l}-S_{l,l}^{\star}\right|\lesssim\mathsf{UB}_{l}\qquad\forall\,l\in\left[d\right].
\]
Lemma \ref{lemma:pca-noise-level-est} also tells us that 
\begin{align*}
\left|\omega_{i}^{2}-\omega_{i}^{\star2}\right| & \lesssim\sqrt{\frac{\log^{2}\left(n+d\right)}{np}}\omega_{i}^{\star2}+\mathsf{UB}_{i}.
\end{align*}
With these basic bounds in mind, we proceed with the following decomposition:
\begin{align*}
\frac{1}{np}\left|\omega_{i}^{2}S_{j,j}-\omega_{i}^{\star2}S_{j,j}^{\star}\right| & \lesssim\frac{1}{np}S_{j,j}\left|\omega_{i}^{2}-\omega_{i}^{\star2}\right|+\frac{1}{np}\omega_{i}^{\star2}\left|S_{j,j}-S_{j,j}^{\star}\right|\\
 & \lesssim\underbrace{\frac{1}{np}S_{j,j}^{\star}\left|\omega_{i}^{2}-\omega_{i}^{\star2}\right|}_{\eqqcolon\gamma_{3,1}}+\underbrace{\frac{1}{np}\omega_{i}^{\star2}\left|S_{j,j}-S_{j,j}^{\star}\right|}_{\eqqcolon\gamma_{3,2}}+\underbrace{\frac{1}{np}\left|S_{j,j}-S_{j,j}^{\star}\right|\left|\omega_{i}^{2}-\omega_{i}^{\star2}\right|}_{\eqqcolon\gamma_{3,3}},
\end{align*}
leaving us with three terms to cope with. 
\begin{itemize}
\item Regarding $\gamma_{3,1}$, we can upper bound 
\begin{align*}
\gamma_{3,1} & \lesssim\underbrace{\frac{1}{np}\left\Vert \bm{U}_{j,\cdot}^{\star}\bm{\Sigma}^{\star}\right\Vert _{2}^{2}\sqrt{\frac{\log^{2}\left(n+d\right)}{np}}\omega_{i}^{\star2}}_{\eqqcolon\gamma_{3,1,1}}+\underbrace{\frac{1}{np}S_{j,j}^{\star}\mathsf{UB}_{i}}_{\eqqcolon\gamma_{3,1,2}}\lesssim\delta v_{i,j}^{\star}.
\end{align*}
Here, the last relation holds since (i) the first term 
\[
\gamma_{3,1,1}\lesssim\sqrt{\frac{\kappa\log^{2}\left(n+d\right)}{n^{3}p^{3}}}\omega_{i}^{\star2}\sigma_{r}^{\star2}\left\Vert \bm{U}_{j,\cdot}^{\star}\right\Vert _{2}^{2}\lesssim\sqrt{\frac{\kappa\log^{2}\left(n+d\right)}{np}}v_{i,j}^{\star}\lesssim\delta v_{i,j}^{\star}
\]
provided that $np\gtrsim\delta^{-2}\kappa\log^{2}(n+d)$; and (ii)
we can also demonstrate that the second term obeys $\gamma_{3,1,2}\lesssim\delta v_{i,j}^{\star}$,
since it is easy seen that $\gamma_{3,1,2}$ admits the same upper
bound as for $\gamma_{1,1}$. 
\item Regarding $\gamma_{3,2}$, it can be seen that 
\begin{align*}
\gamma_{3,2} & \lesssim\underbrace{\frac{1}{np}\omega_{i}^{\star2}\left(\theta+\sqrt{\frac{\kappa^{3}r\log\left(n+d\right)}{n}}\right)\left\Vert \bm{U}_{j,\cdot}^{\star}\bm{\Sigma}^{\star}\right\Vert _{2}^{2}}_{\eqqcolon\gamma_{3,2,1}}+\underbrace{\frac{1}{np}\theta^{2}\omega_{i}^{\star2}\omega_{j}^{\star2}}_{\eqqcolon\gamma_{3,2,2}}\\
 & \quad+\underbrace{\frac{1}{np}\theta\omega_{i}^{\star2}\omega_{j}^{\star}\left\Vert \bm{U}_{j,\cdot}^{\star}\bm{\Sigma}^{\star}\right\Vert _{2}}_{\eqqcolon\gamma_{3,2,3}}+\underbrace{\frac{1}{np}\omega_{i}^{\star2}\zeta_{\mathsf{2nd},j}\sigma_{1}^{\star}\left\Vert \bm{U}_{j,\cdot}^{\star}\bm{\Sigma}^{\star}\right\Vert _{2}}_{\eqqcolon\gamma_{3,2,4}}+\underbrace{\frac{1}{np}\omega_{i}^{\star2}\zeta_{\mathsf{2nd},j}^{2}\sigma_{1}^{\star2}}_{\eqqcolon\gamma_{3,2,5}}\\
 & \lesssim\delta v_{i,j}^{\star},
\end{align*}
where the last line holds true since 
\begin{align*}
\gamma_{3,2,1} & \lesssim\frac{\delta}{np}\omega_{i}^{\star2}\left\Vert \bm{U}_{j,\cdot}^{\star}\bm{\Sigma}^{\star}\right\Vert _{2}^{2}\lesssim\delta v_{i,j}^{\star},\\
\gamma_{3,2,2}+\gamma_{3,2,3} & \lesssim\delta\frac{\sigma_{r}^{\star2}}{np}\left(\omega_{j}^{\star2}\left\Vert \bm{U}_{i,\cdot}^{\star}\right\Vert _{2}^{2}+\omega_{i}^{\star2}\left\Vert \bm{U}_{j,\cdot}^{\star}\right\Vert _{2}^{2}\right)\lesssim\delta v_{i,j}^{\star},\\
\gamma_{3,2,4} & \lesssim\frac{\sigma_{r}^{\star}}{\sqrt{np}}\omega_{i}^{\star}\left\Vert \bm{U}_{j,\cdot}^{\star}\right\Vert _{2}\cdot\sqrt{\frac{\kappa}{np}}\omega_{i}^{\star}\zeta_{\mathsf{2nd},j}\sigma_{1}^{\star}\lesssim v_{i,j}^{\star1/2}\cdot\sqrt{\frac{\kappa}{np}}\omega_{i}^{\star}\zeta_{\mathsf{2nd},j}\sigma_{1}^{\star}\\
 & \overset{\text{(i)}}{\lesssim}v_{i,j}^{\star1/2}\cdot\left(\sqrt{\frac{1}{np}}\varepsilon\frac{\omega_{i}^{\star}\left\Vert \bm{U}_{j,\cdot}^{\star}\bm{\Sigma}^{\star}\right\Vert _{2}+\omega_{i}^{\star}\omega_{j}^{\star}}{\sqrt{\min\left\{ ndp^{2}\kappa,np\right\} }}+\sqrt{\frac{1}{np}}\varepsilon\frac{\omega_{\min}\omega_{i}^{\star}\omega_{j}^{\star}}{\sqrt{np^{2}}\sigma_{1}^{\star}}\right)\\
 & \lesssim v_{i,j}^{\star1/2}\cdot\delta\left(\frac{\sigma_{r}^{\star}}{\sqrt{np}}+\frac{\omega_{\min}}{\sqrt{np^{2}}}\right)\left(\omega_{j}^{\star}\left\Vert \bm{U}_{i,\cdot}^{\star}\right\Vert _{2}+\omega_{i}^{\star}\left\Vert \bm{U}_{j,\cdot}^{\star}\right\Vert _{2}\right)\lesssim\delta v_{i,j}^{\star},\\
\gamma_{3,2,5} & \overset{\text{(ii)}}{\lesssim}\left(\sqrt{\frac{1}{np}}\varepsilon\frac{\omega_{i}^{\star}\left\Vert \bm{U}_{j,\cdot}^{\star}\bm{\Sigma}^{\star}\right\Vert _{2}+\omega_{i}^{\star}\omega_{j}^{\star}}{\sqrt{\min\left\{ ndp^{2}\kappa,np\right\} }}+\sqrt{\frac{1}{np}}\varepsilon\frac{\omega_{\min}\omega_{i}^{\star}\omega_{j}^{\star}}{\sqrt{np^{2}}\sigma_{1}^{\star}}\right)^{2}\\
 & \lesssim\delta\left(\frac{\sigma_{r}^{\star2}}{np}+\frac{\omega_{\min}^{2}}{np^{2}}\right)\left(\omega_{j}^{\star2}\left\Vert \bm{U}_{i,\cdot}^{\star}\right\Vert _{2}^{2}+\omega_{i}^{\star2}\left\Vert \bm{U}_{j,\cdot}^{\star}\right\Vert _{2}^{2}\right)\lesssim\delta v_{i,j}^{\star},
\end{align*}
provided that $\theta\lesssim\delta$, $n\gtrsim\delta^{-2}\kappa^{3}r\log(n+d)$,
$\varepsilon\lesssim1$, $np\gtrsim\delta^{-2}\kappa$, $ndp^{2}\gtrsim\delta^{-2}$
and 
\[
\left\Vert \bm{U}_{i,\cdot}^{\star}\right\Vert _{2}+\left\Vert \bm{U}_{j,\cdot}^{\star}\right\Vert _{2}\gtrsim\delta^{-1}\kappa^{1/2}\theta\frac{\omega_{\max}}{\sigma_{r}^{\star}}+\delta^{-1}\frac{\omega_{\max}}{\sigma_{r}^{\star}}\sqrt{\frac{1}{ndp^{2}\kappa}}+\delta^{-1}\frac{\omega_{\max}}{\sigma_{r}^{\star}}\sqrt{\frac{1}{np}}.
\]
Here, (i) and (ii) arise from (\ref{eq:ce-var-est-inter-2}). 
\item When it comes to $\gamma_{3,3}$, we obtain 
\[
\gamma_{3,3}\lesssim\underbrace{\frac{1}{np}\left|S_{j,j}-S_{j,j}^{\star}\right|\sqrt{\frac{\log^{2}\left(n+d\right)}{np}}\omega_{i}^{\star2}}_{\eqqcolon\gamma_{3,3,1}}+\underbrace{\frac{1}{np}\left|S_{j,j}-S_{j,j}^{\star}\right|\mathsf{UB}_{i}}_{\eqqcolon\gamma_{3,3,2}}\lesssim\delta v_{i,j}^{\star}.
\]
The last relation holds since (i) the first term 
\begin{align*}
\gamma_{3,3,1} & \lesssim\sqrt{\frac{\log^{2}\left(n+d\right)}{np}}\gamma_{3,2}\lesssim\delta v_{i,j}^{\star}
\end{align*}
provided that $np\gtrsim\log^{2}(n+d)$; and (ii) one can easily check
that the second term $\gamma_{3,3,2}$ admits the same upper bound
as for $\gamma_{1,2}$. 
\end{itemize}
With the above results in hand, one can conclude that 
\[
\frac{1}{np}\left|\omega_{i}^{2}S_{j,j}-\omega_{i}^{\star2}S_{j,j}^{\star}\right|\lesssim\gamma_{3,1}+\gamma_{3,2}+\gamma_{3,3}\lesssim\delta v_{i,j}^{\star},
\]
and similarly, 
\[
\frac{1}{np}\left|\omega_{j}^{2}S_{i,i}-\omega_{j}^{\star2}S_{i,i}^{\star}\right|\lesssim\delta v_{i,j}^{\star}.
\]
These allow us to reach 
\[
\left|\alpha_{3}-\beta_{3}\right|\lesssim\frac{1}{np}\left|\omega_{i}^{2}S_{j,j}-\omega_{i}^{\star2}S_{j,j}^{\star}\right|+\frac{1}{np}\left|\omega_{j}^{2}S_{i,i}-\omega_{j}^{\star2}S_{i,i}^{\star}\right|\lesssim\delta v_{i,j}^{\star},
\]
provided that $np\gtrsim\delta^{-2}\kappa\log^{2}(n+d)$, $\theta\lesssim\delta$,
$n\gtrsim\delta^{-2}\kappa^{3}r\log(n+d)$, $\varepsilon\lesssim1$,
and 
\begin{equation}
\left\Vert \bm{U}_{i,\cdot}^{\star}\right\Vert _{2}+\left\Vert \bm{U}_{j,\cdot}^{\star}\right\Vert _{2}\gtrsim\delta^{-1}\kappa^{1/2}\theta\frac{\omega_{\max}}{\sigma_{r}^{\star}}.\label{eq:ce-var-est-inter-11}
\end{equation}
In view of (\ref{eq:theta-exact}), we know that 
\begin{align}
\theta\frac{\omega_{\max}}{\sigma_{r}^{\star}} & \asymp\left(\frac{\omega_{\max}}{\sigma_{r}^{\star}}\sqrt{\frac{d\kappa\log^{2}\left(n+d\right)}{np}}+\frac{\omega_{\max}}{\sigma_{r}^{\star}}\sqrt{\frac{\kappa^{2}\mu r\log^{3}\left(n+d\right)}{np^{2}}}+\frac{\omega_{\max}^{2}}{p\sigma_{r}^{\star2}}\sqrt{\frac{d\kappa\log^{2}\left(n+d\right)}{n}}\right)\cdot\sqrt{\frac{r}{d}}\nonumber \\
 & \asymp\sqrt{\kappa\log^{2}\left(n+d\right)}\left(\frac{\omega_{\max}}{\sigma_{r}^{\star}}\sqrt{\frac{d}{np}}+\frac{\omega_{\max}^{2}}{p\sigma_{r}^{\star2}}\sqrt{\frac{d}{n}}+\frac{\omega_{\max}}{\sigma_{r}^{\star}}\sqrt{\frac{\kappa\mu r\log\left(n+d\right)}{np^{2}}}\right)\cdot\sqrt{\frac{r}{d}}\nonumber \\
 & \asymp\sqrt{\kappa\log^{2}\left(n+d\right)}\left[\frac{\omega_{\max}}{\sigma_{r}^{\star}}\sqrt{\frac{d}{np}}+\frac{\omega_{\max}^{2}}{p\sigma_{r}^{\star2}}\sqrt{\frac{d}{n}}+\frac{\kappa\mu r\log\left(n+d\right)}{\sqrt{ndp^{2}}}\right]\cdot\sqrt{\frac{r}{d}},\label{eq:theta-omega-max}
\end{align}
where we have used the AM-GM inequality in the last line. As a result,
(\ref{eq:ce-var-est-inter-11}) is guaranteed by 
\[
\left\Vert \bm{U}_{i,\cdot}^{\star}\right\Vert _{2}+\left\Vert \bm{U}_{j,\cdot}^{\star}\right\Vert _{2}\gtrsim\delta^{-1}\kappa\log\left(n+d\right)\left[\frac{\omega_{\max}}{\sigma_{r}^{\star}}\sqrt{\frac{d}{np}}+\frac{\omega_{\max}^{2}}{p\sigma_{r}^{\star2}}\sqrt{\frac{d}{n}}+\frac{\kappa\mu r\log\left(n+d\right)}{\sqrt{ndp^{2}}}\right]\cdot\sqrt{\frac{r}{d}}.
\]

\paragraph{Step 4: bounding $\vert\alpha_{4}-\beta_{4}\vert$ and $\vert\alpha_{6}-\beta_{6}\vert$.}

For each $l\in[d]$, let us denote 
\[
\Delta_{l}\coloneqq\big|\omega_{l}^{\star2}+(1-p)S_{l,l}^{\star2}-\omega_{l}^{2}-(1-p)S_{l,l}^{2}\big|.
\]
Lemma \ref{lemma:pca-noise-level-est} tells us that, for each $l\in[d]$,
\begin{equation}
\Delta_{l}\lesssim\sqrt{\frac{\log^{2}\left(n+d\right)}{np}}\omega_{\max}^{2}+\zeta_{\mathsf{1st}}\frac{\sqrt{\kappa^{2}\mu r^{2}\log\left(n+d\right)}}{d}+\sqrt{\frac{\kappa^{2}\mu^{2}r^{3}\log\left(n+d\right)}{nd^{2}}}\sigma_{1}^{\star2}.\label{eq:ce-var-est-inter-5}
\end{equation}
We also know that 
\begin{align}
\Delta_{i} & \lesssim\sqrt{\frac{\log^{2}\left(n+d\right)}{np}}\omega_{i}^{\star2}+\mathsf{UB}_{i}\qquad\text{and}\qquad\Delta_{j}\lesssim\sqrt{\frac{\log^{2}\left(n+d\right)}{np}}\omega_{j}^{\star2}+\mathsf{UB}_{j}.\label{eq:ce-var-est-inter-6}
\end{align}
In addition, for each $l\in[d]$, it holds that 
\begin{equation}
\omega_{l}^{\star2}+\left(1-p\right)S_{l,l}^{\star}\leq\omega_{\max}^{2}+\left\Vert \bm{U}_{l,\cdot}^{\star}\bm{\Sigma}^{\star}\right\Vert _{2}^{2}\lesssim\omega_{\max}^{2}+\frac{\mu r}{d}\sigma_{1}^{\star2}.\label{eq:ce-var-est-inter-7}
\end{equation}
We also have the following bound 
\begin{align}
\left\Vert \bm{U}_{j,\cdot}\right\Vert _{2} & \leq\left\Vert \bm{U}_{j,\cdot}^{\star}\right\Vert _{2}+\left\Vert \left(\bm{U}\bm{R}-\bm{U}^{\star}\right)_{j,\cdot}\right\Vert _{2}\nonumber \\
 & \overset{\text{(i)}}{\lesssim}\left\Vert \bm{U}_{j,\cdot}^{\star}\right\Vert _{2}+\frac{\theta}{\sqrt{\kappa}\sigma_{r}^{\star}}\left(\left\Vert \bm{U}_{j,\cdot}^{\star}\bm{\Sigma}^{\star}\right\Vert _{2}+\omega_{j}^{\star}\right)+\zeta_{\mathsf{2nd},j}\nonumber \\
 & \overset{\text{(ii)}}{\lesssim}\left\Vert \bm{U}_{j,\cdot}^{\star}\right\Vert _{2}+\frac{\theta}{\sqrt{\kappa}\sigma_{r}^{\star}}\omega_{j}^{\star}+\zeta_{\mathsf{2nd},j}\nonumber \\
 & \overset{\text{(iii)}}{\lesssim}\left\Vert \bm{U}_{j,\cdot}^{\star}\right\Vert _{2}+\frac{\theta}{\sqrt{\kappa}\sigma_{r}^{\star}}\omega_{j}^{\star}.\label{eq:ce-var-est-inter-8}
\end{align}
Here (i) follows from Lemma \ref{lemma:pca-1st-err}, (ii) holds provided
that $\theta\lesssim1$, (iii) utilizes (\ref{eq:ce-var-est-inter-2})
and (\ref{eq:theta-exact}) and holds provided that $\varepsilon\lesssim1$.
We can thus decompose 
\begin{align*}
\left|\alpha_{4}-\beta_{4}\right| & \lesssim\underbrace{\frac{1}{np^{2}}\left[\omega_{i}^{\star2}+\left(1-p\right)S_{i,i}^{\star}\right]\sum_{k=1}^{d}\Delta_{k}\left(\bm{U}_{k,\cdot}\bm{U}_{j,\cdot}^{\top}\right)^{2}}_{\eqqcolon\gamma_{4,1}}+\underbrace{\frac{\Delta_{i}}{np^{2}}\sum_{k=1}^{d}\left[\omega_{k}^{2}+\left(1-p\right)S_{k,k}\right]\left(\bm{U}_{k,\cdot}\bm{U}_{j,\cdot}^{\top}\right)^{2}}_{\eqqcolon\gamma_{4,2}}\\
 & \quad+\underbrace{\frac{1}{np^{2}}\left[\omega_{i}^{\star2}+\left(1-p\right)S_{i,i}^{\star}\right]\sum_{k=1}^{d}\left[\omega_{k}^{\star2}+\left(1-p\right)S_{k,k}^{\star}\right]\left|\left(\bm{U}_{k,\cdot}^{\star}\bm{U}_{j,\cdot}^{\star\top}\right)^{2}-\left(\bm{U}_{k,\cdot}\bm{U}_{j,\cdot}^{\top}\right)^{2}\right|}_{\eqqcolon\gamma_{4,3}}.
\end{align*}

\subparagraph{Step 4.1: bounding $\gamma_{4,1}$.}

We have learned from (\ref{eq:ce-var-est-inter-8}) that 
\begin{align*}
\gamma_{4,1} & \lesssim\frac{1}{np^{2}}\left[\omega_{i}^{\star2}+\left(1-p\right)S_{i,i}^{\star}\right]\left\Vert \bm{U}_{j,\cdot}\right\Vert _{2}^{2}\max_{1\leq k\leq d}\Delta_{k}\\
 & \lesssim\underbrace{\frac{1}{np^{2}}\left[\omega_{i}^{\star2}+\left(1-p\right)S_{i,i}^{\star}\right]\left\Vert \bm{U}_{j,\cdot}^{\star}\right\Vert _{2}^{2}\max_{1\leq k\leq d}\Delta_{k}}_{\eqqcolon\gamma_{4,1,1}}+\underbrace{\frac{\theta^{2}\omega_{j}^{\star2}}{np^{2}\sigma_{1}^{\star2}}\left[\omega_{i}^{\star2}+\left(1-p\right)S_{i,i}^{\star}\right]\max_{1\leq k\leq d}\Delta_{k}}_{\eqqcolon\gamma_{4,1,2}}.
\end{align*}
Regarding $\gamma_{4,1,1}$, we can derive 
\begin{align*}
\gamma_{4,1,1} & \lesssim\underbrace{\frac{1}{np^{2}}\omega_{\max}^{2}\omega_{i}^{\star2}\left\Vert \bm{U}_{j,\cdot}^{\star}\right\Vert _{2}^{2}\sqrt{\frac{\log^{2}\left(n+d\right)}{np}}}_{\eqqcolon\gamma_{4,1,1,1}}+\underbrace{\frac{\sqrt{\kappa^{2}\mu r^{2}\log\left(n+d\right)}}{ndp^{2}}\omega_{i}^{\star2}\left\Vert \bm{U}_{j,\cdot}^{\star}\right\Vert _{2}^{2}\zeta_{\mathsf{1st}}}_{\eqqcolon\gamma_{4,1,1,2}}\\
 & \quad+\underbrace{\frac{1}{np^{2}}\omega_{i}^{\star2}\left\Vert \bm{U}_{j,\cdot}^{\star}\right\Vert _{2}^{2}\sqrt{\frac{\kappa^{2}\mu^{2}r^{3}\log\left(n+d\right)}{nd^{2}}}\sigma_{1}^{\star2}}_{\eqqcolon\gamma_{4,1,1,3}}+\underbrace{\frac{1}{np^{2}}\left\Vert \bm{U}_{i,\cdot}^{\star}\bm{\Sigma}^{\star}\right\Vert _{2}^{2}\left\Vert \bm{U}_{j,\cdot}^{\star}\right\Vert _{2}^{2}\sqrt{\frac{\log^{2}\left(n+d\right)}{np}}\omega_{\max}^{2}}_{\eqqcolon\gamma_{4,1,1,4}}\\
 & \quad+\underbrace{\frac{\sqrt{\kappa^{2}\mu r^{2}\log\left(n+d\right)}}{ndp^{2}}\left\Vert \bm{U}_{i,\cdot}^{\star}\bm{\Sigma}^{\star}\right\Vert _{2}^{2}\left\Vert \bm{U}_{j,\cdot}^{\star}\right\Vert _{2}^{2}\zeta_{\mathsf{1st}}}_{\eqqcolon\gamma_{4,1,1,5}}+\underbrace{\frac{1}{np^{2}}\left\Vert \bm{U}_{i,\cdot}^{\star}\bm{\Sigma}^{\star}\right\Vert _{2}^{2}\left\Vert \bm{U}_{j,\cdot}^{\star}\right\Vert _{2}^{2}\sqrt{\frac{\kappa^{2}\mu^{2}r^{3}\log\left(n+d\right)}{nd^{2}}}\sigma_{1}^{\star2}}_{\eqqcolon\gamma_{4,1,1,6}}\\
 & \lesssim\delta v_{i,j}^{\star}.
\end{align*}
Here, the last line holds since 
\begin{align*}
\gamma_{4,1,1,1} & \lesssim\delta\frac{\omega_{\min}^{2}}{np^{2}}\left(\omega_{j}^{\star2}\left\Vert \bm{U}_{i,\cdot}^{\star}\right\Vert _{2}^{2}+\omega_{i}^{\star2}\left\Vert \bm{U}_{j,\cdot}^{\star}\right\Vert _{2}^{2}\right)\lesssim\delta v_{i,j}^{\star},\\
\gamma_{4,1,1,2}+\gamma_{4,1,1,3}+\gamma_{4,1,1,4} & \lesssim\delta\frac{\sigma_{r}^{\star2}}{ndp^{2}}\left(\omega_{j}^{\star2}\left\Vert \bm{U}_{i,\cdot}^{\star}\right\Vert _{2}^{2}+\omega_{i}^{\star2}\left\Vert \bm{U}_{j,\cdot}^{\star}\right\Vert _{2}^{2}\right)\lesssim\delta v_{i,j}^{\star},\\
\gamma_{4,1,1,5}+\gamma_{4,1,1,6} & \lesssim\delta\frac{1}{ndp^{2}\kappa}\left\Vert \bm{U}_{i,\cdot}^{\star}\bm{\Sigma}^{\star}\right\Vert _{2}^{2}\left\Vert \bm{U}_{j,\cdot}^{\star}\bm{\Sigma}^{\star}\right\Vert _{2}^{2}\lesssim\delta v_{i,j}^{\star},
\end{align*}
provided that $np\gtrsim\delta^{-2}\kappa^{2}\mu^{2}r^{2}\log^{2}(n+d)$,
$n\gtrsim\delta^{-2}\kappa^{4}\mu^{2}r^{3}\log(n+d)$, $np\gtrsim\delta^{-2}\kappa_{\omega}^{2}\log^{2}(n+d)$and
\[
\frac{\zeta_{\mathsf{1st}}}{\sigma_{r}^{\star2}}\lesssim\frac{\delta}{\sqrt{\kappa^{2}\mu r^{2}\log\left(n+d\right)}}.
\]
Regarding $\gamma_{4,1,2}$, we can upper bound 
\begin{align*}
\gamma_{4,1,2} & \lesssim\underbrace{\frac{\theta^{2}}{np^{2}}\frac{\omega_{i}^{\star2}\omega_{j}^{\star2}\omega_{\max}^{2}}{\sigma_{1}^{\star2}}\sqrt{\frac{\log^{2}\left(n+d\right)}{np}}}_{\eqqcolon\gamma_{4,1,2,1}}+\underbrace{\frac{\sqrt{\kappa^{2}\mu r^{2}\log\left(n+d\right)}}{ndp^{2}}\frac{\omega_{i}^{\star2}\omega_{j}^{\star2}}{\sigma_{1}^{\star2}}\theta^{2}\zeta_{\mathsf{1st}}}_{\eqqcolon\gamma_{4,1,2,2}}\\
 & \quad+\underbrace{\frac{1}{np^{2}}\omega_{i}^{\star2}\omega_{j}^{\star2}\theta^{2}\sqrt{\frac{\kappa^{2}\mu^{2}r^{3}\log\left(n+d\right)}{nd^{2}}}}_{\eqqcolon\gamma_{4,1,2,3}}+\underbrace{\frac{\theta^{2}}{np^{2}}\left\Vert \bm{U}_{i,\cdot}^{\star}\bm{\Sigma}^{\star}\right\Vert _{2}^{2}\frac{\omega_{\max}^{2}\omega_{j}^{\star2}}{\sigma_{1}^{\star2}}\sqrt{\frac{\log^{2}\left(n+d\right)}{np}}}_{\eqqcolon\gamma_{4,1,2,4}}\\
 & \quad+\underbrace{\frac{\sqrt{\kappa^{2}\mu r^{2}\log\left(n+d\right)}}{ndp^{2}}\left\Vert \bm{U}_{i,\cdot}^{\star}\bm{\Sigma}^{\star}\right\Vert _{2}^{2}\frac{\omega_{j}^{\star2}}{\sigma_{1}^{\star2}}\theta^{2}\zeta_{\mathsf{1st}}}_{\eqqcolon\gamma_{4,1,2,5}}+\underbrace{\frac{1}{np^{2}}\left\Vert \bm{U}_{i,\cdot}^{\star}\bm{\Sigma}^{\star}\right\Vert _{2}^{2}\theta^{2}\omega_{j}^{\star2}\sqrt{\frac{\kappa^{2}\mu^{2}r^{3}\log\left(n+d\right)}{nd^{2}}}}_{\eqqcolon\gamma_{4,1,2,6}}.
\end{align*}
There are six terms on the right-hand side of the above inequality,
which we shall control separately. 
\begin{itemize}
\item With regards to $\gamma_{4,1,2,1}$, we observe that 
\begin{align*}
\gamma_{4,1,2,1} & \asymp\underbrace{\frac{\omega_{i}^{\star2}\omega_{j}^{\star2}\omega_{\max}^{2}}{\sigma_{r}^{\star2}}\frac{r\log^{3}\left(n+d\right)}{n^{2.5}p^{3.5}}}_{\eqqcolon\gamma_{4,1,2,1,1}}+\underbrace{\frac{\omega_{i}^{\star2}\omega_{j}^{\star2}\omega_{\max}^{2}}{\sigma_{r}^{\star2}}\frac{\kappa\mu r^{2}\log^{4}\left(n+d\right)}{n^{2.5}dp^{4.5}}}_{\eqqcolon\gamma_{4,1,2,1,2}}+\underbrace{\frac{\omega_{i}^{\star2}\omega_{j}^{\star2}\omega_{\max}^{4}}{\sigma_{r}^{\star4}}\frac{r\log^{3}\left(n+d\right)}{n^{2.5}p^{4.5}}}_{\eqqcolon\gamma_{4,1,2,1,3}}\lesssim\delta v_{i,j}^{\star}.
\end{align*}
Here, the last relation holds since 
\begin{align*}
\gamma_{4,1,2,1,1}+\gamma_{4,1,2,1,3} & \lesssim\delta\frac{\omega_{\min}^{2}}{np^{2}}\left(\omega_{j}^{\star2}\left\Vert \bm{U}_{i,\cdot}^{\star}\right\Vert _{2}^{2}+\omega_{i}^{\star2}\left\Vert \bm{U}_{j,\cdot}^{\star}\right\Vert _{2}^{2}\right)\lesssim\delta v_{i,j}^{\star},\\
\gamma_{4,1,2,1,2} & \lesssim\delta\frac{\sigma_{r}^{\star2}}{ndp^{2}}\left(\omega_{j}^{\star2}\left\Vert \bm{U}_{i,\cdot}^{\star}\right\Vert _{2}^{2}+\omega_{i}^{\star2}\left\Vert \bm{U}_{j,\cdot}^{\star}\right\Vert _{2}^{2}\right)\lesssim\delta v_{i,j}^{\star},
\end{align*}
provided that 
\[
\left\Vert \bm{U}_{i,\cdot}^{\star}\right\Vert _{2}+\left\Vert \bm{U}_{j,\cdot}^{\star}\right\Vert _{2}\gtrsim\delta^{-1/2}\frac{\sqrt{\kappa_{\omega}\log^{3}\left(n+d\right)}}{n^{1/4}p^{1/4}}\left[\frac{\omega_{\max}}{\sigma_{r}^{\star}}\sqrt{\frac{d}{np}}+\frac{\omega_{\max}^{2}}{p\sigma_{r}^{\star2}}\sqrt{\frac{d}{n}}\right]\sqrt{\frac{r}{d}}.
\]
\item When it comes to $\gamma_{4,1,2,2}$, we obtain 
\begin{align*}
\gamma_{4,1,2,2} & \asymp\underbrace{\frac{\omega_{i}^{\star2}\omega_{j}^{\star2}}{\sigma_{r}^{\star2}}\frac{\kappa\mu^{1/2}r^{2}\log^{5/2}\left(n+d\right)}{n^{2}dp^{3}}\zeta_{\mathsf{1st}}}_{\eqqcolon\gamma_{4,1,2,2,1}}+\underbrace{\frac{\omega_{i}^{\star2}\omega_{j}^{\star2}}{\sigma_{r}^{\star2}}\frac{\kappa^{2}\mu^{3/2}r^{3}\log^{7/2}\left(n+d\right)}{n^{2}d^{2}p^{4}}\zeta_{\mathsf{1st}}}_{\eqqcolon\gamma_{4,1,2,2,2}}\\
 & \quad+\underbrace{\frac{\omega_{i}^{\star2}\omega_{j}^{\star2}\omega_{\max}^{2}}{\sigma_{r}^{\star4}}\frac{\kappa\mu^{1/2}r^{2}\log^{5/2}\left(n+d\right)}{n^{2}dp^{4}}\zeta_{\mathsf{1st}}}_{\eqqcolon\gamma_{4,1,2,2,3}}\lesssim\delta v_{i,j}^{\star}.
\end{align*}
Here, the last relation holds since 
\begin{align*}
\gamma_{4,1,2,2,1}+\gamma_{4,1,2,2,3} & \lesssim\delta\frac{\sigma_{r}^{\star2}}{ndp^{2}}\left(\omega_{j}^{\star2}\left\Vert \bm{U}_{i,\cdot}^{\star}\right\Vert _{2}^{2}+\omega_{i}^{\star2}\left\Vert \bm{U}_{j,\cdot}^{\star}\right\Vert _{2}^{2}\right)\lesssim\delta v_{i,j}^{\star},\\
\gamma_{4,1,2,2,2} & \lesssim\delta\frac{\omega_{\min}^{2}}{np^{2}}\left(\omega_{j}^{\star2}\left\Vert \bm{U}_{i,\cdot}^{\star}\right\Vert _{2}^{2}+\omega_{i}^{\star2}\left\Vert \bm{U}_{j,\cdot}^{\star}\right\Vert _{2}^{2}\right)\lesssim\delta v_{i,j}^{\star},
\end{align*}
provided that 
\begin{align*}
\left\Vert \bm{U}_{i,\cdot}^{\star}\right\Vert _{2}+\left\Vert \bm{U}_{j,\cdot}^{\star}\right\Vert _{2} & \gtrsim\delta^{-1/2}\sqrt{\frac{\zeta_{\mathsf{1st}}}{\sigma_{r}^{\star2}}}\Biggl[\frac{\omega_{\max}}{\sigma_{r}^{\star}}\sqrt{\frac{d\kappa\mu^{1/2}r\log^{5/2}\left(n+d\right)}{np}}+\sqrt{\frac{\kappa^{2}\mu^{3/2}r^{2}\kappa_{\omega}\log^{7/2}\left(n+d\right)}{ndp^{2}}}\\
 & \qquad\qquad\qquad+\frac{\omega_{\max}^{2}}{p\sigma_{r}^{\star2}}\sqrt{\frac{d\kappa\mu^{1/2}r\log^{5/2}\left(n+d\right)}{n}}\Biggr]\sqrt{\frac{r}{d}}.
\end{align*}
\item Regarding $\gamma_{4,1,2,3}$, we can obtain 
\begin{align*}
\gamma_{4,1,2,3} & \asymp\underbrace{\frac{\omega_{i}^{\star2}\omega_{j}^{\star2}}{n^{2.5}dp^{3}}\kappa^{2}\mu r^{5/2}\log^{5/2}\left(n+d\right)}_{\eqqcolon\gamma_{4,1,2,3,1}}+\underbrace{\frac{\omega_{i}^{\star2}\omega_{j}^{\star2}}{n^{2.5}d^{2}p^{4}}\kappa^{3}\mu^{2}r^{7/2}\log^{7/2}\left(n+d\right)}_{\eqqcolon\gamma_{4,1,2,3,2}}\\
 & \quad+\underbrace{\frac{\omega_{\max}^{2}\omega_{i}^{\star2}\omega_{j}^{\star2}}{\sigma_{r}^{\star2}}\frac{\kappa^{2}\mu r^{5/2}\log^{5/2}\left(n+d\right)}{n^{5/2}dp^{4}}}_{\eqqcolon\gamma_{4,1,2,3,3}}\lesssim\delta v_{i,j}^{\star}.
\end{align*}
Here, the last relation holds since 
\begin{align*}
\gamma_{4,1,2,3,1}+\gamma_{4,1,2,3,3} & \lesssim\delta\frac{\sigma_{r}^{\star2}}{ndp^{2}}\left(\omega_{j}^{\star2}\left\Vert \bm{U}_{i,\cdot}^{\star}\right\Vert _{2}^{2}+\omega_{i}^{\star2}\left\Vert \bm{U}_{j,\cdot}^{\star}\right\Vert _{2}^{2}\right)\lesssim\delta v_{i,j}^{\star},\\
\gamma_{4,1,2,3,2} & \lesssim\delta\frac{\omega_{\min}^{2}}{np^{2}}\left(\omega_{j}^{\star2}\left\Vert \bm{U}_{i,\cdot}^{\star}\right\Vert _{2}^{2}+\omega_{i}^{\star2}\left\Vert \bm{U}_{j,\cdot}^{\star}\right\Vert _{2}^{2}\right)\lesssim\delta v_{i,j}^{\star},
\end{align*}
provided that 
\begin{align*}
\left\Vert \bm{U}_{i,\cdot}^{\star}\right\Vert _{2}+\left\Vert \bm{U}_{j,\cdot}^{\star}\right\Vert _{2} & \gtrsim\delta^{-1/2}\frac{\kappa^{1/2}\mu^{1/4}r^{1/4}}{n^{1/4}}\Biggl[\frac{\omega_{\max}}{\sigma_{r}^{\star}}\sqrt{\frac{d\kappa\mu^{1/2}r\log^{5/2}\left(n+d\right)}{np}}+\sqrt{\frac{\kappa^{2}\mu^{3/2}r^{2}\kappa_{\omega}\log^{7/2}\left(n+d\right)}{ndp^{2}}}\\
 & \qquad\qquad\qquad\qquad+\frac{\omega_{\max}^{2}}{p\sigma_{r}^{\star2}}\sqrt{\frac{d\kappa\mu^{1/2}r\log^{5/2}\left(n+d\right)}{n}}\Biggr]\sqrt{\frac{r}{d}}.
\end{align*}
\item For $\gamma_{4,1,2,4}$, $\gamma_{4,1,2,5}$ and $\gamma_{4,1,2,6}$,
it is straightforward to check that 
\begin{align*}
\gamma_{4,1,2,4} & \lesssim\delta\frac{\omega_{\min}^{2}}{np^{2}}\left(\omega_{j}^{\star2}\left\Vert \bm{U}_{i,\cdot}^{\star}\right\Vert _{2}^{2}+\omega_{i}^{\star2}\left\Vert \bm{U}_{j,\cdot}^{\star}\right\Vert _{2}^{2}\right)\lesssim\delta v_{i,j}^{\star},\\
\gamma_{4,1,2,5}+\gamma_{4,1,2,6} & \lesssim\delta\frac{\sigma_{r}^{\star2}}{ndp^{2}}\left(\omega_{j}^{\star2}\left\Vert \bm{U}_{i,\cdot}^{\star}\right\Vert _{2}^{2}+\omega_{i}^{\star2}\left\Vert \bm{U}_{j,\cdot}^{\star}\right\Vert _{2}^{2}\right)\lesssim\delta v_{i,j}^{\star},
\end{align*}
provided that $\theta\lesssim1$, $np\gtrsim\delta^{-2}\kappa_{\omega}^{2}\log^{2}(n+d)$,
$n\gtrsim\delta^{-2}\kappa^{4}\mu^{2}r^{3}\log(n+d)$ and 
\[
\frac{\zeta_{\mathsf{1st}}}{\sigma_{r}^{\star2}}\lesssim\frac{\delta}{\sqrt{\kappa^{2}\mu r^{2}\log\left(n+d\right)}}.
\]
\end{itemize}
Taking the bounds on $\gamma_{4,1,2,1}$ to $\gamma_{4,1,2,6}$ collectively
yields 
\[
\gamma_{4,1,2}\lesssim\gamma_{4,1,2,1}+\gamma_{4,1,2,2}+\gamma_{4,1,2,3}+\gamma_{4,1,2,4}+\gamma_{4,1,2,5}+\gamma_{4,1,2,6}\lesssim\delta v_{i,j}^{\star}.
\]
Then we can combine the bounds on $\gamma_{4,1,1}$ and $\gamma_{4,1,2}$
to arrive at 
\[
\gamma_{4,1}\lesssim\gamma_{4,1,1}+\gamma_{4,1,2}\lesssim\delta v_{i,j}^{\star},
\]
provided that $\theta\lesssim1$, $np\gtrsim\delta^{-2}\kappa^{2}\mu^{2}r^{2}\kappa_{\omega}^{2}\log^{2}(n+d)$,
$n\gtrsim\delta^{-2}\kappa^{4}\mu^{2}r^{3}\log(n+d)$ and 
\[
\frac{\zeta_{\mathsf{1st}}}{\sigma_{r}^{\star2}}\lesssim\frac{\delta}{\sqrt{\kappa^{2}\mu r^{2}\log\left(n+d\right)}},
\]
\[
\left\Vert \bm{U}_{i,\cdot}^{\star}\right\Vert _{2}+\left\Vert \bm{U}_{j,\cdot}^{\star}\right\Vert _{2}\gtrsim\delta^{-1/2}\frac{\sqrt{\kappa_{\omega}\log^{3}\left(n+d\right)}}{n^{1/4}p^{1/4}}\left[\frac{\omega_{\max}}{\sigma_{r}^{\star}}\sqrt{\frac{d}{np}}+\frac{\omega_{\max}^{2}}{p\sigma_{r}^{\star2}}\sqrt{\frac{d}{n}}\right]\sqrt{\frac{r}{d}},
\]
and 
\begin{align*}
\left\Vert \bm{U}_{i,\cdot}^{\star}\right\Vert _{2}+\left\Vert \bm{U}_{j,\cdot}^{\star}\right\Vert _{2} & \gtrsim\delta^{-1/2}\left(\sqrt{\frac{\zeta_{\mathsf{1st}}}{\sigma_{r}^{\star2}}}+\frac{\kappa^{1/2}\mu^{1/4}r^{1/4}}{n^{1/4}}\right)\Biggl[\frac{\omega_{\max}}{\sigma_{r}^{\star}}\sqrt{\frac{d\kappa\mu^{1/2}r\log^{5/2}\left(n+d\right)}{np}}\\
 & \qquad\qquad+\sqrt{\frac{\kappa^{2}\mu^{3/2}r^{2}\kappa_{\omega}\log^{7/2}\left(n+d\right)}{ndp^{2}}}+\frac{\omega_{\max}^{2}}{p\sigma_{r}^{\star2}}\sqrt{\frac{d\kappa\mu^{1/2}r\log^{5/2}\left(n+d\right)}{n}}\Biggr]\sqrt{\frac{r}{d}}.
\end{align*}

\subparagraph{Step 4.2: bounding $\gamma_{4,2}$.}

In view of (\ref{eq:ce-var-est-inter-6}), (\ref{eq:ce-var-est-inter-7})
and (\ref{eq:ce-var-est-inter-8}), we can develop the following upper
bound: 
\begin{align*}
\gamma_{4,2} & \lesssim\frac{\Delta_{i}}{np^{2}}\max_{1\leq k\leq d}\left[\omega_{k}^{2}+\left(1-p\right)S_{k,k}\right]\left\Vert \bm{U}_{j,\cdot}\right\Vert _{2}^{2}\\
 & \lesssim\frac{1}{np^{2}}\left(\sqrt{\frac{\log^{2}\left(n+d\right)}{np}}\omega_{i}^{\star2}+\mathsf{UB}_{i}\right)\left(\omega_{\max}^{2}+\frac{\mu r}{d}\sigma_{1}^{\star2}\right)\left(\left\Vert \bm{U}_{j,\cdot}^{\star}\right\Vert _{2}^{2}+\frac{\theta^{2}}{\kappa\sigma_{r}^{\star2}}\omega_{j}^{\star2}\right)\\
 & \lesssim\delta v_{i,j}^{\star}+\frac{1}{np^{2}}\mathsf{UB}_{i}\left(\omega_{\max}^{2}+\frac{\mu r}{d}\sigma_{1}^{\star2}\right)\left(\left\Vert \bm{U}_{j,\cdot}^{\star}\right\Vert _{2}^{2}+\frac{\theta^{2}}{\kappa\sigma_{r}^{\star2}}\omega_{j}^{\star2}\right)\\
 & \quad+\frac{1}{np^{2}}\sqrt{\frac{\log^{2}\left(n+d\right)}{np}}\omega_{i}^{\star2}\frac{\mu r}{d}\sigma_{1}^{\star2}\left(\left\Vert \bm{U}_{j,\cdot}^{\star}\right\Vert _{2}^{2}+\frac{\theta^{2}}{\kappa\sigma_{r}^{\star2}}\omega_{j}^{\star2}\right)\\
 & \lesssim\delta v_{i,j}^{\star}+\underbrace{\frac{1}{np^{2}}\mathsf{UB}_{i}\omega_{\max}^{2}\left\Vert \bm{U}_{j,\cdot}^{\star}\right\Vert _{2}^{2}}_{\eqqcolon\gamma_{4,2,1}}+\underbrace{\frac{1}{np^{2}}\mathsf{UB}_{i}\omega_{\max}^{2}\omega_{j}^{\star2}\frac{\theta^{2}}{\kappa\sigma_{r}^{\star2}}}_{\eqqcolon\gamma_{4,2,2}}+\underbrace{\frac{\mu r}{ndp^{2}}\mathsf{UB}_{i}\sigma_{1}^{\star2}\left\Vert \bm{U}_{j,\cdot}^{\star}\right\Vert _{2}^{2}}_{\eqqcolon\gamma_{4,2,3}}\\
 & \quad+\underbrace{\frac{\mu r}{ndp^{2}}\mathsf{UB}_{i}\theta^{2}\omega_{j}^{\star2}}_{\eqqcolon\gamma_{4,2,4}}+\underbrace{\frac{\mu r}{ndp^{2}}\sqrt{\frac{\log^{2}\left(n+d\right)}{np}}\omega_{i}^{\star2}\sigma_{1}^{\star2}\left\Vert \bm{U}_{j,\cdot}^{\star}\right\Vert _{2}^{2}}_{\eqqcolon\gamma_{4,2,5}}+\underbrace{\frac{\mu r}{ndp^{2}}\sqrt{\frac{\log^{2}\left(n+d\right)}{np}}\omega_{i}^{\star2}\omega_{j}^{\star2}\theta^{2}}_{\eqqcolon\gamma_{4,2,6}},
\end{align*}
where the penultimate relation holds since 
\[
\frac{1}{np^{2}}\sqrt{\frac{\log^{2}\left(n+d\right)}{np}}\omega_{i}^{\star2}\omega_{\max}^{2}\left(\left\Vert \bm{U}_{j,\cdot}^{\star}\right\Vert _{2}^{2}+\frac{\theta^{2}}{\kappa\sigma_{r}^{\star2}}\omega_{j}^{\star2}\right)\lesssim\gamma_{4,1}\lesssim\delta v_{i,j}^{\star}.
\]
In what follows, we shall bound the six terms from $\gamma_{4,2,1}$
to $\gamma_{4,2,6}$ separately. 
\begin{itemize}
\item Regarding $\gamma_{4,2,1}$, we can derive 
\begin{align*}
\gamma_{4,2,1} & \asymp\underbrace{\frac{1}{np^{2}}\left(\theta+\sqrt{\frac{\kappa^{3}r\log\left(n+d\right)}{n}}\right)\left\Vert \bm{U}_{i,\cdot}^{\star}\bm{\Sigma}^{\star}\right\Vert _{2}^{2}\omega_{\max}^{2}\left\Vert \bm{U}_{j,\cdot}^{\star}\right\Vert _{2}^{2}}_{\eqqcolon\gamma_{4,2,1,1}}+\underbrace{\frac{1}{np^{2}}\theta^{2}\omega_{i}^{\star2}\omega_{\max}^{2}\left\Vert \bm{U}_{j,\cdot}^{\star}\right\Vert _{2}^{2}}_{\eqqcolon\gamma_{4,2,1,2}}\\
 & \quad+\underbrace{\frac{1}{np^{2}}\theta\omega_{i}^{\star}\left\Vert \bm{U}_{i,\cdot}^{\star}\bm{\Sigma}^{\star}\right\Vert _{2}\omega_{\max}^{2}\left\Vert \bm{U}_{j,\cdot}^{\star}\right\Vert _{2}^{2}}_{\eqqcolon\gamma_{4,2,1,3}}+\underbrace{\frac{1}{np^{2}}\zeta_{\mathsf{2nd},i}\sigma_{1}^{\star}\left\Vert \bm{U}_{i,\cdot}^{\star}\bm{\Sigma}^{\star}\right\Vert _{2}\omega_{\max}^{2}\left\Vert \bm{U}_{j,\cdot}^{\star}\right\Vert _{2}^{2}}_{\eqqcolon\gamma_{4,2,1,4}}\\
 & \quad+\underbrace{\frac{1}{np^{2}}\zeta_{\mathsf{2nd},i}^{2}\sigma_{1}^{\star2}\omega_{\max}^{2}\left\Vert \bm{U}_{j,\cdot}^{\star}\right\Vert _{2}^{2}}_{\eqqcolon\gamma_{4,2,1,5}}\lesssim\delta v_{i,j}^{\star}.
\end{align*}
Here, the last relation holds since 
\begin{align*}
\gamma_{4,2,1,1}+\gamma_{4,2,1,3}+\gamma_{4,2,1,4}+\gamma_{4,2,1,5} & \lesssim\delta\frac{\sigma_{r}^{\star2}}{ndp^{2}}\left(\omega_{j}^{\star2}\left\Vert \bm{U}_{i,\cdot}^{\star}\right\Vert _{2}^{2}+\omega_{i}^{\star2}\left\Vert \bm{U}_{j,\cdot}^{\star}\right\Vert _{2}^{2}\right)\lesssim\delta v_{i,j}^{\star},\\
\gamma_{4,2,1,2} & \lesssim\delta\frac{\omega_{\min}^{2}}{np^{2}}\left(\omega_{j}^{\star2}\left\Vert \bm{U}_{i,\cdot}^{\star}\right\Vert _{2}^{2}+\omega_{i}^{\star2}\left\Vert \bm{U}_{j,\cdot}^{\star}\right\Vert _{2}^{2}\right)\lesssim\delta v_{i,j}^{\star},
\end{align*}
provided that $\theta\lesssim\delta/(\kappa\mu r\kappa_{\omega})$,
$n\gtrsim\delta^{-2}\kappa^{5}\mu^{2}r^{3}\kappa_{\omega}^{2}\log(n+d)$,
$\zeta_{\mathsf{2nd},i}\sqrt{d}\lesssim\delta/\sqrt{\kappa^{2}\mu r\kappa_{\omega}^{2}}$,
\begin{equation}
\theta\frac{\omega_{\max}}{\sigma_{r}^{\star}}\sqrt{d}\lesssim\frac{\delta}{\sqrt{\kappa\mu r\kappa_{\omega}}}.\label{eq:theta-omega-sqrt-d}
\end{equation}
Note that in view of (\ref{eq:theta-exact}) and (\ref{eq:pca-1st-err-useful}),
(\ref{eq:theta-omega-sqrt-d}) is guaranteed by $ndp^{2}\gtrsim\delta^{-2}\kappa^{4}\mu^{3}r^{4}\kappa_{\omega}\log^{4}(n+d)$,
\[
\frac{\omega_{\max}}{\sigma_{r}^{\star}}\sqrt{\frac{d}{np}}\lesssim\frac{\delta}{\sqrt{\kappa^{2}\mu r^{2}\kappa_{\omega}\log^{2}\left(n+d\right)}},\qquad\frac{\omega_{\max}^{2}}{p\sigma_{r}^{\star2}}\sqrt{\frac{d}{n}}\lesssim\frac{\delta}{\sqrt{\kappa^{2}\mu r^{2}\kappa_{\omega}\log^{2}\left(n+d\right)}}.
\]
\item Regarding $\gamma_{4,2,2}$, we have 
\begin{align*}
\gamma_{4,2,2} & \asymp\underbrace{\frac{1}{np^{2}}\left(\theta+\sqrt{\frac{\kappa^{3}r\log\left(n+d\right)}{n}}\right)\left\Vert \bm{U}_{i,\cdot}^{\star}\bm{\Sigma}^{\star}\right\Vert _{2}^{2}\omega_{\max}^{2}\omega_{j}^{\star2}\frac{\theta^{2}}{\kappa\sigma_{r}^{\star2}}}_{\eqqcolon\gamma_{4,2,2,1}}+\underbrace{\frac{1}{np^{2}}\theta^{2}\omega_{\max}^{2}\omega_{i}^{\star2}\omega_{j}^{\star2}\frac{\theta^{2}}{\kappa\sigma_{r}^{\star2}}}_{\eqqcolon\gamma_{4,2,2,2}}\\
 & \quad+\underbrace{\frac{1}{np^{2}}\theta\omega_{\max}^{2}\omega_{i}^{\star}\omega_{j}^{\star2}\left\Vert \bm{U}_{i,\cdot}^{\star}\bm{\Sigma}^{\star}\right\Vert _{2}\frac{\theta^{2}}{\kappa\sigma_{r}^{\star2}}}_{\eqqcolon\gamma_{4,2,2,3}}+\underbrace{\frac{1}{np^{2}}\zeta_{\mathsf{2nd},i}\omega_{\max}^{2}\omega_{j}^{\star2}\left\Vert \bm{U}_{i,\cdot}^{\star}\bm{\Sigma}^{\star}\right\Vert _{2}\frac{\theta^{2}}{\sqrt{\kappa}\sigma_{r}^{\star}}}_{\eqqcolon\gamma_{4,2,2,4}}\\
 & \quad+\underbrace{\frac{1}{np^{2}}\zeta_{\mathsf{2nd},i}^{2}\omega_{\max}^{2}\omega_{j}^{\star2}\theta^{2}}_{\eqqcolon\gamma_{4,2,2,5}}\lesssim\delta v_{i,j}^{\star}.
\end{align*}
Here, the last relation holds since 
\begin{align*}
\gamma_{4,2,2,1} & \lesssim\delta\frac{\omega_{\min}^{2}}{np^{2}}\left(\omega_{j}^{\star2}\left\Vert \bm{U}_{i,\cdot}^{\star}\right\Vert _{2}^{2}+\omega_{i}^{\star2}\left\Vert \bm{U}_{j,\cdot}^{\star}\right\Vert _{2}^{2}\right)\lesssim\delta v_{i,j}^{\star},\\
\gamma_{4,2,2,2} & \lesssim\frac{\omega_{\max}^{2}}{\sigma_{r}^{\star2}}\frac{1}{np^{2}\kappa}\cdot\theta^{4}\omega_{i}^{\star2}\omega_{j}^{\star2}\overset{\text{(i)}}{\lesssim}\frac{\omega_{\max}^{2}}{\sigma_{r}^{\star2}}\frac{1}{np^{2}\kappa}\varepsilon^{2}v_{i,j}^{\star}\lesssim\delta v_{i,j}^{\star},\\
\gamma_{4,2,2,3} & \lesssim\theta\frac{\omega_{\max}}{\sqrt{\kappa}\sigma_{r}^{\star}}\sqrt{\frac{\kappa_{\omega}}{np^{2}}}\cdot\frac{\omega_{\min}\omega_{j}^{\star}}{\sqrt{np^{2}}}\left\Vert \bm{U}_{i,\cdot}^{\star}\right\Vert _{2}\cdot\theta^{2}\omega_{i}^{\star}\omega_{j}^{\star}\overset{\text{(ii)}}{\lesssim}\varepsilon\theta\frac{\omega_{\max}}{\sqrt{\kappa}\sigma_{r}^{\star}}\sqrt{\frac{\kappa_{\omega}}{np^{2}}}v_{i,j}^{\star}\lesssim\delta v_{i,j}^{\star},\\
\gamma_{4,2,2,4} & \lesssim\frac{\kappa_{\omega}}{\sqrt{ndp^{2}}}\zeta_{\mathsf{2nd},i}\sqrt{d}\cdot\theta^{2}\omega_{i}^{\star}\omega_{j}^{\star}\cdot\frac{\omega_{\min}\omega_{j}^{\star}}{\sqrt{np^{2}}}\left\Vert \bm{U}_{i,\cdot}^{\star}\right\Vert _{2}\overset{\text{(iii)}}{\lesssim}\frac{\varepsilon\kappa_{\omega}}{\sqrt{ndp^{2}}}\zeta_{\mathsf{2nd},i}\sqrt{d}v_{i,j}^{\star}\lesssim\delta v_{i,j}^{\star},\\
\gamma_{4,2,2,5} & \lesssim\frac{\omega_{\max}^{2}}{\sigma_{1}^{\star2}}\frac{1}{np^{2}}\cdot\zeta_{\mathsf{2nd},i}^{2}\sigma_{1}^{\star2}\cdot\theta^{2}\omega_{\max}^{2}\lesssim\frac{\omega_{\max}^{2}}{\sigma_{1}^{\star2}}\frac{\kappa_{\omega}}{np^{2}}\varepsilon^{2}v_{i,j}^{\star}\lesssim\delta v_{i,j}^{\star},
\end{align*}
provided that $\theta\lesssim\delta/\kappa_{\omega}$, $n\gtrsim\kappa^{3}r\log(n+d)$,
$\varepsilon\lesssim1$, $\zeta_{\mathsf{2nd},i}\sqrt{d}\lesssim\delta$,
$ndp^{2}\gtrsim\kappa_{\omega}^{2}$, 
\[
\frac{\omega_{\max}^{2}}{\sigma_{r}^{\star2}}\frac{1}{np^{2}\kappa}\lesssim\frac{\delta}{\kappa_{\omega}}.
\]
Here, (i)-(iii) rely on (\ref{eq:ce-var-est-inter-3}). In view of
the fact that 
\[
\frac{1}{\sqrt{ndp^{2}}}\cdot\frac{\omega_{\max}^{2}}{p\sigma_{r}^{\star2}}\sqrt{\frac{d}{n}}=\frac{\omega_{\max}^{2}}{\sigma_{r}^{\star2}}\frac{1}{np^{2}},
\]
the last condition above is guaranteed by $ndp^{2}\gtrsim\kappa_{\omega}^{2}$
and 
\[
\frac{\omega_{\max}^{2}}{p\sigma_{r}^{\star2}}\sqrt{\frac{d}{n}}\lesssim\delta.
\]
\item Regarding $\gamma_{4,2,3}$, we have 
\begin{align*}
\gamma_{4,2,3} & \asymp\underbrace{\frac{\mu r}{ndp^{2}}\sigma_{1}^{\star2}\left(\theta+\sqrt{\frac{\kappa^{3}r\log\left(n+d\right)}{n}}\right)\left\Vert \bm{U}_{i,\cdot}^{\star}\bm{\Sigma}^{\star}\right\Vert _{2}^{2}\left\Vert \bm{U}_{j,\cdot}^{\star}\right\Vert _{2}^{2}}_{\eqqcolon\gamma_{4,2,3,1}}+\underbrace{\frac{\mu r}{ndp^{2}}\sigma_{1}^{\star2}\theta^{2}\omega_{i}^{\star2}\left\Vert \bm{U}_{j,\cdot}^{\star}\right\Vert _{2}^{2}}_{\eqqcolon\gamma_{4,2,3,2}}\\
 & \quad+\underbrace{\frac{\mu r}{ndp^{2}}\sigma_{1}^{\star2}\theta\omega_{i}^{\star}\left\Vert \bm{U}_{i,\cdot}^{\star}\bm{\Sigma}^{\star}\right\Vert _{2}\left\Vert \bm{U}_{j,\cdot}^{\star}\right\Vert _{2}^{2}}_{\eqqcolon\gamma_{4,2,3,3}}+\underbrace{\frac{\mu r}{ndp^{2}}\sigma_{1}^{\star3}\zeta_{\mathsf{2nd},i}\left\Vert \bm{U}_{i,\cdot}^{\star}\bm{\Sigma}^{\star}\right\Vert _{2}\left\Vert \bm{U}_{j,\cdot}^{\star}\right\Vert _{2}^{2}}_{\eqqcolon\gamma_{4,2,3,4}}\\
 & \quad+\underbrace{\frac{\mu r}{ndp^{2}}\sigma_{1}^{\star2}\zeta_{\mathsf{2nd},i}^{2}\sigma_{1}^{\star2}\left\Vert \bm{U}_{j,\cdot}^{\star}\right\Vert _{2}^{2}}_{\eqqcolon\gamma_{4,2,3,5}}\lesssim\delta v_{i,j}^{\star}.
\end{align*}
Here, the last relation holds since 
\begin{align*}
\gamma_{4,2,3,1} & \lesssim\frac{\delta}{ndp^{2}\kappa}\left\Vert \bm{U}_{i,\cdot}^{\star}\bm{\Sigma}^{\star}\right\Vert _{2}^{2}\left\Vert \bm{U}_{j,\cdot}^{\star}\bm{\Sigma}^{\star}\right\Vert _{2}^{2}\lesssim\delta v_{i,j}^{\star},\\
\gamma_{4,2,3,2} & \lesssim\delta\frac{\sigma_{r}^{\star2}}{ndp^{2}}\left(\omega_{j}^{\star2}\left\Vert \bm{U}_{i,\cdot}^{\star}\right\Vert _{2}^{2}+\omega_{i}^{\star2}\left\Vert \bm{U}_{j,\cdot}^{\star}\right\Vert _{2}^{2}\right)\lesssim\delta v_{i,j}^{\star},\\
\gamma_{4,2,3,3} & \overset{\text{(i)}}{\lesssim}\delta\frac{1}{ndp^{2}\kappa}\left\Vert \bm{U}_{i,\cdot}^{\star}\bm{\Sigma}^{\star}\right\Vert _{2}^{2}\left\Vert \bm{U}_{j,\cdot}^{\star}\bm{\Sigma}^{\star}\right\Vert _{2}^{2}+\delta\frac{\sigma_{r}^{\star2}}{ndp^{2}}\omega_{i}^{\star2}\left\Vert \bm{U}_{j,\cdot}^{\star}\right\Vert _{2}^{2}\lesssim\delta v_{i,j}^{\star},\\
\gamma_{4,2,3,4} & \lesssim\frac{\kappa\mu r}{\sqrt{ndp^{2}}}\cdot\frac{1}{\sqrt{ndp^{2}\kappa}}\left\Vert \bm{U}_{i,\cdot}^{\star}\bm{\Sigma}^{\star}\right\Vert _{2}\left\Vert \bm{U}_{j,\cdot}^{\star}\bm{\Sigma}^{\star}\right\Vert _{2}\cdot\sigma_{1}^{\star2}\left(\left\Vert \bm{U}_{i,\cdot}^{\star}\right\Vert _{2}\zeta_{\mathsf{2nd},j}+\left\Vert \bm{U}_{j,\cdot}^{\star}\right\Vert _{2}\zeta_{\mathsf{2nd},i}\right)\\
 & \overset{\text{(ii)}}{\lesssim}\frac{\kappa\mu r}{\sqrt{ndp^{2}}}\cdot\frac{1}{\sqrt{ndp^{2}\kappa}}\left\Vert \bm{U}_{i,\cdot}^{\star}\bm{\Sigma}^{\star}\right\Vert _{2}\left\Vert \bm{U}_{j,\cdot}^{\star}\bm{\Sigma}^{\star}\right\Vert _{2}\cdot\varepsilon(v_{i,j}^{\star})^{1/2}\lesssim\frac{\kappa\mu r}{\sqrt{ndp^{2}}}\varepsilon v_{i,j}^{\star}\lesssim\delta v_{i,j}^{\star},\\
\gamma_{4,2,3,5} & \lesssim\frac{\mu r}{ndp^{2}}\left[\sigma_{1}^{\star2}\left(\left\Vert \bm{U}_{i,\cdot}^{\star}\right\Vert _{2}\zeta_{\mathsf{2nd},j}+\left\Vert \bm{U}_{j,\cdot}^{\star}\right\Vert _{2}\zeta_{\mathsf{2nd},i}\right)\right]^{2}\overset{\text{(iii)}}{\lesssim}\frac{\mu r}{ndp^{2}}\varepsilon v_{i,j}^{\star}\lesssim\delta v_{i,j}^{\star},
\end{align*}
provided that $\theta\lesssim\delta/(\kappa^{2}\mu r)$, $n\gtrsim\delta^{-2}\kappa^{7}\mu^{2}r^{3}\log(n+d)$,
$\varepsilon\lesssim1$ and $ndp^{2}\gtrsim\delta^{-2}\kappa^{2}\mu^{2}r^{2}$.
Here, (i) utilizes the AM-GM inequality, whereas (ii) and (iii) utilize
(\ref{eq:ce-var-est-inter-2}). 
\item Regarding $\gamma_{4,2,4}$, it follows that 
\begin{align*}
\gamma_{4,2,4} & \asymp\underbrace{\frac{\mu r}{ndp^{2}}\theta^{2}\omega_{j}^{\star2}\left(\theta+\sqrt{\frac{\kappa^{3}r\log\left(n+d\right)}{n}}\right)\left\Vert \bm{U}_{i,\cdot}^{\star}\bm{\Sigma}^{\star}\right\Vert _{2}^{2}}_{\eqqcolon\gamma_{4,2,4,1}}+\underbrace{\frac{\mu r}{ndp^{2}}\theta^{4}\omega_{i}^{\star2}\omega_{j}^{\star2}}_{\eqqcolon\gamma_{4,2,4,2}}\\
 & \quad+\underbrace{\frac{\mu r}{ndp^{2}}\theta^{3}\omega_{i}^{\star}\omega_{j}^{\star2}\left\Vert \bm{U}_{i,\cdot}^{\star}\bm{\Sigma}^{\star}\right\Vert _{2}}_{\eqqcolon\gamma_{4,2,4,3}}+\underbrace{\frac{\mu r}{ndp^{2}}\theta^{2}\omega_{j}^{\star2}\sigma_{1}^{\star}\zeta_{\mathsf{2nd},i}\left\Vert \bm{U}_{i,\cdot}^{\star}\bm{\Sigma}^{\star}\right\Vert _{2}}_{\eqqcolon\gamma_{4,2,4,4}}\\
 & \quad+\underbrace{\frac{\mu r}{ndp^{2}}\theta^{2}\omega_{j}^{\star2}\zeta_{\mathsf{2nd},i}^{2}\sigma_{1}^{\star2}}_{\eqqcolon\gamma_{4,2,4,5}}\lesssim\delta v_{i,j}^{\star},
\end{align*}
where the last relation holds since 
\begin{align*}
\gamma_{4,2,4,1} & \lesssim\delta\frac{\sigma_{r}^{\star2}}{ndp^{2}}\left(\omega_{j}^{\star2}\left\Vert \bm{U}_{i,\cdot}^{\star}\right\Vert _{2}^{2}+\omega_{i}^{\star2}\left\Vert \bm{U}_{j,\cdot}^{\star}\right\Vert _{2}^{2}\right)\lesssim\delta v_{i,j}^{\star},\\
\gamma_{4,2,4,2} & \overset{\text{(i)}}{\lesssim}\frac{\mu r}{ndp^{2}}\varepsilon^{2}v_{i,j}^{\star}\lesssim\delta v_{i,j}^{\star},\\
\gamma_{4,2,4,3} & \overset{\text{(ii)}}{\lesssim}\frac{\sqrt{\kappa}\mu r}{\sqrt{ndp^{2}}}\theta\cdot\theta^{2}\omega_{i}^{\star}\omega_{j}^{\star}\cdot\frac{\sigma_{r}^{\star}}{\sqrt{ndp^{2}}}\omega_{j}^{\star}\left\Vert \bm{U}_{i,\cdot}^{\star}\right\Vert _{2}\lesssim\frac{\sqrt{\kappa}\mu r}{\sqrt{ndp^{2}}}\varepsilon\theta v_{i,j}^{\star}\lesssim\delta v_{i,j}^{\star},\\
\gamma_{4,2,4,4} & \lesssim\frac{\mu r\kappa_{\omega}^{1/2}}{ndp^{2}}\cdot\theta^{2}\omega_{i}^{\star}\omega_{j}^{\star}\cdot\sigma_{1}^{\star2}\left(\left\Vert \bm{U}_{i,\cdot}^{\star}\right\Vert _{2}\zeta_{\mathsf{2nd},j}+\left\Vert \bm{U}_{j,\cdot}^{\star}\right\Vert _{2}\zeta_{\mathsf{2nd},i}\right)\overset{\text{(iii)}}{\lesssim}\frac{\mu r\kappa_{\omega}^{1/2}\varepsilon^{2}}{ndp^{2}}v_{i,j}^{\star}\lesssim\delta v_{i,j}^{\star},\\
\gamma_{4,2,4,5} & \lesssim\frac{\mu r}{ndp^{2}}\cdot\left(\theta\sigma_{1}^{\star}\omega_{j}^{\star}\zeta_{\mathsf{2nd},i}\right)^{2}\overset{\text{(iv)}}{\lesssim}\frac{\mu r}{ndp^{2}}\varepsilon^{2}v_{i,j}^{\star}\lesssim\delta v_{i,j}^{\star},
\end{align*}
as long as $\theta\lesssim\delta/(\kappa\mu r)$, $n\gtrsim\kappa^{3}r\log(n+d)$,
$ndp^{2}\gtrsim\delta^{-2}\mu r$, $\varepsilon\lesssim1$. Here,
(i), (ii) and (iv) utilize (\ref{eq:ce-var-est-inter-3}); (iii) utilizes
(\ref{eq:ce-var-est-inter-1}). 
\item Regarding $\gamma_{4,2,5}$, we can derive 
\[
\gamma_{4,2,5}\lesssim\delta\frac{\sigma_{r}^{\star2}}{ndp^{2}}\left(\omega_{j}^{\star2}\left\Vert \bm{U}_{i,\cdot}^{\star}\right\Vert _{2}^{2}+\omega_{i}^{\star2}\left\Vert \bm{U}_{j,\cdot}^{\star}\right\Vert _{2}^{2}\right)\lesssim\delta v_{i,j}^{\star}
\]
provided that $np\gtrsim\delta^{-2}\kappa^{2}\mu^{2}r^{2}\log^{2}(n+d)$. 
\item Finally, when it comes to $\gamma_{4,2,6}$, we can bound 
\[
\gamma_{4,2,6}\overset{\text{(i)}}{\lesssim}\theta^{2}\sqrt{\frac{1}{np}}\omega_{i}^{\star2}\omega_{j}^{\star2}\theta^{2}\overset{\text{(ii)}}{\lesssim}\sqrt{\frac{1}{np}}v_{i,j}^{\star}\overset{\text{(iii)}}{\lesssim}\delta v_{i,j}^{\star}.
\]
Here (i) follows from (\ref{eq:theta-exact}), (ii) comes from (\ref{eq:ce-var-est-inter-3}),
whereas (iii) holds provided that $np\gtrsim\delta^{-2}$. 
\end{itemize}
Take the bounds on the terms (from $\gamma_{4,2,1}$ to $\gamma_{4,2,6}$)
together to arrive at 
\[
\gamma_{4,2}\lesssim\gamma_{4,2,1}+\gamma_{4,2,2}+\gamma_{4,2,3}+\gamma_{4,2,4}+\gamma_{4,2,5}+\gamma_{4,2,6}\lesssim\delta v_{i,j}^{\star},
\]
provided that $\theta\lesssim\delta/(\kappa^{2}\mu r\kappa_{\omega})$,
$\varepsilon\lesssim1$, $\zeta_{\mathsf{2nd},i}\sqrt{d}\lesssim\delta/\sqrt{\kappa^{2}\mu r\kappa_{\omega}^{2}}$,
$n\gtrsim\delta^{-2}\kappa^{7}\mu^{2}r^{3}\kappa_{\omega}^{2}\log(n+d)$,
$ndp^{2}\gtrsim\delta^{-2}\kappa^{4}\mu^{3}r^{4}\kappa_{\omega}^{2}\log^{4}(n+d)$,
$np\gtrsim\delta^{-2}\kappa^{2}\mu^{2}r^{2}\log^{2}(n+d)$ 
\[
\frac{\omega_{\max}}{\sigma_{r}^{\star}}\sqrt{\frac{d}{np}}\lesssim\frac{\delta}{\sqrt{\kappa^{2}\mu r^{2}\kappa_{\omega}\log^{2}\left(n+d\right)}}\qquad\text{and}\qquad\frac{\omega_{\max}^{2}}{p\sigma_{r}^{\star2}}\sqrt{\frac{d}{n}}\lesssim\frac{\delta}{\sqrt{\kappa^{2}\mu r^{2}\kappa_{\omega}\log^{2}\left(n+d\right)}}.
\]

\subparagraph{Step 4.3: bounding $\gamma_{4,3}$.}

In view of (\ref{eq:ce-var-est-inter-7}), we can upper bound 
\begin{align*}
\gamma_{4,3} & \lesssim\frac{1}{np^{2}}\left(\omega_{i}^{\star2}+\left\Vert \bm{U}_{i,\cdot}^{\star}\bm{\Sigma}^{\star}\right\Vert _{2}^{2}\right)\left(\omega_{\max}^{2}+\frac{\mu r}{d}\sigma_{1}^{\star2}\right)\sum_{k=1}^{d}\left|\left(\bm{U}_{k,\cdot}^{\star}\bm{U}_{j,\cdot}^{\star\top}\right)^{2}-\left(\bm{U}_{k,\cdot}\bm{U}_{j,\cdot}^{\top}\right)^{2}\right|\\
 & \lesssim\underbrace{\frac{\omega_{i}^{\star2}\omega_{\max}^{2}}{np^{2}}\sum_{k=1}^{d}\left|\left(\bm{U}_{k,\cdot}^{\star}\bm{U}_{j,\cdot}^{\star\top}\right)^{2}-\left(\bm{U}_{k,\cdot}\bm{U}_{j,\cdot}^{\top}\right)^{2}\right|}_{\eqqcolon\gamma_{4,3,1}}+\underbrace{\frac{\mu r}{ndp^{2}}\sigma_{1}^{\star2}\omega_{i}^{\star2}\sum_{k=1}^{d}\left|\left(\bm{U}_{k,\cdot}^{\star}\bm{U}_{j,\cdot}^{\star\top}\right)^{2}-\left(\bm{U}_{k,\cdot}\bm{U}_{j,\cdot}^{\top}\right)^{2}\right|}_{\eqqcolon\gamma_{4,3,2}}\\
 & \quad+\underbrace{\frac{\mu r}{ndp^{2}}\sigma_{1}^{\star2}\left\Vert \bm{U}_{i,\cdot}^{\star}\bm{\Sigma}^{\star}\right\Vert _{2}^{2}\sum_{k=1}^{d}\left|\left(\bm{U}_{k,\cdot}^{\star}\bm{U}_{j,\cdot}^{\star\top}\right)^{2}-\left(\bm{U}_{k,\cdot}\bm{U}_{j,\cdot}^{\top}\right)^{2}\right|}_{\eqqcolon\gamma_{4,3,3}}\\
 & \quad+\underbrace{\frac{1}{np^{2}}\left\Vert \bm{U}_{i,\cdot}^{\star}\bm{\Sigma}^{\star}\right\Vert _{2}^{2}\omega_{\max}^{2}\sum_{k=1}^{d}\left|\left(\bm{U}_{k,\cdot}^{\star}\bm{U}_{j,\cdot}^{\star\top}\right)^{2}-\left(\bm{U}_{k,\cdot}\bm{U}_{j,\cdot}^{\top}\right)^{2}\right|}_{\eqqcolon\gamma_{4,3,4}}.
\end{align*}
Note that for each $k\in[d]$, 
\begin{align*}
\left|\bm{U}_{k,\cdot}^{\star}\bm{U}_{j,\cdot}^{\star\top}-\bm{U}_{k,\cdot}\bm{U}_{j,\cdot}^{\top}\right| & =\left|\bm{U}_{k,\cdot}^{\star}\bm{U}_{j,\cdot}^{\star\top}-\bm{U}_{k,\cdot}\bm{R}\left(\bm{U}_{j,\cdot}\bm{R}\right)^{\top}\right|\\
 & \leq\left|\left(\bm{U}\bm{R}-\bm{U}^{\star}\right)_{k,\cdot}\bm{U}_{j,\cdot}^{\star}\right|+\left|\left(\bm{U}\bm{R}\right)_{k,\cdot}\left(\bm{U}\bm{R}-\bm{U}^{\star}\right)_{j,\cdot}^{\top}\right|\\
 & \leq\left\Vert \bm{U}_{j,\cdot}^{\star}\right\Vert _{2}\left\Vert \bm{U}\bm{R}-\bm{U}^{\star}\right\Vert _{2,\infty}+\left\Vert \bm{U}_{k,\cdot}\right\Vert _{2}\left\Vert \left(\bm{U}\bm{R}-\bm{U}^{\star}\right)_{j,\cdot}\right\Vert _{2}\\
 & \overset{\text{(i)}}{\lesssim}\left\Vert \bm{U}_{j,\cdot}^{\star}\right\Vert _{2}\frac{\zeta_{\mathsf{1st}}}{\sigma_{r}^{\star2}}\sqrt{\frac{r\log\left(n+d\right)}{d}}+\left\Vert \bm{U}_{k,\cdot}\right\Vert _{2}\left[\frac{\theta}{\sqrt{\kappa}\sigma_{r}^{\star}}\left(\left\Vert \bm{U}_{j,\cdot}^{\star}\bm{\Sigma}^{\star}\right\Vert _{2}+\omega_{j}^{\star}\right)\right]\\
 & \overset{\text{(ii)}}{\lesssim}\left\Vert \bm{U}_{j,\cdot}^{\star}\right\Vert _{2}\frac{\zeta_{\mathsf{1st}}}{\sigma_{r}^{\star2}}\sqrt{\frac{r\log\left(n+d\right)}{d}}+\left\Vert \bm{U}_{k,\cdot}\right\Vert _{2}\frac{\theta}{\sqrt{\kappa}\sigma_{r}^{\star}}\omega_{j}^{\star}.
\end{align*}
Here, (i) follows from Lemma \ref{lemma:pca-1st-err} as well as a
direct consequence of Lemma \ref{lemma:pca-1st-err}, (\ref{eq:ce-var-est-inter-2})
and (\ref{eq:theta-exact}): 
\begin{align*}
\left\Vert \left(\bm{U}\bm{R}-\bm{U}^{\star}\right)_{j,\cdot}\right\Vert _{2} & \lesssim\frac{\theta}{\sqrt{\kappa}\sigma_{r}^{\star}}\left(\left\Vert \bm{U}_{j,\cdot}^{\star}\bm{\Sigma}^{\star}\right\Vert _{2}+\omega_{j}^{\star}\right)+\zeta_{\mathsf{2nd},j}\\
 & \lesssim\frac{\theta}{\sqrt{\kappa}\sigma_{r}^{\star}}\left(\left\Vert \bm{U}_{j,\cdot}^{\star}\bm{\Sigma}^{\star}\right\Vert _{2}+\omega_{j}^{\star}\right)+\frac{\varepsilon}{\sqrt{\kappa}}\frac{\left\Vert \bm{U}_{j,\cdot}^{\star}\bm{\Sigma}^{\star}\right\Vert _{2}+\omega_{j}^{\star}}{\sqrt{\min\left\{ ndp^{2}\kappa,np\right\} }\sigma_{1}^{\star}}+\frac{\varepsilon}{\sqrt{\kappa}}\frac{\omega_{\min}\omega_{j}^{\star}}{\sqrt{np^{2}}\sigma_{1}^{\star2}}\\
 & \lesssim\frac{\theta}{\sqrt{\kappa}\sigma_{r}^{\star}}\left(\left\Vert \bm{U}_{j,\cdot}^{\star}\bm{\Sigma}^{\star}\right\Vert _{2}+\omega_{j}^{\star}\right);
\end{align*}
and (ii) utilizes (\ref{eq:theta-zeta-1st}). Therefore, one can derive
\begin{align}
\sum_{k=1}^{d}\left|\left(\bm{U}_{k,\cdot}^{\star}\bm{U}_{j,\cdot}^{\star\top}\right)^{2}-\left(\bm{U}_{k,\cdot}\bm{U}_{j,\cdot}^{\top}\right)^{2}\right| & \lesssim\sum_{k=1}^{d}\left|\bm{U}_{k,\cdot}^{\star}\bm{U}_{j,\cdot}^{\star\top}-\bm{U}_{k,\cdot}\bm{U}_{j,\cdot}^{\top}\right|\left\Vert \bm{U}_{j,\cdot}^{\star}\right\Vert _{2}\left\Vert \bm{U}_{k,\cdot}^{\star}\right\Vert _{2}+\sum_{k=1}^{d}\left|\bm{U}_{k,\cdot}^{\star}\bm{U}_{j,\cdot}^{\star\top}-\bm{U}_{k,\cdot}\bm{U}_{j,\cdot}^{\top}\right|^{2}\nonumber \\
 & \lesssim\sum_{k=1}^{d}\left\Vert \bm{U}_{j,\cdot}^{\star}\right\Vert _{2}^{2}\frac{\zeta_{\mathsf{1st}}}{\sigma_{r}^{\star2}}\sqrt{\frac{r\log\left(n+d\right)}{d}}\left\Vert \bm{U}_{k,\cdot}^{\star}\right\Vert _{2}+\sum_{k=1}^{d}\left\Vert \bm{U}_{k,\cdot}\right\Vert _{2}\left\Vert \bm{U}_{k,\cdot}^{\star}\right\Vert _{2}\frac{\theta}{\sqrt{\kappa}\sigma_{r}^{\star}}\omega_{j}^{\star}\left\Vert \bm{U}_{j,\cdot}^{\star}\right\Vert _{2}\nonumber \\
 & \quad+\sum_{k=1}^{d}\left\Vert \bm{U}_{j,\cdot}^{\star}\right\Vert _{2}^{2}\frac{\zeta_{\mathsf{1st}}^{2}}{\sigma_{r}^{\star4}}\frac{r\log\left(n+d\right)}{d}+\sum_{k=1}^{d}\left\Vert \bm{U}_{k,\cdot}\right\Vert _{2}^{2}\frac{\theta^{2}}{\kappa\sigma_{r}^{\star2}}\omega_{j}^{\star2}\nonumber \\
 & \overset{\text{(i)}}{\lesssim}\left\Vert \bm{U}_{j,\cdot}^{\star}\right\Vert _{2}^{2}\frac{\zeta_{\mathsf{1st}}}{\sigma_{r}^{\star2}}\sqrt{r^{2}\log\left(n+d\right)}+\frac{\theta r}{\sqrt{\kappa}\sigma_{r}^{\star}}\omega_{j}^{\star}\left\Vert \bm{U}_{j,\cdot}^{\star}\right\Vert _{2}\nonumber \\
 & \quad+\left\Vert \bm{U}_{j,\cdot}^{\star}\right\Vert _{2}^{2}\frac{\zeta_{\mathsf{1st}}^{2}}{\sigma_{r}^{\star4}}r^{2}\log\left(n+d\right)+\frac{\theta^{2}r}{\kappa\sigma_{r}^{\star2}}\omega_{j}^{\star2}\nonumber \\
 & \overset{\text{(ii)}}{\lesssim}\left\Vert \bm{U}_{j,\cdot}^{\star}\right\Vert _{2}^{2}\frac{\zeta_{\mathsf{1st}}}{\sigma_{r}^{\star2}}\sqrt{r^{2}\log\left(n+d\right)}+\frac{\theta r}{\sqrt{\kappa}\sigma_{r}^{\star}}\omega_{j}^{\star}\left\Vert \bm{U}_{j,\cdot}^{\star}\right\Vert _{2}+\frac{\theta^{2}r}{\kappa\sigma_{r}^{\star2}}\omega_{j}^{\star2}.\label{eq:ce-var-est-inter-9}
\end{align}
Here, (i) follows from the Cauchy-Schwarz inequality 
\[
\sum_{k=1}^{d}\left\Vert \bm{U}_{k,\cdot}^{\star}\right\Vert _{2}\leq\sqrt{d\sum_{k=1}^{d}\left\Vert \bm{U}_{k,\cdot}^{\star}\right\Vert _{2}^{2}}\leq\sqrt{d\left\Vert \bm{U}^{\star}\right\Vert _{\mathrm{F}}^{2}}\leq\sqrt{dr},
\]
as well as the following bound 
\begin{align*}
\sum_{k=1}^{d}\left\Vert \bm{U}_{k,\cdot}\right\Vert _{2}\left\Vert \bm{U}_{k,\cdot}^{\star}\right\Vert _{2} & \leq\sum_{k=1}^{d}\left\Vert \bm{U}_{k,\cdot}^{\star}\right\Vert _{2}^{2}+\left\Vert \bm{U}\bm{R}-\bm{U}^{\star}\right\Vert _{2,\infty}\sum_{k=1}^{d}\left\Vert \bm{U}_{k,\cdot}^{\star}\right\Vert _{2}\\
 & \lesssim\left\Vert \bm{U}^{\star}\right\Vert _{\mathrm{F}}^{2}+\frac{\zeta_{\mathsf{1st}}}{\sigma_{r}^{\star2}}\sqrt{\frac{r\log\left(n+d\right)}{d}}\cdot\sqrt{dr}\lesssim r,
\end{align*}
which utilizes Lemma \ref{lemma:pca-1st-err} and holds provided that
$\zeta_{\mathsf{1st}}/\sigma_{r}^{\star2}\lesssim1/\sqrt{\log(n+d)}$;
(ii) holds provided that $\zeta_{\mathsf{1st}}/\sigma_{r}^{\star2}\lesssim1/\sqrt{\log(n+d)}$.
Then we shall bound the terms $\gamma_{4,3,1}$, $\gamma_{4,3,2}$,
$\gamma_{4,3,3}$ and $\gamma_{4,3,4}$ respectively. 
\begin{itemize}
\item Regarding $\gamma_{4,3,1}$, we have 
\begin{align*}
\gamma_{4,3,1} & \lesssim\underbrace{\frac{\sqrt{r^{2}\log\left(n+d\right)}}{np^{2}}\omega_{i}^{\star2}\omega_{\max}^{2}\left\Vert \bm{U}_{j,\cdot}^{\star}\right\Vert _{2}^{2}\frac{\zeta_{\mathsf{1st}}}{\sigma_{r}^{\star2}}}_{\eqqcolon\gamma_{4,3,1,1}}+\underbrace{\frac{r}{np^{2}}\omega_{i}^{\star2}\omega_{j}^{\star}\omega_{\max}^{2}\frac{\theta}{\sqrt{\kappa}\sigma_{r}^{\star}}\left\Vert \bm{U}_{j,\cdot}^{\star}\right\Vert _{2}}_{\eqqcolon\gamma_{4,3,1,2}}\\
 & \quad+\underbrace{\frac{r}{np^{2}}\omega_{i}^{\star2}\omega_{j}^{\star2}\omega_{\max}^{2}\frac{\theta^{2}}{\kappa\sigma_{r}^{\star2}}}_{\eqqcolon\gamma_{4,3,1,3}}\lesssim\delta v_{i,j}^{\star},
\end{align*}
where the last relation holds since 
\begin{align*}
\gamma_{4,3,1,1}+\gamma_{4,3,1,2} & \lesssim\delta\frac{\omega_{\min}^{2}}{np^{2}}\left(\omega_{j}^{\star2}\left\Vert \bm{U}_{i,\cdot}^{\star}\right\Vert _{2}^{2}+\omega_{i}^{\star2}\left\Vert \bm{U}_{j,\cdot}^{\star}\right\Vert _{2}^{2}\right)\lesssim\delta v_{i,j}^{\star},\\
\gamma_{4,3,1,3} & \lesssim\theta^{2}\omega_{i}^{\star}\omega_{j}^{\star}\cdot\frac{r}{np^{2}}\frac{\omega_{\max}^{2}\omega_{i}^{\star}\omega_{j}^{\star}}{\kappa\sigma_{r}^{\star2}}\overset{\text{(i)}}{\lesssim}\varepsilon(v_{i,j}^{\star})^{1/2}\cdot\frac{r}{np^{2}}\frac{\omega_{\max}^{2}\omega_{i}^{\star}\omega_{j}^{\star}}{\kappa\sigma_{r}^{\star2}}\\
 & \lesssim\varepsilon(v_{i,j}^{\star})^{1/2}\cdot\delta\frac{\omega_{\min}}{\sqrt{np^{2}}}\left(\omega_{j}^{\star}\left\Vert \bm{U}_{i,\cdot}^{\star}\right\Vert _{2}+\omega_{i}^{\star}\left\Vert \bm{U}_{j,\cdot}^{\star}\right\Vert _{2}\right)\lesssim\delta v_{i,j}^{\star},
\end{align*}
provided that $\zeta_{\mathsf{1st}}/\sigma_{r}^{\star2}\lesssim\delta/\sqrt{r^{2}\kappa_{\omega}^{2}\log(n+d)}$,
$\varepsilon\lesssim1$, 
\begin{align*}
\left\Vert \bm{U}_{i,\cdot}^{\star}\right\Vert _{2}+\left\Vert \bm{U}_{j,\cdot}^{\star}\right\Vert _{2} & \gtrsim\delta^{-1}r\kappa_{\omega}\frac{\omega_{\max}}{\sigma_{1}^{\star}}\theta+\delta^{-1}r\kappa_{\omega}^{1/2}\frac{\omega_{\max}^{2}}{\sigma_{1}^{\star2}}\sqrt{\frac{1}{np^{2}}}.
\end{align*}
Here (i) utilizes (\ref{eq:ce-var-est-inter-3}). 
\item Regarding $\gamma_{4,3,2}$, we have 
\begin{align*}
\gamma_{4,3,2} & \lesssim\underbrace{\frac{\mu r^{2}\sqrt{\log\left(n+d\right)}}{ndp^{2}}\sigma_{1}^{\star2}\omega_{i}^{\star2}\left\Vert \bm{U}_{j,\cdot}^{\star}\right\Vert _{2}^{2}\frac{\zeta_{\mathsf{1st}}}{\sigma_{r}^{\star2}}}_{\eqqcolon\gamma_{4,3,2,1}}+\underbrace{\frac{\mu r^{2}}{ndp^{2}}\sigma_{1}^{\star}\omega_{i}^{\star2}\omega_{j}^{\star}\theta\left\Vert \bm{U}_{j,\cdot}^{\star}\right\Vert _{2}}_{\eqqcolon\gamma_{4,3,2,2}}\\
 & \quad+\underbrace{\frac{\mu r^{2}}{ndp^{2}}\omega_{i}^{\star2}\omega_{j}^{\star2}\theta^{2}}_{\eqqcolon\gamma_{4,3,2,3}}\lesssim\delta v_{i,j}^{\star}.
\end{align*}
Here, the last relation holds since 
\begin{align*}
\gamma_{4,3,2,1}+\gamma_{4,3,2,2} & \lesssim\delta\frac{\sigma_{r}^{\star2}}{ndp^{2}}\left(\omega_{j}^{\star2}\left\Vert \bm{U}_{i,\cdot}^{\star}\right\Vert _{2}^{2}+\omega_{i}^{\star2}\left\Vert \bm{U}_{j,\cdot}^{\star}\right\Vert _{2}^{2}\right)\lesssim\delta v_{i,j}^{\star},\\
\gamma_{4,3,2,3} & \lesssim\frac{\mu r^{2}}{ndp^{2}}\omega_{i}^{\star}\omega_{j}^{\star}\cdot\theta^{2}\omega_{i}^{\star}\omega_{j}^{\star}\overset{\text{(i)}}{\lesssim}\frac{\mu r^{2}}{ndp^{2}}\omega_{i}^{\star}\omega_{j}^{\star}\cdot\varepsilon(v_{i,j}^{\star})^{1/2},\\
 & \lesssim\delta\frac{\omega_{\min}}{\sqrt{np^{2}}}\left(\omega_{j}^{\star}\left\Vert \bm{U}_{i,\cdot}^{\star}\right\Vert _{2}+\omega_{i}^{\star}\left\Vert \bm{U}_{j,\cdot}^{\star}\right\Vert _{2}\right)\cdot\varepsilon(v_{i,j}^{\star})^{1/2}\lesssim\delta v_{i,j}^{\star},
\end{align*}
provided that $\zeta_{\mathsf{1st}}/\sigma_{r}^{\star2}\lesssim\delta/\sqrt{\kappa^{2}\mu^{2}r^{4}\log(n+d)}$
and 
\[
\left\Vert \bm{U}_{i,\cdot}^{\star}\right\Vert _{2}+\left\Vert \bm{U}_{j,\cdot}^{\star}\right\Vert _{2}\gtrsim\delta^{-1}\sqrt{\kappa}\mu r^{2}\frac{\omega_{\max}}{\sigma_{r}^{\star}}\theta+\delta^{-1}\kappa_{\omega}^{1/2}\frac{\mu r^{2}}{\sqrt{ndp^{2}}}\cdot\sqrt{\frac{1}{d}}.
\]
Here (i) follows from (\ref{eq:ce-var-est-inter-3}). 
\item Regarding $\gamma_{4,3,3}$, we have 
\begin{align*}
\gamma_{4,3,3} & \lesssim\underbrace{\frac{\mu r^{2}\sqrt{\log\left(n+d\right)}}{ndp^{2}}\sigma_{1}^{\star2}\left\Vert \bm{U}_{i,\cdot}^{\star}\bm{\Sigma}^{\star}\right\Vert _{2}^{2}\left\Vert \bm{U}_{j,\cdot}^{\star}\right\Vert _{2}^{2}\frac{\zeta_{\mathsf{1st}}}{\sigma_{r}^{\star2}}}_{\eqqcolon\gamma_{4,3,3,1}}+\underbrace{\frac{\mu r^{2}}{ndp^{2}}\sigma_{1}^{\star}\omega_{j}^{\star}\theta\left\Vert \bm{U}_{i,\cdot}^{\star}\bm{\Sigma}^{\star}\right\Vert _{2}^{2}\left\Vert \bm{U}_{j,\cdot}^{\star}\right\Vert _{2}}_{\eqqcolon\gamma_{4,3,3,2}}\\
 & \quad+\underbrace{\frac{\mu r^{2}}{ndp^{2}}\omega_{j}^{\star2}\theta^{2}\left\Vert \bm{U}_{i,\cdot}^{\star}\bm{\Sigma}^{\star}\right\Vert _{2}^{2}}_{\eqqcolon\gamma_{4,3,3,3}}\lesssim\delta v_{i,j}^{\star},
\end{align*}
where the last relation holds since 
\begin{align*}
\gamma_{4,3,3,1} & \lesssim\frac{\delta}{ndp^{2}\kappa}\left\Vert \bm{U}_{i,\cdot}^{\star}\bm{\Sigma}^{\star}\right\Vert _{2}^{2}\left\Vert \bm{U}_{j,\cdot}^{\star}\bm{\Sigma}^{\star}\right\Vert _{2}^{2}\lesssim\delta v_{i,j}^{\star},\\
\gamma_{4,3,3,2} & \overset{\text{(i)}}{\lesssim}\frac{\delta}{ndp^{2}\kappa}\left\Vert \bm{U}_{i,\cdot}^{\star}\bm{\Sigma}^{\star}\right\Vert _{2}^{2}\left\Vert \bm{U}_{j,\cdot}^{\star}\bm{\Sigma}^{\star}\right\Vert _{2}^{2}+\delta\frac{\sigma_{r}^{\star2}}{ndp^{2}}\omega_{j}^{\star2}\left\Vert \bm{U}_{i,\cdot}^{\star}\right\Vert _{2}^{2}\lesssim\delta v_{i,j}^{\star},\\
\gamma_{4,3,3,3} & \lesssim\delta\frac{\sigma_{r}^{\star2}}{ndp^{2}}\omega_{j}^{\star2}\left\Vert \bm{U}_{i,\cdot}^{\star}\right\Vert _{2}^{2}\lesssim\delta v_{i,j}^{\star},
\end{align*}
provided that $\zeta_{\mathsf{1st}}/\sigma_{r}^{\star2}\lesssim\delta/(\kappa^{2}\mu r^{2}\sqrt{\log(n+d)})$
and $\theta\lesssim\delta/(\kappa^{3/2}\mu r^{2})$. Here (i) invokes
the AM-GM inequality. 
\item Regarding $\gamma_{4,3,4}$, we have
\begin{align*}
\gamma_{4,3,4} & \lesssim\underbrace{\frac{1}{np^{2}}\left\Vert \bm{U}_{i,\cdot}^{\star}\bm{\Sigma}^{\star}\right\Vert _{2}^{2}\omega_{\max}^{2}\left\Vert \bm{U}_{j,\cdot}^{\star}\right\Vert _{2}^{2}\frac{\zeta_{\mathsf{1st}}}{\sigma_{r}^{\star2}}\sqrt{r^{2}\log\left(n+d\right)}}_{\eqqcolon\gamma_{4,3,4,1}}+\underbrace{\frac{1}{np^{2}}\left\Vert \bm{U}_{i,\cdot}^{\star}\bm{\Sigma}^{\star}\right\Vert _{2}^{2}\omega_{\max}^{2}\frac{\theta r}{\sigma_{1}^{\star}}\omega_{j}^{\star}\left\Vert \bm{U}_{j,\cdot}^{\star}\right\Vert _{2}}_{\eqqcolon\gamma_{4,3,4,2}}\\
 & \quad+\underbrace{\frac{1}{np^{2}}\left\Vert \bm{U}_{i,\cdot}^{\star}\bm{\Sigma}^{\star}\right\Vert _{2}^{2}\omega_{\max}^{2}\frac{\theta^{2}r}{\kappa\sigma_{r}^{\star2}}\omega_{j}^{\star2}}_{\eqqcolon\gamma_{4,3,4,3}}\lesssim\delta v_{i,j}^{\star},
\end{align*}
where the last relation holds since 
\begin{align*}
\gamma_{4,3,2,1}+\gamma_{4,3,2,2} & \lesssim\delta\frac{\sigma_{r}^{\star2}}{ndp^{2}}\left(\omega_{j}^{\star2}\left\Vert \bm{U}_{i,\cdot}^{\star}\right\Vert _{2}^{2}+\omega_{i}^{\star2}\left\Vert \bm{U}_{j,\cdot}^{\star}\right\Vert _{2}^{2}\right)\lesssim\delta v_{i,j}^{\star},\\
\gamma_{4,3,2,3} & \lesssim\delta\frac{\omega_{\min}^{2}}{np^{2}}\left(\omega_{j}^{\star2}\left\Vert \bm{U}_{i,\cdot}^{\star}\right\Vert _{2}^{2}+\omega_{i}^{\star2}\left\Vert \bm{U}_{j,\cdot}^{\star}\right\Vert _{2}^{2}\right)\lesssim\delta v_{i,j}^{\star},
\end{align*}
provided that $\zeta_{\mathsf{1st}}/\sigma_{r}^{\star2}\lesssim\delta/\sqrt{\kappa^{2}\mu^{2}r^{4}\kappa_{\omega}^{2}\log(n+d)}$,
$\theta\lesssim\sqrt{\delta/(\kappa_{\omega}r)}$, and
\[
\theta\frac{\omega_{\max}}{\sigma_{r}^{\star}}\sqrt{d}\lesssim\frac{\delta}{\sqrt{\kappa_{\omega}\kappa\mu r^{3}}}.
\]
Here, (i) utilizes (\ref{eq:ce-var-est-inter-1}) and (ii) follows
from (\ref{eq:ce-var-est-inter-3}). 
\end{itemize}
Taking the bounds on $\gamma_{4,3,1}$, $\gamma_{4,3,2}$ and $\gamma_{4,3,3}$
collectively yields 
\[
\gamma_{4,3}\lesssim\gamma_{4,3,1}+\gamma_{4,3,2}+\gamma_{4,3,3}\lesssim\delta v_{i,j}^{\star},
\]
provided that $\zeta_{\mathsf{1st}}/\sigma_{r}^{\star2}\lesssim\delta/(\kappa^{2}\mu r^{2}\kappa_{\omega}\sqrt{\log(n+d)})$,
$\theta\lesssim\delta/(\kappa^{3/2}\mu r^{2}\kappa_{\omega}^{1/2})$,
$\varepsilon\lesssim1$, 
\[
\theta\frac{\omega_{\max}}{\sigma_{r}^{\star}}\sqrt{d}\lesssim\frac{\delta}{\sqrt{\kappa_{\omega}\kappa\mu r^{3}}}
\]
and 
\begin{align*}
\left\Vert \bm{U}_{i,\cdot}^{\star}\right\Vert _{2}+\left\Vert \bm{U}_{j,\cdot}^{\star}\right\Vert _{2} & \gtrsim\delta^{-1}\sqrt{\kappa}\mu r^{2}\kappa_{\omega}\frac{\omega_{\max}}{\sigma_{r}^{\star}}\theta+\delta^{-1}r\kappa_{\omega}^{1/2}\frac{\omega_{\max}^{2}}{\sigma_{1}^{\star2}}\sqrt{\frac{1}{np^{2}}}+\delta^{-1}\frac{\mu r^{2}\kappa_{\omega}^{1/2}}{\sqrt{ndp^{2}}}\cdot\sqrt{\frac{1}{d}}.
\end{align*}
In view of (\ref{eq:theta-omega-max}), we know that the last two
conditions above can be guaranteed by $ndp^{2}\gtrsim\delta^{-2}\kappa^{4}\mu^{3}r^{6}\kappa_{\omega}\log^{4}(n+d)$,
\[
\frac{\omega_{\max}}{\sigma_{r}^{\star}}\sqrt{\frac{d}{np}}\lesssim\frac{\delta}{\kappa\mu^{1/2}r^{2}\kappa_{\omega}^{1/2}\log\left(n+d\right)},\qquad\frac{\omega_{\max}^{2}}{p\sigma_{r}^{\star2}}\sqrt{\frac{d}{n}}\lesssim\frac{\delta}{\kappa\mu^{1/2}r^{2}\kappa_{\omega}^{1/2}\log\left(n+d\right)},
\]
and
\begin{align*}
\left\Vert \bm{U}_{i,\cdot}^{\star}\right\Vert _{2}+\left\Vert \bm{U}_{j,\cdot}^{\star}\right\Vert _{2} & \gtrsim\delta^{-1}\kappa\mu r^{2}\kappa_{\omega}\log\left(n+d\right)\left[\frac{\omega_{\max}}{\sigma_{r}^{\star}}\sqrt{\frac{d}{np}}+\frac{\omega_{\max}^{2}}{p\sigma_{r}^{\star2}}\sqrt{\frac{d}{n}}+\frac{\kappa\mu r\log\left(n+d\right)}{\sqrt{ndp^{2}}}\right]\cdot\sqrt{\frac{r}{d}}.
\end{align*}

\subparagraph{Step 4.4: putting the bounds on $\gamma_{4,1}$, $\gamma_{4,2}$
and $\gamma_{4,3}$ together.}

In view of (\ref{eq:theta-zeta-1st}), we can take the bounds on $\gamma_{4,1}$,
$\gamma_{4,2}$ and $\gamma_{4,3}$ together to reach 
\[
\left|\alpha_{4}-\beta_{4}\right|\lesssim\gamma_{4,1}+\gamma_{4,2}+\gamma_{4,3}\lesssim\delta v_{i,j}^{\star},
\]
provided that $\varepsilon\lesssim1$, $n\gtrsim\delta^{-2}\kappa^{7}\mu^{2}r^{3}\kappa_{\omega}^{2}\log(n+d)$,
$ndp^{2}\gtrsim\delta^{-2}\kappa^{4}\mu^{3}r^{6}\kappa_{\omega}^{2}\log^{4}(n+d)$,
$np\gtrsim\delta^{-2}\kappa^{2}\mu^{2}r^{2}\kappa_{\omega}^{2}\log^{2}(n+d)$,
\begin{equation}
\frac{\zeta_{\mathsf{1st}}}{\sigma_{r}^{\star2}}\lesssim\frac{\delta}{\kappa^{2}\mu r^{2}\kappa_{\omega}\sqrt{\log\left(n+d\right)}},\qquad\zeta_{\mathsf{2nd},i}\sqrt{d}\lesssim\frac{\delta}{\sqrt{\kappa^{2}\mu r\kappa_{\omega}^{2}}},\label{eq:ce-var-est-inter-10}
\end{equation}
\[
\frac{\omega_{\max}}{\sigma_{r}^{\star}}\sqrt{\frac{d}{np}}\lesssim\frac{\delta}{\sqrt{\kappa^{2}\mu r^{2}\kappa_{\omega}\log^{2}\left(n+d\right)}},\qquad\frac{\omega_{\max}^{2}}{p\sigma_{r}^{\star2}}\sqrt{\frac{d}{n}}\lesssim\frac{\delta}{\sqrt{\kappa^{2}\mu r^{2}\kappa_{\omega}\log^{2}\left(n+d\right)}},
\]
and 
\[
\left\Vert \bm{U}_{i,\cdot}^{\star}\right\Vert _{2}+\left\Vert \bm{U}_{j,\cdot}^{\star}\right\Vert _{2}\gtrsim\delta^{-1}\kappa\mu r^{2}\kappa_{\omega}\log\left(n+d\right)\left[\frac{\kappa\mu r\log\left(n+d\right)}{\sqrt{nd}p}+\frac{\omega_{\max}^{2}}{p\sigma_{r}^{\star2}}\sqrt{\frac{d}{n}}+\frac{\omega_{\max}}{\sigma_{r}^{\star}}\sqrt{\frac{d}{np}}\right]\sqrt{\frac{r}{d}}.
\]
Let us now take a closer look at (\ref{eq:ce-var-est-inter-10}).
In view of (\ref{eq:pca-1st-err-useful}) and the definition of $\zeta_{\mathsf{2nd},i}$,
we know that (\ref{eq:ce-var-est-inter-10}) is equivalent to 
\begin{align*}
ndp^{2}\gtrsim\delta^{-2}\kappa^{6}\mu^{4}r^{6}\kappa_{\omega}^{2}\log^{5}\left(n+d\right), & \qquad np\gtrsim\delta^{-2}\kappa^{6}\mu^{3}r^{5}\kappa_{\omega}^{2}\log^{3}\left(n+d\right),\\
\frac{\omega_{\max}^{2}}{p\sigma_{r}^{\star2}}\sqrt{\frac{d}{n}}\lesssim\frac{\delta}{\kappa^{2}\mu r^{2}\kappa_{\omega}\log^{3/2}\left(n+d\right)}, & \qquad\frac{\omega_{\max}}{\sigma_{r}^{\star}}\sqrt{\frac{d}{np}}\lesssim\frac{\delta}{\kappa^{5/2}\mu r^{2}\kappa_{\omega}\log\left(n+d\right)},
\end{align*}
provided that $d\gtrsim\kappa^{2}\mu\log(n+d)$. This concludes our
bound on $\left|\alpha_{4}-\beta_{4}\right|$ and the required conditions.

Similarly, we can also prove that (which we omit here for the sake
of brevity) 
\[
\left|\alpha_{6}-\beta_{6}\right|\lesssim\delta v_{i,j}^{\star}
\]
under the above conditions.

\paragraph{Step 5: bounding $\vert\alpha_{5}-\beta_{5}\vert$ and $\vert\alpha_{7}-\beta_{7}\vert$.}

Regarding $\vert\alpha_{5}-\beta_{5}\vert$, we first make the observation
that 
\begin{align*}
\left|\alpha_{5}-\beta_{5}\right| & \lesssim\frac{1}{np^{2}}\left|\sum_{k=1}^{d}S_{i,k}^{\star2}\left(\bm{U}_{k,\cdot}^{\star}\bm{U}_{j,\cdot}^{\star\top}\right)^{2}-\sum_{k=1}^{d}S_{i,k}^{2}\left(\bm{U}_{k,\cdot}\bm{U}_{j,\cdot}^{\top}\right)^{2}\right|\\
 & \lesssim\frac{1}{np^{2}}\left|\sum_{k=1}^{d}\left(S_{i,k}^{2}-S_{i,k}^{\star2}\right)\left(\bm{U}_{k,\cdot}\bm{U}_{j,\cdot}^{\top}\right)^{2}\right|+\frac{1}{np^{2}}\left|\sum_{k=1}^{d}S_{i,k}^{\star2}\left[\left(\bm{U}_{k,\cdot}^{\star}\bm{U}_{j,\cdot}^{\star\top}\right)^{2}-\left(\bm{U}_{k,\cdot}\bm{U}_{j,\cdot}^{\top}\right)^{2}\right]\right|\\
 & \lesssim\frac{1}{np^{2}}\left\Vert \bm{U}_{j,\cdot}\right\Vert _{2}^{2}\max_{k\in[d]}\left|S_{i,k}^{2}-S_{i,k}^{\star2}\right|+\frac{1}{np^{2}}\left\Vert \bm{U}_{i,\cdot}\bm{\Sigma}^{\star}\right\Vert _{2}^{2}\left\Vert \bm{U}^{\star}\bm{\Sigma}^{\star}\right\Vert _{2,\infty}^{2}\sum_{k=1}^{d}\left|\left(\bm{U}_{k,\cdot}^{\star}\bm{U}_{j,\cdot}^{\star\top}\right)^{2}-\left(\bm{U}_{k,\cdot}\bm{U}_{j,\cdot}^{\top}\right)^{2}\right|\\
 & \overset{\text{(i)}}{\lesssim}\frac{1}{np^{2}}\left(\left\Vert \bm{U}_{j,\cdot}^{\star}\right\Vert _{2}^{2}+\frac{\theta^{2}}{\kappa\sigma_{r}^{\star2}}\omega_{j}^{\star2}\right)\max_{k\in[d]}\left|S_{i,k}^{2}-S_{i,k}^{\star2}\right|+\frac{\mu r}{ndp^{2}}\sigma_{1}^{\star2}\left\Vert \bm{U}_{i,\cdot}\bm{\Sigma}^{\star}\right\Vert _{2}^{2}\sum_{k=1}^{d}\left|\left(\bm{U}_{k,\cdot}^{\star}\bm{U}_{j,\cdot}^{\star\top}\right)^{2}-\left(\bm{U}_{k,\cdot}\bm{U}_{j,\cdot}^{\top}\right)^{2}\right|\\
 & \overset{\text{(ii)}}{\lesssim}\frac{1}{np^{2}}\left\Vert \bm{U}_{j,\cdot}^{\star}\right\Vert _{2}^{2}\max_{k\in[d]}\left|S_{i,k}^{2}-S_{i,k}^{\star2}\right|+\frac{1}{np^{2}}\frac{\theta^{2}}{\kappa\sigma_{r}^{\star2}}\omega_{j}^{\star2}\max_{k\in[d]}\left|S_{i,k}^{2}-S_{i,k}^{\star2}\right|+\delta v_{i,j}^{\star}\\
 & \overset{\text{(iii)}}{\lesssim}\underbrace{\frac{1}{np^{2}}\sqrt{\frac{\mu r}{d}}\sigma_{1}^{\star}\left\Vert \bm{U}_{i,\cdot}^{\star}\bm{\Sigma}^{\star}\right\Vert _{2}\left\Vert \bm{U}_{j,\cdot}^{\star}\right\Vert _{2}^{2}\max_{k\in[d]}\left|S_{i,k}-S_{i,k}^{\star}\right|}_{\eqqcolon\gamma_{5,1}}+\underbrace{\frac{1}{np^{2}}\sqrt{\frac{\mu r}{d}}\left\Vert \bm{U}_{i,\cdot}^{\star}\bm{\Sigma}^{\star}\right\Vert _{2}\frac{\theta^{2}}{\sqrt{\kappa}\sigma_{r}^{\star}}\omega_{j}^{\star2}\max_{k\in[d]}\left|S_{i,k}-S_{i,k}^{\star}\right|}_{\eqqcolon\gamma_{5,2}}\\
 & \quad+\underbrace{\frac{1}{np^{2}}\left\Vert \bm{U}_{j,\cdot}^{\star}\right\Vert _{2}^{2}\max_{k\in[d]}\left|S_{i,k}-S_{i,k}^{\star}\right|^{2}}_{\eqqcolon\gamma_{5,3}}+\underbrace{\frac{1}{np^{2}}\frac{\theta^{2}}{\kappa\sigma_{r}^{\star2}}\omega_{j}^{\star2}\max_{k\in[d]}\left|S_{i,k}-S_{i,k}^{\star}\right|^{2}}_{\eqqcolon\gamma_{5,4}}+\delta v_{i,j}^{\star}.
\end{align*}
Here, (i) follows from (\ref{eq:ce-var-est-inter-8}); (ii) holds
since 
\[
\frac{\mu r}{ndp^{2}}\sigma_{1}^{\star2}\left\Vert \bm{U}_{i,\cdot}\bm{\Sigma}^{\star}\right\Vert _{2}^{2}\sum_{k=1}^{d}\left|\left(\bm{U}_{k,\cdot}^{\star}\bm{U}_{j,\cdot}^{\star\top}\right)^{2}-\left(\bm{U}_{k,\cdot}\bm{U}_{j,\cdot}^{\top}\right)^{2}\right|=\gamma_{4,3,3}\lesssim\delta v_{i,j}^{\star};
\]
and (iii) follows from the fact that 
\begin{align*}
\max_{k\in[d]}\left|S_{i,k}^{2}-S_{i,k}^{\star2}\right| & \lesssim\max_{k\in[d]}\left|S_{i,k}-S_{i,k}^{\star}\right|S_{i,k}^{\star}+\max_{k\in[d]}\left|S_{i,k}-S_{i,k}^{\star}\right|^{2}\\
 & \lesssim\left\Vert \bm{U}_{i,\cdot}^{\star}\bm{\Sigma}^{\star}\right\Vert _{2}\left\Vert \bm{U}^{\star}\bm{\Sigma}^{\star}\right\Vert _{2,\infty}\max_{k\in[d]}\left|S_{i,k}-S_{i,k}^{\star}\right|+\max_{k\in[d]}\left|S_{i,k}-S_{i,k}^{\star}\right|^{2}\\
 & \lesssim\sqrt{\frac{\mu r}{d}}\sigma_{1}^{\star}\left\Vert \bm{U}_{i,\cdot}^{\star}\bm{\Sigma}^{\star}\right\Vert _{2}\max_{k\in[d]}\left|S_{i,k}-S_{i,k}^{\star}\right|+\max_{k\in[d]}\left|S_{i,k}-S_{i,k}^{\star}\right|^{2}.
\end{align*}
In view of Lemma \ref{lemma:pca-noise-level-est}, we know that for
each $k\in[d]$, 
\begin{align*}
\left|S_{i,k}-S_{i,k}^{\star}\right| & \lesssim\left(\theta+\sqrt{\frac{\kappa^{3}r\log\left(n+d\right)}{n}}\right)\left\Vert \bm{U}_{i,\cdot}^{\star}\bm{\Sigma}^{\star}\right\Vert _{2}\sqrt{\frac{\mu r}{d}}\sigma_{1}^{\star}+\theta^{2}\omega_{i}^{\star}\omega_{\max}+\theta\omega_{i}^{\star}\sqrt{\frac{\mu r}{d}}\sigma_{1}^{\star}\\
 & \quad+\theta\omega_{\max}\left\Vert \bm{U}_{i,\cdot}^{\star}\bm{\Sigma}^{\star}\right\Vert _{2}+\zeta_{\mathsf{2nd},i}\sqrt{\frac{\mu r}{d}}\sigma_{1}^{\star2}+\zeta_{\mathsf{2nd},k}\left\Vert \bm{U}_{i,\cdot}^{\star}\bm{\Sigma}^{\star}\right\Vert _{2}\sigma_{1}^{\star}+\zeta_{\mathsf{2nd},i}\zeta_{\mathsf{2nd},k}\sigma_{1}^{\star2}\\
 & \lesssim\left(\theta+\sqrt{\frac{\kappa^{3}r\log\left(n+d\right)}{n}}\right)\left\Vert \bm{U}_{i,\cdot}^{\star}\bm{\Sigma}^{\star}\right\Vert _{2}\sqrt{\frac{\mu r}{d}}\sigma_{1}^{\star}+\theta\omega_{i}^{\star}\sqrt{\frac{\mu r}{d}}\sigma_{1}^{\star}\\
 & \quad+\zeta_{\mathsf{2nd},i}\sqrt{\frac{\mu r}{d}}\sigma_{1}^{\star2}+\theta\omega_{\max}\left\Vert \bm{U}_{i,\cdot}^{\star}\bm{\Sigma}^{\star}\right\Vert _{2},
\end{align*}
where the last relation holds due to (\ref{eq:theta-omega-sqrt-d})
and (\ref{eq:ce-var-est-inter-2}), provided that $\varepsilon\lesssim1$,
$np\gtrsim1$ and $ndp^{2}\gtrsim1$. With these preparations in place,
we now proceed to bound the terms $\gamma_{5,1}$, $\gamma_{5,2}$,
$\gamma_{5,3}$ and $\gamma_{5,4}$ separately. 
\begin{itemize}
\item Regarding $\gamma_{5,1}$, we have 
\begin{align*}
\gamma_{5,1} & \lesssim\underbrace{\frac{\mu r}{ndp^{2}}\sigma_{1}^{\star2}\left\Vert \bm{U}_{i,\cdot}^{\star}\bm{\Sigma}^{\star}\right\Vert _{2}^{2}\left\Vert \bm{U}_{j,\cdot}^{\star}\right\Vert _{2}^{2}\left(\theta+\sqrt{\frac{\kappa^{3}r\log\left(n+d\right)}{n}}\right)}_{\eqqcolon\gamma_{5,1,1}}\\
 & \quad+\underbrace{\frac{\mu r}{ndp^{2}}\sigma_{1}^{\star2}\left\Vert \bm{U}_{i,\cdot}^{\star}\bm{\Sigma}^{\star}\right\Vert _{2}\left\Vert \bm{U}_{j,\cdot}^{\star}\right\Vert _{2}^{2}\theta\omega_{i}^{\star}}_{\eqqcolon\gamma_{5,1,2}}+\underbrace{\frac{\mu r}{ndp^{2}}\sigma_{1}^{\star2}\left\Vert \bm{U}_{i,\cdot}^{\star}\bm{\Sigma}^{\star}\right\Vert _{2}\left\Vert \bm{U}_{j,\cdot}^{\star}\right\Vert _{2}^{2}\zeta_{\mathsf{2nd},i}\sigma_{1}^{\star}}_{\eqqcolon\gamma_{5,1,3}}\\
 & \quad+\underbrace{\frac{1}{np^{2}}\sqrt{\frac{\mu r}{d}}\sigma_{1}^{\star}\left\Vert \bm{U}_{i,\cdot}^{\star}\bm{\Sigma}^{\star}\right\Vert _{2}^{2}\left\Vert \bm{U}_{j,\cdot}^{\star}\right\Vert _{2}^{2}\theta\omega_{\max}}_{\eqqcolon\gamma_{5,1,4}}\lesssim\delta v_{i,j}^{\star},
\end{align*}
where the last relation holds since 
\begin{align*}
\gamma_{5,1,1} & \lesssim\frac{\delta}{ndp^{2}\kappa}\left\Vert \bm{U}_{i,\cdot}^{\star}\bm{\Sigma}^{\star}\right\Vert _{2}^{2}\left\Vert \bm{U}_{j,\cdot}^{\star}\bm{\Sigma}^{\star}\right\Vert _{2}^{2}\lesssim\delta v_{i,j}^{\star},\\
\gamma_{5,1,2} & \overset{\text{(i)}}{\lesssim}\frac{\delta}{ndp^{2}\kappa}\left\Vert \bm{U}_{i,\cdot}^{\star}\bm{\Sigma}^{\star}\right\Vert _{2}^{2}\left\Vert \bm{U}_{j,\cdot}^{\star}\bm{\Sigma}^{\star}\right\Vert _{2}^{2}+\delta\frac{\sigma_{r}^{\star2}}{ndp^{2}}\omega_{i}^{\star2}\left\Vert \bm{U}_{j,\cdot}^{\star}\right\Vert _{2}^{2}\lesssim\delta v_{i,j}^{\star},\\
\gamma_{5,1,3} & \lesssim\frac{\kappa\mu r}{\sqrt{ndp^{2}}}\cdot\sqrt{\frac{1}{ndp^{2}\kappa}}\left\Vert \bm{U}_{i,\cdot}^{\star}\bm{\Sigma}^{\star}\right\Vert _{2}\left\Vert \bm{U}_{j,\cdot}^{\star}\bm{\Sigma}^{\star}\right\Vert _{2}\cdot\sigma_{1}^{\star2}\left(\left\Vert \bm{U}_{i,\cdot}^{\star}\right\Vert _{2}\zeta_{\mathsf{2nd},j}+\left\Vert \bm{U}_{j,\cdot}^{\star}\right\Vert _{2}\zeta_{\mathsf{2nd},i}\right)\\
 & \overset{\text{(ii)}}{\lesssim}\frac{\kappa\mu r}{\sqrt{ndp^{2}}}\cdot\sqrt{\frac{1}{ndp^{2}\kappa}}\left\Vert \bm{U}_{i,\cdot}^{\star}\bm{\Sigma}^{\star}\right\Vert _{2}\left\Vert \bm{U}_{j,\cdot}^{\star}\bm{\Sigma}^{\star}\right\Vert _{2}\cdot\varepsilon(v_{i,j}^{\star})^{1/2}\lesssim\frac{\kappa\mu r}{\sqrt{ndp^{2}}}\varepsilon v_{i,j}^{\star}\lesssim\delta v_{i,j}^{\star},\\
\gamma_{5,1,4} & \lesssim\delta\frac{1}{ndp^{2}\kappa}\left\Vert \bm{U}_{i,\cdot}^{\star}\bm{\Sigma}^{\star}\right\Vert _{2}^{2}\left\Vert \bm{U}_{j,\cdot}^{\star}\bm{\Sigma}^{\star}\right\Vert _{2}^{2}\lesssim\delta v_{i,j}^{\star},
\end{align*}
provided that $\theta\lesssim\delta/(\kappa^{2}\mu r)$, $n\gtrsim\delta^{-2}\kappa^{7}\mu^{2}r^{3}\log(n+d)$,
$\varepsilon\lesssim1$, $ndp^{2}\gtrsim\delta^{-2}\kappa^{2}\mu^{2}r^{2}$,
and 
\[
\theta\frac{\omega_{\max}}{\sigma_{r}^{\star}}\sqrt{d}\lesssim\frac{\delta}{\sqrt{\kappa^{3}\mu r}}.
\]
 In the above relations, (i) uses the AM-GM inequality, where (ii)
utilizes (\ref{eq:ce-var-est-inter-1}). 
\item Regarding $\gamma_{5,2}$, we can derive 
\begin{align*}
\gamma_{5,2} & \lesssim\underbrace{\frac{\mu r}{ndp^{2}}\sigma_{1}^{\star}\left\Vert \bm{U}_{i,\cdot}^{\star}\bm{\Sigma}^{\star}\right\Vert _{2}^{2}\frac{\theta^{2}}{\sqrt{\kappa}\sigma_{r}^{\star}}\omega_{j}^{\star2}\left(\theta+\sqrt{\frac{\kappa^{3}r\log\left(n+d\right)}{n}}\right)}_{\eqqcolon\gamma_{5,2,1}}\\
 & \quad+\underbrace{\frac{\mu r}{ndp^{2}}\left\Vert \bm{U}_{i,\cdot}^{\star}\bm{\Sigma}^{\star}\right\Vert _{2}\theta^{3}\omega_{j}^{\star2}\omega_{i}^{\star}}_{\eqqcolon\gamma_{5,2,2}}+\underbrace{\frac{\mu r}{ndp^{2}}\sigma_{1}^{\star}\left\Vert \bm{U}_{i,\cdot}^{\star}\bm{\Sigma}^{\star}\right\Vert _{2}\theta^{2}\omega_{j}^{\star2}\zeta_{\mathsf{2nd},i}}_{\eqqcolon\gamma_{5,2,3}}\\
 & \quad+\underbrace{\frac{1}{np^{2}}\sqrt{\frac{\mu r}{d}}\left\Vert \bm{U}_{i,\cdot}^{\star}\bm{\Sigma}^{\star}\right\Vert _{2}^{2}\frac{\theta^{3}}{\sqrt{\kappa}\sigma_{r}^{\star}}\omega_{j}^{\star2}\omega_{\max}}_{\eqqcolon\gamma_{5,2,4}}\lesssim\delta v_{i,j}^{\star},
\end{align*}
\[
\frac{\mu r}{ndp^{2}}\sigma_{1}^{\star}\left\Vert \bm{U}_{i,\cdot}^{\star}\bm{\Sigma}^{\star}\right\Vert _{2}\theta^{2}\omega_{j}^{\star2}\zeta_{\mathsf{2nd},i}\lesssim\theta\frac{\mu r}{ndp^{2}}\left\Vert \bm{U}_{i,\cdot}^{\star}\bm{\Sigma}^{\star}\right\Vert _{2}\omega_{j}^{\star}\cdot\theta\sigma_{1}^{\star}\left(\omega_{i}^{\star}\zeta_{\mathsf{2nd},j}+\omega_{j}^{\star}\zeta_{\mathsf{2nd},i}\right)
\]
where the last relation holds since 
\begin{align*}
\gamma_{5,2,1} & \lesssim\delta\frac{\sigma_{r}^{\star2}}{ndp^{2}}\left(\omega_{j}^{\star2}\left\Vert \bm{U}_{i,\cdot}^{\star}\right\Vert _{2}^{2}+\omega_{i}^{\star2}\left\Vert \bm{U}_{j,\cdot}^{\star}\right\Vert _{2}^{2}\right)\lesssim\delta v_{i,j}^{\star},\\
\gamma_{5,2,2} & \lesssim\frac{\mu r}{ndp^{2}}\left\Vert \bm{U}_{i,\cdot}^{\star}\bm{\Sigma}^{\star}\right\Vert _{2}\theta^{3}\omega_{j}^{\star2}\omega_{i}^{\star}\lesssim\theta^{2}\gamma_{4,3,2,2}\lesssim\delta v_{i,j}^{\star},\\
\gamma_{5,2,3} & \lesssim\theta\frac{\sqrt{\kappa}\mu r}{\sqrt{ndp^{2}}}\cdot\frac{\sigma_{r}^{\star}}{\sqrt{ndp^{2}}}\omega_{j}^{\star}\left\Vert \bm{U}_{i,\cdot}^{\star}\right\Vert _{2}\cdot\theta\sigma_{1}^{\star}\omega_{j}^{\star}\zeta_{\mathsf{2nd},i}\overset{\text{(i)}}{\lesssim}\theta\frac{\sqrt{\kappa}\mu r}{\sqrt{ndp^{2}}}\varepsilon v_{i,j}^{\star}\lesssim\delta v_{i,j}^{\star},\\
\gamma_{5,2,4} & \lesssim\delta\frac{\sigma_{r}^{\star2}}{ndp^{2}}\left(\omega_{j}^{\star2}\left\Vert \bm{U}_{i,\cdot}^{\star}\right\Vert _{2}^{2}+\omega_{i}^{\star2}\left\Vert \bm{U}_{j,\cdot}^{\star}\right\Vert _{2}^{2}\right)\lesssim\delta v_{i,j}^{\star},
\end{align*}
provided that $\theta\lesssim\delta/(\kappa^{1/2}\mu r)$, $n\gtrsim\kappa^{3}r\log(n+d)$,
$\varepsilon\lesssim1$, $ndp^{2}\gtrsim1$ and 
\[
\theta^{3}\frac{\omega_{\max}}{\sigma_{r}^{\star}}\sqrt{d}\lesssim\frac{\delta}{\sqrt{\kappa\mu r}}.
\]
Here, the relation (i) in the above inequality arises from (\ref{eq:ce-var-est-inter-3}). 
\item With regards to $\gamma_{5,3}$, we make the observation that 
\begin{align*}
\gamma_{5,3} & \lesssim\underbrace{\frac{\mu r}{ndp^{2}}\sigma_{1}^{\star2}\left\Vert \bm{U}_{i,\cdot}^{\star}\bm{\Sigma}^{\star}\right\Vert _{2}^{2}\left\Vert \bm{U}_{j,\cdot}^{\star}\right\Vert _{2}^{2}\left(\theta^{2}+\frac{\kappa^{3}r\log\left(n+d\right)}{n}\right)}_{\eqqcolon\gamma_{5,3,1}}+\underbrace{\frac{\mu r}{ndp^{2}}\sigma_{1}^{\star2}\left\Vert \bm{U}_{j,\cdot}^{\star}\right\Vert _{2}^{2}\theta^{2}\omega_{i}^{\star2}}_{\eqqcolon\gamma_{5,3,2}}\\
 & \quad+\underbrace{\frac{\mu r}{ndp^{2}}\left\Vert \bm{U}_{j,\cdot}^{\star}\right\Vert _{2}^{2}\zeta_{\mathsf{2nd},i}^{2}\sigma_{1}^{\star4}}_{\eqqcolon\gamma_{5,3,3}}+\underbrace{\frac{1}{np^{2}}\left\Vert \bm{U}_{j,\cdot}^{\star}\right\Vert _{2}^{2}\theta^{2}\omega_{\max}^{2}\left\Vert \bm{U}_{i,\cdot}^{\star}\bm{\Sigma}^{\star}\right\Vert _{2}^{2}}_{\eqqcolon\gamma_{5,3,4}}\lesssim\delta v_{i,j}^{\star},
\end{align*}
where the last relation holds since 
\begin{align*}
\gamma_{5,3,1} & \lesssim\left(\theta+\sqrt{\frac{\kappa^{3}r\log\left(n+d\right)}{n}}\right)\gamma_{5,1,1}\lesssim\delta v_{i,j}^{\star},\\
\gamma_{5,3,2} & \lesssim\delta\frac{\sigma_{r}^{\star2}}{ndp^{2}}\left(\omega_{j}^{\star2}\left\Vert \bm{U}_{i,\cdot}^{\star}\right\Vert _{2}^{2}+\omega_{i}^{\star2}\left\Vert \bm{U}_{j,\cdot}^{\star}\right\Vert _{2}^{2}\right)\lesssim\delta v_{i,j}^{\star},\\
\gamma_{5,3,3} & \lesssim\frac{\mu r}{ndp^{2}}\left[\sigma_{1}^{\star2}\left(\left\Vert \bm{U}_{i,\cdot}^{\star}\right\Vert _{2}\zeta_{\mathsf{2nd},j}+\left\Vert \bm{U}_{j,\cdot}^{\star}\right\Vert _{2}\zeta_{\mathsf{2nd},i}\right)\right]^{2}\overset{\text{(i)}}{\lesssim}\frac{\mu r}{ndp^{2}}\varepsilon^{2}v_{i,j}^{\star}\lesssim\delta v_{i,j}^{\star},\\
\gamma_{5,3,4} & \lesssim\delta\frac{1}{ndp^{2}\kappa}\left\Vert \bm{U}_{i,\cdot}^{\star}\bm{\Sigma}^{\star}\right\Vert _{2}^{2}\left\Vert \bm{U}_{j,\cdot}^{\star}\bm{\Sigma}^{\star}\right\Vert _{2}^{2}\lesssim\delta v_{i,j}^{\star},
\end{align*}
provided that $\theta\lesssim\sqrt{\delta/(\kappa\mu r)}$, $n\gtrsim\kappa^{3}r\log(n+d)$,
$\varepsilon\lesssim1$, $ndp^{2}\gtrsim\delta^{-1}\mu r$ and
\[
\theta\frac{\omega_{\max}}{\sigma_{r}^{\star}}\sqrt{d}\lesssim\frac{\delta}{\sqrt{\kappa}}.
\]
Note that the inequality (i) in the above relation results from (\ref{eq:ce-var-est-inter-1}). 
\item When it comes to $\gamma_{5,4}$, we have the following upper bound
\begin{align*}
\gamma_{5,4} & \lesssim\underbrace{\frac{\mu r}{ndp^{2}}\left\Vert \bm{U}_{i,\cdot}^{\star}\bm{\Sigma}^{\star}\right\Vert _{2}^{2}\theta^{2}\omega_{j}^{\star2}\left(\theta^{2}+\frac{\kappa^{3}r\log\left(n+d\right)}{n}\right)}_{\eqqcolon\gamma_{5,4,1}}+\underbrace{\frac{\mu r}{ndp^{2}}\theta^{4}\omega_{j}^{\star2}\omega_{\max}^{2}}_{\eqqcolon\gamma_{5,4,2}}\\
 & \quad+\underbrace{\frac{\mu r}{ndp^{2}}\frac{\theta^{2}}{\kappa\sigma_{r}^{\star2}}\omega_{j}^{\star2}\zeta_{\mathsf{2nd},i}^{2}\sigma_{1}^{\star4}}_{\eqqcolon\gamma_{5,4,3}}+\underbrace{\frac{1}{np^{2}}\frac{\theta^{4}}{\kappa\sigma_{r}^{\star2}}\omega_{j}^{\star2}\omega_{\max}^{2}\left\Vert \bm{U}_{i,\cdot}^{\star}\bm{\Sigma}^{\star}\right\Vert _{2}^{2}}_{\eqqcolon\gamma_{5,4,4}}\lesssim\delta v_{i,j}^{\star},
\end{align*}
where the last inequality holds true since 
\begin{align*}
\gamma_{5,4,1} & \lesssim\delta\frac{\sigma_{r}^{\star2}}{ndp^{2}}\left(\omega_{j}^{\star2}\left\Vert \bm{U}_{i,\cdot}^{\star}\right\Vert _{2}^{2}+\omega_{i}^{\star2}\left\Vert \bm{U}_{j,\cdot}^{\star}\right\Vert _{2}^{2}\right)\lesssim\delta v_{i,j}^{\star},\\
\gamma_{5,4,2} & \overset{\text{(i)}}{\lesssim}\frac{\mu r\kappa_{\omega}}{ndp^{2}}\varepsilon^{2}v_{i,j}^{\star}\lesssim\delta v_{i,j}^{\star},\\
\gamma_{5,4,3} & \overset{\text{(ii)}}{\lesssim}\frac{\mu r\kappa_{\omega}^{1/2}}{ndp^{2}\kappa}\varepsilon^{2}v_{i,j}^{\star}\lesssim\delta v_{i,j}^{\star},\\
\gamma_{5,4,4} & \lesssim\delta\frac{\omega_{\min}^{2}}{np^{2}}\left(\omega_{j}^{\star2}\left\Vert \bm{U}_{i,\cdot}^{\star}\right\Vert _{2}^{2}+\omega_{i}^{\star2}\left\Vert \bm{U}_{j,\cdot}^{\star}\right\Vert _{2}^{2}\right)
\end{align*}
provided that $\theta\lesssim\delta/(\kappa\mu r\kappa_{\omega})$,
$n\gtrsim\kappa^{3}r\log(n+d)$ and $ndp^{2}\gtrsim\delta^{-1}\mu r\kappa_{\omega}$.
Here, the inequalities (i) and (ii) follow from (\ref{eq:ce-var-est-inter-3}). 
\end{itemize}
Combining the above bounds on $\gamma_{5,1}$, $\gamma_{5,2}$, $\gamma_{5,3}$
and $\gamma_{5,4}$ yields 
\[
\left|\alpha_{5}-\beta_{5}\right|\lesssim\gamma_{5,1}+\gamma_{5,2}+\gamma_{5,3}+\gamma_{5,4}\lesssim\delta v_{i,j}^{\star},
\]
provided that $\theta\lesssim\delta/(\kappa^{2}\mu r)$, $n\gtrsim\delta^{-2}\kappa^{7}\mu^{2}r^{3}\log(n+d)$,
$\varepsilon\lesssim1$, $ndp^{2}\gtrsim\delta^{-2}\kappa^{2}\mu^{2}r^{2}$
and
\begin{equation}
\theta\frac{\omega_{\max}}{\sigma_{r}^{\star}}\sqrt{d}\lesssim\frac{\delta}{\sqrt{\kappa^{3}\mu r}}.\label{eq:theta-omega-sqrt-d-again}
\end{equation}
Note that in view of (\ref{eq:theta-exact}) and (\ref{eq:pca-1st-err-useful}),
(\ref{eq:theta-omega-sqrt-d-again}) is guaranteed by $ndp^{2}\gtrsim\delta^{-2}\kappa^{6}\mu^{3}r^{4}\log^{4}(n+d)$,
\[
\frac{\omega_{\max}}{\sigma_{r}^{\star}}\sqrt{\frac{d}{np}}\lesssim\frac{\delta}{\sqrt{\kappa^{4}\mu r^{2}\log^{2}\left(n+d\right)}},\qquad\frac{\omega_{\max}^{2}}{p\sigma_{r}^{\star2}}\sqrt{\frac{d}{n}}\lesssim\frac{\delta}{\sqrt{\kappa^{4}\mu r^{2}\log^{2}\left(n+d\right)}}.
\]

Similarly we can also show that (which we omit here for brevity) 
\[
\left|\alpha_{7}-\beta_{7}\right|\lesssim\delta v_{i,j}^{\star}
\]
holds true under these conditions.

\paragraph{Step 6: putting everything together.}

We are now ready to combine the above bounds on $\vert\alpha_{k}-\beta_{k}\vert$,
$k=1,\ldots,7$, to conclude that 
\[
\left|v_{i,j}-v_{i,j}^{\star}\right|\leq\sum_{k=1}^{7}\left|\alpha_{k}-\beta_{k}\right|\lesssim\delta v_{i,j}^{\star},
\]
provided that $\varepsilon\lesssim1$, $n\gtrsim\delta^{-2}\kappa^{7}\mu^{2}r^{3}\kappa_{\omega}^{2}\log(n+d)$,
$d\gtrsim\kappa^{2}\mu\log(n+d)$, 
\begin{align*}
ndp^{2}\gtrsim\delta^{-2}\kappa^{6}\mu^{4}r^{6}\kappa_{\omega}^{2}\log^{5}\left(n+d\right), & \qquad np\gtrsim\delta^{-2}\kappa^{6}\mu^{3}r^{5}\kappa_{\omega}^{2}\log^{3}\left(n+d\right),\\
\frac{\omega_{\max}^{2}}{p\sigma_{r}^{\star2}}\sqrt{\frac{d}{n}}\lesssim\frac{\delta}{\kappa^{2}\mu r^{2}\kappa_{\omega}\log^{3/2}\left(n+d\right)}, & \qquad\frac{\omega_{\max}}{\sigma_{r}^{\star}}\sqrt{\frac{d}{np}}\lesssim\frac{\delta}{\kappa^{5/2}\mu r^{2}\kappa_{\omega}\log\left(n+d\right)},
\end{align*}
and 
\[
\left\Vert \bm{U}_{i,\cdot}^{\star}\right\Vert _{2}+\left\Vert \bm{U}_{j,\cdot}^{\star}\right\Vert _{2}\gtrsim\delta^{-1}\kappa\mu r^{2}\kappa_{\omega}\log\left(n+d\right)\left[\frac{\kappa\mu r\log\left(n+d\right)}{\sqrt{nd}p}+\frac{\omega_{\max}^{2}}{p\sigma_{r}^{\star2}}\sqrt{\frac{d}{n}}+\frac{\omega_{\max}}{\sigma_{r}^{\star}}\sqrt{\frac{d}{np}}\right]\sqrt{\frac{r}{d}}.
\]
To finish up, we shall take $\varepsilon\asymp1$, and note that (\ref{eq:ce-var-est-inter-1}),
(\ref{eq:ce-var-est-inter-2}) as well as \ref{eq:ce-var-est-inter-3}
are guaranteed by the conditions of Lemma \ref{lemma:ce-normal-approx-1}.

\subsubsection{Proof of Lemma \ref{lemma:ce-CI-validity}\label{appendix:proof-ce-CI-validity}}

To begin with, we know that
\begin{align}
S_{i,j}^{\star}\in\mathsf{CI}_{i,j}^{1-\alpha}\quad & \Longleftrightarrow\quad S_{i,j}^{\star}\in\left[S_{i,j}\pm\Phi^{-1}\left(1-\alpha/2\right)\sqrt{v_{i,j}}\right]\nonumber \\
 & \Longleftrightarrow\quad\frac{S_{i,j}-S_{i,j}^{\star}}{\sqrt{v_{i,j}}}\in\left[-\Phi^{-1}\left(1-\alpha/2\right),\Phi^{-1}\left(1-\alpha/2\right)\right].\label{eq:ce-CI-validity-inter-0}
\end{align}
Note that we have learned from Lemma \ref{lemma:ce-var-est} that
with probability exceeding $1-O((n+d)^{-10})$, 
\[
\left|v_{i,j}^{\star}-v_{i,j}\right|\lesssim\delta v_{i,j}^{\star},
\]
where $\delta$ is the (unspecified) quantity that has appeared in
Lemma \ref{lemma:ce-var-est}. When $\delta\ll1$, an immediate result
is that $v_{i,j}\asymp v_{i,j}^{\star}$, thus indicating that 
\[
\Delta\coloneqq\left|\frac{S_{i,j}-S_{i,j}^{\star}}{\sqrt{v_{i,j}^{\star}}}-\frac{S_{i,j}-S_{i,j}^{\star}}{\sqrt{v_{i,j}}}\right|=\left|S_{i,j}-S_{i,j}^{\star}\right|\left|\frac{v_{i,j}^{\star}-v_{i,j}}{\sqrt{v_{i,j}^{\star}v_{i,j}}\left(\sqrt{v_{i,j}^{\star}}+\sqrt{v_{i,j}}\right)}\right|\lesssim\delta\left|S_{i,j}-S_{i,j}^{\star}\right|/\sqrt{v_{i,j}^{\star}}.
\]
Suppose for the moment that 
\begin{equation}
\Delta\lesssim\frac{1}{\sqrt{\log\left(n+d\right)}}\label{eq:ce-CI-validity-inter-1}
\end{equation}
holds with probability exceeding $1-O((n+d)^{-10})$. From Lemma \ref{lemma:ce-normal-approx-all}
we know that for any $t\in\mathbb{R}$, 
\begin{equation}
\mathbb{P}\left(\frac{S_{i,j}-S_{i,j}^{\star}}{\sqrt{v_{i,j}^{\star}}}\leq t\right)=\Phi\left(t\right)+O\left(\frac{1}{\sqrt{\log\left(n+d\right)}}\right).\label{eq:ce-CI-validity-inter-2}
\end{equation}
Then we know that for any $t\in\mathbb{R}$, 
\begin{align*}
\mathbb{P}\left(\frac{S_{i,j}-S_{i,j}^{\star}}{\sqrt{v_{i,j}}}\leq t\right) & \leq\mathbb{P}\left(\frac{S_{i,j}-S_{i,j}^{\star}}{\sqrt{v_{i,j}}}\leq t,\Delta\lesssim\frac{1}{\sqrt{\log\left(n+d\right)}}\right)+O\left(\left(n+d\right)^{-10}\right)\\
 & \leq\mathbb{P}\left(\frac{S_{i,j}-S_{i,j}^{\star}}{\sqrt{v_{i,j}^{\star}}}\leq t+\frac{1}{\sqrt{\log\left(n+d\right)}}\right)+O\left(\left(n+d\right)^{-10}\right)\\
 & \overset{\text{(i)}}{\leq}\Phi\left(t+\frac{1}{\sqrt{\log\left(n+d\right)}}\right)+O\left(\frac{1}{\sqrt{\log\left(n+d\right)}}\right)\\
 & \overset{\text{(ii)}}{\leq}\Phi\left(t\right)+O\left(\frac{1}{\sqrt{\log\left(n+d\right)}}\right).
\end{align*}
Here (i) follows from (\ref{eq:ce-CI-validity-inter-2}); (ii) holds
since $\Phi(\cdot)$ is a $1/\sqrt{2\pi}$-Lipschitz continuous function.
Similarly we can show that 
\[
\mathbb{P}\left(\frac{S_{i,j}-S_{i,j}^{\star}}{\sqrt{v_{i,j}}}\leq t\right)\geq\Phi\left(t\right)+O\left(\frac{1}{\sqrt{\log\left(n+d\right)}}\right).
\]
Therefore we have 
\begin{equation}
\mathbb{P}\left(\frac{S_{i,j}-S_{i,j}^{\star}}{\sqrt{v_{i,j}}}\leq t\right)=\Phi\left(t\right)+O\left(\frac{1}{\sqrt{\log\left(n+d\right)}}\right).\label{eq:ce-CI-validity-inter-3}
\end{equation}
By taking $t=\Phi^{-1}(1-\alpha/2)$ in (\ref{eq:ce-CI-validity-inter-3}),
we have 
\begin{equation}
\mathbb{P}\left(\frac{S_{i,j}-S_{i,j}^{\star}}{\sqrt{v_{i,j}}}\leq\Phi^{-1}\left(1-\alpha/2\right)\right)=1-\frac{\alpha}{2}+O\left(\frac{1}{\sqrt{\log\left(n+d\right)}}\right).\label{eq:ce-CI-validity-inter-4}
\end{equation}
In addition, for all $\varepsilon\in(0,1/\sqrt{\log(n+d)}]$, by taking
$t=-\Phi^{-1}(1-\alpha/2)-\varepsilon$ in (\ref{eq:ce-CI-validity-inter-3}),
we have 
\begin{align*}
\mathbb{P}\left(\frac{S_{i,j}-S_{i,j}^{\star}}{\sqrt{v_{i,j}}}\leq-\Phi^{-1}\left(1-\alpha/2\right)-\varepsilon\right) & =\Phi\left(-\Phi^{-1}\left(1-\alpha/2\right)-\varepsilon\right)+O\left(\frac{1}{\sqrt{\log\left(n+d\right)}}\right)\\
 & =\frac{\alpha}{2}+O\left(\frac{1}{\sqrt{\log\left(n+d\right)}}\right),
\end{align*}
where the constant hide in $O(\cdot)$ is independent of $\varepsilon$.
Here the last relation holds since $\Phi(\cdot)$ is a $1/\sqrt{2\pi}$-Lipschitz
continuous function and $\varepsilon\lesssim1/\sqrt{\log(n+d)}$.
As a result, 
\begin{align}
\mathbb{P}\left(\frac{S_{i,j}-S_{i,j}^{\star}}{\sqrt{v_{i,j}}}<-\Phi^{-1}\left(1-\alpha/2\right)\right) & =\lim_{\varepsilon\searrow0}\mathbb{P}\left(\frac{S_{i,j}-S_{i,j}^{\star}}{\sqrt{v_{i,j}}}\leq-\Phi^{-1}\left(1-\alpha/2\right)-\varepsilon\right)\nonumber \\
 & =\frac{\alpha}{2}+O\left(\frac{1}{\sqrt{\log\left(n+d\right)}}\right).\label{eq:ce-CI-validity-inter-5}
\end{align}
Taking (\ref{eq:ce-CI-validity-inter-0}), (\ref{eq:ce-CI-validity-inter-3})
and (\ref{eq:ce-CI-validity-inter-5}) collectively yields 
\begin{align*}
\mathbb{P}\left(S_{i,j}^{\star}\in\mathsf{CI}_{i,j}^{1-\alpha}\right) & =\mathbb{P}\left(\frac{S_{i,j}-S_{i,j}^{\star}}{\sqrt{v_{i,j}}}\in\left[-\Phi^{-1}\left(1-\alpha/2\right),\Phi^{-1}\left(1-\alpha/2\right)\right]\right)\\
 & =\mathbb{P}\left(\frac{S_{i,j}-S_{i,j}^{\star}}{\sqrt{v_{i,j}}}\leq\Phi^{-1}\left(1-\alpha/2\right)\right)-\mathbb{P}\left(\frac{S_{i,j}-S_{i,j}^{\star}}{\sqrt{v_{i,j}}}<-\Phi^{-1}\left(1-\alpha/2\right)\right)\\
 & =1-\alpha+O\left(\frac{1}{\sqrt{\log\left(n+d\right)}}\right).
\end{align*}
as long as (\ref{eq:ce-CI-validity-inter-1}) holds with probability
exceeding $1-O((n+d)^{-10})$. With the above calculations in mind,
everything boils down to bounding the quantity $\Delta$ to the desired
level. Towards this, we are in need of the following lemma.

\begin{claim}\label{claim:ce-1st-err} Suppose that the conditions
of Lemma \ref{lemma:ce-2nd-error} hold. Suppose that $np\gtrsim\log^{4}(n+d)$
and $ndp^{2}\gtrsim\mu r\log^{5}(n+d)$. Then with probability exceeding
$1-O((n+d)^{-10})$, we have 
\begin{align*}
\left|S_{i,j}-S_{i,j}^{\star}\right| & \lesssim\zeta_{i,j}+\left\Vert \bm{U}_{j,\cdot}^{\star}\bm{\Sigma}^{\star}\right\Vert _{2}\left\Vert \bm{U}_{j,\cdot}^{\star}\bm{\Sigma}^{\star}\right\Vert _{2}\left(\sqrt{\frac{\kappa\log\left(n+d\right)}{np}}+\sqrt{\frac{\kappa\mu r\log^{2}\left(n+d\right)}{ndp^{2}}}+\frac{\omega_{\max}}{\sigma_{r}^{\star}}\sqrt{\frac{\log\left(n+d\right)}{np^{2}}}\right).\\
 & \quad+\left(\sqrt{\frac{\kappa}{np}}+\sqrt{\frac{\kappa\mu r\log^{2}\left(n+d\right)}{ndp^{2}}}+\frac{\omega_{\max}}{\sigma_{r}^{\star}}\sqrt{\frac{\log^{2}\left(n+d\right)}{np^{2}}}\right)\sigma_{r}^{\star}\left(\omega_{j}^{\star}\left\Vert \bm{U}_{i,\cdot}^{\star}\right\Vert _{2}+\omega_{i}^{\star}\left\Vert \bm{U}_{j,\cdot}^{\star}\right\Vert _{2}\right).
\end{align*}
\end{claim}

Recall from the proof of Lemma \ref{lemma:ce-normal-approx-1} (more
specifically, Step 3 in Appendix \ref{appendix:proof-ce-normal-approx-1})
that 
\begin{equation}
(v_{i,j}^{\star})^{-1/2}\zeta_{i,j}\leq\widetilde{v}_{i,j}^{-1/2}\zeta_{i,j}\lesssim\frac{1}{\sqrt{\log\left(n+d\right)}}\label{eq:ce-CI-validity-inter-6}
\end{equation}
holds under the assumptions of Lemma \ref{lemma:ce-normal-approx-1}.
In addition, recall from Lemma \ref{lemma:ce-variance-concentration}
that 
\begin{equation}
\sqrt{v_{i,j}^{\star}}\gtrsim\frac{1}{\sqrt{\min\left\{ ndp^{2}\kappa,np\right\} }}\left\Vert \bm{U}_{i,\cdot}^{\star}\bm{\Sigma}^{\star}\right\Vert _{2}\left\Vert \bm{U}_{j,\cdot}^{\star}\bm{\Sigma}^{\star}\right\Vert _{2}+\left(\frac{\sigma_{r}^{\star}}{\sqrt{\min\left\{ ndp^{2},np\right\} }}+\frac{\omega_{\min}}{\sqrt{np^{2}}}\right)\left(\omega_{j}^{\star}\left\Vert \bm{U}_{i,\cdot}^{\star}\right\Vert _{2}+\omega_{i}^{\star}\left\Vert \bm{U}_{j,\cdot}^{\star}\right\Vert _{2}\right).\label{eq:ce-CI-validity-inter-7}
\end{equation}
Therefore we can take (\ref{eq:ce-CI-validity-inter-6}), (\ref{eq:ce-CI-validity-inter-7})
and Claim \ref{claim:ce-1st-err} collectively to show that 
\[
\Delta\lesssim\delta\sqrt{\kappa^{2}\mu r\kappa_{\omega}\log^{2}\left(n+d\right)}+\frac{1}{\sqrt{\log\left(n+d\right)}}\lesssim\frac{1}{\sqrt{\log\left(n+d\right)}}
\]
with probability exceeding $1-O((n+d)^{-10})$ as long as we take
\[
\delta\asymp\frac{1}{\kappa\mu^{1/2}r^{1/2}\kappa_{\omega}^{1/2}\log^{3/2}\left(n+d\right)},
\]
where we use the AM-GM inequality that
\begin{align*}
 & \left\Vert \bm{U}_{j,\cdot}^{\star}\bm{\Sigma}^{\star}\right\Vert _{2}\left\Vert \bm{U}_{j,\cdot}^{\star}\bm{\Sigma}^{\star}\right\Vert _{2}\frac{\omega_{\max}}{\sigma_{r}^{\star}}\sqrt{\frac{\log\left(n+d\right)}{np^{2}}}\lesssim\sqrt{\kappa_{\omega}\log\left(n+d\right)}\left\Vert \bm{U}_{j,\cdot}^{\star}\bm{\Sigma}^{\star}\right\Vert _{2}\left\Vert \bm{U}_{j,\cdot}^{\star}\bm{\Sigma}^{\star}\right\Vert _{2}\frac{\omega_{\min}}{\sigma_{r}^{\star}}\sqrt{\frac{1}{np^{2}}}\\
 & \quad\lesssim\sqrt{\kappa\kappa_{\omega}\log\left(n+d\right)}\left[\frac{1}{\sqrt{np^{2}}}\left\Vert \bm{U}_{i,\cdot}^{\star}\bm{\Sigma}^{\star}\right\Vert _{2}\left\Vert \bm{U}_{j,\cdot}^{\star}\bm{\Sigma}^{\star}\right\Vert _{2}\sqrt{\frac{\mu r}{d}}+\frac{\omega_{\min}^{2}}{\sqrt{np^{2}}}\left(\left\Vert \bm{U}_{i,\cdot}^{\star}\right\Vert _{2}+\left\Vert \bm{U}_{j,\cdot}^{\star}\right\Vert _{2}\right)\right]\\
 & \quad\lesssim\sqrt{\kappa^{2}\mu r\kappa_{\omega}\log\left(n+d\right)}\left[\frac{1}{\sqrt{ndp^{2}\kappa}}\left\Vert \bm{U}_{i,\cdot}^{\star}\bm{\Sigma}^{\star}\right\Vert _{2}\left\Vert \bm{U}_{j,\cdot}^{\star}\bm{\Sigma}^{\star}\right\Vert _{2}+\frac{\omega_{\min}}{\sqrt{np^{2}}}\left(\omega_{j}^{\star}\left\Vert \bm{U}_{i,\cdot}^{\star}\right\Vert _{2}+\omega_{i}^{\star}\left\Vert \bm{U}_{j,\cdot}^{\star}\right\Vert _{2}\right)\right]\\
 & \quad\lesssim\sqrt{\kappa^{2}\mu r\kappa_{\omega}\log\left(n+d\right)}\sqrt{v_{i,j}^{\star}}.
\end{align*}
This in turn confirms (\ref{eq:ce-CI-validity-inter-1}).

It remains to specify the conditions of Lemma \ref{lemma:ce-var-est}
when we choose $\delta\asymp\kappa^{-1}\mu^{-1/2}r^{-1/2}\kappa_{\omega}^{-1/2}\log^{-3/2}(n+d)$.
In this case, these conditions read $n\gtrsim\kappa^{9}\mu^{3}r^{4}\kappa_{\omega}^{3}\log^{4}(n+d)$,
$d\gtrsim\kappa^{2}\mu\log(n+d)$, 
\begin{align*}
ndp^{2}\gtrsim\kappa^{8}\mu^{5}r^{7}\kappa_{\omega}^{3}\log^{8}\left(n+d\right), & \qquad np\gtrsim\kappa^{8}\mu^{4}r^{6}\kappa_{\omega}^{3}\log^{6}\left(n+d\right),\\
\frac{\omega_{\max}^{2}}{p\sigma_{r}^{\star2}}\sqrt{\frac{d}{n}}\lesssim\frac{1}{\kappa^{3}\mu^{3/2}r^{5/2}\kappa_{\omega}^{3/2}\log^{3}\left(n+d\right)}, & \qquad\frac{\omega_{\max}}{\sigma_{r}^{\star}}\sqrt{\frac{d}{np}}\lesssim\frac{1}{\kappa^{7/2}\mu^{3/2}r^{5/2}\kappa_{\omega}^{3/2}\log^{5/2}\left(n+d\right)},
\end{align*}
and 
\[
\left\Vert \bm{U}_{i,\cdot}^{\star}\right\Vert _{2}+\left\Vert \bm{U}_{j,\cdot}^{\star}\right\Vert _{2}\gtrsim\kappa^{2}\mu^{3/2}r^{5/2}\kappa_{\omega}^{3/2}\log^{5/2}\left(n+d\right)\left[\frac{\kappa\mu r\log\left(n+d\right)}{\sqrt{nd}p}+\frac{\omega_{\max}^{2}}{p\sigma_{r}^{\star2}}\sqrt{\frac{d}{n}}+\frac{\omega_{\max}}{\sigma_{r}^{\star}}\sqrt{\frac{d}{np}}\right]\sqrt{\frac{r}{d}},
\]
in addition to all the assumptions of Lemma \ref{lemma:ce-normal-approx-1}.
This concludes the proof of this lemma, as long as Claim \ref{claim:ce-1st-err}
can be established.

\begin{proof}[Proof of Claim \ref{claim:ce-1st-err}]Let us start
by decomposing 
\[
X_{i,j}=\sum_{l=1}^{n}\bigg[\underbrace{M_{j,l}^{\natural}E_{i,l}}_{\eqqcolon a_{l}}+\underbrace{M_{i,l}^{\natural}E_{j,l}}_{\eqqcolon b_{l}}+\underbrace{\sum_{k:k\neq i}E_{i,l}E_{k,l}\left(\bm{U}_{k,\cdot}^{\star}\big(\bm{U}_{j,\cdot}^{\star}\big)^{\top}\right)}_{\eqqcolon c_{l}}+\underbrace{\sum_{k:k\neq j}E_{j,l}E_{k,l}\left(\bm{U}_{k,\cdot}^{\star}\big(\bm{U}_{i,\cdot}^{\star}\big)^{\top}\right)}_{\eqqcolon d_{l}}\bigg].
\]
We shall then bound $\sum_{l=1}^{n}a_{l}$, $\sum_{l=1}^{n}b_{l}$,
$\sum_{l=1}^{n}c_{l}$, and $\sum_{l=1}^{n}d_{l}$ separately. 
\begin{itemize}
\item We first bound $\sum_{l=1}^{n}a_{l}$. It is straightforward to calculate
\begin{align*}
L_{a} & \coloneqq\max_{1\leq l\leq n}\left|a_{l}\right|\leq B_{i}\left\Vert \bm{M}_{j,\cdot}^{\natural}\right\Vert _{\infty}\leq B_{i}\left\Vert \bm{U}_{j,\cdot}^{\natural}\right\Vert _{2}\sigma_{1}^{\natural}\left\Vert \bm{V}^{\natural}\right\Vert _{2,\infty},\\
V_{a} & \coloneqq\sum_{l=1}^{n}\mathsf{var}\left(M_{j,l}^{\natural}E_{i,l}|\bm{F}\right)=\sum_{l=1}^{n}M_{j,l}^{\natural2}\sigma_{i,l}^{2}=\sum_{l=1}^{n}M_{j,l}^{\natural2}\sigma_{i,l}^{2}=\sum_{l=1}^{n}M_{j,l}^{\natural2}\left[\frac{1-p}{np}\left(\bm{U}_{i,\cdot}^{\star}\bm{\Sigma}^{\star}\bm{f}_{l}\right)^{2}+\frac{\omega_{i}^{\star2}}{np}\right],
\end{align*}
where$B_{i}$ and $B_{j}$ are two (random) quantities such that 
\[
\max_{l\in[n]}\left|E_{i,l}\right|\leq B_{i}\qquad\text{and}\qquad\max_{l\in[n]}\left|E_{j,l}\right|\leq B_{j}.
\]
On the event $\mathcal{E}_{\mathsf{good}}$, we know that 
\begin{align*}
L_{a} & \overset{\text{(i)}}{\lesssim}\frac{1}{p}\sqrt{\frac{\log\left(n+d\right)}{n}}\left(\left\Vert \bm{U}_{i,\cdot}^{\star}\bm{\Sigma}^{\star}\right\Vert _{2}+\omega_{i}^{\star}\right)\left\Vert \bm{U}_{j,\cdot}^{\star}\right\Vert _{2}\sigma_{1}^{\star}\sqrt{\frac{\log\left(n+d\right)}{n}},\\
 & \lesssim\frac{\log\left(n+d\right)}{np}\left(\left\Vert \bm{U}_{i,\cdot}^{\star}\bm{\Sigma}^{\star}\right\Vert _{2}+\omega_{i}^{\star}\right)\left\Vert \bm{U}_{j,\cdot}^{\star}\right\Vert _{2}\sigma_{1}^{\star},
\end{align*}
and 
\begin{align*}
V_{a} & \overset{\text{(ii)}}{\lesssim}\sum_{l=1}^{n}M_{j,l}^{\natural2}\left[\frac{1-p}{np}\left\Vert \bm{U}_{i,\cdot}^{\star}\bm{\Sigma}^{\star}\right\Vert _{2}^{2}\log\left(n+d\right)+\frac{\omega_{i}^{\star2}}{np}\right]\\
 & =\left\Vert \bm{M}_{j,\cdot}^{\natural}\right\Vert _{2}^{2}\left[\frac{1-p}{np}\left\Vert \bm{U}_{i,\cdot}^{\star}\bm{\Sigma}^{\star}\right\Vert _{2}^{2}\log\left(n+d\right)+\frac{\omega_{i}^{\star2}}{np}\right]\\
 & \lesssim\left\Vert \bm{U}_{j,\cdot}^{\natural}\right\Vert _{2}^{2}\sigma_{1}^{\natural2}\left[\frac{1-p}{np}\left\Vert \bm{U}_{i,\cdot}^{\star}\bm{\Sigma}^{\star}\right\Vert _{2}^{2}\log\left(n+d\right)+\frac{\omega_{i}^{\star2}}{np}\right]\\
 & \overset{\text{(iii)}}{\lesssim}\left[\frac{1-p}{np}\left\Vert \bm{U}_{i,\cdot}^{\star}\bm{\Sigma}^{\star}\right\Vert _{2}^{2}\log\left(n+d\right)+\frac{\omega_{i}^{\star2}}{np}\right]\left\Vert \bm{U}_{j,\cdot}^{\star}\right\Vert _{2}^{2}\sigma_{1}^{\star2}.
\end{align*}
Here (i) follows from (\ref{eq:good-event-sigma-least-largest-hpca}),
(\ref{eq:good-event-V-natural-2-infty-hpca}) and (\ref{eq:good-event-B-hpca});
(ii) follows from (\ref{eq:good-event-Ui-Sigma-f-hpca}); (iii) follows
from (\ref{eq:good-event-sigma-least-largest-hpca}). Therefore in
view of the Bernstein inequality \citep[Theorem 2.8.4]{vershynin2016high},
conditional on $\bm{F}$, with probability exceeding $1-O((n+d)^{-10})$,
\begin{align*}
\sum_{l=1}^{n}a_{l} & \lesssim\sqrt{V_{a}\log\left(n+d\right)}+L_{a}\log\left(n+d\right)\\
 & \lesssim\left[\sqrt{\frac{\log\left(n+d\right)}{np}}\left\Vert \bm{U}_{i,\cdot}^{\star}\bm{\Sigma}^{\star}\right\Vert _{2}+\frac{\omega_{i}^{\star}}{\sqrt{np}}\right]\left\Vert \bm{U}_{j,\cdot}^{\star}\right\Vert _{2}\sigma_{1}^{\star}\frac{\log^{2}\left(n+d\right)}{np}\left(\left\Vert \bm{U}_{i,\cdot}^{\star}\bm{\Sigma}^{\star}\right\Vert _{2}+\omega_{j}^{\star}\right)\left\Vert \bm{U}_{j,\cdot}^{\star}\right\Vert _{2}\sigma_{1}^{\star}\\
 & \lesssim\left[\sqrt{\frac{\log\left(n+d\right)}{np}}\left\Vert \bm{U}_{i,\cdot}^{\star}\bm{\Sigma}^{\star}\right\Vert _{2}+\frac{\omega_{i}^{\star}}{\sqrt{np}}\right]\left\Vert \bm{U}_{j,\cdot}^{\star}\right\Vert _{2}\sigma_{1}^{\star},
\end{align*}
provided that $np\gtrsim\log^{4}(n+d)$. Similarly we can show that
with probability exceeding $1-O((n+d)^{-10})$, 
\[
\sum_{l=1}^{n}b_{l}\lesssim\left[\sqrt{\frac{\log\left(n+d\right)}{np}}\left\Vert \bm{U}_{j,\cdot}^{\star}\bm{\Sigma}^{\star}\right\Vert _{2}+\frac{\omega_{j}^{\star}}{\sqrt{np}}\right]\left\Vert \bm{U}_{i,\cdot}^{\star}\right\Vert _{2}\sigma_{1}^{\star}.
\]
\item We then move on to bound $\sum_{l=1}^{n}c_{l}$. It has been shown
in (\ref{eq:ce-upper-bound-probabilistic}) that with probability
exceeding $1-O((n+d)^{-101})$, 
\[
\max_{1\leq l\leq n}\left|c_{l}\right|\leq\underbrace{\widetilde{C}\sigma_{\mathsf{ub}}B_{i}\sqrt{\log\left(n+d\right)}\left\Vert \bm{U}_{j,\cdot}^{\star}\right\Vert _{2}+\widetilde{C}BB_{i}\left\Vert \bm{U}_{j,\cdot}^{\star}\right\Vert _{2}\sqrt{\frac{\mu r}{d}}\log\left(n+d\right)}_{\eqqcolon C_{\mathsf{prob}}},
\]
where $\widetilde{C}>0$ is some sufficiently large constant. We also
know from (\ref{eq:ce-upper-bound-deterministic}) that $\max_{1\leq l\leq n}\vert c_{l}\vert$
satisfies the following deterministic bound: 
\[
\max_{1\leq l\leq n}\left|c_{l}\right|\leq\underbrace{B_{i}\sqrt{d}B\left\Vert \bm{U}_{j,\cdot}^{\star}\right\Vert _{2}}_{\eqqcolon C_{\mathsf{det}}}.
\]
Then with probability exceeding $1-O((n+d)^{-101})$, one has 
\[
\sum_{l=1}^{n}c_{l}=\sum_{l=1}^{n}c_{l}\ind_{\left|c_{l}\right|\leq C_{\mathsf{prob}}}.
\]
It is then straightforward to calculate that 
\begin{align*}
L_{c} & \coloneqq\max_{1\leq l\leq n}\left|c_{l}\ind_{\left|c_{l}\right|\leq C_{\mathsf{prob}}}\right|\leq C_{\mathsf{prob}},\\
V_{c} & \coloneqq\sum_{l=1}^{n}\mathsf{var}\left(c_{l}\ind_{\left|c_{l}\right|\leq C_{\mathsf{prob}}}\right)\leq\sum_{l=1}^{n}\mathbb{E}\left[c_{l}^{2}\ind_{\left|c_{l}\right|\leq C_{\mathsf{prob}}}\right]\leq\sum_{l=1}^{n}\mathbb{E}\left[c_{l}^{2}\right]=\sum_{l=1}^{n}\mathsf{var}\left(c_{l}\right)\\
 & =\sum_{l=1}^{n}\mathsf{var}\left[\sum_{k:k\neq i}E_{i,l}E_{k,l}\left(\bm{U}_{k,\cdot}^{\star}\bm{U}_{j,\cdot}^{\star\top}\right)\right]\leq\sum_{l=1}^{n}\sum_{k=1}^{d}\sigma_{i,l}^{2}\sigma_{k,l}^{2}\left(\bm{U}_{k,\cdot}^{\star}\bm{U}_{j,\cdot}^{\star\top}\right)^{2},\\
M_{c} & \coloneqq\sum_{l=1}^{n}\mathbb{E}\left[c_{l}\ind_{\left|c_{l}\right|>C_{\mathsf{prob}}}\right]\leq C_{\mathsf{det}}\sum_{l=1}^{n}\mathbb{P}\left(\left|c_{l}\right|>C_{\mathsf{prob}}\right)\lesssim C_{\mathsf{det}}\left(n+d\right)^{-100}.
\end{align*}
On the event $\mathcal{E}_{\mathsf{good}}$, it is seen that 
\begin{align*}
L_{c} & \lesssim\sigma_{\mathsf{ub}}B_{i}\sqrt{\log\left(n+d\right)}\left\Vert \bm{U}_{j,\cdot}^{\star}\right\Vert _{2}+BB_{i}\left\Vert \bm{U}_{j,\cdot}^{\star}\right\Vert _{2}\sqrt{\frac{\mu r}{d}}\log\left(n+d\right),\\
 & \overset{\text{(i)}}{\lesssim}\sigma_{\mathsf{ub}}B_{i}\sqrt{\log\left(n+d\right)}\left\Vert \bm{U}_{j,\cdot}^{\star}\right\Vert _{2}+\sigma_{\mathsf{ub}}\sqrt{\frac{\log\left(n+d\right)}{p}}B_{i}\left\Vert \bm{U}_{j,\cdot}^{\star}\right\Vert _{2}\sqrt{\frac{\mu r}{d}}\log\left(n+d\right),\\
 & \lesssim\sigma_{\mathsf{ub}}B_{i}\left\Vert \bm{U}_{j,\cdot}^{\star}\right\Vert _{2}\sqrt{\log\left(n+d\right)}\left(1+\sqrt{\frac{\mu r}{dp}}\log\left(n+d\right)\right),
\end{align*}
\begin{align*}
V_{c} & \leq\sigma_{\mathsf{ub}}^{2}\sum_{l=1}^{n}\sigma_{i,l}^{2}\sum_{k=1}^{d}\left(\bm{U}_{k,\cdot}^{\star}\bm{U}_{j,\cdot}^{\star\top}\right)^{2}=n\sigma_{\mathsf{ub}}^{2}\left(\frac{1-p}{np}\left(\bm{U}_{i,\cdot}^{\star}\bm{\Sigma}^{\star}\bm{f}_{l}\right)^{2}+\frac{\omega_{i}^{\star2}}{np}\right)\left\Vert \bm{U}_{j,\cdot}^{\star}\right\Vert _{2}^{2},\\
 & \overset{\text{(ii)}}{\lesssim}n\sigma_{\mathsf{ub}}^{2}\left(\frac{1-p}{np}\left\Vert \bm{U}_{i,\cdot}^{\star}\bm{\Sigma}^{\star}\right\Vert _{2}^{2}\log\left(n+d\right)+\frac{\omega_{i}^{\star2}}{np}\right)\left\Vert \bm{U}_{j,\cdot}^{\star}\right\Vert _{2}^{2},
\end{align*}
and 
\begin{align*}
M_{c} & \lesssim B_{i}\sqrt{d}B\left\Vert \bm{U}_{j,\cdot}^{\star}\right\Vert _{2}\left(n+d\right)^{-100}\overset{\text{(iii)}}{\lesssim}\sigma_{\mathsf{ub}}B_{i}\left\Vert \bm{U}_{j,\cdot}^{\star}\right\Vert _{2}\sqrt{\frac{d\log\left(n+d\right)}{p}}\left(n+d\right)^{-100}\\
 & \overset{\text{(iv)}}{\lesssim}\sigma_{\mathsf{ub}}B_{i}\left\Vert \bm{U}_{j,\cdot}^{\star}\right\Vert _{2}\sqrt{nd}\left(n+d\right)^{-100}\lesssim\sigma_{\mathsf{ub}}B_{i}\left\Vert \bm{U}_{j,\cdot}^{\star}\right\Vert _{2}\left(n+d\right)^{-98}.
\end{align*}
Here, (i) and (iii) follow from (\ref{eq:good-event-B-hpca}) and
(\ref{eq:good-event-sigma-hpca}); (ii) arises from (\ref{eq:good-event-Ui-Sigma-f-hpca});
and (iv) holds true provided that $np\gtrsim\log(n+d)$. Therefore,
by virtue of the Bernstein inequality \citep[Theorem 2.8.4]{vershynin2016high},
we see that conditional on $\bm{F}$, 
\begin{align*}
\sum_{l=1}^{n}c_{l} & =\sum_{l=1}^{n}c_{l}\ind_{\left|c_{l}\right|\leq C_{\mathsf{prob}}}\lesssim M_{c}+\sqrt{V_{c}\log\left(n+d\right)}+L_{c}\log\left(n+d\right)\\
 & \lesssim\sigma_{\mathsf{ub}}B_{i}\left\Vert \bm{U}_{j,\cdot}^{\star}\right\Vert _{2}\left(n+d\right)^{-98}+\sigma_{\mathsf{ub}}\left\Vert \bm{U}_{j,\cdot}^{\star}\right\Vert _{2}\sqrt{\frac{\log\left(n+d\right)}{p}}\left(\left\Vert \bm{U}_{i,\cdot}^{\star}\bm{\Sigma}^{\star}\right\Vert _{2}\sqrt{\log\left(n+d\right)}+\omega_{i}^{\star}\right)\\
 & \quad+\sigma_{\mathsf{ub}}B_{i}\left\Vert \bm{U}_{j,\cdot}^{\star}\right\Vert _{2}\log^{3/2}\left(n+d\right)\left(1+\sqrt{\frac{\mu r}{dp}}\log\left(n+d\right)\right)\\
 & \overset{\text{(i)}}{\lesssim}\sigma_{\mathsf{ub}}\left\Vert \bm{U}_{j,\cdot}^{\star}\right\Vert _{2}\left(\left\Vert \bm{U}_{i,\cdot}^{\star}\bm{\Sigma}^{\star}\right\Vert _{2}+\omega_{i}^{\star}\right)\sqrt{\frac{\log\left(n+d\right)}{p}}\\
 & \lesssim\left\Vert \bm{U}_{j,\cdot}^{\star}\right\Vert _{2}\left(\left\Vert \bm{U}_{i,\cdot}^{\star}\bm{\Sigma}^{\star}\right\Vert _{2}+\omega_{i}^{\star}\right)\left(\sqrt{\frac{\mu r\log^{2}\left(n+d\right)}{ndp^{2}}}\sigma_{1}^{\star}+\omega_{\max}\sqrt{\frac{\log\left(n+d\right)}{np^{2}}}\right)
\end{align*}
holds with probability exceeding $1-O((n+d)^{-10})$. Here, (i) holds
true provided that $np\gtrsim\log^{3}(n+d)$ and $ndp^{2}\gtrsim\mu r\log^{5}(n+d)$.
Similarly we can demonstrate that with probability exceeding $1-O((n+d)^{-10})$,
\[
\sum_{l=1}^{n}d_{l}\lesssim\left\Vert \bm{U}_{i,\cdot}^{\star}\right\Vert _{2}\left(\left\Vert \bm{U}_{j,\cdot}^{\star}\bm{\Sigma}^{\star}\right\Vert _{2}+\omega_{j}^{\star}\right)\left(\sqrt{\frac{\mu r\log^{2}\left(n+d\right)}{ndp^{2}}}\sigma_{1}^{\star}+\omega_{\max}\sqrt{\frac{\log\left(n+d\right)}{np^{2}}}\right).
\]
Therefore, the above two bounds taken together lead to 
\begin{align*}
\sum_{l=1}^{n}\left(c_{l}+d_{l}\right) & \lesssim\left[\frac{1}{\sigma_{r}^{\star}}\left\Vert \bm{U}_{j,\cdot}^{\star}\bm{\Sigma}^{\star}\right\Vert _{2}\left\Vert \bm{U}_{j,\cdot}^{\star}\bm{\Sigma}^{\star}\right\Vert _{2}+\left(\omega_{j}^{\star}\left\Vert \bm{U}_{i,\cdot}^{\star}\right\Vert _{2}+\omega_{i}^{\star}\left\Vert \bm{U}_{j,\cdot}^{\star}\right\Vert _{2}\right)\right]\\
 & \quad\cdot\left(\sqrt{\frac{\mu r\log^{2}\left(n+d\right)}{ndp^{2}}}\sigma_{1}^{\star}+\omega_{\max}\sqrt{\frac{\log\left(n+d\right)}{np^{2}}}\right).
\end{align*}
\[
\sum_{l=1}^{n}\left(a_{l}+b_{l}\right)\lesssim\sqrt{\frac{\kappa\log\left(n+d\right)}{np}}\left\Vert \bm{U}_{j,\cdot}^{\star}\bm{\Sigma}^{\star}\right\Vert _{2}\left\Vert \bm{U}_{j,\cdot}^{\star}\bm{\Sigma}^{\star}\right\Vert _{2}+\frac{\sigma_{1}^{\star}}{\sqrt{np}}\left(\omega_{j}^{\star}\left\Vert \bm{U}_{i,\cdot}^{\star}\right\Vert _{2}+\omega_{i}^{\star}\left\Vert \bm{U}_{j,\cdot}^{\star}\right\Vert _{2}\right),
\]
\end{itemize}
Therefore, putting the above bounds together, one can conclude that
\begin{align*}
\left|X_{i,j}\right| & \leq\sum_{l=1}^{n}\left(a_{l}+b_{l}+c_{l}+d_{l}\right)\\
 & \lesssim\left\Vert \bm{U}_{j,\cdot}^{\star}\bm{\Sigma}^{\star}\right\Vert _{2}\left\Vert \bm{U}_{j,\cdot}^{\star}\bm{\Sigma}^{\star}\right\Vert _{2}\left(\sqrt{\frac{\kappa\log\left(n+d\right)}{np}}+\sqrt{\frac{\kappa\mu r\log^{2}\left(n+d\right)}{ndp^{2}}}\right)\\
 & \quad+\sigma_{r}^{\star}\left(\omega_{j}^{\star}\left\Vert \bm{U}_{i,\cdot}^{\star}\right\Vert _{2}+\omega_{i}^{\star}\left\Vert \bm{U}_{j,\cdot}^{\star}\right\Vert _{2}\right)\left(\sqrt{\frac{\kappa}{np}}+\sqrt{\frac{\kappa\mu r\log^{2}\left(n+d\right)}{ndp^{2}}}\right)\\
 & \quad+\left\Vert \bm{U}_{j,\cdot}^{\star}\bm{\Sigma}^{\star}\right\Vert _{2}\left\Vert \bm{U}_{j,\cdot}^{\star}\bm{\Sigma}^{\star}\right\Vert _{2}\frac{\omega_{\max}}{\sigma_{r}^{\star}}\sqrt{\frac{\log\left(n+d\right)}{np^{2}}}+\left(\omega_{j}^{\star}\left\Vert \bm{U}_{i,\cdot}^{\star}\right\Vert _{2}+\omega_{i}^{\star}\left\Vert \bm{U}_{j,\cdot}^{\star}\right\Vert _{2}\right)\omega_{\max}\sqrt{\frac{\log\left(n+d\right)}{np^{2}}}\\
 & \lesssim\left\Vert \bm{U}_{j,\cdot}^{\star}\bm{\Sigma}^{\star}\right\Vert _{2}\left\Vert \bm{U}_{j,\cdot}^{\star}\bm{\Sigma}^{\star}\right\Vert _{2}\left(\sqrt{\frac{\kappa\log\left(n+d\right)}{np}}+\sqrt{\frac{\kappa\mu r\log^{2}\left(n+d\right)}{ndp^{2}}}+\frac{\omega_{\max}}{\sigma_{r}^{\star}}\sqrt{\frac{\log\left(n+d\right)}{np^{2}}}\right)\\
 & \quad+\left[\sqrt{\frac{\kappa}{np}}+\sqrt{\frac{\kappa\mu r\log^{2}\left(n+d\right)}{ndp^{2}}}+\frac{\omega_{\max}}{\sigma_{r}^{\star}}\sqrt{\frac{\log^{2}\left(n+d\right)}{np^{2}}}\right]\sigma_{r}^{\star}\left(\omega_{j}^{\star}\left\Vert \bm{U}_{i,\cdot}^{\star}\right\Vert _{2}+\omega_{i}^{\star}\left\Vert \bm{U}_{j,\cdot}^{\star}\right\Vert _{2}\right),
\end{align*}
Combine the above result with Lemma \ref{lemma:ce-2nd-error} to arrive
at 
\begin{align*}
\left|S_{i,j}-S_{i,j}^{\star}\right| & \lesssim\zeta_{i,j}+\left\Vert \bm{U}_{j,\cdot}^{\star}\bm{\Sigma}^{\star}\right\Vert _{2}\left\Vert \bm{U}_{j,\cdot}^{\star}\bm{\Sigma}^{\star}\right\Vert _{2}\left(\sqrt{\frac{\kappa\log\left(n+d\right)}{np}}+\sqrt{\frac{\kappa\mu r\log^{2}\left(n+d\right)}{ndp^{2}}}+\frac{\omega_{\max}}{\sigma_{r}^{\star}}\sqrt{\frac{\log\left(n+d\right)}{np^{2}}}\right).\\
 & \quad+\left(\sqrt{\frac{\kappa}{np}}+\sqrt{\frac{\kappa\mu r\log^{2}\left(n+d\right)}{ndp^{2}}}+\frac{\omega_{\max}}{\sigma_{r}^{\star}}\sqrt{\frac{\log^{2}\left(n+d\right)}{np^{2}}}\right)\sigma_{r}^{\star}\left(\omega_{j}^{\star}\left\Vert \bm{U}_{i,\cdot}^{\star}\right\Vert _{2}+\omega_{i}^{\star}\left\Vert \bm{U}_{j,\cdot}^{\star}\right\Vert _{2}\right).
\end{align*}
\end{proof} 

\section{Other useful lemmas\label{appendix:technical_lemmas}}

In this section, we gather a few useful results from prior literature
that prove useful for our analysis. The first theorem is a non-asymptotic
version of the Berry-Esseen Theorem, which has been established using
Stein's method; see \citet[Theorem 3.7]{chen2010normal}.

\begin{theorem}\label{thm:berry-esseen}Let $\xi_{1},\ldots,\xi_{n}$
be independent random variables with zero means, satisfying $\sum_{i=1}^{n}\mathsf{var}(\xi_{i})=1$.
Then $W=\sum_{i=1}^{n}\xi_{i}$ satisfies 
\[
\sup_{z\in\mathbb{R}}\big|\,\mathbb{P}\left(W\leq z\right)-\Phi\left(z\right)\big|\leq10\gamma,\qquad\text{where}\ \gamma=\sum_{i=1}^{n}\mathbb{E}\left[\left|\xi_{i}\right|^{3}\right].
\]
\end{theorem}

The next theorem, which we borrow from \citet[Theorem 1.1]{raivc2019multivariate},
generalizes Theorem~\ref{thm:berry-esseen} to the vector case.

\begin{theorem}\label{thm:berry-esseen-multivariate}Let $\bm{\xi}_{1},\ldots,\bm{\xi}_{n}$
be independent, $\mathbb{R}^{d}$-valued random vectors with zero
means. Let $\bm{W}=\sum_{i=1}^{n}\bm{\xi}_{i}$ and assume $\bm{\Sigma}=\mathsf{cov}(\bm{W})$
is invertible. Let $\bm{Z}$ be a $d$-dimensional Gaussian random
vector with zero mean and covariance matrix $\bm{\Sigma}$. Then we
have

\[
\sup_{\mathcal{C}\in\mathscr{C}^{d}}\left|\mathbb{P}\left(\bm{W}\in\mathcal{C}\right)-\mathbb{P}\left(\bm{Z}\in\mathcal{C}\right)\right|\leq\left(42d^{1/4}+16\right)\gamma,\qquad\text{where}\ \gamma=\sum_{i=1}^{n}\mathbb{E}\left[\left\Vert \bm{\Sigma}^{-1/2}\bm{\xi}_{i}\right\Vert _{2}^{3}\right].
\]
Here, $\mathscr{C}^{d}$ represents the set of all convex sets in
$\mathbb{R}^{d}$.\end{theorem}

Moreover, deriving our distributional theory involves some basic results
about the total-variation distance between two Gaussian distributions,
as stated below. This result can be found in \citet[Theorem 1.1]{devroye2018total}.

\begin{theorem}\label{thm:gaussian-TV-distance}Let $\bm{\mu}\in\mathbb{R}^{d}$,
$\bm{\Sigma}_{1}$ and $\bm{\Sigma}_{2}$ be positive definite $d\times d$
matrices. Then the total-variation distance between $\mathcal{N}\left(\bm{\mu},\bm{\Sigma}_{1}\right)$
and $\mathcal{N}\left(\bm{\mu},\bm{\Sigma}_{2}\right)$ --- denoted
by $\mathsf{TV}\left(\mathcal{N}\left(\bm{\mu},\bm{\Sigma}_{1}\right),\mathcal{N}\left(\bm{\mu},\bm{\Sigma}_{2}\right)\right)$
--- satisfies 
\[
\frac{1}{100}\leq\frac{\mathsf{TV}\left(\mathcal{N}\left(\bm{\mu},\bm{\Sigma}_{1}\right),\mathcal{N}\left(\bm{\mu},\bm{\Sigma}_{2}\right)\right)}{\min\left\{ 1,\left\Vert \bm{\Sigma}_{1}^{-1/2}\bm{\Sigma}_{2}\bm{\Sigma}_{1}^{-1/2}-\bm{I}_{d}\right\Vert _{\mathrm{F}}\right\} }\leq\frac{3}{2}.
\]
\end{theorem}

Recall the definition of $\mathcal{C}^{\varepsilon}$ in \eqref{eq:defn-C-epsilon}.
Then for any convex set $\mathcal{C}\in\mathscr{C}^{d}$, let us define
the following quantity related to Gaussian distributions 
\begin{equation}
\gamma\left(\mathcal{C}\right)\coloneqq\sup_{\varepsilon>0}\max\left\{ \frac{1}{\varepsilon}\mathcal{N}\left(\bm{0},\bm{I}_{d}\right)\left\{ \mathcal{C}^{\varepsilon}\setminus\mathcal{C}\right\} ,\frac{1}{\varepsilon}\mathcal{N}\left(\bm{0},\bm{I}_{d}\right)\left\{ \mathcal{C}\setminus\mathcal{C}^{-\varepsilon}\right\} \right\} ,\label{eq:defn-gamma-convex-set}
\end{equation}
and, in addition, 
\begin{equation}
\gamma_{d}\coloneqq\sup_{\mathcal{C}\in\mathscr{C}^{d}}\gamma\left(\mathcal{C}\right).\label{eq:defn-gamma-d-convex}
\end{equation}
The following theorem from \citet[Theorem 1.2]{raivc2019multivariate}
delivers an upper bound on the quantity $\gamma_{d}$.

\begin{theorem} \label{thm:gaussian-perimeter}For all $d\in\mathbb{N}$,
we have 
\[
\gamma_{d}<0.59d^{1/4}+0.21.
\]
\end{theorem}

Finally, we are in need of the following basic lemma in order to translate
results derived for bounded random variables to the ones concerned
with sub-Gaussian random variables.

\begin{lemma}\label{lemma:truncate}There exist two universal constants
$C_{\delta},C_{\sigma}>0$ such that: for any sub-Gaussian random
variable $X$ with $\mathbb{E}[X]=0$, $\mathsf{var}(X)=\sigma^{2}$,
$\Vert X\Vert_{\psi_{2}}\lesssim\sigma$, and any $\delta\in(0,C_{\delta}\sigma)$,
one can construct a random variable $\widetilde{X}$ satisfying the
following properties: 
\begin{enumerate}
\item $\widetilde{X}$ is equal to $X$ with probability at least $1-\delta$; 
\item $\mathbb{E}[\widetilde{X}]=0$; 
\item $\widetilde{X}$ is a bounded random variable: $\vert\widetilde{X}\vert\leq C_{\sigma}\sigma\sqrt{\log(\delta^{-1})}$; 
\item The variance of $\widetilde{X}$ obeys: $\mathsf{var}(\widetilde{X})=(1+O(\sqrt{\delta}))\sigma^{2}$; 
\item $\widetilde{X}$ is a sub-Gaussian random variable obeying $\Vert\widetilde{X}\Vert_{\psi_{2}}\lesssim\sigma$. 
\end{enumerate}
\end{lemma}

\subsection{Proof of Lemma \ref{lemma:truncate}}

\paragraph{Step 1: lower bounding $\mathbb{E}\vert X\vert$.}

For any $t>0$, it is easily seen that 
\[
\mathbb{E}\left[X^{2}\ind_{\left|X\right|>t}\right]\overset{\text{(i)}}{\leq}\left(\mathbb{E}\left[X^{4}\right]\right)^{\frac{1}{2}}\left(\mathbb{P}\left(\left|X\right|>t\right)\right)^{\frac{1}{2}}\overset{\text{(ii)}}{\lesssim}\sigma^{2}\exp\left(-\frac{t^{2}}{C\sigma^{2}}\right),
\]
where $C>0$ is some absolute constant. Here, (i) results from the
Cauchy-Schwarz inequality, whereas (ii) follows from standard properties
of sub-Gaussian random variables \citep[Proposition 2.5.2]{vershynin2016high}.
By taking $t=c_{\mathsf{lb}}^{-1}\sigma/4$ for some sufficiently
small constant $c_{\mathsf{lb}}>0$, we can guarantee that $\mathbb{E}[X^{2}\ind_{\vert X\vert>t}]\leq\sigma^{2}/2$,
which in turn results in 
\[
\mathbb{E}\left[X^{2}\ind_{\left|X\right|\leq t}\right]=\mathbb{E}\left[X^{2}\right]-\mathbb{E}\left[X^{2}\ind_{\left|X\right|>t}\right]\geq\sigma^{2}-\frac{1}{2}\sigma^{2}\geq\frac{1}{2}\sigma^{2}.
\]
As a consequence, we obtain (with the above choice $t=c_{\mathsf{lb}}^{-1}\sigma/4$)
\begin{equation}
\mathbb{E}\left[\left|X\right|\right]\geq\mathbb{E}\left[\left|X\right|\ind_{\left|X\right|\leq t}\right]\geq\frac{\mathbb{E}\left[X^{2}\ind_{\left|X\right|\leq t}\right]}{t}\geq2c_{\mathsf{lb}}\sigma.\label{eq:EX-lower-bound-146}
\end{equation}

\paragraph{Step 2: constructing $\widetilde{X}$ by truncating $X$ randomly.}

For notational simplicity, let us define $X^{+}=X\lor0$ and $X^{-}=(-X)\lor0$.
Given that $\mathbb{E}[X]=0$, we have 
\begin{equation}
\mathbb{E}\left[X^{+}\right]=\mathbb{E}\left[X^{-}\right]=\frac{1}{2}\mathbb{E}\left[\left|X\right|\right]\geq c_{\mathsf{lb}}\sigma,\label{eq:X-moment-lb}
\end{equation}
where the last inequality comes from \eqref{eq:EX-lower-bound-146}.
Define a function $f:\mathbb{R}^{+}\mapsto\mathbb{R}^{+}$ as follows
\[
f\left(x\right)\coloneqq\mathbb{E}\left[X\ind_{X\geq x}\right].
\]
It is straightforward to check that $f(x)$ is monotonically non-increasing
within the domain $x\in[0,\infty)$. In view of the monotone convergence
theorem, we know that $f(x)$ is a left continuous function with $\lim_{x\searrow0}f(x)=\mathbb{E}[X^{+}]$
and $\lim_{x\to+\infty}f(x)=0$. In addition, for any $x\in\mathbb{R}^{+}$,
one has 
\begin{equation}
f\left(x\right)=\mathbb{E}\left[X\ind_{X\geq x}\right]\overset{\text{(i)}}{\leq}\left(\mathbb{E}\left[X^{2}\right]\right)^{\frac{1}{2}}\left(\mathbb{P}\left(X\geq x\right)\right)^{\frac{1}{2}}\overset{\text{(ii)}}{\leq}\sigma\exp\left(-\frac{x^{2}}{C\sigma^{2}}\right),\label{eq:truncate-inter-1}
\end{equation}
where $C>0$ is some absolute constant. Here, (i) comes from the Cauchy-Schwarz
inequality, while (ii) follows from standard properties of sub-Gaussian
random variables \citep[Proposition 2.5.2]{vershynin2016high}.

For any given $\varepsilon\in(0,c_{\mathsf{lb}}\sigma/2)$, we take
\[
x_{\varepsilon}\coloneqq\sup\left\{ x\in\mathbb{R}^{+}:f\left(x\right)\geq\varepsilon\right\} .
\]
Since $f$ is known to be left continuous, we know that 
\[
\lim_{x\nearrow x_{\varepsilon}}f\left(x\right)=f\left(x_{\varepsilon}\right)\geq\varepsilon\geq\lim_{x\searrow x_{\varepsilon}}f\left(x\right).
\]
Taking the definition of $x_{\varepsilon}$ and (\ref{eq:truncate-inter-1})
collectively yields 
\[
\varepsilon\leq f\left(x_{\varepsilon}\right)\leq\sigma\exp\left(-\frac{x_{\varepsilon}^{2}}{C\sigma^{2}}\right),
\]
which further gives an upper bound on $x_{\varepsilon}$ as follows
\begin{equation}
x_{\varepsilon}\leq\sigma\sqrt{C\log\left(\sigma/\varepsilon\right)}.\label{eq:x-varepsilon-ub}
\end{equation}
In addition, we can also lower bound $x_{\varepsilon}$ by observing
that 
\begin{align*}
\mathbb{E}\left[X^{+}\right] & =\mathbb{E}\left[X^{+}\ind_{X\leq x_{\varepsilon}}\right]+\mathbb{E}\left[X\ind_{X>x_{\varepsilon}}\right]=\mathbb{E}\left[X^{+}\ind_{X\leq x_{\varepsilon}}\right]+\lim_{x\to x_{\varepsilon}^{+}}f\left(x\right)\\
 & \leq x_{\varepsilon}+\varepsilon\leq x_{\varepsilon}+\frac{1}{2}c_{\mathsf{lb}}\sigma,
\end{align*}
which taken collectively with (\ref{eq:X-moment-lb}) yields 
\begin{equation}
x_{\varepsilon}\geq\frac{1}{2}c_{\mathsf{lb}}\sigma.\label{eq:x-varepsilon-lb}
\end{equation}
With these calculations in place, we define $\widetilde{X}^{+}$ as
follows: 
\begin{itemize}
\item If $\lim_{x\searrow x_{\varepsilon}}f(x)=f\left(x_{\varepsilon}\right)$,
we immediately know that $f(x_{\varepsilon})=\varepsilon$. Then we
can set 
\[
\widetilde{X}^{+}\coloneqq X^{+}\ind_{X^{+}<x_{\varepsilon}}.
\]
This construction gives $\widetilde{X}^{+}\leq x_{\varepsilon}\leq\sigma\sqrt{C\log(\sigma/\varepsilon)}$
and 
\[
\mathbb{E}\left[\widetilde{X}^{+}\right]=\mathbb{E}\left[X^{+}\right]-f\left(x_{\varepsilon}\right)=\mathbb{E}\left[X^{+}\right]-\varepsilon.
\]
We can also derive from (\ref{eq:x-varepsilon-lb}) that 
\[
\mathbb{P}\left(\widetilde{X}^{+}\neq X^{+}\right)=\mathbb{P}\left(X^{+}\geq x_{\varepsilon}\right)\leq\frac{\mathbb{E}\left[X^{+}\ind_{X^{+}\geq x_{\varepsilon}}\right]}{x_{\varepsilon}}=\frac{f\left(x_{\varepsilon}\right)}{x_{\varepsilon}}=\frac{\varepsilon}{x_{\varepsilon}}\leq\frac{2\varepsilon}{c_{\mathsf{lb}}\sigma}.
\]
\item If $\lim_{x\searrow x_{\varepsilon}}f(x)<f\left(x_{\varepsilon}\right)$,
we know that 
\begin{equation}
\mathbb{E}\left[X^{+}\ind_{X^{+}=x_{\varepsilon}}\right]=f\left(x_{\varepsilon}\right)-\lim_{x\searrow x_{\varepsilon}}>0.\label{eq:truncate-inter-2}
\end{equation}
Then one can set 
\[
\widetilde{X}^{+}\coloneqq X^{+}\ind_{X^{+}<x_{\varepsilon}}+X\ind_{X^{+}=x_{\varepsilon}}Q,
\]
where $Q$ is a Bernoulli random variable (independent of $X$) with
parameter 
\[
q=\frac{f\left(x_{\varepsilon}\right)-\varepsilon}{f\left(x_{\varepsilon}\right)-\lim_{x\searrow x_{\varepsilon}}},
\]
i.e., $\mathbb{P}(Q=1)=1-\mathbb{P}(Q=0)=q$. This construction gives
$\widetilde{X}^{+}\leq x_{\varepsilon}\leq\sigma\sqrt{C\log(\sigma/\varepsilon)}$
and 
\begin{align}
\mathbb{E}\left[\widetilde{X}^{+}\right] & =\mathbb{E}\left[X^{+}\ind_{X^{+}<x_{\varepsilon}}\right]+q\mathbb{E}\left[X\ind_{X^{+}=x_{\varepsilon}}\right]\label{eq:truncate-inter-3}\\
 & =\mathbb{E}\left[X^{+}\right]-f\left(x_{\varepsilon}\right)+q\left[f\left(x_{\varepsilon}\right)-\lim_{x\to x_{\varepsilon}^{+}}f\left(x\right)\right]\nonumber \\
 & =\mathbb{E}\left[X^{+}\right]-\varepsilon.\label{eq:truncate-inter-4}
\end{align}
In addition, we have 
\begin{align*}
\mathbb{P}\left(\widetilde{X}^{+}\neq X^{+}\right) & =\mathbb{P}\left(X^{+}>x_{\varepsilon}\right)+\mathbb{P}\left(X^{+}=x_{\varepsilon},Q=0\right)=\mathbb{P}\left(X^{+}>x_{\varepsilon}\right)+\left(1-q\right)\mathbb{P}\left(X^{+}=x_{\varepsilon}\right)\\
 & \overset{\text{(i)}}{\leq}\frac{\mathbb{E}\left[X^{+}\ind_{X^{+}>x_{\varepsilon}}\right]}{x_{\varepsilon}}+\left(1-q\right)\frac{\mathbb{E}\left[X^{+}\ind_{X^{+}=x_{\varepsilon}}\right]}{x_{\varepsilon}}\\
 & =\frac{\mathbb{E}\left[X^{+}\right]-\mathbb{E}\left[X^{+}\ind_{X^{+}<x_{\varepsilon}}\right]-q\mathbb{E}\left[X^{+}\ind_{X^{+}=x_{\varepsilon}}\right]}{x_{\varepsilon}}\overset{\text{(ii)}}{=}\frac{\mathbb{E}\left[X^{+}\right]-\mathbb{E}\left[\widetilde{X}^{+}\right]}{x_{\varepsilon}}\\
 & \overset{\text{(iii)}}{=}\frac{\varepsilon}{x_{\varepsilon}}\overset{\text{(iv)}}{\leq}\frac{2\varepsilon}{c_{\mathsf{lb}}\sigma}.
\end{align*}
Here, (i) holds since $\mathbb{E}[X^{+}\ind_{X^{+}>x_{\varepsilon}}]\geq x_{\varepsilon}\mathbb{P}(X^{+}>x_{\varepsilon})$
and $\mathbb{E}[X^{+}\ind_{X^{+}=x_{\varepsilon}}]=x_{\varepsilon}\cdot\mathbb{P}(X^{+}=x_{\varepsilon})$;
(ii) follows from (\ref{eq:truncate-inter-3}); (iii) is a consequence
of (\ref{eq:truncate-inter-4}); and (iv) follows from (\ref{eq:x-varepsilon-lb}). 
\end{itemize}
We have thus constructed a random variable $\widetilde{X}^{+}$ that
satisfies: (i) $\widetilde{X}^{+}$ equals either $X^{+}$ or $0$;
(ii) $\mathbb{P}(\widetilde{X}^{+}\neq X^{+})\leq2\varepsilon/(c_{\mathsf{lb}}\sigma)$;
(iii) $\mathbb{E}[\widetilde{X}^{+}]=\mathbb{E}[X^{+}]-\varepsilon$;
and (iv) $0\leq\widetilde{X}^{+}\leq\sigma\sqrt{C\log(\sigma/\varepsilon)}$.
Similarly, we can also construct another random variable $\widetilde{X}^{-}$
satisfying: (i) $\widetilde{X}^{-}$ equals either $X^{-}$ or $0$;
(ii) $\mathbb{P}(\widetilde{X}^{-}\neq X^{-})\leq2\varepsilon/(c_{\mathsf{lb}}\sigma)$;
(iii) $\mathbb{E}[\widetilde{X}^{-}]=\mathbb{E}[X^{-}]-\varepsilon$;
and (iv) $0\leq\widetilde{X}^{-}\leq\sigma\sqrt{C\log(\sigma/\varepsilon)}$.
We shall then construct $\widetilde{X}$ as follows 
\[
\widetilde{X}\coloneqq\widetilde{X}^{+}-\widetilde{X}^{-}.
\]

\paragraph{Step 3: verifying the advertised properties of $\widetilde{X}$.}

To finish up, we can check that the following properties are satisfied: 
\begin{enumerate}
\item $\widetilde{X}$ has mean zero, namely, 
\[
\mathbb{E}\left[\widetilde{X}\right]=\mathbb{E}\left[\widetilde{X}^{+}\right]-\mathbb{E}\left[\widetilde{X}^{-}\right]=\mathbb{E}\left[X^{+}\right]-\varepsilon-\mathbb{E}\left[X^{-}\right]+\varepsilon=\mathbb{E}\left[X\right]=0.
\]
\item $\widetilde{X}$ is identical to $X$ with high probability, namely,
\[
\mathbb{P}\left(X\neq\widetilde{X}\right)=\mathbb{P}\left(X^{+}\neq\widetilde{X}^{+}\right)+\mathbb{P}\left(X^{-}\neq\widetilde{X}^{-}\right)\leq\frac{4\varepsilon}{c_{\mathsf{lb}}\sigma}.
\]
\item $\widetilde{X}$ is a bounded random variable in the sense that $\vert\widetilde{X}\vert\leq\sigma\sqrt{C\log(\sigma/\varepsilon)}\lesssim\sigma\sqrt{\log(\sigma/\varepsilon)}$. 
\item The variance of $\widetilde{X}$ is close to $\sigma^{2}$ in the
sense that 
\[
\mathsf{var}\big(\widetilde{X}\big)=\mathbb{E}\left[\widetilde{X}^{2}\right]=\mathbb{E}\left[X^{2}\right]-\mathbb{E}\left[X^{2}\ind_{X\neq\widetilde{X}}\right]=\left(1+O\left(\sqrt{\varepsilon/\sigma}\right)\right)\sigma^{2},
\]
where the last relation holds due to the following observation 
\[
\mathbb{E}\left[X^{2}\ind_{X\neq\widetilde{X}}\right]\overset{\text{(i)}}{\leq}\left(\mathbb{E}\left[X^{4}\right]\right)^{\frac{1}{2}}\left(\mathbb{P}\left(X\neq\widetilde{X}\right)\right)^{\frac{1}{2}}\overset{\text{(ii)}}{\lesssim}\sigma^{2}\sqrt{\frac{\varepsilon}{\sigma}}.
\]
Here, (i) invokes Cauchy-Schwarz, whereas (ii) is valid due to standard
properties of sub-Gaussian random variables \citep[Proposition 2.5.2]{vershynin2016high}. 
\item By construction, we can see that for any $t\geq0$, 
\[
\mathbb{P}\left(\left|\widetilde{X}\right|\geq t\right)\leq\mathbb{P}\left(\left|X\right|\geq t\right)\leq2\exp\left(-\frac{t^{2}}{\widetilde{c}\sigma^{2}}\right)
\]
holds for some absolute constant $\widetilde{c}>0$, where the last
relation follows from \citet[Proposition 2.5.2]{vershynin2016high}.
By invoking the definition of sub-Gaussian random variables \citep[Definition 2.5.6]{vershynin2016high}
as well as standard properties of sub-Gaussian random variables \citep[Proposition 2.5.2]{vershynin2016high},
we can conclude that $\widetilde{X}$ is sub-Gaussian obeying $\Vert\widetilde{X}\Vert_{\psi_{2}}\lesssim\sigma$. 
\end{enumerate}
By taking $\varepsilon=\delta c_{\mathsf{lb}}\sigma/4$ for any $\delta\in(0,2)$,
we establish the desired result.

\section{Extensions and additional discussions}
In this section, we briefly discuss a couple of potential extensions of our algorithms and theory, 
and provide a few technical remarks about our settings and assumptions.

\subsection{Relaxing the Gaussian design} \label{subsec:subgaussian-design}

In fact, all results in this paper can be generalized to the following setting: 
\begin{itemize}
	\item the samples are generated such that
	\[
	\bm{x}_{j}=\bm{U}^{\star}(\bm{\Lambda}^{\star})^{1/2}\bm{f}_{j},\qquad1\leq j\leq n,
	\]
	where $\bm{f}_{1},\ldots,\bm{f}_{n}$ are independent sub-Gaussian
	random vectors in $\mathbb{R}^{r}$ with independent entries satisfying
	\[
	\mathbb{E}\left[f_{i,j}\right]=0,\qquad\mathbb{E}\left[f_{i,j}^{2}\right]=1,\qquad\mathbb{E}\left[f_{i,j}^{4}\right]=M_{4},\qquad\text{and}\qquad\Vert f_{i,j}\Vert_{\psi_{2}}=O\left(1\right)
	\]
	for any $1\leq i\leq n$ and $1\leq j\leq r$. 
\end{itemize}
Our distributional theory for the principal subspace $\bm{U}^{\star}$ (cf.~Theorem
\ref{thm:pca}) continues to hold with
\begin{align}
	\bm{\Sigma}_{U,l}^{\star} & \coloneqq\left(\frac{1-p}{np}\left\Vert \bm{U}_{l,\cdot}^{\star}\bm{\Sigma}^{\star}\right\Vert _{2}^{2}+\frac{\omega_{l}^{\star2}}{np}\right)\left(\bm{\Sigma}^{\star}\right)^{-2}+\left(M_{4}-3\right)\frac{1-p}{np}\mathsf{diag}\left\{ \bm{U}_{l,\cdot}^{\star\top}\bm{U}_{l,\cdot}^{\star}\right\} \nonumber \\
	& \quad+\frac{2\left(1-p\right)}{np}\bm{U}_{l,\cdot}^{\star\top}\bm{U}_{l,\cdot}^{\star}+\left(\bm{\Sigma}^{\star}\right)^{-2}\bm{U}^{\star\top}\mathsf{diag}\left\{ \left[d_{l,i}^{\star}\right]_{1\leq i\leq d}\right\} \bm{U}^{\star}(\bm{\Sigma}^{\star})^{-2}\label{eq:pca-heteropca-true-covariance-general}
\end{align}
where
\begin{align*}
	d_{l,i}^{\star} & \coloneqq\frac{1}{np^{2}}\left[\omega_{l}^{\star2}+\left(1-p\right)\left\Vert \bm{U}_{l,\cdot}^{\star}\bm{\Sigma}^{\star}\right\Vert _{2}^{2}\right]\left[\omega_{i}^{\star2}+\left(1-p\right)\left\Vert \bm{U}_{i,\cdot}^{\star}\bm{\Sigma}^{\star}\right\Vert _{2}^{2}\right]\\
	& \quad+\frac{2\left(1-p\right)^{2}}{np^{2}}S_{l,i}^{\star2}+\left(M_{4}-3\right)\sum_{s=1}^{r}\sigma_{s}^{\star4}U_{l,s}^{\star2}U_{i,s}^{\star2}.
\end{align*}
In addition, the distributional theory for the spiked covariance matrix
$\bm{S}^{\star}$ (cf.~Theorem \ref{thm:ce}) also holds with
\begin{align}
	v_{i,j}^{\star} & \coloneqq\frac{2-p}{np}S_{i,i}^{\star}S_{j,j}^{\star}+\frac{4-3p}{np}S_{i,j}^{\star2}+\frac{1}{np}\left(\omega_{i}^{\star2}S_{j,j}^{\star}+\omega_{j}^{\star2}S_{i,i}^{\star}\right)+\left(M_{4}-3\right)\frac{2\left(1-p\right)}{np}\sum_{s=1}^{r}\sigma_{s}^{\star4}U_{i,s}^{\star2}U_{j,s}^{\star2}\nonumber \\
	& \quad+\frac{1}{np^{2}}\sum_{k=1}^{d}\left\{ \left[\omega_{i}^{\star2}+\left(1-p\right)S_{i,i}^{\star}\right]\left[\omega_{k}^{\star2}+\left(1-p\right)S_{k,k}^{\star}\right]+2\left(1-p\right)^{2}S_{i,k}^{\star2}\right\} \left(\bm{U}_{k,\cdot}^{\star}\bm{U}_{j,\cdot}^{\star\top}\right)^{2}\nonumber \\
	& \quad+\frac{1}{np^{2}}\sum_{k=1}^{d}\left\{ \left[\omega_{j}^{\star2}+\left(1-p\right)S_{j,j}^{\star}\right]\left[\omega_{k}^{\star2}+\left(1-p\right)S_{k,k}^{\star}\right]+2\left(1-p\right)^{2}S_{j,k}^{\star2}\right\} \left(\bm{U}_{k,\cdot}^{\star}\bm{U}_{i,\cdot}^{\star\top}\right)^{2}\nonumber \\
	& \quad+\left(M_{4}-3\right)\frac{1}{np^{2}}\sum_{k=1}^{d}\sum_{s=1}^{r}\sigma_{s}^{\star4}U_{k,s}^{\star2}\Big[U_{i,s}^{\star2}\left(\bm{U}_{k,\cdot}^{\star}\bm{U}_{j,\cdot}^{\star\top}\right)^{2}+U_{j,s}^{\star2}\left(\bm{U}_{k,\cdot}^{\star}\bm{U}_{i,\cdot}^{\star\top}\right)^{2}\Big]\label{eq:ce-true-variance-i-j-general}
\end{align}
for $i\neq j$, and 
\begin{align}
	v_{i,i}^{\star} & \coloneqq\frac{12-9p}{np}S_{i,i}^{\star2}+\frac{4}{np}\omega_{i}^{\star2}S_{i,i}^{\star}+\left(M_{4}-3\right)\frac{4\left(1-p\right)}{np}\sum_{s=1}^{r}\sigma_{s}^{\star4}U_{i,s}^{\star4}\nonumber \\
	& \quad+\frac{4}{np^{2}}\sum_{k=1}^{d}\left\{ \left[\omega_{i}^{\star2}+\left(1-p\right)S_{i,i}^{\star}\right]\left[\omega_{k}^{\star2}+\left(1-p\right)S_{k,k}^{\star}\right]+2\left(1-p\right)^{2}S_{i,k}^{\star2}\right\} \left(\bm{U}_{k,\cdot}^{\star}\bm{U}_{i,\cdot}^{\star\top}\right)^{2}\nonumber \\
	& \quad+\left(M_{4}-3\right)\frac{4}{np^{2}}\sum_{k=1}^{d}\sum_{s=1}^{r}\sigma_{s}^{\star4}U_{k,s}^{\star2}U_{i,s}^{\star2}\left(\bm{U}_{k,\cdot}^{\star}\bm{U}_{i,\cdot}^{\star\top}\right)^{2}.\label{eq:ce-true-variance-i-general}
\end{align}
It is noteworthy that the fourth-moment $M_4$ plays an important role in the above variance calculation. 
One can then estimate $\bm{\Sigma}_{U,l}^{\star}$ (resp.~$v_{i,j}^{\star}$)
using the plug-in method similar to Algorithm \ref{alg:PCA-HeteroPCA-CR}
(resp.~Algorithm \ref{alg:CE-HeteroPCA-CI}), which naturally leads
to fine-grained confidence regions for $\bm{U}^{\star}$ and entrywise
confidence intervals for $\bm{S}^{\star}$. Note that the theoretical
guarantees in this paper are proved without exploiting the Gaussian
design (we only require basic sub-Gaussian properties and the information about the fourth moment), 
and hence the distributional theory and the validity
of confidence regions/intervals under sub-Gaussian design can be established
in an almost identical manner; for this reason, we omit the details for the sake of brevity.

\subsection{The necessity of incoherence condition} \label{subsec:necessity-incoherence}

Careful readers might wonder whether the incoherence condition (cf.~Assumption
\ref{assumption:incoherence}) is necessary for our algorithm designs
and theoretrial guarantees. 
As has been pointed out in the low-rank matrix estimation literature, 
the incoherence condition plays a crucial role in enabling reliable estimation in the presence of missing data \citep{CanTao10,Se2010Noisy,chi2018nonconvex}. 
To make it more precise, we state below a theorem that captures the fundamental relation between the incoherence
parameter $\mu$ and the information-theoretic sampling limit.

\begin{theorem} \label{thm:sample-complexity}Consider $r=1$ and
	any incoherence parameter $\mu\ll d$. Suppose that $n\geq d$, and
	$p<(1-\varepsilon)\sqrt{\frac{\mu}{nd}}$ for an arbitrarily small constant $0<\varepsilon<1/4$.
	Generate $\Omega$ according to the random sampling model. Then with
	probability at least $0.9$, there exist $(\bm{u}_{1}^{\star},\bm{f}_{1})$
	and $(\bm{u}_{2}^{\star},\bm{f}_{2})$ that satisfy, for both $i=1,2$,
	\begin{enumerate}
		\item $\bm{u}_{i}^{\star}\in\mathbb{R}^{d}$ is a fixed unit vector and
		is $\mu$-incoherent: $\Vert\bm{u}_{i}^{\star}\Vert_{\infty}\leq\sqrt{\mu/d}$,
		\item $\bm{f}_{i}\in\mathcal{N}(\bm{0},\bm{I}_{n})$ is an isotropic Gaussian
		random vector,
	\end{enumerate}
	such that $\min\left\Vert \bm{u}_{1}^{\star}\pm\bm{u}_{2}^{\star}\right\Vert _{2}\asymp1$ 
	but $\mathcal{P}_{\Omega}(\bm{u}_{1}^{\star}\bm{f}_{1}^{\top})=\mathcal{P}_{\Omega}(\bm{u}_{2}^{\star}\bm{f}_{2}^{\top})$.
\end{theorem}

This theorem shows that even in the simplest noiseless setting with $r=1$: 
if $p<(1-\varepsilon)\sqrt{\frac{\mu}{nd}}$,
then with high probability one cannot distinguish $\bm{u}_{1}^{\star}$ and $\bm{u}_{2}^{\star}$ 
given only observations $\mathcal{P}_{\Omega}(\bm{u}_{1}^{\star}\bm{f}_{1}^{\top})$. Consequently, 
in order to ensure identifibility, one needs to require the total sample size
$ndp^{2}$ to exceed the order of $\mu$. 
In other words, this theorem demonstrates how the
incoherence parameter $\mu$ dictates the difficulty of the problem:
when the eigenvectors are spiky, we require more samples to reliably
estimate the principal subspace. 
Additionally, it is worth pointing out that the dependency
w.r.t.~$\mu$ in our theory is very likely to be suboptimal (e.g., Theorem \ref{thm:pca-complete}
requires the sample size to exceed the order of $\mu^{4}$). This
might be improvable via more refined analysis, and we leave it to
future investigation.

\subsubsection{Proof of Theorem \ref{thm:sample-complexity}}

Without loss of generality, assume that $d/\mu$ is an integer. Let
\[
\left(\bm{u}_{1}^{\star}\right)_{i}=\begin{cases}
	\sqrt{\frac{\mu}{d}}, & \text{if }1\leq i\leq\frac{d}{\mu},\\
	0, & \text{otherwise}
\end{cases}\qquad\text{and}\qquad\bm{f}_{1}\sim\mathcal{N}\left(\bm{0},\bm{I}_{n}\right).
\]
Consider a random bipartite graph $\mathcal{G}$ generated by taking
two disjoint vertex sets $\mathcal{U}$ and $\mathcal{V}$ with $|\mathcal{U}|=d/\mu$
and $|\mathcal{V}|=n$, and connecting each $u\in\mathcal{U}$ and
$v\in\mathcal{V}$ independently with probability $p$. One can check
that there is an equivalence between the edge set of $\mathcal{G}$
and a subset $\Omega'$ of the subsampled index set $\Omega\subseteq[d]\times[n]$
defined as
\[
\Omega'\coloneqq\left\{ \left(i,j\right):1\leq i\leq\frac{d}{\mu},\left(i,j\right)\in\Omega\right\} .
\]
More precisely, if one connects the $i$-th vertex of $\mathcal{U}$
and the $j$-th vertex of $\mathcal{V}$ if and only if $(i,j)\in\Omega'$,
then the resulting random graph has the same distribution as $\mathcal{G}$.
In view of \citet[Theorem 6]{johansson2012giant}, when $p<(1-\varepsilon)\sqrt{\frac{\mu}{nd}}$,
for any $k>0$ the probability that there is no connected component
in $\mathcal{G}$ with at least $k$ vertices in $\mathcal{U}$ is
at least
\[
1-\frac{d}{\mu\left(1-\varepsilon\right)}\exp\left(-\frac{1}{6}k\varepsilon^{3}\right). 
\]
By taking $k=12\varepsilon^{-3}\log(d/\mu)$ and recall that $d\gg\mu$
and $\varepsilon<1/4$, we can show that that with probability exceeding
$0.9$, there exists no connected component in $\mathcal{G}$ with at
least $12\varepsilon^{-3}\log(d/\mu)$ vertices in $\mathcal{U}$.
Denote by $\mathcal{C}_{1},\ldots,\mathcal{C}_{K}$ the collection
of connected components in $\mathcal{G}$, and let $\mathcal{U}_{i}$
(resp.~$\mathcal{V}_{i}$) be the set of vertices in $\mathcal{U}$
(resp.~$\mathcal{V}$) that reside in $\mathcal{C}_{i}$. When $d\gg\mu$,
we can always find a subset of $\mathcal{I}\subset[K]$ such that
\[
\frac{1}{2}\left(\frac{d}{\mu}-12\varepsilon^{-3}\log\frac{d}{\mu}\right)\leq\sum_{i\in\mathcal{I}}\left|\,\mathcal{U}_{j}\right|\leq\frac{1}{2}\left(\frac{d}{\mu}+12\varepsilon^{-3}\log\frac{d}{\mu}\right).
\]
Then we can proceed to construct $\bm{u}_{2}^{\star}$ and $\bm{f}_{2}$ in the following manner:
\[
\left(\bm{u}_{2}^{\star}\right)_{i}=\begin{cases}
	\left(\bm{u}_{1}^{\star}\right)_{i}, & \text{if }1\leq i\leq\frac{d}{\mu}\text{ and }i\in\mathcal{U}_{k}\text{ for some }k\in\mathcal{I}\\
	-\left(\bm{u}_{2}^{\star}\right)_{i}, & \text{if }1\leq i\leq\frac{d}{\mu}\text{ and }i\in\mathcal{U}_{k}\text{ for some }k\notin\mathcal{I}\\
	0 & \text{otherwise}
\end{cases}
\]
for any $1\leq i\leq d$, and
\[
\left(\bm{f}_{2}\right)_{j}=\begin{cases}
	\left(\bm{f}_{1}\right)_{j}, & \text{if }j\in\mathcal{V}_{k}\text{ for some }k\in\mathcal{I}\\
	-\left(\bm{f}_{1}\right)_{j}, & \text{if }j\in\mathcal{V}_{k}\text{ for some }k\notin\mathcal{I}
\end{cases}
\]
for any $1\leq j\leq n$. It is straightforward to check that when
$d\gg\mu$, both $\bm{u}_{1}^{\star}$ and $\bm{u}_{2}^{\star}$ are
$\mu$-incoherent obeying $\min\left\Vert \bm{u}_{1}^{\star}\pm\bm{u}_{2}^{\star}\right\Vert _{2}\asymp1$,
and $\bm{f}_{2}\sim\mathcal{N}(\bm{0},\bm{I}_{n})$. This completes
the proof.

\subsection{Other observational models} \label{subsec:general-observation-model}

This subsection discusss the possibility of extending the algorithms
and theory in the current paper to accommodate another model stated below, which finds applications
in medical research and wireless communication. Consider the setting
where, instead of missing all entries outside the sampling set
$\Omega$, we observe pure noise without knowing $\Omega$ in advance; that is, we observe
\[
y_{l,j}=\begin{cases}
	x_{l,j}+\eta_{l,j}, & \text{for}\text{ all }\left(l,j\right)\in\Omega,\\
	\eta_{l,j}, & \text{otherwise},
\end{cases}
\]
or equivalently, $\bm{Y}=\mathcal{P}_{\Omega}(\bm{X})+\bm{N}$ in
the matrix form. Similar to (\ref{eq:PCA-E-gram-decompose}) we can compute
\[
\frac{1}{p^{2}}\mathbb{E}\left[\bm{Y}\bm{Y}^{\top}\,\big|\,\bm{X}\right]=\bm{X}\bm{X}^{\top}+\left(\frac{1}{p}-1\right)\mathcal{P}_{\mathsf{diag}}\left(\bm{X}\bm{X}^{\top}\right)+\frac{n}{p^2}\mathsf{diag}\left\{ \left[\omega_{l}^{\star2}\right]_{1\leq l\leq d}\right\} .
\]
Therefore, it is reasonable to employ \textsf{HeteroPCA} to attempt estimation of 
the principal subspace $\bm{U}^{\star}$ and the spiked covariance
matrix $\bm{S}^{\star}$ under this setting, although we might need to impose 
stronger noise conditions due to the presence of more noise.

Note, however, that it would 
be difficult to perform valid statistical inference (as we did in
Algorithm \ref{alg:PCA-HeteroPCA-CR} and \ref{alg:CE-HeteroPCA-CI}) for this setting; 
the reason is that it becomes fairly difficult to estimate the sampling rate $p$, 
a parameter that 
is required for computing an estimate for the covariance matrix $\bm{\Sigma}_{U,l}^{\star}$
and the variance $v_{i,j}^{\star}$. To further elucidate this point, consider
the setting where the noise components are i.i.d.~$\mathcal{N}(0,\sigma^{2})$
with known $\sigma>0$. Then this observation model can be described
as follows: for each $(i,j)$, 
\[
y_{l,j}\sim\begin{cases}
	\,\mathcal{N}\left(x_{l,j},\sigma^{2}\right), & \text{with probability }p,\qquad\quad\;\;\text{(Case 1)}\\
	\,\mathcal{N}\left(0,\sigma^{2}\right), & \text{with probability }1-p.\qquad\text{(Case 2)}
\end{cases}
\]
The task of estimating $p$ boils down to estimating the portion
of the entries corresponding to Case 1, which, however, is impossible in general.
For instance, suppose that half of the rows of $\bm{U}^{\star}$ are zero.
When $\bm{U}_{l,\cdot}^{\star}=\bm{0}$, we know that $x_{l,j}=0$
for all $j\in[n]$. Therefore for roughly half of the entries, one
cannot distinguish whether they belong to Case 1 or Case 2. This means
that it is in general impossible to estimate $p$ in an accurate manner. 

Nevertheless, 
when stronger noise conditions are met, we shall be able to estimate the sampling rate $p$ reliably using the following scheme. 
Employing similar analysis as in Appendix \ref{appendix:proof-pca-noise-level-est}, we can show that the sum of squares of all observations $y_{l,j}$ concentrates around
\[
\sum_{l=1}^d \sum_{j=1}^n y_{l,j}^2 \approx n p \left\Vert \bm{U}^{\star}\bm{\Sigma}^{\star}\right\Vert _{\mathrm{F}}^{2} + n\sum_{l=1}^d \omega_l^{\star2}
\]
with high probability. 
When the noise level $\omega_{\max}$ is sufficiently small, namely $\omega_{\max} \ll \sqrt{p/d} \sigma_r^\star$, the second term on the right hand side of the above equation becomes negligible, and hence we can faithfully estimate $p$ by means of the following data-driven plug-in estimate: 
\[
\widehat{p} = \frac{1}{n \Vert \bm{U} \bm{\Sigma} \Vert_{\mathrm{F}}^2} \sum_{l=1}^d \sum_{j=1}^n y_{l,j}^2.
\]
This in turn allows us to perform valid statistical inference for this setting, through the same analysis framework introduced in this paper. For example, we can still apply Theorem~\ref{thm:hpca_inference_general} to establish the first- and second-order approximations \eqref{eq:decomposition-U} and \eqref{eq:decomposition-S} with a different effective noise matrix $\bm{E} = p^{-1} \mathcal{P}_{\Omega}(\bm{X}) - \bm{X} + p^{-1} \bm{N}$, then the distributional characterization and data-driven construction of confidence regions/intervals follow naturally from the same analysis. 
Given that this model is not the focus of the current paper, we omit the details for the sake of brevity. 

\subsection{Heteroskedastic noise across rows}

In this paper, we allow the noise levels to vary across different rows, 
but require the noise variance to be identical within each row. Naturally, one might ask whether it is feasible to
further relax such assumptions and accommodate fully heterogeneous noise across all entries. 
Unfortunately, while it is plausible to establish similar estimation guarantees and distributional theory (characterizing the distribution of estimation errors via model parameters) for this extension, 
the design of data-driven inferential procedures --- such as constructing confidence regions for $\bm{U}^{\star}$ or entrywise confidence intervals for $\bm{S}^{\star}$ --- would be substantially more challenging.  The rationale is as follows.
\begin{itemize}
	\item We first notice that,  the setting of Theorem \ref{thm:hpca_inference_general} 
	(i.e., our general second-order subspace perturbation theory for \textsf{HeteroPCA}) 
	allows heteroskedastic random noise across entries (see Assumption
	\ref{assumption:subspace-noise}). This suggests that in the PCA
	context, even if we assume heteroskedastic noise across entries, we
	can still invoke Theorem~\ref{thm:hpca_inference_general} to obtain
	our key error decomposition \eqref{eq:decomposition}, which in turn allows
	us to establish estimation guarantees (cf.~Lemmas \ref{lemma:pca-1st-err}
	and \ref{lemma:pca-noise-level-est}) and distributional theory (cf.~Theorems
	\ref{thm:pca} and \ref{thm:ce}).
	
	\item When it comes to statistical inference, however, we need to reliably estimate
	the covariance matrix $\bm{\Sigma}_{U,l}^{\star}$ and $v_{i,j}^{\star}$
	in a data-driven manner. The key difficulty is that, under our assumptions,
	we cannot hope to provide an entrywise consistent estimate for the
	ground truth data matrix $\bm{X}\in\mathbb{R}^{d\times n}$. We can
	only estimate the underlying spiked covariance matrix $\bm{S}^{\star}\in\mathbb{R}^{d\times d}$
	and its associated eigenspace $\bm{U}^{\star}$ and eigenvalues $\bm{\Lambda}^{\star}$
	reliably. If the noise is heteroskedastic across entries, say, $\mathsf{var}(\eta_{l,j})=\omega_{l,j}^{\star2}$,
	then it is in general impossible to faithfully estimate each noise variance
	$\omega_{l,j}^{\star2}$ or their (nonlinear) functions. This makes
	it hard to construct confidence regions for $\bm{U}^{\star}$
	or entrywise confidence intervals for $\bm{S}^{\star}$ in a data-driven
	manner (as done in Algorithms \ref{alg:PCA-HeteroPCA-CR} and \ref{alg:CE-HeteroPCA-CI}).
\end{itemize}

It is worth noting that several previous studies on PCA \citep{zhu2019high} and random matrix theory \citep{latala2018dimension,alt2021spectral} permit heterogeneous noise across entries. In the realm of statistical inference, however, 
the studies of heteroskedastic noise remained highly limited. For example, \citet{bao2018singular,xia2019normal,koltchinskii2020efficient} requires the noise components to be either i.i.d.~Gaussian or have identical first fourth moments. To the best of our knowledge, tolerating noise that is heterogeneous across rows, as we have in this current study, represents the mildest assumption within the framework of statistical inference for PCA.

\bibliographystyle{apalike}
\bibliography{bibfileNonconvex}

\end{document}